\documentclass[a4paper,11pt]{amsart}
\usepackage{amssymb,amsmath,amsthm,stmaryrd,bm,enumitem,hyperref,mathrsfs,xcolor} 
\usepackage{microtype}

\newcommand\beqn{\begin{equation}}
\newcommand\eeqn{\end{equation}}
\newcommand\beqny{\begin{eqnarray}}
\newcommand\eeqny{\end{eqnarray}}
\newcommand\beqnyn{\begin{eqnarray*}}
\newcommand\eeqnyn{\end{eqnarray*}}

\usepackage{amsmath,amsfonts,amssymb,amsthm}
\newtheorem{theorem}{Theorem}[section]
\newtheorem{lemma}[theorem]{Lemma}
\newtheorem{proposition}[theorem]{Proposition}
\newtheorem{corollary}[theorem]{Corollary}
\newtheorem{definition}[theorem]{Definition}
\newtheorem{remark}[theorem]{Remark} 
\newtheorem{hypothesis}[theorem]{Hypothesis}

\allowdisplaybreaks[3]
\numberwithin{equation}{section}

\newcommand{\op}[1]{\operatorname{\text{\rm #1}}}   
\newcommand{\bs}[1]{\boldsymbol{#1}}

\title[Analysis of singularities of area-minimizing currents, Part~IV]{\textls[-25]{Analysis of singularities of area-minimizing  currents, part~IV: branch points of planar frequency $\geq 2,$ asymptotic normal forms, 
local topology and singular-set structure}} 
\author{Brian Krummel \& Neshan Wickramasekera}

\begin{document}

\begin{abstract}
 This is Part~IV of a five-part series developing a  framework for 
 the local structure of $n$-dimensional area-minimizing rectifiable currents $T$ of codimension $ \geq 2$. 
Part~III analyzed branch points $Z$
of planar frequency ${\mathcal N}_{T, {\rm Pl}}(Z) \neq 2.$ 
This paper treats the 
case ${\mathcal N}_{T, {\rm Pl}}(Z) = 2$, completes the analysis of the singular-set structure,  refines the asymptotics of Part III 
at integer planar frequencies, and extends its topological $n$-disk criterion to integer planar frequencies.
This  completes the Euclidean theory of the series; Part~V 
extends the framework 
to Riemannian ambient spaces. 

Together, Part~III and this paper  establish the program's main result for branch points: 
for any integer $q \geq 2$ and ${\mathcal H}^{n-2}$-a.e.\ density-$q$ branch point $Z$,   
there exist a \emph{branching order} ${\mathcal O}_{T}(Z)$ and an \emph{asymptotic normal form} at $Z$ with quantitative remainder decay. 
The branching order is rational and satisfies 
${\mathcal O}_{T}(Z) \geq {\mathcal N}_{T, {\rm Pl}}(Z)$ with equality 
when ${\mathcal N}_{T, {\rm Pl}}(Z) \not\in {\mathbb N}$. The asymptotic normal form is determined by a unique tangent plane $P_{Z}$ to $T$ at $Z$ and  a unique, nonzero, average-free, 
cylindrical, 
${\mathcal O}_T(Z)$-homogeneous, 
$q$-valued tangent function 
$\varphi^{(Z)} \, : \, P_{Z} \to {\mathcal A}_{q}(P_{Z}^{\perp})$ obtained relative to a center manifold through $Z.$ The center manifold approximates the average deviation of $T$ from $Z + P_{Z}$, and can be taken to be $Z +P_{Z}$ if ${\mathcal N}_{T, {\rm Pl}}(Z) \not\in {\mathbb N}$. 
The branching order 
and the tangent function are intrinsic to $Z$ and independent of 
the center manifold, while  
the center manifold is canonical to  order 
${\mathcal O}_{T}(Z)$ in an $L^{2}$ sense. 
Moreover, the tangent function is a finite superposition of irreducible homogeneous components, and  
there exist coordinates 
$(z, y) \in {\mathbb C} \times {\mathbb R}^{n-2}$ on $P_{Z}$ such that 
each of its nonzero components is  independent of $y$ and has the form 
${\rm Re} \, \left(c^{(Z)} z^{{\mathcal O}_{T}(Z)}\right)$ for some constant $c^{(Z)}\in {\mathbb C}^{n+m}\setminus \{0\}$ with $c^{(Z)} \cdot c^{(Z)} = 0$ and  ${\rm Re} \, (c^{(Z)}),{\rm Im} \, (c^{(Z)}) \in P_{Z}^{\perp};$ in particular, the cross section of the (cylindrical) graph of each component is a holomorphic curve
in ${\mathbb C} \times {\rm span}_{\mathbb R} \,\{{\rm Re} \, (c^{(Z)}),{\rm Im} \, (c^{(Z)})\}$.

A main corollary of the 
normal form theorem  
is a branching-order criterion governing the local topology and regularity of ${\rm spt} \, T$:
if a density-$q$ branch point $Z$ has ${\mathcal O}_{T}(Z) = p/q$ for $p, q$ coprime, and if $Z$ is \emph{$(n-2)$-strongly non-isolated}, i.e.\ at a small scale $\rho>0$ about $Z$ 
density-$q$ branch points with branching order $\geq {\mathcal O}_{T}(Z)$ project fully to cover a $\rho/2$-ball about $Z$ in some 
$(n-2)$-plane through $Z$,
then $Z$ is a \emph{classical branch point}:
near $Z$, 
${\rm spt} \, T$ 
is homeomorphic to an $n$-disk, admits a $C^{1,\mu}$ parameterization, and the singular set of $T$ is an $(n-2)$-dimensional embedded $C^{1, \mu}$ submanifold.  
Algebraic examples show that 
${\rm spt} \, T$ need not be  
an $n$-disk
if the $(n-2)$-strong non-isolatedness of $Z$ is relaxed to requiring full projection of all density $q$
branch points. 
The normal form theorem also leads to our main structure theorem for the singular set of $T$: it is locally the finite union of 
disjoint, locally compact, locally $(n-2)$-rectifiable sets (each having locally finite 
${\mathcal H}^{n-2}$-measure).

Methodologically, for planar frequency~$\neq 2$, Parts~I--III avoid auxiliary center-manifolds  underpinning  
the classical Almgren framework (\cite{Almgren}) and instead use  intrinsic arguments  
that take tangent planes as reference.   
Handling planar frequency $<2$ this way is key to our structural results. 
At  planar frequency~2, 
a center manifold becomes  
necessary, and it has additional 
structure: for any $Z$ with planar frequency 
$\geq 2,$ there exists a single \emph{effective center manifold} supporting a multi-valued normal map $N$ whose critical nodal set contains every point $Z^{\prime}$ sufficiently close to $Z$ that has  planar frequency $\geq 2$ and the same density as $Z.$
The graph of $N$ approximates $T$ 
with quantitatively decaying error  
under rescaling about any such $Z^{\prime}$. This all-scale control on $N$ essentially reduces the analysis of $T$ 
to that of $N.$ We carry out this analysis by adapting 
our  work \cite{KrumWic2} on Dirichlet energy-minimizing multi-valued maps to $N$'s variational structure.  
Central to this is a new   
$L^{\infty}\!-\!L^{2}$ estimate 
that plays the role of Almgren’s  
Hölder 
continuity 
estimate in \cite{KrumWic2}. 
The singular-set structure for $T$ follows 
from the planar-frequency $2$ analysis and Parts~II-III, while the  asymptotic 
normal forms 
require refining Part~III at integer planar frequencies $\geq 3.$ 
This refinement proceeds exactly as our planar frequency 2 analysis. 
The overarching principle of our program is to choose analytic mechanisms according to 
the geometry 
of singularities 
as determined by planar frequency,
rather than 
relying on center manifolds for all branch points.
\end{abstract}

\maketitle

\tableofcontents

\section{Introduction}  
This is Part~IV of a five-part series of papers developing a new framework for studying the interior local structure
of an $n$-dimensional locally area-minimizing rectifiable current $T$ in an $(n+m)$-dimensional Riemannian manifold, $m \geq 2$.  In this paper, as in Parts~I--III (\cite{KrumWica}, \cite{KrumWicb},  \cite{KrumWicc}), 
we shall assume that the ambient manifold is an open subset of ${\mathbb R}^{n+m}$. The present paper completes our program for Euclidean ambient spaces. 
Part~V (\cite{KrumWice}) provides the  technical modifications needed to extend the main results to general Riemannian ambient spaces.

We refer the reader to the introductions 
to Part~I (\cite{KrumWica})  and Part~III (\cite{KrumWicc}) for a comprehensive description of our overall program, including the statements of its main results. In particular, the introduction to Part~III contains the statements of the main theorems, and a detailed overview of our analysis leading to them, concerning   \emph{branch points} (i.e.\ singular points where there is a tangent cone supported on a plane). 

The first part of the branch-point analysis is carried out in Part~III, and this analysis is completed in the present paper. Taken together, 
Part~III and the  present paper yield three main outcomes which we describe below, referring the reader to \cite{KrumWicc} for the precise theorem statements: 
\begin{itemize}
\item[(1)] \cite[Theorem~B]{KrumWicc} establishes that, for  ${\mathcal H}^{n-2}$-a.e.\ branch point $Z,$ there is an intrinsic rational invariant ${\mathcal O}_{T}(Z) >1$ called the \emph{branching order of $T$ at $Z$}, and that at $Z$ the current admits an \emph{asymptotic normal form} of order ${\mathcal O}_{T}(Z)$. If the density of $T$ at $Z$ is $q$ (an integer $\geq 2$), the branching order takes the form ${\mathcal O}_{T}(Z) = \ell_{Z}/q_{Z}$ for coprime integers $\ell_{Z}, q_{Z}$ with $q_{Z} \leq q$ and $\ell_{Z} > q_{Z}$. Moreover, ${\mathcal O}_{T}(Z) \geq {\mathcal N}_{T, {\rm Pl}}(Z)$, with equality if ${\mathcal N}_{T, {\rm Pl}}(Z) \not\in {\mathbb N},$ where ${\mathcal N}_{T, {\rm Pl}}(Z)$ is the planar frequency of $T$ at $Z$ (introduced in Part~I and discussed below).  The asymptotic normal form is determined by: 
\begin{itemize}
\item a unique tangent plane $P_{Z}$ to $T$ at $Z$, and  
\item a unique, $q$-valued, nonzero, ${\mathcal O}_{T}(Z)$-homogeneous, locally Dirichlet energy minimizing, cylindrical, average-free tangent function $\varphi^{(Z)} \, : \, P_{Z} \to {\mathcal A}_{q}(P_{Z}^{\perp})$ to $T$ at $Z.$ 
\end{itemize}
The tangent function $\varphi^{(Z)}$ arises as the blow up of $T$ relative to a  \emph{center manifold} ${\mathcal M}_{Z}$ through $Z$. Such ${\mathcal M}_{Z}$ exists and is given by ${\mathcal M}_{Z} = Z + {\rm graph} \, \psi^{(Z)}$, where 
$\psi^{(Z)} \, : \, P_{Z} \to P_{Z}^{\perp}$ 
is a $C^{3, \kappa}$ function for some $\kappa = \kappa(n, m, q) \in (0, 1),$ with $\psi^{(Z)}(0) = 0$, and $D\psi^{(Z)}(0) = 0$;
the function $\psi^{(Z)}$ is an approximation to the average deviation of $T$ from $Z + P_{Z}$ and may be taken to be the zero-function if ${\mathcal N}_{T, {\rm PL}}(Z) \notin {\mathbb N}$.

Modulo the center manifold, $T$ is quantitatively asymptotic to this normal form; that is, 
the average $L^{2}$ distance of ${\rm spt} \, T$ from the graph of 
$\psi^{(Z)} + \varphi^{(Z)}$ at scale $\rho$ is $O(\rho^{{\mathcal O}_{T}(Z) + \mu_{Z}})$ as $\rho \to 0$, where $\mu_{Z} >0$ is a locally uniform constant.

The branching order ${\mathcal O}_{T}(Z)$ and the tangent function $\varphi^{(Z)}$ are  uniquely determined by $Z$ independently of the choice of the center manifold, and the center manifold is canonical to the order given by 
${\mathcal O}_{T}(Z)$, in the following sense: if $\psi_{1}^{(Z)}, \psi_{2}^{(Z)} \, : \, P_{Z} \to P_{Z}^{\perp}$ define two center manifolds ${\mathcal M}_{Z}^{(1)}, {\mathcal M}_{Z}^{(2)}$ at $Z$ such that $T$ is asymptotic, in the above sense, to the normal forms  determined by the respective tangent functions 
$\varphi^{(Z)}_{1}, \varphi^{(Z)}_{2}$ relative to ${\mathcal M}_{Z}^{(1)}, {\mathcal M}_{Z}^{(2)}$, and if the  homogeneity degrees of these tangent functions are ${\mathcal O}_{T}^{(1)}(Z)$, 
${\mathcal O}_{T}^{(2)}(Z),$  then $\varphi^{(Z)}_{1} = \varphi^{(Z)}_{2}$ 
(in particular ${\mathcal O}_{T}^{(1)}(Z) = 
{\mathcal O}_{T}^{(2)}(Z)$) and, denoting the common homogeneity degree by ${\mathcal O}_{T}(Z)$,  
$\rho^{-n/2}\|\psi_{Z}^{(1)} - \psi_{Z}^{(2)}\|_{L^{2}(B_{\rho}(0))} = O(\rho^{{\mathcal O}_{T}(Z) + \epsilon_{Z}})$
as $\rho \to 0$.

Furthermore, there exist coordinates $(z, y) \in {\mathbb C} \times {\mathbb R}^{n-2}$ relative to an orthonormal basis for $P_{Z}$ such that the tangent function $\varphi^{(Z)}$ is independent of $y$ and is the superposition   
$$\varphi^{(Z)}(z, y) = \sum_{j=1}^{N_{Z}}  m_{j}^{(Z)}\varphi_{j}^{(Z)}(z)$$ 
of a finite number $N_{Z}$ of distinct irreducible components 
$\varphi_{j}^{(Z)}$ with positive integer multiplicities   
$m_{j}^{(Z)},$ for $j=1, 2, \ldots, N_{Z}$. Each non-zero 
component $\varphi_{j}^{(Z)}$ is a $q_{Z}$-valued homogeneous function  
$\varphi^{(Z)}_{j} \, : \, P_{Z} \to {\mathcal A}_{q_{Z}}(P_{Z}^{\perp})$ given by 
$$\varphi_{j}^{(Z)}(z) = {\rm Re} \, \left(c_{j}^{(Z)} z^{{\mathcal O}_{T}(Z)}\right)$$ 
for some constant vector $c^{(Z)}_{j} \in {\mathbb C}^{n+m} \setminus \{0\}$ with  $c^{(Z)}_{j} \cdot c^{(Z)}_{j} = 0$ and ${\rm Re} \, (c^{(Z)}_{j}), {\rm Im} \, (c^{(Z)}_{j}) \in P_{Z}^{\perp}$; in particular, the graph of $\varphi_{j}^{(Z)}$ is a holomorphic curve. (Indeed, the image of 
$\varphi_{j}^{(Z)}$ is contained in $\Sigma \equiv {\rm span}_{\mathbb R} \, \{{\rm Re} \, (c^{(Z)}_{j}), {\rm Im} \, (c^{(Z)}_{j})\}$ and the condition $c^{(Z)}_{j} \cdot c^{(Z)}_{j} = 0$ implies that 
$\Sigma$ is a two-dimensional subspace  of $P_{Z}^{\perp}$  and that there is a complex structure on $\Sigma,$ which is compatible with the Euclidean metric, with respect to which 
$\{(z, \zeta): \zeta \in {\rm spt} \,\varphi_{j}^{(Z)}(z),\, z \in {\mathbb C}\}$ is a holomorphic curve in ${\mathbb C} \times \Sigma$.)

\item[(2)] \cite[Theorem~A]{KrumWicc} establishes a link between  branching-order behaviour and the local topology and regularity of ${\rm spt} \, T$. Specifically, any density-$q$ branch point $Z$ which has branching order  ${\mathcal O}_{T}(Z) = p/q$ for an integer $p >q$ coprime to $q$ and which 
is \emph {$(n-2)$-strongly non-isolated} is a  \emph{classical branch point}; that is, $Z$ has a neighborhood in which ${\rm spt} \, T$ is homeomorphic  to an $n$-disk and admits a $C^{1, \mu}$ parameterization, while the singular set of $T$ in this neighborhood is an $(n-2)$-dimensional $C^{1, \mu}$ submanifold consisting entirely of density $q$ branch points, all having the same branching order ${\mathcal O}_{T}(Z)$. 

The  $(n-2)$-strong non-isolatedness of $Z$ is the requirement that for some sufficiently small $\rho>0$ and some $(n-2)$-dimensional linear subspace $L$ of ${\mathbb R}^{n+m}$, 
the set of density-$q$ branch points $Z^{\prime} \in {\mathbf B}_{\rho}(Z)$  with ${\mathcal O}_{T}(Z^{\prime}) \geq {\mathcal O}_{T}(Z)$ 
``projects fully to $Z + L$ near $Z$,'' in the sense that its orthogonal projection to $Z + L$  contains the ball ${\mathbf B}_{\rho/2}(Z) \cap (Z + L)$.

Explicit classical algebraic surfaces in ${\mathbb C}^{3}$ demonstrate that  if $Z$ fails to be $(n-2)$-strongly non-isolated (for $n=4)$, then ${\rm spt} \, T$ need not be locally homeomorphic to an $n$-disk near $Z$. 
One such example is the complex surface
$$\{(w, x, y) \, : \, w^{2} = x^{3}(x^{2} - y^{3})\} \subset {\mathbb C}^{3};$$  
for $Z = 0$ and $q=2$, this example shows in fact that the requirement in the theorem that $Z$ be $(n-2)$-strongly  non-isolated  
cannot be relaxed to require only that the set of all density-$q$ branch points has the full projection property near $Z$. The stronger condition in the theorem requires this projection property specifically  for those density-$q$ branch points whose  branching order is at least ${\mathcal  O}_{T}(Z)$.

\item[(3)] \cite[Theorem~C]{KrumWicc} establishes that, in any compact subset of the ambient space, the singular set of $T$ is the union of finitely many pairwise disjoint, 
locally compact, locally $(n-2)$-rectifiable sets, each of locally finite ${\mathcal H}^{n-2}$ measure.
\end{itemize}

The topological $n$-disk theorem in (2) is a consequence of our proof  establishing the existence of asymptotic normal forms as in (1). The 
singular-set structure in (3) is a consequence of the asymptotic normal forms together with the results of Part~II (which establish rectifiability conclusions for the set of non-branch point singularities). Thus, asymptotic normal forms lie at the heart of our branch-point analysis.

The methods we employ to establish asymptotic normal forms are dictated by the geometric behavior of $T$ on approach to branch points $Z$. The primary analytic invariant governing this geometry is the \emph{planar frequency} ${\mathcal N}_{T, {\rm Pl}}(Z)$ of $T$ at $Z$. 
The value of  the planar frequency at $Z$, when it exists, determines the order of the current's height deviation from the tangent plane at $Z,$ and whether on approach to $Z$ the maximum separation of its  ``sheets'' from one-another  may be of lower order. The latter can happen only if the planar frequency at $Z$ is an integer. This invariant thus determines the analytic mechanism we use to obtain the asymptotic normal form at $Z$. 

The planar frequency function was introduced in Part~I (\cite{KrumWica}). It is an intrinsically defined frequency function for $T$ that encodes the current's rate-of-decay information when it decays at a point towards a given plane. We  established in \cite[Theorem~4.4]{KrumWica} that if for a point $Z \in {\rm spt} \, T$ there exists a plane $P$ such that over   scales $\rho$ in an interval 
(not necessarily having $0$ as a limit point), 
$T$ remains close to $P$ at $Z$ and decays towards $P$ at a rate $o(\rho^{1 + \alpha})$ for some $\alpha >0$, then the planar frequency function $\rho \mapsto N_{T, P, Z}(\rho)$ associated with $Z$ and $P$ is approximately monotone (non-decreasing) on that interval. Using this property, we further established in \cite{KrumWica} and \cite{KrumWicb} that for each integer $q \geq 2,$ there exists a constant $\alpha = \alpha(n, m, q) \in (0, 1)$ such that at ${\mathcal H}^{n-2}$-a.e.\ density-$q$ branch point $Z$, the current does   decay towards a unique tangent plane $P_{Z}$  at a rate $o(\rho^{1 + \alpha})$ as the scale $\rho \to 0.$ Consequently, at such branch points $Z$---which we term \emph{rapid-decay branch points}---the planar frequency function with respect to $P_{Z}$ is approximately monotone down to arbitrarily small scales $\rho.$ It therefore has a limit  ${\mathcal N}_{T, {\rm Pl}}(Z) \equiv \lim_{\rho \to 0} \, N_{T, P_{Z}, Z}(\rho)$ which we call the planar frequency at $Z,$ which satisfies ${\mathcal N}_{T,{\rm Pl}}(Z) \geq 1 + \alpha.$ Moreover, our analysis in \cite{KrumWica} provides a fundamental locally uniform estimate giving  decay of $T$ at such points $Z$ towards the unique tangent plane $Z + P_{Z}$.  See \cite[Theorem~X]{KrumWicc}, which summarizes these results as well as our conclusions (obtained in \cite{KrumWica} and  \cite{KrumWicb}) for $T$ away from rapid-decay branch points. 

Armed with this decay estimate and the validity of the monotonicity formula for planar frequency at 
${\mathcal H}^{n-2}$-a.e.\ branch point, in Part~III we use intrinsic geometric arguments to establish asymptotic normal forms at ${\mathcal H}^{n-2}$-a.e.\ branch point $Z$ with ${\mathcal N}_{T, {\rm Pl}}(Z) \not\in {\mathbb N}.$  This leads to \cite[Theorem~1.8]{KrumWicc} and  \cite[Theorem~1.12]{KrumWicc}, which are the versions, for points of non-integer planar frequency, 
of the two main theorems in (1) and (2) above, with ${\mathcal O}_{T}(Z) = {\mathcal N}_{T, {\rm Pl}}(Z)$. 
In contrast to the classical Almgren framework (\cite{Almgren}) for the analysis of branch points, our 
analysis in Parts~I--III avoids  auxiliary \emph{center manifolds}, the use of which for all branch points is a technically demanding, foundational aspect of 
the classical approach (discussed in more detail in the next subsection). 

More generally, our planar frequency based intrinsic arguments in Part~III extend to ${\mathcal H}^{n-2}$-a.e.\ point in the set ${\mathcal B}_{q}^{(\neq 2)}$ of density-$q$ branch points $Z$ with ${\mathcal N}_{T, {\rm Pl}}(Z) \neq 2,$ for arbitrary integer $q \geq 2$. This 
center manifold-free analysis suffices to establish the existence of a unique, non-zero, ${\mathcal N}_{T, {\rm Pl}}(Z)$-homogeneous, cylindrical tangent function to $T$ relative to 
the (unique) tangent plane $P_{Z}$ at ${\mathcal H}^{n-2}$-a.e.\ point $Z \in {\mathcal B}_{q}^{(\neq 2)}$ (see \cite[Theorem~1.8]{KrumWicc}) and consequently,  countable 
$(n-2)$-rectifiability of ${\mathcal B}_{q}^{(\neq 2)}$ (see \cite[Theorem~1.10]{KrumWicc} and \cite[Remark~1.11]{KrumWicc}), a conclusion contained in the structure theorem in (3) above. However, unless 
${\mathcal N}_{T, {\rm Pl}}(Z) \not\in {\mathbb N},$
this analysis need not 
establish the more precise asymptotic normal forms described in (1) above.

For branch points 
$Z$ with planar frequency ${\mathcal N}_{T, {\rm Pl}}(Z) \in {\mathbb N}$, 
the asymptotic normal form in (1), and hence also the topological conclusion  in (2), require the use of an Almgren center manifold, but in this case the center manifold arises in a certain ``canonical form'' with crucial additional structure (described in the next subsection) that enables the structural analysis of $T$. Utilizing such a center manifold, which we shall refer to as an \emph{effective center manifold}, 
the present paper establishes  Theorem~\ref{unique tangent function thm}, Theorem~\ref{topology thm} and Theorem~\ref{B rectifiable thm} below, which are the versions, for points with planar frequency at least $2,$ of the three main theorems described in (1), (2) and (3). In this case, the branching order ${\mathcal O}_{T}(Z)$ is the Almgren frequency at $Z$ of the normal map relative the effective center manifold (see the next subsection). Our asymptotic normal form theorem implies that this quantity is intrinsically determined by $Z,$ independently of the choice of center manifold. 

The center manifold-based  analysis in the present paper is the same for all planar-frequency values at least $2$, whether integer or not, 
and thus we do not assume in this paper that the planar frequency is an integer. 
As just mentioned, for non-integer planar frequencies, however, the analysis here is not needed, since it yields (when the planar frequency is a non-integer $>2$) exactly the same asymptotic normal form as that provided by Part~III through its intrinsic, center-manifold-free arguments.

 As discussed in the introduction to \cite{KrumWicc},  outcome (1) above in full generality (i.e.\ \cite[Theorem~B]{KrumWicc}) is then a direct consequence of \cite[Theorem~1.8]{KrumWicc} and Theorem~\ref{unique tangent function thm} below. Similarly, outcome (2) (i.e.\ \cite[Theorem~A]{KrumWicc}) follows directly from  \cite[Theorem~1.12]{KrumWicc} and Theorem~\ref{topology thm} below, while outcome (3) (i.e.\ \cite[Theorem~C]{KrumWicc}) is achieved  by combining \cite[Theorem~1.1]{KrumWicb}, \cite[Theorem~1.10]{KrumWicc} and Theorem~\ref{B rectifiable thm} below.  

Thus, our main task in the present paper is to establish Theorems~\ref{unique tangent function thm}, \ref{topology thm} and \ref{B rectifiable thm} below. 
As mentioned above, these theorems all concern, for arbitrary integer $q \geq 2$, 
density $q$ branch points $Z$ of planar frequency ${\mathcal N}_{T, {\rm Pl}}(Z) \geq 2$. 

\subsection{The role of Almgren's center manifolds in the present framework: reduction to effective center manifolds}
The role played by center manifolds in branch-point analysis  constitutes a fundamental conceptual 
 difference between our framework and the pioneering theory of Almgren,  which established  the sharp Hausdorff dimension upper bound $n-2$ for the singular set of $T$. 
 
 Almgren's monumental work, available as an unpublished preprint since the early 1980's,  was posthumously transcribed into  a book (\cite{Almgren}) by Taylor and Scheffer around the year 2000.
 Subsequently, De~Lellis and Spadaro revisited this theory in a series of papers (\cite{DeLSpa1, DeLSpa2, DeLSpa3}), adhering closely to Almgren's original conceptual architecture while incorporating  technical and expository improvements that enhanced accessibility of the work, particularly concerning the construction and use of center manifolds.

This framework has been adopted in recent work initiated by De~Lellis and Skorobogatova (\cite{DeLSko1, DeLSko2}) and completed by De~Lellis, Minter and Skorobogatova (\cite{DeLMinSko}), which establishes---contemporaneously with and independently of our work---two of the principal results of our program. Specifically,  
their work establishes ${\mathcal H}^{n-2}$-a.e.\ uniqueness of tangent cones relying on a construction of center manifolds at every branch point. In our approach, by contrast, this result is established at an early stage, using the planar frequency function and without recourse to center manifolds (see \cite[Theorem 1.1]{KrumWicb}).

The second main result, namely countable $(n-2)$-rectifiability of the singular set, is obtain in their approach via local upper Minkowski content bounds for a countable decomposition of the branch set. In our program this result is subsumed by the more refined \emph{locally-finite} locally-rectifiable decomposition theorem  stated in (3) above, which is established through asymptotic normal forms.  
We refer the reader to \cite[Section~2]{KrumWica} and \cite[Section~1.9]{KrumWicc} for a comprehensive and detailed comparison of the two approaches.

In \cite{Almgren, DeLSpa1, DeLSpa2, DeLSpa3} as well as in \cite{DeLSko1, DeLSko2}, the key foundational idea in the branch-set analysis is to iteratively construct and use, for an arbitrary branch point $Z$, a sequence of center manifolds  corresponding to a choice of scales (tending to zero) at which the current about $Z$ remains weakly close to a plane. 
In particular, this is done without considering the question of uniqueness of the tangent plane at $Z$.
One then blows up the corresponding sequence of multi-valued, Lipschitz \emph{normal approximation maps} defined over the varying center manifolds. These maps parameterize large portions of $T$ graphically over the center manifolds. However, since the center manifolds are not guaranteed to pass though any branch point, the construction in general does not provide a mechanism to precisely encode the current's branch-set information into the normal maps.

Nonetheless, the center manifolds are constructed so that the normal approximation maps are nearly average free in a quantitative sense; in other words, the center manifolds provide a good approximation to the average of the ``sheets'' of $T$ on approach to $Z$. This property is used to establish approximate monotonicity of the Almgren frequency function associated with the normal approximation maps. Such monotonicity holds over a range of scales at each point of the center manifold but not necessarily down to arbitrarily small scales at \emph{any} point. The nearly-average-free property of the normal maps also ensures that the blow-up limit is average-free. 
A doubling property for the normal map implied by frequency monotonicity guarantees that the blow-up limit is non-trivial. Thus, the  blow-up procedure prevents sheets from all collapsing onto one another to yield a single-valued harmonic function, with some integer multiplicity, as the blow-up limit. It then follows that the top-multiplicity points of the blow-up limit are contained in its singular set. In \cite{Almgren} this allows for the Hausdorff dimension of the branch set of $T$ to be deduced, essentially, from the Hausdorff dimension of the singular set of the blow-up limits 
(with certain exceptional cases handled differently); and in \cite{DeLSko1, DeLSko2}, this provides the starting point for the study of the fine properties of the branch set.

The drawbacks 
of using center manifolds in this generality are two fold. First, the center manifolds are not guaranteed to contain the base point $Z,$ nor indeed any other branch point. 
This presents a substantial obstacle from the perspective of studying the local asymptotic and topological structure of the whole current $T$ near a branch point $Z$---a key goal of our program. 
Second, even for the purpose of establishing the upper bound on the Hausdorff dimension of the singular set, the blow-up argument requires an exceedingly high technical output because, at general branch points, different center manifolds and normal maps must be chosen at different scales.

If, however, uniqueness of the  tangent plane at a branch point together with specific decay information is available a priori, the utility of the center manifold enhances significantly. More precisely, if a density $q$ branch point $Z$ is such that $T$ decays quadratically towards a unique tangent plane $P_{Z}$ at $Z$,  or equivalently, if the planar frequency ${\mathcal N}_{T, {\rm Pl}}(Z)$  exists and is equal to $2$, then the Almgren center manifold at $Z$ becomes an \emph{effective center manifold}: it arises ``canonically'' with  the following properties:
\begin{itemize}
\item[(i)] there is a \emph{single} center manifold ${\mathcal M}_{Z} = Z + {\rm graph} \, \psi^{(Z)}$, where $\psi^{(Z)} \, : \, P_{Z} \to P_{Z}^{\perp}$ is a function of class $C^{3, \kappa}$, $\kappa = \kappa(n, m, q) \in (0, 1)$, such that  ${\mathcal M}_{Z}$ contains $Z$ and all other 
nearby density-$q$ branch points $Z^{\prime}$ with ${\mathcal N}_{T, {\rm Pl}}(Z^{\prime}) \geq 2$; 
\item[(ii)] there is a $q$-valued Lipschitz normal map $N \, : {\mathcal M}_{Z}\, \to {\mathcal A}_{q}({\mathcal M}_{Z}^{\perp})$ which provides an approximate graphical parametrization for $T,$ with quantitative error control down to arbitrarily small scales about $Z$ (in contrast to a range of scales bounded away from zero as is the case for general branch points);
\item[(iii)] ${\mathcal M}_{Z}$ is tangential to $T$ at $Z$ and at all branch points $Z^{\prime}$ as above, and these points are contained in the critical nodal set of $N$; 
\item[(iv)] the Almgren frequency of $N$ at $Z$ and at any branch point $Z^{\prime} \in {\mathcal M}_{Z}$ as above is approximately monotone down to arbitrarily small scales;
\item[(v)] under scaling about $Z$ or any nearby branch point $Z^{\prime}$ as above, the approximation of $T$ by $N$ improves, with error decaying  quantitatively to zero at a sufficiently rapid rate as the scale tends to zero.  
\end{itemize}
More generally,  
all of the above continue to hold when ${\mathcal N}_{T, {\rm Pl}}(Z) \geq 2$, by exactly the same construction as in the case ${\mathcal N}_{T, {\rm Pl}}(Z) = 2$. In our framework the role of the center manifold is confined to this effective setting. The properties (i)--(v) above are central to the local structural analysis of $T$ near branch points $Z$ with ${\mathcal N}_{T, {\rm Pl}}(Z) \geq 2$ that we carry out in the present paper.

The reduction to effective center manifolds is made possible by 
planar-frequency monotonicity. Specifically, planar-frequency monotonicity allows the analysis to focus, at the outset, on establishing uniqueness of tangent cones and, crucially, uniform decay estimates at branch points, prior to introducing center manifolds. Since possible rapid merging of  ``sheets'' of $T$ together on approach to a branch point would facilitate the proof of such decay estimates, it is natural to expect that center manifolds---which serve to rule out this behaviour, namely to separate the sheets---can be avoided altogether in 
establishing planar-decay estimates at branch points. Confirming this expectation, and as discussed above, in Parts~I~and~II (\cite{KrumWica, KrumWicb}) we use the planar frequency function to establish locally uniform decay estimates at ${\mathcal H}^{n-2}$-a.e.\ fixed-density branch point. Then in Part~III we exploit planar-frequency monotonicity once more, using it in its full strength in place of center manifold-based analysis, to  study the asymptotic behaviour of $T$ near branch points of planar frequency $\neq 2$ via intrinsic analysis on the current. 

For branch points $Z$ of planar-frequency $2$ (or more generally, of planar frequency $\geq 2$)---which are the focus of the present paper---properties (i)--(v) above 
 serve essentially to reduce the structural analysis of $T$ near $Z$ to a structural analysis of the $q$-valued normal map $N$; in particular, analyzing the structure of the density-$q$, planar-frequency $\geq 2$ branch set of $T$ near $Z$ amounts to studying the critical nodal set of $N$. In other words, the problem becomes one of studying a multi-valued map subject to variational constraints.
 
We carry out this asymptotic analysis of $N$ by broadly adapting to $N$'s variational structure the strategy developed in our previous work \cite{KrumWic2} for the analysis of Dirichlet energy minimizing $q$-valued maps.  
This theory, in turn, was inspired by the general framework developed in the combined works \cite{Sim93} and \cite{Wic14} for multiplicity 1 classes of minimal submanifolds and the higher-multiplicity setting of stable codimension 1 stationary varifolds, respectively.

In the setting of \cite{KrumWic2}, without loss of generality, the $q$-valued Dirichlet energy minimizing map $v$ is assumed to have zero average (analogous, in the present setting, to $N$ being nearly average free), by subtracting the single-valued harmonic average of $v.$ As is readily verified, the resulting average-free $q$-valued function remains Dirichlet energy minimizing.    Almgren's classical interior H\"older continuity estimate for $v$ then provides uniform convergence of rescalings of 
$v,$ about a multiplicity-$q$ singular point, towards a (homogeneous) $q$-valued, singular tangent function. This uniform convergence guarantees that, at ${\mathcal H}^{n-2}$-a.e.\ singular point, (the rescaled) $v$ is well-separated  in the region away from the spine  of the tangent function. Here, the spine is the $(n-2)$-dimensional subspace of translation invariance, consisting of multiplicity $q$ points of the tangent function. This separation property is of fundamental importance to the asymptotic analysis carried out in \cite{KrumWic2}. 

In the present work, establishing the analogous separation property is considerably more subtle. A continuity estimate for $N$ that would guarantee uniform convergence of blow-up sequences is neither available nor, by itself,  sufficient to imply the required well-separatedness of $T,$ since $T$ need not agree with the graph of $N$ everywhere. Nevertheless, away from the $(n-2)$-dimensional spine of a cylindrical blow-up of $N$ at a branch point $Z^{\prime}$ as above (noting that blow-ups are indeed cylindrical at ${\mathcal H}^{n-2}$-a.e.\ such point $Z^{\prime}$), the whole re-scaled current is well separated. This follows from a new $L^{\infty}-L^{2}$ estimate, established in Theorem~\ref{relharm separation main thm} below. This estimate controls the pointwise distance of ${\rm spt} \, T$ from a disjoint union of graphs of single-valued nearly-harmonic   functions over a center manifold, whose average is zero, in terms of an $L^{2}$ excess quantity involving $N.$

\subsection{Statement of main results in the present paper}  Let $T$ be an $n$-dimensional locally area-minimizing rectifiable current of an open subset $U \subseteq \mathbb{R}^{n+m}$, and let $q$ be an integer $\geq 2.$ Let $Z_{0}$ be a density-$q$ rapid-decay branch point of $T$ (recalling that by~\cite[Theorem~X]{KrumWicc}, $\mathcal{H}^{n-2}$-a.e.~density-$q$ branch point is a rapid-decay branch point) and suppose that the planar frequency of $T$ at $Z_{0}$ is $\geq 2$.  
After translating $Z_0$ to the origin, rotating and rescaling, given a suitably small constant $\varepsilon \in (0,1)$ and a constant $\mu \in (0,1)$ we can assume:

\begin{hypothesis}\label{center mnfd intro hyp}{\rm
$T$ is an $n$-dimensional area-minimizing rectifiable current of $\mathbf{C}_{6\sqrt{n}}(0)$ such that 
\begin{gather*}
	(\partial T) \llcorner \mathbf{C}_{6\sqrt{n}}(0) = 0, \quad \Theta(T,0) = q, \quad 
	\|T\|(\mathbf{C}_{6\sqrt{n}}(0)) \leq q \omega_n (6\sqrt{n})^n + \varepsilon_2^2, \\ 
	\frac{1}{2\omega_n (6\sqrt{n})^n} \int_{\mathbf{C}_{6\sqrt{n}}(Y)} |\vec T - \vec P_0|^2 \,d\|T\| 
		= \inf_{P\in\mathcal{P}} \frac{1}{2\omega_n (6\sqrt{n})^n} \int_{\mathbf{C}_{6\sqrt{n}}(Y)} |\vec T - \vec P|^2 \,d\|T\| \leq \varepsilon_2^2 . \nonumber  
\end{gather*}
Here $\mathcal{P}$ is the set of all $n$-dimensional oriented planes $P$ and $\vec P$ denotes the orientation $n$-vector of $P \in \mathcal{P}$.  $P_0 = \mathbb{R}^n \times \{0\}$ equipped with the orientation $\vec P_0 = e_1 \wedge e_2 \wedge\cdots\wedge e_n$. 
}\end{hypothesis}

\begin{hypothesis}\label{center mnfd branch hyp}{\rm
For $\mathcal{H}^{n-2}$-a.e.~branch point $Z \in \op{sing} T \cap \mathbf{C}_1(0)$ with density $\Theta(T,Z) = q$, $T$ has a unique tangent plane $P_Z$ at $Z$ such that 
\begin{equation}\label{center mnfd hyp6} 
	\sup_{\rho \in (0,1]} \frac{1}{\rho^{n+2+2\mu}} \int_{\mathbf{B}_{2\rho}(Z)} \op{dist}^2(X,Z+P_Z) \,d\|T\|(X) < \infty 
\end{equation}
(cf.~\cite[Hypothesis~1.3]{KrumWicc}).
}\end{hypothesis}

\begin{definition}\label{center mnfd branch defn}{\rm 
Denote by $\mathcal{B}_q = \mathcal{B}_q(T)$ the set of all branch points $Z$ as in Hypothesis~\ref{center mnfd branch hyp}.  Let 
\begin{equation*}
	\mathcal{B}^{(\geq 2)}_q = \mathcal{B}^{(\geq 2)}_q(T) = \{ Z \in \mathcal{B}_q : \mathcal{N}_{T,{\rm Pl}}(Z) \geq 2 \} , 
\end{equation*}
where $\mathcal{N}_{T,{\rm Pl}}(Z)$ denotes the planar frequency of $T$ at $Z$ as defined in~\cite[Definition~4.1 and Definition~5.2]{KrumWica}.
}\end{definition}

\begin{remark}{\rm 
By~\cite[Theorem~3.2(e) and Lemma~3.6]{KrumWicc}, we can equivalently characterize $\mathcal{B}^{(\geq 2)}_q$ as the set of all density-$q$ branch points $Z \in \op{sing} T \cap \mathbf{C}_1(0)$ such that \eqref{center mnfd hyp6} holds true with $\mu = 1$.  Thus $\mathcal{B}^{(\geq 2)}_q$ is well-defined independently of the choice of constant $\mu$.
}\end{remark}

By a standard argument based on the existence of tangent functions relative to the center manifold and the upper semi-continuity of Almgren's frequency function on $\mathcal{B}^{(\geq 2)}_q$, we derive a stratification result, Theorem~\ref{centered stratification thm} below, which in particular gives us the optimal Hausdorff dimension upper bound for $\mathcal{B}^{(\geq 2)}_q$.  Using this together with~\cite[Theorem~1.1]{KrumWica} and~\cite[Theorem~1.5]{KrumWicc}, we recover the optimal Hausdorff dimension upper bound for the singular set as in~\cite{Almgren} (see~\cite[Theorem~1.3]{KrumWica} for further details).

\begin{theorem}\label{centered size sing thm}
Let $q \geq 2$ be an integer.  There exists $\varepsilon_2 = \varepsilon_2(n,m,q) \in (0,1)$ such that if $T$ is as in Hypothesis~\ref{center mnfd intro hyp} and $\mathcal{B}^{(\geq 2)}_q$ is as in Definition~\ref{center mnfd branch defn}, then $\mathcal{B}^{(\geq 2)}_q$ has Hausdorff dimension $\leq n-2$. Hence (by \cite[Theorem~1.1]{KrumWica} and \cite[Theorem~1.5]{KrumWicc}), the singular set of $T$ has Hausdorff dimension $\leq n-2$.
\end{theorem}

Our analysis of the fine structure of $\mathcal{B}_q^{(\geq 2)}$ gives us the following results.

\begin{theorem}\label{B rectifiable thm}
Let $q \geq 2$ be an integer.  There exists $\varepsilon_2 = \varepsilon_2(n,m,q) \in (0,1)$ such that the following holds true.  Let $T$ be as in Hypothesis~\ref{center mnfd intro hyp} and $\mathcal{B}^{(\geq 2)}_q$ be as in Definition~\ref{center mnfd branch defn}.  For each compact set $K \subset \mathbf{C}_1(0)$ 
\begin{equation*}
	\mathcal{B}^{(\geq 2)}_q \cap K = \bigcup_{j=1}^J L_j 
\end{equation*}
for some positive integer $J$ and for some pairwise disjoint, locally compact, locally $(n-2)$-rectifiable sets $L_1,\ldots,L_J$ (and in particular each set $L_j$ has locally finite $\mathcal{H}^{n-2}$-measure in the sense that for each $X \in L_j$ there exists $\delta > 0$ such that $\mathcal{H}^{n-2}(L_j \cap \mathbf{B}_{\delta}(X)) < \infty$).
\end{theorem}

\begin{theorem}\label{unique tangent function thm} 
Let $q \geq 2$ be an integer.  There exists $\varepsilon_2 = \varepsilon_2(n,m,q) \in (0,1)$ such that the following holds true.  Let $T$ be as in Hypothesis~\ref{center mnfd intro hyp} and $\mathcal{B}^{(\geq 2)}_q$ be as in Definition~\ref{center mnfd branch defn}.  For $\mathcal{H}^{n-2}$-a.e.~$Z \in \mathcal{B}^{(\geq 2)}_q$ there exists a radius $\rho_Z \in (0,1/2]$, a (unique) positive rational number $\mathcal{O}_T(Z)$, a $C^{3,\kappa}$-function $\psi^{(Z)} : P_Z \cap \mathbf{B}_{2\rho_Z}(0) \rightarrow P_Z^{\perp}$ (where $\kappa = \kappa(n,m,q) \in (0,1)$), a (unique and possibly zero) homogeneous degree $\mathcal{N}_{T,{\rm Pl}}(Z)$, harmonic single-valued function $h^{(Z)} : P_Z \rightarrow P_Z^{\perp}$ and a (unique) non-zero, homogeneous degree $\mathcal{O}_T(Z)$, average-free, locally Dirichlet energy minimizing $q$-valued function $\varphi^{(Z)} : P_Z \rightarrow \mathcal{A}_q(P_Z^{\perp})$ such that:

\begin{enumerate}[itemsep=2mm,topsep=0mm]
	\item[{\rm (i)}]  $\mathcal{O}_T(Z) = \ell_Z/q_Z$ for some relatively prime positive integers $\ell_Z,q_Z$ with $\ell_Z \geq 2q_Z$ and $q_Z \leq q$; 
	\item[{\rm (ii)}] ${\mathcal M}_{Z} = Z + {\rm graph} \, \psi_{Z}$ is a center manifold for $T$ at $Z$ (as in Subsection~\ref{sec:center mnfd prelims}); 
	\item[{\rm (iii)}] $\varphi^{(Z)}$ is a tangent function of $T$ at $Z$ as constructed in Subsection~\ref{sec:tangent fun} as a blow-up of $T$ relative to ${\mathcal M}_{Z}$;   $\varphi^{(Z)}$ is uniquely determined by $T$ and $Z$; and coordinates $(z, y) \in {\mathbb C} \times {\mathbb R}^{n-2}$ on $P_{Z}$ relative to an orthonormal basis can be chosen with respect to which $\varphi^{(Z)}$ is independent of $y$ and has the form 
	\begin{equation*}
		\varphi^{(Z)} = \sum_{j=1}^{N_Z} m_j^{(Z)} \varphi^{(Z)}_j 
	\end{equation*}
	where $N_Z,m^{(Z)}_1,\ldots,m^{(Z)}_{N_Z}$ are positive integers and $\varphi^{(Z)}_j : P_Z \rightarrow \mathcal{A}_{q^{(Z)}_j}(P_Z^{\perp})$ are distinct functions such that $\sum_{j=1}^{N_Z} m^{(Z)}_j q^{(Z)}_j = q$ and for each $j \in \{1,2,\ldots,N_Z\}$ either:
    \begin{itemize}[itemsep=2mm,topsep=2mm]
    \item $q^{(Z)}_j = 1$ and $\varphi^{(Z)}_j(z,y) = 0$ for all $(z, y) \in {\mathbb C} \times {\mathbb R}^{n-2}$, or 
    \item $q^{(Z)}_j = q_Z$ and $\varphi^{(Z)}_j(z,y) = \op{Re} \, \left(c^{(Z)}_j z^{\ell_Z/q_Z}\right)$ for all $(z, y) \in {\mathbb C} \times {\mathbb R}^{n-2}$ and some $c^{(Z)}_j \in \mathbb{C}^{n+m} \setminus \{0\}$ with $\op{Re} c^{(Z)}_j, \op{Im} c^{(Z)}_j \in P_Z^{\perp};$ 
    \end{itemize}
	\item[{\rm (iv)}]  $T$ is asymptotic to the graph of $\psi^{(Z)} + \varphi^{(Z)}$ at $Z$ in the sense that 
	\begin{align}\label{unique tangent function concl1} 
		&\frac{1}{\omega_n \rho^n} \int_{\mathbf{C}_{\rho}(Z,P_Z)} 
			\op{dist}^2(X, Z + \op{spt} \mathbf{G}_{\psi^{(Z)} + \varphi^{(Z)}}) \,d\|T\|(X) 
		\\&+ \sum_{j=1}^{N_Z} \frac{m^{(Z)}_j}{\omega_n \rho^n} \int_{B_{\rho/2}(0,P_Z) \cap \{ r < \rho/16\}} 
			\sum_{l=1}^{q^{(Z)}_j} \op{dist}^2(Z + X + \psi^{(Z)}(X) + \varphi^{(Z)}_{j,l} (X), \op{spt} T) \,d\mathcal{H}^n(X) \nonumber 
		\\&\leq C_Z  \rho^{2\mathcal{O}_T(Z)} \bigg(\frac{\rho}{\rho_Z}\bigg)^{2\mu_Z} \nonumber  
	\end{align}
	for all $\rho \in (0,\rho_Z]$, where $r(X) = |z|$ is the distance to the linear subspace 
    $\{z=0\} \subset P_{Z}$
    (using the same coordinates $(z, y)$ as in {\rm (iii)}), $\varphi^{(Z)}_j(X) = \sum_{l=1}^{q^{(Z)}_j} \llbracket \varphi^{(Z)}_{j,l}(X) \rrbracket$ for each $X \in P_Z$, $\mathbf{G}_{\psi^{(Z)} + \varphi^{(Z)}}$ is the graph of $\psi^{(Z)} + \varphi^{(Z)}$ as a rectifiable current (as defined in Subsection~\ref{sec:normal graph sec}), and $\mu_Z = \mu(n,m,q,\mathcal{O}_T(Z)) \in (0,\infty)$ and $C_Z = C(n,m,q,\mathcal{O}_T(Z)) \in (0,\infty)$ are constants;
	\item[{\rm (v)}]  $\psi^{(Z)}$ is asymptotic to $h^{(Z)}$ in the sense that 
    \begin{equation}\label{unique tangent function concl2}
        \frac{1}{\rho^n} \int_{B_{\rho}(0,P_Z)} |\psi^{(Z)}(X) - h^{(Z)}(X)|^2 \leq C_Z \rho^{2\mathcal{N}_{T,{\rm Pl}}(Z)} \bigg(\frac{\rho}{\rho_Z}\bigg)^{2\mu_Z}
    \end{equation}
    where $\mu_Z = \mu(n,m,q,\mathcal{O}_T(Z)) \in (0,1)$ and $C_Z = C(n,m,q,\mathcal{O}_T(Z)) \in (0,1)$ are constants;
    \item[{\rm (vi)}] $\mathcal{O}_T(Z) \geq \mathcal{N}_{T,{\rm Pl}}(Z)$ and moreover $\mathcal{O}_T(Z) > \mathcal{N}_{T,{\rm Pl}}(Z)$ $\iff$ $\mathcal{O}_T(Z) \geq \mathcal{N}_{T,{\rm Pl}}(Z)+1/q$ $\iff$ 
    \begin{equation}\label{unique tangent function concl3}
        \lim_{\rho\downarrow 0} \rho^{-n-2\mathcal{N}_{T,{\rm Pl}}(Z) - 2\mu'} \int_{\mathbf{B}_{\rho}(Z)} \op{dist}^2\left(X,Z+\op{spt}\mathbf{G}_{h^{(Z)}}\right) \,d\|T\|(X) = 0 
    \end{equation}
    for some constant $\mu' > 0$; in particular, if $\mathcal{N}_{T,{\rm Pl}}(Z)$ is not an integer then  $\mathcal{O}_T(Z) = \mathcal{N}_{T,{\rm Pl}}(Z),$ $h^{(Z)} = 0$ and $\varphi^{(Z)}$ is also the unique tangent function of $T$ at $Z$ relative to the tangent plane $P_{Z}$; 
    \item[{\rm (vii)}] if $\widetilde{\mathcal O}_T(Z)$ is a positive real number, $\widetilde{\psi}^{(Z)} : P_Z \rightarrow P_Z^{\perp}$ is a $C^3$ single-valued function and $\widetilde{\varphi}^{(Z)} : P_Z \rightarrow \mathcal{A}_q(P_Z^{\perp})$ is a non-zero, average-free, homogeneous degree $\widetilde{\mathcal O}_T(Z)$ $q$-valued such that \eqref{unique tangent function concl1} holds true with $\widetilde{\mathcal O}_T(Z)$, $\widetilde{\psi}^{(Z)}$ and $\widetilde{\varphi}^{(Z)}$ in place of $\mathcal{O}_T(Z)$, $\psi^{(Z)}$ and $\varphi^{(Z)}$ respectively, 
    and possibly with different positive constants in place of $C_{Z}$, $\mu_{Z}$ and $\rho_{Z}$ in \eqref{unique tangent function concl1}, then 
    $\widetilde{\mathcal O}_{T}(Z) = {\mathcal O}_{T}(Z)$, $\varphi^{(Z)} = \widetilde{\varphi}^{(Z)}$ on $P_Z$ and 
	\begin{equation}\label{unique tangent function concl4}
		\frac{1}{\rho^n} \int_{B_{\rho}(0,P_Z)} |\psi^{(Z)} - \widetilde{\psi}^{(Z)}|^2 \leq C \rho^{2\mathcal{O}_T(Z)+2\mu} 
    \end{equation}
    for all sufficiently small $\rho > 0$, where $\mu \in (0,1)$ and $C \in (0,1)$ are constants; moreover, if $\widetilde{h}^{(Z)}$ is a homogeneous degree $\mathcal{N}_{T,{\rm Pl}}(Z)$, single-valued function such that \eqref{unique tangent function concl2} holds true with $\widetilde{\psi}^{(Z)}$ and $\widetilde{h}^{(Z)}$ in place of $\psi^{(Z)}$ and $h^{(Z)}$ respectively, and possibly with different positive constants in place of $C_Z$, $\mu_Z$ and $\rho_Z$ in \eqref{unique tangent function concl2}, then $\widetilde{h}^{(Z)} = h^{(Z)}$ on $P_Z$.
\end{enumerate}
\end{theorem}

\begin{remark}{\rm 
By the upper semi-continuity of $Z \mapsto \mathcal{O}_T(Z)$, in \eqref{unique tangent function concl1} and \eqref{unique tangent function concl2}, $\mu_Z$ and $C_Z$ can be chosen to be locally uniform, i.e.~for $\mathcal{H}^{n-2}$-a.e.~$Z \in \mathcal{B}^{(\geq 2)}_q$ there is a radius $\delta_Z > 0$ such that for $\mathcal{H}^{n-2}$-a.e.~$Z' \in \mathcal{B}^{(\geq 2)}_q \cap \mathbf{B}_{\delta_Z}(Z)$, \eqref{unique tangent function concl1} and \eqref{unique tangent function concl2} hold true with $Z'$ in place of $Z$ and with $\mu_{Z'} = \mu_Z$ and $C_{Z'} = C_Z$ (independent of $Z'$);
}\end{remark}

\begin{remark}
{\rm By~\cite[Theorem~3.2]{MicWhi95}, the constants $c^{(Z)}_j$ in Theorem~\ref{unique tangent function thm}(iii) satisfy $$c^{(Z)}_j \cdot c^{(Z)}_j = 0.$$ Writing $a_{j}^{(Z)} = {\rm Re} \, (c_{j}^{(Z)})$ and $b_{j}^{(Z)} = {\rm Im} \, (c_{j}^{(Z)}),$ this means, as observed in \cite{MicWhi95}, 
that 
$\Sigma^{(Z)}  \equiv {\rm span}_{\mathbb R} \, \{a_{j}^{(Z)}, b_{j}^{(Z)}\}$ is a two-dimensional subspace of $P_{Z}^{\perp}$ on which there is a complex structure $J \, : \, \Sigma^{(Z)} \to \Sigma^{(Z)}$ (defined by $J(a_{j}^{(Z)}) = -b_{j}^{(Z)}$, $J(b_{j}^{(Z)}) = a_{j}^{(Z)}$), which is compatible with the Euclidean metric, with respect to which the (multi-valued) graph of $$z \mapsto \varphi_{j}^{(Z)}(z, 0) = {\rm Re} \, \left(c_{j}^{(Z)}z^{\ell_{Z}/q_{Z}}\right)$$ 
is a holomorphic curve in ${\mathbb C} \times \Sigma^{(Z)}$.  Thus, Theorem~\ref{unique tangent function thm} can be seen as providing a local asymptotic converse, valid near ${\mathcal H}^{n-2}$-a.e.\ branch point of an area-minimizing rectifiable current, to the well-known theorem that holomorphic subvarieties are area-minimizing.}
\end{remark}

Let $\mathcal{M}$ be a center manifold of $T$ and $N : \mathcal{M} \rightarrow \mathcal{A}_q({\rm T}\mathcal{M}^{\perp})$ be a normal approximation map of $T$.  Then the branching order $\mathcal{O}_T(Z)$ is equal to the frequency $\mathbf{I}_{N,Z}(0^+)$ of $N$ (as in \eqref{centered freq}) and we can choose $\psi^{(Z)} : P_Z \cap \mathbf{B}_{2\rho_Z}(0) \rightarrow P_Z^{\perp}$ to be the function whose graph coincides with $\mathcal{M}$, i.e.~$\op{graph} \psi^{(Z)} = \mathcal{M} \cap \mathbf{C}_{2\rho_Z}(Z,P_Z)$.  $\mathcal{M}$ and $N$ are not uniquely determined by their construction; nonetheless, Theorem~\ref{unique tangent function thm} asserts that $\mathcal{O}_T(Z)$ is uniquely determined by $T$ and $Z$ and is independent of the choice of $\mathcal{M}$ and $N$.  Moreover, while $\psi^{(Z)}$ is not unique, 
it is unique up to order $O(|X-Z|^{{\mathcal O}_{T}(Z)})$ (as in \eqref{unique tangent function concl4}).  By standard doubling conditions, $\rho^{-n/2} \|N\|_{L^2(B^{\mathcal M}_{\rho}(Z))}$ decays at a rate between $O(\rho^{\mathcal{O}_T(Z)+\epsilon})$ and $O(\rho^{\mathcal{O}_T(Z)})$ for each $\epsilon > 0$.  The proof of Theorem~\ref{unique tangent function thm} implies the stronger conclusion that $\rho^{-n/2-\mathcal{O}_T(Z)} \|N\|_{L^2(B^{\mathcal M}_{\rho}(Z))}$ converges to a positive (and finite) real number.

As explained in~\cite{KrumWicc}, Theorem~\ref{unique tangent function thm} gives us some preliminary information about the topology and geometric structure of $T$ near $\mathcal{H}^{n-2}$-a.e.~$Z \in \mathcal{B}^{(\geq 2)}_q$.  For instance, if the components $\varphi^{(Z)}_j$ have multiplicity $m^{(Z)}_j = 1$, by the Allard regularity theorem (\cite{Allard}) we can express $T$ as a smooth perturbation of the graph of $\varphi^{(Z)}$ away from the spine.  However, this gives us little information about the topology of $T$ near the spine of $\varphi^{(Z)}$.  Under the additional hypotheses that Almgren's frequency is sufficiently close to $\mathbf{I}_{N,Z}(0^+)$ at some scale (a condition automatically satisfied if the scale is sufficiently small) and $Z$ is \emph{$(n-2)$-strongly non-isolated} at that scale (i.e.\  \eqref{topology hyp2} holds), we show that $T$ is a ``$C^1$-perturbation of the graph of $\varphi^{(Z)}$;'' in particular, ${\rm spt} \, T$ is topologically an $n$-dimensional disk near $Z$ and admits a $C^{1, \mu}$ parameterization. As pointed out in the introductions to \cite{KrumWica} and \cite{KrumWicc}, explicit algebraic examples demonstrate that $(n-2)$-strong non-isolatedness of $Z$ cannot be relaxed to requiring 
$L \cap \mathbf{B}_{\lambda\rho_0}(0) 
		\subseteq \pi_L(\{ X - Z_0 : X \in \mathcal{B}^{(\geq 2)}_q \cap \mathbf{B}_{\rho_0/4}(Z_0)\}).$

\begin{theorem}\label{topology thm} 
Let $\alpha = p/q \geq 2$ where $p,q$ are relatively prime integers.  For each $\lambda \in (0,1/4)$ there exists $\varepsilon_2 = \varepsilon_2(n,m,q) \in (0,1)$ and $\varepsilon = \varepsilon(n,m,q,\alpha,\lambda) \in (0,1)$ such that the following holds true.  Let $T$ be as in Hypothesis~\ref{center mnfd intro hyp} and $\mathcal{B}^{(\geq 2)}_q$ be as in Definition~\ref{center mnfd branch defn}.  Let $\mathcal{M}$ be the center manifold of $T$ and $N : \mathcal{M} \rightarrow {\rm T}\mathcal{M}^{\perp}$ be the normal approximation map of $T$.  Let $Z_0 \in \mathcal{B}^{(\geq 2)}_q$ and $\rho_0 \in (0,1]$ such that 
\begin{equation}\label{topology hyp1}
	\mathbf{I}_{N, Z_{0}}(0^+) \geq  \alpha, \quad \mathbf{I}_{N,Z_0}(\rho_0/2) < \alpha+\varepsilon ,
\end{equation}
where $\mathbf{I}_{N,Z_0}$ is the frequency function of $N$ (as in Definition~\ref{freq fun defn}) and $\mathbf{I}_{N,Z_0}(0^+)$ is the frequency of $N$ (as in \eqref{centered freq}).  Suppose that there exists an $(n-2)$-dimensional linear subspace $L \subset \mathbb{R}^{n+m}$ such that 
\begin{equation}\label{topology hyp2}
	L \cap \mathbf{B}_{\lambda\rho_0}(0) 
		\subseteq \pi_L(\{ X - Z_0 : X \in \mathcal{B}^{(\geq 2)}_q \cap \mathbf{B}_{\rho_0/4}(Z_0) \text{ and } \mathbf{I}_{N,X}(0^+) \geq \alpha \}) , 
\end{equation}
where $\pi_L : \mathbb{R}^{n+m} \rightarrow L$ is the orthogonal projection map onto $L$.  Then there exists a rotation $\Gamma$ of $\mathbb{R}^{n+m}$ with $\Gamma(P_{Z_0}) = \mathbb{R}^n \times \{0\}$, $g \in C^{1,\mu/2}(B^{n-2}_{\lambda\rho_0/4}(0),\mathbb{R}^2)$, $h \in C^{1,\mu/2}(B^{n-2}_{\lambda\rho_0/4}(0),\mathbb{R}^m)$, $\psi \in C^{1,\mu/2}(B^2_{\lambda\rho_0/32}(0) \times B^{n-2}_{\lambda\rho_0/32}(0),\mathbb{R}^m)$, and $f \in C^{1,\mu/2}(B^2_{(\lambda\rho_0/32)^{1/q}}(0) \times B^{n-2}_{\lambda\rho_0/32}(0),\mathbb{R}^m)$ such that: 
\begin{enumerate}[itemsep=2mm,topsep=0mm]
	\item[{\rm (a)}]  $\Gamma(\eta_{Z_0,1}(\op{sing} T)) \cap \mathbf{C}_{\lambda\rho_0/4}(0) = \Gamma(\eta_{Z_0,1}(\{ X \in \mathcal{B}_q^{(\geq 2)} : \mathbf{I}_{N,X}(0^+) = \alpha \})$ and 
	\begin{equation*}
		\Gamma(\eta_{Z_0,1}(\op{sing} T)) \cap \mathbf{C}_{\lambda\rho_0/4}(0) 
		= \{ (g(y), y, h(y)) : y \in B^{n-2}_{\lambda\rho_0/4}(0) \} \cap \mathbf{C}_{\lambda/4}(0) ; 
	\end{equation*}

	\item[{\rm (b)}]  the support of $T$ is parameterized by 
	\begin{equation*}
		\op{spt} \Gamma_{\#} \eta_{Z_0,1\#} T \cap \mathcal{D} \times \mathbb{R}^m 
		= \{ (g(y) + z^q, y, f(z,y)) : z \in B^2_{(\lambda\rho_0/32)^{1/q}}(0),\, y \in B^{n-2}_{\lambda\rho_0/32}(0) \} , 
	\end{equation*}
	where we identify $\mathbb{R}^2 \cong \mathbb{C}$ and $\mathcal{D} = \{ (g(y)+x,y) : x \in B^2_{\lambda\rho_0/32}(0),\, y \in B^{n-2}_{\lambda\rho_0/32}(0) \}$; 

	\item[{\rm (c)}]  there exists a unique tangent function $\varphi$ of $\Gamma_{\#} \eta_{Z_0,1\#} T$ at the origin which is given by $\varphi(x_1,x_2,y) = \op{Re}(c \,(x_1+ix_2)^{p/q})$ for each $(x_1,x_2,y) \in \mathbb{R}^n$ and for some $c \in \mathbb{C}^m \setminus \{0\}$ with $c \cdot c = 0$ and 
	\begin{equation*}
		f(z,y) = \psi(z^q,y) + \op{Re}(cz^p) + e(z,y) 
	\end{equation*} 
	for $z \in B^2_{(\lambda\rho_0/32)^{1/q}}(0)$ and $y \in B^{n-2}_{\lambda\rho_0/32}(0)$, where 
	\begin{gather*}
		|e(z,y)| \leq C \widehat{E} \,|(z^q,y)|^{p/q+\mu/(2nq)} , \quad 
		|D_z e(z,y)| \leq C \widehat{E} \,|(z^q,y)|^{p/q-1+\mu/(2nq)} |z|^{q-1} , \\
		|D_y e(z,y)| \leq C \widehat{E} \,|(z^q,y)|^{p/q-1+\mu/(2nq)} , 
	\end{gather*}
	where $\widehat{E} = \rho_0^{-n/2-1} \|N\|_{L^2(B^{\mathcal M}_{\rho_0}(Z_0))}$ and $C = C(n,m,q,\alpha) \in (0,\infty)$ is a constant; 
	
	\item[{\rm (d)}]  for any fixed $q^{\rm th}$ root of unity $\nu \neq 1$
	\begin{equation*}
		f(z,y) - f(\nu z,y) = \op{Re}(az^p) + \widetilde{e}(z,y) 
	\end{equation*} 
	for $z \in B^2_{(\lambda\rho_0/32)^{1/q}}(0)$ and $y \in B^{n-2}_{\lambda\rho_0/32}(0)$, where $a \in \mathbb{C}^m \setminus \{0\}$ with $a \cdot a = 0$ and 
	\begin{gather*}
		|\widetilde{e}(z,y)| \leq C \widehat{E} \,|(z^q,y)|^{p/q+\mu/(2nq)} , \quad 
		|D_z \widetilde{e}(z,y)| \leq C \widehat{E} \,|(z^q,y)|^{p/q-1+\mu/(2nq)} |z|^{q-1} , \\
		|D_y \widetilde{e}(z,y)| \leq C \widehat{E} \,|(z^q,y)|^{p/q-1+\mu/(2nq)} , 
	\end{gather*}
	where $C = C(n,m,q,\alpha) \in (0,\infty)$ is a constant.
\end{enumerate}
\end{theorem}

As was noted in~\cite[Remark~1.10]{KrumWicc}, it follows from the proof of Theorem~\ref{topology thm} that for each $Z \in \op{sing} T \cap \mathbf{C}_{\lambda\rho_0/32}(Z_0,P_{Z_0})$ there exists a rotation $\Gamma_Z$ of $\mathbb{R}^{n+m}$ and functions $g_Z \in C^{1,\mu/2}(B^{n-2}_8(0),\mathbb{R}^2)$, $\psi_Z \in C^3(B^2_1(0) \times B^{n-2}_1(0),\mathbb{R}^m)$, and $f_Z \in C^{1,\mu/2}(B^2_1(0) \times B^{n-2}_1(0),\mathbb{R}^m)$ such that conclusion~(a)--(d) of Theorem~\ref{topology thm} hold true with $\Gamma_Z,g_Z,\psi_Z,f_Z$ in place of $\Gamma,g,f,\psi$.  It remains open whether there is a \emph{single} $C^1$-parameterization of $T$ near the origin such that the asymptotic expansions as in conclusions~(c) and (d) hold true at every singular point $Z$ near $Z_0$, in particular since it remains open whether the singular set of $T$ is a $C^{2,\mu}$ $(n-2)$-dimensional submanifold of $\mathbb{R}^{n+m}$ so that we can utilize normal coordinates relative to the singular set.  We refer the reader to~\cite[Remark~1.10]{KrumWicc} for a more detailed discussion of this matter.

The proofs of Theorem~\ref{B rectifiable thm}, Theorem~\ref{unique tangent function thm}, and Theorem~\ref{topology thm} as well as the analogous main results of~\cite{KrumWicc} use a blow-up method as developed in~\cite{Sim93} and~\cite{Wic14} and applied to the analysis of branch point singularities in~\cite{KrumWic2} and~\cite{KrumWic1}.  In~\cite{KrumWic2},~\cite{KrumWic1}, and~\cite{KrumWicc}, the blow-up method was implemented using the monotonicity formula for a frequency function relative to a plane---either the fixed plane $\mathbb{R}^n \times \{0\}$ (as in \cite{KrumWic2}) or the varying tangent planes at a branch points (as in \cite{KrumWic1, KrumWicc}).  Here we develop the blow-up method 
for the normal approximation map $N$, using the monotonicity formula for Almgren's frequency function for $N$ and the associated variational properties inherited by $N$ from the current $T$. In contrast to \cite{KrumWic2, KrumWic1, KrumWicc}, the additional subtlety here is that the domain of parameterization is not a plane but the  curved center manifold ${\mathcal M}.$   

Let $T$ be a locally area-minimizing rectifiable current of $\mathbf{C}_{6\sqrt{n}}(0)$ which is weakly close to the plane $P_0 = \mathbb{R}^n \times \{0\}$ with multiplicity $q$ as in Hypothesis~\ref{center mnfd intro hyp}.  Let $\mathcal{M}$ be the center manifold of $T$ and $N : \mathcal{M} \rightarrow \mathcal{A}_q({\rm T}\mathcal{M}^{\perp})$ be the normal approximation map of $T$.  For each $Z \in \mathcal{M} \cap \mathbf{C}_1(0)$, we can find a critical radius $\overline{R}(Z) \in [0,1]$ such that (roughly speaking) the rescaled current $\eta_{Z,\rho\#} T$ decays to a plane at a nearly quadratic (or faster) rate at each scale $\rho \in (\overline{R}(Z),2]$.  The normal approximation map $N$ is constructed to be a Lipschitz approximation of $T$ which satisfies certain estimates in the geodesic ball $B^{\mathcal M}_{\rho}(Z)$ for each $\rho \in (\overline{R}(Z),2]$.  These estimates include bounds on the supremum and Lipschitz constant of $N$, the measure of the set of points where $T$ and the graph of $N$ do not coincide, and the $L^1$-norm of the single-valued average of $N$, the last of which guarantees that $N$ is approximately average-free. In particular, $\mathcal{B}_q^{(\geq 2)}$ is precisely the set of points $Z \in \mathcal{M} \cap \mathbf{C}_1(0)$ for which $\overline{R}(Z) = 0$, so that the center manifold $\mathcal{M}$ passes through $\mathcal{B}_q^{(\geq 2)}$ and the frequency function of $N$ is approximately monotone non-decreasing for all $\rho \in (0,2]$ at each $Z \in \mathcal{B}_q^{(\geq 2)}$.

A key step in developing the blow-up method in the present setting is to adapt~\cite[Theorem~3.15]{KrumWicb} and~\cite[Theorem~4.19]{KrumWicc} to prove local estimates on the uniform distance of the support of $T$ to a finite union of nearly parallel graphs of single-valued nearly harmonic functions defined over the center manifold, whose average is zero.  Since our approach uses the estimates for the normal approximation map $N$ on geodesic balls, we need to carefully keep track of the required size of the critical radii $\overline{R}(X)$ at points $X \in B^{\mathcal M}_{\rho}(Z)$ and show that the conditions on the critical radii are satisfied when obtaining (in Section~\ref{sec:graph rep}) a graphical representation of $T$ relative to a cylindrical multi-valued harmonic function. We will also need to use~\cite[Theorem~4.19]{KrumWicc} to obtain estimates on the distance of $T$ to the union of the graphs of the harmonic functions in $B^{\mathcal M}_{\rho}(Z)$ for $\rho \geq \overline{R}(Z)$ (see Lemma~\ref{relharm critrad1 lemma} and Lemma~\ref{relharm critrad3 lemma}). The prior result \cite[Theorem~3.15]{KrumWicb}, and similarly~\cite[Theorem~4.19]{KrumWicc}, were proved via an inductive argument in which we separate $T$ into a sum of components $T_j$ near the graphs of the harmonic functions and then, in the case that there are two or more components, apply~\cite[Theorem~3.15]{KrumWicb} inductively to bound the distance of each component $T_j$ to the graphs of the harmonic functions.  However, in the present setting, the normal approximation map of each component $T_j$ relative to $\mathcal{M}$ need not be approximately average-free.  Thus to prove the estimates in the present paper, we shall assume that $T$ has already been separated into a sum of $p$ components $T_j$ (see Hypothesis~\ref{relharm separation1 hyp} and Hypothesis~\ref{relharm grand hyp}). Then we prove the estimate via induction, with the inductive hypothesis that the estimate holds true whenever we further decompose $T$ into $p+1$ or more components, and with the base case of induction being $p=q$. Note that $T$ is locally the sum of all of the components $T_j$ and the normal approximation map of $T$ continues to be approximately average-free throughout this argument.  The final theorem (Theorem~\ref{relharm separation main thm}) is the estimate in the case $p = 1$.

\noindent
{\bf AI disclosure:} All mathematical results and their proofs contained in this article and in the other four parts \cite{KrumWica, KrumWicb, KrumWicc, KrumWice} of the series to which this article belongs were developed by the authors without any involvement of generative AI tools. This five-part series, together with the authors' earlier papers \cite{KrumWic0, KrumWic1, KrumWic2}, documents the outcomes of a research program spanning the period 2013--2023 on area-minimizing currents and  related multi-valued harmonic functions and minimal graphs. 

The first versions of the first two parts \cite{KrumWica, KrumWicb} of the five-part series were posted on the arXiv in April, 2023 and contained a brief mention of the remainder of the research program, concerning branch-set analysis via asymptotic normal forms (therein called higher-order expansions) and building on \cite{KrumWic1, KrumWic2}. These two articles were subsequently revised, with expanded introductions. During the preparation of these revised versions, as well as of the manuscript \cite{KrumWicc}(the third part) and the present paper (the fourth part), which together implement the announced methodology for branch-set analysis, publicly available LLMs were used solely for language improvement in portions of the abstracts and  introductions.

\section{Preliminaries: notation and Almgren's center manifold}\label{sec:DeLSpa sec} 

In this section, we introduce the basic terminology and notation that will be used throughout the paper. We also discuss the center manifold and normal approximation map associated with a locally area-minimizing rectifiable current, and the frequency function of the normal approximation map.  These concepts were introduced by Almgren in~\cite{Almgren} and played a fundamental role in the blow-up arguments used in \cite{Almgren} to study the Hausdorff dimension of the singular set of such currents. 
The construction of the center manifold and normal approximation map and their application to the analysis of the singular set were later presented in a more concise, modular and accessible manner by De~Lellis and Spadaro in~\cite{DeLSpa2}\cite{DeLSpa3}.  Here we provide an overview of the theory of the center manifold and normal approximation map as presented in~\cite{DeLSpa2}\cite{DeLSpa3}, providing also the corresponding references to \cite{Almgren} as appropriate for completeness. We refer the reader to~\cite{DeLSpa2}\cite{DeLSpa3} for further details.

\subsection{General notation}\label{sec:overall notation sec}  Throughout the paper we adopt the same notation as in~\cite{KrumWica} and~\cite{KrumWicb} unless otherwise stated.  See~\cite[Section~3]{KrumWica} for a discussion of general notation, as well as an overview of locally area-minimizing rectifiable currents, multi-valued functions, and Dirichlet energy minimizing multi-valued functions.

We let $n \geq 1$ and $m \geq 2$ be fixed integers.  We often let $X$ denote a general point in $\mathbb{R}^{n+m}$.  $e_1,e_2,\ldots,e_{n+m}$ denotes the standard basis of $\mathbb{R}^{n+m}$, with $e_{1}, e_{2}, \ldots, e_{n}$ denoting the standard basis of ${\mathbb R}^{n} \times \{0\}$.  We let $\langle\,\cdot\,,\,\cdot\,\rangle$ denote the standard inner product on $\mathbb{R}^{n+m}$ and the induced inner product on tensors.

Given an $n$-dimensional plane $P$ in $\mathbb{R}^{n+m}$, we let $\pi_P : \mathbb{R}^{n+m} \rightarrow P$ denote the orthogonal projection map onto $P$ and we let $\pi_{P^{\perp}} : \mathbb{R}^{n+m} \rightarrow P^{\perp}$ denote the orthogonal projection map onto the orthogonal complement $P^{\perp}$.  For each $Z \in \mathbb{R}^{n+m}$ and $\rho > 0$, we let  
\begin{align*}
	\mathbf{B}_{\rho}(Z) =\,& \{ X \in \mathbb{R}^{n+m} : |X - Z| < \rho \}, \\
	B_{\rho}(Z,P) =\,& \{ Z + X : X \in P \text{ with } |X| < \rho \}, \\
	\mathbf{C}_{\rho}(Z,P) =\,& \{ Z + X + Y : X \in P \text{ with } |X| < \rho \text{ and } Y \in P^{\perp} \} .
\end{align*}
We say that $P$ is an \emph{oriented $n$-dimensional plane} if we equip $P$ with a simple $n$-vector $\vec P = \tau_1 \wedge \tau_2 \wedge\cdots\wedge \tau_n$, called the \emph{orientation} of $P$, where $\tau_1, \tau_2, \ldots, \tau_n$ is an orthonormal basis for $P$.  In the case $P = \mathbb{R}^n \times \{0\}$, for each $z \in \mathbb{R}^n$ and $\rho > 0$ we let $B_{\rho}(z) = \{ x \in \mathbb{R}^n : |x-z| < \rho \}$ and $\mathbf{C}_{\rho}(z) = B_{\rho}(z) \times \mathbb{R}^m$.  We equip $P = \mathbb{R}^n \times \{0\}$ with the orientation $e_1 \wedge e_2 \wedge\cdots\wedge e_n$.

Let $\mathcal{M}$ be an $n$-dimensional $C^2$-embedded submanifold of $\mathbb{R}^{n+m}$.  For each $X \in \mathcal{M}$, we let ${\rm T}_X \mathcal{M}$ denote the tangent plane to $\mathcal{M}$ at $X$ and we let ${\rm T}_X \mathcal{M}^{\perp}$ denote its orthogonal complement in $\mathbb{R}^{n+m}$.  We let $A_{\mathcal M}$ denote the second fundamental form of $\mathcal{M}$ and $H_{\mathcal M}$ denote the mean curvature vector of $\mathcal{M}$.  We let $\op{dist}_{\mathcal M}(X,Y)$ denote the geodesic distance between two points $X,Y \in \mathcal{M}$.  For each $Z \in \mathcal{M}$ and suitably small $\rho > 0$, we let 
\begin{equation*}
	B^{\mathcal M}_{\rho}(Z) = \{ X \in \mathcal{M} : \op{dist}_{\mathcal M}(X,Z) < \rho \} 
\end{equation*}
be the geodesic ball with center $Y$ and radius $\rho$ in $\mathcal{M}$.

\subsection{Graphs of multi-valued functions}\label{sec:normal graph sec}  Let $q \geq 1$ be an integer and let $\varepsilon = \varepsilon(n,m,q) \in (0,1)$ be suitably small.  Let $\mathcal{M}$ be an $n$-dimensional, oriented, $C^2$-embedded submanifold of an open subset $U \subseteq \mathbb{R}^{n+m}$, $\Omega \subseteq \mathcal{M}$ be an $\mathcal{H}^n$-measurable set and $u : \Omega \rightarrow \mathcal{A}_q(\mathbb{R}^{n+m})$ be a Lipschitz $q$-valued function such that 
\begin{equation}\label{normal graph defn eqn1} 
	\sup_{\mathcal{M}} |A_{\mathcal M}| < \varepsilon, \quad 
    \rho^{-1} \sup_{\Omega} |u| + \op{Lip} u < \varepsilon .
\end{equation}
(Note that the assumption \eqref{normal graph defn eqn1} is not required in the special case that $\mathcal{M}$ is an affine plane and $u$ takes values in the linear subspace orthogonal to $\mathcal{M}$.)  Let $\tau_1,\tau_2,\ldots,\tau_n$ be a positively oriented orthonormal tangent frame field on $\mathcal{M}$ so that $\tau_1(X),\tau_2(X),$ $\ldots,\tau_n(X)$ is an orthonormal basis for ${\rm T}_X \mathcal{M}$ at each $X \in \mathcal{M}$, and $\tau_1 \wedge \tau_2 \wedge\ldots\wedge \tau_n$ is the orientation $n$-vector of $\mathcal{M}$.  Define $F : \Omega \rightarrow \mathcal{A}_q(\mathbb{R}^{n+m})$ by $F(X) = \sum_{l=1}^q \llbracket X + u_l(X) \rrbracket$ for each $X \in \Omega$, where $u(X) = \sum_{l=1}^q \llbracket u_l(X) \rrbracket$.  We associate $u$ with an $n$-dimensional rectifiable current $\mathbf{G}_u$ of $\mathbb{R}^{n+m}$ given by $\mathbf{G}_u 
= F_{\#} (\llbracket \mathcal{M} \rrbracket \llcorner \Omega)$, where the push-forward 
of $F$ is as in~\cite[Definition~1.3]{DeLSpa0}.  In particular, 
\begin{equation*}
	\mathbf{G}_u(\omega) = \int_{\Gamma} \langle \omega,\xi \rangle \,\theta \,d\mathcal{H}^n 
\end{equation*}
for all smooth differential $n$-forms $\omega$ with compact support in $\mathbb{R}^{n+m}$, where 
\begin{gather*}
	\Gamma = \{ X+u_l(X) : X \in \Omega,\, l \in \{1,2,\ldots,q\} \} , \\
	\theta(Y) = \# \{ (X,l) : X \in \Omega,\, l \in \{1,2,\ldots,q\},\, Y = X+u_l(X) \} \text{ for each $Y \in \Gamma$,} \nonumber\\
	\xi(X+u_l(X)) = \dfrac{\bigwedge_{i=1}^n (\tau_i + D_{\tau_i} u_l(X))}{\left| \bigwedge_{i=1}^n (\tau_i + D_{\tau_i} u_l(X)) \right|} 
		\text{ for $\mathcal{H}^n$-a.e.~$X \in \Omega$ and each $l \in \{1,2,\ldots,q\}$,} \nonumber
\end{gather*}
where $u(X) = \sum_{l=1}^q \llbracket u_l(X) \rrbracket$ for each $X \in \Omega$ and $Du(X) = \sum_{l=1}^q \llbracket Du_l(X) \rrbracket$ for $\mathcal{H}^n$-a.e.~$X \in \Omega$ following the convention on~\cite[Subsection~3.5]{KrumWica}.  For each $X \in \mathcal{M}$, $Z \in \mathbb{R}^{n+m}$, and linear map $L : {\rm T}_X \mathcal{M} \rightarrow \mathbb{R}^{n+m}$, we let 
\begin{equation*}
	G_{ij}(X,Z,L) =  \big\langle \tau_i(X) + L(\tau_i(X)), \tau_j(X) + L(\tau_j(X)) \big\rangle , 
\end{equation*}
$[G^{ij}(X,Z,L)] = [G_{ij}(X,Z,L)]^{-1}$, and $G(X,Z,L) = \det [G_{ij}(X,Z,L)]$.  Then 
\begin{equation}\label{normal graph defn eqn5} 
	\mathbf{M}(\mathbf{G}_u) = \int_{\mathcal M} \sum_{l=1}^q \sqrt{G(X,u_l(X),Du_l(X))} \,d\mathcal{H}^n(X) .
\end{equation}
We write $u : \Omega \rightarrow \mathcal{A}_q({\rm T}\mathcal{M}^{\perp})$ to denote that $u : \Omega \rightarrow \mathcal{A}_q(\mathbb{R}^{n+m})$ is a $q$-valued function such that for each $X \in \Omega$ there exist $u_1(X),u_2(X),\ldots,u_q(X) \in {\rm T}_X\mathcal{M}^{\perp}$ such that $u(X) = \sum_{l=1}^q \llbracket u_l(X) \rrbracket$.  In this case, 
\begin{equation}\label{normal graph defn eqn6}
	G_{ij}(X,u_l(X),Du_l(X)) =  \delta_{ij} - 2 \langle A_{ij}(X), u_l(X) \rangle + \langle D_{\tau_i} u_l, D_{\tau_j} u_l(X) \rangle
\end{equation}
for each $X \in \mathcal{M}$, where $A_{ij}(X) = \pi_{{\rm T}_X\mathcal{M}^{\perp}}(D_{\tau_i} \tau_j(X))$ is the second fundamental form.

\subsection{Main hypotheses}  We will need the following notions of oriented tilt excess and height excess of a current $T$ relative to a plane $P$.

\begin{definition}\label{DeLSpa excess and height defn}{\rm (\cite[Definitions~1.1 and 1.2]{DeLSpa2})
Given an $n$-dimensional rectifiable current $T$ of $\mathbf{B}_{\rho}(Y)$ and an $n$-dimensional oriented plane $P$ in $\mathbb{R}^{n+m}$, we define 
\begin{align*}
	\mathbf{E}(T,\mathbf{B}_{\rho}(Y),P) &= \frac{1}{2\omega_n\rho^n} \int_{\mathbf{B}_{\rho}(Y)} |\vec T - \vec P|^2 \,d\|T\|, \\
	\mathbf{h}(T,\mathbf{B}_{\rho}(Y),P) &= \sup_{X_1,X_2 \in \op{spt} T \cap \mathbf{B}_{\rho}(Y)} |\pi_{P^{\perp}}(X_1) - \pi_{P^{\perp}}(X_2)| ,
\end{align*}
where $\vec P$ is the orientation of $P$.  We let 
\begin{equation*}
	\mathbf{E}(T,\mathbf{B}_{\rho}(Y)) = \min_{P \in \mathcal{P}} \mathbf{E}(T,\mathbf{B}_{\rho}(Y),P) ,
\end{equation*}
where the infimum is over the set $\mathcal{P}$ of all $n$-dimensional oriented linear plane $P$ in $\mathbb{R}^{n+m}$, and 
\begin{equation*}
	\mathbf{h}(T,\mathbf{B}_{\rho}(Y)) = \min \{ \mathbf{h}(T,\mathbf{B}_{\rho}(Y),P) : 
		P \in \mathcal{P} \text{ and } \mathbf{E}(T,\mathbf{B}_{\rho}(Y),P) = \mathbf{E}(T,\mathbf{B}_{\rho}(Y)) \} . 
\end{equation*}
}\end{definition}

Throughout the paper, we shall consider parameters $\gamma_1,\beta_2,\delta_2,M_0,N_0,C_e,C_h$ which play a key role in the center manifold construction and satisfy the following requirements:

\begin{hypothesis}\label{DeLSpa constants hypothesis}{\rm
Let $\beta_2 = 4\delta_2 = \min\{1/(2n),\gamma_1/100\}$ where $\gamma_1 = \gamma_1(n,m,q) \in (0,1)$ is as in~\cite[Theorem~1.4]{DeLSpa1}.  Let $M_0 = M_0(n,m,q) \in [4,\infty)$, $N_0 = N_0(n,m,q,M_0) \in \mathbb{N}$, $C_e = C_e(n,m,q,M_0,N_0) \in (0,\infty)$ and $C_h = C_h(n,m,q,M_0,N_0,C_e) \in (0,\infty)$ be suitably large constants and let $\varepsilon_2 = \varepsilon_2(n,m,q,M_0,N_0,C_e,C_h) \in (0,1)$ be a suitably small constant satisfying~\cite[Assumptions~1.8 and 1.9]{DeLSpa2}, and in particular $\sqrt{n} M_0 2^{7-N_0} \leq 1$.
}\end{hypothesis}

Throughout the paper, we shall assume that $T$ is a locally area-minimizing rectifiable current satisfying the requirements of Hypothesis~\ref{center mnfd intro hyp}, which we restate as follows:

\begin{hypothesis}\label{center mnfd hypothesis}{\rm
Let $T$ be an $n$-dimensional area-minimizing rectifiable current of $\mathbf{C}_{6\sqrt{n}}(0)$ such that 
\begin{gather*}
	(\partial T) \llcorner \mathbf{C}_{6\sqrt{n}}(0) = 0, \quad \Theta(T,0) = q, \quad 
	\|T\|(\mathbf{B}_{6\sqrt{n}}(0)) \leq q \omega_n (6\sqrt{n})^n + \varepsilon_2^2, \\ 
	\mathbf{m}_0 = \mathbf{E}(T,\mathbf{B}_{6\sqrt{n}}(0)) = \mathbf{E}(T,P_0,\mathbf{B}_{6\sqrt{n}}(0)) \leq \varepsilon_2^2, \nonumber 
\end{gather*}
where $P_0 = \mathbb{R}^n \times \{0\}$ equipped with the orientation $\vec P_0 = e_1 \wedge e_2 \wedge\cdots\wedge e_n$ (as in~\cite[Assumption~1.3]{DeLSpa2}).  
}\end{hypothesis}

By~\cite[Lemma~1.6]{DeLSpa2}, it follows from Hypothesis~\ref{center mnfd hypothesis} that 
\begin{gather}
	\label{center mnfd hyp3} (\pi_{P_0\#} T) \llcorner \mathbf{C}_{11\sqrt{n}/2}(0) = q \llbracket B_{11\sqrt{n}/2}(0) \rrbracket , \\
	\label{center mnfd hyp4} \mathbf{h}(T,\mathbf{C}_{5\sqrt{n}}(0)) \leq C(n,m,q) \,\mathbf{m}_0^{\frac{1}{2n}} .
\end{gather}

\subsection{General center manifolds and normal approximation maps}\label{sec:center mnfd prelims}   
Consider a cube $L$ of the form 
\begin{equation}\label{DeLSpa cubes eqn}
	L = [a_1,a_1+2\ell] \times [a_2,a_2+2\ell] \times\cdots\times [a_n,a_n+2\ell] \subset \mathbb{R}^n ,
\end{equation}
where $a_1,a_2,\ldots,a_n \in \mathbb{R}^n$ and $\ell > 0$.  We let $x_L = (a_1+\ell,a_2+\ell,\ldots,a_n+\ell)$ denote the center of $L$ and $\ell(L) = \ell$ so that $2\ell(L)$ is the side length of $L$.  For each integer $j \geq 0$, we let $\mathscr{C}^j$ denote the set of all closed cubes $L \subset [-4,4]^n$ which are of the form \eqref{DeLSpa cubes eqn} for some $a_1,a_2,\ldots,a_n \in 2^{1-j} \mathbb{Z}$ and $\ell = 2^{-j}$.  We set $\mathscr{C} = \bigcup_{j=0}^{\infty} \mathscr{C}^j$.  Given $H,L \in \mathscr{C}$, we say that $H$ is an \emph{ancestor} of $L$ if $L \subset H$.

For each cube $L \in \mathscr{C}$, we let $p_L = (x_L,y_L) \in \op{spt} T$ for some $y_L \in \mathbb{R}^m$ (which exists by \eqref{center mnfd hyp3}) and we set $r_L = M_0 \sqrt{n} \,\ell(L)$.  Recalling Definition~\ref{DeLSpa excess and height defn}, we let $P_L$ be an $n$-dimensional oriented linear plane in $\mathbb{R}^{n+m}$ such that 
\begin{equation*}
	\mathbf{E}(T, \mathbf{B}_{64r_L}(p_L),P_L) = \mathbf{E}(T, \mathbf{B}_{64r_L}(p_L)), \quad 
	\mathbf{h}(T, \mathbf{B}_{64r_L}(p_L),P_L) = \mathbf{h}(T, \mathbf{B}_{64r_L}(p_L)) .
\end{equation*}

\begin{definition}\label{DeLSpa whitney defn}{\rm (Whitney decomposition,~\cite[Definition~1.10]{DeLSpa2},~\cite[Section~4.27(21)]{Almgren}) 
Let $T$ be as in Hypothesis~\ref{center mnfd hypothesis}. For each integer $j \geq N_0$, we inductively define collections of cubes $\mathscr{W}^j = \mathscr{W}^j_e \cup \mathscr{W}^j_h \cup \mathscr{W}^j_n \subseteq \mathscr{C}^j$ and $\mathscr{S}^j \subseteq \mathscr{C}^j$ as follows.  For each cube $L \in \mathscr{C}^j$ such that either $j = N_0$ or $j > N_0$ and $L$ does not have an ancestor in $\mathscr{W}^{N_0} \cup \mathscr{W}^{N_0+1} \cup\cdots\cup \mathscr{W}^{j-1}$, we let: 
\begin{enumerate}[itemsep=2mm,topsep=0mm]
	\item[(EX)]  $L \in \mathscr{W}^j_e$ if $\mathbf{E}(T, \mathbf{B}_{64r_L}(p_L)) > C_e \mathbf{m}_0 \ell(L)^{2-2\delta_2}$;
	\item[(HT)]  $L \in \mathscr{W}^j_h$ if $L \not\in \mathscr{W}^j_e$ and $\mathbf{h}(T, \mathbf{B}_{64r_L}(p_L)) > C_h \mathbf{m}_0^{\frac{1}{2n}} \ell(L)^{1+\beta_2}$;
	\item[(NN)]  $L \in \mathscr{W}^j_n$ if $L \not\in \mathscr{W}^j_e \cup \mathscr{W}^j_h$ and $L \cap H \neq \emptyset$ for some $H \in \mathscr{W}^{j-1}$;
	\item[(S)]  $L \in \mathscr{S}^j$ if  $L \not\in \mathscr{W}^j_e \cup \mathscr{W}^j_h \cup \mathscr{W}^j_n$.
\end{enumerate}
We define $\mathscr{W} = \mathscr{W}_e \cup \mathscr{W}_h \cup \mathscr{W}_n$ and $\mathscr{S}$ by $\mathscr{W}_e = \bigcup_{j=N_0}^{\infty} \mathscr{W}^j_e$, $\mathscr{W}_h = \bigcup_{j=N_0}^{\infty} \mathscr{W}^j_h$, $\mathscr{W}_n = \bigcup_{j=N_0}^{\infty} \mathscr{W}^j_n$ and $\mathscr{S} = \bigcup_{j=N_0}^{\infty} \mathscr{S}^j$.  We set 
\begin{equation*}
	\mathbf{\Gamma} = [-4,4]^n \setminus \bigcup_{L \in \mathscr{W}} L = \bigcap_{j=N_0}^{\infty} \bigcup_{L \in \mathscr{S}^j} L . 
\end{equation*}
}\end{definition}

\begin{proposition}\label{domain influence prop} {\rm (}\cite[Corollary~3.2]{DeLSpa2}{\rm )}  
Given any cube $H \in \mathscr{W}_n$ there exists a cube $L \in \mathscr{W}_e$ and a chain $L = L_0, L_1,L_2,\ldots,L_k = H$ such that $L_i \in \mathscr{W}_n$, $L_{i-1} \cap L_i \neq \emptyset$ and $\ell(L_i) = \tfrac{1}{2} \ell(L_{i-1})$ for all $i \in \{1,2,\ldots,k\}$.  We say that $H$ is in the \emph{domain of influence} of $L$.  We have that $H \subset B_{3\sqrt{n}\ell(L)}(x_L)$.  
\end{proposition}

\begin{remark}\label{domain influence rmk}{\rm 
We can enumerate $\mathscr{W}_e = \{L_i\}$ with $\ell(L_i) \geq \ell(L_{i+1})$ for all $i$.  In what follows, we shall take ``$H \in \mathscr{W}_n$ is in the domain of influence of $L_i$'' to also mean that $H$ is not in the domain of influence of $L_j$ for all $j < i$.
}\end{remark}

We shall let $\mathcal{M}$ denote the \textit{center manifold} of $T$ as described in~\cite{DeLSpa2}.  The original construction of the center manifold is in~\cite[Section~4.30]{Almgren}.  Very briefly, to construct $\mathcal{M}$, for each integer $j \geq N_0$ let $\mathscr{P}^j = \bigcup_{i=1}^j \mathscr{W}^i \cup \mathscr{S}^j$.  For each cube $L \in \mathscr{P}^j$, we approximate $T \llcorner \mathbf{C}_{8r_L}(p_L,P_L)$ by the graph of a Lipschitz $q$-valued function $f_L : B_{8r_L}(p_L,P_L) \rightarrow \mathcal{A}_q(P_L^{\perp})$ and smooth the single-valued average of $f_L$ via convolution.  Then we use a smooth partition of unity to glue together the graphs of the smoothed average of $f_L$ (where $L \in \mathscr{P}^j$) to form the graph of a function $\bs{\varphi}_j : (-4,4)^n \to \mathbb{R}^m$.  Since $T$ is locally area-minimizing, the single-valued average of $f_L$ satisfies a variational condition \cite[Proposition~5.2]{DeLSpa2} ~\cite[Proposition~4.3]{Almgren}, from which it follows, through a series of further estimates (cf.~\cite[Sections~4 and 5]{DeLSpa2}, \cite[Sections~4.19, 4.28 and 4.30]{Almgren}), that $\|\bs{\varphi}_j\|_{C^0} \leq C\mathbf{m}_0^{1/(2n)}$ and $\|\bs{\varphi}_j\|_{C^{3,\kappa}} \leq C\mathbf{m}_0^{1/2}$, where $\kappa = \beta_2/4$ and $C = C(n,m,q,M_0,N_0,C_e,C_h) \in (0,\infty)$ is a constant.  The values of $\bs{\varphi}_j$ on Whitney cubes stablizes in the sense that for each $L \in \mathscr{W}^i$ and each $j \geq i+2$, $\bs{\varphi}_j = \bs{\varphi}_{i+2}$ on $L$.  Using this together with the $C^{3,\kappa}$-estimates for the functions $\bs{\varphi}_j$, the whole sequence $(\bs{\varphi}_j)$ converges in the $C^3$-topology, as $j\to\infty$, to some $C^{3,\kappa}$-function $\bs{\varphi} : (-4,4)^n \to \mathbb{R}^m$ such that $\|\bs{\varphi}\|_{C^0} \leq C\mathbf{m}_0^{1/(2n)}$ and $\|\bs{\varphi}\|_{C^{3,\kappa}} \leq C\mathbf{m}_0^{1/2}$, where again $C = C(n,m,q,M_0,N_0,C_e,C_h) \in (0,\infty)$ is a constant.  The center manifold $\mathcal{M}$ is defined by 
$$\mathcal{M} = \{ (x,\bs{\varphi}(x)) : x \in (-4,4)^n \}.$$
By~\cite[Corollary~2.2]{DeLSpa2} (see also~\cite[Theorem~4.33]{Almgren}), the \emph{contact set} $\{ (x,\bs{\varphi}(x)) : x \in \mathbf{\Gamma} \}$ is a subset of the center manifold $\mathcal{M}$ which consists of density $q$ branch points $Z$ of $T$ at which $q \llbracket {\rm T}_Z \mathcal{M} \rrbracket$ is the unique tangent cone to $T$ at $Z$, where ${\rm T}_Z \mathcal{M}$ is the tangent plane to $\mathcal{M}$ at $Z$.  To each cube $L \in \mathscr{W}$, we associate a \emph{Whitney region} 
\begin{equation}\label{Whitney region defn}
    \mathcal{L} = \{ (x,\bs{\varphi}(x)) : x \in H \},
\end{equation}
where $H$ is the cube concentric to $L$ with $\ell(H) = \frac{17}{16} \ell(L)$.

We shall let 
\begin{equation*}
	\mathbf{U} = \{ X+Y : X \in \mathcal{M} \text{ and } Y \in {\rm T}_X \mathcal{M}^{\perp} \cap \mathbf{B}_1(0) \} .
\end{equation*}
Assuming $\varepsilon_2$ is suitably small, there is a well-defined map $\mathbf{p} : \mathbf{U} \rightarrow \mathcal{M}$ such that $\mathbf{p}(X+Y) = X$ for all $X \in \mathcal{M}$ and $Y \in {\rm T}_X \mathcal{M}^{\perp} \cap \mathbf{B}_1(0)$.  In other words, for each $X \in \mathbf{U}$, $\mathbf{p}(X)$ is the unique closest point to $X$ in $\mathcal{M}$.  For each $Z \in \mathcal{M}$ and suitably small $\rho > 0$, we let 
\begin{equation*}
	\mathbf{C}^{\mathcal M}_{\rho}(Z) = \mathbf{p}^{-1}(B^{\mathcal M}_{\rho}(Z)) = \{ X \in \mathbf{U} : \op{dist}_{\mathcal M}(\mathbf{p}(X),Z) < \rho \} , 
\end{equation*}
where we recall that $B^{\mathcal M}_{\rho}(Z)$ is the geodesic ball in $\mathcal{M}$ with center $Z$ and radius $\rho$ and $\op{dist}_{\mathcal M}$ denotes geodesic distance on $\mathcal{M}$.  For each $Z \in \mathcal{M} \cap \mathbf{C}_1(0)$, we let ${\rm exp}_Z : B_3(0,{\rm T}_Z \mathcal{M}) \rightarrow \mathcal{M}$ be the exponential map of $\mathcal{M}$.  Since $\|D\bs{\varphi}\|_{C^{2,\kappa}} \leq C\mathbf{m}_0^{1/2}$,  
\begin{align}\label{exp map est}
	&\sup_{X \in B_3(0,{\rm T}_Z \mathcal{M})} |{\rm exp}_Z(X) - (Z+X)| 
	+ |X| \sup_{X \in B_3(0,{\rm T}_Z \mathcal{M})} |D{\rm exp}_Z(X) - I_{{\rm T}_Z \mathcal{M}}| 
	\\&+ |X|^2 \sup_{X \in B_3(0,{\rm T}_Z \mathcal{M})} |D^2 {\rm exp}_Z(X)| \leq C \mathbf{m}_0^{1/2} |X|^2 , \nonumber 
\end{align}
where $I_{{\rm T}_Z \mathcal{M}}$ is the identity map on ${\rm T}_Z \mathcal{M}$ and $C = C(n,m,q,M_0,N_0,C_e,C_h) \in (0,\infty)$ (see e.g.~\cite[Proposition~A.4]{DeLSpa2}).  

We shall let $N : \mathcal{M} \rightarrow \mathcal{A}_q({\rm T}\mathcal{M}^{\perp})$ denote \emph{normal approximation map} and $\mathcal{K} \subseteq \mathcal{M}$ denote the associated closed set  as constructed in~\cite[Section~6]{DeLSpa2} (cf.~\cite[Theorem~4.33]{Almgren} where $N$ is denoted by $f$ and $\mathcal{K}$ is denoted by $\mathcal{N}^1$).  The map $N$ is a Lipschitz approximation of the current $T$ in the sense that, by~\cite[Theorem~2.4]{DeLSpa2} \cite[Theorem~4.33]{Almgren}, for each $L \in \mathscr{W}$ 
\begin{gather}\label{normal approx est1} 
	T \llcorner \mathbf{p}^{-1}(\mathcal{K}) = \mathbf{G}_N \llcorner \mathbf{p}^{-1}(\mathcal{K}) , \\
	\sup_{\mathcal{L}} |N| \leq C \mathbf{m}_0^{1/(2n)} \ell(L)^{1+\beta_2}, \quad 
	\sup_{\mathcal{L}} |DN| \leq C \mathbf{m}_0^{\gamma_2} \ell(L)^{\gamma_2}, \nonumber \\
	\mathcal{H}^n(\mathcal{L} \setminus \mathcal{K}) + \|T -\mathbf{G}_N\|(\mathbf{p}^{-1}(\mathcal{L})) 
		\leq C \mathbf{m}_0^{1+\gamma_2} \ell(L)^{n+2+\gamma_2}, \nonumber \\
	\int_{\mathcal{L}} |DN|^2 \leq C \mathbf{m}_0 \ell(L)^{n+2-2\delta_2}, \nonumber
\end{gather}
where $\mathcal{L}$ is the Whitney region associated with $L$ (as described above) and $C = C(n,m,q,M_0,N_0,$ $C_e,C_h) \in (0,\infty)$ is a constant.  Furthermore, as described in~\cite[Theorem~2.4]{DeLSpa2}, it follows by applying Young's inequality in (6.13) of~\cite{DeLSpa2} that for each $b > 0$ and each Borel measurable function $\phi : \mathcal{L} \rightarrow [0,\infty)$ the (single-valued) average $N_a$ of $N$ satisfies  
\begin{equation}\label{normal approx est2}
	\int_{\mathcal{L}} \phi |N_a| \leq C \mathbf{m}_0  \ell(L)^{n+3+\beta_2/3} \sup_{\mathcal{L}} \phi + C b^{\frac{1}{1+\gamma_2}} \mathbf{m}_0 \ell(L)^{2+\gamma_2/2} \int_{\mathcal{L}} \phi 
		+ \frac{C}{b} \int_{\mathcal{L}} \phi \,\mathcal{G}(N, q \llbracket N_a \rrbracket)^{2+\gamma_2} , 
\end{equation}
where $C = C(n,m,q,M_0,N_0,C_e,C_h) \in (0,\infty)$ is a constant. By applying Cauchy's inequality instead of Young's inequality in (6.13) of~\cite{DeLSpa2}, we obtain for each $b > 0$ and each Borel measurable function $\phi : \mathcal{L} \rightarrow [0,\infty)$ 
\begin{equation}\label{normal approx est2b}
	\int_{\mathcal{L}} \phi |N_a| \leq C \mathbf{m}_0 \ell(L)^{n+3+\beta_2/3} \sup_{\mathcal{L}} \phi + C b \mathbf{m}_0 \ell(L)^{2+\gamma_2} \int_{\mathcal{L}} \phi  
		+ \frac{C}{b} \int_{\mathcal{L}} \phi\,\mathcal{G}(N, q \llbracket N_a \rrbracket)^2 , 
\end{equation}
where $C = C(n,m,q,M_0,N_0,C_e,C_h) \in (0,\infty)$ is a constant.  (Note that (2.4) of~\cite{DeLSpa2} states that \eqref{normal approx est2} holds true in the special case that $\phi$ is the characteristic function of a Borel set $\mathcal{V} \subset \mathcal{L}$; the same proof works for general non-negative Borel measurable functions $\phi$.  Also, note the difference in the exponent in the first instance of $b$ between \eqref{normal approx est2} and (2.4) of~\cite{DeLSpa2}. 
This has no impact on the overall arguments of~\cite{DeLSpa2} and~\cite{DeLSpa3}, as one can either apply \eqref{normal approx est2} with $b = 1$ or apply \eqref{normal approx est2b}.)  

Following the general conventions of~\cite[Subsection~3.5]{KrumWica}, we write $N(X) = \sum_{l=1}^q \llbracket N_l(X) \rrbracket$ and $N_a(X) = \frac{1}{q} \sum_{l=1}^q N_l(X)$ for each $X \in \mathcal{M}$ and $DN(X) = \sum_{l=1}^q \llbracket DN_l(X) \rrbracket$ for $\mathcal{L}^n$-a.e.~$X \in \mathcal{M}$.

With the above discussion in mind, we shall often assume the following: 

\begin{hypothesis}\label{M and N hypothesis}
Let $T$ and $\mathbf{m}_0$ be as in Hypothesis~\ref{center mnfd hypothesis}.  Let $\mathscr{W}$ be the Whitney decomposition of $T$ as in Definition~\ref{DeLSpa whitney defn}.  Let $\mathcal{M} = \{ (x,\bs{\varphi}(x)) : x \in (-4,4)^n \}$ be the center manifold of $T$ as in~\cite[Theorem~1.17]{DeLSpa2}.  Let $N : \mathcal{M} \rightarrow \mathcal{A}_q({\rm T}\mathcal{M}^{\perp})$ be the normal approximation map of $T$ and $\mathcal{K} \subset \mathcal{M}$ be the corresponding closed set as in~\cite[Theorem~2.4]{DeLSpa2}.  In particular, $N$ and $\mathcal{K}$ satisfy \eqref{normal approx est1}, \eqref{normal approx est2} and \eqref{normal approx est2b} for every $L \in \mathscr{W}$.
\end{hypothesis}

\subsection{Frequency function relative to the center manifold and its conditional monotonicity}

\begin{definition}\label{freq fun defn}{\rm (cf.~\cite[Definition~3.1]{DeLSpa3}, \cite[Section~5.2]{Almgren}) 
Let $\phi : [0,\infty) \rightarrow \mathbb{R}$ is defined by 
\begin{equation}\label{freq phi defn} 
	\phi(s) = \begin{cases} 
		1 &\text{if } 0 \leq s \leq 1/2 \\
		2-2s &\text{if } 1/2 < s < 1 \\
		0 &\text{if } s \geq 1 .
	\end{cases} 
\end{equation}
Let $T, \mathcal{M}, N$ be as in Hypothesis~\ref{M and N hypothesis}.  For each $Z \in B^{\mathcal M}_1(0)$ and $\rho \in (0,2]$ define 
\begin{align*}
	\mathbf{H}_{N,Z}(\rho) &= -\int_{\mathcal{M}} |N|^2 \,\frac{1}{R} \,\phi'(R/\rho) , \\
	\mathbf{D}_{N,Z}(\rho) &= \int_{\mathcal{M}} |DN|^2 \,\phi(R/\rho) , \\
    \mathbf{\Sigma}_{N,Z}(\rho) &= \int |N|^2 \,\phi(R/\rho)
\end{align*}
where $R = R(X) = \op{dist}_{\mathcal M}(X,Z)$ is the geodesic distance of $X$ to $Z$ in $\mathcal{M}$.  We define the frequency function $\mathbf{I}_{N,Z}$ of $N$ at $Z \in B^{\mathcal M}_1(0)$ by 
\begin{equation*}
	\mathbf{I}_{N,Z}(\rho) = \frac{\rho \,\mathbf{D}_{N,Z}(\rho)}{\mathbf{H}_{N,Z}(\rho)}
\end{equation*}
for each $\rho \in (0,2]$ with $\mathbf{H}_{N,Z}(\rho) > 0$. 
}\end{definition}

\begin{definition}\label{crit rad defn}{\rm (Critical radius, cf.~(Stop) and (Go) conditions in~\cite[Subsection~2.1]{DeLSpa3}, \cite[Section~4.32]{Almgren})  Let $T, \mathscr{W}, \mathcal{M}, N$ be as in Hypothesis~\ref{M and N hypothesis}.  For each $Z = (\xi,\bs{\varphi}(\xi)) \in \mathcal{M} \cap \mathbf{C}_1(0)$, the \emph{critical radius} $\overline{R}(Z)$ is the number $\overline{R}(Z) \in [0,2]$ such that for each $\rho \in (\overline{R}(Z),2]$ 
\begin{equation}\label{freq go condition}
	\ell(L) < \frac{\rho}{64\sqrt{n}} \text{ for all } L \in \mathscr{W} \text{ with } L \cap B_{\rho}(\xi) \neq \emptyset 
\end{equation}
and either $\overline{R}(Z) = 0$ or $\overline{R}(Z) > 0$ and 
\begin{equation}\label{freq stop condition}
	\ell(L) = \frac{\overline{R}(Z)}{64\sqrt{n}} \text{ for some } L \in \mathscr{W} \text{ with } L \cap \overline{B_{\overline{R}(Z)}(\xi)} \neq \emptyset .
\end{equation}
}\end{definition}

By \eqref{normal approx est1} and \eqref{freq go condition} \begin{equation}\label{DeLSpa freq mono eqn2}
	\mathbf{D}_{N,Z}(\rho) \leq \int_{B^{\mathcal M}_{\rho}(Z)} |DN|^2 \leq C \mathbf{m}_0 \rho^{n+2-2\delta_2} 
\end{equation}
for all $\rho \in (\overline{R}(Z), 2]$, where $C = C(n,m,q,M_0,N_0,C_e,C_h) \in (0,\infty)$ is a constant.  By~\cite[Lemma~3.6]{DeLSpa3} 
\begin{gather}
	\label{DeLSpa freq mono eqn4} \mathbf{\Sigma}_{N,Z}(\rho) \leq C\rho \,\mathbf{H}_{N,Z}(\rho) + C \rho^2 \mathbf{D}_{N,Z}(\rho), \quad 
		\mathbf{\Sigma}'_{N,Z}(\rho) \leq C \mathbf{H}_{N,Z}(\rho) 
	\end{gather}
for all $\rho \in (0,2]$, where $C = C(n,m,q,M_0,N_0,C_e,C_h) \in (0,\infty)$ are constants.  

By a standard scaling argument~\cite[Subsection~3.1]{DeLSpa3} and \eqref{exp map est} 
\begin{equation}\label{DeLSpa freq mono eqn3}
	\left| \mathbf{H}'_{N,Z}(\rho) - \frac{n-1}{\rho} \,\mathbf{H}_{N,Z}(\rho) + \frac{2}{\rho} \int \sum_{l=1}^q \langle N_l, \partial_{\widehat R} N_l \rangle \,\phi'(R/\rho) \right| \leq C \mathbf{m}_0^{1/2} \mathbf{H}_{N,Z}(\rho) 
\end{equation}
for $\mathcal{L}^1$-a.e.~$\rho \in (0,2]$, where $\partial_{\widehat R}$ denotes the derivative with respect to arclength geodesics in $\mathcal{M}$ starting at $Z$ and $C = C(n,m,q,M_0,N_0,C_e,C_h) \in (0,\infty)$ is a constant.  Using inner and outer variations and estimating error terms as in~\cite[Section~4]{DeLSpa3} using $\|D\bs{\varphi}\|_{C^{2,\kappa}} \leq C\mathbf{m}_0^{1/2}$, \eqref{normal approx est1}, \eqref{normal approx est2} with $b=1$ and~\cite[Lemma~4.5]{DeLSpa3} 
\begin{equation}\label{DeLSpa freq mono eqn16} 
    \left| \mathbf{D}_{N,Z}(\rho) + \frac{1}{\rho} \int \sum_{l=1}^q \langle N_l, \partial_{\widehat R} N_l \rangle \,\phi'(R/\rho) \right| \leq C (\mathbf{\Sigma}_{N,Z}(\rho) + \mathbf{D}_{N,Z}(\rho))^{1+\gamma_3} + C \mathbf{m}_0 \mathbf{\Sigma}_{N,Z}(\rho) 
\end{equation}
for all~$\rho \in (\overline{R}(Z),2]$ and
\begin{align}\label{DeLSpa freq mono eqn17} 
    &\left| \mathbf{D}'_{N,Z}(\rho) - \frac{n-2}{\rho} \,\mathbf{D}_{N,Z}(\rho) + \frac{2}{\rho^2} \int |\partial_{\widehat R} N|^2 \,R \,\phi'(R/\rho) \right| \\&\hspace{15mm} \leq C \mathbf{m}_0 \mathbf{D}_{N,Z}(\rho) + C \mathbf{D}_{N,Z}(\rho)^{\gamma_3} \mathbf{D}'_{N,Z}(\rho) + \frac{C}{\rho} \,\mathbf{D}_{N,Z}(\rho)^{1+\gamma_3} \nonumber 
\end{align}
for $\mathcal{L}^1$-a.e.~$\rho \in (\overline{R}(Z),2]$ with $\mathbf{I}_{N,Z}(\rho) \geq 1$, where $\gamma_3 = \gamma_3(n,m,q) \in (0,1/2]$ and $C = C(n,m,q,M_0,$ $N_0,C_e,C_h) \in (0,\infty)$ are constants.  Note that \eqref{DeLSpa freq mono eqn16} and \eqref{DeLSpa freq mono eqn17} differs from (3.6) and (3.7) of~\cite{DeLSpa3}; here we use the more precise bound $\|A_{\mathcal M}\|_{C^0} \leq C \|D\bs{\varphi}\|_{C^{2,\kappa}} \leq C\mathbf{m}_0^{1/2}$ (rather than $\|A_{\mathcal M}\|_{C^0} \leq C$ for some constant $C$) to estimate the type~2 error terms (as in~\cite[Subsections~3.3 and 4.3]{DeLSpa3}).  Also, in the proof of~\cite[Lemma~4.5]{DeLSpa3}, we apply the argument of~\cite[Lemma~3.6]{DeLSpa3} without using the assumption that $\mathbf{I}_{N,Z}(\rho) \geq 1$; this then requires replacing $\mathbf{D}_{N,Z}(\rho)$ with $\mathbf{\Sigma}_{N,Z}(\rho) + \mathbf{D}_{N,Z}(\rho)$ on the right-hand side of \eqref{DeLSpa freq mono eqn16}.  We shall use \eqref{DeLSpa freq mono eqn16} in this more general form in our proofs below of Lemma~\ref{H positive lemma} (cf.~\eqref{H positive eqn1}) and Lemma~\ref{planar vs alm freq lemma} (cf.~\eqref{planar vs alm freq eqn6}).

Using \eqref{DeLSpa freq mono eqn3}, \eqref{DeLSpa freq mono eqn16} and \eqref{DeLSpa freq mono eqn17} to estimate the derivative of $\mathbf{I}_{N,Z}(\rho)$ at radii $\rho \in (\overline{R}(Z),2]$ for which $\mathbf{I}_{N,Z}(\rho) \geq 1$ gives us the following approximate monotonicity formula (Theorem~\ref{DeLSpa freq mono thm}) for the frequency function $\mathbf{I}_{N,Z}(\cdot)$.  (Note that~\cite{DeLSpa3} considers the frequency function $\mathbf{I}_{N,Z}(\rho)$ and Theorem~\ref{DeLSpa freq mono thm} only in the special case $Z = (0,\bs{\varphi}(0))$, but by the same reasoning Theorem~\ref{DeLSpa freq mono thm} applies for general point $Z \in \mathcal{M}$.)

\begin{theorem}\label{DeLSpa freq mono thm} {\rm (cf.~\cite[Theorem~3.2]{DeLSpa3} \cite[Theorem~5.13]{Almgren})}
Let $M_0,N_0,C_e,C_h,\varepsilon_2$ be as in Hypothesis~\ref{DeLSpa constants hypothesis}.  Let $T, \mathscr{W}, \mathcal{M}, N$ be as in Hypothesis~\ref{M and N hypothesis}.  Let $Z = (\xi,\bs{\varphi}(\xi)) \in \mathcal{M} \cap \mathbf{C}_1(0)$ and let $\overline{R}(Z) \in [0,2]$ be as in Definition~\ref{crit rad defn}.  Then for all $\overline{R}(Z) < \sigma \leq \rho \leq 2$ such that $\mathbf{H}_{N,Z}(\tau) > 0$ for all $\tau \in [\sigma,\rho]$ 
\begin{equation}\label{DeLSpa freq mono concl1}
	\max\{\mathbf{I}_{N,Z}(\sigma),1\} \leq e^{C\mathbf{m}_0^{\gamma_3} \rho^{\gamma_3}} \max\{\mathbf{I}_{N,Z}(\rho),1\}
\end{equation}
and thus 
\begin{equation}\label{DeLSpa freq mono concl2}
	\mathbf{I}_{N,Z}(\sigma) \leq e^{C\mathbf{m}_0^{\gamma_3} \rho^{\gamma_3}} (\mathbf{I}_{N,Z}(\rho) + 1) ,
\end{equation}
where $\gamma_3 = \gamma_3(n,m,q) \in (0,1/2]$ and $C = C(n,m,q,M_0,N_0,C_e,C_h) \in (0,\infty)$ are constants.  
\end{theorem}

\subsection{A critical radius estimate}  We conclude this section with an estimate, Lemma~\ref{relharm critrad2 lemma}, which is a refinement of~\cite[Propositions~3.1 and 3.4]{DeLSpa2} (cf.~\cite[Theorem~4.33(16)(17)]{Almgren}).  In particular, we obtain a lower bound \eqref{relharm critrad2 concl} on the $L^2$-norm of the normal approximation map $N$ at the scale of the critical radius with the constant $C_0$ on the right-hand side independent of $C_e$ and $C_h$.  This will be used to show that Theorem~\ref{graphical_thm} and Theorem~\ref{keyest thm} below hold true without an explicit assumption on the critical radius 
$\overline{R}$ at points in the ball $B^{\mathcal M}_{\rho_{0}}(Z_{0})$ considered in these two results.

\begin{lemma}\label{relharm critrad2 lemma}
Let $M_0,N_0,C_e,C_h$ be as in Hypothesis~\ref{DeLSpa constants hypothesis}. There exist $\varepsilon_2 = \varepsilon_2(n,m,q,M_0,N_0,C_e,C_h) \in (0,1)$ such that the following holds true.  Let $T, \mathscr{W}, \mathcal{M}, N$ be as in Hypothesis~\ref{M and N hypothesis}.  Let $Z = (\xi,\bs{\varphi}(\xi)) \in \mathcal{M} \cap \mathbf{C}_1(0)$ and $\overline{R} = \overline{R}(Z) > 0$ be as in Definition~\ref{crit rad defn}.  Then 
\begin{equation}\label{relharm critrad2 concl} 
	\min\{C_e,C_h\} \cdot \mathbf{m}_0 \overline{R}^2 \leq C_0 \overline{R}^{-n-2+2\delta_2} \int_{B^{\mathcal M}_{\overline{R}}(Z)} |N|^2 , 
\end{equation}
where $\delta_2$ is as in Hypothesis~\ref{DeLSpa constants hypothesis} and $C_0 = C_0(n,m,q,M_0) \in (0,\infty)$ is a constant (independent of $C_e$ and $C_h$; in particular, when we later require $C_e$ and $C_h$ to be sufficiently large, that does not impact the value of $C_0$, nor does this impact the veracity of \eqref{relharm critrad2 concl} as $\overline{R}$ from the center manifold construction becomes smaller).
\end{lemma}

\begin{proof} 
By Definition~\ref{crit rad defn}, there exists $L \in \mathscr{W}$ such that 
\begin{equation}\label{relharm critrad2 eqn1}
	L \cap \overline{B_{\overline{R}}(\xi)} \neq \emptyset \quad\text{and}\quad \ell(L) = \frac{\overline{R}}{64\sqrt{n}}  
\end{equation}
(as in \eqref{freq stop condition}).  Set $\ell = \ell(L)$.  If $L \in \mathscr{W}_h$, then there exists $x_* \in \mathbb{R}^n$ such that $B_{\sqrt{n} \ell/2}(x_*) \subset B_{2\sqrt{n}\ell}(x_L) \cap B_{\overline{R}}(\xi)$.  Assuming $\varepsilon_2$ is suitably small, since $\|D\bs{\varphi}\|_{C^{2,\kappa}} \leq C\mathbf{m}_0^{1/2}$, $\mathcal{M} \cap \mathbf{C}_{\sqrt{n} \ell/4}(x_*) \subset B^{\mathcal M}_{\overline{R}}(Z)$.  By~\cite[Proposition~3.1]{DeLSpa2} 
\begin{equation*}
	\frac{1}{\omega_n \overline{R}^{n+2}} \int_{B^{\mathcal M}_{\overline{R}}(Z)} |N|^2 
		\geq \frac{1}{\omega_n \overline{R}^{n+2}} \int_{\mathcal{M} \cap \mathbf{C}_{\sqrt{n} \ell/4}(x_*)} |N|^2 		
		\geq \frac{1}{4^{4n+14} \,n^2} \,C_h^2 \mathbf{m}_0^{1/n} \overline{R}^{2\beta_2} .
\end{equation*}
Noting that $2\beta_2 \leq 2 - 2\delta_2$, this gives us \eqref{relharm critrad2 concl}.

If $L \in \mathscr{W}_e$, then let $j \geq N_0$ be the integer for which $L \in \mathscr{W}_e^j$.  Let $H \in \mathscr{S}^{j-1}$ and $J \in \mathscr{S}^{j-7}$ be ancestors of $L$.  We argue, via a refinement of the proof of~\cite[Proposition~3.4]{DeLSpa2}, that 
\begin{align}\label{relharm critrad2 eqn7}
	C_e \mathbf{m}_0 \ell^{2-2\delta_2} \leq \frac{C_0}{\ell^{n+2}} \int_{\mathcal{M} \cap \mathbf{C}_{\ell/8}(z,P_H)} |N|^2 
\end{align}
for all $z \in p_J + P_H$ such that $B_{\ell/8}(z,P_H) \subset B_{8r_J}(p_J,P_H)$, where $C_0 = C_0(m,n,q,M_0) \in (0,\infty)$ is a constant.  Let $\lambda = \lambda(n,m,q) \in (0,1)$ be a constant such that $(1+\lambda)^{n+2} < 2^{\delta_2}$.  Let $f : B_{8r_J}(p_J,P_H) \rightarrow \mathcal{A}_q(P_H^{\perp})$ be a Lipschitz approximation of $T$ as in~\cite[Theorem~1.4]{DeLSpa1}.  Let $\overline{\eta} > 0$ and let $w : B_{8r_J}(p_J,P_H) \rightarrow \mathcal{A}_q(P_H^{\perp})$ be a harmonic approximation of $f$ as in~\cite[Theorem~1.6]{DeLSpa1} so that 
\begin{align}\label{relharm critrad2 eqn2}
	&\frac{1}{(8r_J)^{n+2}} \int_{B_{8r_J}(p_H,P_H)} \mathcal{G}(f,w)^2 + \frac{1}{(8r_J)^n} \int_{B_{8r_J}(p_H,P_H)} (|Df| - |Dw|)^2 
	\\&+ \frac{1}{(8r_J)^n} \int_{B_{8r_J}(p_H,P_H)} |Df_{\rm avg} - Dw_{\rm avg}|^2 \leq \overline{\eta} \,\mathbf{E}(T,\mathbf{C}_{32r_J}(p_J,P_H)) . 
	\nonumber 
\end{align}
Here we let $f_{\rm avg}(x) = \frac{1}{q} \sum_{l=1}^q f_l(x)$ and $f_{\rm free}(x) = \sum_{l=1}^q \llbracket f_l(x) - f_{\rm avg}(x) \rrbracket$ for each $x \in B_{\ell/8}(q,P_H)$, where $f(x) = \sum_{l=1}^q \llbracket f_l(x) \rrbracket$.  Similarly, we let $w_{\rm avg}(x) = \frac{1}{q} \sum_{l=1}^q w_l(x)$ and $w_{\rm free}(x) = \sum_{l=1}^q \llbracket w_l(x) - w_{\rm avg}(x) \rrbracket$ for each $x \in B_{\ell/8}(q,P_H)$, where $w(x) = \sum_{l=1}^q \llbracket w_l(x) \rrbracket$.  Let $\rho_* = 32r_H - C \mathbf{m}_0^{1/(2n)} \ell^{1+\beta_2}$ for a suitable constant $C = C(n,m,q,M_0,N_0,C_e,C_h) \in (0,\infty)$.  Recall the first equation on page~551 of~\cite{DeLSpa2} 
\begin{equation}\label{relharm critrad2 eqn3}
	\frac{1}{C} \int_{B_{(1+\lambda)\rho_*}(p_H,P_H)} |Dw|^2 \leq \int_{B_{\ell/8}(z,P_H)} |Dw_{\rm free}|^2 
	\leq \frac{C}{\ell^2} \int_{B_{\ell/8}(z,P_H)} |w_{\rm free}|^2 
\end{equation}
for all $z \in p_J + P_H$ such that $B_{\ell/8}(z,P_H) \subset B_{8r_J}(p_J,P_H)$, where $C = C(n,m,q,\delta_2,M_0) \in (1,\infty)$ is a constant.  By (7.12) in~\cite{DeLSpa2}, $\mathbf{E}(T,\mathbf{C}_{32r_J}(p_J,p_H)) \leq C \mathbf{m}_0 \ell^{2-2\delta_2}$, where $C = C(n,m,q,M_0,N_0,$ $C_e,C_h) > 0$ is a constant.  Arguing as in (7.16) in~\cite{DeLSpa2} except with $A = 0$, by~\cite[Theorem~3.5]{DeLSpa0} and $\mathbf{E}(T,B_{64r_L}(p_L)) > C_e \mathbf{m}_0 \ell^{2-2\delta_2}$ (by condition (EX) of Definition~\ref{DeLSpa whitney defn})
\begin{equation*}
    \frac{1}{2((1+\lambda)\rho_*)^n} \int_{B_{(1+\lambda)\rho_*}(p_H,P_H)} |Df|^2 \geq \mathbf{E}(T,B_{64r_L}(p_L)) - C \mathbf{m}_0^{1+\gamma_1} \ell^{2+\gamma_1/2} , 
\end{equation*}
where $C = C(n,m,q,M_0,N_0,C_e,C_h) \in (0,\infty)$ is a constant.  Hence provided $\varepsilon_2$ is small enough that $C\mathbf{m}_0^{\gamma_1} < C_e/2$,  
\begin{equation}\label{relharm critrad2 eqn4}
    \frac{1}{((1+\lambda)\rho_*)^n} \int_{B_{(1+\lambda)\rho_*}(p_H,P_H)} |Df|^2 \geq \mathbf{E}(T,B_{64r_L}(p_L)) > C_e \mathbf{m}_0 \ell^{2-2\delta_2} .
\end{equation}
Thus assuming $\overline{\eta} = \overline{\eta}(n,m,q,M_0,N_0,C_e,C_h)$ is suitably small, by \eqref{relharm critrad2 eqn2} and \eqref{relharm critrad2 eqn3} 
\begin{equation}\label{relharm critrad2 eqn5}
	\frac{1}{C} \int_{B_{(1+\lambda)\rho_*}(p_H,P_H)} |Df|^2 \leq \int_{B_{\ell/8}(z,P_H)} |Df_{\rm free}|^2 
	\leq \frac{C}{\ell^2} \int_{B_{\ell/8}(z,P_H)} |f_{\rm free}|^2 
\end{equation}
for all $z \in p_J + P_H$ such that $B_{\ell/8}(z,P_H) \subset B_{8r_J}(p_J,P_H)$, where $C = C(n,m,q,\delta_2,M_0) \in (1,\infty)$ is a constant.  

Let $\bs{\varphi}' : B_{8r_J}(p_H,P_H) \rightarrow P_H^{\perp}$ be such that $\op{spt}\mathbf{G}_{\bs{\varphi}'} = \mathcal{M}$ and let $N' : \mathcal{M} \cap \mathbf{C}_{\ell/4}(q,P_H) \rightarrow {\rm T}\mathcal{M}^{\perp}$ be such that $\mathbf{G}_{N'} = \mathbf{G}_f \llcorner \mathbf{C}_{\ell/4}(q,P_H)$.  Notice that by \eqref{normal approx est1}, $N$ and $N'$ agree except on a set of $\mathcal{H}^n$-measure $\leq C\mathbf{m}_0^{1+\gamma_2} \ell^{n+2+\gamma_2}$, where $C = C(n,m,q,M_0,N_0,C_e,C_h) \in (0,\infty)$ is a constant.  By~\cite[Theorem~5.1]{DeLSpa0}, for each $z \in B_{8r_J}(p_H,P_H)$ 
\begin{equation*}
	|N'(z + \bs{\varphi}'(z))| \geq \frac{1}{2\sqrt{q}} \,\mathcal{G}(f(z),q \llbracket \bs{\varphi}'(z) \rrbracket) \geq \frac{1}{4\sqrt{q}} \,|f_{\rm free}(z)| .
\end{equation*}
Therefore, assuming $\varepsilon_2$ is sufficiently small, by \eqref{relharm critrad2 eqn4}, \eqref{relharm critrad2 eqn5}, and \eqref{normal approx est1}
\begin{align}\label{relharm critrad2 eqn6}
	C_e \mathbf{m}_0 \ell^{2-2\delta_2} 
	\leq\,& \frac{C_0}{\ell^{n+2}} \int_{B_{\ell/8}(z,P_H)} |f_{\rm free}|^2 
	\leq \frac{16q C_0}{\ell^{n+2}} \int_{\mathcal{M} \cap \mathbf{C}_{\ell/8}(z,P_H)} |N'|^2
	\\ \leq\,& \frac{16q C_0}{\ell^{n+2}} \int_{\mathcal{M} \cap \mathbf{C}_{\ell/8}(z,P_H)} |N|^2 + C \mathbf{m}_0^{1+\gamma_2+1/n} \ell^{2+2\beta_2+\gamma_2} \nonumber 
\end{align}
for all $z \in P_H$ such that $B_{\ell/8}(z,P_H) \subset B_{8r_J}(p_J,P_H)$, where $C_0 = C_0(m,n,q,M_0) \in (0,\infty)$ and $C = C(n,m,q,M_0,N_0,C_e,C_h) \in (0,\infty)$ are constants.  Thus assuming $\varepsilon_2$ is small enough that $C\mathbf{m}_0^{\gamma_2+1/n} < C_e/2$ (where $C$ is as in \eqref{relharm critrad2 eqn6}), \eqref{relharm critrad2 eqn7} holds true for all $z \in P_H$ such that $B_{\ell/8}(z,P_H) \subset B_{8r_J}(p_J,P_H)$.

Now observe that by \eqref{relharm critrad2 eqn1}, there exists $x_* \in \mathbb{R}^n$ such that $|x_* - x_L| \leq 3 \sqrt{n} \ell/2 < 4\sqrt{n} \ell$ and $B_{\sqrt{n} \ell/2}(x_*) \subset B_{\overline{R}}(\xi)$.  Set $z = p_J + \pi_{P_H}((x_*,\bs{\varphi}(x_*)) - p_J)$.  By the triangle inequality and $|x_* - x_L| \leq 4\sqrt{n} \ell,$
\begin{equation*}
	|x_* - x_J| \leq |x_* - x_L| + |x_L - x_J| \leq 4 \sqrt{n} \ell + 2^7 \sqrt{n} \ell = 132 \sqrt{n} \ell .
\end{equation*}
Assuming $\varepsilon_2$ is sufficiently small, by~\cite[Corollary~2.2]{DeLSpa2} we have that $|y_L - \bs{\varphi}(x_L)| \leq C \mathbf{m}_0^{1/(2n)} \ell^{1+\beta_2} < \ell$ and thus 
\begin{equation*}
	|z - p_J| \leq |(x_*,\bs{\varphi}(x_*)) - p_J| \leq 2 |x_* - x_J| \leq 264 \sqrt{n} \ell < 2^9 M_0 \sqrt{n} \ell = 8r_L .
\end{equation*}
Hence $B_{\ell/8}(z,P_H) \subset B_{8r_J}(p_J,P_H)$ and thus \eqref{relharm critrad2 eqn7} holds true with this choice $z$.  Again assuming $\varepsilon_2$ is sufficiently small, since $\pi_{P_H}(z) = \pi_{P_H}(x_*,\bs{\varphi}(x_*))$ and $B_{\sqrt{n} \ell/2}(x_*) \subset B_{\overline{R}}(\xi)$ 
\begin{equation*}
	\mathcal{M} \cap \mathbf{C}_{\ell/8}(z,P_H) \subset \mathcal{M} \cap \mathbf{C}_{\ell/4}(x_*) \subset B^{\mathcal M}_{\overline{R}}(Z) . 
\end{equation*}
Hence by \eqref{relharm critrad2 eqn7} 
\begin{align}\label{relharm critrad2 eqn8}
	C_e \mathbf{m}_0 \overline{R}^{2-2\delta_2} \leq \frac{C_0}{\overline{R}^{n+2}} \int_{B^{\mathcal M}_{\overline{R}}(Z)} |N|^2 
\end{align}
for some constant $C_0 = C_0(n,m,q,M_0) > 0$.  Therefore, \eqref{relharm critrad2 concl} holds true. 

Finally, if $L \in \mathscr{W}_n$, then by Proposition~\ref{domain influence prop}, $L$ is in the domain of influence of some $M \in \mathscr{W}_e$ and $L \subset B_{3\sqrt{n}\ell(M)}(x_M)$.  Since $L \cap \overline{B_{\overline{R}}(\xi)} \neq \emptyset$, $x_M \in M \cap B_{\overline{R}+3\sqrt{n}\ell(M)}(\xi)$.  Note that $\sqrt{n} \ell(M) \leq \sqrt{n} 2^{-N_0} \leq 1/128$ and $\overline{R} = 64\sqrt{n} \ell \leq 32\sqrt{n} \ell(M) \leq 1/4$ so that $\overline{R} + 3 \sqrt{n} \ell(M) < 1$.  Thus by \eqref{freq go condition} 
\begin{equation*}
	\ell(M) < \frac{\overline{R}+3\sqrt{n}\ell(M)}{64\sqrt{n}}, \quad\text{and therefore } \ell(M) < \frac{\overline{R}}{61\sqrt{n}} . 
\end{equation*}
In particular, $\frac{\overline{R}}{128\sqrt{n}} = \frac{1}{2} \ell(L) \leq \ell(M) < \frac{\overline{R}}{61\sqrt{n}}$.  Since $B_{3\sqrt{n}\ell(M)}(x_M) \cap B_{\overline{R}}(\xi) \neq \emptyset$, there exists $x_* \in \mathbb{R}^n$ such that $|x_* - x_M| \leq 4\sqrt{n} \ell(M)$ and $B_{\sqrt{n}\ell(M)/2}(x_*) \subset B_{\overline{R}}(\xi)$.  Arguing as above with $M$ in place of $L$, we again deduce that \eqref{relharm critrad2 eqn8} holds true, giving us \eqref{relharm critrad2 concl}.
\end{proof}

\section{Branch points of planar frequency 
$\geq 2,$ their effective center manifolds, and unconditional monotonicity of frequency}
\label{sec:quad branch results}  

Let $T$  be as in Hypothesis~\ref{center mnfd hypothesis} and let ${\mathcal M}$ be the center manifold as in Hypothesis~\ref{M and N hypothesis}.  Let $\mathcal{B}^{(\geq 2)}_q$ be as in Definition~\ref{center mnfd branch defn} (with $\mu = 1$).   In this section, we will use planar frequency monotonicity~\cite[Theorem~4.4]{KrumWica} to show that each $Z \in \mathcal{B}^{(\geq 2)}_q$ is contained in the contact set of ${\mathcal M}$ and that for such $Z$, the critical radius $\overline{R}(Z) = 0$ (Lemma~\ref{fast decay lemma});  additionally, the Almgren frequency function satiafies $\mathbf{I}_{N, Z}(\rho) >1$ for all $\rho \in (0, 2]$ (Lemma~\ref{planar vs alm freq lemma}) and thus 
$\mathbf{I}_{N,Z}(\rho)$ satisfies an approximate monotonicity formula for all $\rho \in (0,2]$ unconditionally (Theorem~\ref{freq mono centered thm}), strengthening the conditional monotonicity statement of Theorem~\ref{DeLSpa freq mono thm}.  
In view of these properties, we say that a center manifold ${\mathcal M}$ as in Hypothesis~\ref{M and N hypothesis}  is an \emph{effective center manifold} if the base point $0 \in {\mathcal B}_{q}^{(\geq 2)}$.

As standard consequence of these properties, for $Z \in {\mathcal B}^{(\geq 2)}_{q}$ we obtain growth estimates for $\mathbf{H}_{N,Z}(\rho)$ (Corollary~\ref{doubling centered cor})  and the existence of non-trivial tangent functions for $T$ at $Z$ relative to the center manifold (Subsection~\ref{sec:tangent fun}), followed by a stratification theorem for $\mathcal{B}^{(\geq 2)}_q$ (Theorem~\ref{centered stratification thm}).

\subsection{Some preliminary new estimates for the normal approximation map}\label{sec:normal approx est sec}  In Corollary~\ref{normal approx est lemma}, we are interested in using the estimates \eqref{normal approx est1}, \eqref{normal approx est2} and \eqref{normal approx est2b} for the normal approximation map $N$ on Whitney regions to obtain interior estimates for $N$ 
in geodesic balls $B^{\mathcal M}_{\rho}(Z)$, where $Z \in \mathcal{M}$ and $\rho \in (\overline{R}(Z),2]$ (with $\overline{R}(Z)$ as in Definition~\ref{crit rad defn}).  To this end, as well for later purposes, in Lemma~\ref{normal approx energy lemma} we first prove an energy estimate for $N.$

\begin{lemma}[Energy estimate]\label{normal approx energy lemma} 
Let $M_0,N_0,C_e,C_h$ be as in Hypothesis~\ref{DeLSpa constants hypothesis} and let $\gamma \in (0,1)$.  There exists $\varepsilon_2 = \varepsilon_2(n,m,q,M_0,N_0,C_e,C_h) \in (0,1)$ such that following holds true.  Let $T, \mathscr{W}, \mathcal{M}, N, \mathcal{K}$ be as in Hypothesis~\ref{M and N hypothesis}.  Let $Z \in \mathcal{M} \cap \mathbf{C}_1(0)$ and $\rho \in (\overline{R}(Z),2]$, where $\overline{R}(Z)$ be as in Definition~\ref{crit rad defn}.  Then 
\begin{equation}\label{normal approx energy concl}
	\int_{B^{\mathcal M}_{\gamma\rho}(Z)} |DN|^2 \leq \frac{C}{(1-\gamma)^2 \rho^2} \int_{B^{\mathcal M}_{\rho}(Z)} |N|^2 
\end{equation}
where $C = C(n,m,q,M_0,N_0,C_e,C_h) \in (0,\infty)$ is a constant.
\end{lemma}

\begin{proof} 
Set $\phi_{\rho}(X) = \phi\big(\frac{R(X)}{\rho}\big)$ for each $X \in \mathcal{M}$, where $\phi$ is as in \eqref{freq phi defn} and $R = R(X) = \op{dist}_{\mathcal M}(X,Z)$.  Let $\zeta \in C^1_c(\mathbf{U},\mathbb{R}^{n+m})$ such that $\zeta(X) = \phi_{\rho}^2(\mathbf{p}(X)) \,(X - \mathbf{p}(X))$ for each $X \in \mathbf{U}$ near $\op{spt} T \cup \op{spt} \mathbf{G}_N$.  For $\delta > 0$ sufficiently small and for each $t \in (-\delta,\delta)$, define $\Phi_t : \mathbf{U} \rightarrow \mathbf{U}$ by $\Phi_t(X) = X + t \,\zeta(X)$ for all $X \in \mathbf{U}$.  Notice that $\Phi_{t\#} \mathbf{G}_N$ is the graph of the $q$-valued function $\sum_{l=1}^q \llbracket N_l(X) + t \,\phi_{\rho}^2(X) \,N_l(X) \rrbracket$ for $X \in \mathcal{M}$.  Recalling \eqref{normal graph defn eqn5} and \eqref{normal graph defn eqn6}, 
\begin{align*}
	&\delta\mathbf{G}_N(\zeta) 
	= \frac{d}{dt} \mathbf{M}(\Phi_{t\#} \mathbf{G}_N) \,\Bigg|_{t=0} 
	\\=\,& \int_{\mathcal M} \sum_{l=1}^q \sum_{i,j=1}^n \sqrt{G(X,N_l,DN_l)} \,G^{ij}(X,N_l,DN_l) 
		\,\big( {-\phi_{\rho}^2} \,\langle A_{ij}, N_l \rangle + \langle D_{\tau_i} N_l, D_{\tau_j} (\phi_{\rho}^2 N_l) \rangle \big) , \nonumber 
\end{align*}
where $\delta\mathbf{G}_N$ is the first variation of the integral varifold associated with $\mathbf{G}_N$ (as in~\cite[39.2]{SimonGMT}), $\tau_1,\tau_2,\ldots,\tau_2$ is a positively oriented orthonormal tangent frame field on $\mathcal{M}$, and $G$ and $G^{ij}$ are as in Subsection~\ref{sec:normal graph sec}.  By Taylor's theorem 
\begin{align*}
	\delta\mathbf{G}_N(\zeta) \geq\,& -q \int_{\mathcal M} \phi_{\rho}^2 \,\langle H_{\mathcal M}, N_a \rangle 
		+ \int_{\mathcal M} \sum_{l=1}^q \big( \phi_{\rho}^2 |DN_l|^2 + 2 \phi_{\rho} \,\langle DN_l, N_l \otimes D\phi_{\rho} \rangle \big) 
		\\&- C \int_{\mathcal M} \big( |A_{\mathcal M}| \,|N| + |DN|^2 \big) \big( \phi_{\rho}^2 \,|A_{\mathcal M}| \,|N| + \phi_{\rho}^2 \,|DN|^2 
		+ \phi_{\rho} \,|D\phi_{\rho}| \,|N| \,|DN| \big) \nonumber 
\end{align*}
for some constant $C = C(n,m,q) \in (0,\infty)$.  Thus provided $\varepsilon_2$ is sufficiently small, using Cauchy's inequality, $\|A_{\mathcal M}\| \leq C \mathbf{m}_0^{1/2}$, and \eqref{normal approx est1} 
\begin{equation}\label{normal approx energy eqn1}
	\frac{1}{2} \int_{\mathcal M} \phi_{\rho}^2 |DN|^2 
	\leq q \int_{\mathcal M} \phi_{\rho}^2 \,\langle H_{\mathcal M}, N_a \rangle + \delta\mathbf{G}_N(\zeta) + \frac{C}{\rho^2} \int_{B^{\mathcal M}_{\rho}(Z)} |N|^2 
	 ,  
\end{equation}
where $C = C(n,m,q) \in (0,\infty)$ is a constant. 

To bound the first two terms on the right-hand side of \eqref{normal approx energy eqn1}, we shall apply the strategy from~\cite[Section~4]{DeLSpa3}.  Set $B = \pi_{P_0}(B^{\mathcal M}_{\rho}(Z))$.  Let $\mathcal{T}$ be the set of all cubes $L \in \mathscr{W}_e \cup \mathscr{W}_h$ such that 
\begin{enumerate}[itemsep=3mm,topsep=0mm]
	\item[(i)] $L \in \mathscr{W}_e \cup \mathscr{W}_h$ and $L \cap B \neq \emptyset$ or
	\item[(ii)] $L \in \mathscr{W}_e$ and there exists $H \in \mathscr{W}_n$ in the domain of influence of $L$ (as in Proposition~\ref{domain influence prop} and Remark~\ref{domain influence rmk}) such that $H \cap B \neq \emptyset$.  
\end{enumerate}
For each cube $L \in \mathcal{T}$ set $B^L = B_{s(L)}(x_L)$, where $s(L) = \ell(L)$ if $x_L \in B$ and $s(L) = \op{dist}(x_L,B) + \ell(L)$ otherwise.  By a covering argument (see~\cite[Subsection~4.2]{DeLSpa3}), we can construct a subset $\mathscr{T} \subseteq \mathcal{T}$ such that $\{ B^L : L \in \mathscr{T} \}$ is a pairwise disjoint collection of open balls and we can partition the set of cubes in $\mathscr{W}$ intersecting $B$ as 
\begin{equation*}
	\{ H \in \mathscr{W} : H \cap B \neq \emptyset \} = \bigcup_{L \in \mathscr{T}} \mathscr{W}(L) 
\end{equation*}
(with $\{ \mathscr{W}(L) : L \in \mathscr{T} \}$ pairwise disjoint) where for each $L \in \mathscr{T}$ the set $\mathscr{W}(L)$ consists of cubes $H \in \mathscr{W}$ such that $H \cap B \neq \emptyset$ and either 
\begin{enumerate}[itemsep=3mm,topsep=0mm]
	\item[(i)] $H \in \mathscr{W}_e \cup \mathscr{W}_h$ and $B^H \cap B^L \neq \emptyset$ or 
	\item[(ii)] $H \in \mathscr{W}_n$ and $H$ is in the domain of influence of $J \in \mathscr{W}_e$ such that $B^J \cap B^L \neq \emptyset$. 
\end{enumerate}
For each $L \in \mathscr{T}$, let $B(L)$ be a ball of radius $\ell(L)/4$ such that $B(L) \subset B^L \cap B$ and which lies a distance $\geq \ell(L)/4$ from $\partial B$.  Since $\mathscr{T}$ is a finite or countable set, we can enumerate $\mathscr{T} = \{L_i\}$.  For each $i$, let 
\begin{equation*}
	\ell_i = \ell(L_i), \quad 
	\mathcal{B}^i = \{ (x,\bs{\varphi}(x)) : x \in B(L_i) \}, \quad 
	\mathcal{U}^i = \bigcup_{H \in \mathscr{W}(L_i)} \{ (x,\bs{\varphi}(x)) : x \in H \} \cap B^{\mathcal M}_{\rho}(Z) 
\end{equation*}
and note that $\{\mathcal{B}^i\}$ is a family of pair-wise disjoint subsets of $B^{\mathcal M}_{\rho}(Z)$ and $\bigcup_i \mathcal{U}^i = B^{\mathcal M}_{\rho}(Z)$.  By (4.2) and (4.3) of~\cite{DeLSpa3} 
\begin{equation}\label{DeLSpa freq mono eqn11}
	\frac{\ell_i}{4\rho} \leq \inf_{\mathcal{B}^i} \phi_{\rho} , \quad \sup_{\mathcal{U}^i} \phi_{\rho} \leq C \inf_{\mathcal{B}^i} \phi_{\rho} , 
\end{equation}
where $C = C(n,m,q,M_0,N_0,C_e,C_h) \in (0,\infty)$ is a constant.  
By~\cite[Propositions~3.1 and 3.4]{DeLSpa3}, for each $i$ 
\begin{align}\label{DeLSpa freq mono eqn12}
	\int_{\mathcal{B}^i} |N|^2 \geq\,& c \mathbf{m}_0^{1/n} \ell_i^{n+2+2\beta_2} \quad\text{if } L_i \in \mathscr{W}_h , \\
	\int_{\mathcal{B}^i} |DN|^2 \geq\,& c \mathbf{m}_0 \ell_i^{n+2-2\delta_2} \quad\text{if } L_i \in \mathscr{W}_e , \nonumber 
\end{align}
where $c = c(n,m,q,M_0,N_0,C_e,C_h) > 0$ is a constant.  

Using $8\beta_2 < \gamma_2$, $|H_{\mathcal M}| \leq C\mathbf{m}_0^{1/2}$, \eqref{normal approx est1}, \eqref{normal approx est2} with $b = 1$ and $\phi = \phi_{\rho}^2$, \eqref{DeLSpa freq mono eqn11}, and \eqref{DeLSpa freq mono eqn12} 
\begin{align}\label{normal approx energy eqn2}
	\int_{\mathcal M} \phi_{\rho}^2 \,|H_{\mathcal M}| \,|N_a| 
	\leq\,& C \mathbf{m}_0^{1/2} \sum_i \sup_{\mathcal{U}^i} \phi_{\rho}^2 \left( \mathbf{m}_0 \ell_i^{n+2+\gamma_2/2} 
		+ \int_{\mathcal{U}^i} \phi_{\rho}^2 \,|N|^{2+\gamma_2} \right) 
	\\ \leq\,& C \mathbf{m}_0^{1/2} \rho^{\gamma_2/4} \sum_i \left( \inf_{\mathcal{B}^i} \phi_{\rho}^2 \int_{\mathcal{B}^i} (|N|^2 + |DN|^2) 
		+ \int_{\mathcal{U}^i} \phi_{\rho}^2 \,|N|^2 \right) \nonumber
	\\ \leq\,& C \mathbf{m}_0^{1/2} \rho^{\gamma_2/4} \int_{B^{\mathcal M}_{\rho}(Z)} \phi_{\rho}^2 \,(|N|^2 + |DN|^2) , \nonumber 
\end{align}
where $C = C(n,m,q,M_0,N_0,C_e,C_h) \in (0,\infty)$ are constants.  Notice that 
\begin{equation*}
	|D\zeta| \leq C \left( \phi_{\rho}^2 + \phi_{\rho} \,\frac{|X - \mathbf{p}(X)|}{\rho} \right) 
\end{equation*}
on $\op{spt} T \cup \op{spt} \mathbf{G}_N$ and that by~\cite[Corollary~2.2]{DeLSpa2} and \eqref{normal approx est1}, $|X - \mathbf{p}(X)| \leq C \mathbf{m}_0 \ell_i^{1+\beta_2} < \ell_i$ for each $X \in (\op{spt}T \cup \op{spt}\mathbf{G}_N) \cap \mathbf{p}^{-1}(\mathcal{U}^i)$, where $C = C(n,m,q,M_0,N_0,C_e,C_h) \in (0,\infty)$ is a constant.  Since $T$ is area-stationary, by \eqref{normal approx est1}, \eqref{DeLSpa freq mono eqn11}, and \eqref{DeLSpa freq mono eqn12} 
\begin{align}\label{normal approx energy eqn3}
	|\delta\mathbf{G}_N(\zeta)| 
	=\,& |\delta\mathbf{G}_N(\zeta) - \delta T(\zeta)| 
	\\ \leq\,& C \sum_i \sup_{(\op{spt}T \cup \op{spt}\mathbf{G}_N) \cap \mathbf{p}^{-1}(\mathcal{U}^i)} |D\zeta| 
		\,\|\mathbf{G}_N - T\|(\mathbf{p}^{-1}(\mathcal{U}^i)) \nonumber
	\\ \leq\,& C \sum_i \left( \sup_{\mathcal{U}^i} \phi_{\rho}^2 + \frac{\ell_i}{\rho} \sup_{\mathcal{U}^i} \phi_{\rho} \right) 
		\mathbf{m}_0^{1+\gamma_2} \ell_i^{n+2+\gamma_2} \nonumber
	\\ \leq\,& C (\mathbf{m}_0 \rho)^{3\gamma_2/4} \sum_i \inf_{\mathcal{B}^i} \phi_{\rho}^2 \int_{\mathcal{B}^i} (|N|^2 + |DN|^2) \nonumber
	\\ \leq\,& C (\mathbf{m}_0 \rho)^{3\gamma_2/4} \int_{B^{\mathcal M}_{\rho}(Z)} \phi_{\rho}^2 \,(|N|^2 + |DN|^2) , \nonumber 
\end{align}
where $C = C(n,m,q,M_0,N_0,C_e,C_h) \in (0,\infty)$ are constants.  Bounding the right-hand side of \eqref{normal approx energy eqn1} using \eqref{normal approx energy eqn2} and \eqref{normal approx energy eqn3} 
\begin{equation}\label{normal approx energy eqn4}
	\frac{1}{2} \int_{\mathcal M} \phi_{\rho}^2 |DN|^2 
	\leq \frac{C}{\rho^2} \int_{B^{\mathcal M}_{\rho}(Z)} |N|^2 
		+ C (\mathbf{m}_0 \rho)^{\gamma_2/4} \int_{B^{\mathcal M}_{\rho}(Z)} \phi_{\rho}^2 |DN|^2 , 
\end{equation}
where $C = C(n,m,q,M_0,N_0,C_e,C_h) \in (0,\infty)$ are constants.  Thus provided $\varepsilon_2$ is small enough that $C \mathbf{m}_0^{\gamma_2} \leq C \varepsilon_2^{2\gamma_2} < 1/4$ (where $C$ is as in \eqref{normal approx energy eqn4}), we obtain 
\begin{equation}\label{normal approx energy concl2}
	\int_{\mathcal M} \phi_{\rho}^2 |DN|^2 
	\leq \frac{C}{\rho^2} \int_{B^{\mathcal M}_{\rho}(Z)} |N|^2 , 
\end{equation}
where $C = C(n,m,q,M_0,N_0,C_e,C_h) \in (0,\infty)$ are constants.  Using $\phi_{\rho} \geq 2 (1-\gamma)$ on $B^{\mathcal M}_{\rho}(Z)$ gives us \eqref{normal approx energy concl}.
\end{proof}

\begin{corollary}\label{normal approx est lemma} 
Let $M_0,N_0,C_e,C_h,\varepsilon_2$ be as in Hypothesis~\ref{DeLSpa constants hypothesis} and Lemma~\ref{normal approx energy lemma}.  Let $\gamma \in (0,1)$.  Let $T, \mathscr{W}, \mathcal{M}, N, \mathcal{K}$ be as in Hypothesis~\ref{M and N hypothesis}.  Let $Z = (\xi,\bs{\varphi}(\xi)) \in \mathcal{M} \cap \mathbf{C}_1(0)$ and $\rho \in (\overline{R}(Z),2]$, where $\overline{R}(Z)$ is as in Definition~\ref{crit rad defn}.  Then 
\begin{gather}
	\label{normal approx est3} \mathcal{H}^n(B^{\mathcal M}_{\gamma\rho}(Z) \setminus \mathcal{K}) 
		+ \|T - \mathbf{G}_N\|(\mathbf{C}^{\mathcal M}_{\gamma\rho}(Z)) 
		\leq C \mathbf{m}_0^{3\gamma_2/4} \rho^{3\gamma_2/4-2} \int_{B^{\mathcal M}_{\rho}(Z)} |N|^2 , \\
	\label{normal approx est4} \frac{1}{\rho} \int_{B^{\mathcal M}_{\gamma\rho}(Z)} |N_a| 
		\leq C \rho^{\gamma_2/4-2} \int_{B^{\mathcal M}_{\rho}(Z)} |N|^2 ,  
\end{gather}
where $C = C(n,m,q,M_0,N_0,C_e,C_h,\gamma) \in (0,\infty)$ is a constant.
\end{corollary}

\begin{proof} 
By arguing as in \eqref{normal approx energy eqn2} 
\begin{equation}\label{normal approx energy eqn8}
	\int_{\mathcal M} \phi_{\rho}^2 \,|N_a| 
	\leq C \rho^{\gamma_2/4} \int_{B^{\mathcal M}_{\rho}(Z)} \phi_{\rho}^2 \,(|N|^2 + |DN|^2) ,  
\end{equation}
where $C = C(n,m,q,M_0,N_0,C_e,C_h) \in (0,\infty)$ is a constant.  By arguing as in \eqref{normal approx energy eqn3} using \eqref{normal approx est1}, \eqref{DeLSpa freq mono eqn11} and \eqref{DeLSpa freq mono eqn12}
\begin{align}\label{normal approx energy eqn9}
	\int_{\mathbf{C}^{\mathcal M}_{\rho}(Z)} \phi_{\rho}^2(\mathbf{p}(X)) \,d\|\mathbf{G}_N - T\|(X)
    \leq\,& C \sum_i \sup_{\mathcal{U}^i} \phi_{\rho}^2 \, \|\mathbf{G}_N - T\|(\mathbf{p}^{-1}(\mathcal{U}^i)) 
	\\ \leq\,& C \sum_i \sup_{\mathcal{U}^i} \phi_{\rho}^2 \cdot \mathbf{m}_0^{1+\gamma_2} \ell_i^{n+2+\gamma_2} \nonumber
	\\ \leq\,& C (\mathbf{m}_0 \rho)^{3\gamma_2/4} \int_{B^{\mathcal M}_{\rho}(Z)} \phi_{\rho}^2 \,(|N|^2 + |DN|^2) , \nonumber 
\end{align}
and by similar reasoning 
\begin{equation}\label{normal approx energy eqn10}
	\int_{B^{\mathcal M}_{\rho}(Z) \setminus \mathcal{K}} \phi_{\rho}^2 \,d\mathcal{H}^n
    \leq C (\mathbf{m}_0 \rho)^{3\gamma_2/4} \int_{B^{\mathcal M}_{\rho}(Z)} \phi_{\rho}^2 \,(|N|^2 + |DN|^2) , 
\end{equation}
where $C = C(n,m,q,M_0,N_0,C_e,C_h) \in (0,\infty)$ are constants.  Now \eqref{normal approx est3} and \eqref{normal approx est4} follow from \eqref{normal approx energy eqn8}, \eqref{normal approx energy eqn9} and \eqref{normal approx energy eqn10} by noting that $\phi_{\rho} \geq 2(1 - \gamma)$ on $B^{\mathcal M}_{\gamma\rho}$ and using \eqref{normal approx energy concl2} to bound the right-hand sides.
\end{proof}

\subsection{Unconditional monotonicity of frequency} 
Let $T$ be as in Hypothesis~\ref{center mnfd hypothesis}.  We define the following planar coarse excesses.

\begin{definition}\label{height excess defn}{\rm 
Let $Z \in \mathbf{C}_{11\sqrt{n}/2}(0)$, $\rho > 0$, and $P \subset \mathbb{R}^{n+m}$ be an $n$-dimensional plane.  We define 
\begin{equation*}
	\widehat{E}(T,P,\mathbf{B}_{\rho}(Z)) = \frac{1}{\omega_n\rho^{n+2}} \int_{\mathbf{B}_{\rho}(Z)} \op{dist}^2(X,Z+P) \,d\|T\|(X) .
\end{equation*}
For each $n$-dimensional 
plane $\widetilde{P} \subset \mathbb{R}^{n+m}$, we define 
\begin{equation*}
	\widehat{E}(T,P,\mathbf{C}_{\rho}(Z,\widetilde{P})) = \frac{1}{\omega_n\rho^{n+2}} \int_{\mathbf{C}_{\rho}(Z,\widetilde{P})} \op{dist}^2(X,Z+P) \,d\|T\|(X) .
\end{equation*}
}\end{definition}

\begin{remark}{\rm 
$\widehat{E}(T,P,\mathbf{B}_{\rho}(Z))$ is the $L^2$-distance of $T$ to the affine plane $Z+P$ in $\mathbf{B}_{\rho}(Z)$, in constrast with the oriented tilt-excess $\mathbf{E}(T,\mathbf{B}_{\rho}(Z),P)$ of $T$ relative to $P$ in $\mathbf{B}_{\rho}(Z)$.  Unlike~\cite{KrumWica} and~\cite{KrumWicb}, we use $\widehat{E}(T,P,\mathbf{B}_{\rho}(Z))$ (instead of $E(T,P,\mathbf{B}_{\rho}(Z))$) to denote the coarse excess of $T$ relative to a plane $P$ in $\mathbf{B}_{\rho}(Z)$.
}\end{remark}

Now we will show in Lemma~\ref{fast decay lemma} that $\mathcal{B}^{(\geq 2)}_q$ is contained in the contact set and at each $Z \in \mathcal{B}^{(\geq 2)}_q$ the (Go) condition (as in \eqref{freq go condition}) holds true for all $\rho \in (0,2]$.

\begin{lemma}\label{fast decay lemma}
There are $M_0,N_0,C_e,C_h,\varepsilon_2$ satisfying Hypothesis~\ref{DeLSpa constants hypothesis} such that the following holds true.  Let $T, \mathscr{W}, \mathcal{M}$ be as in Hypothesis~\ref{M and N hypothesis}.  Let $\mathcal{B}^{(\geq 2)}_q$ be as in Definition~\ref{center mnfd branch defn}.  Then for each $Z = (\xi,\zeta) \in \mathcal{B}^{(\geq 2)}_q$ and $\rho \in (0,2]$ 
\begin{equation}\label{fast decay concl}
	\ell(H) < \frac{\rho}{64\sqrt{n}} \text{ for all } H \in \mathscr{W} \text{ with } H \cap B_{\rho}(\xi) \neq \emptyset .
\end{equation}
Moreover, each $Z \in \mathcal{B}^{(\geq 2)}_q$ is in the contact set $\{ (x,\bs{\varphi}(x)) : x \in \mathbf{\Gamma} \}$.
\end{lemma}

\begin{proof} 
Let $Z = (\xi,\zeta) \in \mathcal{B}^{(\geq 2)}_q$ and $L \in \mathscr{W} \cup \mathscr{S}$ such that $|x_L - \xi| \leq 67\sqrt{n} \,\ell(L)$.  Provided $\varepsilon_2$ is sufficiently small, by~\cite[Proposition~4.1(v)]{DeLSpa2}, $p_L = (x_L,y_L)$ satisfies $|y_L - \zeta| \leq C \mathbf{m}_0^{1/(2n)} \ell(L)$ and thus $|p_L - Z| \leq 68\sqrt{n}\,\ell(L)$.  Hence provided $M_0 \geq 68$, $|p_L - Z| \leq M_0 \sqrt{n} \,\ell(L) = r_L$.  

By~\cite[Theorem~2.12]{KrumWicb}, there exists an $n$-dimensional current $\mathbf{P} = \sum_{l=1}^q \llbracket \mathbb{R}^n \times \{a_l\} \rrbracket$, where $a_1,a_2,\ldots,a_q \in \mathbb{R}^m$, such that 
\begin{equation}\label{fast decay eqn1}
	\sup_{X \in \op{spt} T \cap \mathbf{C}_{5\sqrt{n}}(0)} \op{dist}(X,\op{spt}\mathbf{P}) \leq C_0 \mathbf{m}_0^{1/2} , 
\end{equation}
where $C_0 = C_0(n,m,q) \in (0,\infty)$ is a constant.  Let 
\begin{equation*}
	\{ y \in \mathbb{R}^m : \op{dist}(X,\op{spt}\mathbf{P}) \leq C_0 \mathbf{m}_0^{1/2} \} = \bigcup_{i=1}^N (\mathbb{R}^n \times U_i) 
\end{equation*}
for a collection of pairwise disjoint connected open sets $\{U_i\}$ in $\mathbb{R}^m$ (where $C_0$ is as in \eqref{fast decay eqn1}).  Then  
\begin{equation*}
	T \llcorner \mathbf{C}_{5\sqrt{n}}(0) = \sum_{i=1}^N T_i
\end{equation*}
where $T_i = T \llcorner \mathbf{C}_{5\sqrt{n}}(0) \times U_i$ are locally area-minimizing rectifiable currents of $\mathbf{C}_{5\sqrt{n}}(0)$ with $(\partial T_i) \llcorner \mathbf{C}_{5\sqrt{n}}(0) = 0$.  By the constancy theorem and~\cite[Remark~2.18]{KrumWicb}, $\pi_{P_0\#} T_i = q_i \llbracket B_{5\sqrt{n}}(0) \rrbracket$ where $q_i \geq 1$ are integers such that $\sum_{i=1}^N q_i = q$.  Since $0 \in \op{spt} T$ with $\Theta(T,0) = q$, there exists $i \in \{1,2,\ldots,N\}$ such that $0 \in \op{spt} T_i$ and thus by (3.4) from~\cite{KrumWicb} 
\begin{equation*}
	q = \Theta(T,0) = \Theta(T_i,0) \leq \omega_n^{-1} \|T_i\|(\mathbf{B}_1(0)) \leq q_i + C \mathbf{m}_0 , 
\end{equation*}
where $C = C(n) \in (0,\infty)$ is a constant.  Therefore, assuming $\varepsilon_2$ is sufficiently small, we must have that $N = 1$.  Thus by $0 \in \op{spt} T$ and \eqref{fast decay eqn1}
\begin{equation}\label{fast decay eqn3}
	\sup_{X \in \op{spt} T \cap \mathbf{C}_{5\sqrt{n}}(0)} \op{dist}(X,P_0) \leq 2q C_0 \mathbf{m}_0^{1/2} 
\end{equation}
(where $C_0$ is as in \eqref{fast decay eqn1}).  
 
In light of~\cite[Lemma~3.6]{KrumWicc} (applied to $\eta_{0,5\sqrt{n}\#} T$), we can apply~\cite[Corollary~5.3]{KrumWica} and $\mathcal{N}_{T,{\rm Pl}}(Z) \geq 2$ to find $\theta = \theta(n,m,q) \in (0,1)$ such that 
\begin{align}\label{fast decay eqn4}
	\widehat{E}(T,P_Z,\mathbf{C}_{\sigma}(Z,P_Z)) 
	&\leq e^{C (\mathcal{N}_{T,{\rm Pl}}(Z) + 2)} \sigma^{\mathcal{N}_{T,{\rm Pl}}(Z)-1} \widehat{E}(T,P_Z,\mathbf{C}_1(Z,P_Z))
	\\&\leq e^{4C} (e^C \theta)^{\mathcal{N}_{T,{\rm Pl}}(Z)-2} \sigma \,\widehat{E}(T,P_Z,\mathbf{C}_1(Z,P_Z)) \nonumber 
	\\&\leq C \sigma \widehat{E}(T,P_Z,\mathbf{C}_1(Z,P_Z)) \nonumber
\end{align}
for all $\sigma \in (0,\theta]$, where $C = C(n,m,q)$ are constants.  By the triangle inequality, $\Theta(T,Z) \geq q$,~\cite[Lemma~3.4]{KrumWica}, \eqref{fast decay eqn3}, and \eqref{fast decay eqn4} 
\begin{align}\label{fast decay eqn5}
	\|\pi_{P_Z} - \pi_{P_0}\|^2 
	\leq\,& \frac{2}{q \omega_n \sigma^n} \int_{G_n(\mathbf{B}_{\sigma}(Z))} \|\pi_S - \pi_{P_Z}\|^2 \,d|T|(X,S) 
		\\&+ \frac{2}{q \omega_n \sigma^n} \int_{G_n(\mathbf{B}_{\sigma}(Z))} \|\pi_S - \pi_{P_0}\|^2 \,d|T|(X,S) \nonumber 
	\\ \leq\,& C \widehat{E}(T,P_Z,\mathbf{B}_{2\sigma}(Z))^2 + C \sigma^{-n} \mathbf{m}_0 \nonumber 
	\\ \leq\,& C \sigma^2 \widehat{E}(T,P_Z,\mathbf{C}_1(Z,P_Z))^2 + C \sigma^{-n} \mathbf{m}_0 \nonumber 
\end{align}
for all $\sigma \in (0,\theta/2]$, where $C = C(n,m,q) \in (0,\infty)$ are constants.  Thus using $|{\rm dist}(X,Z+P_Z) - {\rm dist}(X,Z+P_0)| \leq \|\pi_{P_Z} - \pi_{P_0}\| \,|X-Z|$ for all $X \in \mathbb{R}^{n+m}$, \eqref{fast decay eqn3}, and \eqref{fast decay eqn5} 
\begin{align*}
	\widehat{E}(T,P_Z,\mathbf{C}_1(Z,P_Z))^2 
	\leq\,& 2 \widehat{E}(T,P_0,\mathbf{C}_1(Z,P_Z))^2 + C \|\pi_{P_Z} - \pi_{P_0}\|^2 \nonumber 
	\\ \leq\,& C \sigma^2 \widehat{E}(T,P_Z,\mathbf{C}_1(Z,P_Z))^2 + C \sigma^{-n} \mathbf{m}_0 
\end{align*}
for all $\sigma \in (0,\theta]$, where $C = C(n,m,q) \in (0,\infty)$ are constants.  Choosing $\sigma = \sigma(n,m,q) \in (0,\theta]$ sufficiently small, 
\begin{equation}\label{fast decay eqn6}
	\widehat{E}(T,P_Z,\mathbf{C}_1(Z,P_Z)) \leq C \mathbf{m}_0^{1/2} ,
\end{equation}
where $C = C(n,m,q) \in (0,\infty)$ is a constant.  By \eqref{fast decay eqn4} and \eqref{fast decay eqn6}  
\begin{equation}\label{fast decay eqn7}
	\widehat{E}(T,P_Z,\mathbf{C}_{\sigma}(Z,P_Z)) \leq C \mathbf{m}_0^{1/2} \sigma 
\end{equation}
for all $\sigma \in (0,\theta]$, where $C = C(n,m,q)$ is a constant.  By \eqref{fast decay eqn3} and \eqref{fast decay eqn5}, we in fact have that \eqref{fast decay eqn7} holds true for all $\sigma \in (0,4\sqrt{n}]$.  Recalling that $\ell(L) \leq 2^{-N_0}$ and assuming $N_0$ is large enough that $130 M_0 \sqrt{n} \,2^{-N_0} \leq 1$, $130r_L \leq 130M_0 \sqrt{n} \,2^{-N_0} \leq 1$.  By \eqref{fast decay eqn7} with $\sigma = 130r_L$ 
\begin{equation}\label{fast decay eqn8}
	\widehat{E}(T,P_Z,\mathbf{B}_{130r_L}(Z)) \leq C \mathbf{m}_0^{1/2} r_L , 
\end{equation}
where $C = C(n,m,q)$ is a constant.  By~\cite[Lemma~3.6]{KrumWica} and \eqref{fast decay eqn8} 
\begin{equation}\label{fast decay eqn9}
	\mathbf{E}(T,\mathbf{B}_{65r_L}(Z),P_Z) 
	\leq C \widehat{E}(T,P_Z,\mathbf{B}_{130r_L}(Z))^2 
	\leq C \mathbf{m}_0 r_L^2 
\end{equation}
and thus using $|p_L - Z| \leq r_L$ 
\begin{equation}\label{fast decay eqn10}
	\mathbf{E}(T,\mathbf{B}_{64r_L}(p_L)) 
	\leq \mathbf{E}(T,\mathbf{B}_{64r_L}(p_L),P_Z) 
	\leq 2^n \mathbf{E}(T,\mathbf{B}_{65r_L}(Z),P_Z) 
	\leq C \mathbf{m}_0 r_L^2 , 
\end{equation}
where $C = C(n,m,q)$ is a constant.  By the triangle inequality, $|p_L - Z| \leq r_L$, \eqref{fast decay eqn9} and \eqref{fast decay eqn10}
\begin{equation}\label{fast decay eqn11}
	|\vec P_Z - \vec P_L| \leq C \mathbf{E}(T,\mathbf{B}_{65r_L}(Z),P_Z)^{1/2} + C \mathbf{E}(T,\mathbf{B}_{64r_L}(p_L))^{1/2} 
	\leq C \mathbf{m}_0^{1/2} r_L ,
\end{equation}
where $C = C(n,m,q)$ are constants.  By $|p_L - Z| \leq r_L$,~\cite[Lemma~3.5]{KrumWica}, \eqref{fast decay eqn8}, and \eqref{fast decay eqn11}  
\begin{align}\label{fast decay eqn12}
	\mathbf{h}(T,\mathbf{B}_{64r_L}(p_L)) 
	\leq\,& \mathbf{h}(T,\mathbf{B}_{65r_L}(Z),P_Z) + C r_L |\vec P_Z - \vec P_L| 
	\\ \leq\,& C r_L \widehat{E}(T,P_Z,\mathbf{B}_{130r_L}(Z)) + C r_L |\vec P_Z - \vec P_L|  
	\leq C \mathbf{m}_0^{1/2} r_L^2 , \nonumber 
\end{align}

Now to see \eqref{fast decay concl}, suppose to the contrary that there exists $Z = (\xi,\zeta) \in \mathcal{B}^{(\geq 2)}_q$, $\rho \in (0,2]$ and $H \in \mathscr{W}$ such that 
\begin{equation}\label{fast decay eqn13}
	H \cap B_{\rho}(\xi) \neq \emptyset, \quad \ell(H) \geq \frac{\rho}{64\sqrt{n}} . 
\end{equation}
We claim that there exists a cube $L \in \mathscr{W}_e \cup \mathscr{W}_h$ such that $|x_L - \xi| \leq 67\sqrt{n}\,\ell(L)$.  If $H \in \mathscr{W}_e \cup \mathscr{W}_h$, we simply set $L = H$ and note that by \eqref{fast decay eqn13} 
\begin{equation*}
	|x_L - \xi| \leq \rho + \sqrt{n}\,\ell(L) \leq 65\sqrt{n}\,\ell(L) .
\end{equation*}
If instead $H \in \mathscr{W}_n$, then by Proposition~\ref{domain influence prop} we know that $H$ is in the domain of influence of some cube $L \in \mathscr{W}_e$ and in particular $\ell(H) < \ell(L)$ and $H \subset B_{3\sqrt{n}\,\ell(L)}(x_L)$.  
Hence by \eqref{fast decay eqn13} 
\begin{equation*}
	\ell(L) \geq \frac{\rho}{64\sqrt{n}}, \quad |x_L - \xi| \leq \rho + 3\sqrt{n}\,\ell(L) \leq 67\sqrt{n}\,\ell(L) . 
\end{equation*}
Since $L$ satisfies $|x_L - \xi| \leq 67\sqrt{n}\,\ell(L)$, we must have that \eqref{fast decay eqn10} and \eqref{fast decay eqn12} hold true.  However, assuming $C_e$ and $C_h$ are sufficiently large, \eqref{fast decay eqn10} and \eqref{fast decay eqn12} contradicts $L \in \mathscr{W}_e \cup \mathscr{W}_h$.

Let $Z = (\xi,\zeta) \in \mathcal{B}^{(\geq 2)}_q$.  If there were a cube $H \in \mathscr{W}$ such that $\xi \in H$, then by $H \cap B_{\rho}(\xi) \neq \emptyset$ for \emph{all} $\rho \in (0,2]$, which since $\ell(H) > 0$ contradicts \eqref{fast decay concl}.  Therefore, $\xi \in \mathbf{\Gamma}$, or equivalently $Z$ is in the contact set $\{ (x,\bs{\varphi}(x)) : x \in \mathbf{\Gamma} \}$.
\end{proof}

\begin{lemma}\label{H positive lemma}
Let $M_0,N_0,C_e,C_h$ be as in Hypothesis~\ref{DeLSpa constants hypothesis} and Lemma~\ref{fast decay lemma}.  There exists $\varepsilon_2 = \varepsilon_2(n,m,q,M_0,N_0,C_e,C_h) \in (0,1)$ such that the following holds true.  Let $T, \mathcal{M}, N, \mathcal{K}$ be as in Hypothesis~\ref{M and N hypothesis}.  Let $\mathcal{B}^{(\geq 2)}_q$ be as in Definition~\ref{center mnfd branch defn}.  For every $Z \in \mathcal{B}^{(\geq 2)}_q$ and $\rho \in (0,2]$, $\mathbf{H}_{N,Z}(\rho) > 0$.
\end{lemma}

\begin{proof}
Let $C \in (0,\infty)$ denote constants depending only on $n,m,q,M_0,N_0,C_e,C_h$.  Suppose to the contrary that $\mathbf{H}_{N,Z}(\rho) = 0$ for some $Z \in \mathcal{B}^{(\geq 2)}_q$ and $\rho \in (0,2]$.  Then $N = q \llbracket 0 \rrbracket$ on $B^{\mathcal M}_{\rho}(Z) \setminus B^{\mathcal M}_{\rho/2}(Z)$.  Moreover, by \eqref{DeLSpa freq mono eqn4} we have that $\mathbf{\Sigma}_{N,Z}(\rho) \leq C \rho^2 \mathbf{D}_{N,Z}(\rho)$.  Note that by Lemma~\ref{fast decay lemma}, $\overline{R}(Z) = 0$.  By \eqref{DeLSpa freq mono eqn16} and $\mathbf{\Sigma}_{N,Z}(\rho) \leq C \rho^2 \mathbf{D}_{N,Z}(\rho)$ 
\begin{equation}\label{H positive eqn1}
	\mathbf{D}_{N,Z}(\rho) \leq C \mathbf{D}_{N,Z}(\rho)^{1+\gamma_3} + C \mathbf{m}_0 \rho^2 \mathbf{D}_{N,Z}(\rho) 
\end{equation}
for some constant $\gamma_3 = \gamma_3(n,m,q) \in (0,1/2]$, which provided $\varepsilon_2$ is sufficiently small gives us $\mathbf{D}_{N,Z}(\rho) = 0$.  Since $\mathbf{\Sigma}_{N,Z}(\rho) \leq C \rho^2 \mathbf{D}_{N,Z}(\rho) = 0$, $\mathbf{\Sigma}_{N,Z}(\rho) = 0$ and thus $N = q \llbracket 0 \rrbracket$ on $B^{\mathcal M}_{\rho}(Z)$.  By \eqref{normal approx est3}, $\mathcal{K} = B^{\mathcal M}_{\rho/2}(Z)$ and thus $T \llcorner \mathbf{C}^{\mathcal M}_{\rho/2}(Z) = q \llbracket B^{\mathcal M}_{\rho/2}(Z) \rrbracket$, contradicting $Z \in \mathcal{B}^{(\geq 2)}_q \subseteq \op{sing} T$.
\end{proof}

In Lemma~\ref{planar vs alm freq lemma} below, for each 
$Z \in {\mathcal B}_{q}^{(\geq 2)}$ we establish a lower bound on Almgren's frequency function $\mathbf{I}_{N,Z}(\rho)$ in terms of the planar frequency $\mathcal{N}_{T,{\rm Pl}}(Z)$;  consequently, $\mathbf{I}_{N,Z}(\rho) > 1$ for all $Z \in \mathcal{B}_q^{(\geq 2)}$ and $\rho \in (0,2]$, and thus in Theorem~\ref{DeLSpa freq mono thm} we can replace $\max\{\mathbf{I}_{N,Z}(\rho),1\}$ with $\mathbf{I}_{N,Z}(\rho).$

\begin{lemma}\label{planar vs alm freq lemma}
Let $M_0,N_0,C_e,C_h$ be as in Hypothesis~\ref{DeLSpa constants hypothesis} and Lemma~\ref{fast decay lemma}.  There exists $\varepsilon_2 = \varepsilon_2(n,m,q,M_0,N_0,C_e,C_h) \in (0,1)$ such that the following holds true.  Let $T, \mathcal{M}, N$ be as in Hypothesis~\ref{M and N hypothesis}.  Let $\mathcal{B}^{(\geq 2)}_q$ be as in Definition~\ref{center mnfd branch defn}.  Then for each $Z \in \mathcal{B}^{(\geq 2)}_q$ 
\begin{equation*}
    \liminf_{\rho\downarrow 0} \mathbf{I}_{N,Z}(\rho) \geq \mathcal{N}_{T,{\rm Pl}}(Z) 
\end{equation*}
(where $\mathcal{N}_{T,{\rm Pl}}(Z)$ is the planar frequency of $T$ at $Z$ as in~\cite[Definition~4.1 and Definition~5.2]{KrumWica}).  
\end{lemma}

\begin{proof} 
To simplify notation, let $\mathbf{H} = \mathbf{H}_{N,Z}$, $\mathbf{D} = \mathbf{D}_{N,Z}$, $\mathbf{\Sigma} = \mathbf{\Sigma}_{N,Z}$ and $\mathbf{I} = \mathbf{I}_{N,Z}$.  Let $C \in (0,\infty)$ denote constants depending only on $n,m,q,M_0,N_0,C_e,C_h,\mathcal{N}_{T,{\rm Pl}}(Z)$.  

Suppose to the contrary that for some $\delta > 0$ and some radius $\rho_0 \in (0,1]$ (which can be chosen to be arbitrarily small), $\mathbf{I}(\rho_0) \leq \mathcal{N}_{T,{\rm Pl}}(Z) - 2\delta$.  That way, provided $\rho_0$ is sufficiently small, by Theorem~\ref{DeLSpa freq mono thm}
\begin{equation}\label{planar vs alm freq eqn1}
    \mathbf{I}(\rho) \leq \mathcal{N}_{T,{\rm Pl}}(Z) - \delta
\end{equation}
for all $\rho \in (0,\rho_0]$.
Thus \eqref{DeLSpa freq mono eqn4} gives us that 
\begin{equation}\label{planar vs alm freq eqn2} 
    \mathbf{\Sigma}(\rho) \leq C\rho \,\mathbf{H}(\rho) 
\end{equation}
for all $\rho \in (0,\rho_0]$.  
Let $\gamma_3 = \gamma_3(n,m,q) \in (0,1)$ be as in \eqref{DeLSpa freq mono eqn16}.  By \eqref{normal approx est1} and \eqref{fast decay concl} 
\begin{equation*}
	\mathbf{H}(\rho) \leq C \rho^{n-1} \sup_{B^{\mathcal M}_{\rho}(Z)} |N|^2 \leq C \mathbf{m}_0^{1/n} \rho^{n+1+2\beta_2} 
\end{equation*}
for all $\rho \in (0,\rho_0]$.  Thus by \eqref{planar vs alm freq eqn2} 
\begin{equation}\label{planar vs alm freq eqn4}
	\mathbf{\Sigma}(\rho)^{1+\gamma_3} \leq (C\rho \,\mathbf{H}(\rho))^{1+\gamma_3} \leq C \mathbf{m}_0^{\gamma_3/n} \rho^{1+\gamma_3 (n+2+2\beta_2)} \mathbf{H}(\rho) 
\end{equation}
for all $\rho \in (0,\rho_0]$.  By \eqref{DeLSpa freq mono eqn2} and \eqref{planar vs alm freq eqn1}
\begin{equation}\label{planar vs alm freq eqn5}
	\mathbf{D}(\rho)^{1+\gamma_3} \leq C \mathbf{m}_0^{\gamma_3} \rho^{\gamma_3 (n+2-2\delta_2)} \mathbf{D}(\rho) 
		\leq \tfrac{7}{4} C \mathbf{m}_0^{\gamma_3} \rho^{\gamma_3 (n+2-2\delta_2)-1} \mathbf{H}(\rho) 
\end{equation}
for all $\rho \in (0,\rho_0]$.  Provided $\gamma_3 \leq 1/(n+2-2\delta_2)$, bounding the right-hand side of \eqref{DeLSpa freq mono eqn16} using \eqref{planar vs alm freq eqn2}, \eqref{planar vs alm freq eqn4} and \eqref{planar vs alm freq eqn5} gives us 
\begin{equation}\label{planar vs alm freq eqn6}
	\left| \mathbf{D}(\rho) + \frac{1}{\rho} \int \phi'(R/\rho) \sum_{l=1}^q \langle N_l, \partial_{\widehat R} N_l \rangle \right| 
		\leq C \mathbf{m}_0^{\gamma_3/n} \rho^{\gamma_3 (n+2-2\delta_2)-1} \mathbf{H}(\rho) 
\end{equation}
for all $\rho \in (0,\rho_0]$.  Combining \eqref{DeLSpa freq mono eqn3} and \eqref{planar vs alm freq eqn6}, 
\begin{equation*}
	\big| \mathbf{H}'(\rho) - \tfrac{n-1}{\rho} \,\mathbf{H}(\rho) - 2 \mathbf{D}(\rho) \big| \leq C \mathbf{m}_0^{\gamma_3/n} \rho^{\gamma_3 (n+2-2\delta_2)-1} \mathbf{H}(\rho)
\end{equation*}
for all $\rho \in (0,\rho_0]$.  Thus provided $\rho_0$ is sufficiently small, by \eqref{planar vs alm freq eqn1}
\begin{equation}\label{planar vs alm freq eqn7}
	\frac{\mathbf{H}'(\rho)}{\mathbf{H}(\rho)} - \frac{n-1}{\rho} \leq \frac{2 \mathbf{I}(\rho)}{\rho} + C \mathbf{m}_0^{\gamma_3/n} \rho^{n\gamma_3-1} \leq \frac{2\mathcal{N}_{T,{\rm Pl}}(Z)-\delta}{\rho} 
\end{equation}
for each $\rho \in (0,\rho_0]$.  By integrating \eqref{planar vs alm freq eqn7} over $[\sigma,\rho]$ 
\begin{equation}\label{planar vs alm freq eqn8}
	\Big(\frac{\sigma}{\rho}\Big)^{n-1+2\mathcal{N}_{T,{\rm Pl}}(Z)-\delta} \mathbf{H}(\rho) \leq \mathbf{H}(\sigma)
\end{equation}
whenever $0 < \sigma \leq \rho \leq \rho_0$.  Noting that $\frac{2}{\rho} \leq \frac{-1}{R} \,\phi'(R/\rho) \leq \frac{4}{\rho}$ on $B^{\mathcal M}_{\rho}(Z) \setminus B^{\mathcal M}_{\rho/2}(Z)$ and summing \eqref{planar vs alm freq eqn8} with $2^{-k} \sigma,2^{-k} \rho$ in place of $\sigma,\rho$ for $k = 0,1,2,\ldots$ gives us 
\begin{equation}\label{planar vs alm freq eqn9}
	\frac{1}{16} \Big(\frac{\sigma}{\rho}\Big)^{n+2\mathcal{N}_{T,{\rm Pl}}(Z)-\delta} \int_{B^{\mathcal M}_{\rho}(Z)} |N|^2 \leq \int_{B^{\mathcal M}_{\sigma}(Z)} |N|^2 
\end{equation}
whenever $0 < \sigma \leq \rho \leq \rho_0$ (see the proof of~\cite[Corollary~5.3]{KrumWica} for further details).

On the other hand, by \eqref{normal approx est1}, Corollary~\ref{normal approx est lemma}, and \eqref{planar vs alm freq eqn9} 
\begin{align*}
    \int_{B^{\mathcal M}_{\rho}(Z)} |N|^2 
    \leq\,& q \int_{B^{\mathcal M}_{\rho}(Z)} |N_a|^2 + \int_{B^{\mathcal M}_{\rho}(Z)} |N_f|^2
    \\ \leq\,& C \mathbf{m}_0^{1/(2n)} \rho^{\gamma_2/4+\beta_2} \int_{B^{\mathcal M}_{2\rho}(Z)} |N|^2 + \int_{B^{\mathcal M}_{\rho}(Z)} |N_f|^2
    \\ \leq\,& 2^{n+2\mathcal{N}_{T,{\rm Pl}}(Z)+5} C \mathbf{m}_0^{1/(2n)} \rho^{\gamma_2/4+\beta_2} \int_{B^{\mathcal M}_{\rho}(Z)} |N|^2 + \int_{B^{\mathcal M}_{\rho}(Z)} |N_f|^2 
\end{align*}
for all $\rho \in (0,\rho_0/2]$, where $N_a(X) = \frac{1}{q} \sum_{l=1}^q N_l(X)$ is the single-valued average of $N$ and $N_f(X) = \sum_{l=1}^q \llbracket N_l(X) - N_a(X) \rrbracket$ is the average-free part of $N$.  Hence provided $\rho_0$ is sufficiently small
\begin{equation*}
    \int_{B^{\mathcal M}_{\rho}(Z)} |N|^2 
    \leq C \int_{B^{\mathcal M}_{\rho}(Z)} |N_f|^2 , 
\end{equation*}
for all $\rho \in (0,\rho_0/2]$.  Recall that $\mathcal{M} = \op{graph}\bs{\varphi}$ where $\|\bs{\varphi}\|_{C^3} \leq C\mathbf{m}_0^{1/(2n)}$ and that by \eqref{fast decay eqn5} and \eqref{fast decay eqn7}, $\op{dist}_{\mathcal H}(P_Z \cap \mathbf{B}_1(0),P_0 \cap \mathbf{B}_1(0)) \leq C\mathbf{m}_0^{1/2}$ (where $P_0 = \mathbb{R}^n \times \{0\}$).  Hence we can express $(Z+P_Z) \cap \mathbf{C}^{\mathcal M}_3(0) = \{ X + \lambda(Z) : X \in B^{\mathcal M}_3(0) \}$ for some $C^3$-function $\lambda : B^{\mathcal M}_1(0) \rightarrow {\rm T}\mathcal{M}^{\perp}$ with $\|\lambda\|_{C^3(B^{\mathcal M}_1(0))} \leq C\mathbf{m}_0^{1/(2n)}$.  Moreover, $|N_l(X) - \lambda(X)| \leq 2 \op{dist}(X+N_l(X),Z+P_Z)$ for all $X \in B^{\mathcal M}_1(0)$ and $l \in \{1,2,\ldots,q\}$.  Thus 
\begin{align*}
    \int_{B^{\mathcal M}_{\rho}(Z)} |N|^2 \leq\,& C \int_{B^{\mathcal M}_{\rho}(Z)} |N_f|^2 \leq C \int_{B^{\mathcal M}_{\rho}(Z)} \mathcal{G}(N,q\llbracket \lambda \rrbracket)^2 \\ \leq\,& C \int_{B^{\mathcal M}_{\rho}(Z)} \sum_{l=1}^q \op{dist}^2(X+N_l(X),Z+P_Z) , 
\end{align*}
for all $\rho \in (0,\rho_0/2]$.  By the area formula and again using \eqref{normal approx est1}, Corollary~\ref{normal approx est lemma}, and \eqref{planar vs alm freq eqn9} 
\begin{equation*}
    \int_{B^{\mathcal M}_{\rho}(Z)} |N|^2 \leq C \int_{\mathbf{C}^{\mathcal M}_{\rho}(Z)} \op{dist}^2(X,P_Z) \,d\|T\|(X) + 2^{n+2\mathcal{N}_{T,{\rm Pl}}(Z)+5} C \mathbf{m}_0^{3\gamma_2/4} \rho^{3\gamma_2/4} \int_{B^{\mathcal M}_{\rho}(Z)} |N|^2 
\end{equation*}
for all $\rho \in (0,\rho_0/2]$.  Therefore, provided $\varepsilon_2$ is sufficiently small 
\begin{equation}\label{planar vs alm freq eqn10}
    \int_{B^{\mathcal M}_{\rho}(Z)} |N|^2 \leq C \int_{\mathbf{C}^{\mathcal M}_{\rho}(Z)} \op{dist}^2(X,Z+P_Z) \,d\|T\|(X) 
\end{equation}
for all $\rho \in (0,\rho_0/2]$.  By comparing both sides of \eqref{planar vs alm freq eqn10} using the growth estimates \eqref{planar vs alm freq eqn9} and~\cite[Corollary~5.3]{KrumWica} 
\begin{align*}
    \frac{1}{4} \Big(\frac{\sigma}{\rho}\Big)^{n+2\mathcal{N}_{T,{\rm Pl}}(Z)-\delta} \int_{B^{\mathcal M}_{\rho}(Z)} |N|^2 \leq\,& \int_{B^{\mathcal M}_{\sigma}(Z)} |N|^2 \leq C \int_{\mathbf{C}^{\mathcal M}_{\sigma}(Z)} \op{dist}^2(X,Z+P_Z) \,d\|T\|(X) \\ \leq\,& C \Big(\frac{\sigma}{\rho}\Big)^{n+2\mathcal{N}_{T,{\rm Pl}}(Z)} \int_{\mathbf{C}^{\mathcal M}_{\rho}(Z)} \op{dist}^2(X,Z+P_Z) \,d\|T\|(X)  
\end{align*}
whenever $0 < \sigma \leq \rho \leq \rho_0/2$. 
Dividing by $\sigma^{n+2\mathcal{N}_{T,{\rm Pl}}(Z)-\delta}$ and letting $\sigma \rightarrow 0^+$ gives us that $N = q \llbracket 0 \rrbracket$ on $B^{\mathcal M}_{\rho_0/2}(Z)$, contradicting Lemma~\ref{H positive lemma}.
\end{proof}

\begin{theorem}\label{freq mono centered thm}
Let $M_0,N_0,C_e,C_h$ be as in Hypothesis~\ref{DeLSpa constants hypothesis} and Lemma~\ref{fast decay lemma}.  There exists $\varepsilon_2 = \varepsilon_2(n,m,q,M_0,N_0,C_e,C_h) \in (0,1)$ such that the following holds true.  Let $T, \mathcal{M}, N$ be as in Hypothesis~\ref{M and N hypothesis}.  Let $\mathcal{B}^{(\geq 2)}_q$ be as in Definition~\ref{center mnfd branch defn}.  Then 
\begin{equation*}
	\mathbf{I}_{N,Z}(\sigma) \leq e^{C\mathbf{m}_0^{\gamma_3} \rho^{\gamma_3}} \mathbf{I}_{N,Z}(\rho)
\end{equation*}
for all $Z \in \mathcal{B}^{(\geq 2)}_q$ and $0 < \sigma \leq \rho \leq 2$, where $\gamma_3 \in (0,1/2]$ and $C = C(n,m,q,M_0,N_0,C_e,C_h) \in (0,\infty)$ are constants.
\end{theorem}

\begin{proof}
Provided $\varepsilon_2$ is sufficiently small, by Lemma~\ref{planar vs alm freq lemma} and Theorem~\ref{DeLSpa freq mono thm}, $\mathbf{I}_{N,Z}(\rho) > 1$ for all $Z \in \mathcal{B}^{(\geq 2)}_q$ and $\rho \in (0,2]$.  Hence $\mathbf{I}_{N,Z}(\rho) = \max\{\mathbf{I}_{N,Z}(\rho),1\}$ for all $Z \in \mathcal{B}^{(\geq 2)}_q$ and $\rho \in (0,2]$.  The  conclusion now follows directly from Lemma~\ref{fast decay lemma} and Theorem~\ref{DeLSpa freq mono thm}. 
\end{proof}

For the remainder of this section, we fix the values of the parameters $M_0,N_0,C_e,C_h,\varepsilon_2$.

\begin{hypothesis}\label{DeLSpa consts freq hypothesis}{\rm
Let $M_0,N_0,C_e,C_h,\varepsilon_2$ be constants satisfying the requirements of Hypothesis~\ref{DeLSpa constants hypothesis} such that Lemma~\ref{normal approx energy lemma}, Corollary~\ref{normal approx est lemma}, Lemma~\ref{fast decay lemma}, Lemma~\ref{H positive lemma}, Lemma~\ref{planar vs alm freq lemma} and Theorem~\ref{freq mono centered thm} all hold true.
}\end{hypothesis}

It now follows from Theorem~\ref{freq mono centered thm} that if $T, \mathcal{M}, N$ are as in Hypothesis~\ref{M and N hypothesis} and if $\mathcal{B}^{(\geq 2)}_q$ is as in Definition~\ref{center mnfd branch defn}, then the \emph{frequency} 
\begin{equation}\label{centered freq}
	\mathbf{I}_{N,Z}(0^+) = \lim_{\rho \rightarrow 0^+} \mathbf{I}_{N,Z}(\rho)
\end{equation}
exists for all $Z \in \mathcal{B}^{(\geq 2)}_q$.  Note that by Lemma~\ref{planar vs alm freq lemma}, $\mathbf{I}_{N,Z}(0^+) \geq \mathcal{N}_{T,{\rm Pl}}(Z)$ for all $Z \in \mathcal{B}^{(\geq 2)}_q$.  As a standard consequence of Theorem~\ref{freq mono centered thm}, the frequency $\mathbf{I}_{N,Z}(0^+)$ is an upper semi-continuous function of $Z \in \mathcal{B}_q^{(\geq 2)}$.

\begin{corollary}\label{centered semicont1 cor}
Let $M_0,N_0,C_e,C_h,\varepsilon_2$ be as in Hypothesis~\ref{DeLSpa consts freq hypothesis}.  
Let $T, \mathcal{M}, N$ be as in Hypothesis~\ref{M and N hypothesis}.  Let $\mathcal{B}^{(\geq 2)}_q$ be as in Definition~\ref{center mnfd branch defn}.  For $\nu = 1,2,3,\ldots$, let $Z_{\nu},Z \in \mathcal{B}^{(\geq 2)}_q$ and $\rho_{\nu},\rho \in (0,2]$ such that $Z_{\nu} \rightarrow Z$ and $\rho_{\nu} \rightarrow \rho$ as $\nu\rightarrow\infty$.  Then 
\begin{equation}\label{centered semicont1 concl1}
	\mathbf{I}_{N,Z}(\rho) = \lim_{\nu \rightarrow \infty} \mathbf{I}_{N,Z_{\nu}}(\rho_{\nu}) . 
\end{equation}
Moreover, 
\begin{equation}\label{centered semicont1 concl2}
	\mathbf{I}_{N,Z}(0^+) \geq \limsup_{k \rightarrow \infty} \mathbf{I}_{N,Z_{\nu}}(0^+) . 
\end{equation}
\end{corollary}

\begin{proof}
To see \eqref{centered semicont1 concl1}, set $R_{\nu}(X) = \op{dist}_{\mathcal{M}}(X,Z_{\nu})$ and $R(X) = \op{dist}_{\mathcal{M}}(X,Z)$ for each $X \in \mathcal{M}$.  By the triangle inequality $|R_{\nu}(X) - R(X)| \leq \op{dist}_{\mathcal{M}}(Z_{\nu},Z)$ and thus $R_{\nu} \rightarrow R$ uniformly on $B^{\mathcal M}_3(0)$ as $\nu \rightarrow \infty$.  In particular, $\phi(R_{\nu}/\rho_{\nu}) \rightarrow \phi(R/\rho)$ uniformly on $B^{\mathcal M}_3(0)$ and $\phi'(R_{\nu}/\rho_{\nu}) \rightarrow \phi'(R/\rho)$ point-wise $\mathcal{H}^n$-a.e.~on $B^{\mathcal M}_3(0)$ as $\nu \rightarrow \infty$.  It then readily follows using the dominated convergence theorem that 
\begin{equation*}
	\mathbf{H}_{N,Z_{\nu}}(\rho) = \lim_{\nu \rightarrow \infty} \mathbf{H}_{N,Z_{\nu}}(\rho_{\nu}) > 0, \quad 
	\mathbf{D}_{N,Z_{\nu}}(\rho) = \lim_{\nu \rightarrow \infty} \mathbf{D}_{N,Z_{\nu}}(\rho_{\nu})
\end{equation*}
and consequently \eqref{centered semicont1 concl1} holds true.  By Theorem~\ref{freq mono centered thm} and \eqref{centered semicont1 concl1} (with $\rho_{\nu} = \rho$), 
\begin{equation}\label{centered semicont1 eqn1}
	e^{C \mathbf{m}_0^{\gamma_3} \rho^{\gamma_3}} \mathbf{I}_{N,Z}(\rho) 
	= \lim_{\nu \rightarrow \infty} e^{C \mathbf{m}_0^{\gamma_3} \rho^{\gamma_3}} \mathbf{I}_{N,Z_{\nu}}(\rho) 
	\geq \limsup_{\nu \rightarrow \infty} \mathbf{I}_{N,Z_{\nu}}(0^+) , 
\end{equation}
for each $\rho \in (0,2]$, where $\gamma_3 = \gamma_3(n,m,q) \in (0,1/2]$ and $C$ are as in Theorem~\ref{freq mono centered thm}.  Letting $\rho \rightarrow 0^+$ in \eqref{centered semicont1 eqn1} gives us \eqref{centered semicont1 concl2}.
\end{proof}

As another standard consequence of Theorem~\ref{freq mono centered thm}, we obtain growth estimates for $\mathbf{H}_{N,Z}(\rho)$.

\begin{corollary}\label{doubling centered cor}
Let $M_0,N_0,C_e,C_h,\varepsilon_2$ be as in Hypothesis~\ref{DeLSpa consts freq hypothesis}.  
Let $T, \mathcal{M}, N$ be as in Hypothesis~\ref{M and N hypothesis}.  Let $\mathcal{B}^{(\geq 2)}_q$ be as in Definition~\ref{center mnfd branch defn} and let $Z \in \mathcal{B}^{(\geq 2)}_q$.  Then 
\begin{align}\label{doubling centered concl1}
	\Big(\frac{\sigma}{\rho}\Big)^{n-1+2e^{C \mathbf{m}_0^{\gamma_3} \rho^{\gamma_3}} \mathbf{I}_{N,Z}(\rho)} &\mathbf{H}_{N,Z}(\rho) 
	\leq \mathbf{H}_{N,Z}(\sigma) \\&\leq e^{2C \mathbf{I}_{N,Z}(\sigma_0) \, \mathbf{m}_0^{\gamma_3}\rho^{\gamma_3}/\gamma_3} 
			\Big(\frac{\sigma}{\rho}\Big)^{n-1+2\mathbf{I}_{N,Z}(\sigma_0)} \mathbf{H}_{N,Z}(\rho) \nonumber 
\end{align}
for all $0 < \sigma_0 \leq \sigma \leq \rho \leq 2$ and 
\begin{align}\label{doubling centered concl2}
	\frac{1}{4} \Big(\frac{\sigma}{\rho}\Big)^{n+2e^{C \mathbf{m}_0^{\gamma_3} \rho^{\gamma_3}} \mathbf{I}_{N,Z}(\rho)} 
		&\int_{B^{\mathcal M}_{\rho}(Z)} |N|^2 
	\leq \int_{B^{\mathcal M}_{\sigma}(Z)} |N|^2 \\&\leq 4 e^{2C \mathbf{I}_{N,Z}(0^+) \, \mathbf{m}_0^{\gamma_3}\rho^{\gamma_3}/\gamma_3} 
			\Big(\frac{\sigma}{\rho}\Big)^{n+2\mathbf{I}_{N,Z}(0^+)} \int_{B^{\mathcal M}_{\rho}(Z)} |N|^2 \nonumber 
\end{align}
for all $0 < \sigma \leq \rho \leq 2$, where $\gamma_3 = \gamma_3(n,m,q) \in (0,1/2]$ and $C = C(n,m,q,M_0,N_0,C_e,C_h) \in (0,\infty)$ are constants.
\end{corollary}

\begin{proof}
Let $0 < \sigma \leq \rho \leq 2$.  
Combining \eqref{DeLSpa freq mono eqn3} and \eqref{DeLSpa freq mono eqn16} 
\begin{align*}
	\big| \mathbf{H}'_{N,Z}(\tau) - \tfrac{n-1}{\tau} \,\mathbf{H}_{N,Z}(\tau) - 2 \mathbf{D}_{N,Z}(\tau) \big| 
	\leq\,& C \mathbf{m}_0^{1/2} \mathbf{H}_{N,Z}(\tau) + C (\mathbf{\Sigma}_{N,Z}(\tau) + \mathbf{D}_{N,Z}(\tau))^{1+\gamma_3} \\&+ C \mathbf{m}_0 \mathbf{\Sigma}_{N,Z}(\tau) 
\end{align*}
for all $\tau \in [\sigma,\rho]$, where $\gamma_3 = \gamma_3(n,m,q) \in (0,1/2]$ and $C = C(n,m,q,M_0,N_0,C_e,C_h) \in (0,\infty)$ are constants.  We can bound the right-hand side using $\mathbf{I}_{N,Z}(\tau) > 1$, \eqref{DeLSpa freq mono eqn2} and \eqref{DeLSpa freq mono eqn4} to obtain 
\begin{equation*}
	\big| \mathbf{H}'_{N,Z}(\tau) - \tfrac{n-1}{\tau} \,\mathbf{H}_{N,Z}(\tau) - 2 \mathbf{D}_{N,Z}(\tau) \big| 
	\leq C \mathbf{m}_0^{\gamma_3} \tau^{\gamma_3} \mathbf{D}_{N,Z}(\tau)  
\end{equation*}
 for all $\tau \in [\sigma,\rho]$, where $C = C(n,m,q,M_0,N_0,C_e,C_h) \in (0,\infty)$ is a constant.  Hence using Theorem~\ref{freq mono centered thm} and $1 + C\mathbf{m}_0^{\gamma_3} \tau^{\gamma_3} \leq e^{C \mathbf{m}_0^{\gamma_3} \tau^{\gamma_3}}$ 
 \begin{equation}\label{doubling centered eqn1}
 	2 (1 - C \mathbf{m}_0^{\gamma_3} \tau^{\gamma_3}) \,\frac{\mathbf{I}_{N,Z}(\sigma_0)}{\tau} 
	\leq \frac{\mathbf{H}'_{N,Z}(\tau)}{\mathbf{H}_{N,Z}(\tau)} - \frac{n-1}{\tau} 
	\leq 2 e^{C \mathbf{m}_0^{\gamma_3} \rho^{\gamma_3}} \,\frac{\mathbf{I}_{N,Z}(\rho)}{\tau}
\end{equation}
for all $\tau \in [\sigma,\rho]$, where $C = C(n,m,q,M_0,N_0,C_e,C_h) \in (0,\infty)$ is a constant.  Integrating \eqref{doubling centered eqn1} over $\tau \in [\sigma,\rho]$, we conclude that \eqref{doubling centered concl1} holds true.  Noting that $\frac{2}{\rho} \leq \frac{-1}{R} \,\phi'(R/\rho) \leq \frac{4}{\rho}$ on $B^{\mathcal M}_{\rho}(Z) \setminus B^{\mathcal M}_{\rho/2}(Z)$ and summing \eqref{doubling centered concl1} with $\sigma_0 \rightarrow 0^+$ and with $2^{-k} \sigma,2^{-k} \rho$ in place of $\sigma,\rho$ for $k = 0,1,2,\ldots$ gives us \eqref{doubling centered concl2} (see the proof of~\cite[Corollary~5.3]{KrumWica}). 
 \end{proof}

\subsection{Coarse blow-ups relative to effective center manifolds}\label{sec:centered blowup sec}  Now we want to define a class of coarse blow-ups of a sequence of locally area-minimizing rectifiable currents relative to the corresponding center manifolds.  Let $M_0,N_0,C_e,C_h,\varepsilon_2$ be as in Hypothesis~\ref{DeLSpa consts freq hypothesis}.  For $\nu = 1,2,3,\ldots$, let $\mathbf{m}_{0,\nu} > 0$ and $T^{(\nu)}$ be an $n$-dimensional locally area-minimizing rectifiable current of $\mathbf{C}_{11\sqrt{n}/2}(0)$ such that Hypothesis~\ref{center mnfd hypothesis} holds true with $\mathbf{m}_{0,\nu},T^{(\nu)}$ in place of $\mathbf{m}_0,T$.  Let $\mathcal{M}^{(\nu)}$ be the center manifold of $T^{(\nu)}$ and $N^{(\nu)}$ be the normal approximation map of $T^{(\nu)}$ as in Hypothesis~\ref{M and N hypothesis} (with $T^{(\nu)},\mathcal{M}^{(\nu)},N^{(\nu)}$ in place of $T,\mathcal{M},N$).  Let $\mathcal{B}^{(\geq 2)}_q(T^{(\nu)})$ be as in Definition~\ref{center mnfd branch defn} (with $T^{(\nu)}$ in place of $T$).  Let $Z_{\nu} \in \mathcal{B}^{(\geq 2)}_q(T^{(\nu)})$, $P_{\nu} = {\rm T}_{Z_{\nu}} \mathcal{M}^{(\nu)}$, $\sigma_0 > 0$, and $\rho_{\nu} > 0$ be such that $0 < \rho_{\nu} \sigma_0 \leq 2$ and $\mathbf{m}_{0,\nu} \rho_{\nu} \rightarrow 0$ as $\nu\rightarrow\infty$.  Set $\widehat{E}_{\nu} = (\rho_{\nu}\sigma_0)^{-n/2-1} \|N^{(\nu)}\|_{L^2(B^{\mathcal{M}^{(\nu)}}_{\rho_{\nu}\sigma_0}(Z_{\nu}))}$ and note that by \eqref{normal approx est1}, $\widehat{E}_{\nu} \leq C\mathbf{m}_{0,\nu}^{1/(2n)} \rho_{\nu}^{\beta_2} \rightarrow 0$.  (More generally, we can let $(\rho_{\nu}\sigma_0)^{-n/2-1} \|N^{(\nu)}\|_{L^2(B^{\mathcal{M}^{(\nu)}}_{\rho_{\nu}\sigma_0}(Z_{\nu}))} \leq C \widehat{E}_{\nu}$ for a suitable constants $C \in (0,\infty)$.)  Let $\mathfrak{q}_{\nu}$ be a rotation of $\mathbb{R}^{n+m}$ such that $\mathfrak{q}_{\nu}(P_{\nu}) = \mathbb{R}^n \times \{0\}$.  Rescale and rotate by setting 
\begin{gather}\label{centered blowup eqn1} 
	\overline{T}^{(\nu)} = \mathfrak{q}_{\nu\#} \eta_{Z_{\nu},\rho_{\nu}\#} T^{(\nu)}, \quad\quad
	\overline{\mathcal M}^{(\nu)} = \mathfrak{q}_{\nu}(\eta_{Z_{\nu},\rho_{\nu}}(\mathcal{M}^{(\nu)})) , \\ 
	\overline{N}^{(\nu)}(X) = \sum_{l=1}^q \llbracket \rho_{\nu}^{-1} \mathfrak{q}_{\nu} N_l^{(\nu)}(Z_{\nu} + \rho_{\nu} \mathfrak{q}_{\nu}^{-1} X) \rrbracket 
		 \text{ for each } X \in \overline{\mathcal M}^{(\nu)} . \nonumber
\end{gather}
Note that each current $\overline{T}^{(\nu)}$ is locally area-minimizing in $\mathbf{C}^{\overline{\mathcal M}^{(\nu)}}_{\sigma_0}(0)$ with $(\partial\overline{T}^{(\nu)}) \llcorner \mathbf{C}^{\overline{\mathcal M}^{(\nu)}}_{\sigma_0}(0) = 0$.  We let $\overline{N}^{(\nu)}(X) = \sum_{l=1}^q \llbracket \overline{N}^{(\nu)}_l(X) \rrbracket$ and $\overline{N}^{(\nu)}_a(X) = \frac{1}{q} \sum_{l=1}^q \overline{N}^{(\nu)}_l(X)$ for all $X \in \overline{\mathcal M}^{(\nu)}$.  We let $\overline{\mathbf e}_{\nu} : B_{\sigma_0}(0) \rightarrow \overline{\mathcal M}^{(\nu)}$ be the exponential map of $\overline{\mathcal M}^{(\nu)}$ at the origin.  Notice that $\|A_{\overline{\mathcal M}^{(\nu)}}\|_{C^1(\overline{\mathcal M}^{(\nu)} \cap \mathbf{C}_{\sigma_0}(0))} \leq C\mathbf{m}_{0,\nu}^{1/2} \rho_{\nu} \rightarrow 0$ as $\nu \rightarrow \infty$ and, by \eqref{exp map est}, $\|D^2 \overline{\mathbf e}_{\nu}\|_{C^0(B_{\sigma_0}(0))} \leq C\mathbf{m}_{0,\nu}^{1/2} \rho_{\nu} \rightarrow 0$ as $\nu \rightarrow \infty$, where $C = C(n,m,q,M_0,N_0,C_e,C_h,\sigma_0) \in (0,\infty)$ are constants.  Hence $\overline{\mathcal M}^{(\nu)} \rightarrow \mathbb{R}^n \times \{0\}$ in the $C^3$-topology in $\mathbf{C}_{\sigma_0}(0)$ and $\overline{\mathbf e}_{\nu}(x) \rightarrow (x,0)$ in the $C^2$-topology in $B_{\sigma_0}(0)$ as $\nu \rightarrow \infty$.  By Lemma~\ref{normal approx energy lemma}, 
\begin{align}\label{centered blowup eqn2} 
	&\frac{1}{\sigma_0^{n+2}} \int_{B^{\overline{\mathcal M}^{(\nu)}}_{\gamma\sigma_0}(0)} |\overline{N}^{(\nu)}|^2 
		+ \frac{1}{\sigma_0^n} \int_{B^{\overline{\mathcal M}^{(\nu)}}_{\gamma\sigma_0}(0)} |D\overline{N}^{(\nu)}|^2 
	\\ \leq\,& \frac{1}{(\rho_{\nu}\sigma_0)^{n+2}} \int_{B^{\mathcal{M}^{(\nu)}}_{\gamma\rho_{\nu}\sigma_0}(Z)} |N^{(\nu)}|^2 
		+ \frac{1}{(\rho_{\nu}\sigma_0)^n} \int_{B^{\mathcal{M}^{(\nu)}}_{\gamma\rho_{\nu}\sigma_0}(Z)} |DN^{(\nu)}|^2 
	\leq C \widehat{E}_{\nu}^2 \nonumber 
\end{align}
for all $\gamma \in (0,1)$ and all positive integers $\nu$, where $C = C(n,m,q,M_0,N_0,C_e,C_h,\gamma) \in (0,\infty)$ is a constant.  Thus by the Rellich compactness theorem~\cite[Proposition 2.11]{DeLSpaDirMin}, after passing to a subsequence there exists a function $\varphi \in W^{1,2}_{\rm loc}(B_{\sigma_0}(0),\mathcal{A}_q(\mathbb{R}^{n+m}))$ such that 
\begin{equation}\label{centered blowup eqn3}
	\frac{\overline{N}^{(\nu)}(\overline{\mathbf e}_{\nu}^{-1}(x))}{\widehat{E}_{\nu}} \rightarrow \varphi(x)
\end{equation}
pointwise $\mathcal{L}^n$-a.e.~in $B_{\sigma_0}(0)$ and in the $L^2$-topology on compact subsets of $B_{\sigma_0}(0)$ as $\nu \rightarrow \infty$.  We call $\varphi$ a \emph{(coarse) blow-up of $T^{(\nu)}$ relative to $\mathcal{M}^{(\nu)}$ in $B^{\mathcal{M}^{(\nu)}}_{\rho_{\nu}\sigma_0}(Z_{\nu})$ by the excess $\widehat{E}_{\nu}$}.  We shall let $\varphi(x) = \sum_{l=1}^q \llbracket \varphi_l(x) \rrbracket$ and $\varphi_a(x) = \frac{1}{q} \sum_{l=1}^q \varphi_l(x)$ for each $x \in B_{\sigma_0}(0)$ and $D\varphi(x) = \sum_{l=1}^q \llbracket D\varphi_l(x) \rrbracket$ for $\mathcal{L}^n$-a.e.~$x \in B_{\sigma_0}(0)$, where we follow the conventions from~\cite[Subsection~3.5]{KrumWica}.  Since $|\pi_{P_0}(\overline{N}^{(\nu)}_i(X))| \leq \|\pi_{{\rm T}_X \overline{\mathcal M}^{(\nu)}} - \pi_{P_0}\| \,|\overline{N}^{(\nu)}_i(X)| \leq C\mathbf{m}_{0,\nu}^{1/2} \rho_{\nu} \sigma_0 |\overline{N}^{(\nu)}_i(X)|$ for all $X \in B^{\overline{\mathcal M}^{(\nu)}}_{\sigma_0}(0)$ 
\begin{equation*} 
	\frac{1}{\sigma_0^{n+2}} \int_{B^{\overline{\mathcal M}^{(\nu)}}_{\gamma\sigma_0}(0)} \sum_{l=1}^q |\pi_{P_0}(\overline{N}^{(\nu)}_l)|^2 
	\leq C \mathbf{m}_{0,\nu} (\rho_{\nu} \sigma_0)^2 \widehat{E}_{\nu}^2, 
\end{equation*}
where $P_0 = \mathbb{R}^n \times \{0\}$ and $C = C(n,m,q,M_0,N_0,C_e,C_h,\gamma) \in (0,\infty)$ are constants.  Hence by \eqref{centered blowup eqn3}, $\varphi_l(x) \in \{0\} \times \mathbb{R}^m$ for all $x \in B_{\sigma_0}(0)$ and $l \in \{1,2,\ldots,q\}$; that is, $\varphi \in W^{1,2}_{\rm loc}(B_{\sigma_0}(0),\mathcal{A}_q(\mathbb{R}^m))$ (where we identify $\mathbb{R}^m \cong \{0\} \times \mathbb{R}^m$).

\begin{lemma}\label{centered blowup lemma}
Let $\varphi$ be the blow-up of $T^{(\nu)}$ relative to $\mathcal{M}^{(\nu)}$ in $B^{\mathcal{M}^{(\nu)}}_{\rho_{\nu}\sigma_0}(Z_{\nu})$ (as in \eqref{centered blowup eqn3}).  Then $\varphi$ is average-free and locally Dirichlet energy minimizing in $B_{\sigma_0}(0)$ and 
\begin{equation}\label{centered blowup eqn4}
	\frac{|D(\overline{N}^{(\nu)} \circ \overline{\mathbf e}_{\nu}^{-1})(X)|}{\widehat{E}_{\nu}} \rightarrow |D\varphi(X)| 
\end{equation}
in the $L^2$-topology on compact subsets of $B_{\sigma_0}(0)$ as $\nu \rightarrow \infty$.
\end{lemma}

\begin{proof}
Recall that $\overline{R}(Z_{\nu}) < \rho_{\nu} \sigma_0 \leq 2$.  By $\beta_2 < \gamma_2$, \eqref{normal approx est1}, \eqref{normal approx est3} and \eqref{normal approx est4} (with $\rho = \rho_{\nu} \sigma_0$) and by \eqref{centered blowup eqn2} 
\begin{gather}
	\label{centered blowup eqn5a} \frac{1}{\sigma_0} \sup_{B^{\overline{\mathcal M}^{(\nu)}}_{\gamma\sigma_0}(0)} |\overline{N}^{(\nu)}| 
		+ \sup_{B^{\overline{\mathcal M}^{(\nu)}}_{\gamma\sigma_0}(0)} |D\overline{N}^{(\nu)}| 
		\leq C (\mathbf{m}_{0,\nu} \rho_{\nu} \sigma_0)^{\beta_2} , \\
	\label{centered blowup eqn5b} \frac{1}{\sigma_0^{n+1}} \int_{B^{\overline{\mathcal M}^{(\nu)}}_{\gamma\sigma_0}(0)} 
		|\overline{N}^{(\nu)}_a| \leq C \widehat{E}_{\nu}^2 , \\
	\label{centered blowup eqn5c} \|\mathbf{G}_{\overline{N}} - \overline{T}\|(\mathbf{C}^{\overline{\mathcal M}^{(\nu)}}_{\gamma\sigma_0}(0)) 
		\leq C (\mathbf{m}_{0,\nu}\rho_{\nu} \sigma_0)^{3\gamma_2/4} \widehat{E}_{\nu}^2 \sigma_0^n
\end{gather}
for all sufficiently large $\nu$, where $C = C(n,m,q,M_0,N_0,C_e,C_h,\gamma) \in (0,\infty)$ are constants.  By \eqref{centered blowup eqn3} and \eqref{centered blowup eqn5b}, $\varphi$ is average-free in $B_{\sigma_0}(0)$.

Let $L = L(m,q) \geq 1$ be an integer and $\bs{\xi} : \mathcal{A}_q(\mathbb{R}^m) \rightarrow \mathbb{R}^L$ be the bi-Lipschitz injection due to B.~White, see~\cite[Definition~1.1(6) and Theorem~1.2]{Almgren}, \cite[p.~706]{Cha88}, \cite[Corollary~2.2]{DeLSpaDirMin}.  By Rellich's compactness lemma, 
\begin{equation*}
	\frac{\bs{\xi}(\overline{N}^{(\nu)}(\overline{\mathbf e}_{\nu}^{-1}(x)))}{\widehat{E}_{\nu}} \rightarrow \bs{\xi}(\varphi(x))
\end{equation*}
weakly in $W^{1,2}$ on compact subsets of $B_{\sigma_0}(0)$ and in particular 
\begin{equation}\label{centered blowup eqn6}
	\int_{B_{\gamma\sigma_0}(0)} |D\varphi|^2 \leq \liminf_{\nu \rightarrow \infty} 
		\frac{1}{\widehat{E}_{\nu}^2} \int_{B^{\overline{\mathcal M}^{(\nu)}}_{\gamma\sigma_0}(0)} |D\overline{N}^{(\nu)}|^2 
\end{equation}
for all $\gamma \in (0,1)$.  Moreover, 
\begin{align}\label{centered blowup eqn7}
	&\limsup_{\nu \rightarrow \infty} \int_{B_{\gamma\sigma_0}(0)} 
		\left( \frac{|D(\overline{N}^{(\nu)} \circ \overline{\mathbf e}_{\nu}^{-1})|}{\widehat{E}_{\nu}} - |D\varphi| \right)^2
	\\=\,& \limsup_{\nu \rightarrow \infty} \int_{B_{\gamma\sigma_0}(0)} 
		\left( \frac{|D(\bs{\xi} \circ \overline{N}^{(\nu)} \circ \overline{\mathbf e}_{\nu}^{-1})|}{\widehat{E}_{\nu}} - |D(\bs{\xi} \circ \varphi)| \right)^2 
		\nonumber
	\\ \leq\,& \limsup_{\nu \rightarrow \infty} \int_{B_{\gamma\sigma_0}(0)} 
		\left| \frac{D(\bs{\xi} \circ \overline{N}^{(\nu)} \circ \overline{\mathbf e}_{\nu}^{-1})}{\widehat{E}_{\nu}} - D(\bs{\xi} \circ \varphi) \right|^2 \nonumber 
	\\=\,& \limsup_{\nu \rightarrow \infty} \int_{B_{\gamma\sigma_0}(0)} 
		\left( \frac{|D(\overline{N}^{(\nu)} \circ \overline{\mathbf e}_{\nu}^{-1})|^2}{\widehat{E}_{\nu}^2} + |D\varphi|^2 - 2 \left\langle 
		 \frac{D(\bs{\xi} \circ \overline{N}^{(\nu)} \circ \overline{\mathbf e}_{\nu}^{-1})}{\widehat{E}_{\nu}}, D(\bs{\xi} \circ \varphi) \right\rangle \right) \nonumber 
	\\=\,&  \limsup_{\nu \rightarrow \infty} \frac{1}{\widehat{E}_{\nu}^2} \int_{B^{\overline{\mathcal M}^{(\nu)}}_{\gamma\sigma_0}(0)} |D\overline{N}^{(\nu)}|^2 
		- \int_{B_{\gamma\sigma_0}(0)} |D\varphi|^2 \nonumber
\end{align}

To show that $\varphi$ is locally Dirichlet energy minimizing, suppose to the contrary that there exist $\delta > 0$, $\gamma \in (0,1)$, and a function $f \in W^{1,2}(B_{\sigma_0}(0),\mathcal{A}_q(\mathbb{R}^m))$ such that 
\begin{gather}\label{centered blowup eqn8} 
	f = \varphi \text{ on } B_{\sigma_0}(0) \setminus B_{\gamma\sigma_0}(0), \\ 
	\frac{1}{\sigma_0^n} \int_{B_{\gamma\sigma_0}(0)} |Df|^2 < \frac{1}{\sigma_0^n} \int_{B_{\gamma\sigma_0}(0)} |D\varphi|^2 - 2\delta. \nonumber
\end{gather}
Since $\varphi$ is average-free in $\mathbb{R}^n$, by (3.4) in~\cite{KrumWica} we may assume that $f$ is average-free in $B_{\sigma_0}(0)$.  Recall that $\overline{\mathbf e}_{\nu}(x) \rightarrow (x,0)$ in the $C^2$-topology in $B_{\sigma_0}(0)$ as $\nu\to\infty$.  Using~\cite[Proposition~4.4]{DeLSpa1} to interpolate between $(\overline{N}^{(\nu)} \circ \overline{\mathbf e}_{\nu}^{-1})/\widehat{E}_{\nu}$ and $f$, we can construct a Lipschitz $q$-valued function $\widehat{N}^{(\nu)} : B^{\overline{\mathcal M}^{(\nu)}}_{\rho_{\nu}\sigma_0}(0) \rightarrow \mathcal{A}_q(\mathbb{R}^{n+m})$ such that  
\begin{gather}
	\label{centered blowup eqn9a} \widehat{N}^{(\nu)} = \overline{N}^{(\nu)} \text{ for each } 
		X \in B^{\overline{\mathcal M}^{(\nu)}}_{\sigma_0}(0) \setminus B^{\overline{\mathcal M}^{(\nu)}}_{(1+\gamma)\sigma_0/2}(0), \\
	\label{centered blowup eqn9b} \frac{1}{\sigma_0} \sup_{B^{\overline{\mathcal M}^{(\nu)}}_{\sigma_0}(0)} |\widehat{N}^{(\nu)}| 
		+ \sup_{B^{\overline{\mathcal M}^{(\nu)}}_{\sigma_0}(0)} |D\widehat{N}^{(\nu)}| \leq C (\mathbf{m}_{0,\nu} \rho_{\nu} \sigma_0)^{\beta_2} , \\
	\label{centered blowup eqn9c} \frac{1}{\sigma_0^{n+1}} \int_{B^{\overline{\mathcal M}^{(\nu)}}_{\sigma_0}(0)} |\widehat{N}^{(\nu)}_a| 
		\leq C \widehat{E}_{\nu}^2 , \\ 
	\label{centered blowup eqn9d} \frac{1}{\sigma_0^{n+2}} \int_{B^{\overline{\mathcal M}^{(\nu)}}_{\sigma_0}(0)} |\widehat{N}^{(\nu)}|^2 
		\leq C \widehat{E}_{\nu}^2 , \\ 
	\label{centered blowup eqn9e} \frac{1}{\sigma_0^n} \int_{B^{\overline{\mathcal M}^{(\nu)}}_{\sigma_0}(0)} |D\widehat{N}^{(\nu)}|^2 
		< \frac{1}{\sigma_0^n} \int_{B^{\overline{\mathcal M}^{(\nu)}}_{\sigma_0}(0)} |D\overline{N}^{(\nu)}|^2 - \frac{\delta}{2} \,\widehat{E}_{\nu}^2 , 
\end{gather}
where $\widehat{N}^{(\nu)}_a(X) = \frac{1}{q} \sum_{i=1}^q \widehat{N}^{(\nu)}_l(X)$ is the single-valued average of $\widehat{N}^{(\nu)}(X) = \sum_{i=1}^q \llbracket \widehat{N}^{(\nu)}_l(X) \rrbracket$ and $C \in (0,\infty)$ is a constant independent of $\nu$.  In particular, \eqref{centered blowup eqn9a} and \eqref{centered blowup eqn9e} follow directly from~\cite[Proposition~4.4]{DeLSpa1} using \eqref{centered blowup eqn3}, \eqref{centered blowup eqn8}, and \eqref{centered blowup eqn6}.  By \eqref{centered blowup eqn5a} we have that $\widehat{E}_{\nu} = \sigma_0^{-n/2-1} \|\overline{N}^{(\nu)}\|_{L^2(B^{\overline{\mathcal M}^{(\nu)}}_{\sigma_0}(0))} \leq C (\mathbf{m}_{0,\nu} \rho_{\nu} \sigma_0)^{\beta_2}$.  Hence \eqref{centered blowup eqn9b} follows from (4.9) of~\cite{DeLSpa1} and \eqref{centered blowup eqn5a} and \eqref{centered blowup eqn9a}.  \eqref{centered blowup eqn9c} follows from (4.11) of~\cite{DeLSpa1} and \eqref{centered blowup eqn5b}.  \eqref{centered blowup eqn9d} follows from (4.10) of~\cite{DeLSpa1}, \eqref{centered blowup eqn2}, and \eqref{centered blowup eqn9e}.  Define $\widetilde{N}^{(\nu)} : B^{\overline{\mathcal M}^{(\nu)}}_{\sigma_0}(0) \rightarrow \mathcal{A}_q(({\rm T}\overline{\mathcal M}^{(\nu)})^{\perp})$ by 
\begin{equation*}
	\widetilde{N}^{(\nu)}(X) = \sum_{l=1}^{q} \llbracket \pi_{({\rm T}_X\overline{\mathcal M}^{(\nu)})^{\perp}} (\widehat{N}^{(\nu)}_l(X)) \rrbracket
\end{equation*}
for each $X \in B^{\mathcal{M}^{(\nu)}}_{\rho_{\nu}}(Z_{\nu})$.  
We can regard $\pi_{{\rm T}_X\overline{\mathcal M}^{(\nu)}}$ and $\pi_{({\rm T}_X\overline{\mathcal M}^{(\nu)})^{\perp}}$ as $(n+m) \times (n+m)$ matrix-valued functions on $X\overline{\mathcal M}^{(\nu)}$ such that $\pi_{{\rm T}_X\overline{\mathcal M}^{(\nu)}} + \pi_{({\rm T}_X\overline{\mathcal M}^{(\nu)})^{\perp}} = I$, where $I$ is the $(n+m) \times (n+m)$ identify matrix.  Using normal coordinates at points of $\overline{\mathcal M}^{(\nu)}$ and the definition of the second fundamental form $A_{\overline{\mathcal M}^{(\nu)}}$, $|D\pi_{({\rm T}_X\overline{\mathcal M}^{(\nu)})^{\perp}}| = |D\pi_{({\rm T}_X\overline{\mathcal M}^{(\nu)})^{\perp}}| \leq C(n,m) \,|A_{\overline{\mathcal M}^{(\nu)}}|$ at each point of $\overline{\mathcal M}^{(\nu)}$.  Hence $|D\widetilde{N}^{(\nu)}| \leq |D\widehat{N}^{(\nu)}| + C(n,m) \,|A_{\overline{\mathcal M}^{(\nu)}}| \,|\widehat{N}^{(\nu)}|$ on $B^{\overline{\mathcal M}^{(\nu)}}_{\sigma_0}(0)$.  Noting that $|A_{\overline{\mathcal M}^{(\nu)}}| \leq C \mathbf{m}_{0,\nu}^{1/2} \rho_{\nu} \rightarrow 0$, it follows from \eqref{centered blowup eqn9a}--\eqref{centered blowup eqn9e} that 
\begin{gather}\label{centered blowup eqn10}
	\widetilde{N}^{(\nu)} = \overline{N}^{(\nu)} \text{ for each } X \in B^{\overline{\mathcal M}^{(\nu)}}_{\sigma_0}(0) \setminus 
		B^{\overline{\mathcal M}^{(\nu)}}_{(1+\gamma)\sigma_0/2}(0), \\
	\frac{1}{\sigma_0} \sup_{B^{\overline{\mathcal M}^{(\nu)}}_{\sigma_0}(0)} |\widetilde{N}^{(\nu)}| 
		+ \sup_{B^{\overline{\mathcal M}^{(\nu)}}_{\sigma_0}(0)} |D\widetilde{N}^{(\nu)}| \leq C (\mathbf{m}_{0,\nu} \rho_{\nu} \sigma_0)^{\beta_2} , 
		\nonumber\\
	\frac{1}{\sigma_0^{n+1}} \int_{B^{\overline{\mathcal M}^{(\nu)}}_{\sigma_0}(0)} |\widetilde{N}^{(\nu)}_a| \leq C \widehat{E}_{\nu}^2 , \quad 
	\frac{1}{\sigma_0^{n+2}} \int_{B^{\overline{\mathcal M}^{(\nu)}}_{\sigma_0}(0)} |\widetilde{N}^{(\nu)}|^2 \leq C \widehat{E}_{\nu}^2 , \nonumber \\ 
	\frac{1}{\sigma_0^n} \int_{B^{\overline{\mathcal M}^{(\nu)}}_{\sigma_0}(0)} |D\widetilde{N}^{(\nu)}|^2 
		\leq \frac{1}{\sigma_0^n} \int_{B^{\overline{\mathcal M}^{(\nu)}}_{\sigma_0}(0)} |D\overline{N}^{(\nu)}|^2 
		- \frac{\delta}{2} \,\widehat{E}_{\nu}^2 , \nonumber
\end{gather}
where $\widetilde{N}^{(\nu)}_a(X) = \frac{1}{q} \sum_{i=1}^q \widetilde{N}^{(\nu)}_l(X)$ is the single-valued average of $\widetilde{N}^{(\nu)}(X) = \sum_{i=1}^q \llbracket \widetilde{N}^{(\nu)}_l(X) \rrbracket$ and $C \in (0,\infty)$ is a constant independent of $\nu$.  By \eqref{centered blowup eqn5c}, after passing to a subsequence there exists $\tau \in \big(\frac{(1+\gamma)\sigma_0}{2},\frac{(3+\gamma)\sigma_0}{4}\big)$ such that the slice $\langle \mathbf{G}_{\overline{N}^{(\nu)}} - \overline{T}^{(\nu)}, R_{\nu} \circ \overline{\mathbf p}_{\nu}, \tau \rangle$ is a rectifiable current 
\begin{equation*}
	\mathbf{M}(\langle \mathbf{G}_{\overline{N}^{(\nu)}} - \overline{T}^{(\nu)}, R_{\nu} \circ \overline{\mathbf p}_{\nu}, \tau \rangle) 
	\leq C (\mathbf{m}_{0,\nu}\rho_{\nu} \sigma_0)^{3\gamma_2/4} \widehat{E}_{\nu}^2 \sigma_0^{n-1} ,
\end{equation*}
where $R_{\nu} = R_{\nu}(X)$ is the geodesic distance from $X$ to the origin in $\overline{\mathcal M}^{(\nu)}$, $\overline{\mathbf p}_{\nu}$ is the closest point projection map onto $\overline{\mathcal M}^{(\nu)}$, and $C = C(n,m,q,M_0,N_0,C_e,C_h,\gamma) \in (0,\infty)$ is a constant.  By the isoperimetric inequality~\cite[Theorem~30.1]{SimonGMT}, there exists an $n$-dimensional rectifiable current $S^{(\nu)}$ such that 
\begin{equation}\label{centered blowup eqn12}
	\partial S^{(\nu)} = \langle \mathbf{G}_{\overline{N}^{(\nu)}} - \overline{T}^{(\nu)}, R_{\nu} \circ \overline{\mathbf p}_{\nu}, \tau \rangle, \quad 
	\mathbf{M}(S^{(\nu)}) \leq C (\mathbf{m}_{0,\nu}\rho_{\nu} \sigma_0)^{3\gamma_2/4} \widehat{E}_{\nu}^2 \sigma_0^n ,
\end{equation}
where $C \in (0,\infty)$ is a constant independent of $\nu$.  Define the competitor $\widetilde{T}^{(\nu)}$ for $\overline{T}^{(\nu)}$ in $\mathbf{C}^{\overline{\mathcal M}^{(\nu)}}_{\sigma_0}(0)$ by 
\begin{equation}\label{centered blowup eqn13}
	\widetilde{T}^{(\nu)} = \mathbf{G}_{\widetilde{N}^{(\nu)}} \llcorner \mathbf{C}^{\overline{\mathcal M}^{(\nu)}}_{\tau}(0) + S^{(\nu)} 
	+ \overline{T}^{(\nu)} \llcorner \mathbf{C}^{\overline{\mathcal M}^{(\nu)}}_{\sigma_0}(0) \setminus \mathbf{C}^{\overline{\mathcal M}^{(\nu)}}_{\tau}(0)
\end{equation}
Assuming $\varepsilon_2$ is sufficiently small, by~\cite[Theorem~3.2]{DeLSpa0}, $\|A_{\overline{\mathcal M}^{(\nu)}}\|_{C^{2,\kappa}} \leq C\mathbf{m}_{0,\nu}^{1/2} \rho_{\nu} \rightarrow 0$, \eqref{centered blowup eqn2}, \eqref{centered blowup eqn5a}--\eqref{centered blowup eqn5c}, \eqref{centered blowup eqn10}, and \eqref{centered blowup eqn12} 
\begin{align*}
	&\mathbf{M}\big( \widetilde{T}^{(\nu)} \llcorner \mathbf{C}^{\overline{\mathcal M}^{(\nu)}}_{\sigma_0}(0) \big) 
		- \mathbf{M}\big( \overline{T}^{(\nu)} \llcorner \mathbf{C}^{\overline{\mathcal M}^{(\nu)}}_{\sigma_0}(0) \big) 
	\\ \leq\,& \mathbf{M}\big( \mathbf{G}_{\widetilde{N}^{(\nu)}} \llcorner \mathbf{C}^{\overline{\mathcal M}^{(\nu)}}_{\tau}(0) \big) 
		- \mathbf{M}\big( \mathbf{G}_{\overline{N}^{(\nu)}} \llcorner \mathbf{C}^{\overline{\mathcal M}^{(\nu)}}_{\tau}(0) \big) 
		+ C (\mathbf{m}_{0,\nu}\rho_{\nu} \sigma_0)^{3\gamma_2/4} \widehat{E}_{\nu}^2 \sigma_0^n 
	\\ \leq\,& q \int_{B^{\overline{\mathcal M}^{(\nu)}}_{\tau}(0)} |H_{\overline{\mathcal M}^{(\nu)}}| \,(|\overline{N}^{(\nu)}_a| + |\widetilde{N}^{(\nu)}_a|) 
		+ \frac{1}{2} \int_{B^{\overline{\mathcal M}^{(\nu)}}_{\tau}(0)} (|D\widetilde{N}^{(\nu)}|^2 - |D\overline{N}^{(\nu)}|^2) 
		\\&+ C \int_{B^{\overline{\mathcal M}^{(\nu)}}_{\tau}(0)} (|A_{\mathcal M}|^2 |\overline{N}^{(\nu)}|^2 + |A_{\mathcal M}|^2 |\widetilde{N}^{(\nu)}|^2 
		+ |D\overline{N}^{(\nu)}|^4 + |D\widetilde{N}^{(\nu)}|^4) 
		\\&+ C (\mathbf{m}_{0,\nu}\rho_{\nu} \sigma_0)^{3\gamma_2/4} \widehat{E}_{\nu}^2 \sigma_0^n 
	\\ \leq\,& - \frac{\delta}{8} \,\widehat{E}_{\nu}^2 \sigma_0^n < 0 
\end{align*}
for all sufficiently large $\nu$ and for some constants $C \in (0,\infty)$ independent of $\nu$, contradicting $\overline{T}^{(\nu)}$ being locally area minimizing.

In light \eqref{centered blowup eqn6} and \eqref{centered blowup eqn7}, to show \eqref{centered tanfun concl2} it suffices to show that for every $\gamma \in (0,1)$ 
\begin{equation*}
	\int_{B_{\gamma\sigma_0}(0)} |D\varphi|^2 \geq \limsup_{\nu\rightarrow \infty} \frac{1}{\widehat{E}_{\nu}^2} 
		\int_{B^{\overline{\mathcal M}^{(\nu)}}_{\gamma\sigma_0}(0)} |D\overline{N}^{(\nu)}|^2 .
\end{equation*}
Suppose to the contrary that after passing to a subsequence there exists $\gamma \in (0,1)$ and $\delta > 0$ such that 
\begin{equation}\label{centered blowup eqn15}
	\frac{1}{\sigma_0^n} \int_{B_{\gamma\sigma_0}(0)} |D\varphi|^2 <  \frac{1}{\widehat{E}_{\nu}^2 \sigma_0^n} 
		\int_{B^{\overline{\mathcal M}^{(\nu)}}_{\gamma\sigma_0}(0)} |D\overline{N}^{(\nu)}|^2 - 2\delta
\end{equation}
for all $\nu$.  Arguing as in the previous paragraph with $\varphi$ in place of $f$, we can use~\cite[Proposition~4.6]{DeLSpa1} to interpolate between $(\overline{N}^{(\nu)} \circ \overline{\mathbf e}_{\nu}^{-1})/\widehat{E}_{\nu}$ and $\varphi$ and rescale to find $\widehat{N}^{(\nu)} : B^{\overline{\mathcal M}^{(\nu)}}_{\sigma_0}(0) \rightarrow \mathcal{A}_q({\rm T}\mathcal{M}^{\perp})$ such that \eqref{centered blowup eqn9a}--\eqref{centered blowup eqn9e}.  Hence arguing as in the previous paragraph, we can construct a competitor $\widetilde{T}^{(\nu)}$ for $\overline{T}^{(\nu)}$ as in \eqref{centered blowup eqn13} and thereby reach a contradiction.
\end{proof}

\begin{lemma}\label{centered semicont2 lemma}
Let $\varphi$ be the blow-up of $T^{(\nu)}$ relative to $\mathcal{M}^{(\nu)}$ in $B^{\mathcal{M}^{(\nu)}}_{\rho_{\nu}\sigma_0}(Z_{\nu})$ (as in \eqref{centered blowup eqn3}).  Let $Y_{\nu} \in \mathcal{B}^{(\geq 2)}_q(T^{(\nu)})$, $\xi \in B_{\sigma_0}(0)$, $\sigma_{\nu} \in (0,\sigma_0 \rho_{\nu} - |Y_{\nu}|)$, and $\sigma \in (0,\sigma_0-|\xi|)$ such that $\mathfrak{q}_{\nu}(Y_{\nu} - Z_{\nu})/\rho_{\nu} \rightarrow (\xi,0)$ and $\sigma_{\nu}/\rho_{\nu} \rightarrow \sigma$ as $\nu\rightarrow\infty$.  Then 
\begin{equation}\label{centered semicont2 concl1}
	N_{\varphi,\xi}(\sigma) = \lim_{\nu \rightarrow \infty} \mathbf{I}_{N^{(\nu)},Y_{\nu}}(\sigma_{\nu}) ,
\end{equation}
where $N_{\varphi,\xi}$ is as in~\cite[Subsection~3.5]{KrumWica}.  Moreover, 
\begin{equation}\label{centered semicont2 concl2}
	\mathcal{N}_{\varphi}(\xi) \geq \limsup_{\nu \rightarrow \infty} \mathbf{I}_{N^{(\nu)},Y_{\nu}}(0^+) , 
\end{equation}
where $\mathcal{N}_{\varphi}(\xi) = \lim_{\rho\rightarrow 0} N_{\varphi,\xi}(\rho)$.
\end{lemma}

\begin{proof}
Set $\overline{Y}_{\nu} = \mathfrak{q}_{\nu}(Y_{\nu} - Z_{\nu})/\rho_{\nu}$ and $\overline{\sigma}_{\nu} = \sigma_{\nu}/\rho_{\nu}$.  For each $X \in B^{\overline{\mathcal M}^{(\nu)}}_{\sigma_0}(0)$, let $R_{\nu} = R_{\nu}(X)$ be the geodesic distance of $X$ from $\overline{Y}_{\nu}$ in $\overline{\mathcal M}^{(\nu)}$.  Recall that $\overline{\mathcal M}^{(\nu)} \rightarrow \mathbb{R}^n \times \{0\}$ in the $C^3$-topology in $\mathbf{C}_{\sigma_0}(0)$ and $\overline{\mathbf e}_{\nu}(x) \rightarrow (x,0)$ in the $C^2$-topology in $B_{\sigma_0}(0)$ as $\nu \rightarrow \infty$.  Hence $R_{\nu}(\overline{\mathbf e}_{\nu}(x)) \rightarrow |x-\xi|$ uniformly on $B_{\sigma_0}(0)$ as $\nu\rightarrow\infty$.  In particular, $\phi((R_{\nu} \circ \overline{\mathbf e}_{\nu})/\overline{\sigma}_{\nu}) \rightarrow \phi(|x-\xi|/\sigma)$ uniformly on $B_{\sigma_0}(0)$ and $\phi'((R_{\nu} \circ \overline{\mathbf e}_{\nu})/\overline{\sigma}_{\nu}) \rightarrow \phi'(|x-\xi|/\sigma)$ point-wise $\mathcal{L}^n$-a.e.~in $B_{\sigma_0}(0)$.  By the dominated convergence theorem, \eqref{centered blowup eqn3} and \eqref{centered blowup eqn4} 
\begin{align*}
	H_{\varphi,\xi}(\sigma) 
		=\,& \lim_{\nu \rightarrow \infty} \frac{-1}{\widehat{E}_{\nu}^2 \overline{\sigma}_{\nu}^{n-1}} 
			\int_{B^{\overline{\mathcal M}^{(\nu)}}_{\overline{\sigma}_{\nu}}(0)} 
			|\overline{N}^{(\nu)}|^2 \,\frac{1}{R_{\nu}} \,\phi'(R_{\nu}/\overline{\sigma}_{\nu}) 
		= \lim_{\nu \rightarrow \infty} \frac{\mathbf{H}_{N^{(\nu)},Y_{\nu}}(\sigma_{\nu})}{\widehat{E}_{\nu}^2 \sigma_{\nu}^{n-1}} > 0, \\
	D_{\varphi,\xi}(\sigma) 
		=\,& \lim_{\nu \rightarrow \infty} \frac{1}{\widehat{E}_{\nu}^2 \overline{\sigma}_{\nu}^{n-2}} 
			\int_{B^{\overline{\mathcal M}^{(\nu)}}_{\overline{\sigma}_{\nu}}(0)} |D\overline{N}^{(\nu)}|^2 \,\phi(R_{\nu}/\overline{\sigma}_{\nu}) 
		= \lim_{\nu \rightarrow \infty} \frac{\mathbf{D}_{N^{(\nu)},Y_{\nu}}(\sigma_{\nu})}{\widehat{E}_{\nu}^2 \sigma_{\nu}^{n-2}} , 
\end{align*}
where $H_{\varphi,\xi}$ and $D_{\varphi,\xi}$ are as in~\cite[Section~3.6]{KrumWica}, and consequently \eqref{centered semicont2 concl1} holds true.  
By Theorem~\ref{freq mono centered thm} and \eqref{centered semicont2 concl1} (with $\sigma_{\nu} = \rho_{\nu} \sigma$), 
\begin{equation}\label{centered semicont2 eqn1}
	N_{\varphi,\xi}(\sigma) 
	= \lim_{\nu \rightarrow \infty} e^{C (\mathbf{m}_{0,\nu} \rho_{\nu} \sigma)^{\gamma_3}} \mathbf{I}_{N^{(\nu)},Y_{\nu}}(\rho_{\nu} \sigma) 
	\geq \limsup_{\nu \rightarrow \infty} \mathbf{I}_{N^{(\nu)},Y_{\nu}}(0^+) , 
\end{equation}
for each $\sigma \in (0,2]$, where $\gamma_3 = \gamma_3(n,m,q) \in (0,1/2]$ and $C$ are as in Theorem~\ref{freq mono centered thm}.  Letting $\sigma \rightarrow 0^+$ in \eqref{centered semicont2 eqn1} gives us \eqref{centered semicont2 concl2}.
\end{proof}

\subsection{Tangent functions relative to an effective center manifold}\label{sec:tangent fun}  Let $M_0,N_0,C_e,C_h,\varepsilon_2$ be as in Hypothesis~\ref{DeLSpa consts freq hypothesis}.  Let $T, \mathcal{M}, N$ be as in Hypothesis~\ref{M and N hypothesis}.  Let $\mathcal{B}^{(\geq 2)}_q$ be as in Definition~\ref{center mnfd branch defn}.  Let $Z \in \mathcal{B}^{(\geq 2)}_q$ and let $(\rho_{\nu})_{\nu=1}^{\infty}$ be a sequence of positive numbers tending to zero.  For each positive integer $\nu$, let $\widehat{E}_{\nu} = \rho_{\nu}^{-n/2-1} \|N\|_{L^2(B^{\mathcal M}_{\rho_{\nu}}(Z))}$.  Recall that $\mathbf{I}_{N,Z}(0^+) \geq 2$.  For each fixed $\sigma \in (1,\infty)$ and every sufficiently large $\nu$, by Lemma~\ref{normal approx energy lemma}, Theorem~\ref{freq mono centered thm}, and Corollary~\ref{doubling centered cor} 
\begin{align}\label{centered tanfun eqn1} 
	\frac{1}{(\rho_{\nu}\sigma)^{n+2}} \int_{B^{\mathcal M}_{\rho_{\nu}\sigma}(Z)} |N|^2 
		+ \frac{1}{(\rho_{\nu}\sigma)^n} \int_{B^{\mathcal M}_{\rho_{\nu}\sigma}(Z)} |DN|^2 
	\leq\,& \frac{C}{(\rho_{\nu}\sigma)^{n+2}} \int_{B^{\mathcal M}_{2\rho_{\nu}\sigma}(Z)} |N|^2
	\\ \leq\,& C (2\sigma)^{2e^{C \mathbf{m}_0^{\gamma_3}} \mathbf{I}_{N,Z}(1) - 2} \widehat{E}_{\nu}^2 , \nonumber 
\end{align}
where $\gamma_3 = \gamma_3(n,m,q) \in (0,1/2]$ and $C = C(n,m,q,M_0,N_0,C_e,C_h) \in (0,\infty)$ are constants.  Thus using the Rellich compactness theorem~\cite[Proposition 2.11]{DeLSpaDirMin} and arguing as in Subsection~\ref{sec:centered blowup sec}, after passing to a subsequence there exists an average-free, locally Dirichlet energy minimizing $q$-valued function $\varphi \in W^{1,2}_{\rm loc}(P_Z,\mathcal{A}_q(P_Z^{\perp}))$ such that 
\begin{equation}\label{centered tanfun concl1}
	\frac{N(\exp_Z(\rho_{\nu} X))}{\rho_{\nu} \widehat{E}_{\nu}} \rightarrow \varphi(X)
\end{equation}
pointwise $\mathcal{L}^n$-a.e.~in $P_Z$ and in $L^2(B_{\sigma}(0,P_Z),\mathcal{A}_q(P_Z^{\perp}))$ for all $\sigma > 0$ as $\nu \rightarrow \infty$ and  
\begin{equation}\label{centered tanfun concl2}
	\frac{|D(N \circ \exp_Z)(\rho_{\nu} X)|}{\widehat{E}_{\nu}} \rightarrow |D\varphi(X)| 
\end{equation}
in $L^2(B_{\sigma}(0,P_Z))$ for all $\sigma > 0$ as $\nu \rightarrow \infty$.  (Note that in Subsection~\ref{sec:centered blowup sec}, we rotated by $\mathfrak{q}_{\nu}$ so that each $\overline{\mathcal M}^{(\nu)}$ is tangent to the same plane $\mathbb{R}^n \times \{0\}$ at the origin.  In this case, since the tangent plane $P_Z$ of $\mathcal{M}$ at $Z$ is fixed, we do not need to rotate.)  By \eqref{centered tanfun concl1}, $\|\varphi\|_{L^2(B_1(0))} = \omega_n^{1/2}$ and thus $\varphi$ is non-zero.  By Lemma~\ref{centered semicont2 lemma} 
\begin{equation*}
	N_{\varphi,0}(\sigma) = \lim_{\nu \rightarrow \infty} \mathbf{I}_{N,Z}(\rho_{\nu} \sigma) = \mathbf{I}_{N,Z}(0^+). 
\end{equation*}
for all $\sigma \in (0,\infty)$.  Therefore, $\varphi$ is homogeneous degree $\mathbf{I}_{N,Z}(0^+)$.  We call any such function $\varphi$ a \emph{tangent function} of $T$ at $Z$ (relative to the center manifold $\mathcal{M}$).  We shall let $\varphi(X) = \sum_{l=1}^q \llbracket \varphi_l(X) \rrbracket$ and $\varphi_a(X) = \frac{1}{q} \sum_{l=1}^q \varphi_l(X)$ for each $X \in P_Z$ and $D\varphi(X) = \sum_{l=1}^q \llbracket D\varphi_l(X) \rrbracket$ for $\mathcal{L}^n$-a.e.~$X \in P_Z$, where we follow the conventions from~\cite[Subsection~3.5]{KrumWica}.

For each $j \in \{0,1,2,\ldots,n-2\}$, we define $S^{(j)} = S^{(j)}(T)$ to be the set of all points $Z \in \mathcal{B}^{(\geq 2)}_q$ such that $\dim\op{spine}\varphi \leq j$ for every tangent function $\varphi$ of $T$ at $Z$.  (Recall from~\cite[Section~3.6]{KrumWica} that $\op{spine}\varphi = \{ X \in P_Z : \mathcal{N}_{\varphi}(X) = \mathcal{N}_{\varphi}(0) \}$.)  Observe that 
\begin{equation*}
	\mathcal{B}^{(\geq 2)}_q = S^{(n-2)} \supseteq S^{(n-3)} \supseteq \cdots \supseteq S^{(1)} \supseteq S^{(0)} .
\end{equation*}

\begin{theorem}\label{centered stratification thm}
Let $M_0,N_0,C_e,C_h,\varepsilon_2$ be as in Hypothesis~\ref{DeLSpa consts freq hypothesis}.  Let $T, \mathcal{M}, N$ be as in Hypothesis~\ref{M and N hypothesis}.  Let $\mathcal{B}^{(\geq 2)}_q$ be as in Definition~\ref{center mnfd branch defn}.  Then:
\begin{enumerate}[itemsep=3mm,topsep=0mm]
	\item[{\rm (i)}]  for each $j \in \{0,1,2,\ldots,n-2\}$, $S^{(j)}$ has Hausdorff dimension at most $j$; 
	\item[{\rm (ii)}]  for each $\alpha > 0$, $\{ Z \in S^{(0)} : \mathbf{I}_{N,Z}(0^+) = \alpha \}$ is discrete;
	\item[{\rm (iii)}]  $\mathcal{B}^{(\geq 2)}_q$ has Hausdorff dimension at most $n-2$. 
\end{enumerate}
In particular, Theorem~\ref{centered size sing thm} from the Introduction holds true.
\end{theorem}

\begin{proof}
Apply the well-known argument~\cite[Theorem~2.26 and Corollary~2.27]{Almgren} using the upper semi-continuity of frequency Corollary~\ref{centered semicont1 cor} and Lemma~\ref{centered semicont2 lemma} and the fact that $\mathbf{I}_{N,Z}(0^+) = \mathcal{N}_{\varphi}(0)$ whenever $\varphi$ is a tangent function to $T$ at $Z$.
\end{proof}

\subsection{Locally uniform weak approximation of area-minimizing currents by homogeneous multi-valued harmonic functions}  The remainder of the paper will be focused on proving Theorem~\ref{B rectifiable thm}, Theorem~\ref{unique tangent function thm}, and Theorem~\ref{topology thm} from the Introduction.  We end this section by proving a weak locally uniform approximation property of locally area-minimizing rectifiable currents $T$ by homogeneous Dirichlet energy minimizing functions $\varphi$ relative to an effective center manifold $\mathcal{M}$ at branch points in $\mathcal{B}^{(\geq 2)}_q$.  This is analogous to \cite[Lemma~2.4]{Sim93} and \cite[Lemma~4.3]{KrumWic2}, and plays an important role in the proofs of Theorem~\ref{unique tangent function thm} and Theorem~\ref{B rectifiable thm}.

\begin{lemma}\label{lemma2_4}
Let $M_0,N_0,C_e,C_h,\varepsilon_2$ be as in Hypothesis~\ref{DeLSpa consts freq hypothesis} and let $\alpha_0 \in [2,\infty)$.  There exists a function $\delta : (0,1) \rightarrow (0,1)$ depending only on $n,m,q,M_0,N_0,C_e,C_h,\varepsilon_2,\alpha_0$ such that the following holds true.  Let $\alpha \in [2,\alpha_0]$ and $\eta \in (0,1)$.  Let $T, \mathcal{M}, N$ be as in Hypothesis~\ref{M and N hypothesis}.  Let $\mathcal{B}^{(\geq 2)}_q$ be as in Definition~\ref{center mnfd branch defn}.  Let $Z_0 \in \mathcal{B}^{(\geq 2)}_q$ and $\rho_0 \in (0,1]$ such that $\mathbf{m}_0 \rho_0 < \delta(\eta)$ and 
\begin{equation}\label{lemma2_4 hyp}
	\mathbf{I}_{N,Z_0}(2\rho_0) < \alpha + \delta(\eta) . 
\end{equation}
If $Z \in \mathcal{B}^{(\geq 2)}_q \cap B^{\mathcal M}_{\delta(\eta) \,\rho_0}(Z_0)$ with $\mathbf{I}_{N,Z}(0^+) \geq \alpha$, then:
\begin{enumerate}[itemsep=3mm,topsep=0mm]
	\item[{\rm (i)}]  $|\mathbf{I}_{N,Z}(\rho) - \alpha| < \eta^2$ for all $\rho \in (0,\rho_0]$; 
	\item[{\rm (ii)}]  for every $\rho \in (0,\delta(\eta) \,\rho_0]$ there exists a non-zero, average-free, locally Dirichlet energy minimizing $q$-valued function $\varphi \in W^{1,2}_{\rm loc}(P_Z,\mathcal{A}_q(P_Z^{\perp}))$ (depending on $\rho$) such that $\varphi$ is homogeneous degree $\mathcal{N}_{\varphi}(0)$, $|\mathcal{N}_{\varphi}(0) - \alpha| < \eta^2$, and 
	\begin{equation}\label{lemma2_4 concl1}
		\int_{B^{\mathcal M}_{\rho}(Z)} \mathcal{G}(N,\varphi\circ\exp_Z^{-1})^2 < \eta^2 \int_{B^{\mathcal M}_{\rho}(Z)} |N|^2
	\end{equation}	
	\item[{\rm (iii)}]  for every $\rho \in (0,\delta(\eta) \,\rho_0]$ either there exists a non-zero, average-free, locally Dirichlet energy minimizing $q$-valued function $\varphi \in W^{1,2}_{\rm loc}(P_Z,\mathcal{A}_q(P_Z^{\perp}))$ (depending on $\rho$) such that $\varphi$ is homogeneous degree $\mathcal{N}_{\varphi}(0)$, $|\mathcal{N}_{\varphi}(0) - \alpha| < \eta^2$, $\dim\op{spine}\varphi = n-2$, and 
	\begin{equation}\label{lemma2_4 concl2}
		\int_{B^{\mathcal M}_{\rho}(Z)} \mathcal{G}(N,\varphi\circ\exp_Z^{-1})^2 < \eta^2 \int_{B^{\mathcal M}_{\rho}(Z)} |N|^2
	\end{equation}
	or there exists an $(n-3)$-dimensional linear subspace $L$ of $P_Z$ (depending on $\rho$) such that 
	\begin{equation}\label{lemma2_4 concl3}
		\{ X \in \mathcal{B}^{(\geq 2)}_q \cap \overline{B^{\mathcal M}_{\rho}(Z)} : \mathbf{I}_{N,X}(0^+) \geq \alpha \} 
			\subset \{ X \in \mathbb{R}^{n+m} : \op{dist}(X,Z+L) < \eta\rho \} .
	\end{equation}
\end{enumerate}
\end{lemma}

\begin{proof}[Proof of Lemma~\ref{lemma2_4}]  \textit{Proof of Conclusion~{\rm (i)}.}  In light of Theorem~\ref{freq mono centered thm}, it suffices to show that $\mathbf{I}_{N,Z}(\rho_0) < \alpha + \eta^2/4$.  We proceed by contradiction.  Fix constants $M_0,N_0,C_e,C_h,\varepsilon_2$ and $\eta \in (0,1)$.  For $\nu = 1,2,3,\ldots$, let $\mathbf{m}_{0,\nu},T^{(\nu)},\mathcal{M}^{(\nu)},N^{(\nu)}$ be as in Hypothesis~\ref{center mnfd hypothesis} with $\mathbf{m}_{0,\nu},T^{(\nu)},\mathcal{M}^{(\nu)},N^{(\nu)}$ in place of $\mathbf{m}_0,T,\mathcal{M},N$.
Suppose to the contrary that for some $\delta_{\nu}\rightarrow 0^+$, $\alpha_{\nu} \in [2,\alpha_0]$, $Z_{0,\nu} \in \mathcal{B}^{(\geq 2)}_q(T^{(\nu)})$, $\rho_{0,\nu} \in (0,1]$, and $Z_{\nu} \in \mathcal{B}^{(\geq 2)}_q(T^{(\nu)}) \cap B^{\mathcal{M}^{(\nu)}}_{\delta_{\nu} \rho_{0,\nu}}(Z_{0,\nu})$ such that $\mathbf{m}_{0,\nu} \rho_{0,\nu} \to 0^+$ and 
\begin{gather}
	\label{lemma2_4 eqn1} \mathbf{I}_{N^{(\nu)},Z_{0,\nu}}(2\rho_{0,\nu}) < \alpha_{\nu} + \delta_{\nu}, \\
	\label{lemma2_4 eqn1b} \mathbf{I}_{N^{(\nu)},Z_{0,\nu}}(\rho_{0,\nu}) \geq \alpha_{\nu} + \eta^2/4 .  
\end{gather}
After passing to a subsequence, let $\alpha_{\nu} \rightarrow \alpha$ in $[2,\alpha_0]$ as $\nu\rightarrow\infty$.  Let $\varphi$ be a blow-up of $T^{(\nu)}$ relative to $\mathcal{M}^{(\nu)}$ in $B^{\mathcal{M}^{(\nu)}}_{2\rho_{0,\nu}}(Z_{0,\nu})$ (as in Subsection~\ref{sec:centered blowup sec} with $\widehat{E}_{\nu} = (2\rho_{0,\nu})^{-n/2-1} \|N^{(\nu)}\|_{L^2(B^{\mathcal{M}^{(\nu)}}_{2\rho_{0,\nu}}(Z_{0,\nu})}$).  By \eqref{lemma2_4 eqn1} and Corollary~\ref{doubling centered cor}, $\|N^{(\nu)}\|_{L^2(B^{\mathcal{M}^{(\nu)}}_{\rho_{0,\nu}}(Z_{0,\nu}))} \geq c \|N^{(\nu)}\|_{L^2(B^{\mathcal{M}^{(\nu)}}_{2\rho_{0,\nu}}(Z_{0,\nu}))}$ for some constant $c = c(n,m,q,M_0,N_0,C_e,C_h,\alpha_0) > 0$, and thus $\|\varphi\|_{L^2(B_2(0))} \geq c 2^{n/2+1}$ so that $\varphi$ is non-zero.  By Lemma~\ref{centered semicont2 lemma}, Theorem~\ref{freq mono centered thm}, and \eqref{lemma2_4 eqn1} 
\begin{gather*}
	N_{\varphi,0}(1) = \lim_{\nu\rightarrow\infty} \mathbf{I}_{N^{(\nu)},Z_{0,\nu}}(\rho_{0,\nu}) 
		\leq \lim_{\nu\rightarrow\infty} \mathbf{I}_{N^{(\nu)},Z_{0,\nu}}(2\rho_{0,\nu}) \leq \alpha, \\
	N_{\varphi,0}(1) = \lim_{\nu\rightarrow\infty} \mathbf{I}_{N^{(\nu)},Z_{\nu}}(\rho_{0,\nu}) \geq \alpha + \eta^2/4 , 
\end{gather*}
giving us a contradiction.

\noindent\textit{Proof of Conclusion~{\rm (ii)}.}  Let $\eta > 0$.  For $\nu = 1,2,3,\ldots$ let $T^{(\nu)}$, $\mathcal{M}^{(\nu)}$, and $N^{(\nu)}$ such that Hypothesis~\ref{M and N hypothesis} holds true with $T^{(\nu)},\mathcal{M}^{(\nu)},N^{(\nu)}$ in place of $T,\mathcal{M},N$.  Recalling conclusion~(i), for $\nu = 1,2,3,\ldots$ let $\delta_{\nu}\rightarrow 0^+$, $\alpha_{\nu} \in [2,\alpha_0]$, $Z_{\nu} \in \mathcal{B}^{(\geq 2)}_q(T^{(\nu)})$, $\rho_{0,\nu} \in (0,1]$, and $\rho_{\nu} \in (0,\delta_{\nu} \rho_{0,\nu}]$ such that $\mathbf{I}_{N^{(\nu)},Z_{\nu}}(0^+) \geq \alpha_{\nu}$ and 
\begin{equation}\label{lemma2_4 eqn2}
	\lim_{\nu\rightarrow\infty} \sup_{\sigma \in (0,\rho_{0,\nu}]} |\mathbf{I}_{N^{(\nu)},Z_{\nu}}(\sigma) - \alpha_{\nu}| = 0.
\end{equation}
Let $P_{\nu} = {\rm T}_{Z_{\nu}} \mathcal{M}^{(\nu)}$.  Let $\mathbf{e}_{\nu} : B_2(0,P_{\nu}) \rightarrow \mathcal{M}^{(\nu)}$ be the exponential map of $\mathcal{M}^{(\nu)}$ at $Z_{\nu}$.  It suffices to show that for infinitely many $\nu$, there exists a non-zero, average-free, homogeneous, locally Dirichlet energy minimizing $q$-valued function $\varphi^{(\nu)} \in W^{1,2}_{\rm loc}(P_{\nu},\mathcal{A}_q(P_{\nu}^{\perp}))$ such that $|\mathcal{N}_{\varphi^{(\nu)}}(0) - \alpha_{\nu}| < \eta^2$ and \eqref{lemma2_4 concl1} holds true with $T^{(\nu)},\mathcal{M}^{(\nu)},N^{(\nu)},\varphi^{(\nu)},Z_{\nu},\rho_{\nu}$ in place of $T,\mathcal{M},N,\varphi,Z,\rho$.

Let $\mathfrak{q}_{\nu}$ be a rotation of $\mathbb{R}^{n+m}$ such that $\mathfrak{q}_{\nu}(P_{\nu}) = \mathbb{R}^n \times \{0\}$.  Rescale, letting $\overline{T}^{(\nu)}$, $\overline{\mathcal M}^{(\nu)}$, and $\overline{N}^{(\nu)}$ be as in \eqref{centered blowup eqn1}.  Let $\overline{\mathbf e}_{\nu} : B_{1/\delta_{\nu}}(0) \rightarrow \overline{\mathcal M}^{(\nu)}$ be the exponential map of $\overline{\mathcal M}^{(\nu)}$ at the origin.  After passing to a subsequence, let $\alpha_{\nu} \rightarrow \alpha$ in $[2,\alpha_0]$ as $\nu\rightarrow\infty$.  By \eqref{lemma2_4 eqn2}, for each fixed $\sigma \in (1,\infty)$ and every sufficiently large $\nu$, $\mathbf{I}_{N^{(\nu)},Z_{\nu}}(2\rho_{\nu} \sigma) \leq \alpha_0+1$ and thus by Lemma~\ref{normal approx energy lemma} and Corollary~\ref{doubling centered cor} 
\begin{align*} 
	&\frac{1}{\sigma^{n+2}} \int_{B^{\overline{\mathcal M}^{(\nu)}}_{\sigma}(0)} |\overline{N}^{(\nu)}|^2 
		+ \frac{1}{\sigma^n} \int_{B^{\overline{\mathcal M}^{(\nu)}}_{\sigma}(0)} |D\overline{N}^{(\nu)}|^2 
	\\=\,& \frac{1}{(\rho_{\nu}\sigma)^{n+2}} \int_{B^{\mathcal{M}^{(\nu)}}_{\rho_{\nu}\sigma}(Z_{\nu})} |N^{(\nu)}|^2 
		+ \frac{1}{(\rho_{\nu}\sigma)^n} \int_{B^{\mathcal{M}^{(\nu)}}_{\rho_{\nu}\sigma}(Z_{\nu})} |DN^{(\nu)}|^2 
	\\ \leq\,& \frac{C}{(\rho_{\nu}\sigma)^{n+2}} \int_{B^{\mathcal{M}^{(\nu)}}_{2\rho_{\nu}\sigma}(Z_{\nu})} |N^{(\nu)}|^2
	\leq C (2\sigma)^{2e^{C \delta_{\nu}^{\gamma_3}} (\alpha_0+1) - 2} \widehat{E}_{\nu}^2 , 
\end{align*}
where $C = C(n,m,q,M_0,N_0,C_e,C_h) \in (0,\infty)$ are constants.  Following the blow-up procedure from Subsection~\ref{sec:centered blowup sec}, after passing to a subsequence there exists an average-free, locally Dirichlet energy minimizing $q$-valued function $\overline{\varphi} \in W^{1,2}_{\rm loc}(\mathbb{R}^n,\mathcal{A}_q(\mathbb{R}^m))$ such that 
\begin{equation}\label{lemma2_4 eqn3}
	\frac{\overline{N}^{(\nu)}(\overline{\mathbf e}_{\nu}^{-1}(x))}{\widehat{E}_{\nu}} \rightarrow \overline{\varphi}(x)
\end{equation}
pointwise $\mathcal{L}^n$-a.e.~in $\mathbb{R}^n$ and in $L^2(B_{\sigma}(0),\mathcal{A}_q(\mathbb{R}^m))$ for all $\sigma > 0$ as $\nu \rightarrow \infty$ 
\begin{equation*}
	\frac{|D(\overline{N}^{(\nu)} \circ \overline{\mathbf e}_{\nu}^{-1})(X)|}{\widehat{E}_{\nu}} \rightarrow |D\overline{\varphi}(X)| 
\end{equation*}
in $L^2(B_{\sigma}(0),\mathcal{A}_q(\mathbb{R}^m))$ for all $\sigma > 0$ as $\nu \rightarrow \infty$.  By Lemma~\ref{centered semicont2 lemma}, it follows from \eqref{lemma2_4 eqn2} that 
\begin{equation*}
	N_{\overline{\varphi},0}(\sigma) = \lim_{\nu\rightarrow\infty} \mathbf{I}_{N^{(\nu)},Z_{\nu}}(\rho_{\nu} \sigma) = \alpha 
\end{equation*}
for all $\sigma \in (0,\infty)$.  Therefore, $\overline{\varphi}$ is homogeneous degree $\alpha$.  Let $\varphi^{(\nu)}(X) = \widehat{E}_{\nu} \rho_{\nu}^{1-\alpha} \,\mathfrak{q}_{\nu}^{-1} \overline{\varphi}(\mathfrak{q}_{\nu} X)$ for each $X \in P_{\nu}$ (where we identify $\mathbb{R}^n \cong \mathbb{R}^n \times \{0\}$).  By \eqref{lemma2_4 eqn3}
\begin{equation*}
	\lim_{\nu\rightarrow\infty} \frac{1}{\widehat{E}_{\nu}^2 \rho_{\nu}^{n+2}} \int_{B^{\mathcal{M}^{(\nu)}}_{\rho_{\nu}}(Z_{\nu})} 
		\mathcal{G}(N^{(\nu)} \circ \mathbf{e}_{\nu}^{-1},\varphi^{(\nu)})^2 
	= \lim_{\nu\rightarrow\infty} \frac{1}{\widehat{E}_{\nu}^2} \int_{B^{\overline{\mathcal M}^{(\nu)}}_1(0)} 
		\mathcal{G}(\overline{N}^{(\nu)} \circ \overline{\mathbf e}_{\nu}^{-1},\overline{\varphi})^2 = 0 .
\end{equation*}
That is, \eqref{lemma2_4 concl1} holds true with $\mathcal{M}^{(\nu)},N^{(\nu)},\varphi^{(\nu)},Z_{\nu},\rho_{\nu}$ in place of $\mathcal{M},N,\varphi,Z,\rho$ for all sufficiently large $\nu$.

\noindent\textit{Proof of Conclusion~{\rm (iii)}.}  We argue exactly as we did for conclusion~(ii), except we observe the following.  If $\dim \op{spine} \overline{\varphi} = n-2$, then we already know that \eqref{lemma2_4 concl2} holds true with $\mathcal{M}^{(\nu)},N^{(\nu)},\varphi^{(\nu)},Z_{\nu},$ $\rho_{\nu}$ in place of $\mathcal{M},N,\varphi,Z,\rho$ for all sufficiently large $\nu$.  If instead $\dim \op{spine} \overline{\varphi} \leq n-3$, let $L$ be an $(n-3)$-dimensional linear subspace such that $\op{spine} \overline{\varphi} \subseteq L$ and set $L_{\nu} = \mathfrak{q}_{\nu}^{-1}(L)$.  Suppose that there exists $X_{\nu} \in \mathcal{B}^{(\geq 2)}_q(T^{(\nu)}) \cap \overline{B^{\mathcal{M}^{(\nu)}}_{\rho_{\nu}}(Z_{\nu})}$ such that $\mathbf{I}_{N^{(\nu)},X_{\nu}}(0^+) \geq \alpha_{\nu}$ and $\op{dist}(X_{\nu}, Z_{\nu} + L_{\nu}) \geq \eta\rho_{\nu}$.  After passing to a subsequence let $\mathfrak{q}_{\nu}(X_{\nu}-Z_{\nu})/\rho_{\nu} \rightarrow (\xi,0)$ as $\nu\rightarrow\infty$ for some $\xi \in \mathbb{R}^n$.  Then $\op{dist}(\xi,L) \geq \eta$, but by Lemma~\ref{centered semicont2 lemma} 
\begin{equation*}
	\mathcal{N}_{\overline{\varphi}}(\xi) \geq \limsup_{\nu\rightarrow\infty} \mathbf{I}_{N^{(\nu)},X_{\nu}}(0^+) \geq \alpha ,
\end{equation*}
giving us a contradiction.  Therefore, \eqref{lemma2_4 concl3} must hold true with $T^{(\nu)},\mathcal{M}^{(\nu)},N^{(\nu)},Z_{\nu},\rho_{\nu},L_{\nu}$ in place of $T,\mathcal{M},N,Z,\rho,L$ for all sufficiently large $\nu$.
\end{proof}

\section{Bounding the distance of an area-minimizing current to a union of non-intersecting single-valued 
$C^2$-graphs}\label{sec:relharm}

In this section, we obtain local estimates on the uniform distance of a locally area minimizing rectifiable current $T$ to the graph of a harmonic multi-valued function $\varphi$, analogous to~\cite[Theorem~2.15]{KrumWicb} and~\cite[Theorem~4.19]{KrumWicc}.  In particular, let $T$ be an $n$-dimensional locally area minimizing rectifiable current of $\mathbf{C}_{11\sqrt{n}/2}(0)$ and $\mathcal{M}$ be the center manifold of $T$ as in Hypothesis~\ref{M and N hypothesis}.  Let $Z \in \mathcal{M} \cap \mathbf{C}_1(0)$ and $\rho \in (0,2]$.  Let $\varphi : B^{\mathcal M}_{\rho}(Z) \rightarrow \mathcal{A}_q({\rm T}\mathcal{M}^{\perp})$ be a $q$-valued function such that $\varphi(X) = \sum_{l=1}^q \llbracket \varphi_l(X) \rrbracket$ for each $X \in B^{\mathcal M}_{\rho}(Z)$ where $\varphi_l : B^{\mathcal M}_{\rho}(Z) \rightarrow {\rm T}\mathcal{M}^{\perp}$ are single-valued $C^{2}$ functions with small $C^{2}$ norm, such that for $i \neq j$, either $\varphi_{i} \equiv \varphi_{j}$ or the graphs of 
$\varphi_{i}$, 
$\varphi_{j}$ are disjoint.  We wish to prove the following local estimate for the uniform distance of $T$ to the graph of $\varphi$:

\begin{theorem}\label{relharm separation main thm} 
Let $M_0,N_0,C_e,C_h$ be as in Hypothesis~\ref{DeLSpa constants hypothesis}.  For every $\kappa \in (0,\infty)$ and $\gamma \in (0,1)$ there exists $\varepsilon_2 = \varepsilon_2(n,m,q,M_0,N_0,C_e,C_h) \in (0,1)$ and $\varepsilon = \varepsilon(n,m,q,M_0,N_0,C_e,C_h,\kappa,\gamma) \in (0,1)$ such that following holds true.  Let $T, \mathscr{W}, \mathcal{M}, N$ be as in Hypothesis~\ref{M and N hypothesis}.  Let $Z \in \mathcal{M} \cap \mathbf{C}_1(0)$ and $\rho \in (0,2]$ such that 
\begin{equation}\label{relharm separation main go condition}
	\sup_{X \in B^{\mathcal M}_{(5+3\gamma)\rho/8}(Z)} \overline{R}(X) < 4^{-q-3} (1-\gamma) \rho , 
\end{equation}
where $\overline{R}$ is as in Definition~\ref{crit rad defn}.  Let $\varphi : B^{\mathcal M}_{\rho}(Z) \rightarrow \mathcal{A}_q(\mathbb{R}^{n+m})$ be given by 
\begin{equation*}
	\varphi(X) = \sum_{l=1}^q \llbracket \varphi_l(X) \rrbracket
\end{equation*} 
for each $X \in B^{\mathcal M}_{\rho}(Z)$, where $\varphi_l \in C^2(B^{\mathcal M}_{\rho}(Z),\mathbb{R}^{n+m})$ are single-valued functions such that 
\begin{gather}
	\label{relharm separation main varphi small} \rho^{-1} \|\varphi_l\|_{C^0(B^{\mathcal M}_{\rho}(Z))} + \|D\varphi_l\|_{C^0(B^{\mathcal M}_{\rho}(Z))} 
		+ \rho \|{\rm Hess}\,\varphi_l\|_{C^0(B^{\mathcal M}_{\rho}(Z))} < \varepsilon , \\
	\label{relharm separation main avgfree} \frac{1}{q} \sum_{l=1}^q \varphi_l(X) = 0 \text{ for all } X \in B^{\mathcal M}_{\rho}(Z), \\	
	\label{relharm separation main normal} \varphi_l(X) \in {\rm T}_X \mathcal{M}^{\perp} \text{ for all } X \in B^{\mathcal M}_{\rho}(Z) 
\end{gather}
for all $l \in \{1,2,\ldots,q\}$ and 
\begin{equation}\label{relharm separation main wellsep} 
	\rho \sup_{B^{\mathcal M}_{\rho}(Z)} |D\varphi_{k} - D\varphi_{l}| \leq \kappa \inf_{B^{\mathcal M}_{\rho}(Z)} |\varphi_{k} - \varphi_{l}| 
\end{equation}
for all $k,l \in \{1,2,\ldots,q\}$.  Then 
\begin{align}\label{relharm separation main concl}
	\sup_{X \in \op{spt} T \cap \mathbf{C}^{\mathcal M}_{\gamma\rho}(Z)} \op{dist}^2(X, \op{spt} \mathbf{G}_{\varphi}) 
	\leq \frac{C}{\rho^n} \int_{B^{\mathcal M}_{\rho}(Z)} \mathcal{G}(N,\varphi)^2 
		+ C \rho^4 \sum_{l=1}^q \|\Delta_{\mathcal M} \varphi_l\|_{C^0(B^{\mathcal M}_{\rho}(Z))}^2 
		&\\+ \frac{C}{\rho^n} \,\|D\varphi\|_{C^0(B^{\mathcal M}_{\rho}(Z))} \int_{B^{\mathcal M}_{\rho}(Z)} |N|^2  
		+ \frac{C}{\rho^{n-\varsigma}} \int_{B^{\mathcal M}_{\rho}(Z)} |N|^2 & \nonumber 
\end{align}
for some constants $\varsigma = \varsigma(n,m,q) \in (0,1)$ and $C = C(n,m,q,M_0,N_0,C_e,C_h,\gamma) \in (0,\infty)$.  Here ${\rm Hess}\,\varphi_l$ is the Hessian of $\varphi_l$ on $\mathcal{M}$ (as in \eqref{app hessian defn}) and $\Delta_{\mathcal M} \varphi_l$ is the Laplacian $\varphi_l$ on $\mathcal{M}$ (as in \eqref{app laplace defn}).
\end{theorem}

The assumption \eqref{relharm separation main go condition} allows us to use  estimates for the normal approximation map $N$, namely \eqref{normal approx est1}, \eqref{normal approx est2}, \eqref{normal approx est3}, and \eqref{normal approx est4}.  Notice that \eqref{normal approx est4} is a quantitative statement that $N$ is approximately average-free in $B^{\mathcal M}_{\gamma\rho}(Z)$ (where $\gamma \in (0,1)$).  Thus we naturally assume in \eqref{relharm separation main avgfree} that $\varphi$ is also average-free.  \eqref{relharm separation main varphi small} asserts that each function $\varphi_l$ has small $C^2$-norm on $B^{\mathcal M}_{\rho}(Z)$.  In \eqref{relharm separation main normal}, we assume that $\varphi(X)$ takes values in ${\rm T}_X \mathcal{M}^{\perp}$ for each $X \in B^{\mathcal M}_{\rho}(Z)$; this condition can be relaxed using a change of coordinates result Lemma~\ref{app lemma}.  \eqref{relharm separation main wellsep} implies that for each $i \neq j$ either $\varphi_i \equiv \varphi_j$ on $B^{\mathcal M}_{\rho}(Z)$ (allowing for higher integer multiplicity) or $\op{spt} \mathbf{G}_{\varphi_i} \cap \op{spt} \mathbf{G}_{\varphi_i} = \emptyset$.  In the latter case, \eqref{relharm separation main wellsep} is a quantitative statement that the graphs of $\varphi_i$ are nearly parallel.  The conclusion \eqref{relharm separation main concl} asserts that the $L^{\infty}$-distance of $T$ to $\op{spt} \mathbf{G}_{\varphi}$ is bounded above by an $L^2$-distance of $N$ to $\varphi$.  The remaining terms on the right-hand side of \eqref{relharm separation main concl} control $\Delta_{\mathcal M} \varphi_l$ as well as error terms arising from the area formula and the estimates for the normal approximation $N$.

Theorem~\ref{relharm separation main thm} will be used in an essential way, together with Almgren's strong Lipschitz approximation theorem, to establish a graphical representation of $T$ relative to homogeneous harmonic functions in Section~\ref{sec:graph rep} below.  Following the strategy employed in~\cite{KrumWicc}, we obtain Theorem~\ref{relharm separation main thm} not as a direct corollary of~\cite[Theorem~2.15]{KrumWicb} but rather by adapting the proof of~\cite[Theorem~2.15]{KrumWicb} to the present context.  There are additional considerations that must be taken into account.  In particular:
\begin{enumerate}[itemsep=2mm,topsep=0mm]
	\item[(i)]  The local estimates from Section~\ref{sec:normal approx est sec} for the normal approximation map $N$ in $B^{\mathcal M}_{\rho}(Z)$ require that $\rho > \overline{R}(Z)$, where $\overline{R}$ is as in Definition~\ref{crit rad defn}.  To apply these estimates, we must keep track of the size of $\overline{R}$ thoroughout the proof of Theorem~\ref{relharm separation main thm}, and this is facilitated by assumption \eqref{relharm separation main go condition} above.
	
	\item[(ii)]  Recall that \eqref{normal approx est4} is a quantitative statement that $N$ is approximately average-free and that \eqref{relharm separation main avgfree} asserts that $\varphi$ is average free. The averages of the $q$-valued function $N$ and $\varphi$ are computed with their multiplicities as weights.  Unlike~\cite{KrumWicb}\cite{KrumWicc}, we do not consider the $L^2$-distance $\int_{\mathbf{C}^{\mathcal M}_{\rho}(Z)} \op{dist}^2(X,\op{spt}\mathbf{G}_{\varphi}) \,d\|T\|(X)$, which disregards the multiplicity of the components of $\varphi$ and leads to ambiguity in the value of $\varphi_a$.  Rather, we consider the $L^2$-distance $\int_{B^{\mathcal M}_{\rho}(Z)} \mathcal{G}(N,\varphi)^2$ which accounts for the multiplicities.
	
	\item[(iii)]  \cite[Theorem~3.15]{KrumWicb} uses an inductive argument in which we decompose $T$ into components near the graphs of the functions $\varphi_i$ and then apply inductively~\cite[Theorem~3.15]{KrumWicb} to each of the components of $T$.  We do not expect this inductive argument to preserve \eqref{normal approx est4} in the sense that we do not expect the normal approximation maps of the components of $T$ to be approximately average-free.  We proceed by  \emph{inductively assuming} $T$ separates into finitely many locally area minimizing rectifiable currents $T_j$  and $\varphi$ separates into corresponding $q_j$-valued functions $\varphi_j$ to which we compare $T_j$, see Hypothesis~\ref{relharm separation1 hyp} and Hypothesis~\ref{relharm grand hyp} below; the induction is on the number $p$ of such pieces $T_{j}$ that make up $T$, with the base case of induction being $p=q$ and the theorem corresponding to $p=1$. In Lemma~\ref{relharm normal approx lemma}, we construct normal approximation maps $N_j$ of $T_j$ with estimates analogous to \eqref{normal approx est1} and \eqref{normal approx est3}, as required for the proof of the theorem.
\end{enumerate}

\subsection{Normal approximation maps of components of area-minimizing currents}\label{sec:relharm energy sec}  Let $M_0,N_0,$ $C_e,C_h,\varepsilon_2$ be as in Hypothesis~\ref{DeLSpa constants hypothesis} and Lemma~\ref{normal approx energy lemma}.  In Sections~\ref{sec:relharm energy sec}--\ref{sec:relharm separation2 sec}, given a suitably small constant $\varepsilon \in (0,1)$ we shall assume the following: 

\begin{hypothesis}\label{relharm separation1 hyp} 
Let $\mathbf{m}_0,T, \mathscr{W}, \mathcal{M}, N, \mathcal{K}$ be as in Hypothesis~\ref{M and N hypothesis}.  Let $Z \in \mathcal{M} \cap \mathbf{C}_1(0)$ and $\rho \in (0,2]$ such that $\mathbf{m}_0 \rho < \varepsilon$.  Let $p,q_1,q_2,\ldots,q_p$ be positive integers with $\sum_{j=1}^p q_j = q$.  For each $j \in \{1,2,\ldots,p\}$, let $T_j$ be $n$-dimensional locally area minimizing rectifiable currents of $\mathbf{C}^{\mathcal M}_{\rho}(Z)$ such that 
\begin{gather}\label{relharm separation1 hyp eqn1} 
	(\partial T_j) \llcorner \mathbf{C}^{\mathcal M}_{\rho}(Z) = 0, \quad \mathbf{p}_{\#} T_j = q_j \llbracket B^{\mathcal M}_{\rho}(Z) \rrbracket , \\
	T \llcorner \mathbf{C}^{\mathcal M}_{\rho}(Z) = \sum_{j=1}^p T_j , \quad 
		\|T\| \llcorner \mathbf{C}^{\mathcal M}_{\rho}(Z) = \sum_{j=1}^p \|T_j\| . \nonumber 
\end{gather}
For each $j \in \{1,2,\ldots,p\}$, let $\varphi_j \in C^2(B^{\mathcal M}_{\rho}(Z),\mathbb{R}^{n+m})$ be single-valued functions such that 
\begin{gather}
	\label{relharm separation1 small norm} \rho^{-1} \|\varphi_j\|_{C^0(B^{\mathcal M}_{\rho}(Z))} + \|D\varphi_j\|_{C^0(B^{\mathcal M}_{\rho}(Z))} 
		+ \rho \|{\rm Hess}\,\varphi_j\|_{C^0(B^{\mathcal M}_{\rho}(Z))} < \varepsilon , \\
	\label{relharm separation1 avgfree} \frac{1}{q} \sum_{j=1}^p q_j \varphi_j(X) = 0 \text{ for each } X \in B^{\mathcal M}_{\rho}(Z) . 
\end{gather}
(Here we do not require $\varphi(X)$ to take values in ${\rm T}_X \mathcal{M}^{\perp}$ for $X \in B^{\mathcal M}_{\rho}(Z)$.)
\end{hypothesis}

In Lemma~\ref{relharm normal approx lemma} we introduce normal approximation maps $N_j$ of the currents $T_j$.  Note that the maps $N_j$ are defined on $B^{\mathcal M}_{\gamma\rho}(Z)$ where $0 < \gamma < 1$, and we must bear this in mind when defining the excesses $\mathcal{E}_{Z,\gamma\rho,\rho}$ and $E_{Z,\gamma\rho,\rho}$ in \eqref{relharm excess defn1}, \eqref{relharm excess defn2}, and \eqref{relharm excess defn3} below. 

\begin{lemma}\label{relharm normal approx lemma} 
Let $M_0,N_0,C_e,C_h,\varepsilon_2$ be as in Hypothesis~\ref{DeLSpa constants hypothesis} and Lemma~\ref{normal approx energy lemma}.  Let $\gamma \in (0,1)$, $Z \in \mathcal{M} \cap \mathbf{C}_1(0)$, and $\rho \in (\overline{R}(Z),2]$, where $\overline{R}(Z)$ be as in Definition~\ref{crit rad defn}.  Let $p,q_1,q_2,\ldots,q_p$ be positive integers such that $\sum_{j=1}^p q_j = q$.  For each $j \in \{1,2,\ldots,p\}$, let $T_j$ be $n$-dimensional locally area minimizing rectifiable currents of $\mathbf{C}^{\mathcal M}_{\rho}(Z)$ such that \eqref{relharm separation1 hyp eqn1} holds true.  Then for each $j \in \{1,2,\ldots,J\}$ there exist a closed set $\mathcal{K}_j \subseteq B^{\mathcal M}_{\gamma\rho}(Z)$ and a Lipschitz $q_j$-valued function $N_j : B^{\mathcal M}_{\gamma\rho}(Z) \rightarrow \mathcal{A}_{q_j}({\rm T}\mathcal{M}^{\perp})$ such that 
\begin{gather}\label{relharm normal approx concl1}
	T_j \llcorner \mathbf{p}^{-1}(\mathcal{K}_j) = \mathbf{G}_{N_j} \llcorner \mathbf{p}^{-1}(\mathcal{K}_j) , \quad 
	N(X) = \sum_{j=1}^p N_j(X) \text{ for all } X \in \mathcal{K}_j , \\
	N_j(x,\bs{\varphi}(x)) = q_j \llbracket 0 \rrbracket \text{ for all } x \in \bs{\Gamma} , \nonumber
\end{gather}
and for each $L \in \mathscr{W}$ such that the Whitney region $\mathcal{L}$ associated with $L$ (as in \eqref{Whitney region defn}) satisfies $B^{\mathcal M}_{\gamma\rho}(Z) \cap \mathcal{L} \neq \emptyset$
\begin{gather}\label{relharm normal approx concl2} 
	\sup_{B^{\mathcal M}_{\gamma\rho}(Z) \cap \mathcal{L}} |N_j| \leq C \mathbf{m}_0^{1/(2n)} \ell(L)^{1+\beta_2}, \quad 
	\sup_{B^{\mathcal M}_{\gamma\rho}(Z) \cap \mathcal{L}} |DN_j| \leq C \mathbf{m}_0^{\gamma_2} \ell(L)^{\gamma_2}, \\
	\mathcal{H}^n(B^{\mathcal M}_{\gamma\rho}(Z) \cap \mathcal{L} \setminus \mathcal{K}_j) 
		+ \|T_j -\mathbf{G}_{N_j}\|(\mathbf{p}^{-1}(B^{\mathcal M}_{\gamma\rho}(Z) \cap \mathcal{L})) 
		\leq C \mathbf{m}_0^{1+\gamma_2} \ell(L)^{n+2+\gamma_2}, \nonumber \\
	\int_{B^{\mathcal M}_{\gamma\rho}(Z) \cap \mathcal{L}} |DN_j|^2 \leq C \mathbf{m}_0 \ell(L)^{n+2-2\delta_2} , \nonumber
\end{gather}
where $C = C(n,m,q,M_0,N_0,C_e,C_h,\gamma) \in (0,\infty)$ is a constant.
\end{lemma}

\begin{remark}{\rm 
For each $j \in \{1,2,\ldots,p\}$, let $\widetilde{N}_j : B^{\mathcal M}_{\gamma\rho}(Z) \cap \mathcal{K} \to \mathcal{A}_{q_j}({\rm T}\mathcal{M})$ be the multi-valued function whose graph is $T_j \llcorner \cap \mathbf{p}^{-1}(B^{\mathcal M}_{\gamma\rho}(Z) \cap \mathcal{K})$.  One plausible approach to proving Lemma~\ref{relharm normal approx lemma} is to extend $\widetilde{N}_j$ to a suitable multi-valued function $N_j$ on $B^{\mathcal M}_{\gamma\rho}(Z)$.  Notice that any such extension would have to be constructed in such a way that the local estimate \eqref{relharm normal approx concl2} on each Whitney region $\mathcal{L}$ intersecting $B^{\mathcal M}_{\gamma\rho}(Z)$ holds true.  
The main difficulty with this approach is verifying that $\widetilde{N}_j$ are Lipschitz multi-valued functions.  If the components $T_j$ of $T \llcorner \mathbf{C}^{\mathcal M}_{\rho}(Z)$ were sufficiently well-separated to guarantee that 
\begin{equation}\label{relharm normal approx rmk eqn1}
    \mathcal{G}(N(X),N(Y))^2 = \sum_{j=1}^p \mathcal{G}(\widetilde{N}_j(X),\widetilde{N}_j(Y))^2 \text{ for all } X,Y \in B^{\mathcal M}_{\gamma\rho}(Z) \cap \mathcal{K} ,
\end{equation}
then it would follow that $\op{Lip} \widetilde{N}_j \leq \op{Lip} N$; however, \eqref{relharm normal approx rmk eqn1} is not evident.  Instead, we proceed by modifying the construction of the normal approximation map (\cite[Section~6]{DeLSpa2}, \cite[Theorem~4.33]{Almgren}).
}\end{remark}

\begin{proof}
Let $L \in \mathscr{W}$ such that the Whitney region $\mathcal{L}$ associated with $L$ (as in \eqref{Whitney region defn}) satisfies $B^{\mathcal M}_{\gamma\rho}(Z) \cap \mathcal{L} \neq \emptyset$.  Notice that $\ell(L)$ might be large enough that $\mathbf{B}_{64r_L}(p_L) \not\subset \mathbf{C}^{\mathcal M}_{\rho}(Z)$.  Nonetheless, we can subdivide $L$ into $2^{nk}$ congruent cubes $H$ of side length $2^{1-k} \ell(L)$, where $k = k(L)$ is the smallest non-negative number such that $2^{7-k} r_L \leq (1-\gamma)\,\rho/4$ (and in particular $k = 0$ if $2^7 r_L \leq (1-\gamma)\,\rho/4$).  By \eqref{freq go condition}, $r_L = M_0 \sqrt{n} \ell(L) < M_0 \rho/64$ and thus $0 \leq k \leq k_0$, where $k_0 = k_0(M_0,\gamma)$ is the smallest positive integer such that $2^{1-k_0} M_0 \leq (1-\gamma)/4$.  Let $\widetilde{\mathscr W}_0$ denote the set of the cubes $H$ obtained by subdividing a cube $L \in \mathscr{W}$ in this way.  For each cube $H \in \widetilde{\mathscr W}_0$, let \begin{equation}\label{relharm normal approx eqn01}
    \mathcal{H}' = \{ (x,\bs{\varphi}(x)) : x \in \widetilde{J} \} 
\end{equation}
where $\widetilde{J}$ is the cube concentric to $H$ with $\ell(\widetilde{J}) = \frac{9}{8} \ell(H)$.  We let $\widetilde{\mathcal W}$ be the set of all cubes $H$ obtained by subdividing a cube $L \in \mathscr{W}$ such that the set $\mathcal{H}'$ associated with $H$ as in \eqref{relharm normal approx eqn01} satisfies $\mathcal{H}' \cap B^{\mathcal M}_{(1+\gamma)\rho/2}(Z) \neq \emptyset$.

Let $H \in \widetilde{\mathscr W}$.  Let $\mathcal{H}'$ be as in \eqref{relharm normal approx eqn01} and let $L \in \mathscr{W}$ be the (unique) cube with $H \subseteq L$.  Let $j \in \{1,2,\ldots,p\}$.  Let $x_H$ be the center of $H$, $p_{j,H} = (x_H,y_{j,H}) \in \op{spt} T_j$, and $r_H = M_0 \sqrt{n} \ell(H)$.  By~\cite[Proposition~2.2(ii)]{DeLSpa2}, $\mathcal{H}' \cap B^{\mathcal M}_{(1+\gamma)\rho/2}(Z) \neq \emptyset$ and $2^7 r_H \leq (1-\gamma)\,\rho/2$
\begin{equation*}
	\op{spt} T \cap \mathbf{C}_{32 r_H}(p_{j,H},P_L) 
	\subset \mathbf{C}^{\mathcal M}_{32 r_H + C \mathbf{m}_0^{\frac{1}{2n}} r_H^{1+\beta_2}}(x_H,\bs{\varphi}(x_H))
	\subset \mathbf{C}^{\mathcal M}_{(1+\gamma)\rho/2 + 64 r_H}(Z)
	\subset \mathbf{C}^{\mathcal M}_{\rho}(Z) , 
\end{equation*}
where here and throughout this proof $C = C(n,m,q,M_0,N_0,C_e,C_h,\gamma) \in (0,\infty)$ are constants.  Moreover, by~\cite[Proposition~4.1]{DeLSpa2} 
\begin{equation*}
    \op{spt} T_j \cap \mathbf{C}_{32 r_H}(p_{j,H},P_L) \subset \mathbf{B}_{32r_H + \sqrt{n} \ell(L) + C\mathbf{m}_0^{\frac{1}{2n}} r_L^{1+\beta_2}}(p_L) \subset \mathbf{B}_{64r_L}(p_L) .
\end{equation*}
Thus by \eqref{relharm separation1 hyp eqn1} and the construction of $\mathscr{W}$ in Definition~\ref{DeLSpa whitney defn} 
\begin{align}\label{relharm normal approx eqn1}
	\mathbf{E}(T_j,\mathbf{C}_{32 r_H}(p_{j,H},P_L),P_L) \leq\,& C \mathbf{m}_0 \ell(L)^{2-2\delta_2} , \\ 
	\mathbf{h}(T_j,\mathbf{C}_{32 r_H}(p_{j,H},P_L),P_L) \leq\,& C \mathbf{m}_0^{\frac{1}{2n}} \ell(L)^{1+\beta_2} . \nonumber 
\end{align}
Thus arguing as in~\cite[Proposition~4.2]{DeLSpa2} there exists a Lipschitz approximation $f_{j,H} : B_{8r_H}(p_{j,H},P_L) \rightarrow \mathcal{A}_{q_j}(P_L^{\perp})$ and a set $K_{j,H} \subseteq B_{8r_H}(p_{j,H},P_L)$ such that 
\begin{gather}\label{relharm normal approx eqn2}
	T_j \llcorner \pi_{P_L}^{-1}(K_{j,H}) = \mathbf{G}_{f_{j,H}} \llcorner \pi_{P_L}^{-1}(K_{j,H}) , \quad 
	\op{Lip} f_{j,H} \leq C \mathbf{m}_0^{\gamma_1} \ell(L)^{2\gamma_1 (1-\delta_2)} , \\
	\mathcal{H}^n(B_{8r_H}(p_{j,H},P_L) \setminus K_{j,H}) + \|T_j\|(\pi_{P_L}^{-1}(B_{8r_H}(p_{j,H},P_L) \setminus K_{j,H})) 
		\leq C \mathbf{m}_0^{1+\gamma_1} \ell(L)^{n+2(1+\gamma_1)(1-\delta_2)} \nonumber 
\end{gather}
Let 
\begin{equation*}
	\mathscr{D}_j(H) = (\op{spt} T_j \cup \op{spt} \mathbf{G}_{f_{j,H}}) \cap \pi_{P_L}^{-1}(B_{8r_H}(p_{j,H},P_L) \setminus K_{j,H})
\end{equation*}
be the set of points where the supports of $T_j$ and the graph of $f_{j,H}$ disagree. There is a Lipschitz function $N_{j,H} : \mathcal{H}' \rightarrow \mathcal{A}_{q_j}({\rm T}\mathcal{M}^{\perp})$ such that $\mathbf{G}_{N_{j,H}} = \mathbf{G}_{f_{j,H}} \llcorner \mathbf{p}^{-1}(\mathcal{H}')$ (see e.g.\ ~\cite[Theorem~5.1]{DeLSpa0}).  

Recall that $\{(x,\bs{\varphi}(x)) : x \in \mathbf{\Gamma} \}$ is the contact set.  Set 
\begin{equation*}
	\mathcal{K}_j = B^{\mathcal M}_{\gamma\rho}(Z) \cap \mathcal{K} \setminus \left( \bigcup_{H \in \widetilde{\mathscr W}} \Big( \mathcal{H}' 
		\cap \bigcup \{ \mathbf{p}(\mathscr{D}_j(M)) : M \in \widetilde{\mathscr W} \text{ with } M \cap H \neq \emptyset \} \Big) \right) ,
\end{equation*}
where $\mathcal{H}'$ is associated with $H$ as in \eqref{relharm normal approx eqn01}.  That way $B^{\mathcal M}_{\gamma\rho}(Z) \cap \{(x,\bs{\varphi}(x)) : x \in \mathbf{\Gamma} \} \subseteq \mathcal{K}_j$ and $\mathbf{G}_{N_{j,H}} \llcorner \mathbf{p}^{-1}(\mathcal{K}_j \cap \mathcal{H}') = T \llcorner \mathbf{p}^{-1}(\mathcal{K}_j \cap \mathcal{H}')$ for all $H \in \widetilde{\mathscr W}$.  
Let $L \in \mathcal{W}$ and let $\mathcal{L}$ be the Whitney region associated with $L$ (as in \eqref{Whitney region defn}).  Suppose that $H,M \in \widetilde{\mathcal W}$ such that $\mathcal{H}' \cap \mathcal{L} \neq \emptyset$ and $M \cap H \neq \emptyset$, where $\mathcal{H}'$ is as in \eqref{relharm normal approx eqn01}.  Then there exists unique cubes $L_1,L_2 \in \mathcal{W}$ such that $H \subseteq L_1$ and $M \subseteq L_2$.  Moreover, $L \cap L_1 \neq \emptyset$ and $L_1 \cap L_2 \neq \emptyset$, from which it follows that $2^{-k_0-2} \ell(L) \leq 2^{-k_0}\ell(L_2) \leq \ell(M) \leq \ell(L_2) \leq 4 \ell(L)$.  Hence by \eqref{relharm normal approx eqn2}
\begin{align*}
	\mathcal{H}^n(B^{\mathcal M}_{\gamma\rho}(Z) \cap \mathcal{L} \setminus \mathcal{K}_j)
	\leq\,& \mathcal{H}^n(B^{\mathcal M}_{\gamma\rho}(Z) \cap \mathcal{L} \setminus \mathcal{K}) \\&+ \sum_{H \in \widetilde{\mathscr W}}\,\sum_{\{ M \in \widetilde{\mathscr W} : 
		M \cap H \neq \emptyset \}} \mathcal{H}^n(B^{\mathcal M}_{\gamma\rho}(Z) \cap \mathcal{L} \cap \mathcal{H}' \cap \mathbf{p}(\mathscr{D}_j(M)))
	\\ \leq\,& C \mathbf{m}_0 \ell(L)^{n+2-2\delta_2} 
\end{align*}
and similarly $\|T_j -\mathbf{G}_{N_j}\|(\mathbf{p}^{-1}(B^{\mathcal M}_{\gamma\rho}(Z) \cap \mathcal{L})) \leq C \mathbf{m}_0 \ell(L)^{n+2-2\delta_2}$.  Define $N_j |_{\mathcal{K}_j}$ by $N_j = q_j \llbracket 0 \rrbracket$ on $B^{\mathcal M}_{\gamma\rho}(Z) \cap \{(x,\bs{\varphi}(x)) : x \in \mathbf{\Gamma} \}$ and $N_j = N_{j,H}$ on $\mathcal{K}_j \cap \mathcal{H}'$ for all $H \in \widetilde{\mathscr W}$, where $\mathcal{H}'$ is as in \eqref{relharm normal approx eqn01}.  It follows from \eqref{relharm normal approx eqn1} and \eqref{relharm normal approx eqn2}  that 
\begin{equation*}
	\sup_{\mathcal{K}_j \cap \mathcal{H}} |N_j| \leq C \mathbf{m}_0 \ell(L)^{2-2\delta_2} , \quad 
	\op{Lip} (N_j |_{\mathcal{K}_j \cap \mathcal{H}}) \leq C \mathbf{m}_0^{\gamma_2} \ell(L)^{\gamma_2} 
\end{equation*}
for all $H \in \widetilde{\mathscr W}$, where $L \in \mathcal{W}$ is the unique cube with $H \subseteq L$ and $\mathcal{H}$ be the graph of $\bs{\varphi}$ over the cube $J$ with the same center as $H$ and $\ell(J) = \frac{17}{16} \ell(H)$.  By inductively applying the Lipschitz extension theorem~\cite[Theorem~1.7]{DeLSpaDirMin} (\cite[Theorem~1.3(2)]{Almgren}) on tubular neighborhoods of the $k$-skeleton of $\widetilde{\mathscr W}$ for $k = 0,1,2,\ldots,n$, construct a Lipschitz function $\overline{N}_j : B^{\mathcal M}_{\gamma\rho}(Z) \rightarrow \mathcal{A}_{q_j}(\mathbb{R}^{n+m})$ such that $\overline{N}_j = N_j$ on $\mathcal{K}_j$ and \eqref{relharm normal approx concl2} holds true with $\overline{N}_j$ in place of $N_j$.  Define $N_j : B^{\mathcal M}_{\gamma\rho}(Z) \rightarrow \mathcal{A}_{q_j}({\rm T}\mathcal{M}^{\perp})$ by $N_j(X) = \sum_{l=1}^q \llbracket \pi_{{\rm T}_X \mathcal{M}^{\perp}}(\overline{N}_{j,l}(X)) \rrbracket$ for each $X \in B^{\mathcal M}_{\gamma\rho}(Z)$, where $\overline{N}_j(X) = \sum_{l=1}^q \llbracket \bar N_{j,l}(X) \rrbracket$.  Note that this procedure does not change the values of $N_j |_{\mathcal{K}_j}$.  It is readily verified that $N_j$ satisfies \eqref{relharm normal approx concl1} and \eqref{relharm normal approx concl2} (with possibly a worse constant $C$).
\end{proof}

\begin{remark}\label{relharm norm scaling rmk}{\rm 
Suppose that the hypotheses of Lemma~\ref{relharm normal approx lemma} hold true and $N_j$ satisfies the conclusion of Lemma~\ref{relharm normal approx lemma}.  Then for each $Z' \in \mathcal{M}$, $\rho' \in (\overline{R}(Z'),2]$, and $\gamma' \in (0,1)$ such that $B^{\mathcal M}_{\rho'}(Z') \subseteq B^{\mathcal M}_{\rho}(Z)$ and $B^{\mathcal M}_{\gamma'\rho'}(Z') \subseteq B^{\mathcal M}_{\gamma\rho}(Z)$, the conclusion of Lemma~\ref{relharm normal approx lemma} obviously holds true with $B^{\mathcal M}_{\gamma'\rho'}(Z')$, $\mathcal{K}_j \cap B^{\mathcal M}_{\gamma'\rho'}(Z')$, and $N_j |_{B^{\mathcal M}_{\gamma'\rho'}(Z')}$ in place of $B^{\mathcal M}_{\gamma\rho}(Z)$, $\mathcal{K}_j$, and $N_j$ respectively.
}\end{remark}

Given $N_j : B^{\mathcal M}_{\gamma\rho}(Z) \rightarrow \mathcal{A}_{q_j}({\rm T}\mathcal{M}^{\perp})$ satisfying \eqref{relharm normal approx concl1} and \eqref{relharm normal approx concl2}, we will write $N_j(X) = \sum_{l=1}^{q_j} \llbracket N_{j,l}(X) \rrbracket$ and $N_{j,a}(X) = \frac{1}{q_j} \sum_{l=1}^{q_j} N_{j,l}(X)$ for each $X \in B^{\mathcal M}_{\gamma\rho}(Z)$.  We write $DN_j(X) = \sum_{l=1}^{q_j} \llbracket DN_{j,l}(X) \rrbracket$ at $\mathcal{H}^n$-a.e.~$X \in B^{\mathcal M}_{\gamma\rho}(Z)$, following the conventions from~\cite[Subsection~3.5]{KrumWica}.

Suppose that the hypotheses of Lemma~\ref{relharm normal approx lemma} hold true.  Let $N_j : B^{\mathcal M}_{\gamma\rho}(Z) \rightarrow \mathcal{A}_{q_j}({\rm T}\mathcal{M}^{\perp})$ satisfy \eqref{relharm normal approx concl1} and \eqref{relharm normal approx concl2}.  Arguing as we did in Section~\ref{sec:normal approx est sec} to obtain \eqref{normal approx est3}, it follows using \eqref{relharm normal approx concl2}, \eqref{DeLSpa freq mono eqn11}, \eqref{DeLSpa freq mono eqn12} and \eqref{normal approx energy concl2} that
\begin{gather}
	\label{relharm norm est rmk eqn1} \mathcal{H}^n(B^{\mathcal M}_{\gamma\rho}(Z) \setminus \mathcal{K}_j) 
		+ \|T_j - \mathbf{G}_{N_j}\|(\mathbf{C}^{\mathcal M}_{\gamma\rho}(Z)) 
		\leq C \mathbf{m}_0^{3\gamma_2/4} \rho^{3\gamma_2/4-2} \int_{B^{\mathcal M}_{\rho}(Z)} |N|^2 , 
	\end{gather}
where $C =  C(n,m,q,M_0,N_0,C_e,C_h,\gamma) \in (0,\infty)$ are constants.  By \eqref{normal approx est1}, \eqref{normal approx est4}, \eqref{relharm normal approx concl1}, \eqref{relharm normal approx concl2}, and \eqref{relharm norm est rmk eqn1} 
\begin{align}\label{relharm norm est rmk eqn3}
	\frac{1}{q\rho} \int_{B^{\mathcal M}_{\gamma\rho}(Z)} \left| \sum_{j=1}^J q_j N_{j,a} \right| 
	\leq\,& \frac{1}{\rho} \int_{B^{\mathcal M}_{\gamma\rho}(Z)} |N_a| 
		+ C \mathbf{m}_0^{1/(2n)} {\rho^{\beta_2}} \mathcal{H}^n(B^{\mathcal M}_{\gamma\rho}(Z) \setminus \mathcal{K}_j)
	\\ \leq\,& C \rho^{\gamma_2/4-2} 
    \int_{B^{\mathcal M}_{\rho}(Z)} |N|^2 \nonumber 
\end{align}
where $C =  C(n,m,q,M_0,N_0,C_e,C_h,\gamma) \in (0,\infty)$ are constants.

\begin{lemma}[Energy estimate]\label{relharm energy lemma} 
Let $M_0,N_0,C_e,C_h,\varepsilon_2$ be as in Hypothesis~\ref{DeLSpa constants hypothesis} and Lemma~\ref{normal approx energy lemma}.  For every integer $p \in \{1,2,\ldots,q\}$ and every $\gamma \in (0,1)$ there exists $\varepsilon = \varepsilon(n,m,q,p,M_0,N_0,C_e,C_h,\gamma) \in (0,1)$ such that following holds true.  Assume that Hypothesis~\ref{relharm separation1 hyp} holds true and that $\rho > \overline{R}(Z)$, where $\overline{R}(Z)$ is as in Definition~\ref{crit rad defn}.  For each $j \in \{1,2,\ldots,p\}$ let $\mathcal{K}_j \subseteq B^{\mathcal M}_{(1+\gamma)\rho/2}(Z)$ and $N_j : B^{\mathcal M}_{(1+\gamma)\rho/2}(Z) \rightarrow \mathcal{A}_{q_j}({\rm T}\mathcal{M}^{\perp})$ satisfy \eqref{relharm normal approx concl1} and \eqref{relharm normal approx concl2} with $(1+\gamma)/2$ in place of $\gamma$.  Then 
\begin{align}\label{relharm energy concl}
	&\sum_{j=1}^p \int_{B^{\mathcal M}_{\gamma\rho}(Z)} \mathcal{G}(DN_j, q \llbracket D\varphi_j \rrbracket)^2 
	\leq \frac{C}{\rho^2} \sum_{j=1}^p \int_{B^{\mathcal M}_{(1+\gamma)\rho/2}(Z)} \mathcal{G}(N_j, q \llbracket \varphi_j \rrbracket)^2 
		\\&\hspace{15mm} + C \rho^2 \sum_{j=1}^p \|\Delta_{\mathcal M} \varphi_j\|_{L^2(B^{\mathcal M}_{(1+\gamma)\rho/2}(Z))}^2 
		+ C \mathbf{m}_0^{\gamma_2/4} \rho^{\gamma_2/4-2} \int_{B^{\mathcal M}_{\rho}(Z)} |N|^2 , \nonumber 
\end{align}
where $C = C(n,m,q,M_0,N_0,C_e,C_h,\gamma) \in (0,\infty)$ is a constant.
\end{lemma}

\begin{proof} 
Let $\eta \in C^1_c(\mathcal{M})$ such that $0 \leq \eta \leq 1$ on $\mathcal{M}$, $\eta = 1$ on $B^{\mathcal M}_{\gamma\rho}(Z)$, $\eta = 0$ on $\mathcal{M} \setminus B^{\mathcal M}_{(1+\gamma)\rho/2}(Z)$, and $|D\eta| \leq \frac{6}{1-\gamma}$.  Let $\zeta_j \in C^1_c(\mathbf{C}^{\mathcal M}_{\rho}(Z),\mathbb{R}^{n+m})$ such that $\zeta_j(X) = \eta(\mathbf{p}(X))^2 \,(X - \mathbf{p}(X) - \varphi_j(\mathbf{p}(X)))$ for each $X \in \mathbf{C}^{\mathcal M}_{\rho}(Z)$ near $\op{spt} T_j \cup \op{spt} \mathbf{G}_{N_j}$.  For $\delta > 0$ sufficiently small and for each $t \in (-\delta,\delta)$, define $\Phi_{j,t} : \mathbf{C}^{\mathcal M}_{\rho}(Z) \rightarrow \mathbf{C}^{\mathcal M}_{\rho}(Z)$ by $\Phi_{j,t}(X) = X + t \,\zeta_j(X)$ for all $X \in \mathbf{U}$.  Notice that $\Phi_{j,t\#} \mathbf{G}_{N_j}$ is the graph of the $q_j$-valued function $\sum_{l=1}^{q_j} \llbracket N_{j,l}(X) + t \,\eta(X)^2 \,(N_{j,l}(X) - \varphi_j(X)) \rrbracket$ for $X \in \mathcal{M}$.  Recalling \eqref{normal graph defn eqn5}, 
\begin{align}\label{relharm energy eqn1}
	&\delta\mathbf{G}_{N_j}(\zeta_j) 
	= \frac{d}{dt} \mathbf{M}(\Phi_{j,t\#} \mathbf{G}_{N_j}) \,\Bigg|_{t=0} 
	\\=\,& \int_{\mathcal M} \sum_{l=1}^{q_j} \sum_{i,k=1}^n \sqrt{G(X,N_{j,l},DN_{j,l})} \,G^{ik}(X,N_{j,l},DN_{j,l}) 
		\,\big\langle \tau_i + D_{\tau_i} N_{j,l}, D_{\tau_k} (\eta^2 (N_{j,l} - \varphi_j)) \big\rangle , \nonumber 
\end{align}
where $\delta\mathbf{G}_{N_j}$ is the first variation of the integral varifold associated with $\mathbf{G}_{N_j}$ (as in~\cite[39.2]{SimonGMT}), $\tau_1,\tau_2,\ldots,\tau_2$ is a positively oriented orthonormal tangent frame field on $\mathcal{M}$, and $G$ and $G^{ik}$ are as in Subsection~\ref{sec:normal graph sec}.  Recall \eqref{normal graph defn eqn6} and notice that $\langle \tau_i, D_{\tau_k} (\eta^2 N_{j,l}) \rangle = -\eta^2 \,\langle A_{ik}, N_{j,l} \rangle$, where $A_{ik} = (D_{\tau_i} \tau_k)^{\perp}$ is the second fundamental form of $\mathcal{M}$.  Using this together with Taylor's theorem 
\begin{align}\label{relharm energy eqn2}
	\delta\mathbf{G}_{N_j}(\zeta_j) \geq\,& -q_j \int_{\mathcal M} \big( \eta^2 \langle H_{\mathcal M}, N_{j,a} \rangle 
		+ \op{div}_{\mathcal M}(\eta^2 \varphi_j) \big) 
		\\&+ \int_{\mathcal M} \sum_{l=1}^{q_j} \big( \eta^2 \langle DN_{j,l}, DN_{j,l} - D\varphi_j \rangle 
		+ 2 \eta \,\langle DN_{j,l}, (N_{j,l} - \varphi_j) \otimes D\eta \rangle \big) \nonumber 
		\\&- C \int_{\mathcal M} \big( |A_{\mathcal M}| \,|N_j| + |DN_j|^2 \big) \big( \eta^2 \,\mathcal{G}(DN_j,q_j\llbracket D\varphi_j \rrbracket) 
		+ \eta \,|D\eta| \,\mathcal{G}(N_j,q\llbracket \varphi_j \rrbracket) \big) \nonumber 
\end{align}
for some constant $C = C(n,m,q) \in (0,\infty)$.  By the divergence theorem $\int_{\mathcal M} \op{div}_{\mathcal M} \xi \,d\mathcal{H}^n = 0$ with $\xi = \sum_{i=1}^n \eta^2 \,\langle N_{j,a} - \varphi_j, D_{\tau_i} \varphi_j \rangle \,\tau_i$ 
\begin{equation}\label{relharm energy eqn3}
	\int_{\mathcal M} \big( \eta^2 \langle D\varphi_j, DN_{j,a} - D\varphi_j \rangle + 2 \eta \,\langle D\varphi_j, (N_{j,a} - \varphi_j) \otimes D\eta \rangle \big) 
	= -\int_{\mathcal M} \eta^2 \langle N_{j,a} - \varphi_j, \Delta_{\mathcal M} \varphi_j \rangle .  
\end{equation}
Subtracting \eqref{relharm energy eqn2} and \eqref{relharm energy eqn3} and using Cauchy's inequality 
\begin{align}\label{relharm energy eqn4}
	\frac{1}{2} \int_{\mathcal M} \eta^2 \mathcal{G}(DN_j, q_j \llbracket D\varphi_j \rrbracket)^2 
	\leq\,& q_j \int_{\mathcal M} \big( \eta^2 \langle H_{\mathcal M}, N_{j,a} \rangle + \op{div}_{\mathcal M}(\eta^2 \varphi_j) \big) 
	\\&+ \frac{C}{\rho^2} \int_{B^{\mathcal M}_{(1+\gamma)\rho/2}(Z)} \mathcal{G}(N_j, q \llbracket \varphi_j \rrbracket)^2 
	+ C \rho^2 \int_{\mathcal M} \eta^2 |\Delta_{\mathcal M} \varphi_j|^2 \nonumber
	\\&+ C \int_{\mathcal M} \eta^2 (|A_{\mathcal M}|^2 |N_j|^2 + |DN_j|^4) + \delta\mathbf{G}_{N_j}(\zeta_j) , \nonumber 
\end{align}
where $C = C(n,m,q,\gamma) \in (0,\infty)$ is a constant.  Recall that $\varphi_j$ satisfies \eqref{relharm separation1 small norm} and by~\cite[Corollary~2.2]{DeLSpa2}, $\sup_{X \in \op{spt} T \cap B^{\mathcal M}_{\rho}(Z)} \op{dist}(X,\mathcal{M}) \leq C \mathbf{m}_0^{1/(2n)} \rho^{1+\beta_2}$.  Hence $|D\zeta_j| \leq C$ near $\op{spt} T \cup \op{spt}\mathbf{G}_N$.  Thus since $T_j$ is area-stationary in $\mathbf{C}^{\mathcal M}_{\rho}(Z)$ and using \eqref{relharm norm est rmk eqn1}
\begin{align}\label{relharm energy eqn5}
	|\delta\mathbf{G}_{N_j}(\zeta_j)| =\,& |\delta\mathbf{G}_{N_j}(\zeta_j) - \delta T_j(\zeta_j)| 
	\leq \|D\zeta_j\|_{C^0} \,\|\mathbf{G}_{N_j} - T_j\|(\mathbf{C}^{\mathcal M}_{(1+\gamma)\rho/2}(Z)) 
	\\ \leq\,& C \mathbf{m}_0^{3\gamma_2/4} \rho^{3\gamma_2/4-2} \int_{B^{\mathcal M}_{\rho}(Z)} |N|^2 , \nonumber 
\end{align}
where $\delta T_j$ denotes the first variation of the integral varifold associated with $T_j$.  Now summing \eqref{relharm energy eqn4} over $j \in \{1,2,\ldots,J\}$ and bounding the terms on the right-hand side using $\|A_{\mathcal M}\|_{C^0} \leq C \mathbf{m}_0^{1/2}$, \eqref{relharm normal approx concl2}, \eqref{relharm norm est rmk eqn3}, \eqref{relharm separation1 avgfree}, \eqref{relharm norm est rmk eqn1}, Lemma~\ref{normal approx energy lemma}, and \eqref{relharm energy eqn5} gives us \eqref{relharm energy concl}.
\end{proof}

\subsection{Harmonic approximation of an area-minimizing current relative to $C^2$-graphs}\label{sec:relharm harm approx}  Let $\gamma,\varepsilon,\varsigma \in (0,1)$.  Assume that Hypothesis~\ref{relharm separation1 hyp} holds true.  For each $j \in \{1,2,\ldots,p\}$, let $N_j : B^{\mathcal M}_{\gamma\rho/2}(Z) \rightarrow \mathcal{A}_{q_j}({\rm T}\mathcal{M}^{\perp})$ satisfy \eqref{relharm normal approx concl1} and \eqref{relharm normal approx concl2}.  We define 
\begin{align}\label{relharm excess defn1} 
	\mathcal{E}_{Z,\gamma\rho,\rho}^2 =\,& \frac{1}{\rho^n} \sum_{j=1}^p \int_{B^{\mathcal M}_{\gamma\rho}(Z)} 
		\mathcal{G}(DN_j, q_j \llbracket D\varphi_j \rrbracket)^2 
		+ \frac{\rho^2}{\varepsilon^2} \sum_{j=1}^p \|\Delta_{\mathcal M} \varphi_j\|_{C^0(B^{\mathcal M}_{\rho}(Z))}^2  
		\\&+ \sum_{j=1}^p \|D\varphi_j\|_{C^0(B^{\mathcal M}_{\rho}(Z))} \cdot \frac{1}{\rho^n} \int_{B^{\mathcal M}_{\rho}(Z)} |DN|^2  
		+ \frac{1}{\rho^{n+2-\varsigma}} \int_{B^{\mathcal M}_{\rho}(Z)} |N|^2 . \nonumber 
\end{align}
The first term on the right-hand side of \eqref{relharm excess defn1} is a tilt excess of $T_j$ with respect to the graph of $\varphi_j$.  The second term $\varepsilon^{-1} \rho \|\Delta_{\mathcal M} \varphi_j\|_{C^0(B^{\mathcal M}_{\rho}(Z))}$ is zero if $\varphi_j$ is harmonic on $B^{\mathcal M}_{\rho}(Z)$ and otherwise guarantees that $\varphi_j$ is nearly-harmonic in the sense that $\rho \|\Delta_{\mathcal M} \varphi_j\|_{C^0(B^{\mathcal M}_{\rho}(Z))} \leq \varepsilon \,\mathcal{E}_{Z,\gamma\rho,\rho}$.  The last term on the right-hand side of \eqref{relharm excess defn1} controls error terms arising from the area formula and the estimates for the normal approximation maps $N$ and $N_j$.

\begin{definition}\label{parallel harmonic defn} 
Let $\mathcal{M}$ be as in Hypothesis~\ref{M and N hypothesis} and let $\Omega \subseteq \mathcal{M}$.  Let $\varphi : \Omega \rightarrow \mathbb{R}^{n+m}$ be a single-valued function and let $w : \Omega \rightarrow \mathcal{A}_q(\mathbb{R}^{n+m})$ be a $q$-valued function.  We define $(\varphi+w) : \Omega \rightarrow \mathcal{A}_q(\mathbb{R}^{n+m})$
\begin{equation*}
	(\varphi+w)(X) = \sum_{l=1}^q \llbracket \varphi(X) + w_l(X) \rrbracket
\end{equation*}
for each $X \in \Omega$, where $w(X) = \sum_{l=1}^q \llbracket w_l(X) \rrbracket$.  In particular, given $a = \sum_{l=1}^q \llbracket a_l \rrbracket \in \mathcal{A}_q(\mathbb{R}^{n+m})$, we define $(\varphi+a) : \Omega \rightarrow \mathcal{A}_q(\mathbb{R}^{n+m})$
\begin{equation*}
	(\varphi+a)(X) = \sum_{l=1}^q \llbracket \varphi(X) + a_l \rrbracket
\end{equation*}
for each $X \in \Omega$.
\end{definition}

\begin{definition}\label{a sep width defn} 
Let $a = \sum_{i=1}^N q_i \llbracket a_i \rrbracket \in \mathcal{A}_q(\mathbb{R}^{n+m})$, where $N,q_1,\ldots,q_N$ are positive integers and $a_1,\ldots,a_N \in \mathbb{R}^{n+m}$ are distinct.  If $N \geq 2$, we define the size $\op{size} a$, separation $\op{sep} a$, and width $\op{width} a$ of $a$ by 
\begin{equation*}
	\op{size} a = N , \quad \op{sep} a = \min_{i \neq j} |a_i - a_j|, \quad \op{width} a = \max_{i \neq j} |a_i - a_j| .
\end{equation*}
If instead $N = 1$, we define $\op{size} a = 1$, $\op{sep} a = \infty$, and $\op{width} a = 0$.
\end{definition}

\begin{lemma}[Harmonic approximation]\label{relharm harmonic approx lemma}
Let $M_0,N_0,C_e,C_h,\varepsilon_2$ be as in Hypothesis~\ref{DeLSpa constants hypothesis} and Lemma~\ref{normal approx energy lemma}.  For every integer $p \in \{1,2,\ldots,q\}$ and every $\eta,\gamma \in (0,1)$ there exists a constant 
$\varepsilon = \varepsilon(n,m,q,p,M_0,N_0,C_e,C_h,\eta,\gamma) \in (0,1)$ such that the following holds true.  Assume that Hypothesis~\ref{relharm separation1 hyp} holds true and that $\rho > \overline{R}(Z)$, where $\overline{R}(Z)$ is as in Definition~\ref{crit rad defn}.  For each $j \in \{1,2,\ldots,p\}$ let $\mathcal{K}_j \subseteq B^{\mathcal M}_{(1+\gamma)\rho/2}(Z)$ and $N_j : B^{\mathcal M}_{(1+\gamma)\rho/2}(Z) \rightarrow \mathcal{A}_{q_j}({\rm T}\mathcal{M}^{\perp})$  satisfy \eqref{relharm normal approx concl1} and \eqref{relharm normal approx concl2} with $(1+\gamma)/2$ in place of $\gamma$.  Then there exist Dirichlet energy minimizing $q$-valued functions $w_j : B_{\gamma\rho}(0,{\rm T}_Z \mathcal{M}) \rightarrow \mathcal{A}_{q_j}(\mathbb{R}^{n+m})$ such that 
\begin{align}\label{relharm harmonic approx concl}
	&\frac{1}{\rho^{n+2}} \sum_{j=1}^p \int_{B^{\mathcal M}_{\gamma\rho}(Z)} \mathcal{G}(N_j,\varphi_j+w_j \circ \exp_Z^{-1})^2 
	\\&+ \frac{1}{\rho^n} \sum_{j=1}^p \int_{B^{\mathcal M}_{\gamma\rho}(Z)} (|D(N_j - \varphi_j)| - |D(w_j \circ \exp_Z^{-1})|)^2 
	\leq \eta^2 \mathcal{E}_{Z,(1+\gamma)\rho/2,\rho}^2 , \nonumber 
\end{align}
where $\mathcal{E}_{Z,(1+\gamma)\rho/2,\rho}$ is as in \eqref{relharm excess defn1} with $\varsigma = 2\beta_2$ (and with $(1+\gamma)/2$ in place of $\gamma$).
\end{lemma}

\begin{proof} 
First observe that if $\mathcal{E}_{Z,(1+\gamma)\rho/2,\rho} = 0$, then $N = q \llbracket 0 \rrbracket$ on $B^{\mathcal M}_{\rho}(Z)$ and $DN_j = q_j \llbracket D\varphi_j \rrbracket$ on $B^{\mathcal M}_{(1+\gamma)\rho/2}(Z)$.  By \eqref{normal approx est3} and \eqref{relharm norm est rmk eqn1}, $N = q \llbracket 0 \rrbracket$ on $B^{\mathcal M}_{\rho}(Z)$ implies that $\mathcal{K} \cap B^{\mathcal M}_{(1+\gamma)\rho/2}(Z) = \mathcal{K}_j = B^{\mathcal M}_{(1+\gamma)\rho/2}(Z)$ and thus $N_j = q_j \llbracket 0 \rrbracket$ on $B^{\mathcal M}_{(1+\gamma)\rho/2}(Z)$.  Hence $D\varphi_j = 0$ on $B^{\mathcal M}_{(1+\gamma)\rho/2}(Z)$; that is, there exists a constant $a_j \in {\rm T}_Z\mathcal{M}$ such that $\varphi_j(X) = a_j$ for all $X \in B^{\mathcal M}_{(1+\gamma)\rho/2}(Z)$.  Set $w_j = q_j \llbracket -a_j \rrbracket$ on $B_{\gamma\rho}(0,{\rm T}_Z \mathcal{M})$.  Since $w_j$ is a constant $q_j$-valued function, $w_j$ is Dirichlet energy minimizing on $B_{\gamma\rho}(0,{\rm T}_Z \mathcal{M})$.  Since $N_j - \varphi_j = q_j \llbracket -a_j \rrbracket = w_j \circ \exp_Z^{-1}$ on $B^{\mathcal M}_{\gamma\rho}(Z)$, \eqref{relharm harmonic approx concl} holds true.  

Now let us suppose that $\mathcal{E}_{Z,(1+\gamma)\rho/2,\rho} > 0$.  Fix positive integers $p,q_1,q_2,\ldots,q_p$ such that $\sum_{j=1}^p q_j = q$.  For $\nu = 1,2,3,\ldots$, suppose that $\varepsilon_{\nu} \rightarrow 0^+$, $\mathbf{m}_{0,\nu}$, $\mathcal{M}^{(\nu)}$, $Z_{\nu}$, $\rho_{\nu}$, $T^{(\nu)}$, $T^{(\nu)}_j$, $\mathcal{K}^{(\nu)}$, $N^{(\nu)}$, $\mathcal{K}^{(\nu)}_j$, $N^{(\nu)}_j$, and $\varphi^{(\nu)}_j$ satisfy the hypotheses of Lemma~\ref{relharm harmonic approx lemma} in place of $\varepsilon$, $\mathbf{m}_0$, $\mathcal{M}$, $Z$, $\rho$, $T$, $T_j$, $\mathcal{K}$, $N$, $\mathcal{K}_j$, $N_j$, and $\varphi_j$.  Let 
\begin{align*}
	\mathcal{E}_{\nu}^2 =\,& \frac{1}{\rho_{\nu}^n} \sum_{j=1}^p \int_{B^{\mathcal{M}^{(\nu)}}_{(1+\gamma)\rho_{\nu}/2}(Z_{\nu})} 
		\mathcal{G}(DN^{(\nu)}_j, q_j \llbracket D\varphi^{(\nu)}_j \rrbracket)^2 
		+ \frac{\rho_{\nu}^2}{\varepsilon_{\nu}^2} \sum_{j=1}^p \|\Delta_{\mathcal{M}^{(\nu)}} \varphi^{(\nu)}_j 
		\|_{C^0(B^{\mathcal{M}^{(\nu)}}_{\rho_{\nu}}(Z_{\nu}))}^2  
		\\&+ \sum_{j=1}^p \|D\varphi^{(\nu)}_j\|_{C^0(B^{\mathcal{M}^{(\nu)}}_{\rho_{\nu}}(Z_{\nu}))} 
		\cdot \frac{1}{\rho_{\nu}^n} \int_{B^{\mathcal{M}^{(\nu)}}_{\rho_{\nu}}(Z_{\nu})} |DN^{(\nu)}|^2 
		+ \frac{1}{\rho_{\nu}^{n+2-\varsigma}} \int_{B^{\mathcal{M}^{(\nu)}}_{\rho_{\nu}}(Z_{\nu})} |N^{(\nu)}|^2 \nonumber
\end{align*}
(as in \eqref{relharm excess defn1}).  Let $\mathbf{e}_{\nu} : B_2(0,{\rm T}_{Z_{\nu}}\mathcal{M}^{(\nu)}) \rightarrow \mathcal{M}^{(\nu)}$ be the exponential map of $\mathcal{M}^{(\nu)}$ at $Z_{\nu}$.  It suffices to show that for infinitely many $\nu$ there exists a Dirichlet energy minimizing $q$-valued function $w^{(\nu)}_j : B_{\gamma\rho_{\nu}}(0,{\rm T}_{Z_{\nu}}\mathcal{M}^{(\nu)}) \rightarrow \mathcal{A}_{q_j}(\mathbb{R}^{n+m})$ such that 
\begin{align}\label{relharm harmonic approx eqn2}
	&\frac{1}{\rho_{\nu}^{n+2}} \sum_{j=1}^p \int_{B^{\mathcal{M}^{(\nu)}}_{\gamma\rho_{\nu}}(Z_{\nu})} 
		\mathcal{G}(N^{(\nu)}_j,\varphi^{(\nu)}_j+w^{(\nu)}_j \circ \mathbf{e}_{\nu}^{-1})^2 
	\\&+ \frac{1}{\rho_{\nu}^n} \sum_{j=1}^p \int_{B^{\mathcal{M}^{(\nu)}}_{\gamma\rho_{\nu}}(Z_{\nu})} 
		(|D(N^{(\nu)}_j - \varphi^{(\nu)}_j)| - |D(w^{(\nu)}_j \circ \mathbf{e}_{\nu}^{-1})|)^2 
	\leq \eta^2 \mathcal{E}_{\nu}^2 \nonumber 
\end{align}
(as in \eqref{relharm harmonic approx concl}).

\noindent\textit{Step~1.  Construction of a fine blow-up.}  Let $\mathfrak{q}_{\nu}$ be a rotation of $\mathbb{R}^{n+m}$ such that $\mathfrak{q}_{\nu}({\rm T}_{Z_{\nu}}\mathcal{M}^{(\nu)}) = \mathbb{R}^n \times \{0\}$.  Rescale and rotate by setting 
\begin{gather*}
	\overline{T}^{(\nu)} = \mathfrak{q}_{\nu\#} \eta_{Z_{\nu},\rho_{\nu}\#} T^{(\nu)}, \quad 
	\overline{T}^{(\nu)}_j = \mathfrak{q}_{\nu\#} \eta_{Z_{\nu},\rho_{\nu}\#} T^{(\nu)}_j, \quad
	\overline{\mathcal M}^{(\nu)} = \mathfrak{q}_{\nu}(\eta_{Z_{\nu},\rho_{\nu}}(\mathcal{M}^{(\nu)})) , \\ 
	\overline{N}^{(\nu)}_j(X) = \sum_{l=1}^{q_j} \llbracket \rho_{\nu}^{-1} \mathfrak{q}_{\nu} N_{j,l}^{(\nu)}(Z_{\nu} + \rho_{\nu} \mathfrak{q}_{\nu}^{-1} X) 
		\rrbracket \text{ for each } X \in \overline{\mathcal M}^{(\nu)} , \\
	\overline{\varphi}^{(\nu)}(X) = \rho_{\nu}^{-1} \mathfrak{q}_{\nu} \varphi^{(\nu)}(Z_{\nu} + \rho_{\nu} \mathfrak{q}_{\nu}^{-1} X) 
		\text{ for each } X \in \overline{\mathcal M}^{(\nu)} ,
\end{gather*}
where $N_j^{(\nu)}(X) = \sum_{l=1}^{q_j} \llbracket N^{(\nu)}_{j,l}(X) \rrbracket$ for all $X \in \mathcal{M}^{(\nu)}$.  Note that each current $\overline{T}^{(\nu)}$ and $\overline{T}^{(\nu)}_j$ is locally area-minimizing in $\mathbf{C}^{\overline{\mathcal M}^{(\nu)}}_1(0)$.  Let $\overline{\mathbf e}_{\nu} : B_1(0) \rightarrow \overline{\mathcal M}^{(\nu)}$ be the exponential map of $\overline{\mathcal M}^{(\nu)}$ at the origin.  For each $j \in \{1,2,\ldots,p\}$, define $g_j^{(\nu)} : B_{(1+\gamma)/2}(0) \rightarrow \mathcal{A}_q(\mathbb{R}^{n+m})$ by  
\begin{equation*}
	g_j^{(\nu)}(x) = \sum_{l=1}^q \left\llbracket \frac{\overline{N}^{(\nu)}_{j,l}(\overline{\mathbf e}_{\nu}(x)) 
		- \overline{\varphi}^{(\nu)}_j(\overline{\mathbf e}_{\nu}(x))}{\mathcal{E}_{\nu}} \right\rrbracket
\end{equation*}
for each $x \in B_{(1+\gamma)/2}(0)$.  Notice that 
\begin{equation}\label{relharm harmonic approx eqn3}
	\sum_{j=1}^p \int_{B_{(1+\gamma)/2}(0)} |Dg_j^{(\nu)}|^2 \leq \frac{4}{\mathcal{E}_{\nu}^2 \rho_{\nu}^n} \sum_{j=1}^p 
		\int_{B^{\mathcal{M}^{(\nu)}}_{(1+\gamma)\rho_{\nu}/2}(Z_{\nu})} \mathcal{G}(DN^{(\nu)}_j, q_j \llbracket D\varphi^{(\nu)}_j \rrbracket)^2 \leq 4 
\end{equation}
for all large $\nu$.  Hence by~\cite[Proposition~4.3]{DeLSpa1} after passing to a subsequence for each $j \in \{1,2,\ldots,p\}$ there exists 
\begin{enumerate}[itemsep=2mm,topsep=0mm]
	\item[(i)]  positive integers $p_j,q_{j,1},\ldots,q_{j,p_j}$ with $\sum_{k=1}^{p_j} q_{j,k} = q_j$; 
	\item[(ii)]  $y^{(\nu)}_{j,k} \in \mathbb{R}^{n+m}$ with $\lim_{\nu\rightarrow\infty} |y^{(\nu)}_{j,k} - y^{(\nu)}_{j,k'}| = \infty$ for all $k \neq k'$; 
	\item[(iii)]  $q_{j,k}$-valued functions $w_{j,k} \in W^{1,2}(B_{(1+\gamma)/2}(0),\mathcal{A}_{q_{j,k}}(\mathbb{R}^{n+m}))$; 
	\item[(iv)]  a sequence of translating sheets $(h_j^{(\nu)})$ in $W^{1,2}(B_{(1+\gamma)/2}(0),\mathcal{A}_q(\mathbb{R}^{n+m}))$ given by $h_j^{(\nu)}(x) = \sum_{k=1}^{p_j} \sum_{l=1}^{q_{j,k}}$ $\llbracket y^{(\nu)}_{j,k} + w_{j,k,l}(x) \rrbracket$ for each $x \in B_{(1+\gamma)/2}(0)$, where $w_{j,k}(x) = \sum_{l=1}^{q_{j,k}} \llbracket w_{j,k,l}(x) \rrbracket$;
\end{enumerate}
such that 
\begin{gather}
	\label{relharm harmonic approx eqn4a} \lim_{\nu\rightarrow\infty} \int_{B_{(1+\gamma)/2}(0)} \mathcal{G}(g_j^{(\nu)}, h_j^{(\nu)})^2 = 0 , \\
	\label{relharm harmonic approx eqn4b} \sum_{k=1}^{p_j} \int_{B_{\gamma}(0)} |Dw_{j,k}|^2 
		\leq \liminf_{\nu\rightarrow\infty} \int_{B_{\gamma}(0)} |Dg_j^{(\nu)}|^2 , \\
	\label{relharm harmonic approx eqn4c} \limsup_{\nu\rightarrow\infty} \int_{B_{\gamma}(0)} (|Dg_j^{(\nu)}| - |Dh_j^{(\nu)}|)^2 
		\leq \limsup_{\nu\rightarrow\infty} \int_{B_{\gamma}(0)} |Dg_j^{(\nu)}|^2 - \sum_{k=1}^{p_j} \int_{B_{\gamma}(0)} |Dw_{j,k}|^2 .
\end{gather}
Note that $|Dh_j^{(\nu)}|^2 = \sum_{k=1}^{p_j} |Dw_{j,k}|^2$ on $B_{(1+\gamma)/2}(0)$ for all $\nu$.  

In light \eqref{relharm harmonic approx eqn4a} and \eqref{relharm harmonic approx eqn4c}, to complete the proof it suffices to show that for each $j \in \{1,2,\ldots,p\}$, $h^{(\nu)}_j$ is locally Dirichlet energy minimizing and 
\begin{equation}\label{relharm harmonic approx eqn20} 
	\lim_{\nu\rightarrow\infty} \int_{B_{\gamma}(0)} |Dg_j^{(\nu)}|^2 = \sum_{k=1}^{p_j} \int_{B_{\gamma}(0)} |Dw_{j,k}|^2 .
\end{equation}
Then it will follow that \eqref{relharm harmonic approx eqn2} holds true with $w^{(\nu)}_j(X) = \mathcal{E}_{\nu} \rho_{\nu} \mathfrak{q}_{\nu}^{-1} h^{(\nu)}_j(\rho_{\nu}^{-1} \mathfrak{q}_{\nu} X)$ for each $X \in B_{\gamma\rho_{\nu}}(0,{\rm T}_{Z_{\nu}}\mathcal{M}^{(\nu)})$.

\noindent\textit{Step~2.  Average-free condition.}  We claim that we can assume that 
\begin{gather}\label{relharm harmonic approx eqn8}
	\frac{1}{q} \sum_{j=1}^p \sum_{k=1}^{p_j} q_{j,k} w_{j,k,a} = 0 \text{ on } B_{(1+\gamma)/2}(0), \quad 
	\frac{1}{q} \sum_{j=1}^p \sum_{k=1}^{p_j} q_{j,k} y^{(\nu)}_{j,k} = 0 , \\ 
	\frac{1}{q} \sum_{j=1}^p q_j h^{(\nu)}_{j,a} = 0 \text{ on } B_{(1+\gamma)/2}(0) , \nonumber 
\end{gather}
where $w_{j,k,a}(x) = \frac{1}{q_{j,k}} \sum_{l=1}^{q_{j,k}} w_{j,k,l}(x)$ is the single-valued average of $w_{j,k}$ and $h^{(\nu)}_{j,a}(x) = \frac{1}{q_j} \sum_{l=1}^{q_j} h^{(\nu)}_{j,l}(x)$ is the single-valued average of $h^{(\nu)}_j(x) = \sum_{l=1}^{q_j} \llbracket h^{(\nu)}_{j,l}(x) \rrbracket$.  
By $\varsigma \leq \gamma_2/4$, \eqref{relharm norm est rmk eqn3}, and \eqref{relharm separation1 avgfree} 
\begin{equation}\label{relharm harmonic approx eqn6}
	\frac{1}{q} \int_{B_{(1+\gamma)/2}(0)} \left| \sum_{j=1}^p q_j g^{(\nu)}_{j,a} \right| 
	\leq \frac{C}{q \,\mathcal{E}_{\nu}} \int_{B^{\mathcal{M}^{(\nu)}}_{(1+\gamma)\rho_{\nu}/2}(Z_{\nu})} 
		\left| \sum_{j=1}^p   q_j (N^{(\nu)}_{j,a} - \varphi^{(\nu)}_j) \right| J\mathbf{e}_{\nu}^{-1} 
	\leq C \mathcal{E}_{\nu} , 
\end{equation}
where $J\mathbf{e}_{\nu}^{-1}$ is the Jacobian of the inverse function of the exponential map $\mathbf{e}_{\nu}$ of $\mathcal{M}^{(\nu)}$ at $Z_{\nu}$, which is bounded, and $C = C(n,m,q,M_0,N_0,C_e,C_h) \in (0,\infty)$ are constants.  Hence by \eqref{relharm harmonic approx eqn4a} 
\begin{equation}\label{relharm harmonic approx eqn7}
	\lim_{\nu\rightarrow\infty} \frac{1}{q} \int_{B_{(1+\gamma)/2}(0)} \left| \sum_{j=1}^p q_j h^{(\nu)}_{j,a} \right|
	= \lim_{\nu\rightarrow\infty} \frac{1}{q} \int_{B_{(1+\gamma)/2}(0)} \left| \sum_{j=1}^p \sum_{k=1}^{p_j} q_j (y^{(\nu)}_{j,k} + w_{j,k,a}) \right| = 0 .
\end{equation}
By \eqref{relharm harmonic approx eqn7}, $\frac{1}{q} \sum_{j=1}^p q_j w_{j,k,a}$ is a constant function on $B_{(1+\gamma)/2}(0)$.  Thus after adding to $y^{(\nu)}_{j,k}$ and $w_{j,k,a}$ appropriate constant vectors in $\mathbb{R}^{n+m}$, we may assume that \eqref{relharm harmonic approx eqn8} holds true while also guaranteeing that \eqref{relharm harmonic approx eqn4a} holds true and without changing \eqref{relharm harmonic approx eqn4b} and \eqref{relharm harmonic approx eqn4c}.

\noindent\textit{Step~3.  Show that $h_j^{(\nu)}$ and $w_{j,k}$ are locally Dirichlet energy minimizing.}  To show that $h_j^{(\nu)}$ is locally Dirichlet energy minimizing in $B_{\gamma}(0)$ for all sufficiently large $\nu$, it suffices to show that $w_{j,k}$ is locally Dirichlet energy minimizing in $B_{\gamma}(0)$ for all $k \in \{1,2,\ldots,p_j\}$ (see the proof of~\cite[Lemma~4.7]{KrumWicc}).  Suppose to the contrary that at least one $w_{j,k}$ is not locally Dirichlet energy minimizing in $B_{\gamma}(0)$ and thus there exist $\delta > 0$ and $\widetilde{w}_{j,k} \in W^{1,2}(B_{\gamma}(0),\mathcal{A}_{q_{j,k}}(\mathbb{R}^{n+m}))$ such that 
\begin{gather}
	\label{relharm harmonic approx eqn9a} \widetilde{w}_{j,k} = w_{j,k} \text{ on } B_{(1+\gamma)/2}(0) \setminus B_{\gamma}(0), \\ 
	\label{relharm harmonic approx eqn9b} \sum_{j=1}^p \sum_{k=1}^{p_j} \int_{B_{\gamma}(0)} |D\widetilde{w}_{j,k}|^2 
		\leq \sum_{j=1}^p \sum_{k=1}^{p_j} \int_{B_{\gamma}(0)} |Dw_{j,k}|^2 - \delta .  
\end{gather}
By \eqref{relharm harmonic approx eqn8} and \eqref{relharm harmonic approx eqn9a}, $\frac{1}{q} \sum_{j=1}^p \sum_{k=1}^{p_j} q_{j,k} \widetilde{w}_{j,k,a} = 0$ on $B_{(1+\gamma)/2}(0) \setminus B_{\gamma}(0)$.  Thus, by possibly reducing $\sum_{j=1}^p \sum_{k=1}^{p_j} \int_{B_{\gamma}(0)} |D\widetilde{w}_{j,k}|^2$, we may assume that 
\begin{equation}\label{relharm harmonic approx eqn9c} 
	\frac{1}{q} \sum_{j=1}^p \sum_{k=1}^{p_j} q_{j,k} \widetilde{w}_{j,k,a} = 0 \text{ on } B_{(1+\gamma)/2}(0) . 
\end{equation}

Set $\widetilde{h}^{(\nu)}_j(x) = \sum_{k=1}^{p_j} \sum_{l=1}^{q_{j,k}} \llbracket y^{(\nu)}_{j,k} + \widetilde{w}_{j,k,l}(x) \rrbracket$ for each $x \in B_{(1+\gamma)/2}(0)$, where $\widetilde{w}_{j,k}(x) = \sum_{l=1}^{q_{j,k}} \llbracket \widetilde{w}_{j,k,l}(x) \rrbracket$.  Note that by \eqref{relharm harmonic approx eqn9a}, \eqref{relharm harmonic approx eqn9b}, \eqref{relharm harmonic approx eqn9c}, and \eqref{relharm harmonic approx eqn8} 
\begin{gather}
	\label{relharm harmonic approx eqn10a} \widetilde{h}^{(\nu)}_j = h^{(\nu)}_j \text{ on } B_{(1+\gamma)/2}(0) \setminus B_{\gamma}(0), \\ 
	\label{relharm harmonic approx eqn10b} \frac{1}{q} \sum_{j=1}^p q_j \,\widetilde{h}^{(\nu)}_{j,a} = 0 \text{ on } B_{(1+\gamma)/2}(0) , \\
	\label{relharm harmonic approx eqn10c} \sum_{j=1}^p \int_{B_{\gamma}(0)} |D\widetilde{h}^{(\nu)}_j|^2 
		\leq \sum_{j=1}^p \int_{B_{\gamma}(0)} |Dh^{(\nu)}_j|^2 - \delta , 
\end{gather}
where $\widetilde{h}^{(\nu)}_{j,a}(x) = \frac{1}{q_j} \sum_{k=1}^{p_j} \sum_{l=1}^{q_{j,k}} (y^{(\nu)}_{j,k} + \widetilde{w}_{j,k,l}(x))$ is the single-valued average for $\widetilde{h}^{(\nu)}_{j,a}$.  Using~\cite[Proposition~4.4]{DeLSpa1} to interpolate between $g^{(\nu)}_j$ and $\widetilde{h}_j^{(\nu)}$, we can construct a Lipschitz $q$-valued function $\widehat{h}_j^{(\nu)} : B_{(1+\gamma)/2}(0) \rightarrow \mathcal{A}_q(\mathbb{R}^{n+m})$ such that for all sufficiently large $\nu$ 
\begin{gather} 
	\label{relharm harmonic approx eqn11a} \widehat{h}_j^{(\nu)} = g_j^{(\nu)} \text{ on } B_{(1+\gamma)/2}(0) \setminus B_{(1+3\gamma)/4}(0), \\
	\label{relharm harmonic approx eqn11b} \op{Lip} \widehat{h}_j^{(\nu)} \leq C \big( \op{Lip} g_j^{(\nu)} + 1 \big) , \\
	\label{relharm harmonic approx eqn11c} \frac{1}{q} \int_{B_{(1+\gamma)/2}(0)} \left| \sum_{j=1}^{p_j} q_j \widehat{h}^{(\nu)}_{j,a} \right| 
		\leq C \mathcal{E}_{\nu} \text{ on } B_{(1+\gamma)/2}(0), \\
	\label{relharm harmonic approx eqn11d} \int_{B_{(1+3\gamma)/4}(0)} \mathcal{G}(g_j^{(\nu)},\widehat{h}_j^{(\nu)})^2 
		\leq C \int_{B_{(1+3\gamma)/4}(0)} (|Dg_j^{(\nu)}|^2 + |D\widehat{h}_j^{(\nu)}|^2) + C \leq 9C , \\
	\label{relharm harmonic approx eqn11e} \sum_{j=1}^p \int_{B_{(1+3\gamma)/4}(0)} |D\widehat{h}_j^{(\nu)}|^2 
		\leq \sum_{j=1}^p \int_{B_{(1+3\gamma)/4}(0)} |Dg_j^{(\nu)}|^2 - \frac{3\delta}{4} , 
\end{gather}
where $\widehat{h}^{(\nu)}_{j,a}(x) = \frac{1}{q_j} \sum_{l=1}^{q_j} \widehat{h}^{(\nu)}_{j,l}(x)$ is the single-valued average of $\widehat{h}^{(\nu)}_j(x) = \sum_{l=1}^{q_j} \llbracket \widehat{h}^{(\nu)}_{j,l}(x) \rrbracket$ and $C \in (0,\infty)$ are constants independent of $\nu$.  In particular, \eqref{relharm harmonic approx eqn11a} and \eqref{relharm harmonic approx eqn11e} follow directly from~\cite[Proposition~4.4]{DeLSpa1} using \eqref{relharm harmonic approx eqn4a}, \eqref{relharm harmonic approx eqn10a}, \eqref{relharm harmonic approx eqn10c}, and \eqref{relharm harmonic approx eqn4b}. The estimate  \eqref{relharm harmonic approx eqn11b} is (4.9) of~\cite{DeLSpa1},  and \eqref{relharm harmonic approx eqn11c} follows from (4.11) of~\cite{DeLSpa1}, \eqref{relharm harmonic approx eqn6}, and \eqref{relharm harmonic approx eqn10b}.  The bound \eqref{relharm harmonic approx eqn11d} follows from (4.10) of~\cite{DeLSpa1}, \eqref{relharm harmonic approx eqn4a}, \eqref{relharm harmonic approx eqn3}, and \eqref{relharm harmonic approx eqn11e}.

Define $\widehat{N}^{(\nu)}_j : B^{\overline{\mathcal M}^{(\nu)}}_{(1+\gamma)/2}(0) \rightarrow \mathcal{A}_{q_j}(\mathbb{R}^{n+m})$ by 
\begin{equation*}
	\widehat{N}^{(\nu)}_j(X) = \sum_{l=1}^{q_j} \llbracket \overline{\varphi}^{(\nu)}_j(X) 
		+ {\mathcal E}_{\nu} \widehat{h}^{(\nu)}_{j,l}(\overline{\mathbf e}_{\nu}^{-1}(X)) \rrbracket
\end{equation*}
for each $X \in B^{\overline{\mathcal M}^{(\nu)}}_{(1+\gamma)/2}(0)$, where $\widehat{h}^{(\nu)}_j(x) = \sum_{l=1}^{q_j} \llbracket \widehat{h}^{(\nu)}_{j,l}(x) \rrbracket$ for each $x \in B_{(1+\gamma)/2}(0)$.  By \eqref{relharm harmonic approx eqn11a}, \eqref{relharm separation1 avgfree}, \eqref{relharm harmonic approx eqn11c}, and \eqref{relharm harmonic approx eqn11d} 
\begin{gather}
	\label{relharm harmonic approx eqn12a} \widehat{N}^{(\nu)}_j = \overline{N}^{(\nu)}_j 
		\text{ on } B^{\overline{\mathcal M}^{(\nu)}}_{(1+\gamma)/2}(0) \setminus B^{\overline{\mathcal M}^{(\nu)}}_{(1+3\gamma)/4}(0) , \\
	\label{relharm harmonic approx eqn12b} \frac{1}{q} \int_{B^{\overline{\mathcal M}^{(\nu)}}_{(1+\gamma)/2}(0)} 
		\left| \sum_{j=1}^p q_j \widehat{N}^{(\nu)}_{j,a} \right| \leq C \mathcal{E}_{\nu}^2 , \\
	\label{relharm harmonic approx eqn12c} \int_{B^{\overline{\mathcal M}^{(\nu)}}_{(1+\gamma)/2}(0)} 
		\mathcal{G}(\overline{N}^{(\nu)}_j,\widehat{N}^{(\nu)}_j)^2 \leq C \mathcal{E}_{\nu}^2 
\end{gather}
where $\widehat{N}^{(\nu)}_{j,a}(X) = \frac{1}{q_j} \sum_{l=1}^{q_j} \widehat{N}^{(\nu)}_{j,l}(X)$ is the single-valued average of $\widehat{N}^{(\nu)}_j(X) = \sum_{l=1}^{q_j} \llbracket \widehat{N}^{(\nu)}_{j,l}(X) \rrbracket$ and where $C \in (0,\infty)$ are constants independent of $\nu$.  By \eqref{relharm normal approx concl2}, $|\overline{N}^{(\nu)}_j| \leq C \mathbf{m}_{0,\nu}^{1/(2n)} \rho_{\nu}^{\beta_2}$ on $B^{\overline{\mathcal M}^{(\nu)}}_{(1+\gamma)/2}(0)$ and $\op{Lip} \overline{N}^{(\nu)}_j \leq C \mathbf{m}_{0,\nu}^{\gamma_2} \rho_{\nu}^{\gamma_2}$.  Thus using $\beta_2 < \gamma_2$, $\op{Lip} g^{(\nu)}_j \leq \mathcal{E}_{\nu}^{-1} (\sqrt{q} \op{Lip} \overline{\varphi}^{(\nu)}_j + \op{Lip} \overline{N}^{(\nu)}_j)$, \eqref{relharm harmonic approx eqn11b}, and \eqref{relharm harmonic approx eqn12a} 
\begin{equation}\label{relharm harmonic approx eqn12d} 
	\sup_{B^{\overline{\mathcal M}^{(\nu)}}_{(1+\gamma)/2}(0)} |\widehat{N}^{(\nu)}_j| 
		+ \sup_{B^{\overline{\mathcal M}^{(\nu)}}_{(1+\gamma)/2}(0)} |D\widehat{N}^{(\nu)}_j| 
	\leq C \op{Lip} \varphi^{(\nu)}_j + C \mathbf{m}_{0,\nu}^{\beta_2} \rho_{\nu}^{\beta_2} + C {\mathcal E}_{\nu} , 
\end{equation}
where $C \in (0,\infty)$ is a constant independent of $\nu$.  Observe that $\|\Delta_{\overline{\mathcal M}^{(\nu)}}\overline{\varphi}^{(\nu)}_j\|_{C^0(B^{\overline{\mathcal M}^{(\nu)}}_1(0))} = \rho_{\nu} \|\Delta_{\mathcal{M}^{(\nu)}} \varphi^{(\nu)}_j\|_{C^0(B^{\mathcal{M}^{(\nu)}}_{\rho_{\nu}}(Z_{\nu}))} \leq \varepsilon_{\nu} \mathcal{E}_{\nu}$.  Hence using \eqref{relharm harmonic approx eqn11a} and \eqref{relharm harmonic approx eqn11d} 
\begin{align*}
	&\int_{B^{\overline{\mathcal M}^{(\nu)}}_{(1+3\gamma)/4}(0)} \big\langle D\overline{\varphi}^{(\nu)}_j(X), D((\widehat{h}^{(\nu)}_{j,a} - g^{(\nu)}_{j,a}) 
		\circ \overline{\mathbf e}_{\nu}^{-1})(X) \big\rangle \,d\mathcal{H}^n(X) 
	\\=\,& -\int_{B^{\overline{\mathcal M}^{(\nu)}}_{(1+3\gamma)/4}(0)} \big\langle \Delta_{\overline{\mathcal M}^{(\nu)}}\overline{\varphi}^{(\nu)}_j(X), 
		\widehat{h}^{(\nu)}_{j,a}(\overline{\mathbf e}_{\nu}^{-1}(X)) - g^{(\nu)}_{j,a}(\overline{\mathbf e}_{\nu}^{-1}(X)) \big\rangle \,d\mathcal{H}^n(X) 
	\leq C \varepsilon_{\nu} \mathcal{E}_{\nu} , 
\end{align*}
where $C = C(n,\gamma) \in (0,\infty)$ is a constant.  Hence by \eqref{relharm harmonic approx eqn3}, \eqref{relharm harmonic approx eqn11e}, and the fact that $\|D\overline{\mathbf e}_{\nu} - I\|_{C^0(B_1(0))} \leq C\mathbf{m}_{0,\nu}^{1/2} \rho_{\nu} < C\varepsilon_{\nu}^{1/2}$
\begin{align}\label{relharm harmonic approx eqn12e} 
	&\int_{B^{\overline{\mathcal M}^{(\nu)}}_{(1+3\gamma)/4}(0)} |D\widehat{N}^{(\nu)}_j|^2 
		- \int_{B^{\overline{\mathcal M}^{(\nu)}}_{(1+3\gamma)/4}(0)} |D\overline{N}^{(\nu)}_j|^2 
	\\ \leq\,& 2 {\mathcal E}_{\nu} \int_{B^{\overline{\mathcal M}^{(\nu)}}_{(1+3\gamma)/4}(0)} \big\langle D\overline{\varphi}^{(\nu)}_j(X), 
		D((\widehat{h}^{(\nu)}_{j,a} - g^{(\nu)}_{j,a}) \circ \overline{\mathbf e}_{\nu}^{-1})(X) \big\rangle \,d\mathcal{H}^n(X) \nonumber
		\\&+ {\mathcal E}_{\nu}^2 \int_{B_{(1+3\gamma)/4}(0)} (|D\widehat{h}^{(\nu)}_j|^2 - |Dg^{(\nu)}_j|^2) \,J\overline{\mathbf e}_{\nu} \nonumber 
	\\ \leq\,& -\frac{\delta}{2} {\mathcal E}_{\nu}^2 \nonumber
\end{align}
for all sufficiently large $\nu$, where $C = C(n,m,q,M_0,N_0,C_e,C_h) \in (0,\infty)$ are constants.

Define $\widetilde{N}^{(\nu)}_j : B^{\overline{\mathcal M}^{(\nu)}}_{(1+\gamma)/2}(0) \rightarrow \mathcal{A}_{q_j}(({\rm T}\overline{\mathcal M}^{(\nu)})^{\perp})$ by 
\begin{equation*}
	\widetilde{N}^{(\nu)}_j(X) = \sum_{l=1}^{q_j} \llbracket \pi_{({\rm T}_X\overline{\mathcal M}^{(\nu)})^{\perp}} (\widehat{N}^{(\nu)}_{j,l}(X)) \rrbracket
\end{equation*}
for each $X \in B^{\overline{\mathcal M}^{(\nu)}}_{(1+\gamma)/2}(0)$, where $\widehat{N}^{(\nu)}_j(X) = \sum_{l=1}^{q} \llbracket \widehat{N}^{(\nu)}_{j,l}(X) \rrbracket$.  Recall from the proof of Lemma~\ref{centered blowup lemma} that $|A_{\overline{\mathcal M}^{(\nu)}}| \leq C \mathbf{m}_{0,\nu}^{1/2} \rho_{\nu}$ and $|D\widetilde{N}^{(\nu)}_j| \leq |D\widehat{N}^{(\nu)}_j| + C(n,m) \,|A_{\overline{\mathcal M}^{(\nu)}}| \,|\widehat{N}^{(\nu)}_j|$ on $B^{\overline{\mathcal M}^{(\nu)}}_{(1+\gamma)/2}(0)$.  Hence by \eqref{relharm harmonic approx eqn12a}, \eqref{relharm harmonic approx eqn12b}, \eqref{relharm harmonic approx eqn12c}, \eqref{relharm harmonic approx eqn12d}, and \eqref{relharm harmonic approx eqn12e} 
\begin{gather}
	\label{relharm harmonic approx eqn13a} \widetilde{N}^{(\nu)}_j = \overline{N}^{(\nu)}_j 
		\text{ on } B^{\overline{\mathcal M}^{(\nu)}}_{(1+\gamma)/2}(0)) \setminus B^{\overline{\mathcal M}^{(\nu)}}_{(1+3\gamma)/4}(0) , \\
	\label{relharm harmonic approx eqn13b} \frac{1}{q} \int_{B^{\overline{\mathcal M}^{(\nu)}}_{(1+\gamma)/2}(0)} 
		\left| \sum_{j=1}^p \widetilde{N}^{(\nu)}_{j,a} \right| \leq C \mathcal{E}_{\nu}^2 , \\
	\label{relharm harmonic approx eqn13c} \int_{B^{\overline{\mathcal M}^{(\nu)}}_{(1+\gamma)/2}(0)} 
		\mathcal{G}(\overline{N}^{(\nu)}_j,\widetilde{N}^{(\nu)}_j)^2 \leq C \mathcal{E}_{\nu}^2 , \\
	\label{relharm harmonic approx eqn13d} \sup_{B^{\overline{\mathcal M}^{(\nu)}}_{(1+\gamma)/2}(0)} |\widetilde{N}^{(\nu)}_j| 
		+ \sup_{B^{\overline{\mathcal M}^{(\nu)}}_{(1+\gamma)/2}(0)} |D\widetilde{N}^{(\nu)}_j| 
		\leq C \op{Lip} \varphi^{(\nu)}_j + C \mathbf{m}_{0,\nu}^{\beta_2} \rho_{\nu}^{\beta_2} + C {\mathcal E}_{\nu} , \\
	\label{relharm harmonic approx eqn13e} \int_{B^{\overline{\mathcal M}^{(\nu)}}_{(1+3\gamma)/4}(0)} |D\widetilde{N}^{(\nu)}_j|^2
		\leq \int_{B^{\overline{\mathcal M}^{(\nu)}}_{(1+3\gamma)/4}(0)} |D\overline{N}^{(\nu)}_j|^2 - \frac{3\delta}{8} {\mathcal E}_{\nu}^2 
\end{gather}
where $\widetilde{N}^{(\nu)}_{j,a}(X) = \frac{1}{q_j} \sum_{l=1}^{q_j} \widetilde{N}^{(\nu)}_{j,l}(X)$ is the single-valued average of $\widetilde{N}^{(\nu)}_j(X) = \sum_{l=1}^{q_j} \llbracket \widetilde{N}^{(\nu)}_{j,l}(X) \rrbracket$ and where $C \in (0,\infty)$ are constants independent of $\nu$. 

By $\varsigma \leq \gamma_2/4$ and \eqref{relharm norm est rmk eqn1}
\begin{equation*}
	\|\overline{T}^{(\nu)}_j - \mathbf{G}_{\overline{N}^{(\nu)}_j}\|(\mathbf{C}^{\overline{\mathcal M}^{(\nu)}}_{(1+\gamma)/2}(0)) \leq C \mathbf{m}_{0,\nu}^{3\gamma_2/4} \rho_{\nu}^{-n-2+3\gamma_2/2} \int_{B^{\mathcal{M}^{(\nu)}}_{\rho_{\nu}}(Z_{\nu})} |N^{(\nu)}|^2 \leq C \mathbf{m}_{0,\nu}^{3\gamma_2/4} \rho_{\nu}^{\gamma_2/2} \mathcal{E}_{\nu}^2 , 
\end{equation*}
where $C = C(n,m,q,M_0,N_0,C_e,C_h,\gamma) \in (0,\infty)$ are constants.  By the coarea formula and Fatou's lemma 
\begin{equation*}
    \int_{(1+3\gamma)/4}^{(1+\gamma)/2} \liminf_{\nu\to\infty} \frac{\mathbf{M}(\langle \overline{T}^{(\nu)}_j - \mathbf{G}_{\overline{N}^{(\nu)}_j}, R_{\nu} \circ \overline{\mathbf p}_{\nu}, s \rangle)}{\mathbf{m}_{0,\nu}^{3\gamma_2/4} \rho_{\nu}^{\gamma_2/2} \mathcal{E}_{\nu}^2} \,ds \leq \liminf_{\nu\to\infty} \frac{C \|\overline{T}^{(\nu)}_j - \mathbf{G}_{\overline{N}^{(\nu)}_j}\|(\mathbf{C}^{\overline{\mathcal M}^{(\nu)}}_{(1+\gamma)/2}(0))}{\mathbf{m}_{0,\nu}^{3\gamma_2/4} \rho_{\nu}^{\gamma_2/2} \mathcal{E}_{\nu}^2} \leq C , 
\end{equation*}
where $R_{\nu} = R_{\nu}(X)$ is the geodesic distance from $X$ to the origin in $\overline{\mathcal M}^{(\nu)}$, $\overline{\mathbf p}_{\nu}$ is the closest point projection map onto $\overline{\mathcal M}^{(\nu)}$, and $C = C(n,m,q,M_0,N_0,C_e,C_h,\gamma) \in (0,\infty)$ are constants.  Hence after passing to a subsequence, there exists $\overline{\gamma} \in \big(\frac{1+3\gamma}{4},\frac{1+\gamma}{2}\big)$ such that the slice $\langle \overline{T}^{(\nu)}_j - \mathbf{G}_{\overline{N}^{(\nu)}_j}, R_{\nu} \circ \overline{\mathbf p}_{\nu}, \overline{\gamma} \rangle$ exists and is an $(n-1)$-dimensional rectifiable current with 
\begin{equation*}
    \mathbf{M}(\langle \overline{T}^{(\nu)}_j - \mathbf{G}_{\overline{N}^{(\nu)}_j}, R_{\nu} \circ \overline{\mathbf p}_{\nu}, \overline{\gamma} \rangle) 
	\leq C \mathbf{m}_{0,\nu}^{3\gamma_2/4} \rho_{\nu}^{\gamma_2/2} \mathcal{E}_{\nu}^2 , 
\end{equation*}
where $C = C(n,m,q,M_0,N_0,C_e,C_h,\gamma) \in (0,\infty)$ is a constant.  
Thus by the isoperimetric inequality, there exists an $n$-dimensional area-minimizing rectifiable current $S^{(\nu)}_j$ of $\mathbb{R}^{n+m}$ such that 
\begin{equation}\label{relharm harmonic approx eqn17} 
	\partial S^{(\nu)}_j = \langle \overline{T}^{(\nu)}_j - \mathbf{G}_{\overline{N}^{(\nu)}_j}, R_{\nu} \circ \overline{\mathbf p}_{\nu}, \overline{\gamma} \rangle , \quad 
	\mathbf{M}(S^{(\nu)}_j) \leq C \mathbf{m}_{0,\nu}^{3\gamma_2/4} \rho_{\nu}^{\gamma_2/2} \mathcal{E}_{\nu}^2 ,  
\end{equation}
where $C = C(n,m,q,M_0,N_0,C_e,C_h,\gamma) \in (0,\infty)$ is a constant.  Let $\widetilde{T}^{(\nu)}_j$ be a competitor for $\overline{T}^{(\nu)}_j \llcorner \mathbf{C}^{\overline{\mathcal M}^{(\nu)}}_1(0)$ defined by 
\begin{equation*}
	\widetilde{T}^{(\nu)}_j = \mathbf{G}_{\widetilde{N}^{(\nu)}_j} \llcorner \mathbf{C}^{\overline{\mathcal M}^{(\nu)}}_{\overline{\gamma}}(0)) + S^{(\nu)}_j 
		+ \overline{T}^{(\nu)}_j \llcorner \mathbf{C}^{\overline{\mathcal M}^{(\nu)}}_1(0) \setminus 
		\mathbf{C}^{\overline{\mathcal M}^{(\nu)}}_{\overline{\gamma}}(0) . 
\end{equation*}
Using \eqref{relharm harmonic approx eqn13a} and \eqref{relharm harmonic approx eqn17}, it is readily verified that $(\partial \widetilde{T}^{(\nu)}_j) \llcorner \mathbf{C}^{\overline{\mathcal M}^{(\nu)}}_1(0) = 0$ and $\widetilde{T}^{(\nu)}_j \llcorner \mathbf{C}^{\overline{\mathcal M}^{(\nu)}}_1(0) \setminus \mathbf{C}^{\overline{\mathcal M}^{(\nu)}}_{\overline{\gamma}}(0) = T^{(\nu)}_j \llcorner \mathbf{C}^{\overline{\mathcal M}^{(\nu)}}_1(0) \setminus \mathbf{C}^{\overline{\mathcal M}^{(\nu)}}_{\overline{\gamma}}(0)$.  By \eqref{relharm norm est rmk eqn1}, \eqref{relharm harmonic approx eqn13a}, \eqref{relharm harmonic approx eqn17}, and~\cite[Theorem~3.2]{DeLSpa0} 
\begin{align*}
	&\mathbf{M}(\widetilde{T}^{(\nu)}_j) - \mathbf{M}(\overline{T}^{(\nu)}_j) 
	\\ \leq\,& \mathbf{M}(\mathbf{G}_{\widetilde{N}^{(\nu)}_j} \llcorner \mathbf{C}^{\overline{\mathcal M}^{(\nu)}}_{\overline{\gamma}}(0)) 
		- \mathbf{M}(\mathbf{G}_{\overline{N}^{(\nu)}_j} \llcorner \mathbf{C}^{\overline{\mathcal M}^{(\nu)}}_{\overline{\gamma}}(0)) 
		+ C \mathbf{m}_{0,\nu}^{3\gamma_2/4} \rho_{\nu}^{\gamma_2/2} {\mathcal E}_{\nu}^2 \nonumber 
	\\ \leq\,& -q_j \int_{B^{\overline{\mathcal M}^{(\nu)}}_{\overline{\gamma}}(0)} \langle H_{\overline{\mathcal M}^{(\nu)}}, \widetilde{N}^{(\nu)}_{j,a} 
		- \overline{N}^{(\nu)}_{j,a} \rangle  + \frac{1}{2} \int_{B^{\overline{\mathcal M}^{(\nu)}}_{\overline{\gamma}}(0)} |D\widetilde{N}^{(\nu)}_j|^2 
			- \frac{1}{2} \int_{B^{\overline{\mathcal M}^{(\nu)}}_{\overline{\gamma}}(0)} |D\overline{N}^{(\nu)}_j|^2 \nonumber
		\\&+ C \int_{B^{\overline{\mathcal M}^{(\nu)}}_{\overline{\gamma}}(0)} (|A_{\overline{\mathcal M}^{(\nu)}}|^2 |\overline{N}^{(\nu)}_j|^2 
			+ |A_{\overline{\mathcal M}^{(\nu)}}|^2 |\widetilde{N}^{(\nu)}_j|^2 + |D\overline{N}^{(\nu)}_j|^4 + |D\widetilde{N}^{(\nu)}_j|^4) 
			+ C \mathbf{m}_{0,\nu}^{3\gamma_2/4} \rho_{\nu}^{\gamma_2/2} {\mathcal E}_{\nu}^2 , \nonumber
\end{align*}
where $C = C(n,m,q,M_0,N_0,C_e,C_h,\gamma) \in (0,\infty)$ is a constant.  By $\varsigma = 2\beta_2$, $\mathbf{m}_{0,\nu} \rho_{\nu} < \varepsilon_{\nu}$, $\|A_{\overline{\mathcal M}^{(\nu)}}\|_{C^0(\overline{\mathcal M}^{(\nu)})} \leq C \mathbf{m}_{0,\nu}^{1/2} \rho_{\nu}$, $\op{Lip} \overline{N}^{(\nu)}_j \leq C \mathbf{m}_{0,\nu}^{\gamma_2} \rho_{\nu}^{\gamma_2}$, Lemma~\ref{normal approx energy lemma}, \eqref{relharm normal approx concl2}, \eqref{relharm norm est rmk eqn3}, and \eqref{relharm harmonic approx eqn13b}--\eqref{relharm harmonic approx eqn13e} 
\begin{equation*}
	\sum_{j=1}^p (\mathbf{M}(\widetilde{T}^{(\nu)}_j) - \mathbf{M}(\overline{T}^{(\nu)}_j)) \leq -\frac{\delta}{8} {\mathcal E}_{\nu}^2 < 0
\end{equation*}
for all equation large $\nu$.  This contradicts $\overline{T}^{(\nu)}_j$ being locally area-minimizing in $\mathbf{C}^{\overline{\mathcal M}^{(\nu)}}_1(0)$ for all $j \in \{1,2,\ldots,p\}$.  

\noindent\textit{Step~4.  Convergence of the Dirichlet energy of $g_j^{(\nu)}$.}  Finally, we want to show that  we claim that \eqref{relharm harmonic approx eqn20} holds true for all $j \in \{1,2,\ldots,p\}$.  In light of \eqref{relharm harmonic approx eqn4b}, it suffices to assume to the contrary that after passing to a subsequence there exists $\delta > 0$ such that 
\begin{equation}\label{relharm harmonic approx eqn21} 
	\sum_{j=1}^p \sum_{k=1}^{p_j} \int_{B_{\gamma}(0)} |Dw_{j,k}|^2 \leq \sum_{j=1}^p \int_{B_{\gamma}(0)} |Dg_j^{(\nu)}|^2 - \delta . 
\end{equation}
In light of \eqref{relharm harmonic approx eqn21}, we can apply~\cite[Proposition~4.4]{DeLSpa1} to interpolate between $h^{(\nu)}$ and $g^{(\nu)}$ and find a Lipschitz $q$-valued function $\widehat{h}^{(\nu)} : B_{(1+\gamma)/2}(0) \rightarrow \mathcal{A}_q(\mathbb{R}^{n+m})$ such that \eqref{relharm harmonic approx eqn11a}--\eqref{relharm harmonic approx eqn11d} hold true.  Repeating the above argument, we again derive a contradiction with $\overline{T}^{(\nu)}$ being locally area-minimizing in $\mathbf{C}^{\overline{\mathcal M}^{(\nu)}}_1(0)$.  Therefore, \eqref{relharm harmonic approx eqn20} holds true.  
\end{proof}

\subsection{$L^{\infty}$-distance estimate in terms of tilt excess relative to a single-valued $C^2$-graph}\label{sec:relharm separation1 sec}  For $j \in \{1,2,\ldots,p\}$, let $T_j$ and $\varphi_j$ be as in Hypothesis~\ref{relharm separation1 hyp}.  In this subsection, we want to prove Theorem~\ref{relharm separation1 thm}, in which we bound the $L^{\infty}$-distance of $T_j \llcorner \mathbf{C}^{\mathcal M}_{\gamma\rho}(Z)$ to some union of graphs parallel to the graph of $\varphi_j$ in terms of the excess $\mathcal{E}_{Z,(1+\gamma)\rho/2,\rho}$ (as in \eqref{relharm excess defn1}).  This is analogous to~\cite[Theorem~2.12]{KrumWicb}.  

\begin{lemma}[Controlled growth of excess]\label{relharm tilt growth lemma}
Let $M_0,N_0,C_e,C_h,\varepsilon_2$ be as in Hypothesis~\ref{DeLSpa constants hypothesis} and Lemma~\ref{normal approx energy lemma}.  For every integer $p \in \{1,2,\ldots,q\}$ and every $\theta \in (0,1/16)$ there exist $\varepsilon = \varepsilon(n,m,q,p,$ $M_0,N_0,C_e,C_h,\theta) \in (0,1)$ such that the following holds true.  Assume that Hypothesis~\ref{relharm separation1 hyp} holds true and that $\rho > \overline{R}(Z)$, where $\overline{R}(Z)$ is as in Definition~\ref{crit rad defn}.  For each $j \in \{1,2,\ldots,p\}$, let $N_j : B^{\mathcal M}_{5\rho/8}(Z) \rightarrow \mathcal{A}_{q_j}({\rm T}\mathcal{M}^{\perp})$ satisfy \eqref{relharm normal approx concl1} and \eqref{relharm normal approx concl2} with $\gamma = 5/8$.  
Then 
\begin{equation}\label{relharm tilt growth concl} 
	\mathcal{E}_{Z,5\theta\rho/8,\theta\rho} \leq C \theta^{\mu-1} \mathcal{E}_{Z,5\rho/8,\rho} ,
\end{equation}
where $\mathcal{E}_{Z,5\theta\rho/8,\theta\rho}$ and $\mathcal{E}_{Z,5\rho/8,\rho}$ are as in \eqref{relharm excess defn1} with $\gamma = 5/8$ and $\varsigma = 2\beta_2$ and where $\mu = \mu(n,m,q) \in (0,1)$ and $C = C(n,m,q) \in (0,\infty)$ are constants.
\end{lemma}

\begin{proof} 
Let $\eta = \eta(n,m,q,p,M_0,N_0,C_e,C_h,\theta) \in (0,1)$ to be later determined.  By Lemma~\ref{relharm harmonic approx lemma}, for each $j \in \{1,2,\ldots,p\}$ there exists a Dirichlet energy minimizing $q_j$-valued function $w_j : B_{\rho/4}(0) \rightarrow \mathcal{A}_{q_j}(\mathbb{R}^{n+m})$ such that \eqref{relharm harmonic approx concl} holds true with $\gamma = 1/4$.  By \cite[(3.20)]{DeLSpaDirMin} (\cite[Theorem~2.13]{Almgren}), there exists a constant $\mu = \mu(n,m,q) \in (0,1)$ such that 
\begin{equation}\label{relharm tilt growth eqn1}
	\frac{1}{(5\theta/8)^{n-2}} \int_{B_{5\theta\rho/8}(0)} |Dw_j|^2 \leq 4^{n-2+2\mu} (5\theta/8)^{2\mu} \int_{B_{\rho/4}(0)} |Dw_j|^2 .
\end{equation}
Thus by \eqref{relharm harmonic approx concl} 
\begin{align*}
	\frac{1}{(\theta\rho)^n} \int_{B^{\mathcal M}_{5\theta\rho/8}(Z)} \mathcal{G}(DN_j,q\llbracket D\varphi_j \rrbracket)^2 
	\leq\,& \frac{4 (5/2)^n \theta^{2\mu-2}}{\rho^n} \int_{\mathbf{C}^{\mathcal M}_{\rho/4}(Z)} \mathcal{G}(DN_j,q\llbracket D\varphi_j \rrbracket)^2
		\\&+ C \eta^2 (\theta^{2\mu-2} + \theta^{-n}) \,\mathcal{E}_{Z,5\rho/8,\rho}^2 , 
\end{align*}
for each $j \in \{1,2,\ldots,p\}$, where $C = C(n,m,q,M_0,N_0,C_e,C_h) \in (0,\infty)$ is a constant.  By choosing $\eta$ sufficiently small 
\begin{equation}\label{relharm tilt growth eqn2}
	\frac{1}{(\theta\rho)^n} \int_{B^{\mathcal M}_{5\theta\rho/8}(Z))} \mathcal{G}(DN_j,q\llbracket D\varphi_j \rrbracket)^2 
	\leq C \theta^{2\mu-2} \mathcal{E}_{Z,5\rho/8,\rho}^2 ,  
\end{equation}
where $C = 4 (5/2)^n + 1$.  

Let 
\begin{equation*}
	\widehat{E}_{Z,\rho}^2 = \frac{1}{\rho^n} \int_{B^{\mathcal M}_{5\rho/8}(Z)} |DN|^2 
		+ \frac{1}{\rho^{n+2-\varsigma}} \int_{B^{\mathcal M}_{\rho}(Z)} |N|^2 . 
\end{equation*}
By Lemma~\ref{relharm harmonic approx lemma} (in the special case where $p = 1$ and $\varphi_1$ is identically zero), there exists a Dirichlet energy minimizing $q$-valued function $w : B_{\rho/4}(0) \rightarrow \mathcal{A}_q(\mathbb{R}^{n+m})$ such that 
\begin{equation}\label{relharm tilt growth eqn3}
	\frac{1}{\rho^{n+2}} \int_{B^{\mathcal M}_{\rho/4}(Z)} \mathcal{G}(N,w \circ \mathbf{e}^{-1})^2 
	+ \frac{1}{\rho^n} \int_{B^{\mathcal M}_{\rho/4}(Z)} (|DN| - |D(w \circ \mathbf{e}^{-1})|)^2 
	\leq \eta^2 \widehat{E}_{Z,\rho}^2 . 
\end{equation}
By $L^2$-growth estimates for Dirichlet energy minimizing multi-valued functions \cite[Corollary~3.18]{DeLSpaDirMin} (\cite[Theorem~2.6(8)]{Almgren}) 
\begin{equation}\label{relharm tilt growth eqn4}
	\frac{1}{\theta^n} \int_{B_{\theta\rho}(0)} |w|^2 \leq 4^n \int_{B_{\rho/4}(0)} |w|^2 .
\end{equation}
By \cite[(3.20)]{DeLSpaDirMin} (\cite[Theorem~2.13]{Almgren}), there exists a (possibly smaller) constant $\mu = \mu(n,m,q) \in (0,1)$ such that 
\begin{equation}\label{relharm tilt growth eqn5}
	\theta^{2-n} \int_{B_{\theta\rho}(0)} |Dw|^2 \leq 4^{n-2+2\mu} \theta^{2\mu} \int_{B_{\rho/4}(0)} |Dw|^2 
\end{equation}
(cf.~\eqref{relharm tilt growth eqn1}).  By \eqref{relharm tilt growth eqn3}, \eqref{relharm tilt growth eqn4}, and \eqref{relharm tilt growth eqn5} 
\begin{align*}
	\frac{1}{(\theta\rho)^{n+2}} \int_{B^{\mathcal M}_{\theta\rho}(Z)} |N|^2 
		\leq\,& \frac{4^{n+1}}{\theta^2 \rho^{n+2}} \int_{B^{\mathcal M}_{\rho/4}(Z)} |N|^2 
		+ C \eta^2 (\theta^{-n-2} + \theta^{-2}) \,\widehat{E}_{Z,\rho}^2 , \\
	\frac{1}{(\theta\rho)^n} \int_{B^{\mathcal M}_{\theta\rho}(Z)} |DN|^2 
		\leq\,& \frac{4^{n+1} \theta^{2\mu-2}}{\rho^n} \int_{B^{\mathcal M}_{\rho/4}(Z)} |DN|^2 
		+ C \eta^2 (\theta^{-n} + \theta^{2\mu-2}) \,\widehat{E}_{Z,\rho}^2
\end{align*}
where $C = C(n,m,q,p,M_0,N_0,C_e,C_h) \in (0,\infty)$ is a constant.  Note that by Lemma~\ref{normal approx energy lemma}, $\widehat{E}_{Z,\rho}^2 \leq C \rho^{-n-2} \int_{B^{\mathcal M}_{\rho}(Z)} |N|^2$, where $C = C(n,m,q,p,M_0,N_0,C_e,C_h) \in (0,\infty)$ is a constant.  Hence by taking $\eta$ to be sufficiently small 
\begin{align}
	\label{relharm tilt growth eqn6} \frac{1}{(\theta\rho)^{n+2}} \int_{B^{\mathcal M}_{\theta\rho}(Z)} |N|^2 
		\leq\,& \frac{C}{\theta^2 \rho^{n+2}} \int_{B^{\mathcal M}_{\rho}(Z)} |N|^2 , \\
	\label{relharm tilt growth eqn7} \frac{1}{(\theta\rho)^n} \int_{B^{\mathcal M}_{\theta\rho}(Z)} |DN|^2 
		\leq\,& \frac{C \theta^{2\mu-2}}{\rho^n} \int_{B^{\mathcal M}_{\rho}(Z)} |DN|^2 
		+ \frac{\theta^{2\mu-2}}{\rho^{n+2-\varsigma}} \int_{B^{\mathcal M}_{\rho}(Z)} |N|^2
\end{align}
where $C = 4^{n+1} + 1$.  Combining \eqref{relharm tilt growth eqn2}, \eqref{relharm tilt growth eqn6}, and \eqref{relharm tilt growth eqn7} gives us \eqref{relharm tilt growth concl} with $\min\{\mu,\varsigma/2\}$ in place of $\mu$.
\end{proof}

\begin{lemma}[Poincar\'e-type inequality]\label{relharm poincare lemma} 
Let $M_0,N_0,C_e,C_h,\varepsilon_2$ be as in Hypothesis~\ref{DeLSpa constants hypothesis} and Lemma~\ref{normal approx energy lemma}.  For every integer $p \in \{1,2,\ldots,q\}$ there exist $\varepsilon = \varepsilon(n,m,q,p,M_0,N_0,C_e,C_h) \in (0,1)$ such that the following holds true.  Assume that Hypothesis~\ref{relharm separation1 hyp} holds true and that $\rho > \overline{R}(Z)$, where $\overline{R}(Z)$ is as in Definition~\ref{crit rad defn}.  For each $j \in \{1,2,\ldots,p\}$, let $\mathcal{K}_j \subseteq B^{\mathcal M}_{(1+\gamma)\rho/2}(Z)$ and $N_j : B^{\mathcal M}_{\gamma\rho}(Z) \rightarrow \mathcal{A}_{q_j}({\rm T}\mathcal{M}^{\perp})$ satisfy \eqref{relharm normal approx concl1} and \eqref{relharm normal approx concl2}.  Then for each $j \in \{1,2,\ldots,p\}$ there exists $a_j \in \mathcal{A}_{q_j}(\mathbb{R}^{n+m})$ such that $|a_j| < C\varepsilon^{\beta_2}\rho$ and 
\begin{align}\label{relharm poincare concl}
	\frac{1}{\rho^{n+2}} \int_{\mathbf{C}^{\mathcal M}_{\gamma\rho}(Z)} \op{dist}^2(X, \op{spt} \mathbf{G}_{\varphi_j+a_j}) \,d\|T_j\|(X) & 
	\\+ \frac{1}{\rho^{n+2}} \int_{B^{\mathcal M}_{\gamma\rho}(Z)} \mathcal{G}(N_j,\varphi_j+a_j)^2  
	&\leq C \mathcal{E}_{Z,\gamma\rho,\rho}^2 , \nonumber
\end{align}
where $\mathcal{E}_{Z,\gamma\rho,\rho}$ is as in \eqref{relharm excess defn1} with $\varsigma = 3\gamma_2/4$ and where $C = C(n,m,q,M_0,N_0,C_e,C_h,\gamma) \in (0,\infty)$ is a constant.
\end{lemma}

\begin{proof} 
By the Poincar\'e inequality~\cite[Proposition~2.11]{DeLSpaDirMin} applied to $(N_j - \varphi_j) \circ \exp_Z$ (where $\exp_Z$ is the exponential map of $\mathcal{M}$ at $Z$), there exists $a_j \in \mathcal{A}_{q_j}(\mathbb{R}^{n+m})$ such that 
\begin{equation}\label{relharm poincare eqn1}
	\frac{1}{\rho^{n+2}} \int_{B^{\mathcal M}_{\gamma\rho}(Z)} \mathcal{G}(N_j,\varphi_j+a_j)^2 
	\leq \frac{C}{\rho^n} \int_{B^{\mathcal M}_{\gamma\rho}(Z)} \mathcal{G}(DN_j, q_j \llbracket D\varphi_j \rrbracket)^2 
	\leq 4C \mathcal{E}_{Z,\gamma\rho,\rho}^2 , 
\end{equation}
where $C = C(n,m,q) \in (0,\infty)$ is a constant.  Assuming that we chose $a_j$ to minimize $\int_{B^{\mathcal M}_{\gamma\rho}(Z)} \mathcal{G}(N_j,\varphi_j+a_j)^2$, by \eqref{relharm separation1 small norm} and \eqref{relharm normal approx concl2} we have that $|a_j| < C\mathbf{m}_0^{1/(2n)} \rho^{1+\beta_2} + C\varepsilon\rho < C\varepsilon^{\beta_2}\rho$, where $C = C(n,m,q,M_0,$ $N_0,C_e,C_h,\gamma) \in (0,\infty)$ are constants.  Recall that by~\cite[Corollary~2.2]{DeLSpa2}, 
\begin{equation*}
	\sup_{X \in \op{spt} T \cap \mathbf{C}^{\mathcal M}_{\rho}(Z)} \op{dist}(X,\mathcal{M}) \leq C\mathbf{m}_0^{1/(2n)} \rho^{1+\beta_2} , 
\end{equation*}
where $C = C(n,m,q,M_0,N_0,C_e,C_h) \in (0,\infty)$ is a constant.  Hence using \eqref{relharm separation1 small norm}, \eqref{relharm normal approx concl1}, \eqref{relharm normal approx concl2}, \eqref{relharm norm est rmk eqn1}, and \eqref{relharm poincare eqn1} 
\begin{align*}
	&\frac{1}{\rho^{n+2}} \int_{\mathbf{C}^{\mathcal M}_{\gamma\rho}(Z)} \op{dist}^2(X, \op{spt} \mathbf{G}_{\varphi_j+a_j}) \,d\|T_j\|(X) 
	\\ \leq\,& \frac{C}{\rho^{n+2}} \sum_{j=1}^p \int_{B^{\mathcal M}_{\gamma\rho}(Z)} \mathcal{G}(N_j,\varphi_j+a_j)^2 
		+ C \mathbf{m}_0^{3\gamma_2/4} \rho^{-n-2+3\gamma_2/4} \int_{B^{\mathcal M}_{\gamma\rho}(Z)} |N|^2 
	\leq C \mathcal{E}_{Z,\gamma\rho,\rho}^2
\end{align*}
where $C = C(n,m,q,M_0,N_0,C_e,C_h,\gamma) \in (0,\infty)$ are constants, completing the proof of \eqref{relharm poincare concl}.
\end{proof}

The next result, Lemma~\ref{relharm critrad1 lemma}, together with Lemma~\ref{relharm critrad2 lemma} allows us to bound the $L^{\infty}$-distance of $T_j$ to a union of graphs parallel to the graph of $\varphi_j$ at the scale $\rho = \overline{R}(Z)$ where \eqref{freq go condition} fails.  In Lemma~\ref{relharm critrad1 lemma}, we express $\varphi_j$ as graphs of functions over the tangent plane to $\mathcal{M}$ at $Z$ and then use~\cite[Theorem~4.19]{KrumWicc} to bound the uniform distance of $T_j$ to $\varphi_j$.  The combined use of Lemma~\ref{relharm critrad1 lemma} and Lemma~\ref{relharm critrad2 lemma} will play a role analogous to~\cite[Lemma~3.6]{KrumWicb} in the proof of~\cite[Theorem~3.10]{KrumWicb}.

\begin{lemma}\label{relharm critrad1 lemma} 
Let $M_0,N_0,C_e,C_h,\varepsilon_2$ be as in Hypothesis~\ref{DeLSpa constants hypothesis} and Lemma~\ref{normal approx energy lemma}.  For every integer $p \in \{1,2,\ldots,q\}$ there exist $\varepsilon = \varepsilon(n,m,q,p,M_0,N_0,C_e,C_h) \in (0,1)$ such that the following holds true.  Assume that Hypothesis~\ref{relharm separation1 hyp} holds true, 
\begin{equation}\label{relharm separation1 normal} 
	\varphi_j(X) \in {\rm T}_X \mathcal{M}^{\perp} 
\end{equation}
for all $X \in B^{\mathcal M}_{\rho}(Z)$, and $\rho \geq \overline{R}(Z)$ where $\overline{R}(Z)$ is as in Definition~\ref{crit rad defn}.  For each $j \in \{1,2,\ldots,p\}$, let $a_j \in \mathcal{A}_{q_j}({\rm T}_Z \mathcal{M}^{\perp})$ with $|a_j| < \varepsilon\rho$.  Then for each $j \in \{1,2,\ldots,p\}$  
\begin{align}\label{relharm critrad1 concl}
	&\sup_{X \in \op{spt} T_j \cap \mathbf{C}^{\mathcal M}_{\rho/2}(Z)} \op{dist}^2(X, \op{spt} \mathbf{G}_{\varphi_j+a_j}) 
	\leq \frac{C}{\rho^n} \int_{\mathbf{C}^{\mathcal M}_{\rho}(Z)} \op{dist}^2(X, \op{spt} \mathbf{G}_{\varphi_j+a_j}) \,d\|T_j\|(X) 
	\\&\hspace{15mm} + C \rho^4 \sum_{j=1}^p \|\Delta_{\mathcal M} \varphi_j\|_{C^0(B^{\mathcal M}_{\rho}(Z))}^2  
	+ C \rho^2 \sum_{j=1}^p \|D\varphi_j\|_{C^0(B^{\mathcal M}_{\rho}(Z))} \cdot \frac{1}{\rho^n} \int_{B^{\mathcal M}_{\rho}(Z)} |DN|^2 \nonumber 
	\\&\hspace{15mm} + \frac{1}{\rho^{n-\varsigma}} \int_{B^{\mathcal M}_{\rho}(Z)} |N|^2 + C \mathbf{m}_0 \rho^4 , 
	\nonumber 
\end{align}
where $\varsigma = \varsigma(n,m,q) \in (0,1)$ and $C = C(n,m,q,M_0,N_0,C_e,C_h) \in (0,\infty)$ are constants.
\end{lemma}

\begin{proof} 
Set $P = {\rm T}_Z \mathcal{M}$.  By Lemma~\ref{app lemma}, for each $j \in \{1,2,\ldots,p\}$ there exists a (unique) single-valued function $\psi_j : B_{3\rho/4}(Z,P) \rightarrow P^{\perp}$ such that 
\begin{equation}\label{relharm critrad1 eqn1}
	\mathbf{G}_{\psi_j} = \mathbf{G}_{\varphi_j} \llcorner \mathbf{C}_{3\rho/4}(Z,P) . 
\end{equation} 
Since $a_j \in \mathcal{A}_q({\rm T}_Z\mathcal{M}^{\perp})$ 
\begin{equation}\label{relharm critrad1 eqn2}
	\mathbf{G}_{\psi_j+a_j} = \mathbf{G}_{\varphi_j+a_j} \llcorner \mathbf{C}_{3\rho/4}(Z,P) . 
\end{equation} 
We want to apply Lemma~\ref{app lemma} to obtain $C^2$-estimates for $\psi_j$.  First let's bound the various terms on the right-hand side of \eqref{app concl3}--\eqref{app concl6}.  Note that the right-hand side of \eqref{app concl3}--\eqref{app concl6} simplify since $P$ has zero second fundamental form.  For each $j \in \{1,2,\ldots,p\}$, define $\xi_j : B_{3\rho/4}(Z,P) \rightarrow \mathcal{M}$ be the unique injective map such that 
\begin{equation}\label{relharm critrad1 eqn3}
	X + \psi_j(X) = \xi_j(X) + \varphi_j(\xi_j(X)) 
\end{equation}
for each $X \in B_{3\rho/4}(Z,P)$.  Since $P = {\rm T}_Z \mathcal{M}$ and $\|D\bs{\varphi}\|_{C^{2,\kappa}} \leq C \mathbf{m}_0^{1/2}$ 
\begin{equation}\label{relharm critrad1 eqn4}
	\|\pi_{{\rm T}_X \mathcal{M}} - \pi_P\| \leq C \mathbf{m}_0^{1/2} \rho 
\end{equation}
for all $X \in B^{\mathcal M}_{\rho}(Z)$, where $C = C(n,m,q,M_0,N_0,C_e,C_h) \in (0,\infty)$ are constants.  Since $\varphi_j(X) \in {\rm T}_X\mathcal{M}^{\perp}$ for each $X \in B^{\mathcal M}_{\rho}(Z)$, it follows from the definition of $A_{\mathcal M}$ and \eqref{relharm separation1 small norm} that 
\begin{equation}\label{relharm critrad1 eqn5}
	|\pi_{{\rm T}_X \mathcal{M}}(D_{\vec\tau} \varphi_j(X))| \leq C |A_{\mathcal M}| \,|\varphi_j(X)| \leq C \mathbf{m}_0^{1/2} \rho 
\end{equation}
for all $X \in B^{\mathcal M}_{\rho}(Z)$ and $\vec\tau \in {\rm T}_X \mathcal{M}$, where $C = C(n,m,q,M_0,N_0,C_e,C_h) \in (0,\infty)$ is a constant.  Therefore, by Lemma~\ref{app lemma} together with $\|A_{\mathcal M}\|_{C^0(\mathcal{M})} \leq C \mathbf{m}_0^{1/2}$, \eqref{relharm separation1 small norm}, \eqref{relharm critrad1 eqn4}, and \eqref{relharm critrad1 eqn5} 
\begin{align}
	\label{relharm critrad1 eqn6} \|D\xi_j - I_P\|_{C^0(B_{3\rho/4}(Z,P))} 
		\leq\,& C \sup_{X \in B^{\mathcal M}_{\rho}(Z)} \|\pi_{{\rm T}_X \mathcal{M}} - \pi_P\| 
			+ C \|\pi_{{\rm T}\mathcal{M}}(D\varphi_j)\|_{C^0(B^{\mathcal M}_{\rho}(Z))} 
		\\ \leq\,& C \mathbf{m}_0^{1/2} \rho , \nonumber\\
	\label{relharm critrad1 eqn7} \|D\psi_j\|_{C^0(B_{3\rho/4}(Z,P))} 
		\leq\,& C \sup_{X \in B^{\mathcal M}_{\rho}(Z)} \|\pi_{{\rm T}_X \mathcal{M}} - \pi_P\| 
			+ C \|D\varphi_j\|_{C^0(B^{\mathcal M}_{\rho}(Z))}
		\\ \leq\,& C\mathbf{m}_0^{1/2} \rho + C \|D\varphi_j\|_{C^0(B^{\mathcal M}_{\rho}(Z))} , \nonumber\\
	\label{relharm critrad1 eqn9} \|\Delta\psi_j\|_{C^0(B_{3\rho/4}(Z,P))} 
		\leq\,& C \|H_{\mathcal M}\|_{C^0(B^{\mathcal M}_{\rho}(Z))} + C \|\Delta_{\mathcal M}\varphi_j\|_{C^0(B^{\mathcal M}_{\rho}(Z))} 
			\\&+ C \mathbf{m}_0^{1/2} \rho \,(\|A_{\mathcal M}\|_{C^0(B^{\mathcal M}_{\rho}(Z))} 
			+ \|{\rm Hess}\,\varphi_j\|_{C^0(B^{\mathcal M}_{\rho}(Z))}) \nonumber 
		\\ \leq\,& C \mathbf{m}_0^{1/2} + C \|\Delta_{\mathcal M}\varphi_j\|_{C^0(B^{\mathcal M}_{\rho}(Z))} , \nonumber
\end{align}
where $I_P$ is the identity map on $P$ and $C = C(n,m,q,M_0,N_0,C_e,C_h) \in (0,\infty)$ is a constant.

Let $\overrightarrow{\mathcal M}$ denote the orientation on $\mathcal{M}$ and let $\vec P = \overrightarrow{\mathcal M}(Z)$ be the orientation on $P$ induced by $\mathcal{M}$.  Let 
\begin{equation*}
	\widehat{\mathcal E}(T,P,\mathbf{C}_{3\rho/4}(Z,P))^2 
		= \frac{1}{2\omega_n(3\rho/4)^n} \int_{\mathbf{C}_{3\rho/4}(Z,P)} |\vec T - \vec P|^2 \,d\|T\|
\end{equation*}
(as in (2.4) and (2.5) in~\cite{KrumWicc}).  By the triangle inequality,~\cite[Proposition~3.4]{DeLSpa0}, $\|\bs{\varphi}\|_{C^2} \leq C\mathbf{m}_0^{1/2}$ and \eqref{normal approx est1} 
\begin{align}\label{relharm critrad1 eqn10} 
	&\widehat{\mathcal E}(T,P,\mathbf{C}_{3\rho/4}(Z,P))^2 
	\\ \leq\,& \frac{1}{\omega_n (3\rho/4)^n} \int_{\mathbf{C}_{3\rho/4}(Z,P)} |\vec T - \overrightarrow{\mathcal M}|^2 \,d\|T\| 
		+ \frac{1}{\omega_n(3\rho/4)^n} \int_{\mathbf{C}_{3\rho/4}(Z,P)} |\overrightarrow{\mathcal M} - \vec P|^2 \,d\|T\| \nonumber
	\\ \leq\,& \frac{1}{\omega_n(3\rho/4)^n} \int_{B^{\mathcal M}_{3\rho/4}(Z)} |DN|^2 
		+ \frac{C}{\rho^n} \int_{B^{\mathcal M}_{3\rho/4}(Z)} (|A_{\mathcal M}|^2 |N|^2 + |DN|^4) 
		+ C\mathbf{m}_0 \rho^2 \nonumber
	\\ \leq\,& \frac{C}{\rho^n} \int_{B^{\mathcal M}_{3\rho/4}(Z)} |DN|^2 + C\mathbf{m}_0 \rho^2 , \nonumber 
\end{align}
where $C = C(n,m,q,M_0,N_0,C_e,C_h) \in (0,\infty)$ are constants.  

Fix $j \in \{1,2,\ldots,p\}$.  By \eqref{relharm separation1 hyp eqn1}, we have that $(\partial T_j) \llcorner \mathbf{C}^{\mathcal M}_{\rho}(Z) = 0$.  By $\|\bs{\varphi}\|_{C^3} \leq C\mathbf{m}_0^{1/(2n)}$ and~\cite[Corollary~2.2]{DeLSpa0}, the uniform height of $T_j \llcorner \mathbf{C}_{3\rho/4}(Z,P)$ relative to $Z+P$ is $\leq C\mathbf{m}_0^{1/(2n)} \rho^{1+\beta_2} < C \varepsilon^{\beta_2} \rho$ and moreover $\op{spt} T_j \cap \mathbf{C}_{3\rho/4}(Z,P) \subset \mathbf{C}^{\mathcal M}_{3\rho/4 + C\mathbf{m}_0^{1/(2n)} \rho^{1+\beta}}(Z)$, where $C = C(n,m,q,M_0,N_0,C_e,C_h) \in (0,\infty)$ are constants.  By the constancy theorem, $(\pi_{P\#}T_j) \llcorner \mathbf{C}_{3\rho/4}(Z,P)$ is a constant multiple of $\llbracket P \rrbracket \llcorner \mathbf{C}_{3\rho/4}(Z,P)$.  Since $\mathbf{m}_0 \rho < \varepsilon$ and $\|D\bs{\varphi}\|_{C^2} \leq C\mathbf{m}_0^{1/2}$ and $\mathcal{M}$ is tangent to $P$ at $Z$, $\eta_{Z,\rho\#} \mathcal{M}$ is $C\varepsilon^{1/2}$-close to $P$ in the $C^3$-topology in $\mathbf{C}_{3/4}(0,P)$.  Thus by \eqref{relharm separation1 hyp eqn1} and the continuity of push-forwards in the weak topology 
\begin{gather}\label{relharm critrad1 eqn11} 
	(\pi_{P\#} T_j) \llcorner \mathbf{C}_{3\rho/4}(Z,P) = q_j \llbracket P \rrbracket \llcorner \mathbf{C}_{3\rho/4}(Z,P) .
\end{gather}
By \eqref{relharm critrad1 eqn10} and \eqref{relharm normal approx concl2}, $\widehat{\mathcal E}(T,P,\mathbf{C}_{3\rho/4}(Z,P)) \leq C \mathbf{m}_0^{1/2} \rho^{1-\delta_2}$.  Having shown $(\partial T_j) \llcorner \mathbf{C}^{\mathcal M}_{\rho}(Z) = 0$, \eqref{relharm critrad1 eqn11}, and $\widehat{\mathcal E}(T,P,\mathbf{C}_{3\rho/4}(Z,P)) \leq C \mathbf{m}_0 \rho^{1-\delta_2}$, we can apply~\cite[Theorem~4.19]{KrumWicc} to obtain 
\begin{align}\label{relharm critrad1 eqn12}
	&\sup_{X \in \op{spt} T_j \cap \mathbf{C}_{5\rho/8}(Z,P)} \op{dist}^2(X, \op{spt} \mathbf{G}_{\psi_j+a_j}) 
	\leq \frac{C}{\rho^n} \int_{\mathbf{C}_{3\rho/4}(Z,P)} \op{dist}^2(X, \op{spt} \mathbf{G}_{\psi_j+a_j}) \,d\|T_j\|(X) 
		\\&\hspace{15mm} + C \rho^2 \|D\psi_j\|_{C^0(B_{3\rho/4}(Z,P))} 
		\,\widehat{\mathcal E}(T,P,\mathbf{C}_{3\rho/4}(Z,P))^2 \nonumber
		\\&\hspace{15mm} + C\rho^4 \|\Delta\psi_j\|_{C^0(B_{3\rho/4}(Z,P))}^2 
		+ C \rho^2 \widehat{\mathcal E}(T,P,\mathbf{C}_{3\rho/4}(Z,P))^{2+\overline{\varsigma}} , \nonumber
\end{align}
where $\overline{\varsigma} = \overline{\varsigma}(n,m,q) \in (0,1)$ and $C = C(n,m,q) \in (0,\infty)$ are constants.  By \eqref{relharm critrad1 eqn7}, \eqref{relharm critrad1 eqn10}, and Lemma~\ref{normal approx energy lemma} 
\begin{align}\label{relharm critrad1 eqn13}
	& \|D\psi_j\|_{C^0(B_{3\rho/4}(Z,P))} \,\widehat{\mathcal E}(T,P,\mathbf{C}_{3\rho/4}(Z,P))^2 
	\\ \leq\,& C \Big( \mathbf{m}_0^{1/2} \rho + \|D\varphi_j\|_{C^0(B^{\mathcal M}_{\rho}(Z))} \Big)
			\cdot \frac{1}{\rho^n} \int_{B^{\mathcal M}_{3\rho/4}(Z)} |DN|^2 + C\mathbf{m}_0 \rho^2 \nonumber 
	\\ \leq\,& \frac{C \mathbf{m}_0^{1/2}}{\rho^{n+1}} \int_{B^{\mathcal M}_{\rho}(Z)} |N|^2 
		+ C \|D\varphi_j\|_{C^0(B^{\mathcal M}_{\rho}(Z))} \cdot \frac{1}{\rho^n} \int_{B^{\mathcal M}_{3\rho/4}(Z)} |DN|^2 
		+ C\mathbf{m}_0 \rho^2 \nonumber 
\end{align}
where $C = C(n,m,q,M_0,N_0,C_e,C_h) \in (0,\infty)$ are constants.   By \eqref{relharm critrad1 eqn10}, \eqref{normal approx est1}, and Lemma~\ref{normal approx energy lemma} 
\begin{align}\label{relharm critrad1 eqn14}
	\widehat{\mathcal E}(T,P,\mathbf{C}_{3\rho/4}(Z,P))^{2+\overline{\varsigma}}
	\leq\,& C \left( \frac{1}{\rho^n} \int_{B^{\mathcal M}_{3\rho/4}(Z)} |DN|^2 \right)^{1+\overline{\varsigma}/2} 
		+ C\mathbf{m}_0^{1+\overline{\varsigma}/2} \rho^{2+\overline{\varsigma}}  
	\\ \leq\,& C (\mathbf{m}_0 \rho)^{\overline{\varsigma}\gamma_2} \cdot \frac{1}{\rho^n} \int_{B^{\mathcal M}_{3\rho/4}(Z)} |DN|^2 
		+ C\mathbf{m}_0^{1+\overline{\varsigma}/2} \rho^{2+\overline{\varsigma}} \nonumber 
	\\ \leq\,& \frac{C \mathbf{m}_0^{\overline{\varsigma}\gamma_2}}{\rho^{n+2-\overline{\varsigma}\gamma_2}} \int_{B^{\mathcal M}_{\rho}(Z)} |N|^2
		+ C\mathbf{m}_0^{1+\overline{\varsigma}/2} \rho^{2+\overline{\varsigma}} . \nonumber 
\end{align}
Assume that $\varsigma \leq \overline{\varsigma} \gamma_2$.  Recalling \eqref{relharm critrad1 eqn2} and bounding terms on the right-hand side of \eqref{relharm critrad1 eqn12} using \eqref{relharm critrad1 eqn9}, \eqref{relharm critrad1 eqn13}, and \eqref{relharm critrad1 eqn14}, we deduce that \eqref{relharm critrad1 concl} holds true.
\end{proof}

\begin{theorem}\label{relharm separation1 thm}
Let $M_0,N_0,C_e,C_h,\varepsilon_2$ be as in Hypothesis~\ref{DeLSpa constants hypothesis}, Lemma~\ref{normal approx energy lemma} and Lemma~\ref{relharm critrad2 lemma}.  For every integer $p \in \{1,2,\ldots,q\}$ and every $\gamma \in (0,1)$ there exist $\varepsilon = \varepsilon(n,m,q,p,M_0,N_0,C_e,$ $C_h,\gamma) \in (0,1)$ such that the following holds true.  Suppose that Hypothesis~\ref{relharm separation1 hyp} and \eqref{relharm separation1 normal} hold true.  Suppose that 
\begin{equation}\label{relharm separation1 go condition}
	\sup_{X \in B^{\mathcal M}_{\gamma\rho}(Z)} \overline{R}(X) < 5(1-\gamma)\rho/32 ,
\end{equation}
where $\overline{R}$ is as in Definition~\ref{crit rad defn}.  For each $j \in \{1,2,\ldots,p\}$, let $N_j : B^{\mathcal M}_{(1+\gamma)\rho/2}(Z) \rightarrow \mathcal{A}_{q_j}({\rm T}\mathcal{M}^{\perp})$ satisfy \eqref{relharm normal approx concl1} and \eqref{relharm normal approx concl2} with $(1+\gamma)/2$ in place of $\gamma$.  Then for each $j \in \{1,2,\ldots,p\}$ there exists $a_j \in \mathcal{A}_{q_j}(\mathbb{R}^{n+m})$ such that 
\begin{equation}\label{relharm separation1 concl}
	\sup_{X \in \op{spt} T_j \cap \mathbf{C}^{\mathcal M}_{\gamma\rho}(Z)} \op{dist}(X, \op{spt} \mathbf{G}_{\varphi_j+a_j}) 
	\leq C \mathcal{E}_{Z,(1+\gamma)\rho/2,\rho} \,\rho ,
\end{equation}
where $\mathcal{E}_{Z,(1+\gamma)\rho/2,\rho}$ is as in \eqref{relharm excess defn1} with $\varsigma = \varsigma(n,m,q) \in (0,1)$ (and with $(1+\gamma)/2$ in place of $\gamma$) and where $C = C(n,m,q,M_0,N_0,C_e,C_h,\gamma) \in (0,\infty)$ is a constant.
\end{theorem}

\begin{proof}
The proof is similar to that of~\cite[Theorem~2.12]{KrumWicb} and we provide details for completion.  Let $\varepsilon$ be the minimum of its values in Lemma~\ref{relharm tilt growth lemma}, Lemma~\ref{relharm poincare lemma}, Lemma~\ref{relharm critrad1 lemma}, and Lemma~\ref{relharm critrad2 lemma}.  Choose $\theta = \theta(n,m,q) \in (0,1/16)$ so that $C\theta^{\mu/2} < 1$, where $\mu$ and $C$ are as in Lemma~\ref{relharm tilt growth lemma}.  Let $\varsigma = \varsigma(n,m,q) \in (0,1)$ be the minimum of $2\delta_2$ and its value in Lemma~\ref{relharm critrad1 lemma}.  Let $\mathcal{K}_j \subseteq B^{\mathcal M}_{(1+\gamma)\rho/2}(Z)$ and $N_j : B^{\mathcal M}_{(1+\gamma)\rho/2}(Z) \rightarrow \mathcal{A}_{q_j}({\rm T}\mathcal{M}^{\perp})$ satisfy \eqref{relharm normal approx concl1} and \eqref{relharm normal approx concl2} with $(1+\gamma)/2$ in place of $\gamma$.  By Lemma~\ref{relharm poincare lemma} for each $j \in \{1,2,\ldots,p\}$ there exists $a_j \in \mathcal{A}_{q_j}(\mathbb{R}^{n+m})$ such that $|a_j| < C_1 \varepsilon^{\beta_2} \rho$
\begin{align}\label{relharm separation1 eqn1}
	\frac{1}{\rho^{n+2}} \int_{\mathbf{C}^{\mathcal M}_{(1+\gamma)\rho/2}(Z)} \op{dist}^2(X, \op{spt} \mathbf{G}_{\varphi_j+a_j}) \,d\|T_j\|(X)&
	\\+ \frac{1}{\rho^{n+2}} \int_{B^{\mathcal M}_{(1+\gamma)\rho/2}(Z)} \mathcal{G}(N_j,\varphi_j+a_j)^2 
	&\leq C_1 \mathcal{E}_{Z,(1+\gamma)\rho/2,\rho}^2 , \nonumber 
\end{align}
where $C_1 = C_1(n,m,q,M_0,N_0,C_e,C_h,\gamma) \in (0,\infty)$ is a constant.  

Let $Z' \in B^{\mathcal M}_{\gamma\rho}(Z)$.  By \eqref{relharm separation1 go condition}, $\overline{R}(Z') < 5(1-\gamma)\rho/32$.  If $Z'$ is not in the contact set $\{(x,\bs{\varphi}(x)) : x \in \mathbf{\Gamma} \}$ so that $\overline{R}(Z') > 0$, let $i_0 \geq 0$ such that 
\begin{equation}\label{relharm separation1 eqn2}
	5\theta^{i_0+1} (1-\gamma)\rho/32 < \overline{R}(Z') \leq 5\theta^{i_0} (1-\gamma)\rho/32 . 
\end{equation}
If $Z'$ is in the contact set, let $i_0 = \infty$.  Recalling Remark~\ref{relharm norm scaling rmk}, for each integer $i$ with $1 \leq i \leq i_0$ let $\mathcal{E}_{Z',5\theta^i (1-\gamma)\rho/32,\theta^i (1-\gamma)\rho/4}$ be as in \eqref{relharm excess defn1} with $5/8,Z',\theta^i (1-\gamma)\rho/4$ in place of $\gamma,Z,\rho$ and with the normal approximation maps $N_j$.  By inductively applying Lemma~\ref{relharm tilt growth lemma} 
\begin{equation}\label{relharm separation1 eqn3}
	\mathcal{E}_{Z',5\theta^i (1-\gamma)\rho/32,\theta^i (1-\gamma)\rho/4} \leq C \theta^{i(\mu/2-1)} \mathcal{E}_{Z,(1+\gamma)\rho/2,\rho} 
\end{equation}
for all integers $i$ with $1 \leq i \leq i_0$, where $C = C(n,\gamma) \in (0,\infty)$ is a constant.  By Lemma~\ref{relharm poincare lemma} and \eqref{relharm separation1 eqn3}, for each integer $i$ with $1 \leq i \leq i_0$ and each $j \in \{1,2,\ldots,p\}$ there exists $a^{(i)}_j \in \mathcal{A}_{q_j}(\mathbb{R}^{n+m})$ such that $|a^{(i)}_j| < C_1 \varepsilon^{\beta_2} \theta^i \rho$ and 
\begin{align}\label{relharm separation1 eqn4}
	&\frac{1}{(\theta^i \rho)^{n+2}} \int_{\mathbf{C}^{\mathcal M}_{5\theta^i (1-\gamma)\rho/32}(Z')} \op{dist}^2(X, \op{spt} \mathbf{G}_{\varphi_j+a^{(i)}_j}) 
	\,d\|T_j\|(X) \\&+ \frac{1}{(\theta^i \rho)^{n+2}} \int_{B^{\mathcal M}_{5\theta^i (1-\gamma)\rho/32}(Z')} \mathcal{G}(N_j,\varphi_j+a^{(i)}_j)^2 
	\nonumber\\ \leq\,& C \mathcal{E}_{Z',5\theta^i (1-\gamma)\rho/32,\theta^i (1-\gamma)\rho/2}^2 
	\leq C \theta^{2i (\mu/2-1)} \mathcal{E}_{Z,(1+\gamma)\rho/2,\rho}^2 , \nonumber 
\nonumber 
\end{align}
where $C = C(n,m,q,M_0,N_0,C_e,C_h,\gamma) \in (0,\infty)$ are constants.  Let $a^{(0)}_j = a_j$, where $a_j$ is as in \eqref{relharm separation1 eqn1}.  By the triangle inequality, \eqref{relharm separation1 eqn1}, and \eqref{relharm separation1 eqn4} 
\begin{equation}\label{relharm separation1 eqn5}
	\mathcal{G}(a^{(i)}_j,a^{(i+1)}_j) \leq C \theta^{i\mu/2} \mathcal{E}_{Z,(1+\gamma)\rho/2,\rho} \,\rho
\end{equation}
for all integers $i$ with $0 \leq i < i_0$ and all $j \in \{1,2,\ldots,p\}$, where $C = C(n,m,q,M_0,N_0,C_e,C_h,\gamma) \in (0,\infty)$ is a constant.  Hence by the triangle inequality 
\begin{equation}\label{relharm separation1 eqn6}
	\mathcal{G}(a^{(i)}_j,a_j) \leq C \mathcal{E}_{Z,(1+\gamma)\rho/2,\rho} \,\rho
\end{equation}
for all integers $i$ with $0 \leq i \leq i_0$ and $j \in \{1,2,\ldots,p\}$, where $C = C(n,m,q,M_0,N_0,C_e,C_h,\gamma) \in (0,\infty)$ is a constant.  

Suppose that $Z'$ is in the contact set and recall that $\op{spt} T_j \cap \mathbf{p}^{-1}(Z') = \{Z'\}$.  Fix $j \in \{1,2,\ldots,p\}$.  By \eqref{relharm separation1 eqn5}, $(a^{(i)}_j)$ is a Cauchy sequence in $\mathcal{A}_q(\mathbb{R}^{n+m})$ and thus there exists $\overline{a}_j \in \mathcal{A}_{q_j}(\mathbb{R}^{n+m})$ such that $a^{(i)}_j \rightarrow \overline{a}_j$ in $\mathcal{A}_{q_j}(\mathbb{R}^{n+m})$ as $i \rightarrow \infty$.  By \eqref{relharm separation1 eqn5}, $\mathcal{G}(a^{(i)}_j,\overline{a}_j) \leq C \theta^{i\mu/2} \mathcal{E}_{Z,(1+\gamma)\rho/2,\rho} \,\rho$ for each integer $i \geq 0$, and in particular $\mathcal{G}(a_j,\overline{a}_j) \leq C \mathcal{E}_{Z,(1+\gamma)\rho/2,\rho} \,\rho$, where $C = C(n,m,q,M_0,N_0,C_e,C_h,\gamma) \in (0,\infty)$ is a constant.  Hence by the triangle inequality and \eqref{relharm separation1 eqn4} 
\begin{equation*}
	\frac{1}{(\theta^i \rho)^{n+2}} \int_{B^{\mathcal M}_{5\theta^i (1-\gamma)\rho/32}(Z')} \mathcal{G}(N_j,\varphi_j+\overline{a}_j)^2 
	\leq C \theta^{2i (\mu/2-1)} \mathcal{E}_{Z,(1+\gamma)\rho/2,\rho}^2 
\end{equation*}
for all integers $i \geq 1$, where $C = C(n,m,q,M_0,N_0,C_e,C_h,\gamma) \in (0,\infty)$ is a constant.  Letting $i\rightarrow\infty$, we deduce that $N_j(Z') = q_j \llbracket 0 \rrbracket = \varphi_j(Z') + \overline{a}_j$ so that $Z' \in \op{spt}\mathbf{G}_{\varphi_j+\overline{a}_j}$.  Therefore, by $\mathcal{G}(a_j,\overline{a}_j) \leq C \mathcal{E}_{Z,(1+\gamma)\rho/2,\rho} \,\rho$ 
\begin{equation}\label{relharm separation1 eqn7}
	\op{dist}(Z', \op{spt} \mathbf{G}_{\varphi_j+a_j}) \leq C \mathcal{E}_{Z,(1+\gamma)\rho/2,\rho} \,\rho , 
\end{equation}
where $C = C(n,m,q,M_0,N_0,C_e,C_h,\gamma) \in (0,\infty)$ is a constant.

Now suppose that $Z'$ is not in the contact set so that $i_0 < \infty$.  For each $j \in \{1,2,\ldots,p\}$, let $\widetilde{a}_j = \sum_{l=1}^{q_j} \llbracket \pi_{{\rm T}_{Z'} \mathcal{M}^{\perp}}(a^{(i_0)}_{j,l}) \rrbracket$, where $a^{(i_0)}_j = \sum_{l=1}^{q_j} \llbracket a^{(i_0)}_{j,l} \rrbracket$.  Since $N^{(i_0)}_j$ and $\varphi_j$ take values orthogonal to $\mathcal{M}$  
\begin{align}\label{relharm separation1 eqn8}
	&\frac{1}{(\theta^{i_0} \rho)^{n+2}} \int_{B^{\mathcal M}_{5\theta^{i_0}(1-\gamma)\rho/32}(Z')} \sum_{l=1}^{q_j} 
		|\pi_{{\rm T}_X \mathcal{M}}(a^{(i_0)}_{j,l})|^2 
	\\ \leq\,& \frac{1}{(\theta^{i_0} \rho)^{n+2}} \int_{B^{\mathcal M}_{5\theta^{i_0}(1-\gamma)\rho/32}(Z')} \mathcal{G}(N^{(i_0)}_j, \varphi_j + a^{(i_0)}_j)^2 
	\leq \mathcal{E}_{Z',5\theta^{i_0}(1-\gamma)\rho/32,\theta^{i_0} (1-\gamma)\rho/4}^2 . \nonumber 
\end{align}
Thus using $\|\pi_{{\rm T}_X \mathcal{M}} - \pi_{{\rm T}_{Z'} \mathcal{M}}\| \leq C \mathbf{m}_0^{1/2} |X - Z'|$, $|a^{(i_0)}_j| < \theta^{i_0} \rho$, and \eqref{relharm separation1 eqn8}
\begin{align}\label{relharm separation1 eqn9}
	\mathcal{G}(a^{(i_0)}_j,\widetilde{a}_j)^2 
	\leq\,& \sum_{l=1}^{q_j} |\pi_{{\rm T}_{Z'} \mathcal{M}}(a^{(i_0)}_{j,l})|^2 
	\\ \leq\,& \frac{C}{(\theta^{i_0} \rho)^n} \int_{B^{\mathcal M}_{5\theta^{i_0}(1-\gamma)\rho/32}(Z')} \sum_{l=1}^{q_j} |\pi_{{\rm T}_X \mathcal{M}}(a^{(i_0)}_{j,l})|^2 + C \mathbf{m}_0 \theta^{2i_0} \rho^2 |a^{(i_0)}_j|^2 \nonumber 
	\\ \leq\,& C \mathcal{E}_{Z',5\theta^{i_0}(1-\gamma)\rho/32,\theta^{i_0} (1-\gamma)\rho/4}^2 \theta^{2i_0} \rho^2 + C \mathbf{m}_0 \theta^{4i_0} \rho^4 , \nonumber
\end{align}
where $C = C(n,m,q,M_0,N_0,C_e,C_h,\gamma) \in (0,\infty)$ are constants.  By \eqref{relharm separation1 eqn2},  Lemma~\ref{relharm critrad2 lemma}, and \eqref{relharm separation1 eqn3}
\begin{equation}\label{relharm separation1 eqn10}
    \mathbf{m}_0 \theta^{4i_0} \rho^4 \leq C \mathcal{E}_{Z',5\theta^{i_0}(1-\gamma)\rho/32,\theta^{i_0} (1-\gamma)\rho/4}^2 \theta^{2i_0} \rho^2 \leq C \theta^{i_0 \mu} \mathcal{E}_{Z,(1+\gamma)\rho/2,\rho}^2 \rho^2
\end{equation}
for some constant $C = C(n,m,q,M_0,N_0,C_e,C_h,\gamma) \in (0,\infty)$.  Thus by \eqref{relharm separation1 eqn9}, \eqref{relharm separation1 eqn10}, and \eqref{relharm separation1 eqn3}
\begin{align}\label{relharm separation1 eqn11}
	\mathcal{G}(a^{(i_0)}_j,\widetilde{a}_j)^2 
	\leq C \mathcal{E}_{Z',5\theta^{i_0}(1-\gamma)\rho/32,\theta^{i_0} (1-\gamma)\rho/4}^2 \theta^{2i_0} \rho^2 
	\leq C \theta^{i_0 \mu} \mathcal{E}_{Z,(1+\gamma)\rho/2,\rho}^2 \rho^2 , 
\end{align}
where $C = C(n,m,q,M_0,N_0,C_e,C_h,\gamma) \in (0,\infty)$ are constants.  By \eqref{relharm separation1 eqn4} and \eqref{relharm separation1 eqn11} 
\begin{align}\label{relharm separation1 eqn12}
	&\frac{1}{(\theta^{i_0} \rho)^{n+2}} \int_{\mathbf{C}^{\mathcal M}_{5\theta^{i_0}(1-\gamma)\rho/32}(Z')} 
		\op{dist}^2(X, \op{spt} \mathbf{G}_{\varphi_j+\widetilde{a}_j}) \,d\|T_j\|(X) 
		\\&+ \frac{1}{(\theta^{i_0} \rho)^{n+2}} \int_{B^{\mathcal M}_{5\theta^{i_0} (1-\gamma)\rho/32}(Z')} 
		\mathcal{G}(N^{(i_0)}_j,\varphi_j+\widetilde{a}_j)^2 
	\leq C \theta^{2i_0 (\mu/2-1)} \mathcal{E}_{Z,(1+\gamma)\rho/2,\rho}^2 , \nonumber 
\end{align}
where $C = C(n,m,q,M_0,N_0,C_e,C_h,\gamma) \in (0,\infty)$ is a constant.  By Lemma~\ref{relharm critrad1 lemma}, \eqref{relharm separation1 eqn12}, \eqref{relharm separation1 eqn10}, and \eqref{relharm separation1 eqn3}
\begin{equation}\label{relharm separation1 eqn13}
	\sup_{X \in \op{spt} T_j \cap \mathbf{C}^{\mathcal M}_{\overline{R}(Z')/2}(Z')} \op{dist}(X, \op{spt} \mathbf{G}_{\varphi_j+\widetilde{a}_j}) 
	\leq C \theta^{i_0 \mu/2} \mathcal{E}_{Z,(1+\gamma)\rho/2,\rho} \,\rho , 
\end{equation}
where $C = C(n,m,q,M_0,N_0,C_e,C_h) \in (0,\infty)$ are constants.  By \eqref{relharm separation1 eqn6} and \eqref{relharm separation1 eqn11}, $\mathcal{G}(a_j,\widetilde{a}_j) \leq C \mathcal{E}_{Z,(1+\gamma)\rho/2,\rho} \,\rho$, and thus \eqref{relharm separation1 eqn13} gives us 
\begin{equation}\label{relharm separation1 eqn14}
	\sup_{X \in \op{spt} T_j \cap \mathbf{C}^{\mathcal M}_{\overline{R}(Z')/2}(Z')} \op{dist}(X, \op{spt} \mathbf{G}_{\varphi_j+a_j}) 
	\leq C \mathcal{E}_{Z,(1+\gamma)\rho/2,\rho} \,\rho , 
\end{equation}
where $C = C(n,m,q,M_0,N_0,C_e,C_h,\gamma) \in (0,\infty)$ is a constant.  Since either \eqref{relharm separation1 eqn7} or \eqref{relharm separation1 eqn14} holds true, we conclude that \eqref{relharm separation1 concl} holds true.
\end{proof}

\begin{remark}\label{relharm separation1 rmk}{\rm 
Notice that if $\mathcal{E}_{Z,(1+\gamma)\rho/2,\rho} = 0$, then $\op{spt} T_j \subseteq \op{spt} \mathbf{G}_{\varphi_j + a_j}$.  Assume that $\mathcal{E}_{Z,(1+\gamma)\rho/2,\rho} > 0$.  Recall from the proof of Theorem~\ref{relharm separation1 rmk} that we chose $a_j \in \mathcal{A}_q(\mathbb{R}^{n+m})$ so that \eqref{relharm separation1 eqn1} holds true, and we proved that \eqref{relharm separation1 concl} holds true.  Fix $j \in \{1,2,\ldots,p\}$ and let  
\begin{equation*}
	\{ X \in \mathbf{C}^{\mathcal M}_{\gamma\rho}(Z) : \op{dist}(X, \op{spt}\op{graph}(\varphi_j+a_j)) < 2C \mathcal{E}_{Z,(1+\gamma)\rho/2,\rho} \,\rho \} 
	= \bigcup_{i=1}^{M_j} U_{j,i}
\end{equation*}
where $\{U_{j,i}\}$ is a collection of mutually disjoint, connected, open subsets of $\mathbf{C}^{\mathcal M}_{\gamma\rho}(Z)$ and $C$ is as in \eqref{relharm separation1 concl}.  Then 
\begin{equation*}
	T_j \llcorner \mathbf{C}_{\gamma}(0) = \sum_{i=1}^{M_j} T_{j,i}  \quad\text{where}\quad T_{j,i} = T \llcorner U_{j,i} 
\end{equation*}
and $T_{j,i}$ are locally area-minimizing rectifiable currents of $\mathbf{C}_{\gamma}(0)$ with $(\partial T_i) \llcorner \mathbf{C}^{\mathcal M}_{\gamma\rho}(Z) = 0$.  By the constancy theorem 
\begin{equation*}
	\mathbf{p}_{\#} T_{j,i} = q_{j,i} \llbracket B^{\mathcal M}_{\gamma\rho}(Z) \rrbracket 
\end{equation*}
for some integers $q_{j,i}$ with $\sum_{i=1}^{M_j} q_{j,i} = q_j$. Arguing as we did in the proof of Lemma~\ref{relharm critrad1 lemma} to justify the hypotheses of~\cite[Theorem~4.19]{KrumWicc} 
\begin{gather*}
	\sup_{X \in \op{spt} T_{j,i} \cap \mathbf{C}^{\mathcal M}_{\gamma\rho}(Z)} \op{dist}(X,\mathcal{M}) 
		\leq C \mathbf{m}_0^{\frac{1}{2n}} \rho^{1+\beta_2} , \\ 
	\pi_{P\#}(T_{j,i} \llcorner \mathbf{C}_{3\gamma\rho/4}(Z,{\rm T}_Z \mathcal{M})) 
		= q_{j,i} \llbracket B_{3\gamma\rho/4}(Z,{\rm T}_Z \mathcal{M})) \rrbracket , \\
	\mathcal{E}(T_{j,l},{\rm T}_Z \mathcal{M},\mathbf{C}_{3\gamma\rho/4}(Z,{\rm T}_Z \mathcal{M})) \leq C\mathbf{m}_0^{1/2} \rho^{1-\delta_2}, 
\end{gather*}
where $C = C(n,m,q,M_0,N_0,C_e,C_h) \in (0,\infty)$ is a constant.  Assuming $\varepsilon_2$ is sufficiently small, by~\cite[Lemma~2.7]{KrumWicb} we have that $q_{j,i} \geq 0$ for all $i \in \{1,2,\ldots,M_j\}$.  By~\cite[Lemma~2.6]{KrumWicb}, $q_{j,i} = 0$ if and only if $\op{spt} T_{j,i} \cap \mathbf{C}^{\mathcal M}_{\gamma\rho/2}(Z) = \emptyset$.  We claim that by arguing as in~\cite[Remark~2.18]{KrumWicb}, we can assume that $q_{j,i} > 0$ for all $i \in \{1,2,\ldots,M_j\}$.  Observe that by \eqref{relharm separation1 eqn1} and \eqref{normal approx est1}, there exists a set $\Omega_j \subset \mathcal{K}_j \cap B^{\mathcal M}_{\gamma\rho/2}(Z)$ such that $\mathcal{H}^n(\Omega_j) \geq \frac{1}{2} \omega_n (\gamma\rho/2)^n$ and $\mathcal{G}(N_j(X),\varphi_j(X)+a_j) \leq C_2 \mathcal{E}_{Z,(1+\gamma)\rho/2,\rho} \,\rho$ for all $X \in \Omega_j$, where $C_2 = \sqrt{\frac{4C_1}{\omega_n(\gamma/2)^n}}$ (for $C_1$ as in \eqref{relharm separation1 eqn1}).  Take any $i \in \{1,2,\ldots,M_j\}$ and $X \in \Omega$.  Let $a_j = \sum_{l=1}^{q_j} \llbracket a_{j,l} \rrbracket$ and $N_j(X) = \sum_{l=1}^{q_j} \llbracket N_{j,l}(X) \rrbracket$, where $a_{j,l} \in \mathbb{R}^{n+m}$ and $N_{j,l}(X) \in {\rm T}_X \mathcal{M}^{\perp}$.  There exists $k(i) \in \{1,2,\ldots,q_j\}$ such that $\op{graph}(\varphi_j+a_{j,k(i)}) \cap \mathbf{C}_{\gamma}(0) \subset U_{j,i}$, and there exists $l(i) \in \{1,2,\ldots,q_j\}$ such that $X + N_{j,l(i)}(X) \in \op{spt} T_j$ and 
\begin{align*}
	\op{dist}(X + N_{j,l(i)}(X),\op{graph}(\varphi_j+a_{j,k(i)})) \leq\,& |N_{j,l(i)}(X) - \varphi_j(X) - a_{j,k(i)}| \\ \leq\,& \mathcal{G}(N_j(X),\varphi_j(X)+a_j) 
		\leq C_2 \mathcal{E}_{Z,(1+\gamma)\rho/2,\rho} \,\rho .
\end{align*}
Hence provided we take $C_1$ in \eqref{relharm separation1 concl} large enough that $C_2 < 2C_1$, $(X,N_{j,l(i)}(X)) \in \op{spt} T_{j,i} \cap \mathbf{C}^{\mathcal M}_{\gamma\rho/2}(Z)$.  Therefore, $q_{j,i} > 0$ for all $i \in \{1,2,\ldots,M_j\}$.
}\end{remark}

\subsection{$L^{\infty}$-distance estimate in terms $L^2$-distance relative to parallel $C^2$-graphs}\label{sec:relharm separation2 sec}  Let $\gamma,\varepsilon,\varsigma \in (0,1)$.  Assume that Hypothesis~\ref{relharm separation1 hyp} holds true.  For each $j \in \{1,2,\ldots,p\}$, let $N_j : B^{\mathcal M}_{(1+\gamma)\rho/2}(Z) \rightarrow \mathcal{A}_{q_j}({\rm T}\mathcal{M}^{\perp})$ satisfy \eqref{relharm normal approx concl1} and \eqref{relharm normal approx concl2}.  For each $j \in \{1,2,\ldots,p\}$, let $a_j \in \mathcal{A}_{q_j}(\mathbb{R}^{n+m})$ with $|a_j| < \varepsilon\rho$.  We define 
\begin{align}\label{relharm excess defn2}
	E_{Z,\gamma\rho,\rho}^2 =\,& \frac{1}{\rho^{n+2}} \sum_{j=1}^p \sum_{l=1}^{q_j} \int_{B^{\mathcal M}_{\gamma\rho}(Z)} 
		\mathcal{G}(N_j, \varphi_j + a_j)^2 + \frac{\rho^2}{\varepsilon^2} \sum_{j=1}^p \|\Delta_{\mathcal M} \varphi_j\|_{C^0(B^{\mathcal M}_{\rho}(Z))}^2  
		\\&+ \left( \sum_{j=1}^p \|D\varphi_j\|_{C^0(B^{\mathcal M}_{\rho}(Z))} + \rho^{\varsigma} \right) 
		\frac{1}{\rho^{n+2}} \int_{B^{\mathcal M}_{\rho}(Z)} |N|^2 \nonumber 
\end{align}
The first term on the right-hand side of \eqref{relharm excess defn2} is an $L^2$-distance of $T_j$ to the graph of $\varphi_j+a_j$.  Much like the case with \eqref{relharm excess defn1}, the remaining terms on the right-hand side of \eqref{relharm excess defn2} control $\Delta_{\mathcal M} \varphi_j$ and error terms arising from the area formula and estimates for the normal approximation maps $N_j$ and $N$.

In this subsection, we want to bound the $L^{\infty}$-distance of $T_j \llcorner \mathbf{C}^{\mathcal M}_{\gamma\rho}(Z)$ to the graph of $\varphi_j + a_j$ in terms of $E_{Z,(1+\gamma)\rho/2,\rho}$.  First, let us consider the special case where $\op{size} a_j = 1$ for all $j$ (where $\op{size} a_j$ is as in Definition~\ref{a sep width defn}).  This is analogous to~\cite[Lemma~3.18]{KrumWicb}; however, we prove the result using the approach of~\cite[Lemma~4.13]{KrumWicc}.

\begin{lemma}\label{relharm separation0 lemma} 
Let $M_0,N_0,C_e,C_h,\varepsilon_2$ be as in Hypothesis~\ref{DeLSpa constants hypothesis}, Lemma~\ref{normal approx energy lemma} and Lemma~\ref{relharm critrad2 lemma}.  For every integer $p \in \{1,2,\ldots,q\}$ and every $\gamma \in (0,1)$ there exist $\varepsilon = \varepsilon(n,m,q,p,M_0,N_0,$ $C_e,C_h,\gamma) \in (0,1)$ such that the following holds true.  Suppose that Hypothesis~\ref{relharm separation1 hyp} and \eqref{relharm separation1 normal} hold hold true.  Suppose that 
\begin{equation}\label{relharm separation0 go condition}
	\sup_{X \in B^{\mathcal M}_{\gamma\rho}(Z)} \overline{R}(X) < 5(1-\gamma)\rho/64 ,
\end{equation}
where $\overline{R}$ is as in Definition~\ref{crit rad defn}.  For each $j \in \{1,2,\ldots,p\}$, let $N_j : B^{\mathcal M}_{(1+\gamma)\rho/2}(Z) \rightarrow \mathcal{A}_{q_j}({\rm T}\mathcal{M}^{\perp})$ satisfy \eqref{relharm normal approx concl1} and \eqref{relharm normal approx concl2} with $(1+\gamma)/2$ in place of $\gamma$.  For each $j \in \{1,2,\ldots,p\}$, let $a_j = q_j \llbracket a_{j,1} \rrbracket \in \mathcal{A}_{q_j}(\mathbb{R}^{n+m})$ such that $|a_j| < \varepsilon\rho$ and $\frac{1}{q} \sum_{j=1}^p q_j a_{j,1} = 0$.  Then for each $j \in \{1,2,\ldots,p\}$  
\begin{equation}\label{relharm separation0 concl}
	\sup_{X \in \op{spt} T_j \cap \mathbf{C}^{\mathcal M}_{\gamma\rho}(Z)} \op{dist}(X, \op{spt} \mathbf{G}_{\varphi_j+a_j}) 
	\leq C E_{Z,(1+\gamma)\rho/2,\rho} \,\rho
\end{equation}
where $E_{Z,(1+\gamma)\rho/2,\rho}$ is as in \eqref{relharm excess defn2} with $\varsigma = \varsigma(n,m,q) \in (0,1)$ (and with $(1+\gamma)/2$ in place of $\gamma$) and where $C = C(n,m,q,M_0,N_0,C_e,C_h,\gamma) \in (0,\infty)$ is a constant.
\end{lemma}

\begin{proof} 
Let $\varsigma$ be the minimum of $\gamma_2/4$ and its value in Theorem~\ref{relharm separation1 thm}.  Note that by \eqref{relharm separation0 go condition}, $(3+\gamma)\rho/4 > \overline{R}(Z)$.  By Lemma~\ref{relharm energy lemma} (with $\frac{(3+\gamma)\rho}{4},\frac{1+3\gamma}{3+\gamma}$ in place of $\rho,\gamma$) 
\begin{equation}\label{relharm separation0 eqn1}
	\frac{1}{\rho^n} \int_{B^{\mathcal M}_{(1+3\gamma)\rho/4}(Z)} \mathcal{G}(DN_j,q_j \llbracket D\varphi_j \rrbracket)^2 
	\leq C E_{Z,(1+\gamma)\rho/2,\rho}^2 
\end{equation}
and thus by Lemma~\ref{normal approx energy lemma}, $\mathcal{E}_{Z,(1+3\gamma)\rho/4,(1+\gamma)\rho/2} \leq C E_{Z,(1+\gamma)\rho/2,\rho}$ where $C = C(n,m,q,M_0,N_0,C_e,C_h,$ $\gamma) \in (0,\infty)$ are constants.  Recall Remark~\ref{relharm norm scaling rmk} and note that by \eqref{relharm separation0 go condition}, \eqref{relharm separation1 go condition} holds true with $\frac{(1+\gamma)\rho}{2}, \frac{2\gamma}{1+\gamma}$ in place of $\rho,\gamma$.  Thus using Lemma~\ref{relharm poincare lemma} and \eqref{relharm separation0 eqn1}, for each $j \in \{1,2,\ldots,p\}$ there exists $b_j \in \mathcal{A}_{q_j}(\mathbb{R}^{n+m})$ such that 
\begin{align}\label{relharm separation0 eqn2}
	\frac{1}{\rho^{n+2}} \int_{\mathbf{C}^{\mathcal M}_{(1+3\gamma)\rho/4}(Z)} \op{dist}^2(X, \op{spt} \mathbf{G}_{\varphi_j+b_j}) \,d\|T_j\|(X)&
	\\+ \frac{1}{\rho^{n+2}} \int_{B^{\mathcal M}_{(1+3\gamma)\rho/4}(Z)} \mathcal{G}(N_j,\varphi_j+b_j)^2 &\leq C E_{Z,(1+\gamma)\rho/2,\rho}^2 , \nonumber 
\end{align}
where $C = C(n,m,q,M_0,N_0,C_e,C_h) \in (0,\infty)$ is a constant.  Then by Theorem~\ref{relharm separation1 thm} (with $\frac{(1+\gamma)\rho}{2}, \frac{2\gamma}{1+\gamma}$ in place of $\rho,\gamma$)
\begin{equation}\label{relharm separation0 eqn3}
	\sup_{X \in \op{spt} T_j \cap \mathbf{C}^{\mathcal M}_{\gamma\rho}(Z)} \op{dist}(X, \op{spt} \mathbf{G}_{\varphi_j+b_j}) 
	\leq C E_{Z,(1+\gamma)\rho/2,\rho} \,\rho , 
\end{equation}
where $C = C(n,m,q,M_0,N_0,C_e,C_h,\gamma) \in (0,\infty)$ is a constant.  
By the triangle inequality and \eqref{relharm separation0 eqn2}  
\begin{align}\label{relharm separation0 eqn4}
	\mathcal{G}(b_j, a_j)^2 
	\leq\,& \frac{C}{\rho^n} \int_{B^{\mathcal M}_{(1+3\gamma)\rho/4}(Z)} \mathcal{G}(N_j,\varphi_j+b_j)^2 
		\\&+ \frac{C}{\rho^n} \int_{B^{\mathcal M}_{(1+3\gamma)\rho/4}(Z)} \mathcal{G}(N_j,\varphi_j+a_j)^2
	\leq C E_{Z,(1+\gamma)\rho/2,\rho}^2 \,\rho^2 , \nonumber 
\end{align}
where $C = C(n,m,q,M_0,N_0,C_e,C_h,\gamma) \in (0,\infty)$ are constants.  Therefore, by the triangle inequality, \eqref{relharm separation0 eqn3}, and \eqref{relharm separation0 eqn4}, we conclude that \eqref{relharm separation0 concl} holds true.
\end{proof}

Next in Lemma~\ref{relharm separation2 lemma}, we want to bound the $L^{\infty}$-distance of $T_j$ to the graph of $\varphi_j + a_j$ in terms of $E_{Z,(1+\gamma)\rho/2,\rho}$ in the setting where $\op{size} a_j > 1$ for some $j$, but with the additional hypotheses \eqref{relharm separation2 hyp1} and \eqref{relharm separation2 hyp2}.  This is analogous to~\cite[Lemma~3.19]{KrumWicb}.  The proof of Lemma~\ref{relharm separation2 lemma} will use and follow after Lemma~\ref{relharm wellsep2 lemma} below. 

\begin{lemma}\label{relharm separation2 lemma} 
Let $M_0,N_0,C_e,C_h,\varepsilon_2$ be as in Hypothesis~\ref{DeLSpa constants hypothesis}, Lemma~\ref{normal approx energy lemma} and Lemma~\ref{relharm critrad2 lemma}.  For all $p \in \{1,2,\ldots,q\}$, $\gamma \in (0,1)$, and $\Lambda \in (0,\infty)$ there exist $\varepsilon, \beta \in (0,1)$ depending only on $n,m,q,p,s,M_0,N_0,C_e,C_h,\gamma,\Lambda$ such that the following holds true.  Suppose that Hypothesis~\ref{relharm separation1 hyp} and \eqref{relharm separation1 normal} hold true.  Suppose that 
\begin{equation}\label{relharm separation2 go condition}
	\sup_{X \in B^{\mathcal M}_{(1+7\gamma)\rho/8}(Z)} \overline{R}(X) < 5(1-\gamma)\rho/512 , 
\end{equation}
where $\overline{R}$ is as in Definition~\ref{crit rad defn}.  For each $j \in \{1,2,\ldots,p\}$, let $N_j : B^{\mathcal M}_{(1+\gamma)\rho/2}(Z) \rightarrow \mathcal{A}_{q_j}({\rm T}\mathcal{M}^{\perp})$ satisfy \eqref{relharm normal approx concl1} and \eqref{relharm normal approx concl2} with $(1+\gamma)/2$ in place of $\gamma$.  For each $j \in \{1,2,\ldots,p\}$, let $a_j = \sum_{l=1}^{q_j} \llbracket a_{j,l} \rrbracket \in \mathcal{A}_{q_j}(\mathbb{R}^{n+m})$ such that $|a_j| < \varepsilon\rho$ and $\frac{1}{q} \sum_{j=1}^p q_j a_{j,a} = 0$, where $a_{j,a} = \frac{1}{q_j} \sum_{l=1}^{q_j} a_{j,l}$.  Suppose that $\op{size} a_j > 1$ for some $j \in \{1,2,\ldots,p\}$ and 
\begin{align}
	\label{relharm separation2 hyp1} \max_{1 \leq j \leq p} \op{width} a_j &\leq \Lambda \min_{1 \leq j \leq p} \op{sep} a_j , \\
	\label{relharm separation2 hyp2} E_{Z,(1+\gamma)\rho/2,\rho} \,\rho &\leq \beta \min_{1 \leq j \leq p} \op{sep} a_j , 
\end{align}
Then for each $j \in \{1,2,\ldots,p\}$  
\begin{equation}\label{relharm separation2 concl}
	\sup_{X \in \op{spt} T_j \cap \mathbf{C}^{\mathcal M}_{\gamma\rho}(Z)} \op{dist}(X, \op{spt} \mathbf{G}_{\varphi_j+a_j}) 
	\leq C E_{Z,(1+\gamma)\rho/2,\rho} \,\rho 
\end{equation}
for some constant $C = C(n,m,q,M_0,N_0,C_e,C_h,\gamma) \in (0,\infty)$.  Here $\op{sep} a_j$ and $\op{width} a_j$ are as in Definition~\ref{a sep width defn} with $a_j$ in place of $a$.  $E_{Z,(1+\gamma)\rho/2,\rho}$ is as in \eqref{relharm excess defn2} with $\varsigma = \varsigma(n,m,q) \in (0,1)$ (and $(1+\gamma)/2$ in place of $\gamma$).
\end{lemma}

\begin{lemma}\label{relharm wellsep2 lemma} 
Let $M_0,N_0,C_e,C_h,\varepsilon_2$ be as in Hypothesis~\ref{DeLSpa constants hypothesis}, Lemma~\ref{normal approx energy lemma} and Lemma~\ref{relharm critrad2 lemma}.  For all $p \in \{1,2,\ldots,q\}$, $\gamma \in (0,1)$, $\Lambda \in (0,\infty)$, and $\eta \in (0,1)$ there exist $\varepsilon, \beta \in (0,1)$ depending only on $n,m,q,p,M_0,N_0,C_e,C_h,\gamma,\Lambda,\eta$ such that the following holds true.  Assume that Hypothesis~\ref{relharm separation1 hyp} holds true 
and that $\rho > \frac{2}{1+\gamma}\overline{R}(Z)$, where $\overline{R}$ is as in Definition~\ref{crit rad defn}.  For each $j \in \{1,2,\ldots,p\}$, let $N_j : B^{\mathcal M}_{(1+\gamma)\rho/2}(Z) \rightarrow \mathcal{A}_{q_j}({\rm T}\mathcal{M}^{\perp})$ satisfy \eqref{relharm normal approx concl1} and \eqref{relharm normal approx concl2} with $(1+\gamma)/2$ in place of $\gamma$.  For each $j \in \{1,2,\ldots,p\}$, let $a_j = \sum_{l=1}^{q_j} \llbracket a_{j,l} \rrbracket \in \mathcal{A}_{q_j}(\mathbb{R}^{n+m})$ such that $|a_j| < \varepsilon\rho$ and $\frac{1}{q} \sum_{j=1}^p q_j a_{j,a} = 0$.  Suppose that $\op{size} a_j > 1$ for some $j \in \{1,2,\ldots,p\}$ and that \eqref{relharm separation2 hyp1} and \eqref{relharm separation2 hyp2} hold true.  Then 
\begin{equation}\label{relharm wellsep2 concl}
	\frac{1}{\rho^{n-2}} \sum_{j=1}^p \int_{B^{\mathcal M}_{(1+3\gamma)\rho/4}(Z)} \mathcal{G}(DN_j, q_j \llbracket D\varphi_j \rrbracket)^2 
	\leq \eta^2 \min_{1 \leq j \leq p} (\op{sep} a_j)^2 .
\end{equation}
Here $E_{Z,(1+\gamma)\rho/2,\rho}$ is as in \eqref{relharm excess defn2} with $\varsigma = 2\beta_2$ (and $(1+\gamma)/2$ in place of $\gamma$).
\end{lemma} 

\begin{proof} 
The proof is similar to that of~\cite[Lemma~2.20]{KrumWicb} and we sketch the proof for completion.  Fix $\gamma \in (0,1)$, $\Lambda \in (0,\infty)$, and $\eta \in (0,1)$.  For $\nu = 1,2,3,\ldots$, suppose that $\varepsilon_{\nu} \rightarrow 0^+$, $\beta_{\nu} \rightarrow 0^+$, $\mathbf{m}_{0,\nu}$, $\mathcal{M}^{(\nu)}$, $Z_{\nu}$, $\rho_{\nu}$, $T^{(\nu)}$, $T^{(\nu)}_j$, $\mathcal{K}_{\nu}$, $N^{(\nu)}$, $N^{(\nu)}_j$, $\varphi^{(\nu)}_j$, and $a^{(\nu)}_j$ satisfy the hypotheses of Lemma~\ref{relharm wellsep2 lemma} (and in particular \eqref{relharm separation2 hyp1} and \eqref{relharm separation2 hyp2}) in place of $\varepsilon$, $\beta$, $\mathbf{m}_0$, $\mathcal{M}$, $Z$, $\rho$, $T$, $T_j$, $\mathcal{K}$, $N$, $N_j$, $\varphi_j$, and $a_j$.  For each $j \in \{1,2,\ldots,p\}$, let $a^{(\nu)}_j = \sum_{l=1}^{q_j} \llbracket a^{(\nu)}_{j,l} \rrbracket$ where $a^{(\nu)}_{j,l} \in \mathbb{R}^{n+m}$ and let $a^{(\nu)}_{j,a} = \frac{1}{q_j} \sum_{l=1}^{q_j} a^{(\nu)}_{j,l}$.  Suppose to the contrary that 
\begin{equation}\label{relharm wellsep2 eqn1}
	\frac{1}{\rho_{\nu}^{n-2}} \sum_{j=1}^p \int_{B^{\mathcal{M}^{(\nu)}}_{(1+3\gamma)\rho_{\nu}/4}(Z_{\nu})} 
		\mathcal{G}(DN^{(\nu)}_j, q_j \llbracket D\varphi^{(\nu)}_j \rrbracket)^2 
	> \eta^2 \min_{1 \leq j \leq p} (\op{sep} a^{(\nu)}_j)^2 . 
\end{equation}
After relabeling, assume that $\op{sep} a^{(\nu)}_1 = \max_{1 \leq j \leq p} \op{sep} a^{(\nu)}_j$.  Recall that $|a^{(\nu)}_j| < \varepsilon_{\nu} \rightarrow 0$ and in particular $\op{sep} a^{(\nu)}_1 \rightarrow 0$.  

By \eqref{relharm separation2 hyp1} 
\begin{equation}\label{relharm wellsep2 eqn2}
	|a^{(\nu)}_{j,l} - a^{(\nu)}_{j,a}| \leq \Lambda \op{sep} a^{(\nu)}_1
\end{equation}
for all $j \in \{1,2,\ldots,p\}$ and $l \in \{1,2,\ldots,q_j\}$.  By the triangle inequality, \eqref{relharm separation2 hyp2}, and \eqref{relharm wellsep2 eqn2} 
\begin{equation}\label{relharm wellsep2 eqn3}
	\frac{1}{\rho_{\nu}^n} \sum_{j=1}^p \int_{B^{\mathcal{M}^{(\nu)}}_{(1+\gamma)\rho_{\nu}/2}(Z_{\nu})} 
		\mathcal{G}(N^{(\nu)}_j, q_j \llbracket \varphi^{(\nu)}_j + a^{(\nu)}_{j,a} \rrbracket)^2 \leq C (\op{sep} a^{(\nu)}_1)^2 , 
\end{equation}
where $C = C(n,m,q,M_0,N_0,C_e,C_h,\Lambda,\gamma) \in (0,\infty)$ is a constant.  Note that by \eqref{relharm separation2 hyp2} and Lemma~\ref{normal approx energy lemma} 
\begin{equation}\label{relharm wellsep2 eqn4}
	\frac{\rho_{\nu}^4}{\varepsilon_{\nu}^2} \sum_{j=1}^p \|\Delta_{\mathcal{M}^{(\nu)}} \varphi^{(\nu)}_j\|_{C^0(B^{\mathcal{M}^{(\nu)}}_{\rho_{\nu}}(Z_{\nu}))}^2 
	\leq \beta_{\nu}^2 \,(\op{sep} a^{(\nu)}_1)^2
\end{equation}
and 
\begin{align}\label{relharm wellsep2 eqn5}
	&\sum_{j=1}^p \|D\varphi^{(\nu)}_j\|_{C^0(B^{\mathcal{M}^{(\nu)}}_{\rho_{\nu}}(Z_{\nu}))} 
		\cdot \frac{1}{\rho_{\nu}^{n-2}} \int_{B^{\mathcal{M}^{(\nu)}}_{(1+\gamma)\rho_{\nu}/2}(Z_{\nu})} |DN^{(\nu)}|^2 
		+ \frac{1}{\rho_{\nu}^{n-\varsigma}} \int_{B^{\mathcal{M}^{(\nu)}}_{\rho_{\nu}}(Z_{\nu})} |N^{(\nu)}|^2 
	\\ \leq\,& \left( \sum_{j=1}^p \|D\varphi^{(\nu)}_j\|_{C^0(B^{\mathcal{M}^{(\nu)}}_{\rho_{\nu}}(Z_{\nu}))} + \rho_{\nu}^{\varsigma} \right) 
		\frac{1}{\rho_{\nu}^n} \int_{B^{\mathcal{M}^{(\nu)}}_{\rho_{\nu}}(Z_{\nu})} |N^{(\nu)}|^2 \leq \beta_{\nu}^2 \,(\op{sep} a^{(\nu)}_1)^2 . 
		\nonumber 
\end{align}
Recall Remark~\ref{relharm norm scaling rmk} and that $\rho_{\nu} > \frac{2}{1+\gamma} \overline{R}(Z_{\nu})$, where $\overline{R}$ is as in Definition~\ref{crit rad defn} with $T^{(\nu)}$ in place of $T$.  By Lemma~\ref{relharm energy lemma} (with $\frac{(1+\gamma)\rho_{\nu}}{2}, \frac{5+7\gamma}{6+6\gamma}$ in place of $\rho,\gamma$), \eqref{relharm wellsep2 eqn3}, \eqref{relharm wellsep2 eqn4}, and \eqref{relharm wellsep2 eqn5} 
\begin{equation}\label{relharm wellsep2 eqn6}
	\frac{1}{\rho_{\nu}^{n-2}} \sum_{j=1}^p \int_{B^{\mathcal{M}^{(\nu)}}_{(5+7\gamma)\rho_{\nu}/12}(Z_{\nu})} 
		\mathcal{G}(DN^{(\nu)}_j, q_j \llbracket D\varphi^{(\nu)}_j \rrbracket)^2 \leq C (\op{sep} a^{(\nu)}_1)^2 , 
\end{equation}
where $C = C(n,m,q,M_0,N_0,C_e,C_h,\Lambda,\gamma) \in (0,\infty)$ is a constant.  Let $\mathfrak{q}_{\nu}$ be a rotation of $\mathbb{R}^{n+m}$ such that $\mathfrak{q}_{\nu}({\rm T}_{Z_{\nu}} \mathcal{M}^{(\nu)}) = \mathbb{R}^n \times \{0\}$.  Let $\mathbf{e}_{\nu} : B_1(0) \cong B_1(0,{\rm T}_{Z_{\nu}} \mathcal{M}^{(\nu)}) \rightarrow \mathcal{M}^{(\nu)}$ be the composition of $\mathfrak{q}_{\nu}^{-1}$ and the exponential map of $\mathcal{M}^{(\nu)}$ at $Z_{\nu}$.  By Lemma~\ref{relharm harmonic approx lemma} (with $\frac{(1+\gamma)\rho_{\nu}}{2}, \frac{2+4\gamma}{3+3\gamma}$ in place of $\rho,\gamma$), \eqref{relharm wellsep2 eqn4}, \eqref{relharm wellsep2 eqn5}, and \eqref{relharm wellsep2 eqn6}, for each $j \in \{1,2,\ldots,p\}$ there exists a locally Dirichlet energy minimizing $q_j$-valued function $w^{(\nu)} : B^{\mathcal{M}^{(\nu)}}_{(1+2\gamma)\rho_{\nu}/3}(Z_{\nu}) \rightarrow \mathcal{A}_q(\mathbb{R}^{n+m})$ such that 
\begin{align}\label{relharm wellsep2 eqn7}
	&\frac{1}{\rho_{\nu}^n} \sum_{j=1}^p \int_{B^{\mathcal{M}^{(\nu)}}_{(1+2\gamma)\rho_{\nu}/3}(Z_{\nu})} 
		\mathcal{G}(N^{(\nu)}_j,\varphi^{(\nu)}_j+w^{(\nu)}_j \circ \mathbf{e}_{\nu}^{-1})^2 
	\\&+ \frac{1}{\rho_{\nu}^{n-2}} \sum_{j=1}^p \int_{B^{\mathcal{M}^{(\nu)}}_{(1+2\gamma)\rho_{\nu}/3}(Z_{\nu})} 
		(|D(N^{(\nu)}_j-\varphi^{(\nu)}_j)| - |D(w^{(\nu)}_j \circ \mathbf{e}_{\nu}^{-1})|)^2 
	\leq \frac{1}{\nu^2} \,(\op{sep} a^{(\nu)}_1)^2 . \nonumber 
\end{align}
By \eqref{relharm wellsep2 eqn2}, after passing to a subsequence 
\begin{equation*}
	\frac{a^{(\nu)} - a^{(\nu)}_{j,a}}{\op{sep} a^{(\nu)}_1} \rightarrow \widetilde{a}_j 
\end{equation*} 
in $\mathcal{A}_{q_j}(\mathbb{R}^m)$ as $\nu\rightarrow\infty$.  By \eqref{relharm wellsep2 eqn3} and \eqref{relharm wellsep2 eqn7}, 
\begin{equation*}
	\rho_{\nu}^{-n/2} \|w^{(\nu)}_j - a^{(\nu)}_{j,a}\|_{L^2(B_{(1+2\gamma)\rho_{\nu}/3}(0))} \leq C \op{sep} a^{(\nu)}_1, 
\end{equation*}
where $C = C(n,m,q,M_0,N_0,C_e,C_h,\Lambda,\gamma) \in (0,\infty)$ is a constant.  Hence after passing to a subsequence for each $j \in \{1,2,\ldots,p\}$ there exists a Dirichlet energy minimizing $q$-valued function $\widetilde{w}_j : B_{(1+3\gamma)/4}(0) \rightarrow \mathcal{A}_{q_j}(\mathbb{R}^{n+m})$ such that 
\begin{equation}\label{relharm wellsep2 eqn8}
	\frac{w^{(\nu)}_j(\rho_{\nu} x) - a^{(\nu)}_{j,a}}{\op{sep} a^{(\nu)}_1} \rightarrow \widetilde{w}_j(x)
\end{equation} 
uniformly on $B_{(1+3\gamma)/4}(0)$ as $\nu\rightarrow\infty$.  It follows from \eqref{relharm separation2 hyp2} and \eqref{relharm wellsep2 eqn7} that $\op{spt} \widetilde{w}_j(x) \subseteq \op{spt} \widetilde{a}_j$ for all $x \in B_{(1+3\gamma)/4}(0)$.  Hence by the continuity of Dirichlet energy minimizing $q$-valued functions, $\widetilde{w}_j$ is a constant $q$-valued function.  Therefore, by \eqref{relharm wellsep2 eqn7}, \eqref{relharm wellsep2 eqn8}, and~\cite[Proposition~3.20]{DeLSpaDirMin} (\cite[Theorem~2.15]{Almgren}) 
\begin{align*}
	&\lim_{\nu\rightarrow\infty} \frac{1}{\rho_{\nu}^{n-2} (\op{sep} a^{(\nu)}_1)^2} \sum_{j=1}^p 
		\int_{B^{\mathcal{M}^{(\nu)}}_{(1+3\gamma)\rho_{\nu}/4}(Z_{\nu})} \mathcal{G}(DN^{(\nu)}_j,q_j \llbracket D\varphi^{(\nu)}_j \rrbracket)^2 
	\\=\,& \lim_{\nu\rightarrow\infty} \frac{1}{\rho_{\nu}^{n-2} (\op{sep} a^{(\nu)}_1)^2} \sum_{j=1}^p \int_{B_{(1+3\gamma)\rho_{\nu}/4}(0)} |Dw^{(\nu)}_j|^2 
	= \sum_{j=1}^p \int_{B_{(1+3\gamma)/4}(0)} |D\widetilde{w}_j|^2 = 0, 
\end{align*}
contradicting \eqref{relharm wellsep2 eqn1}.
\end{proof}

\begin{proof}[Proof of Lemma~\ref{relharm separation2 lemma}] 
Let $\varsigma$ be the minimum of $2\beta_2$ and its value in Theorem~\ref{relharm separation1 thm}.  Let $\eta = \eta(n,m,q,p,s,M_0,N_0,C_e,C_h,\gamma) \in (0,1)$ to be later determined.  Recall Remark~\ref{relharm norm scaling rmk} and note that by \eqref{relharm separation2 go condition}, $\rho > \frac{2}{1-\gamma}\,\overline{R}(Z)$.  Assuming $\beta = \beta(n,m,q,p,s,M_0,$ $N_0,C_e,C_h,\gamma,\eta)$ is sufficiently small, by Lemma~\ref{relharm wellsep2 lemma} 
\begin{equation}\label{relharm separation2 eqn1}
	\frac{1}{\rho^{n-2}} \sum_{j=1}^p \int_{B^{\mathcal M}_{(1+3\gamma)\rho/4}(Z)} \mathcal{G}(DN_j, q_j \llbracket D\varphi_j \rrbracket)^2 
	\leq \eta^2 \min_{1 \leq j \leq p} (\op{sep} a_j)^2 .
\end{equation}
By \eqref{relharm separation2 hyp2} and Lemma~\ref{normal approx energy lemma} 
\begin{align}\label{relharm separation2 eqn2}
	&\frac{\rho^4}{\varepsilon^2} \sum_{j=1}^p \|\Delta_{\mathcal M} \varphi_j\|_{C^0(B^{\mathcal M}_{\rho}(Z))}^2  
	+ \sum_{j=1}^p \|D\varphi_j\|_{C^0(B^{\mathcal M}_{\rho}(Z))} \cdot \frac{1}{\rho^{n-2}} \int_{B^{\mathcal M}_{(3+5\gamma)\rho/8}(Z)} |DN|^2 
		\\&+ \frac{1}{\rho^{n-\varsigma}} \int_{B^{\mathcal M}_{\rho}(Z)} |N|^2 
	\leq \beta^2 \min_{1 \leq j \leq p} (\op{sep} a_j)^2 \nonumber 
\end{align}
By \eqref{relharm separation2 go condition}, \eqref{relharm separation1 go condition} holds true with $\frac{(3+5\gamma)\rho}{8}, \frac{1+7\gamma}{3+5\gamma}$ in place of $\rho,\gamma$.  Thus using Lemma~\ref{relharm poincare lemma}, \eqref{relharm separation2 eqn1}, and \eqref{relharm separation2 eqn2}, for each $j \in \{1,2,\ldots,p\}$ there exists $b_j \in \mathcal{A}_{q_j}(\mathbb{R}^{n+m})$ such that 
\begin{align}\label{relharm separation2 eqn3}
	\frac{1}{\rho^n} \int_{\mathbf{C}^{\mathcal M}_{(1+3\gamma)\rho/4}(Z)} \op{dist}^2(X, \op{spt} \mathbf{G}_{\varphi_j+b_j}) \,d\|T_j\|(X)&
	\\+ \frac{1}{\rho^n} \int_{B^{\mathcal M}_{(1+3\gamma)\rho/4}(Z)} \mathcal{G}(N_j,\varphi_j+b_j)^2 
	&\leq C (\beta+\eta)^2 (\op{sep} a_j)^2 , \nonumber 
\end{align}
where $C = C(n,m,q,M_0,N_0,C_e,C_h) \in (0,\infty)$ is a constant.  Then by Theorem~\ref{relharm separation1 thm} (with $\frac{(3+5\gamma)\rho}{8}, \frac{1+7\gamma}{3+5\gamma}$ in place of $\rho,\gamma$), 
\begin{equation}\label{relharm separation2 eqn4}
	\sup_{X \in \op{spt} T_j \cap \mathbf{C}^{\mathcal M}_{(1+7\gamma)\rho/8}(Z)} \op{dist}(X, \op{spt}\mathbf{G}_{\varphi_j+b_j}) 
	< C (\beta+\eta) \op{sep} a_j ,
\end{equation}
where $C = C(n,m,q,M_0,N_0,C_e,C_h,\gamma) \in (0,\infty)$ is a constant.  
Now by the triangle inequality, \eqref{relharm separation2 hyp2}, and \eqref{relharm separation2 eqn3} 
\begin{align}\label{relharm separation2 eqn5}
	\mathcal{G}(b_j,a_j)^2 
	\leq\,& \frac{C}{\rho^n} \int_{B^{\mathcal M}_{(1+3\gamma)\rho/4}(Z)} \mathcal{G}(N_j,\varphi_j+b_j)^2 
		\\&+ \frac{C}{\rho^n} \int_{B^{\mathcal M}_{(1+3\gamma)\rho/4}(Z)} \sum_{l=1}^{q_j} \mathcal{G}(N_j,\varphi_j+a_j)^2  \nonumber
	\leq C (\beta+\eta)^2 (\op{sep} a_j)^2 ,  
\end{align}
where $C = C(n,m,q,M_0,N_0,C_e,C_h,\gamma) \in (0,\infty)$ are constants.  Thus by \eqref{relharm separation2 eqn4} and \eqref{relharm separation2 eqn5}
\begin{equation}\label{relharm separation2 eqn6}
	\sup_{X \in \op{spt} T_j \cap \mathbf{C}^{\mathcal M}_{(1+7\gamma)\rho/8}(Z)} \op{dist}(X, \op{spt}\mathbf{G}_{\varphi_j+a_j}) 
	\leq C (\beta+\eta) \op{sep} a_j ,
\end{equation}
where $C = C(n,m,q,M_0,N_0,C_e,C_h,\gamma) \in (0,\infty)$ are constants.  Choosing $\eta$ and $\beta$ small enough that $C (\beta+\eta) < 1/3$ (where $C$ is as in \eqref{relharm separation2 eqn6}), we deduce that 
\begin{equation*}
	T_j = \sum_{l=1}^{q_j} T_l \quad\text{where}\quad 
	T_{j,l} = T_j \llcorner \Big\{ X \in \mathbf{C}^{\mathcal M}_{(1+7\gamma)\rho/8}(Z) : 
		\op{dist}(X,\op{spt}\mathbf{G}_{\varphi_j+a_{j,l}}) < \tfrac{1}{3} \op{sep} a_j \Big\} 
\end{equation*}
and $T_{j,l}$ are locally area-minimizing rectifiable currents of $\mathbf{C}^{\mathcal M}_{(1+3\gamma)\rho/4}(Z)$ with $(\partial T_{j,l}) \llcorner \mathbf{C}^{\mathcal M}_{(1+7\gamma)\rho/8}(Z) = 0$.  $\varphi_j+a_{j,l}$ in place of $\varphi_j$.  By \eqref{relharm separation2 go condition}, \eqref{relharm separation0 go condition} holds true with $\frac{(1+7\gamma)\rho}{8}, \frac{8\gamma}{1+7\gamma}$ in place of $\rho,\gamma$.  Thus we can apply Lemma~\ref{relharm separation0 lemma} with $T_{j,l},\varphi_j+a_{j,l}$ in place of $T_j,\varphi_j+a_j$ (and with $\frac{(1+7\gamma)\rho}{8}, \frac{8\gamma}{1+7\gamma}$ in place of $\rho,\gamma$) to conclude that \eqref{relharm separation2 concl} holds true.
\end{proof}

\subsection{$L^{\infty}$-distance estimate in terms $L^2$-distance relative to disjoint, not-necessarily parallel $C^2$-graphs}  In this subsection, we prove Theorem~\ref{relharm separation main thm}.  As discussed at the beginning of this section, we shall proceed by induction by assuming Hypothesis~\ref{relharm grand hyp} that $T$ separates into two or more locally area minimizing rectifiable currents $T_j$ and $\varphi$ separating into corresponding $q_j$-valued functions $\varphi_j$.  Then we prove Theorem~\ref{relharm separation3 thm} in which we bound the $L^{\infty}$-distance of $T_j$ to the graph of $\varphi_j$.  We deduce that Theorem~\ref{relharm separation main thm} holds true as Theorem~\ref{relharm separation3 thm} in the special case $p = 1$.  

Let $M_0,N_0,C_e,C_h,\varepsilon_2$ be as in Hypothesis~\ref{DeLSpa constants hypothesis}, Lemma~\ref{normal approx energy lemma} and Lemma~\ref{relharm critrad2 lemma}.  Let $1 \leq p \leq q$ be integers.  Let $\varepsilon, \beta \in (0,1)$, $\gamma \in (0,1)$, and $\kappa \in (0,\infty)$ be suitable constants.  We shall assume that:

\begin{hypothesis}\label{relharm grand hyp} 
Let $T, \mathscr{W}, \mathcal{M}, N, \mathcal{K}$ be as in Hypothesis~\ref{M and N hypothesis}.  Let $Z \in \mathcal{M} \cap \mathbf{C}_1(0)$ and $\rho \in (0,2]$.  Let $p,q_1,q_2,\ldots,q_p$ be positive integers with $\sum_{j=1}^p q_j = q$ and for each $j \in \{1,2,\ldots,p\}$ let $T_j$ be $n$-dimensional locally area minimizing rectifiable currents of $\mathbf{C}^{\mathcal M}_{\rho}(Z)$ such that \eqref{relharm separation1 hyp eqn1} holds true. For each $j \in \{1,2,\ldots,p\}$, let $\varphi_j : B^{\mathcal M}_{\rho}(Z) \rightarrow \mathcal{A}_{q_j}(\mathbb{R}^{n+m})$ such that 
\begin{equation*}
	\varphi_j(X) = \sum_{l=1}^{q_j} \llbracket \varphi_{j,l}(X) \rrbracket 
\end{equation*}
for each $X \in B^{\mathcal M}_{\rho}(Z)$, where $\varphi_{j,1},\ldots,\varphi_{j,q_j} \in C^2(B^{\mathcal M}_{\rho}(Z),\mathbb{R}^{n+m})$ are single-valued functions such that 
\begin{gather}
	\label{relharm grand hyp2} \rho^{-1} \|\varphi_{j,l}\|_{C^0(B^{\mathcal M}_{\rho}(Z))} + \|D\varphi_{j,l}\|_{C^0(B^{\mathcal M}_{\rho}(Z))} 
		+ \rho \|D^2 \varphi_{j,l}\|_{C^0(B^{\mathcal M}_{\rho}(Z))} < \varepsilon , \\
	\label{relharm separation3 avgfree} \frac{1}{q} \sum_{j=1}^p \sum_{l=1}^{q_j} \varphi_{j,l}(X) = 0 \text{ for all } X \in B^{\mathcal M}_{\rho}(Z) , \\
	\label{relharm separation3 normal} \varphi_{j,l}(X) \in {\rm T}_X \mathcal{M}^{\perp} \text{ for all } X \in B^{\mathcal M}_{\rho}(Z) 
\end{gather}
for all $j \in \{1,2,\ldots,p\}$ and $l \in \{1,2,\ldots,q_j\}$ and 
\begin{equation}\label{relharm separation3 wellsep} 
	\rho \sup_{B^{\mathcal M}_{\rho}(Z)} |D\varphi_{j,k} - D\varphi_{j,l}| \leq \kappa \inf_{B^{\mathcal M}_{\rho}(Z)} |\varphi_{j,k} - \varphi_{j,l}| 
\end{equation}
for all $j \in \{1,2,\ldots,p\}$ and $k,l \in \{1,2,\ldots,q_j\}$.
\end{hypothesis}

Let $\varepsilon,\varsigma \in (0,1)$.  Assume that Hypothesis~\ref{relharm grand hyp} holds true.  For each $j \in \{1,2,\ldots,p\}$, let $N_j : B^{\mathcal M}_{(1+\gamma)\rho/2}(Z) \rightarrow \mathcal{A}_{q_j}({\rm T}\mathcal{M}^{\perp})$ satisfy \eqref{relharm normal approx concl1} and \eqref{relharm normal approx concl2}.  As a slight abuse of notation, we define 
\begin{align}\label{relharm excess defn3}
	E_{Z,\gamma\rho,\rho}^2 =\,& \frac{1}{\rho^{n+2}} \sum_{j=1}^p \int_{B^{\mathcal M}_{\gamma\rho}(Z)} \mathcal{G}(N_j,\varphi_j)^2 
		+ \frac{\rho^2}{\varepsilon^2} \sum_{j=1}^p \sum_{l=1}^{q_j} \|\Delta_{\mathcal M} \varphi_{j,l}\|_{C^0(B^{\mathcal M}_{\rho}(Z))}^2  		
		\\&+ \left( \sum_{j=1}^p \|D\varphi_j\|_{C^0(B^{\mathcal M}_{\rho}(Z))} + \rho^{\varsigma} \right) 
		\frac{1}{\rho^{n+2}} \int_{B^{\mathcal M}_{\rho}(Z)} |N|^2 \nonumber 
\end{align}
$E_{Z,\gamma\rho,\rho}$ as defined in \eqref{relharm excess defn3} is an $L^2$-distance of $T_j$ to the graph of $\varphi_j$ together with additional terms to control $\Delta_{\mathcal M} \varphi_{j,l}$ and error terms arising from the area formula and estimates for the normal approximation maps $N_j$ and $N$ (cf.~\eqref{relharm excess defn2}).

\begin{theorem}\label{relharm separation3 thm} 
Let $M_0,N_0,C_e,C_h,\varepsilon_2$ be as in Hypothesis~\ref{DeLSpa constants hypothesis}, Lemma~\ref{normal approx energy lemma} and Lemma~\ref{relharm critrad2 lemma}.  For all integers $1 \leq p \leq q$ and $s \geq p$ and all $\kappa \in (0,\infty)$ and $\gamma \in (0,1)$ there exists $\varepsilon, \beta \in (0,1)$ depending only on $n,m,q,p,s,M_0,N_0,C_e,C_h,\gamma,\kappa$ such that the following holds true.  Suppose that Hypothesis~\ref{relharm grand hyp} holds true.  For each $j \in \{1,2,\ldots,p\}$, let $N_j : B^{\mathcal M}_{(1+\gamma)\rho/2}(Z) \rightarrow \mathcal{A}_{q_j}({\rm T}\mathcal{M}^{\perp})$ satisfy \eqref{relharm normal approx concl1} and \eqref{relharm normal approx concl2} with $(1+\gamma)/2$ in place of $\gamma$.  Suppose that 
\begin{equation}\label{relharm separation3 go condition} 
	\sup_{X \in B^{\mathcal M}_{(5+3\gamma)\rho/8}(Z)} \overline{R}(X) < 4^{-s+p-4} (1-\gamma) \rho , 
\end{equation}
where $\overline{R}$ is as in Definition~\ref{crit rad defn}, and that $\sum_{j=1}^p \op{size}\varphi_j(X) = s$ for all $X \in B^{\mathcal M}_{\rho}(Z)$ (where $\op{size}$ is as in Definition~\ref{a sep width defn}).  Then for all $j \in \{1,2,\ldots,p\}$ 
\begin{align}\label{relharm separation3 concl}
	\sup_{X \in \op{spt} T_j \cap \mathbf{C}^{\mathcal M}_{\gamma\rho}(Z)} \op{dist}(X, \op{spt} \mathbf{G}_{\varphi_j}) 
	\leq C E_{Z,(1+\gamma)\rho/2,\rho} \rho 
\end{align}
where $E_{Z,(1+\gamma)\rho/2,\rho}$ is as in \eqref{relharm excess defn3} with $\varsigma = \varsigma(n,m,q) \in (0,1)$ (and $(1+\gamma)/2$ in place of $\gamma$) and where $C = C(n,m,q,p,s,M_0,N_0,C_e,C_h,\gamma) \in (0,\infty)$ are constants.
\end{theorem}

\begin{definition}\label{relharm minsep width}  
For each $j \in \{1,2,\ldots,p\}$, let $\varphi_j : B^{\mathcal M}_{\rho}(Z) \rightarrow \mathcal{A}_{q_j}(\mathbb{R}^{n+m})$ be as in Hypotheses~\ref{relharm grand hyp} holds true.  For each $j \in \{1,2,\ldots,p\}$, if $\op{size} \varphi_j > 1$ on $B^{\mathcal M}_{\rho}(X)$, then we define 
\begin{align*}
	\op{minsep} \varphi_j =\,& \inf_{X \in B^{\mathcal M}_{\rho}(Z)} \op{sep} \varphi_j(X) , \\
	\op{width} \varphi_j =\,& \sup_{X \in B^{\mathcal M}_{\rho}(Z)} \op{width} \varphi_j(X) .
\end{align*}
If instead $\op{size} \varphi_j = 1$ on $B^{\mathcal M}_{\rho}(X)$, then we define $\op{minsep} \varphi_j = \infty$ and $\op{width} \varphi_j = 0$.  (Here $\op{size}$, $\op{sep}$, and $\op{width}$ are as in Definition~\ref{a sep width defn}.)
\end{definition}

\begin{remark}{\rm 
By \eqref{relharm separation3 wellsep}, $\op{size}\varphi_j(X)$ is a constant function of $X \in B^{\mathcal M}_{\rho}(Z)$.  Moreover, if $\op{size} \varphi_j > 1$ on $B^{\mathcal M}_{\rho}(Z)$, then $0 < \op{minsep} \varphi_j < \infty$ and $0 < \op{width} \varphi_j < \infty$.
}\end{remark}

Let $\varphi_j : B^{\mathcal M}_{\rho}(Z) \rightarrow \mathcal{A}_{q_j}(\mathbb{R}^{n+m})$ be as in Hypotheses~\ref{relharm grand hyp}.  For all $\sigma \in (0,\rho]$, $j \in \{1,2,\ldots,p\}$, and $k,l \in \{1,2,\ldots,q_j\}$, by \eqref{relharm separation3 wellsep} 
\begin{align}\label{relharm dist to sep}
	\sup_{B^{\mathcal M}_{\rho}(Z)} |\varphi_{j,k} - \varphi_{j,l}| 
	\leq\,& \inf_{B^{\mathcal M}_{\sigma}(Z)} |\varphi_{j,k} - \varphi_{j,l}| + 2\rho \sup_{B^{\mathcal M}_{\rho}(Z)} |D\varphi_{j,k} - D\varphi_{j,l}| 
	\\ \leq\,& (1 + 2\kappa) \inf_{B^{\mathcal M}_{\sigma}(Z)} |\varphi_{j,k} - \varphi_{j,l}| . \nonumber 
\end{align}

\begin{lemma}\label{relharm critrad3 lemma} 
Let $M_0,N_0,C_e,C_h,\varepsilon_2$ be as in Hypothesis~\ref{DeLSpa constants hypothesis}, Lemma~\ref{normal approx energy lemma} and Lemma~\ref{relharm critrad2 lemma}.  For every integer $p \in \{1,2,\ldots,q\}$ there exist $\varepsilon = \varepsilon(n,m,q,p,M_0,N_0,C_e,C_h) \in (0,1)$ such that if Hypothesis~\ref{relharm grand hyp} holds true and $\rho \geq \overline{R}(Z)$ where $\overline{R}(Z)$ is as in Definition~\ref{crit rad defn}, then for each $j \in \{1,2,\ldots,p\}$  
\begin{align}\label{relharm critrad3 concl}
	&\sup_{X \in \op{spt} T_j \cap \mathbf{C}^{\mathcal M}_{\rho/2}(Z)} \op{dist}^2(X, \op{spt} \mathbf{G}_{\varphi_j}) 
	\leq \frac{C}{\rho^n} \int_{\mathbf{C}^{\mathcal M}_{\rho}(Z)} \op{dist}^2(X, \op{spt} \mathbf{G}_{\varphi_j}) \,d\|T_j\|(X) 
	\\&+ \rho^4 \sum_{j=1}^p \sum_{l=1}^{q_j} \|\Delta_{\mathcal M} \varphi_{j,l}\|_{C^0(B^{\mathcal M}_{\rho}(Z))}^2  		
	+ \left( \sum_{j=1}^p \sum_{l=1}^{q_j} \|D\varphi_{j,l}\|_{C^0(B^{\mathcal M}_{\rho}(Z))} + \rho^{\varsigma} \right) 
	\frac{1}{\rho^n} \int_{B^{\mathcal M}_{\rho}(Z)} |N|^2 \nonumber \\&+ C \mathbf{m}_0 \rho^4 , \nonumber 
\end{align}
where $\varsigma = \varsigma(n,m,q) \in (0,1)$ and $C = C(n,m,q,M_0,N_0,C_e,C_h) \in (0,\infty)$ are constants.
\end{lemma}

\begin{proof} 
The proof is nearly identical to Lemma~\ref{relharm critrad1 lemma}.  Set $P = {\rm T}_Z \mathcal{M}$.  By Lemma~\ref{app lemma}, for each $j \in \{1,2,\ldots,p\}$ and $l \in \{1,2,\ldots,q_j\}$ there exists a (unique) single-valued function $\psi_{j,l} : B_{3\rho/4}(Z,P) \rightarrow P^{\perp}$ such that \begin{gather*}
    \mathbf{G}_{\psi_{j,l}} = \mathbf{G}_{\varphi_{j,l}} \llcorner \mathbf{C}_{3\rho/4}(Z,P) , \\
    \|D\psi_{j,l}\|_{C^0(B_{3\rho/4}(Z,P))} \leq C \mathbf{m}_0^{1/2} \rho + \|D\varphi_{j,l}\|_{C^0(B^{\mathcal M}_{\rho}(Z))}, \\ 
    \|\Delta\psi_{j,l}\|_{C^0(B_{3\rho/4}(Z,P))} \leq C \mathbf{m}_0^{1/2} + C \|\Delta\varphi_{j,l}\|_{C^0(B^{\mathcal M}_{\rho}(Z))} ,
\end{gather*} 
where $C = C(n,m,q,M_0,N_0,C_e,C_h) \in (0,\infty)$ is a constant (cf.~\eqref{relharm critrad1 eqn1}, \eqref{relharm critrad1 eqn2}, \eqref{relharm critrad1 eqn7}, \eqref{relharm critrad1 eqn9}).  By \eqref{relharm critrad1 eqn4}
\begin{equation*}
    |\pi_{{\rm T}_X \mathcal{M}}(\varphi_{j,l}(X) - \varphi_{j',l'}(X))| \leq C \mathbf{m}_0^{1/2} \rho \,|\varphi_{j,l}(X) - \varphi_{j',l'}(X)|
\end{equation*}
for all $X \in B^{\mathcal M}_{\rho}(Z)$, $1 \leq j,j' \leq p$, $1 \leq l \leq q_j$, and $1 \leq l' \leq q_{j'}$, where $C = C(n,m,q,M_0,N_0,C_e,C_h) \in (0,\infty)$ is a constant.  Hence provided $\varepsilon$ is sufficiently small, by Lemma~\ref{app lemma2}
\begin{equation*}
    \rho \sup_{B_{3\rho/4}(Z,P)} |D\psi_{j,l} - D\psi_{j',l'}| \leq 40\kappa \inf_{B_{3\rho/4}(Z,P)} |\psi_{j,l} - \psi_{j',l'}| 
\end{equation*}
for all $1 \leq j,j' \leq p$, $1 \leq l \leq q_j$, and $1 \leq l' \leq q_{j'}$ (cf. (4.69) of~\cite{KrumWicc}).  Thus to obtain \eqref{relharm critrad3 concl}, we can apply~\cite[Theorem~4.19]{KrumWicc} as in the proof of Lemma~\ref{relharm critrad1 lemma}, using the functions $\psi_{j,l}$ in place of the multi-valued $\psi_j+a_j$ and using \eqref{relharm critrad1 eqn10}, \eqref{relharm critrad1 eqn14}, and the direct analog of \eqref{relharm critrad1 eqn13} (i.e.~with $\varphi_{j,l},\psi_{j,l}$ in place of $\varphi_j,\psi_j$ and using Lemma~\ref{normal approx energy lemma}).
\end{proof}

\begin{lemma}\label{relharm wellsep3 lemma}
Let $M_0,N_0,C_e,C_h,\varepsilon_2$ be as in Hypothesis~\ref{DeLSpa constants hypothesis}, Lemma~\ref{normal approx energy lemma} and Lemma~\ref{relharm critrad2 lemma}.  For all integers $1 \leq p \leq q$ and $s \geq p$ and all $\kappa,\lambda \in (0,\infty)$ and $\gamma \in (0,1)$ there exists $\varepsilon, \beta \in (0,1)$ depending only on $n,m,q,p,s,M_0,N_0,C_e,C_h,\kappa,\lambda,\gamma$ such that the following holds true.  Suppose that Hypothesis~\ref{relharm grand hyp} holds true.  For each $j \in \{1,2,\ldots,p\}$, let $N_j : B^{\mathcal M}_{(1+\gamma)\rho/2}(Z) \rightarrow \mathcal{A}_{q_j}({\rm T}\mathcal{M}^{\perp})$ satisfy \eqref{relharm normal approx concl1} and \eqref{relharm normal approx concl2} with $(1+\gamma)/2$ in place of $\gamma$.  Suppose that 
\begin{equation}\label{relharm wellsep3 go condition}
	\sup_{X \in B^{\mathcal M}_{(5+3\gamma)\rho/8}(Z)} \overline{R}(X) < (1-\gamma)\rho/1024 , 
\end{equation}
where $\overline{R}$ is as in Definition~\ref{crit rad defn}, and that 
\begin{align}
	\label{relharm wellsep3 hyp1} \max_{1 \leq j \leq p} \op{width} \varphi_j &\leq \lambda \min_{1 \leq j \leq p} \op{minsep} \varphi_j , \\
	\label{relharm wellsep3 hyp2} E_{Z,(1+\gamma)\rho/2,\rho} \,\rho &\leq \beta \min_{1 \leq j \leq p} \op{minsep} \varphi_j ,
\end{align}
where $E_{Z,(1+\gamma)\rho/2,\rho}$ is as in \eqref{relharm excess defn3} with $\varsigma = \varsigma(n,m,q) \in (0,1)$ (and with $(1+\gamma)/2$ in place of $\gamma$).  Then for each $j \in \{1,2,\ldots,p\}$ 
\begin{equation}\label{relharm wellsep3 concl}
	\sup_{X \in \op{spt} T \cap \mathbf{C}^{\mathcal M}_{\gamma\rho}(Z)} \op{dist}(X, \op{spt} \mathbf{G}_{\varphi_j}) 
	\leq C E_{Z,(1+\gamma)\rho/2,\rho} \,\rho , 
\end{equation}
where $C = C(n,m,q,p,s,M_0,N_0,C_e,C_h,\gamma) \in (0,\infty)$ is a constant.
\end{lemma}

\begin{proof} 
The proof is similar to that of~\cite[Lemma~2.21]{KrumWicb}.  Let $\varsigma$ be the minimum of $2\delta_2$ and its values in Lemma~\ref{relharm critrad3 lemma}, Lemma~\ref{relharm separation0 lemma}, and Lemma~\ref{relharm separation2 lemma}.  Let $\sigma = \sigma(n,m,q,p,M_0,N_0,C_e,C_h,\kappa,\lambda,\gamma) \in (0,(1-\gamma)/8)$ be a constant to be later determined.  For each $X \in B^{\mathcal M}_{(1+\gamma)\rho/2}(Z)$, let 
\begin{equation*}
	\overline{\tau}(X) = \inf \bigg\{ \tau \in [2\sigma\rho, (1-\gamma)\rho/8] : \sup_{Y \in B^{\mathcal M}_{\rho'}(X)} \overline{R}(Y) < \rho'/128 
		\text{ for all } \rho' \in [\tau,(1-\gamma)\rho/8] \bigg\} 
\end{equation*}
By \eqref{relharm wellsep3 go condition}, when $\rho' = (1-\gamma)\rho/8$ we have that $\sup_{Y \in B^{\mathcal M}_{(1-\gamma)\rho/8}(X)} \overline{R}(Y) < (1-\gamma)\rho/1024$.  Noting that $\overline{\tau}(X)/2 \geq \sigma\rho$ for all $X \in B^{\mathcal M}_{(1+\gamma)\rho/2}(Z)$, we can cover $B^{\mathcal M}_{(1+\gamma)\rho/2}(Z)$ by a finite collection of balls $\{B^{\mathcal M}_{\rho_k}(X_k) : k = 1,2,\ldots,K\}$, where $X_k \in B^{\mathcal M}_{(1+\gamma)\rho/2}(Z)$ and $\rho_k = \overline{\tau}(X_k)/2$ for all $k \in \{1,2,\ldots,K\}$.  We may now consider the following two cases:
\begin{enumerate}[itemsep=2mm,topsep=0mm]
	\item[(i)]  $\rho_k = \sigma\rho$ and $\sup_{Y \in B^{\mathcal M}_{\rho_k}(X_k)} \overline{R}(Y) < \rho_k/128$; 
	\item[(ii)]  $\sup_{Y \in B^{\mathcal M}_{\rho_k}(X_k)} \overline{R}(Y) = \rho_k/128$.
\end{enumerate}

\noindent\textit{Case~(i).}  Suppose that for some fixed $k \in \{1,2,\ldots,K\}$, $\rho_k = \sigma\rho$ and $\sup_{Y \in B^{\mathcal M}_{11\rho_k/8}(X_k)} \overline{R}(Y) < \rho_k/128$.  Let $\varphi_{j,a}(X) = \frac{1}{q_j} \sum_{l=1}^{q_j} \varphi_{j,l}(X)$ for all $X \in B^{\mathcal M}_{\rho}(Z)$ and all $j \in \{1,2,\ldots,p\}$.  For each $j \in \{1,2,\ldots,p\}$ and $l \in \{1,2,\ldots,q_j\}$, define $\widetilde{\varphi}_{k,j,l} : B^{\mathcal M}_{2\rho_k}(X_k) \rightarrow \mathcal{A}_{q_j}(\mathbb{R}^{n+m})$ by 
\begin{equation*}
	\widetilde{\varphi}_{k,j,l}(X) = \varphi_{j,a}(X) + \varphi_{j,l}(X_k) - \varphi_{j,a}(X_k) 
\end{equation*}
for each $X \in B^{\mathcal M}_{2\rho_k}(X_k)$.  For each $j \in \{1,2,\ldots,p\}$, define $\widetilde{\varphi}_{k,j} : B^{\mathcal M}_{2\rho_k}(X_k) \rightarrow \mathcal{A}_{q_j}(\mathbb{R}^{n+m})$ by $\widetilde{\varphi}_{k,j}(X) = \sum_{l=1}^{q_j} \llbracket \widetilde{\varphi}_{k,j,l}(X) \rrbracket$ for each $X \in B^{\mathcal M}_{2\rho_k}(X_k)$.  Let 
\begin{align}\label{relharm wellsep3 eqn2}
	\widetilde{E}_k^2 =\,& \frac{1}{\rho_k^{n+2}} \int_{B^{\mathcal M}_{2\rho_k}(X_k)} \mathcal{G}(N_j,\widetilde{\varphi}_{k,j})^2 
		+ \rho_k^2 \sum_{j=1}^p \|\Delta_{\mathcal M} \varphi_{j,a}\|_{C^0(B^{\mathcal M}_{2\rho_k}(X_k))}^2 
		\\&+ \left( \sum_{j=1}^p \|D\varphi_{j,a}\|_{C^0(B^{\mathcal M}_{2\rho_k}(Z))} + \rho_k^{\varsigma} \right) 
		\frac{1}{\rho_k^{n+2}} \int_{B^{\mathcal M}_{4\rho_k}(X_k)} |N|^2 \nonumber
\end{align}
Recall that \eqref{relharm separation3 avgfree} holds true and notice that $\frac{1}{q_j} \sum_{l=1}^{q_j} (\varphi_{j,l}(X_k) - \varphi_{j,a}(X_k)) = 0$ for all $j \in \{1,2,\ldots,p\}$.  $\op{sep} \widetilde{\varphi}_{k,j}(X)$ is a constant function of $X \in B^{\mathcal M}_{2\rho_k}(X_k)$, whose value we can denote by $\op{sep} \widetilde{\varphi}_{k,j}$.  (In fact, $\op{sep} \widetilde{\varphi}_{k,j}(X) = \op{sep} \varphi_j(X_k)$ for all $X \in B^{\mathcal M}_{2\rho_k}(X_k)$.)  Since $\widetilde{\varphi}_{k,j}(X_k) = \varphi_j(X_k)$, we have that $\op{minsep} \varphi_j \leq \op{sep} \widetilde{\varphi}_j$ and by \eqref{relharm separation3 wellsep} and \eqref{relharm wellsep3 hyp1} 
\begin{align}\label{relharm wellsep3 eqn3}
	\sup_{B^{\mathcal M}_{2\rho_k}(X_k)} |\widetilde{\varphi}_{k,j,l} - \varphi_{j,l}| 
	\leq\,& 2\rho_k \sup_{B^{\mathcal M}_{2\rho_k}(X_k)} |D\widetilde{\varphi}_{k,j,l} - D\varphi_{j,l}| 
	= 2\sigma\rho \sup_{B^{\mathcal M}_{2\rho_k}(X_k)} |D\varphi_{j,a} - D\varphi_{j,l}| 
	\\ \leq\,& \frac{2\sigma\rho}{q_j} \sum_{i=1}^{q_j} \sup_{B^{\mathcal M}_{\rho}(Z)} |D\varphi_{j,i} - D\varphi_{j,l}| 
	\leq \frac{2\kappa\sigma}{q_j} \sum_{i=1}^{q_j} \inf_{B^{\mathcal M}_{\rho}(Z)} |\varphi_{j,i} - \varphi_{j,l}| \nonumber 
	\\ \leq\,& 2\kappa\sigma \op{width} \varphi_j 
	\leq 2\kappa\lambda\sigma \min_{1 \leq j' \leq p} \op{minsep} \varphi_{j'} . \nonumber 
\end{align}
Arguing as in~\cite{KrumWicb} using the triangle inequality, \eqref{relharm wellsep3 hyp2}, and \eqref{relharm wellsep3 eqn3} 
\begin{equation}\label{relharm wellsep3 eqn4}
	\widetilde{E}_k \rho_k \leq C \left( \frac{\beta}{\sigma^{n/2}} + \sigma \right) \min_{1 \leq j \leq p} \op{minsep} \varphi_j 
\end{equation}
for some constant $C = C(n,q,\gamma,\kappa,\lambda) \in (0,\infty)$ (cf.~(2.81) in~\cite{KrumWicb}).  Since $\sup_{Y \in B^{\mathcal M}_{\rho_k}(X_k)} \overline{R}(Y) < \rho_k/128$, we can apply Lemma~\ref{relharm separation2 lemma} (with $X_k,2\rho_k,1/2$ in place of $Z,\rho,\gamma$) to deduce that for each $j \in \{1,2,\ldots,p\}$ 
\begin{equation}\label{relharm wellsep3 eqn5}
	\sup_{X \in \op{spt} T_j \cap \mathbf{C}^{\mathcal M}_{\rho_k}(X_k)} \op{dist}(X, \op{spt} \mathbf{G}_{\widetilde{\varphi}_{k,j}}) 
	\leq C \left( \frac{\beta}{\sigma^{n/2}} + \sigma \right) \op{minsep} \varphi_j , 
\end{equation}
where $C = C(n,m,q,p,M_0,N_0,C_e,C_h,\kappa,\lambda,\gamma) \in (0,\infty)$ is a constant. 

\noindent\textit{Case~(ii).}  Suppose that for some fixed $k \in \{1,2,\ldots,K\}$, $\sup_{Y \in B^{\mathcal M}_{\rho_k}(X_k)} \overline{R}(Y) = \rho_k/128$.  For each $j \in \{1,2,\ldots,p\}$, set $\widetilde{\varphi}_{k,j}(X) = \varphi_j(X)$ for each $X \in B^{\mathcal M}_{2\rho_k}(X_k)$.  Let $\widetilde{E}_k$ is as in \eqref{relharm wellsep3 eqn2} with $\widetilde{\varphi}_{k,j} = \varphi_j$.  Notice that by the area formula, \eqref{relharm separation1 small norm}, \eqref{relharm normal approx concl1}, \eqref{relharm normal approx concl2}, and \eqref{relharm norm est rmk eqn1} 
\begin{align*}
	\int_{\mathbf{C}^{\mathcal M}_{2\rho_k}(X_k)} \op{dist}^2(X, \op{spt} \mathbf{G}_{\varphi_j}) \,d\|T_j\|(X) 
	\leq\,& C \int_{B^{\mathcal M}_{2\rho_k}(X_k)} \sum_{l=1}^{q_j} \op{dist}^2(N_{j,l}(X), \op{spt} \varphi_j(X)) \,d\mathcal{H}^n(X) 
		\\&+ C \mathbf{m}_0^{3\gamma_2/4} \rho^{3\gamma_2/4} \int_{B^{\mathcal M}_{4\rho_k}(X_k)} |N|^2 
	\\ \leq\,& C \widetilde{E}_k^2 \rho_k^{n+2} , 
\end{align*}
where $C = C(n,m,q,p,M_0,N_0,C_e,C_h,\kappa,\lambda,\gamma) \in (0,\infty)$ is a constant.  Thus by Lemma~\ref{relharm critrad3 lemma}

\begin{equation}\label{relharm wellsep3 eqn6}
	\sup_{X \in \op{spt} T_j \cap \mathbf{C}^{\mathcal M}_{\rho_k}(X_k)} \op{dist}(X, \op{spt} \mathbf{G}_{\varphi_j}) 
		\leq C \widetilde{E}_k \rho_k + C \mathbf{m}_0^{1/2} \rho_k^2 
\end{equation}
where $C = C(n,m,q,p,M_0,N_0,C_e,C_h,\kappa,\lambda,\gamma) \in (0,\infty)$ is a constant.  Observe that there exists $Y_k \in B^{\mathcal M}_{\rho_k}(X_k)$ such that $\rho_k/256 < \overline{R}(Y_k) \leq \rho_k/128$.  Thus provided $\varsigma \leq 2\delta_2$, by Lemma~\ref{relharm critrad2 lemma} (with $Y_k$ in place of $Z$) 
\begin{equation}\label{relharm wellsep3 eqn7}
	\mathbf{m}_0 \rho_k^2 
	\leq \frac{C}{\rho_k^{n-\varsigma}} \int_{B^{\mathcal M}_{\rho_k/128}(Y_k)} |N|^2 
	\leq \frac{C}{\rho_k^{n-\varsigma}} \int_{B^{\mathcal M}_{2\rho_k}(X_k)} |N|^2 , 
\end{equation}
where $C = C(n,m,q,M_0,N_0,C_e,C_h) \in (0,\infty)$ is a constant.  Combining \eqref{relharm wellsep3 eqn6}, \eqref{relharm wellsep3 eqn7}, \eqref{relharm wellsep3 hyp2}, and $\sigma\rho < \rho_k$ gives us 
\begin{equation*}
	\sup_{X \in \op{spt} T_j \cap \mathbf{C}^{\mathcal M}_{\rho_k}(X_k)} \op{dist}(X, \op{spt} \mathbf{G}_{\varphi_j}) 
	\leq C \widetilde{E}_k \rho_k \leq C\beta \bigg(\frac{\rho}{\rho_k}\bigg)^{n/2} \op{minsep} \varphi_j 
	\leq \frac{C\beta}{\sigma^{n/2}} \op{minsep} \varphi_j , 
\end{equation*}
where $C = C(n,m,q,p,M_0,N_0,C_e,C_h,\kappa,\lambda,\gamma) \in (0,\infty)$ are constants, cf.~\eqref{relharm wellsep3 eqn5}.  

\noindent\textit{Separation of $T_j$ into components $T_{j,l}$ and bounding the distance to the graph of $\varphi_j$.}  Now we have shown that \eqref{relharm wellsep3 eqn5} holds true for all $k \in \{1,2,\ldots,K\}$.  Assuming that $\sigma$ is small enough that $C\sigma < 1/9$ and $\beta$ is small enough (depending on $\sigma$) that $C \beta \sigma^{-n/2} < 1/9$, it follows from \eqref{relharm wellsep3 eqn5} using $\op{minsep} \varphi_j \leq \op{minsep} \widetilde{\varphi}_{k,j}$ that 
\begin{equation}\label{relharm wellsep3 eqn8}
	\sup_{X \in \op{spt} T_j \cap \mathbf{C}^{\mathcal M}_{\rho_k}(X_k)} \op{dist}(X, \op{spt} \mathbf{G}_{\widetilde{\varphi}_{k,j}}) 
	< \frac{2}{9} \op{minsep} \varphi_j \leq \frac{2}{9} \op{minsep} \widetilde{\varphi}_{k,j} .
\end{equation}
By \eqref{relharm wellsep3 eqn8}, for each $k \in \{1,2,\ldots,K\}$ and $j \in \{1,2,\ldots,p\}$ 
\begin{equation*}
	T_j \llcorner \mathbf{C}^{\mathcal M}_{\rho_k}(X_k) = \sum_{i=1}^s T_{k,j,l} \text{ where } 
	T_{k,j,l} = T_j \llcorner \Big\{ X \in \mathbf{C}^{\mathcal M}_{\rho_k}(X_k) : 
		\op{dist}(X, \op{spt} \mathbf{G}_{\widetilde{\varphi}_{k,j,l}}) < \tfrac{2}{9} \op{minsep} \varphi_j \Big\} 
\end{equation*}
and $T_{k,j,l}$ are locally area-minimizing rectifiable currents of $\mathbf{C}^{\mathcal M}_{\rho_k}(X_k)$ with $(\partial T_{k,j,l}) \llcorner \mathbf{C}^{\mathcal M}_{\rho_k}(X_k) = 0$.  Assuming that $2\kappa\lambda\sigma < 1/9$, by the triangle inequality, \eqref{relharm wellsep3 eqn3}, and \eqref{relharm wellsep3 eqn8} 
\begin{equation*}
	\sup_{X \in \op{spt} T_{k,j,l}} \op{dist}(X, \op{spt} \mathbf{G}_{\varphi_{j,l}}) < \tfrac{1}{3} \op{minsep} \varphi_j .
\end{equation*}
That is, $\op{dist}(X, \op{spt} \mathbf{G}_{\varphi_j}) < \frac{1}{3} \op{minsep} \varphi_j$ for all $X \in \op{spt} T_j \cap \mathbf{C}_{(1+\gamma)\rho/2}(Z)$.  Hence for each $j \in \{1,2,\ldots,p\}$ 
\begin{equation*}
	T_j \llcorner \mathbf{C}^{\mathcal M}_{(1+\gamma)\rho/2}(Z) = \sum_{i=1}^{q_j} T_{j,l} 
\end{equation*}
where 
\begin{equation*}
	T_{j,l} = T_j \llcorner \Big\{ X \in \mathbf{C}^{\mathcal M}_{(1+\gamma)\rho/2}(Z) : 
		\op{dist}(X, \op{spt} \mathbf{G}_{\varphi_{j,l}}) < \tfrac{1}{3} \op{minsep} \varphi_j \Big\} 
\end{equation*}
and $T_{j,l}$ are locally area-minimizing rectifiable currents of $\mathbf{C}^{\mathcal M}_{(1+\gamma)\rho/2}(Z)$ with $(\partial T_{j,l}) \llcorner \mathbf{C}^{\mathcal M}_{(1+\gamma)\rho/2}(Z) = 0$.  (In fact, $T_{j,l} \llcorner (\mathbf{C}^{\mathcal M}_{(1+\gamma)\rho/2}(Z) \cap \mathbf{C}^{\mathcal M}_{\rho_k}(X_k)) = T_{k,j,l} \llcorner (\mathbf{C}^{\mathcal M}_{(1+\gamma)\rho/2}(Z) \cap \mathbf{C}^{\mathcal M}_{\rho_k}(X_k))$ for all $k,j,l$.)  By \eqref{relharm wellsep3 go condition}, \eqref{relharm separation0 go condition} holds true with $\frac{1+\gamma}{2}\rho,\frac{2\gamma}{1+\gamma}$ in place of $\rho,\gamma$.  Thus we can apply Lemma~\ref{relharm separation0 lemma} with $\frac{1+\gamma}{2}\rho,\frac{2\gamma}{1+\gamma},T_{j,l},\varphi_{j,l},0$ in place of $\rho,\gamma,T_j,\varphi_j,a_j$ to obtain \eqref{relharm wellsep3 concl}.
\end{proof}

\begin{proof}[Proof of Theorem~\ref{relharm separation3 thm}] 
We will prove Theorem~\ref{relharm separation3 thm} by the inductive argument from~\cite[Lemma~2.13]{KrumWicb} with some modification.  Let $\varsigma$ to be the minimum of $2\beta_2$ and its values in Lemma~\ref{relharm separation0 lemma} and Lemma~\ref{relharm wellsep3 lemma}.  In the case $p = s$, Theorem~\ref{relharm separation3 thm} follows from Lemma~\ref{relharm separation0 lemma}.  Suppose that $p_0, s_0$ are positive integers such that $1 \leq p_0 < s_0 \leq q$ and: 
\begin{enumerate}[itemsep=3mm,topsep=0mm]
	\item[(H1)]  Theorem~\ref{relharm separation3 thm} holds true whenever $p_0 < p \leq q$;
	\item[(H2)]  Theorem~\ref{relharm separation3 thm} holds true whenever $p = p_0$ and $p_0 \leq s < s_0$.
\end{enumerate}
As a slight change of notation, for each $j \in \{1,2,\ldots,p_0\}$ let's rewrite $\varphi_j$ as 
\begin{equation*}
	\varphi_j(X) = \sum_{l=1}^{s_j} m_{j,l} \llbracket \varphi_{j,l}(X) \rrbracket , 
\end{equation*}
where $s_j, m_{j,1}, m_{j,2},\ldots, m_{j,s_j}$ are positive integers such that $\sum_{j=1}^{p_0} s_j = s_0$ and $\sum_{l=1}^{s_j} m_{j,l} = q_j$ and $\varphi_{j,1},\varphi_{j,2},\ldots,\varphi_{j,s_j} \in C^2(B^{\mathcal M}_{\rho}(Z),\mathbb{R}^{n+m})$ are distinct single-valued functions such that \eqref{relharm grand hyp2}--\eqref{relharm separation3 normal}, and \eqref{relharm separation3 wellsep} hold true.  Assume that $s_1 > 1$ and 
\begin{equation*}
	\min_{1 \leq j \leq p_0} \op{minsep} \varphi_j = \inf_{B^{\mathcal M}_{\rho}(Z)} |\varphi_{1,s_1-1} - \varphi_{1,s_1}| .
\end{equation*}
Now define $\widetilde{\varphi}_1 : B^{\mathcal M}_{\rho}(Z) \rightarrow \mathcal{A}_{q_1}({\rm T}\mathcal{M}^{\perp})$ by collapsing $\varphi_{1,s_1}$ to $\varphi_{1,s_1-1}$ and then subtracting off the average of the resulting function; more precisely, let $\widetilde{\varphi}_1 = \sum_{l=1}^{\widetilde{s}_1} \widetilde{m}_{1,l} \llbracket \widetilde{\varphi}_{1,l} \rrbracket$, where $\widetilde{s}_1 = s_1 - 1$, $\widetilde{m}_{1,l} = m_{1,l}$ for all $l \in \{1,2,\ldots,s_1-2\}$, $\widetilde{m}_{1,\widetilde{s}_1} = m_{1,s_1-1} + m_{1,s_1}$, and $\widetilde{\varphi}_{1,l} = \varphi_{1,l} - \frac{m_{1,s_1}}{q_1} (\varphi_{1,s_1-1} - \varphi_{1,s_1})$ for all $l \in \{1,2,\ldots,\widetilde{s}_1\}$.  For all $j \in \{2,3,\ldots,p_0\}$, let $\widetilde{\varphi}_j = \varphi_j$ and write $\widetilde{\varphi}_j = \sum_{l=1}^{\widetilde{s}_j} \widetilde{m}_{j,l} \llbracket \widetilde{\varphi}_{j,l} \rrbracket$ where $\widetilde{s}_j = s_j$ and $\widetilde{m}_{j,l} = m_{j,l}$ and $\widetilde{\varphi}_{j,l} = \varphi_{j,l}$ for all $l \in \{1,2,\ldots,\widetilde{s}_j\}$.  Since $\widetilde{\varphi}_1$ has fewer distinct components than $\varphi_1$, for each $X \in B^{\mathcal M}_{\rho}(Z)$ there exists $l,l' \in \{1,2,\ldots,s_1\}$ and $k \in \{1,2,\ldots,\widetilde{s}_1\}$ such that 
\begin{equation*}
	|\varphi_{1,l}(X) - \widetilde{\varphi}_{1,k}(X)|^2 + |\varphi_{1,l'}(X) - \widetilde{\varphi}_{1,k}(X)|^2 
	\leq \mathcal{G}(\varphi_1(X),\widetilde{\varphi}_1(X))^2
\end{equation*}
and thus 
\begin{align*}
	\op{minsep} \varphi_1 \leq\,& |\varphi_{1,l}(X) - \varphi_{1,l'}(X)| 
	\leq |\varphi_{1,l}(X) - \widetilde{\varphi}_{1,k}(X)| + |\varphi_{1,l'}(X) - \widetilde{\varphi}_{1,k}(X)| 
	\\ \leq\,& 2 \mathcal{G}(\varphi_1(X),\widetilde{\varphi}_1(X)) .
\end{align*}
That is, 
\begin{equation}\label{relharm separation3 eqn1}
	\min_{1 \leq j \leq p_0} \op{minsep} \varphi_j \leq 2 \inf_{B^{\mathcal M}_{\rho}(Z)} \mathcal{G}(\varphi_1,\widetilde{\varphi}_1) .
\end{equation}
On the other hand, by \eqref{relharm dist to sep} for each $l \in \{1,2,\ldots,\widetilde{s}_1\}$ 
\begin{align}\label{relharm separation3 eqn2}
	\sup_{B^{\mathcal M}_{\rho}(Z)} |\widetilde{\varphi}_{1,l} - \varphi_{1,l}| 
	\leq\,& \sup_{B^{\mathcal M}_{\rho}(Z)} |\varphi_{1,s_1-1} - \varphi_{1,s_1}| 
	\\ \leq\,& (1+2\kappa) \inf_{B^{\mathcal M}_{\rho}(Z)} |\varphi_{1,s_1-1} - \varphi_{1,s_1}| \nonumber 
	\\=\,& (1+2\kappa) \min_{1 \leq j \leq p_0} \op{minsep} \varphi_j . \nonumber 
\end{align}
and similarly 
\begin{equation}\label{relharm separation3 eqn03}
	\sup_{B^{\mathcal M}_{\rho}(Z)} |\widetilde{\varphi}_{1,\widetilde{s}_1} - \varphi_{1,s_1}| 
	\leq (1+2\kappa) \min_{1 \leq j \leq p_0} \op{minsep} \varphi_j .  
\end{equation}

Let 
\begin{align*}
	\widetilde{H}^2 =\,& \frac{1}{\rho^{n+2}} \sum_{j=1}^{p_0} \int_{B^{\mathcal M}_{(1+\gamma)\rho/2}(Z)} \mathcal{G}(N_j,\widetilde{\varphi}_j)^2 
		+ \rho^2 \sum_{j=1}^{p_0} \sum_{l=1}^{\widetilde{s}_j} \|\Delta_{\mathcal M} \widetilde{\varphi}_{j,l}\|_{C^0(B^{\mathcal M}_{\rho}(Z))}^2  		
		\\&+ \left( \sum_{j=1}^{p_0} \sum_{l=1}^{\widetilde{s}_j} \|D\widetilde{\varphi}_{j,l}\|_{C^0(B^{\mathcal M}_{\rho}(Z))} + \rho^{\varsigma} \right) 
		\frac{1}{\rho^{n+2}} \int_{B^{\mathcal M}_{\rho}(Z)} |N|^2 . \nonumber 
\end{align*}
(as in \eqref{relharm excess defn3} with $(1+\gamma)/2,\widetilde{\varphi}_j$ in place of $\gamma,\varphi_j$).  We may assume that $\widetilde{H} > 0$ as otherwise arguing as in the the beginning of the proof of Lemma~\ref{relharm harmonic approx lemma} we have that $\op{spt} T_j \cap \mathbf{C}^{\mathcal M}_{\gamma\rho}(Z) \subseteq \op{spt}\mathbf{G}_{\varphi_j}$ for all $j \in \{1,2,\ldots,p_0\}$ and thus \eqref{relharm separation3 concl} holds true.  Notice that \eqref{relharm separation3 go condition} with $p_0,s_0$ in place of $p,s$ implies that \eqref{relharm separation3 go condition} holds true with $p_0,s_0-1,\frac{(1+\gamma)\rho}{2},\frac{1+3\gamma}{2+2\gamma}$ in place of $p,s,\rho,\gamma$.  Thus by (H2), we can apply Theorem~\ref{relharm separation3 thm} (with $\frac{(1+\gamma)\rho}{2},\frac{1+3\gamma}{2+2\gamma}$ in place of $\rho,\gamma$) to obtain 
\begin{equation}\label{relharm separation3 eqn3}
	\sup_{X \in \op{spt} T \cap \mathbf{C}^{\mathcal M}_{(1+3\gamma)\rho/4}(Z)} \op{dist}(X,\op{spt}\mathbf{G}_{\widetilde{\varphi}_j}) 
	< C_0 \widetilde{H} \rho 
\end{equation}
for some constant $C_0 = C_0(n,m,q,p_0,s_0,M_0,N_0,C_e,C_h,\gamma) \in (0,\infty)$.  Let $\beta = \beta(n,m,q,p_0,s_0,M_0,$ $N_0,C_e,C_h,\gamma) \in (0,1)$ to be later determined.  We may assume that 
\begin{equation}\label{relharm separation3 eqn4}
	E_{Z,(1+\gamma)\rho/2,\rho} \,\rho \leq \beta \min_{1 \leq j \leq p_0} \op{minsep} \varphi_j ; 
\end{equation}
otherwise, by the triangle inequality, \eqref{relharm separation3 eqn2}, and \eqref{relharm separation3 eqn03} 
\begin{equation*}
	\widetilde{H} \leq E_{Z,(1+\gamma)\rho/2,\rho} + \frac{C}{\rho} \min_{1 \leq j \leq p_0} \op{minsep}\varphi_j 
	\leq \left( 1 + \frac{C}{\beta} \right) E_{Z,(1+\gamma)\rho/2,\rho} ,
\end{equation*}
where $C = C(n,q,\kappa) \in (0,\infty)$ is a constant.  Hence by \eqref{relharm separation3 eqn3} and $\op{spt}\mathbf{G}_{\widetilde{\varphi}_j} \subseteq \op{spt}\mathbf{G}_{\varphi_j}$, \eqref{relharm separation3 concl} holds true for some constant $C = C(n,m,q,p_0,s_0,$ $M_0,N_0,C_e,C_h,\gamma,\beta) \in (0,\infty)$.

Fix any $j \in \{1,2,\ldots,p_0\}$ and let  
\begin{equation*}
	\{ X \in \mathbf{C}^{\mathcal M}_{(1+3\gamma)\rho/4}(Z) : \op{dist}(X,\op{spt}\mathbf{G}_{\widetilde{\varphi}_j}) < C_0 \widetilde{H} \rho \} 
	= \bigcup_{i=1}^{\widetilde{M}_j} \widetilde{U}_{j,i} 
\end{equation*}
for some collection of mutually disjoint connected open subsets $\{\widetilde{U}_{j,i}\}$ of $\mathbf{C}^{\mathcal M}_{(1+3\gamma)\rho/4}(Z)$ (where $C_0$ is as in \eqref{relharm separation3 eqn3}).  By \eqref{relharm separation3 eqn3}, we can express 
\begin{equation*}
	T_j = \sum_{i=1}^{\widetilde{M}_j} \widetilde{T}_{j,i} \quad\text{where}\quad \widetilde{T}_{j,i} = T_j \llcorner \widetilde{U}_{j,i} 
\end{equation*}
and $\widetilde{T}_{j,i}$ are locally area-minimizing rectifiable currents of $\mathbf{C}^{\mathcal M}_{(1+3\gamma)\rho/4}(Z)$ with $(\partial\widetilde{T}_{j,i}) \llcorner \mathbf{C}^{\mathcal M}_{(1+3\gamma)\rho/4}(Z) = 0$.  By the constancy theorem, $\mathbf{p}_{\#} \widetilde{T}_{j,i} = \widetilde{q}_{j,i} \llbracket B^{\mathcal M}_{(1+3\gamma)\rho/4}(Z) \rrbracket$ for some integers $\widetilde{q}_{j,i}$ with $\sum_{i=1}^{\widetilde{M}_j} \widetilde{q}_{j,i} = q_j$.  Since $\mathcal{M}$ is close to $P_Z$ in the $C^3$-topology in $\mathbf{C}_{\gamma\rho}(Z,P)$, by the continuity of push-forwards in the weak topology we have that $\pi_{P_Z \#} \widetilde{T}_{j,i} = \widetilde{q}_{j,i} \llbracket B_{\gamma\rho}(Z,P_Z) \rrbracket$.  Thus~\cite[Lemma~2.7]{KrumWicb}, $\widetilde{q}_{j,i} \geq 0$ for all $i \in \{1,2,\ldots,\widetilde{M}_j\}$.  Let $\widetilde{\varphi}_j = \sum_{l=1}^{\widetilde{M}_j} \widetilde{\psi}_{j,l}$ for multi-valued functions $\widetilde{\psi}_{j,l}$ on $B^{\mathcal M}_{(1+3\gamma)\rho/4}(Z)$ such that $\mathbf{G}_{\widetilde{\psi}_{j,i}} = \mathbf{G}_{\widetilde{\varphi}_j} \llcorner \widetilde{U}_{j,i}$.  Provided $C_0 \geq \sqrt{\frac{4}{\omega_n} \big(\frac{4}{1+3\gamma}\big)^n}$ (where $C_0$ is as in \eqref{relharm separation3 eqn3}), we can argue as in Remark~\ref{relharm separation1 rmk} that there exists a set $\Omega_j \subseteq \mathcal{K}_j \cap B^{\mathcal M}_{(1+3\gamma)\rho/4}(Z)$ such that $\mathcal{H}^n(\Omega_j) \geq \frac{1}{2} \omega_n \big(\frac{(1+3\gamma)\rho}{4}\big)^n$ and $\mathcal{G}(N_j,\varphi_j) < C_0 \widetilde{H} \rho$ on $\Omega_j$.  For each $X \in \Omega_j$, the slice $\langle \widetilde{T}_{j,i}, \mathbf{p}, X \rangle$ is the sum of precisely the values of $N_j(X)$ (counting multiplicity) which are within a distance $< C_0 \widetilde{H} \rho$ of the values of $\widetilde{\psi}_{j,i}(X)$.  Therefore, $\widetilde{q}_{j,i} > 0$ and $\widetilde{\psi}_{j,i}$ is a $\widetilde{q}_{j,i}$-valued function.  If $j = 1$, by the triangle inequality, \eqref{relharm separation3 eqn1}, and \eqref{relharm separation3 eqn4} 
\begin{align*}
	\Big(\frac{1+\gamma}{2}\Big)^{n/2} \min_{1 \leq j' \leq p_0} \op{minsep} \varphi_{j'} 
	\leq\,& 2 \rho^{-n/2} \|\mathcal{G}(\widetilde{\varphi}_1,\varphi_1)\|_{L^2(B^{\mathcal M}_{(1+\gamma)\rho/2}(Z))} 
	\leq E_{Z,(1+\gamma)\rho/2,\rho} \rho + \widetilde{H} \rho 
    \\ \leq\,& \beta \min_{1 \leq j' \leq p_0} \op{minsep} \varphi_{j'} + \widetilde{H} \rho 
\end{align*}
so that $\min_{1 \leq j' \leq p_0} \op{minsep} \varphi_{j'} \leq 2 \big(\frac{1+\gamma}{2}\big)^{-n/2} \widetilde{H}\rho$.  Thus by \eqref{relharm separation3 eqn2}, for each $l \in \{1,2,\ldots,s_1-1\}$ 
\begin{equation*}
	\sup_{B^{\mathcal M}_{\rho}(Z)} |\widetilde{\varphi}_{1,l} - \varphi_{1,l}| 
	\leq 2 (1+2\kappa) \,\Big(\frac{1+\gamma}{2}\Big)^{-n/2} \widetilde{H} \rho
\end{equation*}
Hence provided $C_0 > 2 (1+2\kappa) \,\big(\frac{1+\gamma}{2}\big)^{-n/2}$, there exists $i(l) \in \{1,2,\ldots,\widetilde{M}_1\}$ such that $\op{spt}\mathbf{G}_{\widetilde{\varphi}_{1,l}} \subset \widetilde{U}_{1,i(l)}$ and thus $\op{spt}\mathbf{G}_{\varphi_{1,l}} \subset \widetilde{U}_{1,i(l)}$.  By similar reasoning using \eqref{relharm separation3 eqn03}, $\op{spt}\mathbf{G}_{\widetilde{\varphi}_{1,s_1-1}} \subset \widetilde{U}_{1,i(s_1-1)}$ and thus $\op{spt}\mathbf{G}_{\varphi_{1,s_1}} \subset \widetilde{U}_{1,i(s_1-1)}$.  Therefore, there exists $\widetilde{q}_{1,i}$-valued functions $\psi_{1,i}$ of $B^{\mathcal M}_{(1+3\gamma)\rho/4}(Z)$ such that $\mathbf{G}_{\psi_{1,i}} = \mathbf{G}_{\varphi_1} \llcorner \widetilde{U}_{1,i}$.  If instead $j \geq 2$, set $\psi_{j,i} = \widetilde{\psi}_{j,i}$ for all $i \in \{1,2,\ldots,\widetilde{M}_j\}$ so that $\psi_{j,i}$ is a $\widetilde{q}_{j,i}$-valued function on $B^{\mathcal M}_{(1+3\gamma)\rho/4}(Z)$ such that $\mathbf{G}_{\psi_{j,i}} = \mathbf{G}_{\varphi_j} \llcorner \widetilde{U}_{j,i}$.  Notice that \eqref{relharm separation3 go condition} with $p_0,s_0$ in place of $p,s$ implies that \eqref{relharm separation3 go condition} holds true with $p_0,s_0-1,\frac{(1+3\gamma)\rho}{4},\frac{4\gamma}{1+3\gamma}$ in place of $p,s,\rho,\gamma$.  In the case that $\widetilde{M}_{j_1} > 1$ for some $j_1 \in \{1,2,\ldots,p_0\}$, by (H1) we can apply Theorem~\ref{relharm separation3 thm} with $\widetilde{T}_{j,i}$ and $\psi_{j,i}$ in place of $T_j$ and $\varphi_j$ (and with $\frac{(1+3\gamma)\rho}{4},\frac{4\gamma}{1+3\gamma}$ in place of $\rho,\gamma$) to deduce that \eqref{relharm separation3 concl} holds true.  

Now suppose that $\widetilde{M}_j = 1$ for all $j \in \{1,2,\ldots,p_0\}$.  By the triangle inequality and \eqref{relharm separation3 eqn4}, $\widetilde{H} \leq C(n,q) \min_{1 \leq j \leq p_0} \op{minsep} \varphi_j$.  For each $j \in \{1,2,\ldots,p_0\}$, recalling that $\widetilde{M}_j = 1$ and assuming $\widetilde{s}_j > 1$, we can relabel the components $\widetilde{\varphi}_{j,l}$ of $\widetilde{\varphi}_j$ so that for each $k \in \{2,3,\ldots,\widetilde{s}_j\}$ there exists $l(k) \in \{1,2,\ldots,k-1\}$ such that 
\begin{equation*}
	\inf_{B^{\mathcal M}_{(1+3\gamma)\rho/4}(Z)} |\widetilde{\varphi}_{j,k} - \widetilde{\varphi}_{j,l(k)}| \leq 2C_0 \widetilde{H} \rho .
\end{equation*}
By \eqref{relharm dist to sep} 
\begin{equation*}
	\sup_{B^{\mathcal M}_{\rho}(Z)} |\widetilde{\varphi}_{j,k} - \widetilde{\varphi}_{j,l(k)}| \leq 2 (1+2\kappa) \,C_0 \widetilde{H} \rho .
\end{equation*}
It follows that 
\begin{equation}\label{relharm separation3 eqn5}
	\max_{1 \leq j \leq p_0} \op{width}\widetilde{\varphi}_j \leq 2q_0 (1+2\kappa) \,C_0 \widetilde{H} \rho .
\end{equation}
(Note that in the case $\widetilde{s}_j = 1$, $\widetilde{\varphi}_j$ has exactly one component and thus $\widetilde{\varphi}_j$ has width zero.)  By \eqref{relharm separation3 eqn2}, \eqref{relharm separation3 eqn03}, \eqref{relharm separation3 eqn5}, and $\widetilde{H} \leq C(n,q) \min_{1 \leq j \leq p_0} \op{minsep} \varphi_j$ 
\begin{equation}\label{relharm separation3 eqn6}
	\max_{1 \leq j \leq p_0} \op{width} \varphi_j \leq C \min_{1 \leq j \leq p_0} \op{minsep} \varphi_j
\end{equation}
for some constant $C = C(n,m,q,p_0,s_0,M_0,N_0,C_e,C_h,\gamma,\kappa) \in (0,\infty)$.  Note that by \eqref{relharm separation3 go condition} with $p_0,s_0$ in place of $p,s$, \eqref{relharm wellsep3 go condition} holds true.  In light of \eqref{relharm wellsep3 go condition}, \eqref{relharm separation3 eqn4}, and \eqref{relharm separation3 eqn6}, we can apply Lemma~\ref{relharm wellsep3 lemma} to deduce that provided $\beta$ is sufficiently small, \eqref{relharm separation3 concl} holds true.
\end{proof}

\begin{proof}[Proof of Theorem~\ref{relharm separation main thm}] 
Theorem~\ref{relharm separation main thm} is simply Theorem~\ref{relharm separation3 thm} in the special case $p = 1$.  Note that \eqref{relharm normal approx concl1} and \eqref{relharm normal approx concl2} hold true with $(1+\gamma)/2$ in place of $\gamma$ and with $N_j = N$ on $B^{\mathcal M}_{(1+\gamma)\rho/2}(Z)$ and $\mathcal{K}_j = \mathcal{K} \cap B^{\mathcal M}_{(1+\gamma)\rho/2}(Z)$ for $j = 1$.
\end{proof}

\section{A graphical representation and blow-up procedure}\label{sec:graph rep} 

Let $T$ be an $n$-dimensional locally area minimizing rectifiable current of $\mathbf{C}_{11\sqrt{n}/2}(0)$ as in Hypothesis~\ref{center mnfd hypothesis} and $\mathcal{B}^{(\geq 2)}_q$ be as in Definition~\ref{center mnfd branch defn}.  Let $\mathcal{M}$ be the center manifold of $T$ as in Hypothesis~\ref{M and N hypothesis}.  Let $Z_0 \in \mathcal{B}^{(\geq 2)}_q$, $\rho_0 \in (0,1]$, and $P_{Z_0} = {\rm T}_{Z_0} \mathcal{M}$.  Recall that $P_{Z_0}$ is the common tangent plane to $\op{spt} T$ and $\mathcal{M}$ at $Z_0$.  Let $\varphi : P_{Z_0} \rightarrow \mathcal{A}_q(P_{Z_0}^{\perp})$ be a non-zero, homogeneous degree $\alpha$, harmonic function which is translation invariant along an $(n-2)$-dimensional linear subspace (as in Definition~\ref{Phi-p_defn} below).  In the case $\varphi$ has two or more components, assume that $T$ is significantly closer in $L^2$ to the graph of $\varphi$ than to the graph of any other homogeneous harmonic function $\varphi'$ having fewer components than $\varphi$ (as in Hypothesis~$(\star\star)$ in Section~\ref{sec:graph_rep_results} below).  Following an approach introduced in~\cite{KrumWicb}, we use Theorem~\ref{relharm separation main thm} to show that away from $\{0\} \times \mathbb{R}^{n-2}$, $T$ separates as the sum of locally area-minimizing rectifiable currents $T_j$ such that $T_j$ is uniformly close to the graph of a component $\varphi_j$ of $\varphi$.  Moreover, we use Almgren's Strong Lipschitz Approximation to approximate $T_j$ by the graph of a Lipschitz multi-valued function $v_j$ over a domain in the graph of $\varphi_j$.  These results will allow us to produce ``fine blow-ups'' of a sequence $(T^{(\nu)})$ of locally area-minimizing rectifiable currents relative to a sequence of non-zero, homogeneous degree $\alpha$, harmonic functions $(\varphi^{(\nu)})$.

\subsection{Notation and terminology}\label{sec:graphical_rep_notation}  Let $P$ be an $n$-dimensional plane in $\mathbb{R}^{n+m}$ and let $S$ be an $(n-2)$-dimensional linear subspace of $P$.  (We will often take $P = P_{Z_0} = {\rm T}_{Z_0} \mathcal{M}$ as the tangent plane to a center manifold $\mathcal{M}$ of a locally area minimizing rectifiable current at a branch point $Z_0 \in \mathcal{M}$.)  
We shall use the following notation, which is  in line with the notation in  ~\cite[Subsection~6.1]{KrumWic2} and in ~\cite[Subsection~5.1]{KrumWicc}.

\begin{definition}\label{Phi-p_defn}
Let $q \geq 2$ and $p \geq 0$ be integers.  Let $\alpha = \ell_0/q_0 \geq 2$ for relatively prime positive integers $\ell_0$ and $q_0$ with $q_0 \leq q$.  We define $\Phi_{\alpha,q,p}(P,S)$ to be the set of all average-free $q$-valued functions $\varphi : P \rightarrow \mathcal{A}_q(P^{\perp})$ of the form 
\begin{equation*}
	\varphi = \sum_{j=1}^J m_j \varphi_j 
\end{equation*}
where 
\begin{enumerate}[itemsep=3mm,topsep=0mm]
	\item[{\rm (i)}]  either $\varphi_j : P \rightarrow P^{\perp}$ is the single-valued zero function $\varphi_j(\,\cdot\,) = 0$ or $\varphi_j : P \rightarrow \mathcal{A}_{q_0}(P^{\perp})$ is the $q_0$-valued function given by 
\begin{equation*}
	\varphi_j(X) = \op{Re}(c_j (x_1+ix_2)^{\alpha}) 
\end{equation*}
for each $X \in P$, where $c_j \in \mathbb{C}^{n+m} \setminus \{0\}$ with $\op{Re}(c_j),\op{Im}(c_j) \in P^{\perp}$, $i = \sqrt{-1}$, and $(x_1,x_2,\ldots,x_n)$ are coordinates of $P$ with respect to an orthonormal basis for $P$ such that $S = \{ x_1 = x_2 = 0 \}$; 

	\item[{\rm (ii)}]  $\varphi_j \neq \varphi_{j'}$ whenever $j \neq j'$;
	
	\item[{\rm (iii)}]  $J$ and $m_j$ are positive integers such that $J = p+1$ if $\varphi_{j_1}(\,\cdot\,) = 0$ for some $j_1 \in \{1,2,\ldots,J\}$ and $J = p$ otherwise.
\end{enumerate}
\end{definition}

\begin{remark}{\rm (1)  Throughout we let $q_j = 1$ if $\varphi_j$ identically zero and $q_j = q_0$ otherwise.  That way, each component $\varphi_j$ is a $q_j$-valued function and the multiplicities $m_j$ satisfy $\sum_{j=1}^J q_j m_j = q$.

\noindent
(2)  We have that $q_0 \sum_{j=1}^p m_j = q$ in the case that $\varphi_j \not\equiv 0$ for all $j$ and $m_{j_1} + q_0 \sum_{j \neq j_1} m_j = q$ in the case that $\varphi_{j_1} \equiv 0$ for some (unique) $j_1$.  In particular, $q_0 p \leq q$ and thus $1 \leq p \leq \lfloor q/q_0 \rfloor$.

\noindent 
(3)  We do not assume that the functions $\varphi \in \Phi_{\alpha,q, p}(P,S)$ are locally energy minimizing.
} \end{remark}

\begin{definition}\label{Phi_defn}
Let $q$ be a positive integer.  Let $\alpha = \ell_0/q_0$ for relatively prime positive integers $\ell_0$ and $q_0$ with $q_0 \leq q$.  We let $\Phi_{\alpha,q}(P,S) = \bigcup_{p=0}^{\lfloor q/q_0 \rfloor} \Phi_{\alpha,q,p}(P,S)$.
\end{definition}

For $p \in \{0,1,\ldots,\lfloor q/q_0 \rfloor\}$, let $\varphi \in \Phi_{\alpha,q,p}(P,S)$ be as in Definition~\ref{Phi-p_defn}.  For each $j \in \{1,2,\ldots,J\}$, the component $\varphi_j$ is a $q_j$-valued function with integer multiplicity $m_j$ and the graph of the component $\varphi_j$ is a smoothly immersed submanifold of $\mathbb{R}^{n+m} \setminus \pi_P^{-1}(S)$.  Thus corresponding to any ball $B \subset P \setminus S$, there are single-valued harmonic functions $\varphi_{j,l} : B \rightarrow P^{\perp}$ such that 
\begin{equation}\label{varphi_localized}
	\varphi_{j}(X) = \sum_{l=1}^{q_j} \llbracket \varphi_{j,l}(X) \rrbracket 
\end{equation}
for each $X \in B$, where we recall that $q_j = 1$ if $\varphi_j \equiv 0$ and $q_j = q_0$ otherwise.

Let $T$ be an $n$-dimensional locally area minimizing rectifiable current of $\mathbf{C}_{11\sqrt{n}/2}(0)$ as in Hypothesis~\ref{center mnfd hypothesis} and $\mathcal{B}^{(\geq 2)}_q$ be as in Definition~\ref{center mnfd branch defn}.  Let $\mathcal{M}$ be the center manifold of $T$ as in Hypothesis~\ref{M and N hypothesis}.  Let $Z_0 \in \mathcal{B}^{(\geq 2)}_q$, $P_{Z_0} = {\rm T}_{Z_0}\mathcal{M}$, and $S$ be an $(n-2)$-dimensional linear subspace of $P_{Z_0}$.  Recall that $\exp_{Z_0} : B_3(0,P_{Z_0}) \rightarrow B^{\mathcal M}_1(Z_0)$ is the exponential map of $\mathcal{M}$ at $Z_0$.  Let $p \in \{0,1,\ldots,\lfloor q/q_0 \rfloor\}$ and let $\varphi \in \Phi_{\alpha,q,p}(P_{Z_0},S)$ be as in Definition~\ref{Phi-p_defn}.  Let $\Omega \subset P_{Z_0} \setminus S$ be an open set and for each $j \in \{1,2,\ldots,J\}$ let $v_j : \op{spt}\mathbf{G}_{\varphi_j \circ \exp_{Z_0}^{-1}} \cap \mathbf{p}^{-1}(\exp_{Z_0}^{-1}(\Omega)) \rightarrow \mathcal{A}_{m_j}(\mathbb{R}^{n+m})$ be Lipschitz $m_j$-valued functions.  Let $B \subset \Omega$ be an open ball and let $\varphi_{j,l}$ is as in \eqref{varphi_localized}.  For each $j \in \{1,2,\ldots,J\}$ and $l \in \{1, \ldots, q_j\}$ we define $v_{j,l} : B \rightarrow {\mathcal A}_{m_{j, k}}(\mathbb{R}^{n+m})$ by 
\begin{equation}\label{v_localized}
	v_{j,l}(X) = v_j(\exp_{Z_0}(X) + \varphi_{j,l}(X)) = \sum_{h=1}^{m_j} \llbracket v_{j,l,h}(X) \rrbracket
\end{equation}
for each $X \in B$, where $v_{j,l,h}(X) \in {\mathbb R}^{n+m}$ for each $h \in \{1,2,\ldots,m_j\}$.  For each $j \in \{1,2,\ldots,J\}$, we define a $q_j m_j$-valued function $N_j : \exp_{Z_0}^{-1}(\Omega) \rightarrow \mathcal{A}_{q_j m_j}(\mathbb{R}^{n+m})$ by 
\begin{equation} \label{N of v defn} 
	N_j(\exp_{Z_0}(X)) = \sum_{l=1}^{q_j} \sum_{h=1}^{m_j} \llbracket N_{j,l,h}(X) \rrbracket 
		\quad\text{where}\quad N_{j,l,h}(X) = \varphi_{j,l}(X) + v_{j,l,h}(X)
\end{equation}
for every ball $B \subset \Omega$ and every $X \in B$.  Hence $N = \sum_{j=1}^J N_j$ is a $q$-valued function associated with $v = (v_j)$.

\begin{definition}\label{v comp minimizing defn}
Let $Z_0 \in \mathcal{B}^{(\geq 2)}_q$, $P_{Z_0} = {\rm T}_{Z_0} \mathcal{M}$, and $S$ be an $(n-2)$-dimensional linear subspace of $P_{Z_0}$.  Let $p \in \{0,1,\ldots,\lfloor q/q_0 \rfloor\}$ and let $\varphi \in \Phi_{\alpha,q,p}(P_{Z_0},S)$ be as in Definition~\ref{Phi-p_defn}.  Let $\Omega \subset P_{Z_0} \setminus S$ be an open set and for each $j \in \{1,2,\ldots,J\}$ let $v_j : \op{spt}\mathbf{G}_{\varphi_j \circ \exp_{Z_0}^{-1}} \cap \mathbf{p}^{-1}(\exp_{Z_0}^{-1}(\Omega)) \rightarrow \mathcal{A}_{m_j}(\mathbb{R}^{n+m})$ be Lipschitz $m_j$-valued functions.  We say that $v = (v_j)$ is \emph{component-wise Dirichlet energy minimizing in $\Omega$} if for every ball $B \subset\subset \Omega$, each function $v_{j,l} : B \rightarrow \mathcal{A}_{m_j}(\mathbb{R}^{n+m})$ (as in \eqref{v_localized}) is in $W^{1,2}(B;\mathcal{A}_{m_j}(\mathbb{R}^{n+m}))$ and is locally Dirichlet energy minimizing in $B$.
\end{definition}

Let $T$ be an $n$-dimensional locally area minimizing rectifiable current of $\mathbf{C}_{11\sqrt{n}/2}(0)$ as in Hypothesis~\ref{center mnfd hypothesis} and $\mathcal{B}^{(\geq 2)}_q$ be as in Definition~\ref{center mnfd branch defn}.  Let $\mathcal{M}$ be the center manifold of $T$ as in Hypothesis~\ref{M and N hypothesis}.  Let $Z_0 \in \mathcal{B}^{(\geq 2)}_q$, $P_{Z_0} = {\rm T}_{Z_0}\mathcal{M}$, and $S$ be an $(n-2)$-dimensional linear subspace of $P_{Z_0}$.  Let $\Omega \subset P_{Z_0} \setminus S$ be a connected open set which is rotationally symmetric about $S$, i.e.~given $X,\widetilde{X} \in P_{Z_0} \cap S^{\perp}$ and $Y \in S$, if $X+Y \in \Omega$ and $|X| = |\widetilde{X}|$ then $\widetilde{X}+Y \in \Omega$.  Let $p \in \{0,1,\ldots,\lfloor q/q_0 \rfloor\}$ and let $\varphi \in \Phi_{\alpha,q,p}(P_{Z_0},S)$ be as in Definition~\ref{Phi-p_defn} (with $P = P_{Z_0}$).  Assume that for $\varepsilon = \varepsilon(n,m,q,\Omega) \in (0,\infty)$ suitably small, $\|\varphi\|_{L^2(\Omega)} < \varepsilon$.  Suppose that  
\begin{equation*}
	\sup_{X \in \op{spt} T \cap \mathbf{p}^{-1} (\exp_{Z_0}(\Omega))} \op{dist}(X,\op{spt} \mathbf{G}_{\varphi \circ \exp_{Z_0}^{-1}}) 
	< \tfrac{1}{3} \inf_{X \in \Omega} \op{sep} \varphi(X) ,
\end{equation*}
where the separation $\op{sep} \varphi(X)$ is as in Definition~\ref{a sep width defn}.  Then there exists (unique) $n$-dimensional locally area-minimizing rectifiable currents $T_j$ of $\mathbf{p}^{-1} (\exp_{Z_0}(\Omega))$ such that 
\begin{gather*}
	T \llcorner \mathbf{p}^{-1} (\exp_{Z_0}(\Omega)) = \sum_{j=1}^J T_j , \quad\quad
	(\partial T_j) \llcorner \mathbf{p}^{-1} (\exp_{Z_0}(\Omega)) = 0, \\
	\sup_{X \in \op{spt} T_j} \op{dist}(X, \op{spt} \mathbf{G}_{\varphi_j \circ \exp_{Z_0}^{-1}}) < \tfrac{1}{3} \inf_{X \in \Omega} \op{sep} \varphi(X) . \nonumber 
\end{gather*}
Moreover, corresponding to any ball $B \subset \Omega$ there are $n$-dimensional locally area-minimizing rectifiable currents $T_{j,l}$ of $\mathbf{p}^{-1} (\exp_{Z_0}^{-1}(B))$ such that 
\begin{gather}\label{T_separated eqn4}
	T_j \llcorner \mathbf{p}^{-1} (\exp_{Z_0}(B)) = \sum_{l=1}^{q_j} T_{j,l} , \quad\quad
	(\partial T_{j,l}) \llcorner \mathbf{p}^{-1} (\exp_{Z_0}(B)) = 0, \\
	\sup_{X \in \op{spt} T_{j,l}} \op{dist}(X, \op{spt} \mathbf{G}_{\varphi_{j,l} \circ \exp_{Z_0}^{-1}}) < \tfrac{1}{3} \inf_{X \in \Omega} \op{sep} \varphi(X) , \nonumber 
\end{gather}
where $\varphi_{j,l}$ is as in \eqref{varphi_localized}.  By the constancy theorem, for each $j \in \{1,\ldots,J\}$ there exists an integer $\widehat{m}_j$ such that 
\begin{equation*}
	\mathbf{p}_{\#} T_{j,l} = \widehat{m}_j \llbracket \exp_{Z_0}(B) \rrbracket 
\end{equation*}
for every ball $B \subset \Omega$ and $l \in \{1,\ldots,q_j\}$.  Since $\Omega$ is rotationally symmetric about $S$, $\op{graph} \varphi_j |_{\Omega}$ is connected and thus it follows that $\widehat{m}_j$ is independent of the ball $B$ and index $l$.  In particular, 
\begin{equation*}
	\mathbf{p}_{\#} T_j = q_j \widehat{m}_j \llbracket \exp_{Z_0}(\Omega) \rrbracket . 
\end{equation*}

\subsection{Graphical representation of area-minimizing currents}\label{sec:graph_rep_results}  Let $\varsigma,\varepsilon \in (0,1)$.  Let $T$, $\mathcal{M}$, and $N$ be as in Hypothesis~\ref{M and N hypothesis}.  Let $Z_0 \in \mathcal{B}^{(\geq 2)}_q$, $\rho_0 \in (0,1]$, $P_{Z_0} = {\rm T}_{Z_0} \mathcal{M}$, and $S$ be an $(n-2)$-dimensional linear subspace of $P_{Z_0}$.  Let $p \in \{0,1,\ldots,\lfloor q/q_0 \rfloor\}$ and let $\varphi \in \Phi_{\alpha,q,p}(P_{Z_0},S)$ be as in Definition~\ref{Phi-p_defn} (with $P = P_{Z_0}$).  We define 
\begin{equation}\label{graphical E defn}
	E_{\varsigma,\varepsilon}(N,\varphi,B^{\mathcal M}_{\rho_0}(Z_0))^2 = \frac{1}{\omega_n \rho_0^{n+2}} \int_{B^{\mathcal M}_{\rho_0}(Z_0)} 
		\mathcal{G}(N, \varphi \circ \exp_{Z_0}^{-1})^2 + \frac{1}{\varepsilon \omega_n \rho_0^{n+2-\varsigma}} \int_{B^{\mathcal M}_{\rho_0}(Z_0)} |N|^2 . 
\end{equation}
When $\varepsilon = 1$, we let $E_{\varsigma}(N,\varphi,B^{\mathcal M}_{\rho_0}(Z_0)) = E_{\varsigma,1}(N,\varphi,B^{\mathcal M}_{\rho_0}(Z_0))$

Observe that we can express each $X \in B^{\mathcal M}_1(Z_0)$ as $X = \exp_{Z_0}(x+y)$ for (unique) $x \in P_{Z_0} \cap S^{\perp}$ and $y \in S$.  Bearing this in mind, we define $r : B^{\mathcal M}_1(Z_0) \rightarrow \mathbb{R}$ by 
\begin{equation}\label{dist to S defn}
    r(X) = \op{dist}(\exp_{Z_0}^{-1}(X),S) = |x|
\end{equation}
for each $X \in B^{\mathcal M}_1(Z_0)$.  We extend $r$ to a function on $\mathbf{C}^{\mathcal M}_1(Z_0)$ defined by $r(X) = r(\mathbf{p}(X))$ for each $X \in \mathbf{C}^{\mathcal M}_1(Z_0)$.

Let $\alpha = \ell_0/q_0 \geq 2$ for some relatively prime positive integers $\ell_0$ and $q_0 \leq q$.  Let $M_0,N_0,C_e,C_h,\varepsilon_2$ as in Hypothesis~\ref{DeLSpa constants hypothesis} and $\varepsilon_0, \beta_0 \in (0,1)$ be a suitably small constants.  In Theorem~\ref{graphical_thm} and subsequent results throughout the paper, we shall assume:

\noindent\textbf{Hypothesis~$\bs{(\star)}$.}  \textit{Let $\varsigma = \varsigma(m,n,q) \in (0,1)$ be as in Theorem~\ref{relharm separation main thm}.  Let $T$, $\mathcal{M}$, $N$, and $\mathcal{K}$ be as in Hypothesis~\ref{M and N hypothesis} and let $\mathcal{B}^{(\geq 2)}_q$ be as in Definition~\ref{center mnfd branch defn}.  Let $Z_0 \in \mathcal{B}^{(\geq 2)}_q$ and $\rho_0 \in (0,1]$ such that $\mathbf{m}_0 \rho_0 < \varepsilon_0$.  Let $P_{Z_0} = {\rm T}_{Z_0} \mathcal{M}$ and $S$ be an $(n-2)$-dimensional linear subspace of $P_{Z_0}$.  Let $p \in \{1,2,\ldots, \lfloor q/q_0 \rfloor \}$ and let $\varphi \in \Phi_{\alpha,q,p}(P_{Z_0},S)$ be as in Definition~\ref{Phi-p_defn} (with $P = P_{Z_0}$).}

\noindent\textbf{Hypothesis~$\bs{(\star\star)}$.}  For $p \in \{1,2,\ldots, \lfloor q/q_0 \rfloor \}$, $T$, and $\varphi \in \Phi_{\alpha,q,p}(P_{Z_0},S)$ as in Hypothesis~$(\star)$ 
\begin{equation*}
	E_{\varsigma/2}(N,\varphi,B^{\mathcal M}_{\rho_0}(Z_0)) 
		\leq \beta_0 \inf_{\varphi' \in \bigcup_{p'=0}^{p-1} \Phi_{\alpha,q,p'}(P_{Z_0},S)} E_{\varsigma/2}(N,\varphi',B^{\mathcal M}_{\rho_0}(Z_0)) .
\end{equation*}

\begin{theorem}\label{graphical_thm}
Let $\alpha = \ell_0/q_0 \geq 2$ for some relatively prime positive integers $\ell_0$ and $q_0 \leq q$.  Given $p \in \{1,2,\ldots, \lfloor q/q_0 \rfloor \}$ and each $\tau \in (0,1/2]$ there exists $M_0,N_0,C_e,C_h,\varepsilon_2$ as in Hypothesis~\ref{DeLSpa constants hypothesis} (independent of $p,\alpha,\tau$) and there exists $\varepsilon_0,\beta_0 \in (0,1)$ depending only on $n,m,q,p,M_0,N_0,C_e,C_h,\alpha,\tau$ such that if Hypothesis~$(\star)$ and Hypothesis~$(\star\star)$ hold true, then: 
\begin{enumerate}[itemsep=3mm,topsep=0mm]
	\item[(A)]  
	\begin{equation}\label{graphical conclA1} 
		\inf_{x \in \mathbb{S}^{n+m-1} \cap P_{Z_0} \cap S^{\perp}} \op{sep} \varphi(x) 
			\geq c \rho_0 \inf_{\varphi' \in \bigcup_{p'=0}^{p-1} \Phi_{\alpha,q,p'}(P_{Z_0},S)} E_{\varsigma/2}(N,\varphi',B^{\mathcal M}_{\rho_0}(Z_0)) , 
	\end{equation}
	where the separation $\op{sep} \varphi(x)$ is as in Definition~\ref{a sep width defn} and $c = c(n,m,q,p,M_0,N_0,C_e,C_h,\alpha) > 0$ is a constant, and 
	\begin{equation}\label{graphical conclA2}
		\|\varphi_{j,l} - \varphi_{j',l'}\|_{L^2(B_{1/2}(x_0,P_{Z_0}))} \leq C \inf_{B_{1/2}(x_0,P_{Z_0})} |\varphi_{j,l} - \varphi_{j',l'}|
	\end{equation}
	for all $x_0 \in \mathbb{S}^{n+m-1} \cap P_{Z_0} \cap S^{\perp}$, $j,j' \in \{1,\ldots,J\}$, $l \in \{1,\ldots,q_j\}$, and $l' \in \{1,\ldots,q_{j'}\}$, where $\varphi_{j,l}$ and $\varphi_{j',l'}$ are as in \eqref{varphi_localized} with $B = B_{1/2}(x_0,P_{Z_0})$ and where $C = C(n,m,q,p,M_0,N_0,C_e,$ $C_h,\alpha) \in (0,\infty)$ is a constant; 

	\item[(B)]  there exist $n$-dimensional locally area-minimizing rectifiable currents $T_j$ of $\mathbf{C}^{\mathcal M}_{15\rho_0/16}(Z_0) \cap \{r > \tau\rho_0/2\}$ such that 
	\begin{gather*}
		T \llcorner \mathbf{C}^{\mathcal M}_{15\rho_0/16}(Z_0) \cap \{r > \tau\rho_0/2\} = \sum_{j=1}^J T_j , \\
		(\partial T_j) \llcorner \mathbf{C}^{\mathcal M}_{15\rho_0/16}(Z_0) \cap \{r > \tau\rho_0/2\} = 0, \nonumber \\
		\mathbf{p}_{\#} T_j = q_j m_j \llbracket B^{\mathcal M}_{15\rho_0/16}(Z_0) \cap \{r > \tau\rho_0/2\} \rrbracket , \nonumber \\ 
		\sup_{X \in \op{spt} T_j \cap \{r > \sigma\rho_0\}} \op{dist}(X, \op{spt} \mathbf{G}_{\varphi_j \circ \exp_{Z_0}^{-1}}) 
			\leq C_{\sigma} \rho_0 \,E_{\varsigma/2}(N,\varphi,B^{\mathcal M}_{\rho_0}(Z_0)) \text{ for all } \sigma \in [\tau/2,1/2], \nonumber 
	\end{gather*}
	where $C_{\sigma} = C_{\sigma}(n,m,q,p,M_0,N_0,C_e,C_h,\alpha,$ $\sigma) \in (0,\infty)$ is a constant;

	\item[(C)]  for each $j \in \{1,2,\ldots,J\}$ there exists a Lipschitz $m_j$-valued function $v_j : \op{spt} \mathbf{G}_{\varphi_j \circ \exp_{Z_0}^{-1}} \cap \mathbf{C}^{\mathcal M}_{7\rho_0/8}(Z_0) \cap \{r > \tau\rho_0\} \rightarrow \mathcal{A}_{m_j}(\mathbb{R}^{n+m})$ and a closed set $\mathcal{K}_j \subseteq B^{\mathcal M}_{7\rho_0/8}(Z_0) \cap \{r > \tau\rho_0\}$ such that, with $N_{j}$ as in \eqref{N of v defn}, we have: 
	\begin{equation*}
		T_j \llcorner \mathbf{p}^{-1}(\mathcal{K}_j) = \mathbf{G}_{N_j} \llcorner \mathbf{p}^{-1}(\mathcal{K}_j), 
	\end{equation*}
	for each $L \in \mathscr{W}$ such that the Whitney region $\mathcal{L}$ associated with $\mathcal{L}$ (as in \eqref{Whitney region defn} satisfies $\mathcal{L} \cap B^{\mathcal M}_{7\rho_0/8}(Z_0) \cap \{r > \tau\rho_0\} \neq \emptyset$ 
	\begin{gather*}
		\sup_{\mathcal{L} \cap B^{\mathcal M}_{7\rho_0/8}(Z_0) \cap \{r > \tau\rho_0\}} |N_j| 
			\leq C \mathbf{m}_0^{1/(2n)} \ell(L)^{1+\beta_2} , \\ 
		\sup_{\mathcal{L} \cap B^{\mathcal M}_{7\rho_0/8}(Z_0) \cap \{r > \tau\rho_0\}} |DN_j| 
			\leq C \mathbf{m}_0^{\gamma_2} \ell(L)^{\gamma_2} , \\
		\mathcal{H}^n(\mathcal{L} \cap B^{\mathcal M}_{7\rho_0/8}(Z_0) \cap \{r > \tau\rho_0\} \setminus \mathcal{K}_j) 
			\leq C \mathbf{m}_0^{1+\gamma_2} \ell(L)^{n+2+\gamma_2} , \\
		\|T_j\|(\mathbf{p}^{-1}(\mathcal{L} \cap B^{\mathcal M}_{7\rho_0/8}(Z_0) \cap \{r > \tau\rho_0\} \setminus \mathcal{K}_j)) 
			\leq C \mathbf{m}_0^{1+\gamma_2} \ell(L)^{n+2+\gamma_2} , \\
		\int_{\mathcal{L} \cap B^{\mathcal M}_{7\rho_0/8}(Z_0) \cap \{r > \tau\rho_0\}} |DN_j|^2 
			\leq C \mathbf{m}_0 \ell(L)^{n+2-2\delta_2} , 
	\end{gather*}
	and 
	\begin{align}\label{graphical conclC2}	
		&\frac{1}{\rho_0^2} \sup_{B^{\mathcal M}_{7\rho_0/8}(Z_0) \cap \{r > \sigma\rho_0\}} \sum_{j=1}^J \sum_{l=1}^{q_j} |v_{j,l}|^2 
			\\&\hspace{10mm} + \frac{1}{\rho_0^n} \int_{B^{\mathcal M}_{7\rho_0/8}(Z_0) \cap \{r > \sigma\rho_0\}} 
			\sum_{j=1}^J \sum_{l=1}^{q_j} |Dv_{j,l}|^2 
		\leq C_{\sigma} E_{\varsigma/2}(N,\varphi,B^{\mathcal M}_{7\rho_0/8}(Z_0))^2 , \nonumber 
	\end{align}
	for all $\sigma \in [\tau,7/8)$, where $v_{j,l}$ is as in \eqref{v_localized}, and $C = C(n,m,q,p,M_0,N_0,$ $C_e,C_h,\alpha) \in (0,\infty)$ and $C_{\sigma} = C_{\sigma}(n,m,q,p,M_0,N_0,C_e,C_h,\alpha,\sigma) \in (0,\infty)$ is a constant. 
\end{enumerate}
\end{theorem} 

\begin{remark}\label{graphical beta0 rmk}{\rm 
A close examination of the use of Hypothesis~$(\star\star)$ in the proof of Theorem~\ref{graphical_thm} below shows that we may take $\beta_0 = \beta_{\star} \tau^{n/2+\alpha}$ for some $\beta_{\star} = \beta_{\star}(n,m,q,p,M_0,N_0,C_e,C_h,\alpha) \in (0,1)$ (independent of $\tau$).
}\end{remark}

In Lemma~\ref{graphical_lemma} and subsequently, we shall use the following notation.  Let $\varsigma \in (0,1)$.  Let $T$ be an $n$-dimensional locally area minimizing rectifiable current of $\mathbf{C}_{11\sqrt{n}/2}(0)$ as in Hypothesis~\ref{center mnfd hypothesis}.  Let $\mathcal{M}$ be the center manifold of $T$ and let $N$ be the normal approximation map of $T$ as in Hypothesis~\ref{M and N hypothesis}.  Fix $Z_0 \in \mathcal{B}^{(\geq 2)}_q$, $\rho_0 \in (0,1]$, $P_{Z_0} = {\rm T}_{Z_0} \mathcal{M}$, and an $(n-2)$-dimensional linear subspace $S$ of $P_{Z_0}$.  Given $\kappa \in (0,1]$, $\zeta \in S$, and $\rho > 0$ with $\rho^2 + |\zeta|^2 \leq (15\rho_0/16)^2$, we define the annulus $A_{\rho,\kappa}(\zeta) \subset P_{Z_0}$ by 
\begin{equation} \label{annuli_defn}
	A_{\rho,\kappa}(\zeta) = \{ x+y : x  \in P_{Z_0} \cap S^{\perp} \text{, } y \in S \text{, } (|x| - \rho)^2 + |y-\zeta|^2 < (\kappa\rho/64)^2 \} .  
\end{equation} 
Let $A^{\mathcal M}_{\rho,\kappa}(\zeta) = \exp_{Z_0}(A_{\rho,\kappa}(\zeta)) \subset \mathcal{M}$ and $\mathbf{A}^{\mathcal M}_{\rho,\kappa}(\zeta) = \{ X \in \mathbf{U} : \mathbf{p}(X) \in A^{\mathcal M}_{\rho,\kappa}(\zeta) \}$.  For each $\varsigma \in (0,1]$ and $\varphi \in \Phi_{\alpha,q,p}(P_{Z_0},S)$ we define 
\begin{equation}\label{graphical annuli E defn}
	E_{\varsigma}(N,\varphi,A^{\mathcal M}_{\rho,\kappa}(\zeta))^2 
		= \frac{1}{\rho^{n+2}} \int_{A^{\mathcal M}_{\rho,\kappa}(\zeta)} \mathcal{G}(N,\varphi \circ \exp_{Z_0}^{-1})^2 
		+ \frac{\rho_0^{\varsigma}}{\rho^{n+2}} \int_{A^{\mathcal M}_{\rho,\kappa}(\zeta)} |N|^2 .
\end{equation}

\begin{lemma}\label{graphical_lemma} 
Let $\alpha = \ell_0/q_0 \geq 2$ for some relatively prime positive integers $\ell_0$ and $q_0 \leq q$.  Let $M_0,N_0,C_e,C_h$ be as in Hypothesis~\ref{DeLSpa constants hypothesis}.  For all $p \in \{0,1,2,\ldots, \lfloor q/q_0 \rfloor \}$ and $\kappa \in (0,1)$ there exists $\varepsilon_2 = \varepsilon_2(n,m,q,M_0,N_0,C_e,C_h) \in (0,1)$ and $\overline{\varepsilon}, \overline{\beta} \in (0,1)$ depending only on $n,m,q,p,M_0,N_0,C_e,C_h,$ $\alpha,\kappa$ such that such that if Hypothesis~$(\star)$ holds true with $\overline{\varepsilon}$ in place of $\varepsilon_0$ and if $\zeta \in S$ and $\rho > 0$ such that $\rho^2 + |\zeta|^2 \leq (15\rho_0/16)^2$,  
\begin{equation}\label{graphical annuli go condition}
	\sup_{X \in A^{\mathcal M}_{\rho,(3+\kappa)/4}(\zeta)} \overline{R}(X) < 2^{-2q-15} (1-\kappa) \,\rho , 
\end{equation}
and either 
\begin{enumerate}[itemsep=3mm,topsep=0mm]
	\item[(i)]  $p = 0$ or 
	\item[(ii)]  $p > 0$ and 
	\begin{equation}\label{graphical annuli hyp} 
		E_{\varsigma/2}(N,\varphi,A^{\mathcal M}_{\rho,1}(\zeta)) \leq \overline{\beta} \inf_{\varphi' \in \bigcup_{p'=0}^{p-1} \Phi_{\alpha,q,p'}(P_{Z_0},S)} 
			E_{\varsigma/2}(N,\varphi',A^{\mathcal M}_{\rho,1}(\zeta)) .
	\end{equation}
\end{enumerate}
Then:
\begin{enumerate}[itemsep=3mm,topsep=0mm]
	\item[(A)]  when $p > 0$ 
	\begin{equation}\label{graphical annuli conclA1} 
		\inf_{x \in \mathbb{S}^{n+m-1} \cap P_{Z_0} \cap S^{\perp}} \op{sep} \varphi(x) 
			\geq c \rho \inf_{\varphi' \in \bigcup_{p'=0}^{p-1} \Phi_{\alpha,q,p'}(P_{Z_0},S)} E_{\varsigma/2}(N,\varphi',A^{\mathcal M}_{\rho,1}(\zeta)) ,
	\end{equation}
	where the separation $\op{sep} \varphi(x)$ is as in Definition~\ref{a sep width defn} and $c = c(n,m,q,p,M_0,N_0,C_e,C_h,\alpha) > 0$ is a constant, and 
	\begin{equation}\label{graphical annuli conclA2} 
		\|\varphi_{j,l} - \varphi_{j',l'}\|_{L^2(B_{1/2}(x_0,P_{Z_0}))} \leq C \inf_{B_{1/2}(x_0,P_{Z_0})} |\varphi_{j,l} - \varphi_{j',l'}|
	\end{equation}
	for all $x_0 \in \mathbb{S}^{n+m-1} \cap P_{Z_0} \cap S^{\perp}$, $j,j' \in \{1,\ldots,p\}$, $l \in \{1,\ldots,q_j\}$, and $l' \in \{1,\ldots,q_{j'}\}$, where $\varphi_{j,l}$ and $\varphi_{j',l'}$ are as in \eqref{varphi_localized} with $B = B_{1/2}(x_0,P_{Z_0})$ and where $C = C(n,m,q,p,M_0,N_0,C_e,$ $C_h,\alpha) \in (0,\infty)$ is a constant;
	
	\item[(B)]  there exists $n$-dimensional locally area-minimizing rectifiable currents $T_j$ of $\mathbf{A}^{\mathcal M}_{\rho,(1+\kappa)/2}(\zeta)$ such that 
	\begin{gather}
		\label{graphical annuli conclB1} T \llcorner \mathbf{A}^{\mathcal M}_{\rho,(1+\kappa)/2}(\zeta) = \sum_{j=1}^J T_j , \quad 
		(\partial T_j) \llcorner \mathbf{A}^{\mathcal M}_{\rho,(1+\kappa)/2}(\zeta) = 0, \\ 
		\label{graphical annuli conclB2} \mathbf{p}_{\#} T_j = q_j m_j \llbracket A^{\mathcal M}_{\rho,(1+\kappa)/2}(\zeta) \rrbracket , \\
		\label{graphical annuli conclB3} \sup_{X \in \op{spt} T_j} \op{dist}(X, \op{spt} \mathbf{G}_{\varphi_j \circ \exp_{Z_0}^{-1}}) 
			\leq C \rho \,E_{\varsigma/2}(N,\varphi,A^{\mathcal M}_{\rho,1}(\zeta)) ,  
	\end{gather}
	where $C = C(n,m,q,p,\alpha,M_0,N_0,C_e,$ $C_h,\kappa) \in (0,\infty)$ is a constant;

	\item[(C)]  for each $j \in \{1,2,\ldots,J\}$ there exists a Lipschitz $\widehat{m}_j$-valued function $v_j : \op{spt} \mathbf{G}_{\varphi_j \circ \exp_{Z_0}^{-1}} \cap \mathbf{A}^{\mathcal M}_{\rho,\kappa}(\zeta) \rightarrow \mathcal{A}_{\widehat{m}_j}(\mathbb{R}^{n+m})$ and a closed set $\mathcal{K}_j \subseteq A^{\mathcal M}_{\rho,\kappa}(\zeta)$ such that, with $N_{j}$ as in \eqref{N of v defn}, we have:
	\begin{equation}\label{graphical annuli conclC1}
		T_j \llcorner \mathbf{p}^{-1}(\mathcal{K}_j) = \mathbf{G}_{N_j} \llcorner \mathbf{p}^{-1}(\mathcal{K}_j) , 
	\end{equation}
	for each $L \in \mathscr{W}$ such that the Whitney region $\mathcal{L}$ associated with $\mathcal{L}$ (as in \eqref{Whitney region defn} satisfies $\mathcal{L} \cap A^{\mathcal M}_{\rho,\kappa}(\zeta) \neq \emptyset$ 
	\begin{gather}\label{graphical annuli conclC2}
		\sup_{A^{\mathcal M}_{\rho,\kappa}(\zeta) \cap \mathcal{L}} |N_j| \leq C \mathbf{m}_0^{1/(2n)} \ell(L)^{1+\beta_2}, \quad 
		\sup_{A^{\mathcal M}_{\rho,\kappa}(\zeta) \cap \mathcal{L}} |DN_j| \leq C \mathbf{m}_0^{\gamma_2} \ell(L)^{\gamma_2} , \\
		\mathcal{H}^n(A^{\mathcal M}_{\rho,\kappa}(\zeta) \cap \mathcal{L} \setminus \mathcal{K}_j) 
			+ \|T_j\|(\mathbf{p}^{-1}(A^{\mathcal M}_{\rho,\kappa}(\zeta) \cap \mathcal{L} \setminus \mathcal{K}_j)) 
			\leq C \mathbf{m}_0^{1+\gamma_2} \ell(L)^{n+2+\gamma_2}, \nonumber\\
		\int_{A^{\mathcal M}_{\rho,\kappa}(\zeta) \cap \mathcal{L}} |DN_j|^2 \leq C \mathbf{m}_0 \ell(L)^{n+2-2\delta_2} , \nonumber
	\end{gather}
	and 
	\begin{equation}\label{graphical annuli conclC3}
		\frac{1}{\rho^2} \sup_{A^{\mathcal M}_{\rho,\kappa}(\zeta)} \sum_{j=1}^J \sum_{l=1}^{q_j} |v_{j,l}|^2 
			+ \frac{1}{\rho^n} \int_{A^{\mathcal M}_{\rho,\kappa}(\zeta)} \sum_{j=1}^J \sum_{l=1}^{q_j} |Dv_{j,l}|^2  
		\leq C E_{\varsigma/2}(N,\varphi,A^{\mathcal M}_{\rho,1}(\zeta))^2 , 
	\end{equation}
	where $v_j$ is as in \eqref{v_localized} and $C = C(n,m,q,p,M_0,N_0,C_e,C_h,\alpha,\kappa) \in (0,\infty)$ is a constant.
\end{enumerate}
\end{lemma}

\begin{remark}\label{graphical rmk}{\rm 
Let $\alpha = \ell_0/q_0 \geq 2$, where $\ell_0$ and $q_0 \leq q$ are relatively prime positive integers.  Let $M_0,N_0,C_e,C_h$ be as in Hypothesis~\ref{DeLSpa constants hypothesis} and let $\varepsilon_2 = \varepsilon_2(n,m,q,M_0,N_0,C_e,C_h) \in (0,1)$ and $\overline{\varepsilon} = \overline{\varepsilon}(n,m,q,M_0,N_0,C_e,C_h,\alpha) \in (0,1)$ be suitably small.  Assume that Hypothesis~$(\star)$ holds true with $\overline{\varepsilon}$ in place of $\varepsilon_0$.  Consider $A^{\mathcal M}_{\rho,1}(\zeta)$, where $\zeta \in S$ and $\rho > 0$ such that $\rho^2 + |\zeta|^2 \leq (15\rho_0/16)^2$.  (Note that the same conclusions hold true when we replace $\rho$, $A^{\mathcal M}_{\rho,1}(\zeta)$, and \eqref{graphical annuli hyp} with $\rho_0$, $B^{\mathcal M}_{\rho_0}(Z_0)$, and Hypothesis~$(\star\star)$ respectively.)

\noindent (1)  Suppose that for some constant $\overline{\beta} \in (0,1)$ such that \eqref{graphical annuli hyp} holds true.  By the compactness of $\Phi_{\alpha,q}(P_{Z_0},S)$, there exists (one or more) $\widetilde{p} \in \{0,1,2,\ldots,p-1\}$ and $\widetilde{\varphi} \in \Phi_{\alpha,q,\widetilde{p}}(P_{Z_0},S)$ such that 
\begin{equation*} 
	E_{\varsigma/2}(N,\widetilde{\varphi},A^{\mathcal M}_{\rho,1}(\zeta)) 
	= \inf_{\varphi' \in \bigcup_{p'=0}^{p-1} \Phi_{\alpha,q,p'}(P_{Z_0},S)} E_{\varsigma/2}(N,\varphi',A^{\mathcal M}_{\rho,1}(\zeta)) . 
\end{equation*}
By the triangle inequality 
\begin{equation}\label{graphical rmk eqn1} 
	\rho^{\alpha-1} \|\widetilde{\varphi}\|_{L^2(B_1(0))} \leq C E_{\varsigma/2}(N,\widetilde{\varphi},A^{\mathcal M}_{\rho,1}(\zeta)) + C \widehat{E} \leq 2C \widehat{E} , 
\end{equation}
where $\widehat{E} = \rho^{-n/2-1} \|N\|_{L^2(A^{\mathcal M}_{\rho,1}(\zeta))}$ and $C = C(n,\alpha) \in (0,\infty)$ is a constant.  Thus \eqref{graphical annuli hyp} is equivalent to 
\begin{equation*}
	E_{\varsigma/2}(N,\varphi,A^{\mathcal M}_{\rho,1}(\zeta))) \leq \overline{\beta} \inf_{\varphi' \in \bigcup_{p'=0}^{p-1} \Phi_{\alpha,q,p'}(P_{Z_0},S) 
		\text{ with } \rho^{\alpha-1} \|\varphi'\|_{L^2(B_1(0))} \leq 2C \widehat{E}} E_{\varsigma/2}(N,\varphi',A^{\mathcal M}_{\rho,1}(\zeta))
\end{equation*}
where $C$ is as in \eqref{graphical rmk eqn1}.

\noindent (2)  \textit{For $\overline{\varepsilon},\overline{\beta} \in (0,1)$ sufficiently small, if }
\begin{equation}\label{graphical rmk eqn2} 
	E_{\varsigma/2}(N,\varphi,A^{\mathcal M}_{\rho,1}(\zeta)) \leq \overline{\beta} \rho^{-n/2-1} \|N\|_{L^2(A^{\mathcal M}_{\rho,1}(\zeta)} 
\end{equation}
\textit{(cf.~\eqref{graphical annuli hyp}), then}
\begin{equation}\label{graphical rmk eqn3} 
	\frac{1}{C} \,\rho^{-n/2-\alpha} \|N\|_{L^2(A^{\mathcal M}_{\rho,1}(\zeta))} \leq \|\varphi\|_{L^2(B_1(0,P_{Z_0}))} 
		\leq C \rho^{-n/2-\alpha} \|N\|_{L^2(A^{\mathcal M}_{\rho,1}(\zeta))},
\end{equation}
\textit{where $C = C(n,\alpha) \in (1,\infty)$ is a constant.}

We shall prove the first inequality, as the second inequality is similar.  Notice that by the triangle inequality, $A^{\mathcal M}_{\rho,1}(\zeta) \subset B^{\mathcal M}_{2\rho}(\exp_{Z_0}(\zeta))$, $\|D\bs{\varphi}\|_{C^0(B^{\mathcal M}_{\rho_0}(Z_0))} \leq C\mathbf{m}_0^{1/2} \rho_0 <C\overline{\varepsilon}$ and \eqref{graphical rmk eqn2} 
\begin{align*}
	\rho^{-n/2-1} \|N\|_{L^2(A^{\mathcal M}_{\rho,1}(\zeta))}
	\leq\,& E_{\varsigma/2}(N,\varphi,A^{\mathcal M}_{\rho,1}(\zeta)) + 2 \rho^{-n/2-1} \|\varphi\|_{L^2(A_{\rho,1}(\zeta))} 
	\\ \leq\,& \overline{\beta} \rho^{-n/2-1} \|N\|_{L^2(A^{\mathcal M}_{\rho,1}(\zeta))} + 2^{\alpha} \rho^{2\alpha-2} \|\varphi\|_{L^2(B_1(0))} .
\end{align*}
Thus provided $\overline{\beta} \leq 1/2$, we deduce that the first inequality in \eqref{graphical rmk eqn3} holds true.

\noindent (3)  \textit{There exists $\overline{\beta} \in (0,1)$ such that if \eqref{graphical rmk eqn2} holds true, then}
\begin{equation}\label{graphical rmk eqn4} 
	\rho \sup_{A^{\mathcal M}_{\rho,1}(\zeta)} |\Delta_{\mathcal M} (\varphi \circ \exp_{Z_0}^{-1})| 
	\leq C \mathbf{m}_0^{1/2} \rho_0 \rho^{-n/2-1} \|N\|_{L^2(A^{\mathcal M}_{\rho,1}(\zeta))} , 
\end{equation}
\textit{where $C = C(n,m,q,p,M_0,N_0,C_e,C_h,\alpha,\kappa) \in (1,\infty)$ is a constant.} 

To see \eqref{graphical rmk eqn4}, let $B \subset A_{\rho,1}(\zeta)$ be any open ball and let $\varphi_{j,l}$ be as in \eqref{varphi_localized}.  
We compute that 
\begin{equation*}
	\Delta_{\mathcal M} (\varphi_{j,l} \circ \exp_{Z_0}^{-1})(X) 
	= \sum_{i=1}^n D_{D_{\vec\tau_i} \exp_{Z_0}^{-1}(X)} D_{D_{\vec\tau_i} \exp_{Z_0}^{-1}} \varphi_{j,l}(\exp_{Z_0}^{-1}(X)) 
		- D_{\Delta_{\mathcal M} \exp_{Z_0}^{-1}(X)} \varphi_{j,l}(\exp_{Z_0}^{-1}(X)) 
\end{equation*}
for every $X \in \exp_{Z_0}(B)$ and orthonormal basis $\vec\tau_1,\vec\tau_2,\ldots,\vec\tau_n$ for ${\rm T}_X \mathcal{M}$.  Hence using \eqref{exp map est} and \eqref{graphical rmk eqn3} 
\begin{equation*}
	\rho \|\Delta_{\mathcal M} (\varphi \circ \exp_{Z_0}^{-1})\|_{C^0(A_{\rho,1}(\zeta))} 
	\leq C \mathbf{m}_0^{1/2} \rho_0 \rho^{\alpha-1} \|\varphi\|_{L^2(B_1(0,P_{Z_0}))} 
	\leq C \mathbf{m}_0^{1/2} \rho_0 \rho^{-n/2-1} \|N\|_{L^2(A^{\mathcal M}_{\rho,1}(\zeta))}
\end{equation*}
where $C = C(n,m,q,p,M_0,N_0,C_e,C_h,\alpha,\kappa) \in (0,\infty)$ are constants.

\noindent (4)  \textit{There exists $\overline{\varepsilon},\overline{\beta} \in (0,1)$ such that if $p > 0$ and \eqref{graphical annuli hyp} holds true, then 
\begin{equation}\label{graphical rmk eqn5}
	\left( \int_{B_1(0,P_{Z_0})} \mathcal{G}(\varphi_j,\varphi_{j'})^2 \right)^{1/2} 
		\geq c \rho^{1-\alpha} \inf_{\varphi' \in \bigcup_{p'=1}^{p-1} \Phi_{\alpha,q,p'}(P_{Z_0},S)} 
			E_{\varsigma/2}(N,\varphi',A^{\mathcal M}_{\rho,1}(\zeta)) 
\end{equation}
whenever $j \neq j'$ and $\varphi_j, \varphi_{j'}$ are both non-zero and 
\begin{equation}\label{graphical rmk eqn6}
	\left( \int_{B_1(0,P_{Z_0})} |\varphi_j|^2 \right)^{1/2} 
		\geq c \rho^{1-\alpha} \inf_{\varphi' \in \bigcup_{p'=1}^{p-1} \Phi_{\alpha,q,p'}(P_{Z_0},S)} 
		E_{\varsigma/2}(N, \varphi',A^{\mathcal M}_{\rho,1}(\zeta))  
\end{equation}
whenever $\varphi_j$ is non-zero, where $c = c(n,m,q,p,M_0,N_0,C_e,C_h,\alpha,\kappa) > 0$ is a constant.}  

Notice that \eqref{graphical rmk eqn5} and \eqref{graphical rmk eqn6} are weaker conclusions than \eqref{graphical annuli conclA1} and \eqref{graphical annuli conclA2}, since assuming the former it is still possible that the graph of $\varphi$ has self-intersections (or is close to self-intersecting) along an $(n-1)$-dimensional half-plane.  To see \eqref{graphical rmk eqn5}, let $\widetilde{\varphi} = \sum_{k \neq j'} m_k \varphi_k + m_{j'} \varphi_j$ (i.e.~collapse the components $\varphi_j,\varphi_{j'}$ to the same component $\varphi_j$) and let $\varphi' \in \bigcup_{p'=1}^{p-1} \Phi_{\alpha,q,p'}$ be the average-free part of $\widetilde{\varphi}$.  Notice that 
\begin{equation}\label{graphical rmk eqn7}
	\sup_{B_1(0,P_{Z_0})} \mathcal{G}(\varphi',\varphi) 
	\leq \sup_{B_1(0,P_{Z_0})} \mathcal{G}(\widetilde{\varphi},\varphi) 
	\leq \sup_{B_1(0,P_{Z_0})} \mathcal{G}(\varphi_j,\varphi_{j'}) .
\end{equation}
By the triangle inequality, \eqref{graphical annuli hyp},\eqref{graphical rmk eqn7}, and~\cite[Lemma~A.5]{KrumWic2} 
\begin{align*}
	E_{\varsigma/2}(N,\varphi',A^{\mathcal M}_{\rho,1}(\zeta)) 
	\leq\,& E_{\varsigma/2}(N,\varphi,A^{\mathcal M}_{\rho,1}(\zeta))^2 
		+ C \rho^{-1} \sup_{X \in A_{\rho,1}(\zeta)} \mathcal{G}(\varphi,\varphi') 
	\\ \leq\,& \overline{\beta} \,E_{\varsigma/2}(N,\varphi',A^{\mathcal M}_{\rho,1}(\zeta))^2 
		+ C \rho^{\alpha-1} \|\mathcal{G}(\varphi_j,\varphi_{j'})\|_{L^2(B_1(0,P_{Z_0}))}
\end{align*}
where $C = C(n,m,q,p,M_0,N_0,C_e,C_h,\alpha,\kappa) \in (0,\infty)$ are constants.  Taking $\overline{\beta} \leq 1/2$ gives us \eqref{graphical rmk eqn5}.  \eqref{graphical rmk eqn6} follows by a similar argument in which we let $\varphi'$ be the average-free part of $\widetilde{\varphi} =  m_j \llbracket 0 \rrbracket + \sum_{k \neq j} m_k \varphi_k$ (i.e.~collapse the component $\varphi_j$ to zero).
}\end{remark}

We shall give the proofs of Theorem~\ref{graphical_thm} and Lemma~\ref{graphical_lemma} at the end of this section.  The proof of Lemma~\ref{graphical_lemma} will proceed by induction on $p$ and will make use of a blow-up procedure that inductively uses the lemma itself.  We shall next discuss this blow-up procedure, which will also play an essential role in the proof of the excess decay lemma, Theorem~\ref{main excess decay thm}.

\subsection{Fine blow-ups of area-minimizing currents relative to cylindrical harmonic functions}\label{sec:blowup subsec}  We follow a blow-up procedure based on that of~\cite[Subsection~6.3]{KrumWic2}.  Let $M_0,N_0,C_e,C_h$ be suitable constants and in particular satisfy the requirements of Hypothesis~\ref{DeLSpa constants hypothesis} and let 

.  Let $p \in \{1,2,\ldots,\lfloor q/q_0 \rfloor\}$ and let $\varsigma = \varsigma(n,m,q) \in (0,1)$ be as in Theorem~\ref{relharm separation main thm}.  For $\nu = 1,2,3,\ldots$, let $\mathbf{m}_{0,\nu}$, $T^{(\nu)}$, $\mathcal{M}^{(\nu)}$, and $N^{(\nu)}$ satisfy Hypothesis~\ref{M and N hypothesis} with $\mathbf{m}_{0,\nu},T^{(\nu)},\mathcal{M}^{(\nu)},N^{(\nu)}$ in place of $\mathbf{m}_0,T,\mathcal{M},N$.  Let $Z_{\nu} \in \mathcal{B}^{(\geq 2)}_q(T^{(\nu)})$, $P_{\nu} = {\rm T}_{Z_{\nu}} \mathcal{M}^{(\nu)}$, and $S_{\nu}$ be $(n-2)$-dimensional linear subspace of $P_{\nu}$.  Let $\varphi^{(\nu)} \in \Phi_{\alpha,q,p}(P_{\nu},S_{\nu})$.  We want to define the following two types of fine blow-ups: 
\begin{enumerate}[itemsep=3mm,topsep=0mm]
	\item[(I)]  the fine blow-up of $T^{(\nu)}$ relative to a graph of $\varphi^{(\nu)}$ in $\mathbf{C}^{\mathcal{M}^{(\nu)}}_{\rho_{\nu}}(Z_{\nu})$, where $\rho_{\nu} \in (0,1]$ with $\mathbf{m}_{0,\nu} \rho_{\nu} \to 0$; 
	\item[(II)]  the fine blow-up of $T^{(\nu)}$ relative to a graph of $\varphi^{(\nu)}$ in $\mathbf{A}^{\mathcal{M}^{(\nu)}}_{\rho_{\nu},1}(\zeta_{\nu})$, where $\zeta_{\nu} \in S_{\nu}$, $\rho_{0,\nu} \in (0,1]$, and $\rho_{\nu} > 0$ such that with $\mathbf{m}_{0,\nu} \rho_{0,\nu} \to 0$ and $|\zeta_{\nu}|^2 + \rho_{\nu}^2 < (3+\gamma)^2 \rho_{0,\nu}^2/16$.  (Here we let $A^{\mathcal{M}^{(\nu)}}_{\rho_{\nu},1}(\zeta_{\nu}) = \mathbf{e}_{\nu}(\{ x+y : x  \in P_{\nu} \cap S_{\nu}^{\perp} \text{, } y \in S_{\nu} \text{, } (|x| - \rho)^2 + |y-\zeta_{\nu}|^2 < (\rho_{\nu}/64)^2 \})$ and $\mathbf{A}^{\mathcal{M}^{(\nu)}}_{\rho_{\nu},1}(\zeta_{\nu}) = \mathbf{p}_{\nu}^{-1}(A^{\mathcal{M}^{(\nu)}}_{\rho_{\nu},1}(\zeta_{\nu}))$ where $\mathbf{e}_{\nu}$ is the exponential map of $\mathcal{M}^{(\nu)}$ at $Z_{\nu}$ and $\mathbf{p}_{\nu}$ is the closest point projection map onto $\mathcal{M}^{(\nu)}$.)
\end{enumerate}

Let $\mathfrak{q}_{\nu}$ be a rotation of $\mathbb{R}^{n+m}$ such that $\mathfrak{q}_{\nu}(P_{\nu}) = \mathbb{R}^n \times \{0\}$ and $\mathfrak{q}_{\nu}(S_{\nu}) = \{0\}^2 \times \mathbb{R}^{n-2} \times \{0\}^m$.  Rescale and rotate by setting 
\begin{gather}\label{blowup eqn1}
	\overline{T}^{(\nu)} = \mathfrak{q}_{\nu\#} \eta_{Z_{\nu},\rho_{\nu}\#} T^{(\nu)}, \quad\quad
	\overline{\mathcal M}^{(\nu)} = \mathfrak{q}_{\nu}(\eta_{Z_{\nu},\rho_{\nu}}(\mathcal{M}^{(\nu)})) , \\ 
	\overline{N}^{(\nu)}(X) = \sum_{l=1}^q \llbracket \rho_{\nu}^{-1} \mathfrak{q}_{\nu} N_l^{(\nu)}(Z_{\nu} + \rho_{\nu} \mathfrak{q}_{\nu}^{-1} X) \rrbracket 
		 \text{ for each } X \in \overline{\mathcal M}^{(\nu)} , \nonumber \\
	\overline{\varphi}^{(\nu)}(X) = \sum_{l=1}^q \llbracket \rho_{\nu}^{\alpha-1} \mathfrak{q}_{\nu} \varphi_l^{(\nu)}(\mathfrak{q}_{\nu}^{-1} X) \rrbracket 
		\text{ for each } X \in \mathbb{R}^n , \nonumber 
\end{gather}
where $N^{(\nu)}(X) = \sum_{l=1}^q \llbracket N_l^{(\nu)}(X) \rrbracket$ for each $X \in \mathcal{M}^{(\nu)}$.  Express 
\begin{equation*}
	\overline{\varphi}^{(\nu)} = \sum_{j=1}^J m_j \overline{\varphi}^{(\nu)}_j 
\end{equation*}
on $\mathbb{R}^n$, where $\overline{\varphi}^{(\nu)}_j$ are the distinct $q_j$-valued components of $\overline{\varphi}^{(\nu)}$ with multiplicity $m_j$ (as in Definition~\ref{Phi-p_defn}).  Let $\overline{\mathbf p}_{\nu}$ be the closest point projection map onto $\overline{\mathcal M}^{(\nu)}$.  Now in cases~(I) and (II), we apply Theorem~\ref{graphical_thm} and Lemma~\ref{graphical_lemma} as follows.

\noindent\textbf{Case~(I).}  Let $\widehat{E}_{\nu} = \widehat{E} = \rho_{\nu}^{-n/2-1} \|N^{(\nu)}\|_{L^2(B^{\mathcal{M}^{(\nu)}}_{\rho_{\nu}}(Z_{\nu}))}$ and $E_{\nu} = E_{\varsigma/2}(N^{(\nu)},\varphi^{(\nu)},B^{\mathcal{M}^{(\nu)}}_{\rho_{\nu}}(Z_{\nu}))$.  For $\beta_{\nu} \rightarrow 0^+$, we shall assume that 
\begin{equation}\label{blowup eqn3}
	E_{\nu} \leq \beta_{\nu} \inf_{\varphi' \in \bigcup_{p'=0}^{p-1} \Phi_{\alpha,q,p'}(P_{\nu},S_{\nu})} 
		E_{\varsigma/2}(N^{(\nu)},\varphi',B^{\mathcal{M}^{(\nu)}}_{\rho_{\nu}}(Z_{\nu})) . 
\end{equation}
We let $\Omega = B_{3/4}(0)$ and $\Omega_{\nu} = B_{3/4}(0) \setminus (B^2_{\tau_{\nu}}(0) \times \mathbb{R}^{n-2})$ for some suitable sequence $\tau_{\nu} \rightarrow 0^+$.  Let $\overline{\mathbf e}_{\nu} = \exp_{\overline{\mathcal M}^{(\nu)}} : B_1(0) \rightarrow \overline{\mathcal M}^{(\nu)}$ be the exponential map of $\overline{\mathcal M}^{(\nu)}$ at the origin.  (That way, $\overline{\mathbf e}_{\nu}(\Omega) = B^{\overline{\mathcal M}^{(\nu)}}_{3/4}(0)$ and $\overline{\mathbf e}_{\nu}(\Omega_{\nu}) = B^{\overline{\mathcal M}^{(\nu)}}_{3/4}(0) \cap \{r_{\nu} > \tau_{\nu}\}$, where $r_{\nu}(X) = \op{dist}(\overline{\mathbf e}_{\nu}^{-1}(X),\{0\}^2 \times \mathbb{R}^{n-2} \times \{0\}^m)$.)  By applying Theorem~\ref{graphical_thm} with $T^{(\nu)}, \mathcal{M}^{(\nu)}, N^{(\nu)}, \varphi^{(\nu)}$ in place of $T,\mathcal{M},N,\varphi$ and then rotating and rescaling, for some suitable sequence of positive numbers $(\tau_{\nu})$ there exists: 
\begin{enumerate}[itemsep=3mm,topsep=0mm]
	\item[(i)]  $n$-dimensional locally area-minimizing rectifiable currents $\overline{T}^{(\nu)}_j$ in $\overline{\mathbf p}_{\nu}^{-1}(\overline{\mathbf e}_{\nu}(\Omega_{\nu}))$ such that 
	\begin{gather}\label{blowup eqn5}
		\overline{T}^{(\nu)} \llcorner \overline{\mathbf p}_{\nu}^{-1}(\overline{\mathbf e}_{\nu}(\Omega_{\nu})) = \sum_{j=1}^J \overline{T}^{(\nu)}_j , \\
		(\partial \overline{T}^{(\nu)}_j) \llcorner \overline{\mathbf p}_{\nu}^{-1}(\overline{\mathbf e}_{\nu}(\Omega_{\nu})) = 0 , \quad 
		\overline{\mathbf p}_{\nu \#} \overline{T}^{(\nu)}_j = q_j m_j \llbracket \overline{\mathbf e}_{\nu}(\Omega_{\nu}) \rrbracket , \nonumber \\
		\sup_{X \in \op{spt} \overline{T}^{(\nu)}_j \cap \overline{\mathbf p}_{\nu}^{-1}(\overline{\mathbf e}_{\nu}(\Omega'))} 
			\op{dist}(X, \op{spt} \mathbf{G}_{\overline{\varphi}^{(\nu)}_j \circ \overline{\mathbf e}_{\nu}^{-1}}) \leq C E_{\nu} \nonumber 
	\end{gather}
    for each (fixed) set $\Omega' \subset\subset \Omega \setminus \{0\} \times \mathbb{R}^{n-2}$ such that $\Omega' \subseteq \Omega_{\nu}$, where $C = C(n,m,q,M_0,N_0,C_e,$ $C_h,\alpha,\Omega') \in (0,\infty)$ is a constant independent of $\nu$; 
	
	\item[(ii)]  Lipschitz $m_j$-valued functions $v^{(\nu)}_j : \mathbf{G}_{\overline{\varphi}^{(\nu)}_j \circ \overline{\mathbf e}_{\nu}^{-1}} \cap \overline{\mathbf p}_{\nu}^{-1}(\overline{\mathbf e}_{\nu}(\Omega_{\nu})) \rightarrow \mathcal{A}_{m_j}(\mathbb{R}^{n+m})$ and closed sets $K^{(\nu)}_j \subseteq \overline{\mathbf e}_{\nu}(\Omega_{\nu})$ such that 
	\begin{gather}\label{blowup eqn6}
		\overline{T}^{(\nu)}_j \llcorner \overline{\mathbf p}_{\nu}^{-1}(K^{(\nu)}_j) 
			= \mathbf{G}_{\overline{N}^{(\nu)}_j} \llcorner \overline{\mathbf p}_{\nu}^{-1}(K^{(\nu)}_j) , \\
		\mathcal{H}^n(\overline{\mathbf e}_{\nu}(\Omega') \setminus K^{(\nu)}_j) 
			+ \|\overline{T}^{(\nu)}_j\|(\overline{\mathbf p}^{-1}(\overline{\mathbf e}_{\nu}(\Omega') \setminus K^{(\nu)}_j)) 
			\leq C \widehat{E}_{\nu}^{2+\varsigma}, \nonumber\\
		\sup_{\overline{\mathbf e}_{\nu}(\Omega')} |\overline{N}^{(\nu)}_j| \leq C \widehat{E}_{\nu}, \quad\quad 
			\sup_{\overline{\mathbf e}_{\nu}(\Omega')} |D\overline{N}^{(\nu)}_j| 
			\leq C \widehat{E}_{\nu}^{\varsigma} \nonumber 
	\end{gather}
	and 
	\begin{equation}\label{blowup eqn7}
		\sup_{\overline{\mathbf e}_{\nu}(\Omega')} \sum_{j=1}^J \sum_{l=1}^{q_j} |v^{(\nu)}_{j,l}|^2 
		+ \int_{\overline{\mathbf e}_{\nu}(\Omega')} \sum_{j=1}^J \sum_{l=1}^{q_j} |Dv^{(\nu)}_{j,l}|^2 
		\leq C E_{\nu}^2
	\end{equation}
	for each (fixed) set $\Omega' \subset\subset \Omega \setminus \{0\} \times \mathbb{R}^{n-2}$ such that $\Omega' \subseteq \Omega_{\nu}$, where $v^{(\nu)}_{j,l}$ and $\overline{N}^{(\nu)}_j$ are as in \eqref{v_localized} and \eqref{N of v defn} with $\overline{\varphi}^{(\nu)},\overline{N}^{(\nu)}_j,v^{(\nu)}_{j,l}$ in place of $\varphi,N_j,v_j$ and $C = C(n,m,q,M_0,N_0,C_e,C_h,\alpha,\Omega') \in (0,\infty)$ is a constant independent of $\nu$.
\end{enumerate}

\noindent\textbf{Case~(II).}  Let $\widehat{E}_{\nu} = \rho_{\nu}^{-n/2-1} \|N^{(\nu)}\|_{L^2(A^{\mathcal{M}^{(\nu)}}_{\rho_{\nu},1}(\zeta_{\nu}))}$ and let $E_{\nu} = E_{\varsigma/2}(N^{(\nu)},\varphi^{(\nu)},A^{\mathcal{M}^{(\nu)}}_{\rho_{\nu},1}(\zeta_{\nu}))$.  We shall assume that 
\begin{equation*} 
	E_{\nu} \leq \overline{\beta} \inf_{\varphi' \in \bigcup_{p'=0}^{p-1} \Phi_{\alpha,q,p'}(P_{\nu},S_{\nu})} 
		E_{\varsigma/2}(N^{(\nu)},\varphi',A^{\mathcal{M}^{(\nu)}}_{\rho_{\nu},1}(\zeta_{\nu})) ,
\end{equation*}
where $\overline{\beta}$ is as in Lemma~\ref{graphical_lemma}.  (Here $E_{\varsigma/2}(N^{(\nu)},\varphi^{(\nu)},A^{\mathcal{M}^{(\nu)}}_{\rho_{\nu},1}(\zeta_{\nu}))$ and $E_{\varsigma/2}(N^{(\nu)},\varphi',A^{\mathcal{M}^{(\nu)}}_{\rho_{\nu},1}(\zeta_{\nu}))$ are as in \eqref{graphical annuli E defn} with $Z_{\nu},\rho_{0,\nu}$ in place of $Z_0,\rho_0$.)  Upon rescaling, we also set $\overline{\zeta}_{\nu} = \mathfrak{q}_{\nu}(\eta_{0,\rho_{\nu}}(\zeta_{\nu}))$.  We let $\Omega = A_{1,1}(0)$ and $\Omega_{\nu} = A_{1,\kappa_{\nu}}(0)$ for some suitable sequence $\kappa_{\nu} \rightarrow 1^-$, where for each $\kappa \in (0,1]$ we let $A_{1,\kappa}(0) = \{ (x,y) \in \mathbb{R}^2 \times \mathbb{R}^{n-2} : (|x| - \rho)^2 + |y|^2 < (1-\gamma)^2 \rho^2/64 \}$.  We let $\overline{\mathbf e}_{\nu} : A_{1,1}(0) \rightarrow \overline{\mathcal M}^{(\nu)}$ such that $\overline{\mathbf e}_{\nu}(X) = \exp_{\overline{\mathcal M}^{(\nu)}}(\overline{\zeta}_{\nu} + X)$ for each $X \in A_{1,1}(0)$, where again $\exp_{\overline{\mathcal M}^{(\nu)}}$ is the exponential map of $\overline{\mathcal M}^{(\nu)}$ at the origin.  (That way, $\overline{\mathbf e}_{\nu}(\Omega) = A^{\overline{\mathcal M}^{(\nu)}}_{1,1}(\overline{\zeta}_{\nu})$ and $\overline{\mathbf e}_{\nu}(\Omega_{\nu}) = A^{\overline{\mathcal M}^{(\nu)}}_{1,\kappa_{\nu}}(\overline{\zeta}_{\nu})$.)  By applying Lemma~\ref{graphical_lemma} with $T^{(\nu)}, \mathcal{M}^{(\nu)}, N^{(\nu)}, \varphi^{(\nu)}$ in place of $T,\mathcal{M},N,\varphi$ and then rotating and rescaling, for some suitable sequence of numbers $(\kappa_{\nu})$ there exists $n$-dimensional locally area-minimizing rectifiable currents $\overline{T}^{(\nu)}_j$ in $\overline{\mathbf p}_{\nu}^{-1}(\overline{\mathbf e}_{\nu}(\Omega_{\nu}))$ such that \eqref{blowup eqn5} holds true and there exists Lipschitz $m_j$-valued functions $v^{(\nu)}_j : \mathbf{G}_{\overline{\varphi}^{(\nu)}_j \circ \overline{\mathbf e}_{\nu}^{-1}} \cap \overline{\mathbf p}_{\nu}^{-1}(\overline{\mathbf e}_{\nu}(\Omega_{\nu})) \rightarrow \mathcal{A}_{m_j}(\mathbb{R}^{n+m})$ and closed sets $K^{(\nu)}_j \subseteq \overline{\mathbf e}_{\nu}(\Omega_{\nu})$ such that \eqref{blowup eqn6} and \eqref{blowup eqn7} hold true.

Now after passing to a subsequence of $(\nu)$, we may assume that for each $j \in \{1,\ldots,J\}$ one of the following two possibilities holds:
\begin{enumerate}[itemsep=3mm,topsep=0mm]
	\item[(a)]  $\overline{\varphi}^{(\nu)}_j$ is not identically zero for all $\nu = 1,2,3,\ldots$; 
	\item[(b)]  $\overline{\varphi}^{(\nu)}_j \equiv 0$ for all $\nu = 1,2,3,\ldots$.
\end{enumerate}
We shall construct cylindrical functions $\varphi^{(0)}_j$ and $\varphi^{(\infty)}_j$ and functions (fine blow-ups) 
\begin{equation*}
	w_j : \op{spt} \mathbf{G}_{\varphi^{(\infty)}_j} \cap (\Omega \setminus \{0\} \times \mathbb{R}^{n-2}) \times \mathbb{R}^m
		\rightarrow \mathcal{A}_{m_j}(\mathbb{R}^{n+m})
\end{equation*}
by considering these two cases as follows:

\noindent\textbf{Case~(a):}  Let $(x,y) = (re^{i\theta},y)$ denote cylindrical coordinates on $\mathbb{R}^n$, where $r > 0$, $\theta \in \mathbb{R}$, and $y \in \mathbb{R}^{n-2}$.  Let $\overline{\varphi}^{(\nu)}_j(re^{i\theta},y) = \op{Re}(c^{(\nu)}_j r^{\alpha} e^{i\alpha\theta})$ for $c^{(\nu)}_j \in \mathbb{C}^m \setminus \{0\}$.  Recalling Remark~\ref{graphical rmk}(2), $|c^{(\nu)}_j| \leq C \widehat{E}_{\nu}$, where $C = C(n,\alpha) \in (0,\infty)$ is a constant.  Hence after passing to a subsequence of $(\nu)$ let 
\begin{equation*}
	c^{(0)}_j = \lim_{\nu\rightarrow\infty} \frac{c^{(\nu)}_j}{\widehat{E}_{\nu}} , \quad\quad 
	 c^{(\infty)}_j = \lim_{\nu\rightarrow\infty} \frac{c^{(\nu)}_j}{|c^{(\nu)}_j|} . 
\end{equation*}
(Note that $c^{(0)}_j$ may or may not equal zero, whereas $c^{(\infty)}_j$ is a unit vector.)  Define $\varphi^{(0)}_j : \mathbb{R}^n \rightarrow \mathcal{A}_{q_0}(\mathbb{R}^m)$ and $\varphi^{(\infty)}_j : \mathbb{R}^n \rightarrow \mathcal{A}_{q_0}(\mathbb{R}^m)$ by 
\begin{equation*}
	\varphi^{(0)}_j(x_1,x_2,y) = \op{Re}(c^{(0)}_j (x_1+ix_2)^{\alpha}) , \quad 
	\varphi^{(\infty)}_j(x_1,x_2,y) = \op{Re}(c^{(\infty)}_j (x_1+ix_2)^{\alpha}) 
\end{equation*}
for each $(x_1,x_2,y) \in \mathbb{R}^n$.  Observe that $\op{Re}(c^{(\nu)}_{j,k} r^{\alpha} e^{i\alpha\theta})$, $\op{Re}(c^{(0)}_j r^{\alpha} e^{i\alpha\theta})$, and $\op{Re}(c^{(\infty)}_j r^{\alpha} e^{i\alpha\theta})$ are well-defined single-valued functions of $r > 0$, $\theta \in \mathbb{R}$, and $y \in \mathbb{R}^{n-2}$ which are $2\pi q_0$-periodic as functions of $\theta$.  In particular, for each $r > 0$, $\theta \in \mathbb{R}$, and $y \in \mathbb{R}^{n-2}$, $(re^{i\theta},y,\op{Re}(c^{(\nu)}_j r^{\alpha} e^{i\alpha\theta}))$ and $(re^{i\theta},y,\op{Re}(c^{(\infty)}_j r^{\alpha} e^{i\alpha\theta}))$ are well-defined points on $\op{graph} \varphi^{(\nu)}$ and $\op{graph} \varphi^{(\infty)}$ respectively.  Define $w^{(\nu)}_j : \op{spt} \mathbf{G}_{\varphi^{(\infty)}_j} \cap \overline{\mathbf p}_{\nu}^{-1}(\overline{\mathbf e}_{\nu}(\Omega_{\nu})) \rightarrow \mathcal{A}_{m_j}(\mathbb{R}^{n+m})$ by 
\begin{equation}\label{blowup eqn10}
	w^{(\nu)}_j(re^{i\theta},y,\op{Re}(c^{(\infty)}_j r^{\alpha} e^{i\alpha\theta})) 
	= \frac{v^{(\nu)}_j(\overline{\mathbf e}_{\nu}(re^{i\theta},y),\op{Re}(c^{(\nu)}_j r^{\alpha} e^{i\alpha\theta}))}{E_{\nu}}
\end{equation}
for every $r > 0$, $\theta \in \mathbb{R}$, and $y \in \mathbb{R}^{n-2}$ such that $(re^{i\theta},y) \in \Omega_{\nu}$.  By \eqref{blowup eqn7} and the Rellich compactness lemma for multi-valued Sobolev functions~\cite[Proposition~2.11]{DeLSpaDirMin}, after passing to a subsequence of $(\nu)$ there exists a function $w_j : \op{spt} \mathbf{G}_{\varphi^{(\infty)}_j} \cap (\Omega \setminus \{0\} \times \mathbb{R}^{n-2}) \times \mathbb{R}^m \rightarrow \mathcal{A}_{m_j}(\mathbb{R}^{n+m})$ such that $w^{(\nu)}_j \rightarrow w_j$ in $L^2(\op{spt} \mathbf{G}_{\varphi^{(\infty)}_j} \cap \Omega' \times \mathbb{R}^m)$ for each $\Omega' \subset\subset \Omega \setminus \{0\}\times\mathbb{R}^{n-2}$.  Note that by Remark~\ref{graphical rmk}(3), $\|\Delta_{\overline{\mathcal M}^{(\nu)}} \overline{\varphi}^{(\nu)}_j\|_{C^0(\overline{\mathbf e}_{\nu}(\Omega))} \leq C \mathbf{m}_0^{1/2} E_{\nu}$.  Thus arguing as in Lemma~\ref{relharm harmonic approx lemma}, $w_j$ is component-wise Dirichlet energy minimizing in $\Omega\setminus\{0\}\times\mathbb{R}^{n-2}$ (as in Definition~\ref{v comp minimizing defn}) and $|Dw_j^{(\nu)}| \rightarrow |Dw_j|$ in $L^2(\op{spt} \mathbf{G}_{\varphi^{(\infty)}_j} \cap \Omega' \times \mathbb{R}^m)$ for each $\Omega' \subset\subset \Omega \setminus \{0\}\times\mathbb{R}^{n-2}$.

\noindent\textbf{Case~(b):}  Define $\varphi^{(0)} : \mathbb{R}^n \rightarrow \mathbb{R}^m$ and $\varphi^{(\infty)} : \mathbb{R}^n \rightarrow \mathbb{R}^m$ by $\varphi^{(0)}_j(x,y) = \varphi^{(\infty)}_j(x,y) = 0$ for all $(x,y) \in \mathbb{R}^n$.  Define $w^{(\nu)}_j : \op{spt} \mathbf{G}_{\varphi^{(\infty)}_j} \cap \overline{\mathbf p}_{\nu}^{-1}(\overline{\mathbf e}_{\nu}(\Omega_{\nu})) \rightarrow \mathcal{A}_{m_j}(\mathbb{R}^{n+m})$ by 
\begin{equation}\label{blowup eqn11}
	w^{(\nu)}_j(x,y,0) = \frac{v^{(\nu)}_j(\overline{\mathbf e}_{\nu}(x,y),0)}{E_{\nu}} 
\end{equation}
for all $(x,y) \in \Omega_{\nu}$.  By \eqref{blowup eqn7} and the Rellich compactness lemma for multi-valued Sobolev functions~\cite[Proposition~2.11]{DeLSpaDirMin}, after passing to a subsequence of $(\nu)$ there exists a function $w_j : \op{spt} \mathbf{G}_{\varphi^{(\infty)}_j} \cap (\Omega \setminus \{0\} \times \mathbb{R}^{n-2}) \times \mathbb{R}^m \rightarrow \mathcal{A}_{m_j}(\mathbb{R}^{n+m})$ such that $w^{(\nu)}_j \rightarrow w_j$ in $L^2(\op{spt} \mathbf{G}_{\varphi^{(\infty)}_j} \cap \Omega' \times \mathbb{R}^m)$ for each $\Omega' \subset\subset \Omega \setminus \{0\}\times\mathbb{R}^{n-2}$.  Arguing as in Lemma~\ref{relharm harmonic approx lemma}, $w_j$ is component-wise Dirichlet energy minimizing in $\Omega\setminus\{0\}\times\mathbb{R}^{n-2}$ and $|Dw_j^{(\nu)}| \rightarrow |Dw_j|$ in $L^2(\op{spt} \mathbf{G}_{\varphi^{(\infty)}_j} \cap \Omega' \times \mathbb{R}^m)$ for each $\Omega' \subset\subset \Omega \setminus \{0\}\times\mathbb{R}^{n-2}$.

Let $\varphi^{(0)} = \sum_{j=1}^J m_j \varphi^{(0)}_j$.  In Case~(I), by $\varphi^{(\nu)}_j/\widehat{E}_{\nu} \rightarrow \varphi^{(0)}_j$ and \eqref{blowup eqn3}, it follows that $\varphi^{(0)}$ is the coarse blow-up of $T^{(\nu)}$ relative to $\mathcal{M}^{(\nu)}$ in $B^{\mathcal{M}^{(\nu)}}_{\rho_{\nu}}(Z_{\nu})$ by excess $\widehat{E}_{\nu}$ (as in Subsection~\ref{sec:centered blowup sec}).  In case~(I), we say that $w = (w_j)$ is a \emph{(fine) blow-up of $T^{(\nu)}$ relative to $\varphi^{(\nu)}$ in $\mathbf{C}^{\mathcal{M}^{(\nu)}}_{\rho_{\nu}}(Z_{\nu})$ by excess $E_{\nu}$}.  Similarly, in case~(II), we say that $w = (w_j)$ is a \emph{(fine) blow-up of $T^{(\nu)}$ relative to $\varphi^{(\nu)}$ in $\mathbf{A}^{\mathcal{M}^{(\nu)}}_{\rho_{\nu},1}(\zeta_{\nu})$ by excess $E_{\nu}$}.

\begin{remark}\label{cylindrical blowups rmk}{\rm 
Let $\beta = \beta(n,m,q,\alpha) \in (0,1)$ be a small constant.  For $\nu = 1,2,3,\ldots$, let $\varepsilon_{\nu} \rightarrow 0^+$, $P_{\nu}$ be $n$-dimensional linear subspaces of $\mathbb{R}^{n+m}$, $S_{\nu}$ be $(n-2)$-dimensional linear subspaces of $P_{\nu}$, $p \in \{1,2,\ldots,\lfloor q/q_0 \rfloor\}$, $\varphi^{(\nu)} \in \Phi_{\alpha,q,p}(P_{\nu},S_{\nu})$, and $u^{(\nu)} \in \Phi_{\alpha,q}(P_{\nu},S_{\nu})$ such that 
\begin{gather}
	\|\varphi^{(\nu)}\|_{L^2(B_1(0))} < \varepsilon_{\nu}, \quad \|u^{(\nu)}\|_{L^2(B_1(0))} < \varepsilon_{\nu}, \nonumber \\  
	\label{blowup eqn12} \int_{B_1(0)} \mathcal{G}(u^{(\nu)}, \varphi^{(\nu)})^2 
		\leq \beta \inf_{\varphi' \in \bigcup_{p'=0}^{p-1} \Phi_{\alpha,q,p'}(P_{\nu},S_{\nu})} \int_{B_1(0)} \mathcal{G}(u^{(\nu)}, \varphi')^2 .
\end{gather}
By modifying the approach of~\cite[Remark~5.11]{KrumWicc}, we can define the fine blow-up of $u^{(\nu)}$ relative to $\varphi^{(\nu)}$ as follows.  Let $E_{\nu} = \left( \int_{B_1(0)} \mathcal{G}(u^{(\nu)},\varphi^{(\nu)})^2 \right)^{1/2}$.  Rotate so that $P_{\nu} = \mathbb{R}^n \times \{0\}$ and $S_{\nu} = \{0\} \times \mathbb{R}^{n-2} \times \{0\}$.  Let $\varphi^{(\nu)} = \sum_{j=1}^J m_j \varphi^{(\nu)}_j$ on $\mathbb{R}^n$, where $\varphi^{(\nu)}_j$ are the distinct $q_j$-valued components of $\varphi^{(\nu)}$ with multiplicity $m_j$ (as in Definition~\ref{Phi-p_defn}).  Assuming $\beta$ is sufficiently small, by~\ref{blowup eqn12} and~\cite[Lemma~A.5 and Remark~6.6(b)]{KrumWic2}, after passing to a subsequence 
\begin{equation*}
	u^{(\nu)} = \sum_{j=1}^J \sum_{k=1}^{s_j} m_{j,k} u^{(\nu)}_{j,k} , 
\end{equation*}
where $u^{(\nu)}_{j,k}$ are distinct $q_{j,k}$-valued components of $u^{(\nu)}$ such that either 
\begin{enumerate}[itemsep=3mm,topsep=0mm]
	\item[(i)]  $\varphi^{(\nu)}_j$ and $u^{(\nu)}_{j,k}$ are both non-zero and $\sup_{B_1(0)} \mathcal{G}(u^{(\nu)}_{j,k},\varphi^{(\nu)}_{j,k}) \leq C(n,m,q,\alpha) \,E_{\nu}$, 
	\item[(ii)]  $\varphi^{(\nu)}_j$ and $u^{(\nu)}_{j,k}$ are both identically zero, or 
	\item[(iii)]  $\varphi^{(\nu)}_j$ is identically zero, $u^{(\nu)}_{j,k}$ is non-zero, and $|u^{(\nu)}_{j,k}| \leq C(n,m,q,\alpha) \,E_{\nu}$
\end{enumerate}
and $s_j$ and $m_{j,k}$ are positive integers independent of $\nu$ such that $\sum_{j=1}^J \sum_{k=1}^{s_j} m_{j,k} q_{j,k} = q$.  Set $\widehat{m}_j = \sum_{k=1}^{s_j} m_{j,k}$ if $\varphi^{(\nu)}_j$ is non-zero and $\widehat{m}_j = \sum_{k=1}^{s_j} m_{j,k} q_{j,k}$ if $\varphi^{(\nu)}_j$ is identically zero.  Set $u^{(\nu)}_j = \sum_{k=1}^{s_j} m_{j,k} u^{(\nu)}_{j,k}$.  Since $u^{(\nu)}_{j,k}$ is uniformly close to $\varphi^{(\nu)}_j$, there exists a unique continuous function $v^{(\nu)}_j : \op{spt} \mathbf{G}_{\varphi^{(\nu)}_j} \rightarrow \mathcal{A}_{\widehat{m}_j}(\mathbb{R}^m)$ such that $\sup_{\op{spt} \mathbf{G}_{\varphi^{(\nu)}_j} \cap \mathbf{C}_1(0)} |v^{(\nu)}_j| \leq CE_{\nu}$ and $u^{(\nu)}_j(X) = \sum_{l=1}^{q_j} \sum_{h=1}^{\widehat{m}_j} \llbracket \varphi^{(\nu)}_{j,l}(X) + v^{(\nu)}_{j,l,h}(X) \rrbracket$ for each $X \in \mathbb{R}^n$, where $\varphi^{(\nu)}_j(X) = \sum_{l=1}^{q_j} \llbracket \varphi^{(\nu)}_{j,l}(X) \rrbracket$ and $v^{(\nu)}_j(X,\varphi^{(\nu)}_{j,l}(X)) = \sum_{h=1}^{\widehat{m}_j} \llbracket v^{(\nu)}_{j,l,h}(X) \rrbracket$.  Following the procedure in cases~(a) and (b) above, we produce a fine blow-up $w = (w^{(\nu)}_j)$ of $u^{(\nu)}$ relative to $\varphi^{(\nu)}$ by the excess $E_{\nu}$ in $\mathbb{R}^n$.

Note that if we are additionally given a sequence of radii $\rho_{\nu} \in (0,1]$, we can rotate so that $P_{\nu} = \mathbb{R}^n \times \{0\}$ and $S_{\nu} = \{0\} \times \mathbb{R}^{n-2} \times \{0\}$ and rescale by setting $\overline{\varphi}^{(\nu)}(X) = \rho_{\nu}^{-1} \varphi^{(\nu)}(\rho_{\nu} X) = \rho_{\nu}^{\alpha-1} \varphi^{(\nu)}(X)$ and $\overline{u}^{(\nu)}(X) = \rho_{\nu}^{-1} u^{(\nu)}(\rho_{\nu} X) = \rho_{\nu}^{\alpha-1} u^{(\nu)}(X)$ for each $X \in \mathbb{R}^n$.  Set $\overline{E}_{\nu} = \left( \int_{B_1(0)} \mathcal{G}(\overline{u}^{(\nu)},\overline{\varphi}^{(\nu)})^2 \right)^{1/2} = \left( \rho_{\nu}^{-n-2\alpha} \int_{B_{\rho_{\nu}}(0)} \mathcal{G}(u^{(\nu)},\varphi^{(\nu)})^2 \right)^{1/2}$.  Using the above procedure we can construct a fine blow-up $w = (w^{(\nu)}_j)$ of $\overline{u}^{(\nu)}$ relative to $\overline{\varphi}^{(\nu)}$ by the excess $\overline{E}_{\nu}$ in $\mathbb{R}^n$.
}\end{remark}

\subsection{Proofs of the graphical representation results}

\begin{proof}[Proof of Lemma~\ref{graphical_lemma}]
The proof proceeds by induction on $p$.  The base case $p = 0$, in which $\varphi = q \llbracket 0 \rrbracket$, is a special case of~\cite[Lemma~4.6]{KrumWicb}.  Suppose that for some $p_0 \in \{1,2,\ldots,\lfloor q/q_0 \rfloor\}$, Lemma~\ref{graphical_lemma} holds true for all $p \in \{0,1,2,\ldots,p_0-1\}$.  We want to show that Lemma~\ref{graphical_lemma} holds true when $p = p_0$.

\noindent\textit{Proof of Conclusion~(A).}  The proof is nearly identical to that of~\cite[Lemma~5.9(A)]{KrumWicc} and we sketch the argument for completion.  To see \eqref{graphical annuli conclA1}, suppose that for each $\nu = 1,2,3,\ldots$ there exists $\varepsilon_{\nu} \rightarrow 0^+$, $\beta_{\nu} \rightarrow 0^+$, $\mathbf{m}_{0,\nu}$, $T^{(\nu)}$, $\mathcal{M}^{(\nu)}$, $N^{(\nu)}$, $\mathcal{K}^{(\nu)}$, $\varphi^{(\nu)}$, $Z_{0,\nu}$, $\rho_{0,\nu}$, $P_{\nu}$, $S_{\nu}$, $\zeta_{\nu}$, and $\rho_{\nu}$ which satisfy the hypotheses of Lemma~\ref{graphical_lemma} in place of $\overline{\varepsilon}$, $\overline{\beta}$, $\mathbf{m}_0$, $T$, $\mathcal{M}$, $N$, $\mathcal{K}$, $Z_0$, $\varphi$, $\rho_0$, $P_{Z_0}$, $S$, $\zeta$, and $\rho$ but 
\begin{equation}\label{graphical annuli eqnA1}
	\inf_{x \in \mathbb{S}^{n+m-1} \cap P_{\nu} \cap S_{\nu}^{\perp}} \op{sep} \varphi^{(\nu)}(x) 
	< \frac{1}{\nu} \inf_{\varphi' \in \bigcup_{p'=0}^{p-1} \Phi_{\alpha,q,p'}(P_{\nu},S_{\nu})} 
		E_{\varsigma/2}(N^{(\nu)},\varphi',A^{\mathcal M}_{\rho_{\nu},1}(\zeta_{\nu})) . 
\end{equation}
(Here $E_{\varsigma/2}(N^{(\nu)},\varphi^{(\nu)},A^{\mathcal M}_{\rho_{\nu},1}(\zeta_{\nu}))$ and $E_{\varsigma/2}(N^{(\nu)},\varphi',A^{\mathcal M}_{\rho_{\nu},1}(\zeta_{\nu}))$ is as in \eqref{graphical annuli E defn} with $Z_{0,\nu},\rho_{0,\nu}$ in place of $Z_0,\rho_0$.)  Let $\widetilde{\beta} = \widetilde{\beta}(n,m,q,p_0,M_0,N_0,C_e,C_h,\alpha) \in (0,1)$ be a suitably small constant.  Choose $\widetilde{p} \in \{0,1,2,\ldots,p_0-1\}$ and $\widetilde{\varphi}^{(\nu)} \in \Phi_{\alpha,q,\widetilde{p}}(P_{\nu},S_{\nu})$ such that 
\begin{equation*}
	E_{\varsigma/2}(N^{(\nu)},\widetilde{\varphi}^{(\nu)},A^{\mathcal M}_{\rho_{\nu},1}(\zeta_{\nu}))
	\leq 2^{p_0} \widetilde{\beta}^{1-p_0} \inf_{\varphi' \in \bigcup_{p'=0}^{p_0-1} \Phi_{\alpha,q,p'}(P_{\nu},S_{\nu})} 
		E_{\varsigma/2}(N^{(\nu)},\varphi',A^{\mathcal M}_{\rho_{\nu},1}(\zeta_{\nu})) 
\end{equation*}
and either $\widetilde{p} = 0$ or $\widetilde{p} > 0$ and 
\begin{equation}\label{graphical annuli eqnA2}
	E_{\varsigma/2}(N^{(\nu)},\widetilde{\varphi}^{(\nu)},A^{\mathcal M}_{\rho_{\nu},1}(\zeta_{\nu}))
	\leq \widetilde{\beta} \inf_{\varphi' \in \bigcup_{p'=0}^{\widetilde{p}-1} \Phi_{\alpha,q,p'}(P_{\nu},S_{\nu})} 
		E_{\varsigma/2}(N^{(\nu)},\varphi',A^{\mathcal M}_{\rho_{\nu},1}(\zeta_{\nu})) . 
\end{equation}
Let $\widetilde{\varphi}^{(\nu)} = \sum_{j=1}^{\widetilde{p}} \widetilde{m}_j \widetilde{\varphi}^{(\nu)}_j$ on $\mathbb{R}^n$, where $\widetilde{\varphi}^{(\nu)}_j$ are the distinct $\widetilde{q}_j$-valued components of $\widetilde{\varphi}^{(\nu)}$ with multiplicity $\widetilde{m}_j$ (as in Definition~\ref{Phi-p_defn}).  After passing to a subsequence, we may take $\widetilde{p}$, $\widetilde{m}_j$, and $\widetilde{q}_j$ to be independent of $\nu$.  Set $\widetilde{E}_{\nu} = E_{\varsigma/2}(N^{(\nu)},\widetilde{\varphi}^{(\nu)},A^{\mathcal M}_{\rho_{\nu},1}(\zeta_{\nu}))$.  By the inductive hypothesis and \eqref{graphical annuli eqnA2}, we can apply the blow-up procedure from Subsection~\ref{sec:blowup subsec} to construct a blow-up $w = (w_j)$ of $T^{(\nu)}$ relative to $\widetilde{\varphi}^{(\nu)}$ by excess $\widetilde{E}_{\nu}.$  Moreover, we can apply the blow-up procedure from Remark~\ref{cylindrical blowups rmk} to construct a blow-up $\psi = (\psi_j)$ of $\rho_{\nu}^{\alpha-1} \varphi^{(\nu)}$ relative to $\rho_{\nu}^{\alpha-1} \widetilde{\varphi}^{(\nu)}$ by excess $\widetilde{E}_{\nu}$.  By \eqref{graphical annuli hyp}, $w_j = \psi_j$ on $\op{spt} \mathbf{G}_{\varphi^{(\infty)}_j} \cap \mathbf{A}_{1,1/2}(0)$.  (By contrast, in~\cite[Lemma~5.9]{KrumWicc} we only get $\op{spt} w_j(X) = \op{spt} \psi_j(X)$ at each $X \in \op{spt} \mathbf{G}_{\varphi^{(\infty)}_j} \cap \mathbf{A}_{1,1/2}(0)$.)  Hence $\psi$ is homogeneous degree $\alpha$, translation invariant along $\{0\} \times \mathbb{R}^{n-2}$, and  component-wise Dirichlet energy minimizing in $A_{1,1}(0)$ (as in Definition~\ref{v comp minimizing defn}).  By the inductive hypothesis and \eqref{graphical annuli eqnA2}, \eqref{graphical annuli conclA1} holds true with $N^{(\nu)},\widetilde{\varphi}^{(\nu)},P_{\nu},S_{\nu},\zeta_{\nu},\rho_{\nu}$ in place of $N,\varphi,P_{Z_0},S,\zeta,\rho$.  Moreover, by the triangle inequality and \eqref{graphical annuli hyp} 
\begin{align*}
	\rho^{\alpha-1} \|\mathcal{G}(\varphi^{(\nu)},\widetilde{\varphi}^{(\nu)})\|_{L^2(B_1(0,P_{\nu}))} 
	\leq\,& C \|\mathcal{G}(N^{(\nu)},\varphi^{(\nu)})\|_{L^2(A^{\mathcal M}_{\rho,1}(\zeta))} 
		+ C \|\mathcal{G}(N^{(\nu)},\widetilde{\varphi}^{(\nu)})\|_{L^2(A^{\mathcal M}_{\rho,1}(\zeta))} 
	\\ \leq\,& 2C \widetilde{E}_{\nu} , 
\end{align*}
where $C = C(n,q,\alpha,\kappa) \in (0,\infty)$ is a constant.  Thus we can write $\varphi^{(\nu)} = \sum_{j=1}^{p_0} m^{(\nu)}_{j,k} \varphi^{(\nu)}_{j,k}$ where $\varphi^{(\nu)}_{j,k}$ are the distinct $q^{(\nu)}_{j,k}$-valued components of $\varphi^{(\nu)}$ with multiplicity $m^{(\nu)}_{j,k}$ (as in Definition~\ref{Phi-p_defn}) and $\|\mathcal{G}(\varphi^{(\nu)}_{j,k},\widetilde{\varphi}^{(\nu)}_j)\|_{L^{\infty}(B_1(0,P_{\nu}))} \leq C \widetilde{E}_{\nu}$ if $\widetilde{\varphi}^{(\nu)}_j$ is non-zero and $\|\varphi^{(\nu)}_{j,k}\|_{L^{\infty}(B_1(0,P_{\nu}))} \leq C \widetilde{E}_{\nu}$ if $\widetilde{\varphi}^{(\nu)}_j$ is identically zero, where $C = C(n,m,q,\alpha,\kappa) \in (0,\infty)$ is a constant.  It follows from \eqref{graphical annuli eqnA1} that upon blowing-up, there exists $j_1 \in \{1,2,\ldots,J\}$ such that the graph of $\psi_{j_1}$ has a self-intersection along an $(n-1)$-dimensional half-plane.  By Remark~\ref{graphical rmk}(4), $\psi_{j_1}$ is not a single-valued function with constant integer multiplicity.  Therefore, $\psi$ is component-wise Dirichlet energy minimizing in $A_{1,1}(0)$ whose singular set contains an $(n-1)$-dimensional half-plane, giving us a contradiction.  Therefore, \eqref{graphical annuli conclA1} must hold true.  

To see \eqref{graphical annuli conclA2}, let $\widetilde{\beta} = \widetilde{\beta}(n,m,q,p_0,M_0,N_0,C_e,C_h,\alpha) \in (0,1)$ be a suitably small constant.  Choose $\widetilde{p} \in \{0,1,2,\ldots,p_0-1\}$ and $\widetilde{\varphi} \in \Phi_{\alpha,q,\widetilde{p}}(P_{Z_0},S)$ such that 
\begin{equation}\label{graphical annuli eqnA3}
	E_{\varsigma/2}(N,\widetilde{\varphi},A^{\mathcal M}_{\rho,1}(\zeta))
	\leq 2^{p_0} \widetilde{\beta}^{1-p_0} \inf_{\varphi' \in \bigcup_{p'=0}^{p_0-1} \Phi_{\alpha,q,p'}(P_{Z_0},S)} 
		E_{\varsigma/2}(N,\varphi',A^{\mathcal M}_{\rho,1}(\zeta)) 
\end{equation}
and either $\widetilde{p} = 0$ or $\widetilde{p} > 0$ and 
\begin{equation*}
	E_{\varsigma/2}(N,\widetilde{\varphi},A^{\mathcal M}_{\rho,1}(\zeta))
	\leq \widetilde{\beta} \inf_{\varphi' \in \bigcup_{p'=0}^{\widetilde{p}-1} \Phi_{\alpha,q,p'}(P_{Z_0},S)} 
		E_{\varsigma/2}(N,\varphi',A^{\mathcal M}_{\rho,1}(\zeta)) . 
\end{equation*}
Let $\widetilde{\varphi} = \sum_{j=1}^{\widetilde{p}} \widetilde{m}_j \widetilde{\varphi}_j$, where $\widetilde{\varphi}_j$ are the distinct components of $\widetilde{\varphi}_j$ with multiplicity $\widetilde{m}_j$ (as in Definition~\ref{Phi-p_defn}).  By the inductive hypothesis and \eqref{graphical annuli eqnA3}, \eqref{graphical annuli conclA2} holds true with $\widetilde{\varphi}$ in place of $\varphi$.  By the triangle inequality and \eqref{graphical annuli hyp}, $\|\mathcal{G}(\varphi,\widetilde{\varphi})\|_{L^2(B_1(0,P_{Z_0}))} \leq C E_{\varsigma/2}(N,\widetilde{\varphi},A^{\mathcal M}_{\rho,1}(\zeta))$, where $C = C(n,q,\alpha,\kappa) \in (0,\infty)$ is a constant.  Thus we can write $\varphi = \sum_{j=1}^{p_0} m_{j,k} \varphi_{j,k}$ where $\varphi_{j,k}$ are the distinct components of $\varphi$ with multiplicity $m_{j,k}$ (as in Definition~\ref{Phi-p_defn}) and $\|\mathcal{G}(\varphi_{j,k},\widetilde{\varphi}_j)\|_{L^{\infty}(B_1(0,P_{Z_0}))} \leq C E_{\varsigma/2}(N,\widetilde{\varphi},A^{\mathcal M}_{\rho,1}(\zeta))$ if $\widetilde{\varphi}_j$ is non-zero and $\|\varphi_{j,k}\|_{L^{\infty}(B_1(0,P_{Z_0}))} \leq C E_{\varsigma/2}(N,\widetilde{\varphi},A^{\mathcal M}_{\rho,1}(\zeta))$ if $\widetilde{\varphi}_j$ is identically zero, where $C = C(n,m,q,\alpha,\kappa) \in (0,\infty)$ is a constant.  Thus using \eqref{graphical annuli conclA1}, we can argue as in~\cite[Lemma~5.9]{KrumWicc} that \eqref{graphical annuli conclA2} holds true.

\noindent\textit{Proof of Conclusion~(B).}  Let $X_0 \in A^{\mathcal M}_{\rho,(1+\kappa)/2}(\zeta)$ and set $\rho' = (1-\kappa)\rho/256$.  Let $\varphi_{j,l}$ be as in \eqref{varphi_localized} with $B = B^{\mathcal M}_{\rho'}(X_0)$.  Since $\varphi$ is homogeneous degree $\alpha$ and harmonic and $\zeta$ is in the spine of $\varphi$, by Remark~\ref{graphical rmk}(2)
\begin{align}\label{graphical annuli eqnB1}
	&\rho^{-1} \|\varphi\|_{C^0(A_{\rho,1}(Y))} + \|D\varphi\|_{C^0(A_{\rho,1}(Y))} + \rho \|D^2\varphi\|_{C^0(A_{\rho,1}(Y))} 
	\\ \leq\,& C \rho^{-n/2-1} \|\varphi\|_{L^2(A_{\rho,1}(Y))} 
	\leq C \rho^{\alpha-1} \|\varphi\|_{L^2(B_1(0,P_{Z_0}))} 
	\leq C \widehat{E} , \nonumber 
\end{align} 
where $\widehat{E} = \rho^{-n/2-1} \|N\|_{L^2(B^{\mathcal M}_{\rho}(Z_0))}$ and $C = C(n,m,q,\alpha) \in (0,\infty)$ are constants.  By Lemma~\ref{app lemma}, for each $j \in \{1,2,\ldots,J\}$ and $l \in \{1,2,\ldots,q_j\}$ there exists a (unique) single-valued function $\psi_{j,l} : B^{\mathcal M}_{3\rho'/4}(X_0) \rightarrow {\rm T}\mathcal{M}^{\perp}$ such that 
\begin{equation}\label{graphical annuli eqnB2}
	\mathbf{G}_{\psi_{j,l}} = \mathbf{G}_{\varphi_{j,l} \circ \exp_{Z_0}^{-1}} \,\llcorner\, \mathbf{C}^{\mathcal M}_{3\rho'/4}(X_0) . 
\end{equation} 

We want to apply Lemma~\ref{app lemma} to obtain $C^2$-estimates for $\psi_{j,l}$.  For each $j \in \{1,2,\ldots,J\}$ and $l \in \{1,2,\ldots,q_j\}$, define $\xi_{j,l} : B^{\mathcal M}_{3\rho'/4}(X_0) \rightarrow \mathcal{M}$ be the unique injective map such that 
\begin{equation}\label{graphical annuli eqnB3}
	X + \psi_{j,l}(X) = \xi_{j,l}(X) + \varphi_{j,l}(\exp_{Z_0}^{-1}(\xi_{j,l}(X))) 
\end{equation}
for each $X \in B^{\mathcal M}_{3\rho'/4}(X_0)$.  Since $\|D\bs{\varphi}\|_{C^{2,\kappa}} \leq C \mathbf{m}_0^{1/2}$ 
\begin{equation}\label{graphical annuli eqnB4}
	\|\pi_{{\rm T}_X \mathcal{M}} - \pi_{{\rm T}_Y \mathcal{M}}\| \leq C \mathbf{m}_0^{1/2} |X - Y|  
\end{equation}
for all $X,Y \in B^{\mathcal M}_{\rho_0}(Z_0)$, where $C = C(n,m,q,M_0,N_0,C_e,C_h) \in (0,\infty)$ are constants.  Let $X \in B^{\mathcal M}_{\rho_0}(Z_0)$.  Since $\xi_{j,l}(X) \in \mathcal{M}$ and $X$ is the closest point in $\mathcal{M}$ to $X+\psi_{j,l}(X)$ 
\begin{equation}\label{graphical annuli eqnB5}
    |\psi_{j,l}(X)| \leq |\varphi_{j,l}(\exp_{Z_0}^{-1}(\xi_{j,l}(X)))| \leq C \widehat{E} \rho , 
\end{equation} 
where $C = C(n,m,q,\alpha) \in (0,\infty)$.  Hence by \eqref{graphical annuli eqnB3}
\begin{equation}\label{graphical annuli eqnB6}
	|X - \xi_{j,l}(X)| \leq |\psi_{j,l}(X)| + |\varphi_{j,l}(\xi_{j,l}(X))| \leq 2 |\varphi_{j,l}(\exp_{Z_0}^{-1}(\xi_{j,l}(X)))| \leq 2C \widehat{E} \rho .
\end{equation}
Since $\varphi_{j,l}$ takes values $P_{Z_0}^{\perp}$, by \eqref{graphical annuli eqnB4} and \eqref{graphical annuli eqnB1} 
\begin{equation}\label{graphical annuli eqnB7}
	|\pi_{{\rm T}_X \mathcal{M}}(D_{\vec\tau} (\varphi_{j,l}\circ\exp_{Z_0}^{-1})(X))| 
	\leq \|\pi_{{\rm T}_X \mathcal{M}} - \pi_{P_{Z_0}}\| \,|D(\varphi_{j,l}\circ\exp_{Z_0}^{-1})(X)| 
	\leq C \mathbf{m}_0^{1/2} \rho_0 \widehat{E} 
\end{equation}
for all $X \in B^{\mathcal M}_{\rho'}(X_0)$ and $\vec\tau \in {\rm T}_X \mathcal{M}$, where $C = C(n,m,q,M_0,N_0,C_e,C_h,\alpha) \in (0,\infty)$ is a constant.  Therefore, by Lemma~\ref{app lemma} together with $\|D\bs{\varphi}\|_{C^2(\mathcal{M})} \leq C \mathbf{m}_0^{1/2}$, \eqref{graphical annuli eqnB1}, \eqref{graphical annuli eqnB4}, \eqref{graphical annuli eqnB6}, \eqref{graphical annuli eqnB5}, and \eqref{graphical annuli eqnB7} 
\begin{align}\label{graphical annuli eqnB8} 
	&\sup_{X \in B^{\mathcal M}_{3\rho'/4}(X_0)} |D\xi_{j,l}(X) - I_{{\rm T}_X \mathcal{M}}| 
		\leq C \sup_{X \in B^{\mathcal M}_{3\rho'/4}(X_0)} \|\pi_{{\rm T}_{\xi_{j,l}(X)} \mathcal{M}} - \pi_{{\rm T}_X \mathcal{M}}\| 
		\\&+ \|A_{\mathcal M}\|_{C^0(\mathcal{M})} \|\psi_{j,l}\|_{C^0(B^{\mathcal M}_{3\rho'/4}(X_0))}  + C \|\pi_{{\rm T}\mathcal{M}}(D\varphi_{j,l})\|_{C^0(B^{\mathcal M}_{\rho'}(X_0))} 
        \leq C \mathbf{m}_0^{1/2} \rho_0 \widehat{E} \nonumber 
\end{align}
and 
\begin{align}\label{graphical annuli eqnB9} 
    &\|D\psi_{j,l}\|_{C^0(B^{\mathcal M}_{3\rho'/4}(X_0))} 
    \leq C \sup_{X \in B^{\mathcal M}_{3\rho'/4}(X_0)} \|\pi_{{\rm T}_{\xi_{j,l}(X)} \mathcal{M}} - \pi_{{\rm T}_X \mathcal{M}}\| 
	\\&+ \|A_{\mathcal M}\|_{C^0(\mathcal{M})} \,\|\psi_{j,l}\|_{C^0(B^{\mathcal M}_{3\rho'/4}(X_0))} + C \|D\varphi_{j,l}\|_{C^0(B^{\mathcal M}_{\rho'}(X_0))} 
    \leq  C \widehat{E} \leq C \mathbf{m}_0^{1/(2n)} \rho^{\beta_2} , \nonumber
\end{align}
where $C = C(n,m,q,M_0,N_0,C_e,C_h,\alpha,\kappa) \in (0,\infty)$ are constants.  The last step in \eqref{graphical annuli eqnB9} follows by using \eqref{normal approx est1} to obtain $\widehat{E} = \rho^{-n/2-1} \|N\|_{L^2(B^{\mathcal M}_{\rho}(Z))} \leq C\mathbf{m}_0^{1/(2n)} \rho^{\beta_2}$.  By Remark~\ref{graphical rmk}(3) 
\begin{equation}\label{graphical annuli eqnB10}
	\rho \|\Delta_{\mathcal M} (\varphi_{j,l} \circ \exp_{Z_0}^{-1})\|_{C^0(B^{\mathcal M}_{3\rho'/4}(X_0))} 
	\leq C \mathbf{m}_0^{1/2} \rho_0 \widehat{E} , 
\end{equation}
where $C = C(n,m,q,M_0,N_0,C_e,C_h,\alpha,\kappa) \in (0,\infty)$ is a constant.  Therefore, by Lemma~\ref{app lemma} together with $\|D\bs{\varphi}\|_{C^2(\mathcal{M})} \leq C \mathbf{m}_0^{1/2}$, \eqref{graphical annuli eqnB1}, \eqref{graphical annuli eqnB4}, \eqref{graphical annuli eqnB6}, \eqref{graphical annuli eqnB5}, \eqref{graphical annuli eqnB8}, \eqref{graphical annuli eqnB9}, and \eqref{graphical annuli eqnB10}
\begin{align}\label{graphical annuli eqnB11} 
    &\rho \|\Delta_{\mathcal M} \psi_{j,l}\|_{C^0(B^{\mathcal M}_{3\rho'/4}(X_0))} 
		\leq C \rho \|A_{\mathcal M}\|_{C^1(\mathcal{M})} \,\|\xi_{j,l} - {\rm Id}_{\mathcal M}\|_{C^0(B^{\mathcal M}_{3\rho'/4}(X_0))} 
		\\&+ C \rho \|A_{\mathcal M}\|_{C^1(\mathcal{M})} \,\|\psi_{j,l}\|_{C^1(B^{\mathcal M}_{3\rho'/4}(X_0))} + C \rho \|\Delta_{\mathcal M} (\varphi_{j,l} \circ \exp_{Z_0}^{-1})\|_{C^0(B^{\mathcal M}_{3\rho'/4}(X_0))} \nonumber
		\\&+ C \sup_{X \in B^{\mathcal M}_{3\rho'/4}(X_0)} |D\xi_{j,l}(X) - I_{{\rm T}_X \mathcal{M}}| \leq C \mathbf{m}_0^{1/2} \rho_0 \widehat{E} , \nonumber 
\end{align}
where $C = C(n,m,q,M_0,N_0,C_e,C_h,\alpha,\kappa) \in (0,\infty)$ is a constant.  Since each $\varphi_{j,l}$ takes values in $P_{Z_0}^{\perp}$, by \eqref{graphical annuli eqnB4}  
\begin{align*}
    &|\pi_{{\rm T}_X \mathcal{M}}(D_{\vec\tau} (\varphi_{j,l}\circ\exp_{Z_0}^{-1})(X) - D_{\vec\tau} (\varphi_{j',l'}\circ\exp_{Z_0}^{-1})(X))| 
    \\ \leq\,& \|\pi_{{\rm T}_X \mathcal{M}} - \pi_{P_{Z_0}}\| \,|D_{\vec\tau} (\varphi_{j,l}\circ\exp_{Z_0}^{-1})(X) - D_{\vec\tau} (\varphi_{j',l'}\circ\exp_{Z_0}^{-1})(X)|
    \\ \leq\,& C \mathbf{m}_0^{1/2} \rho_0 |D_{\vec\tau} (\varphi_{j,l}\circ\exp_{Z_0}^{-1})(X) - D_{\vec\tau} (\varphi_{j',l'}\circ\exp_{Z_0}^{-1})(X)|
\end{align*}
for all $X \in B^{\mathcal M}_{\rho}(X_0)$, $\vec\tau \in {\rm T}_X \mathcal{M}$, $j,j' \in \{1,2,\ldots,J\}$, $l \in \{1,2,\ldots,q_j\}$, and $l' \in \{1,2,\ldots,q_{j'}\}$, where $C = C(n,m,q,M_0,N_0,C_e,C_h,\alpha,\kappa)$ is a constant.  Hence by Lemma~\ref{app lemma2} and \eqref{graphical annuli conclA2} 
\begin{equation}\label{graphical annuli eqnB20}
	\rho \,\|D\psi_{j,k} - D\psi_{j,l}\|_{C^0(B^{\mathcal M}_{3\rho'/4}(X_0))} 
	\leq C \inf_{B^{\mathcal M}_{3\rho'/4}(X_0)} |\psi_{j,l} - \psi_{j',l'}| , 
\end{equation}
where $C = C(n,m,q,p_0,M_0,N_0,C_e,C_h,\alpha,\kappa) \in (0,\infty)$ is a constant.

Observe that for each $X \in A^{\mathcal M}_{\rho,1}(\zeta)$, $N(X) = \sum_{j=1}^p \sum_{l=1}^{q_j} \sum_{h=1}^{m_j} \llbracket N_{j,l,h}(X) \rrbracket$ for some $N_{j,l,h}(X) \in {\rm T}_X \mathcal{M}^{\perp}$ such that 
\begin{equation*}
    \mathcal{G}(N(X),\varphi(\exp_{Z_0}^{-1}(X)))^2 = \sum_{j=1}^p \sum_{l=1}^{q_j} \sum_{h=1}^{m_j} |N_{j,l,h}(X) - \varphi_{j,l}(\exp_{Z_0}^{-1}(X))|^2 .
\end{equation*}
Provided $\varepsilon$ is sufficiently small, since $\varphi_{j,l}$ and $\psi_{j,l}$ satisfy \eqref{graphical annuli eqnB2} and both have small Lipschitz constants ($\leq C\widehat{E} \leq C\varepsilon^{\beta_2}$)
\begin{align}\label{graphical annuli eqnB21}
    \mathcal{G}(N(X),\psi(X))^2 \leq\,& \sum_{j=1}^p \sum_{l=1}^{q_j} \sum_{h=1}^{m_j} |N_{j,l,h}(X) - \psi_{j,l}(X))|^2 \\ \leq\,& 4 \sum_{j=1}^p \sum_{l=1}^{q_j} \sum_{h=1}^{m_j} |N_{j,l,h}(X) - \varphi_{j,l}((\exp_{Z_0}^{-1}(X))|^2 \nonumber \\=\,& 4 \,\mathcal{G}(N(X),\varphi(\exp_{Z_0}^{-1}(X)))^2 . \nonumber 
\end{align}

Since $X_0 \in A^{\mathcal M}_{\rho,(1+\kappa)/2}(\zeta)$, $B^{\mathcal M}_{\rho'}(X_0) \subset A^{\mathcal M}_{\rho,(3+\kappa)/4}(\zeta)$.  Hence by \eqref{graphical annuli go condition}, \eqref{relharm separation main go condition} holds true with $X_0,\rho',1/2$ in place of $Z, \rho, \gamma$.  In light of \eqref{graphical annuli eqnB20}, we can apply Theorem~\ref{relharm separation main thm} (with $X_0, \rho', 1/2$ in place of $Z, \rho, \gamma$) together with \eqref{graphical annuli eqnB9}, \eqref{graphical annuli eqnB11}, and \eqref{graphical annuli eqnB21} to obtain 
\begin{align}\label{graphical annuli eqnB22}
	\sup_{X \in \op{spt} T \cap \mathbf{C}^{\mathcal M}_{\rho'/2}(X_0)} \op{dist}(X, \op{spt} \mathbf{G}_{\varphi \circ \exp_{Z_0}^{-1}}) 
	\leq C \rho \,E_{\varsigma/2}(N,\varphi,A^{\mathcal M}_{\rho,1}(\zeta)) 
\end{align}
where $C = C(n,m,q,p_0,M_0,N_0,C_e,C_h,\alpha,\kappa) \in (0,\infty)$ is a constant.  That is, 
\begin{equation}\label{graphical annuli eqnB23} 
	\sup_{X \in \op{spt} T \cap \mathbf{A}^{\mathcal M}_{\rho,(1+\kappa)/2}(\zeta)} \op{dist}(X,\op{spt} \mathbf{G}_{\varphi \circ \exp_{Z_0}^{-1}}) 
	\leq C \rho \,E_{\varsigma/2}(N,\varphi,A^{\mathcal M}_{\rho,1}(\zeta)) 
\end{equation}
(where $C$ is as in \eqref{graphical annuli eqnB22}).  By \eqref{graphical annuli hyp} and \eqref{graphical annuli conclA1}, $\rho \,E_{\varsigma/2}(N,\varphi,A^{\mathcal M}_{\rho,1}(\zeta)) \leq c^{-1} \overline{\beta} \inf_{A_{\rho,(1+\kappa)/2}(\zeta)} \op{sep} \varphi$ (where $c$ is as in \eqref{graphical annuli conclA1}).  Thus provided $\overline{\beta}$ is sufficiently small, by \eqref{graphical annuli eqnB23}, there exists (unique) $n$-dimensional locally area-minimizing rectifiable currents $T_j$ of $\mathbf{A}^{\mathcal M}_{\rho,(1+\kappa)/2}(\zeta)$ such that \eqref{graphical annuli conclB1} and \eqref{graphical annuli conclB3} hold true.  Moreover, we can argue as in Remark~\ref{relharm separation1 rmk} that there exists a set $\Omega_j \subseteq \mathcal{K} \cap A^{\mathcal M}_{\rho,(1+\kappa)/2}(\zeta)$ such that $\mathcal{H}^n(\mathcal{K} \cap A^{\mathcal M}_{\rho,(1+\kappa)/2}(\zeta) \setminus \Omega_j) \leq c^{-2} \overline{\beta}^2 \rho^n$ (where $c$ is as in \eqref{graphical annuli conclA1}) and $\mathcal{G}(N_j,\varphi_j) < \frac{1}{3} \op{sep} \varphi(X)$ on $\Omega_j$.  For each $X \in \Omega_j$, the slice $\langle T_j, \mathbf{p}, X \rangle$ is the sum of precisely the $q_j m_j$ values of $N(X)$ (counting multiplicity) which are within a distance $<  \frac{1}{3} \op{sep} \varphi(X)$ from the values of $\varphi_j(X)$.  It follows using the constancy theorem that \eqref{graphical annuli conclB2} holds true.

\noindent\textit{Proof of Conclusion~(C).}  Let $X_0 \in A^{\mathcal M}_{\rho,(1+\kappa)/2}(\zeta)$ and set $\rho' = (1-\kappa)\rho/256$.  Let $T_{j,l}$ be as in \eqref{T_separated eqn4} where $B = B^{\mathcal M}_{\rho'/2}(Z)$.  By Lemma~\ref{relharm normal approx lemma}, there exist Lipschitz $m_j$-valued function $N_{X_0,j,l} : B^{\mathcal M}_{\rho'/4}(Z) \rightarrow \mathcal{A}_{m_j}({\rm T}\mathcal{M}^{\perp})$ and closed sets $\mathcal{K}_{X_0,j,l} \subseteq B^{\mathcal M}_{\rho'/4}(Z)$ such that \eqref{relharm normal approx concl1} and \eqref{relharm normal approx concl2} hold true with $1/4,\rho',N_{X_0,j,l},\mathcal{K}_{X_0,j,l}$ in place of $\gamma,\rho,N_j,\mathcal{K}_j$.  By replacing $\mathcal{K}_{X_0,j,l}$ with $\mathcal{K}_{X_0,j} = \bigcap_{l'=1}^{q_j} \mathcal{K}_{X_0,j,l'}$, assume that $\mathcal{K}_{X_0,j,l} = \mathcal{K}_{X_0,j}$ independent of $l$.  Define $v_{X_0,j,l} : \op{spt} \mathbf{G}_{\varphi_{j,l} \circ \exp_{Z_0}^{-1}} \cap \mathbf{C}^{\mathcal M}_{\rho'/4}(Z) \rightarrow \mathcal{A}_{m_j}(\mathbb{R}^{n+m})$ by 
\begin{equation*}
	v_{X_0,j,l}(X,\varphi_{j,l}(\exp_{Z_0}^{-1}(X))) = \sum_{h=1}^{m_j} \llbracket N_{X_0,j,l,h}(X) - \varphi_{j,l}(\exp_{Z_0}^{-1}(X)) \rrbracket 
\end{equation*}
for each $X \in B^{\mathcal M}_{\rho'/4}(Z)$, where $N_{X_0,j,l}(X) = \sum_{h=1}^{m_j} \llbracket N_{X_0,j,l,h}(X) \rrbracket$.  Using a partition of unity argument (like in the proof of~\cite[Theorem~3.8]{KrumWicb}), we can construct $v_j : \op{spt} \mathbf{G}_{\varphi_{j,l} \circ \exp_{Z_0}^{-1}} \cap \mathbf{A}^{\mathcal M}_{\rho,\kappa}(\zeta) \rightarrow \mathcal{A}_{m_j}(\mathbb{R}^{n+m})$ and a closed set $\mathcal{K}_j \subseteq \mathbf{A}^{\mathcal M}_{\rho,\kappa}(\zeta)$ such that \eqref{graphical annuli conclC1} and \eqref{graphical annuli conclC2} holds true, where $N_j$ is as in \eqref{N of v defn}.  By \eqref{graphical annuli eqnB23} and truncation, we can guarantee the supremum estimate for $v_j$ in \eqref{graphical annuli conclC3}.  By Lemma~\ref{relharm energy lemma} and \eqref{graphical annuli eqnB10}, we obtain the gradient $L^2$-estimate for $v_j$ in \eqref{graphical annuli conclC3}.
\end{proof}

\begin{proof}[Proof of Theorem~\ref{graphical_thm}] 
Let $\tau \in (0,1/2]$.  Let $Z_0 \in \mathcal{B}^{(\geq 2)}_q$ and $\rho_0 \in (0,1]$.  We claim that for all $\sigma \in (0,2]$, $\sigma^{-n/2-1} \|N\|_{L^2(B^{\mathcal M}_{\sigma}(Z_0))} \leq C\sigma$ where $C = C(n,m,q) \in (0,\infty)$ is a constant (independent of $C_e$ and $C_h$).  To see this, notice that by \eqref{normal approx est1} and \eqref{normal approx est2} with $b = 1$ and $\phi$ equal to the characteristic function on $B^{\mathcal M}_2(Z_0)$
\begin{equation*}
	\int_{B^{\mathcal M}_2(Z_0)} |N|^2 \leq q \int_{B^{\mathcal M}_2(Z_0)} |N_a|^2 + \int_{B^{\mathcal M}_2(Z_0)} |N_f|^2 \leq C\mathbf{m}_0^{1+1/(2n)} + \int_{B^{\mathcal M}_2(Z_0)} |N_f|^2 , 
\end{equation*}
where $N_a$ is the single-valued average of $N$, $N_f$ is the average-free part of $N$, and $C = C(n,m,q,M_0,$ $N_0,C_e,C_h) \in (0,\infty)$ is a constant.  Express $P_0 \cap B^{\mathcal M}_3(0) = \{ X + \lambda(Z) : X \in B^{\mathcal M}_3(0) \}$ for some $C^3$-function $\lambda : B^{\mathcal M}_1(0) \rightarrow {\rm T}\mathcal{M}^{\perp}$ with $\|\lambda\|_{C^3(B^{\mathcal M}_1(0))} \leq C\mathbf{m}_0^{1/(2n)}$ (where $P_0 = \mathbb{R}^n \times \{0\}$).  Arguing as in the proof of Lemma~\ref{planar vs alm freq lemma} 
\begin{equation*}
	\int_{B^{\mathcal M}_2(Z_0)} |N_f|^2 \leq \int_{B^{\mathcal M}_2(Z_0)} \mathcal{G}(N,q \llbracket \lambda \rrbracket)^2 \leq 4 \int_{B^{\mathcal M}_2(Z_0)} \sum_{l=1}^q \op{dist}^2(X+N_l(X),P_0) , 
\end{equation*}
and so 
\begin{equation*}
	\int_{B^{\mathcal M}_2(Z_0)} |N|^2 \leq C\mathbf{m}_0^{1+1/(2n)} + 4 \int_{B^{\mathcal M}_2(Z_0)} \sum_{l=1}^q \op{dist}^2(X+N_l(X),P_0) , 
\end{equation*}
where $C = C(n,m,q,M_0,N_0,C_e,C_h) \in (0,\infty)$ is a constant.  Assuming $\varepsilon_2$ is sufficiently small, by the area formula, \eqref{normal approx est1}, and \eqref{fast decay eqn3}
\begin{equation}\label{graphical thm eqn01}
	\int_{B^{\mathcal M}_2(Z_0)} |N|^2 \leq C\mathbf{m}_0^{1+\min\{1/(2n),\gamma_2\}} + 8 \int_{\mathbf{C}^{\mathcal M}_2(Z_0)} \sum_{l=1}^q \op{dist}^2(X,P_0) \,d\|T\|(X) \leq C_1 \mathbf{m}_0 , 
\end{equation}
where $C = C(n,m,q,M_0,N_0,C_e,C_h) \in (0,\infty)$ and $C_1 = C_1(n,m,q) \in (0,\infty)$ are constants.  By Lemma~\ref{planar vs alm freq lemma}, $\mathbf{I}_{N,Z}(0^+) \geq 2$.  Hence provided $\varepsilon_2$ is sufficiently small, by Corollary~\ref{doubling centered cor} and \eqref{graphical thm eqn01}
\begin{equation}\label{graphical thm eqn4}
	\frac{1}{\sigma^{n+2}} \int_{B^{\mathcal M}_{\sigma}(Z_0)} |N|^2 
	\leq C \sigma^2 \int_{B^{\mathcal M}_2(Z_0)} |N|^2 \leq C \mathbf{m}_0 \sigma^2
\end{equation}
for all $\sigma \in (0,2]$, where $C = C(n,m,q)$ are constants.  

Now let $\rho \geq \tau\rho_0/4$ and $\zeta \in S$ such that $\rho^2 + |\zeta|^2 \leq (15\rho_0/16)^2$.  We want to verify that the hypotheses \eqref{graphical annuli go condition} (with $\kappa = 1/2$) and \eqref{graphical annuli hyp} of Lemma~\ref{graphical_lemma} hold true.  To see that \eqref{graphical annuli go condition} holds true (with $\kappa = 1/2$), suppose to the contrary that there exists $Y \in A^{\mathcal M}_{\rho,7/8}(\zeta)$ such that $\overline{R}(Y) \geq 4^{-q-8} \rho$.  By Lemma~\ref{relharm critrad2 lemma} (with $Y$ in place of $Z$) 
\begin{equation}\label{graphical thm eqn5} 
	\min\{C_e,C_h\} \cdot \mathbf{m}_0 \overline{R}^{2-2\delta_2} \leq \frac{C}{\overline{R}^{n+2}} \int_{B^{\mathcal M}_{\overline{R}}(Y)} |N|^2 , 
\end{equation}
where $\overline{R} = \overline{R}(Y)$ and $C = C(n,m,q,M_0) \in (0,\infty)$ is a constant.  Note that $Y \in B^{\mathcal M}_{61\rho_0/64}(Z_0)$.  Provided $C_e$ and $C_h$ are sufficiently large, by \eqref{graphical thm eqn4} with $\sigma = 61\rho_0/64 + \overline{R}$ we have that $\overline{R} < \rho_0/64$ and in particular $B^{\mathcal M}_{\overline{R}}(Y) \subset B^{\mathcal M}_{\rho_0}(Z_0)$.  Using the triangle inequality and Hypothesis~$(\star\star)$ as in Remark~\ref{graphical rmk}(2), we obtain $\|\varphi\|_{L^2(B_1(0,P_{Z_0}))} \leq 4 \rho_0^{-n/2-\alpha} \|N\|_{L^2(B^{\mathcal M}_{\rho_0}(Z_0))}$.  Since $\varphi$ is homogeneity degree $\alpha$ and translation invariant along $S$, 
\begin{align*}
	\frac{1}{\overline{R}^{n+2}} \int_{B_{\overline{R}}(\exp_{Z_0}^{-1}(Y),P_{Z_0})} |\varphi|^2 
	\leq\,& \frac{C}{\overline{R}^{n+2}} \int_{B_{2\rho+\overline{R}}(\zeta,P_{Z_0})} |\varphi|^2 
	= C (2\rho+\overline{R})^{2\alpha-2} \int_{B_1(0,P_{Z_0})} |\varphi|^2 
	\\ \leq\,& 8 C \bigg(\frac{2\rho+\overline{R}}{\rho_0}\bigg)^{2\alpha-2} \cdot \frac{1}{\rho_0^{n+2}} \int_{B^{\mathcal M}_{\rho_0}(Z_0)} |N|^2 , \nonumber 
\end{align*}
where $C = C(n,q) \in (0,\infty)$ is a constant.  By triangle inequality, $\rho \geq \tau\rho_0/4$, and Hypothesis~$(\star\star)$
\begin{align*} 
	\frac{1}{\overline{R}^{n+2}} \int_{B^{\mathcal M}_{\overline{R}}(Y)} |N|^2 
	\leq\,& \frac{C}{\overline{R}^{n+2}} \int_{B_{\overline{R}}(\exp_{Z_0}^{-1}(Y),P_{Z_0})} |\varphi|^2 
		+ \frac{C}{\rho^{n+2}} \int_{B^{\mathcal M}_{\rho_0}(Z_0)} \mathcal{G}(N,\varphi \circ \exp_{Z_0}^{-1})^2 
	\\ \leq\,& C \left( \bigg(\frac{2\rho+\overline{R}}{\rho_0}\bigg)^{2\alpha-2}  
		+ \frac{\beta_0^2}{\tau^{n+2\alpha}} \bigg(\frac{\rho}{\rho_0}\bigg)^{2\alpha-2} \right) 
		\frac{1}{\rho_0^{n+2}} \int_{B^{\mathcal M}_{\rho_0}(Z_0)} |N|^2 , 
\end{align*}
where $C = C(n,q) \in (0,\infty)$ are constants.  Assuming $\beta_0$ is small enough that $\beta_0^2 \tau^{-n-2\alpha} \leq 1$ 
\begin{align*} 
	\frac{1}{\overline{R}^{n+2}} \int_{B^{\mathcal M}_{\overline{R}}(Y)} |N|^2 
	\leq 2C \bigg(\frac{2\rho+\overline{R}}{\rho_0}\bigg)^{2\alpha-2} \cdot \frac{1}{\rho_0^{n+2}} \int_{B^{\mathcal M}_{\rho_0}(Z_0)} |N|^2 .
\end{align*}
By $\alpha \geq 2$ and \eqref{graphical thm eqn4} 
\begin{align}\label{graphical thm eqn7} 
	\frac{1}{\overline{R}^{n+2}} \int_{B^{\mathcal M}_{\overline{R}}(Y)} |N|^2 \leq C \mathbf{m}_0 \overline{R}^2
\end{align}
where $C = C(n,q) \in (0,\infty)$ are constants.  Assuming $C_e$ and $C_h$ are sufficiently large, \eqref{graphical thm eqn5} and \eqref{graphical thm eqn7} give us a contradiction. 

To see \eqref{graphical annuli hyp}, let $p' \in \{0,1,2,\ldots,p-1\}$ and $\varphi' \in \Phi_{\alpha,q,p'}(P_{Z_0},S)$.  By the triangle inequality and Hypothesis~$(\star\star)$ 
\begin{align*} 
	E_{\varsigma/2}(N,\varphi',B^{\mathcal M}_{\rho_0}(Z_0)) 
	\leq\,& E_{\varsigma/2}(N,\varphi,B^{\mathcal M}_{\rho_0}(Z_0)) + C \rho_0^{\alpha-1} \|\mathcal{G}(\varphi,\varphi')\|_{L^2(B_1(0,P_{Z_0}))} 
	\\ \leq\,& \beta_0 E_{\varsigma/2}(N,\varphi',B^{\mathcal M}_{\rho_0}(Z_0)) + C \rho_0^{\alpha-1} \|\mathcal{G}(\varphi,\varphi')\|_{L^2(B_1(0,P_{Z_0}))} 
\end{align*}
and thus provided $\beta_0 < 1/2$ 
\begin{equation}\label{graphical thm eqn8} 
	E_{\varsigma/2}(N,\varphi',B^{\mathcal M}_{\rho_0}(Z_0)) \leq 2C \rho_0^{\alpha-1} \|\mathcal{G}(\varphi,\varphi')\|_{L^2(B_1(0,P_{Z_0}))} , 
\end{equation}
where $C = C(n,m,q,p,M_0,N_0,C_e,C_h) \in (0,\infty)$ is a constant.  On the other hand, by again using the triangle inequality 
\begin{equation}\label{graphical thm eqn9} 
	\|\mathcal{G}(\varphi,\varphi')\|_{L^2(B_1(0,P_{Z_0}))} 
	\leq C E_{\varsigma/2}(N,\varphi,A^{\mathcal M}_{\rho,1}(\zeta)) + C E_{\varsigma/2}(N,\varphi',A^{\mathcal M}_{\rho,1}(\zeta)) , 
\end{equation}
where $C = C(n,m,q,p,M_0,N_0,C_e,C_h) \in (0,\infty)$ is a constant.  By $\rho \geq \tau\rho_0/4$, $A^{\mathcal M}_{\rho,1}(\zeta) \subset B^{\mathcal M}_{\rho_0}(Z_0)$, Hypothesis~$(\star\star)$, \eqref{graphical thm eqn8}, and \eqref{graphical thm eqn9}  
\begin{align*} 
	E_{\varsigma/2}(N,\varphi,A^{\mathcal M}_{\rho,1}(\zeta)) 
	\leq\,& C \tau^{-n/2-1} E_{\varsigma/2}(N,\varphi,B^{\mathcal M}_{\rho_0}(Z_0)) 
	\\ \leq\,& C \beta_0 \tau^{-n/2-1} E_{\varsigma/2}(N,\varphi',B^{\mathcal M}_{\rho_0}(Z_0)) 
	\\ \leq\,& C \beta_0 \tau^{-n/2-1} \|\mathcal{G}(\varphi,\varphi')\|_{L^2(B_1(0,P_{Z_0}))} 
	\\ \leq\,& C \beta_0 \tau^{-n/2-1} (E_{\varsigma/2}(N,\varphi,A^{\mathcal M}_{\rho,1}(\zeta)) + E_{\varsigma/2}(N,\varphi',A^{\mathcal M}_{\rho,1}(\zeta))) 
\end{align*}
where $C = C(n,m,q,p,M_0,N_0,C_e,C_h) \in (0,\infty)$ are constants.  Assuming $\beta_0$ is small enough that $C \beta_0 \tau^{-n/2-1} < 1/2$, we conclude that 
\begin{equation}\label{graphical thm eqn10} 
	E_{\varsigma/2}(N,\varphi,A^{\mathcal M}_{\rho,1}(\zeta)) \leq 2C \beta_0 \tau^{-n/2-1} 
		\inf_{\varphi' \in \bigcup_{p'=0}^{p-1} \Phi_{\alpha,q,p'}(P_{Z_0},S)} E_{\varsigma/2}(N,\varphi',A^{\mathcal M}_{\rho,1}(\zeta)) 
\end{equation}
as in \eqref{graphical annuli hyp}.  Note that provide $\beta_0$ is small enough that $C \beta_0 \tau^{-n/2-1} < 1/2$ (where $C$ is as in \eqref{graphical thm eqn10}), it follows from \eqref{graphical thm eqn9} and \eqref{graphical thm eqn10} that 
\begin{equation}\label{graphical thm eqn11} 
	\|\mathcal{G}(\varphi,\varphi')\|_{L^2(B_1(0,P_{Z_0}))} 
	\leq C \inf_{\varphi' \in \bigcup_{p'=0}^{p-1} \Phi_{\alpha,q,p'}(P_{Z_0},S)} E_{\varsigma/2}(N,\varphi',A^{\mathcal M}_{\rho,1}(\zeta)) ,
\end{equation}
where $C = C(n,m,q,p,M_0,N_0,C_e,C_h) \in (0,\infty)$ is a constant.  

Now having verified \eqref{graphical annuli go condition} (with $\kappa = 1/2$) and \eqref{graphical annuli hyp}, we can argue as in the proof of~\cite[Theorem~5.6]{KrumWicc} using Lemma~\ref{graphical_lemma} as follows.  By applying Lemma~\ref{graphical_lemma}(A) with $\zeta = 0$ and $\rho = \rho_0/4$ and using \eqref{graphical thm eqn8} and \eqref{graphical thm eqn11}, we deduce that Theorem~\ref{graphical_thm}(A) holds true.  By the Vitali covering lemma there is a finite collection of $\rho_k \geq \tau\rho_0/2$ and $\zeta_k \in S$ such that $\rho_k^2 +|\zeta_k|^2 \leq (15\rho_0/16)^2$, $\rho_{k+1} \leq \rho_k$, $\{A^{\mathcal M}_{\rho_k,1/8}(\zeta_k)\}$ covers $B^{\mathcal M}_{15\rho_0/16}(Z_0) \cap \{r > \tau\rho_0/2\}$, and $\{A^{\mathcal M}_{\rho_k,1/40}(\zeta_k)\}$ are pairwise disjoint.  Arguing as in the proof of~\cite[Theorem~5.6]{KrumWicc}, we can apply Theorem~\ref{graphical_thm}(B) with $\zeta = \zeta_k$ and $\rho = \rho_0$ to obtain locally area minimizing rectifiable currents $T_{k,j}$ of $\mathbf{A}^{\mathcal M}_{\rho_k,3/4}(\zeta_k)$ such that $T \llcorner \mathbf{A}^{\mathcal M}_{\rho_k,3/4}(\zeta_k) = \sum_{j=1}^J T_{k,j}$, $(\partial T_{k,j}) \llcorner \mathbf{A}^{\mathcal M}_{\rho_k,3/4}(\zeta_k) = 0$, and uniform distance of $T_{k,j}$ to the graph of $\varphi_j$ is $\leq C \rho_k^{-(n+2)(4+\varsigma)/8} E_{\varsigma/2}(N,\varphi,\mathbf{C}^{\mathcal M}_{\rho_0}(Z_0))$.  Then we can ``glue'' the currents $T_{k,j}$ together to form currents $T_j$ of $B^{\mathcal M}_{15\rho_0/16}(Z_0) \cap \{r > \tau\rho_0/2\}$ satisfying Theorem~\ref{graphical_thm}(B).  Similarly, we can apply Theorem~\ref{graphical_thm}(C) with $\zeta_0 = \zeta_k$ and $\rho = \rho_0$ to approximate $T_{k,j}$ by the graph of Lipschitz functions $v_{k,j} : \op{spt} \mathbf{G}_{\varphi_j \circ \exp_{Z_0}^{-1}} \cap \mathbf{A}^{\mathcal M}_{\rho_k,1/2}(\zeta_k) \rightarrow \mathcal{A}_{m_j}(\mathbb{R}^{n+m})$.  Using a partition of unity argument, we can ``glue'' together the functions $v_{k,j}$ to find functions $v_j : \op{spt} \mathbf{G}_{\varphi_j \circ \exp_{Z_0}^{-1}} \cap \mathbf{C}^{\mathcal M}_{7\rho_0/8}(Z_0) \cap \{r > \tau\rho_0\} \rightarrow \mathcal{A}_{m_j}(\mathbb{R}^{n+m})$ such that Theorem~\ref{graphical_thm}(C) holds true.
\end{proof}

\section{A priori estimates: Part~I}\label{sec:key_est_sec}

In this section and the next, we prove a number of integral estimates for $N$ inspired by~\cite{Sim93} and  \cite{Wic14}. These estimates play a fundamental role in the proof of the main excess decay estimate, Theorem~\ref{main excess decay thm}.  In particular, Theorem~\ref{keyest thm} below is the main estimate; it implies the remaining estimates for $N,$ and its direct consequence for the fine blow-ups (as constructed in Subsection~\ref{sec:blowup subsec}) is a key ingredient in the proof of an asymptotic decay estimate for the fine blow-ups (Theorem~\ref{blowup decay thm}).  The proof of Theorem~\ref{keyest thm} is based on an adaptation of the arguments from the authors' previous work~\cite{KrumWic2} on Dirichlet energy minimizing multi-valued functions; see also ~\cite{KrumWicc}. 

The integral estimates of this section and Section~\ref{sec:nonconcentration} hold true under the following assumption which is weaker than Hypothesis~$(\star\star)$ 
in Section~\ref{sec:graph_rep_results}.  For a suitably small constant $\beta \in (0,1)$: 

\noindent\textbf{Hypothesis~$\bs{(\dagger)}$.}  For $p \in \{1,2,\ldots, \lfloor q/q_0 \rfloor \}$, $T$, and $\varphi \in \Phi_{\alpha,q,p}(P_{Z_0},S)$ as in Hypothesis~$(\star)$ 
\begin{equation*}
	E_{\varsigma/2}(N,\varphi,B^{\mathcal M}_{\rho_0}(Z_0)) \leq \beta \rho_0^{-n/2-1} \|N\|_{L^2(B^{\mathcal M}_{\rho_0}(Z_0))} . 
\end{equation*}

To state and prove the results of this section and the next, we need the following notation.  Assume that Hypothesis~$(\star)$ holds true and in particular that $Z_0 \in \mathcal{B}^{(\geq 2)}_q$, $P_{Z_0} = {\rm T}_{Z_0}\mathcal{M}$, and $S$ be an $(n-2)$-dimensional linear subspace of $P_{Z_0}$.  Throughout this section, we shall use the following local coordinates on $\mathcal{M}$.  Let $\vec\tau_1,\vec\tau_2,\ldots,\vec\tau_n$ be an orthonormal basis for $P_{Z_0} = {\rm T}_{Z_0} \mathcal{M}$ such that $\vec\tau_3,\ldots,\vec\tau_n$ is a basis for $S$.  Parameterize $B^{\mathcal M}_1(Z_0)$ by 
\begin{equation}\label{normal coord eqn1}
	\exp_{Z_0}(x_1 \vec\tau_1 + x_2 \vec\tau_2 +\cdots + x_n \vec\tau_n) 
\end{equation}
so that $(x_1,x_2,\ldots,x_n)$ are normal coordinates of $\mathcal{M}$ at $Z_0$.  We can express each $X \in B^{\mathcal M}_1(Z_0)$ as $X = \exp_{Z_0}(x+y)$ where $x = x_1 \vec\tau_1 + x_2 \vec\tau_2$ and $y = x_3 \vec\tau_3 +\cdots + x_n \vec\tau_n$.  We define $r : B^{\mathcal M}_1(Z_0) \to \mathbb{R}$ by $r(X) = \op{dist}(\exp_{Z_0}^{-1}(X),S) = |x|$ for each $X \in B^{\mathcal M}_1(Z_0)$ (as in \eqref{dist to S defn}) and we define $R : B^{\mathcal M}_1(Z_0) \to \mathbb{R}$ so that $R(X)$ is the geodesic distance of $X \in B^{\mathcal M}_1(Z_0)$ to $Z_0$.  In normal coordinates 
\begin{equation}\label{normal coord eqn2}
	r = r(X) = \sqrt{x_1^2 + x_2^2}, \quad\quad R = R(X) = \sqrt{x_1^2+x_2^2+\cdots+x_n^2} .
\end{equation}
We shall let $\frac{\partial}{\partial x_i}$ denote partial derivatives with respect to the coordinate $x_i$.  We shall let $\widehat{r}$ and $\widehat{R}$ denote the gradients (on $\mathcal{M}$) of $r$ and $R$ respectively.  We shall let $\partial_{\widehat r}$ and $\partial_{\widehat R}$ be derivatives in the directions of the gradients of $r$ and $R$ respectively.  (Note that the gradient of $R$ is a unit vector but the gradient of $r$ need not be a unit vector.)  We denote the induced metric $\langle\,\cdot\,,\,\cdot\,\rangle$ on $\mathcal{M}$ in $(x_1,x_2,\ldots,x_n)$ coordinates by 
\begin{equation}\label{normal coord eqn3}
	\sum_{i,j=1}^n g_{ij} \,dx^i \,dx^j \quad\text{where}\quad g_{ij} = \left\langle \frac{\partial}{\partial x_i}, \frac{\partial}{\partial x_j} \right\rangle . 
\end{equation}
We shall let $[g^{ij}]$ denote the inverse matrix of $[g_{ij}]$ and $g$ denote the determinant of $[g_{ij}]$.  By \eqref{exp map est} 
\begin{equation}\label{normal coord eqn4}
	|g_{ij} - \delta_{ij}| + R \,|Dg_{ij}| \leq C\mathbf{m}_0^{1/2} R 
\end{equation} 
on $B^{\mathcal M}_1(Z_0)$, where $C = C(n,m,q,M_0,N_0,C_e,C_h) \in (0,\infty)$ is a constant.

Recall that for each $Z \in \mathcal{B}^{(\geq 2)}_q$, $\mathbf{H}_{N,Z}$, $\mathbf{D}_{N,Z}$, $\mathbf{\Sigma}_{N,Z}$, and $\mathbf{I}_{N,Z}$ are as in Definition~\ref{freq fun defn} and that $\mathbf{I}_{N,Z}(0^+) = \lim_{\rho \rightarrow 0^+} \mathbf{I}_{N,Z}(\rho)$.

\begin{lemma}\label{weiss lemma}
Let $M_0,N_0,C_e,C_h,\varepsilon_2$ be as in Hypothesis~\ref{DeLSpa consts freq hypothesis}.  Let $\varsigma = \varsigma(m,n,q) \in (0,1)$ be as in Theorem~\ref{relharm separation main thm}.  For all $\alpha \geq 2$ there exists there exist $\varepsilon_2 = \varepsilon_2(n,m,q,M_0,N_0,C_e,C_h) \in (0,1)$ such that the following holds true.  Let $T$ be as in Hypothesis~\ref{center mnfd hypothesis} and $\mathcal{B}^{(\geq 2)}_q$ be as in Definition~\ref{center mnfd branch defn}.  Let $\mathcal{M}$ and $N$ be as in Hypothesis~\ref{M and N hypothesis}.  Then 
\begin{align}\label{weiss concl}
	&\frac{d}{d\rho} \Big( \rho^{2-n-2\alpha} \,\mathbf{D}_{N,Z_0}(\rho) - \alpha \rho^{1-n-2\alpha} \,\mathbf{H}_{N,Z_0}(\rho) \Big) 
	\\ \geq\,& -2\rho^{-n-2\alpha} \int_{\mathcal M} R^{2\alpha+1} \,|\partial_{\widehat R} (N/R^{\alpha})|^2 \phi'(R/\rho) 
		- C \mathbf{m}_0^{\gamma_3} \rho^{-n-1-2\alpha+\gamma_3} \int_{B^{\mathcal M}_{16\rho/15}(Z_0)} |N|^2 \nonumber
\end{align}
for all $Z_0 \in \mathcal{B}^{(\geq 2)}_q$ and $\rho \in (0,15/16]$, where $C = C(n,m,q,M_0,N_0,C_e,C_h,\alpha) \in (0,\infty)$ is a constant.
\end{lemma}

\begin{proof} 
Fix $Z_0 \in \mathcal{B}^{(\geq 2)}_q$ and let $\mathbf{H}(\rho) = \mathbf{H}_{N,Z_0}(\rho)$,  $\mathbf{D}(\rho) = \mathbf{D}_{N,Z_0}(\rho)$ and $\mathbf{\Sigma}(\rho) = \mathbf{\Sigma}_{N,Z_0}(\rho)$ for each $\rho \in (0,1]$.  Let $\widehat{E}(\rho) = \int_{B^{\mathcal M}_{16\rho/15}(Z_0)} |N|^2$.  Clearly $\mathbf{H}(\rho) \leq 4\rho^{-1} \widehat{E}(\rho)$ and $\mathbf{\Sigma}(\rho) \leq \widehat{E}(\rho)$ for each $\rho \in (0,1]$.  By Lemma~\ref{normal approx energy lemma}, $\mathbf{D}(\rho) + \rho \,\mathbf{D}'(\rho) \leq C \rho^{-2} \widehat{E}(\rho)$ for each $\rho \in (0,1]$.  Hence \eqref{DeLSpa freq mono eqn2}, \eqref{DeLSpa freq mono eqn4}, \eqref{DeLSpa freq mono eqn3}, \eqref{DeLSpa freq mono eqn16}, and \eqref{DeLSpa freq mono eqn17} give us 
\begin{align*}
	\left| \mathbf{H}'(\rho) - \frac{n-1}{\rho} \,\mathbf{H}(\rho) + \frac{2}{\rho} \int_{\mathcal M} \sum_{l=1}^q \langle N_l, \partial_{\widehat R} N_l \rangle 
		\,\phi'(R/\rho) \right| \leq\,& C \mathbf{m}_0^{1/2} \rho^{-1} \widehat{E}(\rho) , \\
	\left| \mathbf{D}(\rho) + \frac{1}{\rho} \int_{\mathcal M} \sum_{l=1}^q \langle N_l, \partial_{\widehat R} N_l \rangle \,\phi'(R/\rho) \right| 
		\leq\,& C \mathbf{m}_0^{\gamma_3} \rho^{-2+\gamma_3} \widehat{E}(\rho) , \\
	\left| \mathbf{D}'(\rho) - \frac{n-2}{\rho} \,\mathbf{D}(\rho) + \frac{2}{\rho^2} \int_{\mathcal M} |\partial_{\widehat R} N|^2 \,R \,\phi'(R/\rho) \right| 
		\leq\,& C \mathbf{m}_0^{\gamma_3} \rho^{-3+\gamma_3} \widehat{E}(\rho) 
\end{align*}
for each $\rho \in (0,1]$, where $C = C(n,m,q,M_0,N_0,C_e,C_h) \in (0,\infty)$ are constants.  Hence 
\begin{align*}
	&\frac{d}{d\rho} \Big( \rho^{2-n-2\alpha} \mathbf{D}(\rho) - \alpha \rho^{1-n-2\alpha} \mathbf{H}(\rho) \Big) 
	\\=\,& \rho^{2-n-2\alpha} \mathbf{D}'(\rho) - \alpha \rho^{1-n-2\alpha} \mathbf{H}'(\rho) 
		- (2\alpha+n-2) \,\rho^{1-n-2\alpha} \mathbf{D}(\rho) \\&+ \alpha (2\alpha+n-1) \,\rho^{-n-2\alpha} \mathbf{H}(\rho) 
	\\ \geq\,& -2\rho^{2-n-2\alpha} \int_{\mathcal M} |\partial_{\widehat R} N|^2 \,R \,\phi'(R/\rho) 
		+ 4\alpha \rho^{-n-2\alpha} \int_{\mathcal M} \sum_{l=1}^q \langle N_l, \partial_{\widehat R} N_l \rangle \,\phi'(R/\rho) 
		\\&- 2\alpha^2 \rho^{-n-2\alpha} \int_{\mathcal M} |N|^2 \,\frac{1}{R} \,\phi'(R/\rho) 
		- C \mathbf{m}_0^{\gamma_3} \rho^{-n-1-2\alpha+\gamma_3} \widehat{E}(\rho)
	\\=\,& -2\rho^{2-n-2\alpha} \int_{\mathcal M} |R \,\partial_{\widehat R} N - \alpha \,N|^2 \,\frac{1}{R} \,\phi'(R/\rho) 
		- C \mathbf{m}_0^{\gamma_3} \rho^{-n-1-2\alpha+\gamma_3} \widehat{E}(\rho) 
\end{align*}
for some constant $C = C(n,m,q,M_0,N_0,C_e,C_h,\alpha) \in (0,\infty)$, thereby proving \eqref{weiss concl}.
\end{proof}

\begin{theorem}\label{keyest thm} 
Let $\alpha = \ell_0/q_0 \geq 2$ for some relatively prime positive integers $\ell_0$ and $q_0 \leq q$.  Let $\varsigma = \varsigma(m,n,q) \in (0,1)$ be as in Theorem~\ref{relharm separation main thm}.  For every $p \in \{1,2,\ldots, \lfloor q/q_0 \rfloor \}$ there exists $M_0,N_0,C_e,C_h,\varepsilon_2$ as in Hypothesis~\ref{DeLSpa constants hypothesis} (independent of $p$ and $\alpha$) and there exists $\varepsilon_0,\beta \in (0,1)$ depending only on $n,m,q,p,M_0,N_0,C_e,C_h,\alpha$ such that if Hypothesis~$(\star)$ and Hypothesis~$(\dagger)$ hold true and $\mathbf{I}_{N,Z_0}(0^+) \geq \alpha$, then 
\begin{equation}\label{keyest concl}
	\frac{1}{\rho_0^n} \int_{B^{\mathcal M}_{3\rho_0/4}(Z_0)} \sum_{i=3}^n \bigg|\frac{\partial N}{\partial x_i}\bigg|^2 
	+ \rho_0^{2\alpha-2} \int_{B^{\mathcal M}_{3\rho_0/4}(Z_0)} R^{2-n} \,\bigg| \partial_{\hat R} \bigg(\frac{N}{R^{\alpha}}\bigg) \bigg|^2
	\leq C E_{\varsigma/2}(N,\varphi,B^{\mathcal M}_{\rho_0}(Z_0))^2 ,  
\end{equation}
where $C = C(n,m,q,M_0,N_0,C_e,C_h,\alpha) \in (0,\infty)$ is a constant.
\end{theorem}

\begin{proof}  \textit{Step~1.  Integral estimates based on a modification of the monotonicity formula for Almgren's frequency function.}  By Corollary~\ref{doubling centered cor}, either $\mathbf{I}_{N,Z_0}(0^+) = \alpha$ and $\limsup_{\rho\rightarrow 0^+} \rho^{-\alpha} \mathbf{H}_{N,Z_0}(\rho) < \infty$ or $\mathbf{I}_{N,Z_0}(0^+) > \alpha$ and $\lim_{\rho\rightarrow 0^+} \mathbf{H}_{N,Z_0}(\rho) = 0$.  Hence 
\begin{equation*}
	\lim_{\rho\rightarrow 0^+} (\rho^{2-n-2\alpha} \,\mathbf{D}_{N,Z_0}(\rho) - \alpha \rho^{1-n-2\alpha} \,\mathbf{H}_{N,Z_0}(\rho)) 
	= \lim_{\rho\rightarrow 0^+} \rho^{1-n-2\alpha} \mathbf{H}_{N,Z_0}(\rho) \,(\mathbf{I}_{N,Z_0}(\rho) - \alpha) = 0.
\end{equation*}
Hence integrating \eqref{weiss concl}  
\begin{align}\label{keyest eqn2}
	&\rho^{2-n-2\alpha} \,\mathbf{D}_{N,Z_0}(\rho) - \alpha \rho^{1-n-2\alpha} \,\mathbf{H}_{N,Z_0}(\rho)) 
	\\ \geq\,& -2 \int_{0}^{\rho} \tau^{-n-2\alpha} \int_{\mathcal M} R^{2\alpha+1} \,|\partial_{\widehat R} (N/R^{\alpha})|^2 
		\phi'(R/\tau) \,d\mathcal{H}^n(X) \,d\tau \nonumber\\& - C \mathbf{m}_0^{\gamma_3} \int_{0}^{\rho} \tau^{-n-1-2\alpha+\gamma_3} 
		\int_{B^{\mathcal M}_{16\tau/15}(Z_0)} |N|^2 \,d\mathcal{H}^n(X) \,d\tau \nonumber 
\end{align}
for all $\rho \in (0,15\rho_0/16]$, where $C = C(n,m,q,M_0,N_0,C_e,C_h,\alpha) \in (0,\infty)$ are constants.  For each $\rho \in (0,15\rho_0/16]$ and each $X \in B^{\mathcal M}_{24\rho/29}(Z_0)$, using the definition of $\phi$ and $R < 29R/24 \leq \min\{\rho,2R\}$ 
\begin{equation*}
	-2 \int_{0}^{\rho} \tau^{-n-2\alpha} \,\phi'(R/\tau) \,d\tau 
	= 4 \int_{R}^{\min\{\rho,2R\}} \tau^{-n-2\alpha} \,d\tau 
	\geq c R^{1-n-2\alpha} , 
\end{equation*}
where $c = c(n,\alpha) > 0$ is a constant.  Thus by Fubini's theorem 
\begin{align}\label{keyest eqn3}
	&-2 \int_{0}^{\rho} \tau^{-n-2\alpha} \int_{\mathcal M} R^{2\alpha+1} \,|\partial_{\widehat R} (N/R^{\alpha})|^2 \phi'(R/\tau) \,d\mathcal{H}^n(X) \,d\tau
	\\&\geq c \int_{B^{\mathcal M}_{24\rho/29}(Z_0)} R^{2-n} \,|\partial_{\widehat R} (N/R^{\alpha})|^2 \,d\mathcal{H}^n(X) , \nonumber 
\end{align}
where $c = c(n,\alpha) > 0$ is a constant.  By Corollary~\ref{doubling centered cor} 
\begin{equation}\label{keyest eqn4}
	\int_{B^{\mathcal M}_{16\tau/15}(Z_0)} |N|^2 \leq C \Big(\frac{\tau}{\rho_0}\Big)^{n+2\alpha} \int_{B^{\mathcal M}_{\rho_0}(Z_0)} |N|^2
\end{equation}
for all $\tau \in (0,15\rho_0/16]$, where $C = C(n,m,q,M_0,N_0,C_e,C_h,\alpha) \in (0,\infty)$ is a constant.  Multiplying both sides of \eqref{keyest eqn4} by $\tau^{-n-1-2\alpha+\gamma_3}$ and integrating over $\tau \in (0,\rho]$ 
\begin{align}\label{keyest eqn5}
	&\int_{0}^{\rho} \tau^{-n-1-2\alpha+\gamma_3} \int_{B^{\mathcal M}_{16\tau/15}(Z_0)} |N|^2 \,d\mathcal{H}^n(X) \,d\tau 
	\\ \leq\,& C \rho_0^{-n-2\alpha} \int_{0}^{\rho} \tau^{-1+\gamma_3} \,d\tau \int_{B^{\mathcal M}_{\rho_0}(Z_0)} |N|^2 \,d\mathcal{H}^n(X) \nonumber
	\\=\,& \frac{C}{\gamma_3} \rho_0^{-n-2\alpha+\gamma_3} \int_{B^{\mathcal M}_{\rho_0}(Z_0)} |N|^2 \,d\mathcal{H}^n(X) \nonumber
\end{align}
for all $\rho \in (0,\rho_0]$, where $C = C(n,m,q,M_0,N_0,C_e,C_h,\alpha) \in (0,\infty)$ is a constant.  Bounding the right-hand side of \eqref{keyest eqn2} from below using \eqref{keyest eqn3} and \eqref{keyest eqn5}
\begin{align}\label{keyest eqn6}
	&\rho^{2-n-2\alpha} \,\mathbf{D}_{N,Z_0}(\rho) - \alpha \rho^{1-n-2\alpha} \,\mathbf{H}_{N,Z_0}(\rho)) 
	\\ \geq\,& c \int_{B^{\mathcal M}_{24\rho/29}(Z_0)} R^{2-n} \,|\partial_{\widehat R} (N/R^{\alpha})|^2 
		- C \mathbf{m}_0^{\gamma_3} \rho^{-n-2\alpha+\gamma_3} \int_{B^{\mathcal M}_{\rho_0}(Z_0)} |N|^2 \nonumber 
\end{align}
for all $\rho \in (0,15\rho_0/16]$, where $c = c(n,\alpha) > 0$ and $C = C(n,m,q,M_0,N_0,C_e,C_h,\alpha) \in (0,\infty)$ are constants.  On the other hand, by the definitions of $\mathbf{D}_{N,Z_0}(\rho)$ and $\mathbf{H}_{N,Z_0}(\rho)$ 
\begin{equation}\label{keyest eqn7}
	\mathbf{D}_{N,Z_0}(\rho) - \alpha \rho^{-1} \,\mathbf{H}_{N,Z_0}(\rho)) 
	= \int_{\mathcal M} \left( |DN|^2 \,\phi(R/\rho) + \alpha \,|N|^2 \,\frac{1}{\rho R} \,\phi'(R/\rho) \right)  
\end{equation}
for all $\rho \in (0,15\rho_0/16]$.  Thus by combining \eqref{keyest eqn6} and \eqref{keyest eqn7}
\begin{align}\label{keyest eqn8}
	\int_{B^{\mathcal M}_{24\rho/29}(Z_0)} R^{2-n} \,|\partial_{\widehat R} (N/R^{\alpha})|^2 
	\leq\,& C\rho^{2-n-2\alpha} \int_{\mathcal M} \left( |DN|^2 \,\phi(R/\rho) + \alpha \,|N|^2 \,\frac{1}{\rho R} \,\phi'(R/\rho) \right) 
		\\&+ C \mathbf{m}_0^{\gamma_3} \rho^{-n-2\alpha+\gamma_3} \int_{B^{\mathcal M}_{\rho_0}(Z_0)} |N|^2 \nonumber 
\end{align}
for all $\rho \in (0,15\rho_0/16]$, where $C = C(n,m,q,M_0,N_0,C_e,C_h,\alpha) \in (0,\infty)$ is a constant.  Define $\psi : [0,\infty) \rightarrow \mathbb{R}$ by 
\begin{equation*}
	\psi(s) = \frac{16}{\rho_0} \int_{7\rho_0/8}^{15\rho_0/16} \phi(s/\rho) \,d\rho 
\end{equation*}
for all $s \in [0,\infty)$.  Since $\phi$ is the Lipschitz function given by \eqref{freq phi defn}, $\psi \in C^{1,1}([0,\infty))$ with $\psi(s) = 1$ if $0 \leq s \leq 7\rho_0/16$, $\psi(s) = 0$ if $s \geq 15\rho_0/16$, and $\psi'(s) = \frac{16}{\rho_0} \int_{7\rho_0/8}^{15\rho_0/16} \frac{1}{\rho} \,\phi'(s/\rho) \,d\rho$ for all $s \in [0,\infty)$.  Note that for each $\rho \in [29\rho_0/32,15\rho_0/16]$ we have $24\rho/29 \geq 3\rho_0/4$.  Multiplying both sides of \eqref{keyest eqn8} by $\rho^{n-2+2\alpha}$ and integrating over $\rho \in [7\rho_0/8,15\rho_0/16]$, 
\begin{align}\label{keyest eqn9}
	\rho_0^{n-2+2\alpha} \int_{B^{\mathcal M}_{3\rho_0/4}(Z_0)} R^{2-n} \,|\partial_{\widehat R} (N/R^{\alpha})|^2 
	\leq\,& C \int_{\mathcal M} \left( |DN|^2 \,\psi(R) + \alpha \,|N|^2 \,\frac{1}{R} \,\psi'(R) \right) 
	\\&+ C \mathbf{m}_0^{\gamma_3} \rho_0^{\gamma_3-2} \int_{B^{\mathcal M}_{\rho_0}(Z_0)} |N|^2 , \nonumber 
\end{align}
where $C = C(n,m,q,M_0,N_0,C_e,C_h,\alpha) \in (0,\infty)$ is a constant.

\noindent\textit{Step~2.  Integral estimates based on the outer and inner variation formula.}  Let $\zeta_{\rm out} \in C^1_c(\mathbf{U})$ be a vector field such that $\zeta_{\rm out}(X) = \psi(R(\mathbf{p}(X))) \,(X - \mathbf{p}(X))$ on an open neighborhood of 
$\op{spt} T \,\cup \,\op{spt} \mathbf{G}_N$.  By~\cite[Theorem~4.2]{DeLSpa0} (cf.\ \cite[Theorem~5.6]{Almgren})
\begin{align}\label{keyest eqn10}
	&\int_{\mathcal M} \left( |DN|^2 \,\psi(R) + \sum_{l=1}^q \langle N_l, \partial_{\widehat R} N_l \rangle \,\psi'(R) \right) 
	\leq \delta\mathbf{G}_N(\zeta_{\rm out}) + q \int_{\mathcal M} \langle H_{\mathcal M}, N_a \rangle \,\psi(R) 
		\\&+ C \int_{\mathcal M} (|A_{\mathcal M}| \,|N| + |DN|^2)^2 \,\psi(R) 
		+ C \int_{\mathcal M} (|A_{\mathcal M}| \,|N| + |DN|^2) \,|N| \,|DN| \,|\psi'(R)| , \nonumber 
\end{align}
where $C = C(n,m,q,M_0,N_0,C_e,C_h) \in (0,\infty)$ is a constant.  We claim that we can bound the terms on the right-hand side of \eqref{keyest eqn10} to obtain 
\begin{equation}\label{keyest eqn11}
	\int_{\mathcal M} \left( |DN|^2 \,\psi(R) + \sum_{l=1}^q \langle N_l, \partial_{\widehat R} N_l \rangle \,\psi'(R) \right) 
	\leq C \mathbf{m}_0^{\gamma_2/4} \rho_0^{\gamma_2/4-2} \int_{B^{\mathcal M}_{\rho_0}(Z_0)} |N|^2
\end{equation}
where $C = C(n,m,q,M_0,N_0,C_e,C_h) \in (0,\infty)$ is a constant.  To bound the the first term on the right-hand side of \eqref{keyest eqn10}, notice that by~\cite[Corollary~2.2]{DeLSpa2} and \eqref{normal approx est1} that $|D\zeta_{\rm out}| \leq C$ on $\op{spt} T \cup \op{spt} \mathbf{G}_N$.  Thus by \eqref{normal approx est3} 
\begin{align}\label{keyest eqn12}
	|\delta\mathbf{G}_N(\zeta_{\rm out})| 
	=\,& |\delta\mathbf{G}_N(\zeta_{\rm out}) - \delta T(\zeta_{\rm out})| 
	\leq C \|\mathbf{G}_N - T\|(\mathcal{U}^i \setminus \mathcal{K}) 
	\\ \leq\,& C \mathbf{m}_0^{3\gamma_2/4} \rho_0^{3\gamma_2/4-2} \int_{B^{\mathcal M}_{\rho_0}(Z_0)} |N|^2 \nonumber 
\end{align}
where $C = C(n,m,q,M_0,N_0,C_e,C_h) \in (0,\infty)$ are constants.  We can bound second term on the right-hand side of \eqref{keyest eqn10} using $\|H_{\mathcal M}\|_{C^0} \leq C\mathbf{m}_0^{1/2}$ and \eqref{normal approx est4}.  We can bound third term on the right-hand side of \eqref{keyest eqn10} by noting that $(|A_{\mathcal M}| \,|N| + |DN|^2)^2 \leq 2 |A_{\mathcal M}|^2 \,|N|^2 + 2 |DN|^4$ and then using $\|A_{\mathcal M}\|_{C^0} \leq C\mathbf{m}_0^{1/2}$, \eqref{normal approx est1}, and Lemma~\ref{normal approx energy lemma}.  To bound the the last term on the right-hand side of \eqref{keyest eqn10}, we use $\|A_{\mathcal M}\|_{C^0} \leq C\mathbf{m}_0^{1/2}$, $|\psi'(R)| \leq C/\rho$, \eqref{normal approx est1}, and Lemma~\ref{normal approx energy lemma} to obtain 
\begin{align*}
	&\int_{\mathcal M} (|A_{\mathcal M}| \,|N| + |DN|^2) \,|N| \,|DN| \,|\psi'(R)|
	\\ \leq\,& C \mathbf{m}_0^{1/2+\gamma_2} \rho_0^{\gamma_2-1} \int_{B^{\mathcal M}_{15\rho_0/16}(Z_0)} |N|^2 
		+ C \mathbf{m}_0^{1/(2n)+\gamma_2} \rho_0^{\beta_2+\gamma_2} \int_{B^{\mathcal M}_{15\rho_0/16}(Z_0)} |DN|^2 
	\\ \leq\,& C \mathbf{m}_0^{1/(2n)+\gamma_2} \rho_0^{\gamma_2-2} \int_{B^{\mathcal M}_{\rho_0}(Z_0)} |N|^2 
\end{align*}
where $C = C(n,m,q,M_0,N_0,C_e,C_h) \in (0,\infty)$ are constants.  

Next recall that $(x_1,x_2,\ldots,x_n)$ denotes normal coordinates on $\mathcal{M}$ on $Z_0$ as in \eqref{normal coord eqn1} and that $r$ is as in \eqref{normal coord eqn2}.  Let $\zeta_{\rm in} \in C^1_c(\mathcal{M})$ be a vector field such that $\zeta_{\rm in}(X) = \psi(R) \,r \,\widehat{r}$ on $\mathcal{M}$ and extend $\zeta_{\rm in}$ to a vector field in $C^1_c(\mathbf{U})$ such that $\zeta_{\rm in}(X) = \zeta_{\rm in}(\mathbf{p}(X))$ on an open neighborhood of $\op{spt} T \cup \op{spt} \mathbf{G}_N$.  By~\cite[Theorem~4.3]{DeLSpa0} with $Y = \zeta_{\rm in}$ (see also~\cite[Theorem~5.7]{Almgren})
\begin{align}\label{keyest eqn13}
	&\int_{\mathcal M} \left( \tfrac{1}{2} |DN|^2 \,(1 + r\,\Delta_{\mathcal M} r) - |\partial_{\widehat r} N|^2 
		- \sum_{l=1}^q r \,\langle DN_l, DN_l \cdot \op{Hess} r \rangle \right) \psi(R) 
	\\&+ \int_{\mathcal M} \left( \tfrac{1}{2} |DN|^2 \,r \,\langle \widehat{r}, \widehat{R} \rangle 
		- \sum_{l=1}^q r \,\langle \partial_{\widehat r} N_l, \partial_{\widehat R} N_l \rangle \right) \psi'(R) \nonumber \\ 
	\leq\,& \delta\mathbf{G}_N(\zeta_{\rm in}) + q \int_{\mathcal M} (|H_{\mathcal M}| \,|N_a| \,|D\zeta_{\rm in}| 
		+ |DH_{\mathcal M}| \,|N_a| \,|\zeta_{\rm in}|) \nonumber\\&+ C \int_{\mathcal M} (|A_{\mathcal M}|^2 \,|N| \,|DN| + |A_{\mathcal M}| \,|N| \,|DN|^2 
		+ |A_{\mathcal M}| \,|DN|^3) \,|\zeta_{\rm in}| \nonumber\\&+ C \int_{\mathcal M} (|A_{\mathcal M}|^2 \,|N|^2 + |A_{\mathcal M}| |N|^2 |DN| + |DN|^4) 
		\,|D\zeta_{\rm in}| , \nonumber 
\end{align}
where $C = C(n,m,q,M_0,N_0,C_e,C_h) \in (0,\infty)$ is a constant.  We claim that we can bound the terms on the right-hand side of \eqref{keyest eqn13} to obtain 
\begin{align}\label{keyest eqn14}
	&\int_{\mathcal M} \left( \tfrac{1}{2} |DN|^2 \,(1 + r\,\Delta_{\mathcal M} r) - |\partial_{\widehat r} N|^2 
		- \sum_{l=1}^q r \,\langle DN_l, DN_l \cdot \op{Hess} r \rangle \right) \psi(R) 
	\\&+ \int_{\mathcal M} \left( \tfrac{1}{2} |DN|^2 \,r \,\langle \widehat{r}, \widehat{R} \rangle 
		- \sum_{l=1}^q r \,\langle \partial_{\widehat r} N_l, \partial_{\widehat R} N_l \rangle \right) \psi'(R) \nonumber \\ 
	\leq\,& C \mathbf{m}_0^{\gamma_2/4} \rho_0^{\gamma_2/4-2} \int_{B^{\mathcal M}_{\rho_0}(Z_0)} |N|^2 , \nonumber 
\end{align}
where $C = C(n,m,q,M_0,N_0,C_e,C_h) \in (0,\infty)$ is a constant.  To bound the first term on the right-hand side of \eqref{keyest eqn13}, we notice that $|D\zeta_{\rm in}| \leq C$ on $\op{spt} T \cup \op{spt} \mathbf{G}_N$ and then argue as in \eqref{keyest eqn12}.  We can bound second term on the right-hand side of \eqref{keyest eqn12} using $\|H_{\mathcal M}\|_{C^1} \leq C\mathbf{m}_0^{1/2}$, $|D\zeta_{\rm in}| \leq C$ on $\op{spt} \mathbf{G}_N$, and \eqref{normal approx est4}.  To bound the third term on the right-hand side of \eqref{keyest eqn13}, we use $\|A_{\mathcal M}\|_{C^0} \leq C\mathbf{m}_0^{1/2}$, $|\zeta_{\rm in}| \leq r \leq \rho_0$ on $\op{spt} \mathbf{G}_N$, \eqref{normal approx est1}, and Lemma~\ref{normal approx energy lemma} to obtain 
\begin{align}\label{keyest eqn15}
	&\int_{\mathcal M} (|A_{\mathcal M}|^2 \,|N| \,|DN| + |A_{\mathcal M}| \,|N| \,|DN|^2 + |A_{\mathcal M}| \,|DN|^3) \,|\zeta_{\rm in}| 
	\\ \leq\,& C \mathbf{m}_0 \rho_0 \int_{B^{\mathcal M}_{15\rho_0/16}(Z_0)} |N| \,|DN| 
		+ C \mathbf{m}_0^{1/2} \rho_0^{1+\gamma_2} \int_{B^{\mathcal M}_{15\rho_0/16}(Z_0)} |DN|^2 \nonumber 
	\\ \leq\,& C \mathbf{m}_0 \int_{B^{\mathcal M}_{15\rho_0/16}(Z_0)} |N|^2 
		+ C \mathbf{m}_0^{1/2} \rho_0^{1+\gamma_2} \int_{B^{\mathcal M}_{15\rho_0/16}(Z_0)} |DN|^2 \nonumber
	\\ \leq\,& C \mathbf{m}_0^{1/2} \rho_0^{\gamma_2-1} \int_{B^{\mathcal M}_{\rho_0}(Z_0)} |N|^2 \nonumber 
\end{align}
where $C = C(n,m,q,M_0,N_0,C_e,C_h) \in (0,\infty)$ are constants.  To bound the last term on the right-hand side of \eqref{keyest eqn13}, we use $\|A_{\mathcal M}\|_{C^0} \leq C\mathbf{m}_0^{1/2}$, $|D\zeta_{\rm in}| \leq C$, \eqref{normal approx est1}, and Lemma~\ref{normal approx energy lemma}.  

Having shown \eqref{keyest eqn14}, observe that by \eqref{normal coord eqn4} and Gauss's lemma
\begin{gather*}
	\widehat{r} = \sum_{i=1}^2 \sum_{j=1}^n g^{ij} \,\frac{x_i}{r} \,\frac{\partial}{\partial x_j} , \quad 
	\widehat{R} = \sum_{i,j=1}^n g^{ij} \,\frac{x_i}{R} \,\frac{\partial}{\partial x_j} = \sum_{i=1}^n \frac{x_i}{R} \,\frac{\partial}{\partial x_i} , \\ 
	\left| \op{Hess} r - \sum_{i=1}^2 \left( \frac{\delta_{ij}}{r} - \frac{x_i x_j}{r^3} \right) dx_i \otimes dx_j \right| \leq C\mathbf{m}_0^{1/2} , 
\end{gather*}
where $C = C(n,m,q,M_0,N_0,C_e,C_h) \in (0,\infty)$ is a constant.  Hence using \eqref{normal coord eqn4} and Lemma~\ref{normal approx energy lemma}, \eqref{keyest eqn14} gives us 
\begin{align}\label{keyest eqn16}
	&\int_{\mathcal M} \sum_{i=3}^n \bigg|\frac{\partial N}{\partial x_i}\bigg|^2 \,\psi(R) 
		+ \int_{\mathcal M} \left( \frac{r^2}{2R} \,|DN|^2 - \sum_{l=1}^q r \,\langle \partial_{\widehat r} N_l, \partial_{\widehat R} N_l \rangle \right) \psi'(R)  
	\\ \leq\,& C\mathbf{m}_0^{1/2} \int_{B^{\mathcal M}_{15\rho_0/16}(Z_0)} R \,|DN|^2
		+ C \mathbf{m}_0^{\gamma_2/4} \rho_0^{\gamma_2/4-2} \int_{B^{\mathcal M}_{\rho_0}(Z_0)} |N|^2 , \nonumber 
	\\ \leq\,& C \mathbf{m}_0^{\gamma_2/4} \rho_0^{\gamma_2/4-2} \int_{B^{\mathcal M}_{\rho_0}(Z_0)} |N|^2 , \nonumber 
\end{align}
where $C = C(n,m,q,M_0,N_0,C_e,C_h) \in (0,\infty)$ are constants.  

Now computing $\alpha\,\cdot\,$\eqref{keyest eqn11}$\,+\,$\eqref{keyest eqn16} and then adding $\alpha^2 \int_{\mathcal M} |N|^2 \,R^{-1} \psi'(R)$ to both sides 
\begin{align}\label{keyest eqn17}
	&\int_{\mathcal M} \left( \alpha |DN|^2 + \sum_{i=3}^n \bigg|\frac{\partial N}{\partial x_i}\bigg|^2 \right) \psi(R) 
		+ \alpha^2 \int_{\mathcal M} |N|^2 \,\frac{1}{R} \,\psi'(R)
	\\ \leq\,& -\int_{\mathcal M} \left( \frac{r^2}{2R} \,|DN|^2 - \sum_{l=1}^q \langle r \partial_{\widehat r} N_l - \alpha N_l , \partial_{\widehat R} N_l \rangle 
		- \alpha^2 |N|^2 \,\frac{1}{R} \right) \psi'(R) \nonumber 
	\\&+ C \mathbf{m}_0^{\gamma_2/4} \rho_0^{\gamma_2/4-2} \int_{B^{\mathcal M}_{\rho_0}(Z_0)} |N|^2 \nonumber 
\end{align}
where $C = C(n,m,q,M_0,N_0,C_e,C_h,\alpha) \in (0,\infty)$ is a constant.  Observe that 
that 
\begin{equation*}
	\partial_{\widehat R} N_l = \sum_{i,j=1}^n g^{ij} \frac{x_i}{R} \,\frac{\partial N_l}{\partial x_j} 
		= \frac{r}{R} \,\partial_{\widehat r} N_l + \sum_{i=3}^n \sum_{j=1}^n  g^{ij} \frac{x_i}{R} \,\frac{\partial N_l}{\partial x_j}
\end{equation*}
at $\mathcal{H}^n$-a.e.~point of $B^{\mathcal M}_{\rho_0}(Z_0)$ and for $l = 1,2,\ldots,q$ (where $g_{ij}$ is as in \eqref{normal coord eqn3}).  Hence using \eqref{normal coord eqn4} and Cauchy's inequality 
\begin{align}\label{keyest eqn18}
    &-\langle r \partial_{\widehat r} N_l - \alpha N_l , \partial_{\widehat R} N_l \rangle \,\psi'(R) \leq -\langle r \partial_{\widehat r} N_l - \alpha N_l , r \partial_{\widehat r} N_l \rangle\,\frac{1}{R}\,\psi'(R)
    \\&\hspace{15mm}+ \frac{1}{2} \sum_{i=3}^n \left|\frac{\partial N}{\partial x_i}\right|^2 \psi(R) + C \mathbf{m}_0 R^2 |DN|^2 \psi(R) + C |r \partial_{\widehat r} N_l - \alpha N_l|^2 \,\frac{\psi'(R)^2}{\psi(R)} \nonumber  
\end{align}
at $\mathcal{H}^n$-a.e.~point of $B^{\mathcal M}_{\rho_0}(Z_0)$ and for $l = 1,2,\ldots,q$, where $C = C(n,m,q,M_0,N_0,C_e,C_h) \in (0,\infty)$ is a constant.  We have that $|\psi'(s)|^2/\psi(s) \leq 2 \|\psi''\|_{L^{\infty}([0,\infty))} \leq C$ for all $s \in [0,15\rho_0/16)$ for some constant $C \in (0,\infty)$ (to see this, one can for instance smooth $\psi$ via mollification and argue as in~\cite[Lemma~6.6]{IlmanenReg}).  Now by \eqref{keyest eqn17}, \eqref{keyest eqn18}, and Lemma~\ref{normal approx energy lemma}
\begin{align}\label{keyest eqn19}
	&\frac{1}{2} \int_{\mathcal M} \sum_{i=3}^n \bigg|\frac{\partial N}{\partial x_i}\bigg|^2 \,\psi(R) + \alpha \int_{\mathcal M} |DN|^2 \,\psi(R) 
		+ \alpha^2 \int_{\mathcal M} |N|^2 \,\frac{1}{R} \,\psi'(R)
	\\ \leq\,& -\int_{\mathcal M} \left( \tfrac{1}{2} r^2 |DN|^2 - \sum_{l=1}^q \langle r \partial_{\widehat r} N_l - \alpha N_l , r \partial_{\widehat r} N_l \rangle 
		- \alpha^2 |N|^2 \right) \frac{1}{R} \,\psi'(R) \nonumber 
	\\&+ C \int_{B^{\mathcal M}_{15\rho_0/16}(Z_0)} \sum_{l=1}^q |r \partial_{\widehat r} N_l - \alpha N_l|^2 
		+ C \mathbf{m}_0^{\gamma_2/4} \rho_0^{\gamma_2/4-2} \int_{B^{\mathcal M}_{\rho_0}(Z_0)} |N|^2 , \nonumber 
\end{align}
where $C = C(n,m,q,M_0,N_0,C_e,C_h) \in (0,\infty)$ is a constant.  Combining \eqref{keyest eqn9} and \eqref{keyest eqn19} 
\begin{align}\label{keyest eqn20}
	&\int_{\mathcal M} \sum_{i=3}^n \bigg|\frac{\partial N}{\partial x_i}\bigg|^2 \,\psi(R) 
		+ \frac{1}{C} \,\rho_0^{n-2+2\alpha} \int_{B^{\mathcal M}_{3\rho_0/4}(Z_0)} R^{2-n} \,|\partial_{\widehat R} (N/R^{\alpha})|^2 
	\\ \leq\,& -\int_{\mathcal M} \left( \tfrac{1}{2} r^2 |DN|^2 - \sum_{l=1}^q \langle r \partial_{\widehat r} N_l - \alpha N_l , r \partial_{\widehat r} N_l \rangle 
		- \alpha^2 |N|^2 \right) \frac{1}{R} \,\psi'(R) \nonumber 
	\\&+ C \int_{B^{\mathcal M}_{15\rho_0/16}(Z_0)} \sum_{l=1}^q |r \partial_{\widehat r} N_l - \alpha N_l|^2 
		+ C \mathbf{m}_0^{\gamma_2/4} \rho_0^{\gamma_2/4-2} \int_{B^{\mathcal M}_{\rho_0}(Z_0)} |N|^2 \nonumber 
\end{align}
where $C = C(n,m,q,M_0,N_0,C_e,C_h,\alpha) \in (1,\infty)$ is a constant.

\noindent\textit{Step~3.  Definition of the set $U$ and local estimates on $U$.}  Let $U$ be the union of all sets $B^{\mathcal M}_{\rho_0}(Z_0) \cap A^{\mathcal M}_{\rho,1/2}(\zeta)$, where $\rho > 0$ and $\zeta \in S$ such that $\rho^2 + |\zeta|^2 < (15\rho_0/16)^2$ and either $\rho \leq \rho_0/160$ or 
\begin{equation}\label{keyest eqn21}
	\sup_{X \in B^{\mathcal M}_{2\rho}(\exp_{Z_0}(\zeta))} \overline{R}(X) < 4^{-q-8} \rho . 
\end{equation}
Recall that $(x_1,x_2,\ldots,x_n)$ are normal coordinates on $\mathcal{M}$ on $Z_0$ as in \eqref{normal coord eqn1} and that $r$ is as in \eqref{normal coord eqn2}.  Define $d : B^{\mathcal M}_{\rho_0}(Z_0) \rightarrow \mathbb{R}$ by letting $d(X)$ be the Euclidean distance between $\exp_{Z_0}^{-1}(X)$ and $\exp_{Z_0}^{-1}(B^{\mathcal M}_{\rho_0}(Z_0) \setminus \overline{U})$ for each $X \in B^{\mathcal M}_{\rho_0}(Z_0)$ (i.e.~$d(X)$ is the Euclidean distance between the normal coordinates of $X$ and the normal coordinates of points of $B^{\mathcal M}_{\rho_0}(Z_0) \setminus \overline{U}$).  Define the locally Lipschitz cutoff function $\chi : B^{\mathcal M}_{\rho_0}(Z_0) \setminus \exp_{Z_0}(S) \rightarrow \mathbb{R}$ by 
\begin{equation*}
	\chi(X) = \begin{cases} 
		1 &\text{if } d(X) \geq r(X)/2 \\
		4 d(X)/r(X) - 1 &\text{if } r(X)/4 < d(X) < r(X)/2 \\
		0 &\text{if } d(X) \leq r(X)/4 .
	\end{cases} 
\end{equation*}
Observe that $0 \leq \chi \leq 1$ on $B^{\mathcal M}_{\rho_0}(Z_0) \setminus \exp_{Z_0}(S)$.  Moreover, if $x,\widetilde{x} \in P_{Z_0} \cap S^{\perp}$ and $y \in S$ with $|x| = |\widetilde{x}|$ and $|x|^2 + |y|^2 < (15\rho_0/16)^2$, then $\chi(\exp_{Z_0}(x+y)) = \chi(\exp_{Z_0}(\widetilde{x}+y))$.

Using the the normal coordinates $(x_1,x_2,\ldots,x_n)$ from \eqref{normal coord eqn1}, we want to cover $U$ by a countable collection of annuli $\{A^{\mathcal M}_{\rho,1/4}(\zeta)\}_{(\rho,\zeta)\in\mathcal{I}}$ satisfying certain properties to be stated below.  Notice that $B_{\rho_0}(0)$ is the set of all normal coordinates of points in $B^{\mathcal M}_{\rho_0}(Z_0)$ and for each $\rho > 0$, $\zeta \in S$, and $\kappa \in (0,1]$, $\{ (r\omega,y) : (r,y) \in B^{n-1}_{\kappa\rho/64}(\rho,\widetilde{\zeta})$ is the set of all normal coordinates of points in $A^{\mathcal M}_{\rho,\kappa}(\zeta)$, where $\widetilde{\zeta} = (\widetilde{\zeta}_1,\widetilde{\zeta}_2,\ldots,\widetilde{\zeta}_{n-2}) \in \mathbb{R}^{n-2}$ are the coordinates of $\zeta$ with respect to the basis $\vec\tau_3,\vec\tau_4,\ldots,\vec\tau_n$ of $S$ so that $\zeta = \sum_{i=1}^{n-2} \widetilde{\zeta}_i \vec\tau_{2+i}$.  By applying the Vitali covering lemma, there exists a countable collection $\mathcal{I}$ of $(\rho,\zeta)$ with $\rho > 0$, $\zeta \in S$, and $\rho^2 + |\zeta|^2 < (15\rho_0/16)^2$ such that either $\rho \leq \rho_0/160$ or \eqref{keyest eqn21} holds true for all $(\rho,\zeta) \in \mathcal{I}$, $\{A^{\mathcal M}_{\rho,1/20}(\zeta)\}_{(\rho,\zeta)\in\mathcal{I}}$ is a collection of pairwise disjoint annuli and $\{A^{\mathcal M}_{\rho,1/4}(\zeta)\}_{(\rho,\zeta)\in\mathcal{I}}$ covers $U$.  Let $\zeta,\zeta' \in S$ and let $\widetilde{\zeta},\widetilde{\zeta}'$ be the respective coordinates of $\zeta,\zeta'$ with respect to the basis $\vec\tau_3,\vec\tau_4,\ldots,\vec\tau_n$ of $S$.  Observe that if $A^{\mathcal M}_{\rho,1}(\zeta) \cap A^{\mathcal M}_{\rho',1}(\zeta') \neq \emptyset$ then $|(\rho,\widetilde{\zeta}) - (\rho',\widetilde{\zeta}')| < (\rho+\rho')/64$ and thus $\frac{63}{65} \,\rho < \rho' < \frac{65}{63} \,\rho$.  Hence 
\begin{equation*}
	B^{n-1}_{\frac{63\rho}{20 \cdot 64 \cdot 65}}(\rho',\widetilde{\zeta}') 
	\subset B^{n-1}_{\frac{\rho'}{20 \cdot 64}}(\rho',\widetilde{\zeta}') 
	\subset B^{n-1}_{\frac{\rho}{20 \cdot 64} \big( 20 + \frac{21 \cdot 65}{63}\big)}(\rho,\widetilde{\zeta}) .
\end{equation*}
Since $\{B^{n-1}_{\rho'/(20 \cdot 64)}(\rho',\widetilde{\zeta}')\}_{(\rho',\zeta')\in\mathcal{I}}$ is a pairwise disjoint collection of balls implies that for each $(\rho,\zeta) \in \mathcal{I}$ 
\begin{equation*}
	\#\{ (\rho',\zeta') \in \mathcal{I} : A^{\mathcal M}_{\rho,1}(\zeta) \cap A^{\mathcal M}_{\rho',1}(\zeta') \neq \emptyset \} \leq C(n) .
\end{equation*}
Thus there exists an integer $N \leq C(n)$ and pairwise disjoint sets $\mathcal{I}_1, \mathcal{I}_2, \ldots, \mathcal{I}_N \subset \mathcal{I}$ such that $\mathcal{I} = \bigcup_{j=1}^N \mathcal{I}_j$ and $\{A^{\mathcal M}_{\rho,1}(\zeta)\}_{(\rho,\zeta) \in \mathcal{I}_j}$ is a collection of pairwise disjoint annuli for each $j \in \{1,2,\ldots,N\}$.  Let $\{\eta_{(\rho,\zeta)}\}_{(\rho,\zeta) \in \mathcal{I}}$ is a smooth partition of unity subordinate to $\{A^{\mathcal M}_{\rho,1/2}(\zeta)\}_{(\rho,\zeta) \in \mathcal{I}}$ such that if $x,\widetilde{x} \in P_{Z_0} \cap S^{\perp}$ and $y \in S$ with $|x| = |\widetilde{x}|$ and $x+y \in A_{\rho,1/2}(\zeta)$, then $\eta_{(\rho,\zeta)}(\exp_{Z_0}(x+y)) = \eta_{(\rho,\zeta)}(\exp_{Z_0}(\widetilde{x}+y))$ and 
\begin{equation*}
	\op{spt} \eta_{(\rho,\zeta)} \subset A^{\mathcal M}_{\rho,1/2}(\zeta), \quad\quad 
	\sum_{(\rho,\zeta) \in \mathcal{I}} \eta_{(\rho,\zeta)} = 1 \text{ on } U . 
\end{equation*}
We claim that for each $(\rho,\zeta) \in \mathcal{I}$ 
\begin{align}\label{keyest eqn23} 
	&-\int_{\mathcal M} \left( \tfrac{1}{2} r^2 |DN|^2 - \sum_{l=1}^q \langle r \partial_{\widehat r} N_l - \alpha N_l , r \partial_{\widehat r} N_l \rangle 
		- \alpha^2 |N|^2 \right) \frac{1}{R} \,\psi'(R) \,\chi \,\eta_{(\rho,\zeta)}
	\\&\leq \frac{C}{\rho_0^2} \int_{A^{\mathcal M}_{\rho,1}(\zeta)} \mathcal{G}(N,\varphi\circ\exp_{Z_0}^{-1})^2 
		+ C \rho_0^{\varsigma/2-2} \int_{A^{\mathcal M}_{\rho,1}(\zeta)} |N|^2 \nonumber
\end{align}
and 
\begin{equation}\label{keyest eqn24} 
	\int_{\mathcal M} \sum_{l=1}^q |r \partial_{\widehat r} N_l - \alpha N_l|^2 \,\chi \,\eta_{(\rho,\zeta)}
	\leq \frac{C}{\rho_0^2} \int_{A^{\mathcal M}_{\rho,1}(\zeta)} \mathcal{G}(N,\varphi\circ\exp_{Z_0}^{-1})^2 
		+ C \rho_0^{\varsigma/2-2} \int_{A^{\mathcal M}_{\rho,1}(\zeta)} |N|^2 ,  
\end{equation}
where $C = C(n,m,q,M_0,N_0,C_e,C_h,\alpha) \in (0,\infty)$ are constants.  Then by summing \eqref{keyest eqn23} and \eqref{keyest eqn24} over $(\rho,\zeta) \in \mathcal{I} = \bigcup_{j=1}^N \mathcal{I}_j$ and keeping in mind that $\op{spt} \psi \subset B^{\mathcal M}_{15\rho_0/16}(Z_0)$
\begin{align}\label{keyest eqn25} 
	&-\int_{\mathcal M} \left( \tfrac{1}{2} r^2 |DN|^2 - \sum_{l=1}^q \langle r \partial_{\widehat r} N_l - \alpha N_l , r \partial_{\widehat r} N_l \rangle 
		- \alpha^2 |N|^2 \right) \frac{1}{R} \,\psi'(R) \,\chi 
	\\&\leq C \rho_0^n \,E_{\varsigma/2}(N,\varphi,B^{\mathcal M}_{\rho_0}(Z_0)) , \nonumber
\end{align}
and 
\begin{equation}\label{keyest eqn26} 
	\int_{B^{\mathcal M}_{15\rho_0/16}(Z_0)} \sum_{l=1}^q |r \partial_{\widehat r} N_l - \alpha N_l|^2 \,\chi 
	\leq C \rho_0^n \,E_{\varsigma/2}(N,\varphi,B^{\mathcal M}_{\rho_0}(Z_0)) ,  
\end{equation}
where $C = C(n,m,q,M_0,N_0,C_e,C_h,\alpha) \in (0,\infty)$ are constants.

To see \eqref{keyest eqn23} and \eqref{keyest eqn24}, fix $(\rho,\zeta) \in \mathcal{I}$.  Suppose that $\rho \geq \rho_0/160$ and let $Y \in A^{\mathcal M}_{\rho,7/8}(\zeta) \subset B^{\mathcal M}_{61\rho_0/64}(0)$.  Provided $C_e$ and $C_h$ are sufficiently small, we can argue as in the proof of Theorem~\ref{graphical_thm} that if $\overline{R}(Y) \geq 4^{-q-8} \rho$, then \eqref{graphical thm eqn4} with $\sigma = 61\rho_0/64 + \overline{R}(Y)$ contradicts Lemma~\ref{relharm critrad2 lemma}.  Therefore, \eqref{graphical annuli go condition} holds true with $\kappa = 1/2$.  On the other hand, if $\rho < \rho_0/160$, it follows from \eqref{keyest eqn21} that \eqref{graphical annuli go condition} holds true with $\kappa = 1/2$.  Let $\overline{\beta}$ be the minimum of the values $\overline{\beta}(n,m,q,p',M_0,$ $N_0,C_e,C_h,\alpha,1/2)$ in Lemma~\ref{graphical_lemma} (with $\kappa = 1/2$) over $p' \in \{0,1,2,\ldots,\lfloor q/q_0 \rfloor \}$.  Suppose that \eqref{graphical annuli hyp} hold true.  Then there exists $n$-dimensional locally area-minimizing rectifiable currents $T_j$ of $\mathbf{A}^{\mathcal M}_{\rho,3/4}(\zeta)$ such that \eqref{graphical annuli conclB1}--\eqref{graphical annuli conclB3} hold true with $\kappa = 1/2$.  For each $j \in \{1,2,\ldots,J\}$ there exists a Lipschitz $m_j$-valued function $v_j : \op{spt} \mathbf{G}_{\varphi_j \circ \exp_{Z_0}^{-1}} \cap \mathbf{A}^{\mathcal M}_{\rho,1/2}(\zeta) \rightarrow \mathcal{A}_{m_j}(\mathbb{R}^{n+m})$ and a closed set $\mathcal{K}_j \subseteq A^{\mathcal M}_{\rho,1/2}(\zeta)$ such that \eqref{graphical annuli conclC1}--\eqref{graphical annuli conclC3} hold true.  Here $N_j(X) = \sum_{l=1}^{q_j} \llbracket N_{j,l}(X) \rrbracket$ for $X \in \mathbf{A}^{\mathcal M}_{\rho,1/2}(\zeta)$.  (Note that there is ambiguity in the notation $N_l(X)$ for the values of $N(X)$ and $N_j$ for the $q_j$-valued function as in \eqref{N of v defn}.  To reduce the ambiguity, whenever possible we shall write out the values $N_l(X)$ and $N_{j,l}(X)$ of the multi-valued function $N$ and $N_j$.)  By \eqref{graphical annuli go condition} (with $\kappa = 1/2$), \eqref{normal approx est1}, and \eqref{graphical annuli conclC1}, $|N_{j,l}(X)| \leq C \mathbf{m}_0^{1/(2n)} \rho^{1+\beta_2}$ and $|N_l(X)| \leq C \mathbf{m}_0^{1/(2n)} \rho^{1+\beta_2}$ for each $X \in A^{\mathcal M}_{\rho,1/2}(\zeta)$.  By \eqref{graphical annuli conclC1}, \eqref{normal approx est1}, and \eqref{relharm norm est rmk eqn1} and a covering argument, $N = \sum_{j=1}^J N_j$ on except on a set of $\mathcal{H}^n$-measure $\leq C \mathbf{m}_0^{3\gamma_2/4} \rho^{3\gamma_2/4-2} \int_{A^{\mathcal M}_{\rho,1}(\zeta)} |N|^2$.  Hence 
\begin{align*}
	&-\int_{\mathcal M} \sum_{l=1}^q \left( \tfrac{1}{2} r^2 |DN_l|^2 - \sum_{l=1}^q \langle r \partial_{\widehat r} N_l - \alpha N_l , 
		r \partial_{\widehat r} N_l \rangle - \alpha^2 |N_l|^2 \right) \frac{1}{R} \,\psi'(R) \,\chi \,\eta_{(\rho,\zeta)} \,d\mathcal{H}^n 
	\\ \leq\,& -\int_{\mathcal M} \sum_{j=1}^J \sum_{l=1}^{q_j} \left( \tfrac{1}{2} r^2 |DN_{j,l}|^2 - \langle r \partial_{\widehat r} N_{j,l} - \alpha N_{j,l} , 
		r \partial_{\widehat r} N_{j,l} \rangle - \alpha^2 |N_{j,l}|^2 \right) \frac{1}{R} \,\psi'(R) \,\chi \,\eta_{(\rho,\zeta)} \,d\mathcal{H}^n \nonumber 
		\\&+ C \mathbf{m}_0^{3\gamma_2/4} \rho^{3\gamma_2/4} \rho_0^{-2} \int_{A^{\mathcal M}_{\rho,1}(\zeta)} |N|^2 \,d\mathcal{H}^n , \nonumber 
\end{align*}
where $C = C(n,m,q,\alpha,M_0,N_0,C_e,C_h) \in (0,\infty)$ is a constant.  Now observe 
\begin{equation*}
	|DN_{j,l}|^2 = \sum_{i,k=1}^n g^{ik} \,\frac{\partial N_{j,l}}{\partial x_i} \,\frac{\partial N_{j,l}}{\partial x_k}, \quad 
	\partial_{\hat r} N_{j,l} = \sum_{i=1}^2 \sum_{k=1}^n g^{ik} \,\frac{x_i}{r} \,\frac{\partial N_{j,l}}{\partial x_k} , \quad
	d\mathcal{H}^n = \sqrt{g} \,dx
\end{equation*}
where $g_{ij}$ is as in \eqref{normal coord eqn3} and $dx = dx_1 \,dx_2 \dots dx_n$.  Hence using \eqref{normal coord eqn4} and Lemma~\ref{normal approx energy lemma} 
\begin{align*} 
	&-\int_{\mathcal M} \sum_{l=1}^q \left( \tfrac{1}{2} r^2 |DN_l|^2 - \sum_{l=1}^q \langle r \partial_{\widehat r} N_l - \alpha N_l , 
		r \partial_{\widehat r} N_l \rangle - \alpha^2 |N_l|^2 \right) \frac{1}{R} \,\psi'(R) \,\chi \,\eta_{(\rho,\zeta)} \,d\mathcal{H}^n 
	\\ \leq\,& -\int_{\mathcal M} \sum_{j=1}^J \sum_{l=1}^{q_j} \left( \frac{1}{2} r^2 \sum_{i=1}^n \left|\frac{\partial N_{j,l}}{\partial x_i}\right|^2 
		- \left\langle r \,\frac{\partial N_{j,l}}{\partial r} - \alpha N_{j,l}, r \,\frac{\partial N_{j,l}}{\partial r} \right\rangle - \alpha^2 |N_{j,l}|^2 \right) 
		\frac{1}{R} \,\psi'(R) \,\chi \,\eta_{(\rho,\zeta)} \,dx \nonumber 
		\\&+ C \mathbf{m}_0^{3\gamma_2/4} \rho_0^{3\gamma_2/4-2} \int_{A^{\mathcal M}_{\rho,1}(\zeta)} |N|^2 \,d\mathcal{H}^n , \nonumber 
\end{align*}
where $\frac{\partial}{\partial r} = \frac{x_1}{r} \,\frac{\partial}{\partial x_1} + \frac{x_2}{r} \,\frac{\partial}{\partial x_2}$ and $C = C(n,m,q,\alpha,M_0,N_0,C_e,C_h) \in (0,\infty)$ is a constant.  Representing $N_j$ as in \eqref{N of v defn} and using \eqref{graphical annuli conclC3} 
\begin{align}\label{keyest eqn27} 
	&-\int_{\mathcal M} \sum_{l=1}^q \left( \tfrac{1}{2} r^2 |DN_l|^2 - \sum_{l=1}^q \langle r \partial_{\widehat r} N_l - \alpha N_l , 
		r \partial_{\widehat r} N_l \rangle - \alpha^2 |N_l|^2 \right) \frac{1}{R} \,\psi'(R) \,\chi \,\eta_{(\rho,\zeta)} \,d\mathcal{H}^n 
	\\ \leq\,& -\sum_{j=1}^J q_j \int_{P_{Z_0}} \sum_{l=1}^{q_j} \left( \tfrac{1}{2} r^2 |D\varphi_{j,l}|^2 - \alpha^2 |\varphi_{j,l}|^2 \right) 
		\frac{1}{R} \,\psi'(R) \,\chi \,\eta_{(\rho,\zeta)} \nonumber
		\\&- \int_{P_{Z_0}} \sum_{j=1}^J \sum_{l=1}^{q_j} \sum_{h=1}^{m_j} \left( r^2 \langle D\varphi_{j,l}, Dv_{j,l,h} \rangle 
		- \alpha \left\langle r \frac{\partial v_{j,l,h}}{\partial r} , \varphi_{j,l,h} \right\rangle \right) 
		\frac{1}{R} \,\psi'(R) \,\chi \,\eta_{(\rho,\zeta)} \nonumber 
		\\&+ \alpha^2 \int_{P_{Z_0}} \sum_{j=1}^J \sum_{l=1}^{q_j} \sum_{h=1}^{m_j} \langle\varphi_{j,l}, v_{j,l,h} \rangle 
		\,\frac{1}{R} \,\psi'(R) \,\chi \,\eta_{(\rho,\zeta)}	
		+ \frac{C}{\rho_0^2} \int_{A^{\mathcal M}_{\rho,1}(\zeta)} \mathcal{G}(N,\varphi)^2 \nonumber
		\\&+ C \rho_0^{\varsigma/2-2} \int_{A^{\mathcal M}_{\rho,1}(\zeta)} |N|^2 , \nonumber
\end{align}
where $v_{j,l,h}$ is as in \eqref{v_localized} and $C = C(n,m,q,\alpha,M_0,N_0,C_e,C_h) \in (0,\infty)$ is a constant.  Here as a slight abuse of notation, we identify normal coordinates $(x_1,x_2,\ldots,x_n)$ on $\mathcal{M}$ at $Z_0$ with coordinates with respect to the basis $\vec\tau_1,\vec\tau_2,\ldots,\vec\tau_n$ on $P_{Z_0}$, and we thereby use $r$, $R$, and $\eta_{(\rho,\zeta)}$ to denote $r \circ \exp_{Z_0}$, $R \circ \exp_{Z_0}$, and $\eta_{(\rho,\zeta)} \circ \exp_{Z_0}$.  Hence $r(X) = \op{dist}(X,S)$ and $R(X) = |X|$ for all $X \in B_{\rho_0}(Z_0,P_{Z_0})$.  
Write $x_1 = r \cos\theta$ and $x_2 = r \sin\theta$ in polar coordinates, where $r > 0$ and $\theta \in \mathbb{R}$.  Since $\frac{\partial^2 \varphi_{j,l}}{\partial\theta^2} + \alpha^2 \varphi_{j,l} = 0$ in $P_{Z_0}$, by integration by parts in the $\theta$ variable  
\begin{align*} 
	&\int_{P_{Z_0}} \sum_{l=1}^{q_j} \left( \tfrac{1}{2} r^2 |D\varphi_{j,l}|^2 - \alpha^2 |\varphi_{j,l}|^2 \right) 
		\frac{1}{R} \,\psi'(R) \,\chi \,\eta_{(\rho,\zeta)} 
	\\=\,& \frac{1}{2} \int_{P_{Z_0}} \sum_{l=1}^{q_j} \left( \left|\frac{\partial\varphi_{j,l}}{\partial\theta}\right|^2 - \alpha^2 |\varphi_{j,l}|^2 \right) 
		\frac{1}{R} \,\psi'(R) \,\chi \,\eta_{(\rho,\zeta)} 
	\\=\,& \frac{-1}{2} \int_{P_{Z_0}} \sum_{l=1}^{q_j} \left\langle \varphi_{j,l}, \frac{\partial^2 \varphi_{j,l}}{\partial\theta^2} + \alpha^2 \varphi_{j,l} \right\rangle 
		\frac{1}{R} \,\psi'(R) \,\chi \,\eta_{(\rho,\zeta)} 
	= 0 
\end{align*}
and similarly 
\begin{align*} 
	&\int_{P_{Z_0}} \sum_{l=1}^{q_j} \sum_{h=1}^{m_j} \left( r^2 \langle D\varphi_{j,l}, Dv_{j,l,h} \rangle 
		- \alpha \left\langle r \,\frac{\partial v_{j,l,h}}{\partial r} , \varphi_{j,l} \right\rangle - \alpha^2 \langle v_{j,l,h}, \varphi_{j,l} \rangle \right) 
		\frac{1}{R} \,\psi'(R) \,\chi \,\eta_{(\rho,\zeta)} 
	\\=\,& \int_{P_{Z_0}} \sum_{l=1}^{q_j} \sum_{h=1}^{m_j} \left( \left\langle \frac{\partial\varphi_{j,l}}{\partial\theta}, \frac{\partial v_{j,l,h}}{\partial\theta} 
		\right\rangle - \alpha^2 \langle v_{j,l,h}, \varphi_{j,l} \rangle \right) 
		\frac{1}{R} \,\psi'(R) \,\chi \,\eta_{(\rho,\zeta)} 
	\\=\,& -\int_{P_{Z_0}} \sum_{l=1}^{q_j} \sum_{h=1}^{m_j} \left\langle \frac{\partial^2 \varphi_{j,l}}{\partial\theta^2} + \alpha^2 \varphi_{j,l} , v_{j,l,h} 
		\right\rangle \frac{1}{R} \,\psi'(R) \,\chi \,\eta_{(\rho,\zeta)} = 0 .
\end{align*}
Therefore, \eqref{keyest eqn27} gives us that \eqref{keyest eqn23} holds true.  By a similar argument, we obtain \eqref{keyest eqn24}.

Suppose that $(\rho,\zeta) \in \mathcal{I}$ such that \eqref{graphical annuli hyp} does not hold true.  Then we can find $\widetilde{p} \in \{0,1,2,\ldots,p-1\}$ and $\widetilde{\varphi} \in \Phi_{\alpha,q,\widetilde{p}}$ such that 
\begin{align}
	\label{keyest eqn28} E_{\varsigma/2}(N,\widetilde{\varphi},A^{\mathcal M}_{\rho,1}(\zeta)) 
		&\leq 2^p \overline{\beta}^{-p} E_{\varsigma/2}(N,\varphi,A^{\mathcal M}_{\rho,1}(\zeta)) , \\
	\label{keyest eqn29} E_{\varsigma/2}(N,\widetilde{\varphi},A^{\mathcal M}_{\rho,1}(\zeta)) 
		&\leq \overline{\beta} \inf_{\varphi' \in \bigcup_{p'=0}^{\widetilde{p}-1} \Phi_{\alpha,q,p'}} 
		E_{\varsigma/2}(N,\varphi',A^{\mathcal M}_{\rho,1}(\zeta)) .
\end{align}
Since \eqref{graphical annuli go condition} (with $\kappa = 1/2$) and \eqref{keyest eqn29} hold true, we can argue as before that \eqref{keyest eqn23} and \eqref{keyest eqn24} hold true with $\widetilde{\varphi}$ in place of $\varphi$.  Then by \eqref{keyest eqn28}, we obtain \eqref{keyest eqn23} and \eqref{keyest eqn24} with the original $\varphi$.

\noindent\textit{Step~4.  Local estimates on $B^{\mathcal M}_{15\rho_0/16}(Z_0) \setminus U$.}  Let $X \in B^{\mathcal M}_{15\rho_0/16}(Z_0)$ with $d(X) \leq r(X)/2$ and let $X = \exp_{Z_0}(x+y)$ for unique $x \in P_{Z_0} \cap S^{\perp}$ and $y \in S$.  Then there exists $Z \in B^{\mathcal M}_{15\rho_0/16}(Z_0) \cap \partial U$ such that $Z = \exp_{Z_0}(\xi+\zeta)$ for unique $\xi \in P_{Z_0} \cap S^{\perp}$ and $\zeta \in S$ and $|(x,y) - (\xi,\zeta)| \leq 3 |x|/5$.  Thus $|x| \leq 5 |\xi|/2$ and so $|(x,y) - (\xi,\zeta)| \leq 3 |\xi|/2$.  Hence 
\begin{equation*}
	\{ X \in B^{\mathcal M}_{15\rho_0/16}(Z_0) : d(X) < r(X)/2 \} 
	\subset \bigcup_{Z \in B^{\mathcal M}_{15\rho_0/16}(Z_0) \cap \partial U} B^{\mathcal M}_{3 |\xi|/2}(\exp_{Z_0}(\zeta)) .
\end{equation*}

Let $Z \in B^{\mathcal M}_{15\rho_0/16}(Z_0) \cap \partial U$ and let $Z = \exp_{Z_0}(\xi+\zeta)$, where $\xi \in P_{Z_0} \cap S^{\perp}$ and $\zeta \in S$.  Then \eqref{keyest eqn21} (with $\rho = |\xi|$) does not hold true; this is, there exists $Y \in B^{\mathcal M}_{3|\xi|/2}(\exp_{Z_0}(\zeta))$ such that $\overline{R}(Y) = 4^{-q-8} |\xi|$.  Recall that $|\xi| < \rho_0/160$.  By applying the Vitali covering lemma, there exists a countable collection $\mathcal{J}$ of $(|\xi|,\zeta)$ such that $\xi \in P_{Z_0} \cap S^{\perp}$ and $\zeta \in S$ with $Z = \exp_{Z_0}(\xi+\zeta) \in B^{\mathcal M}_{15\rho_0/16}(Z_0) \cap \partial U$, $\{B^{\mathcal M}_{2|\xi|}(\exp_{Z_0}(\zeta))\}_{(|\xi|,\zeta) \in \mathcal{J}}$ is pairwise disjoint, $\{B^{\mathcal M}_{10|\xi|}(\exp_{Z_0}(\zeta))\}_{(|\xi|,\zeta) \in \mathcal{J}}$ covers $\{ X \in B^{\mathcal M}_{15\rho_0/16}(Z_0) : d(X) < r(X)/2 \}$.  

Let $(\rho,\zeta) \in \mathcal{J}$.  Let $Y \in A^{\mathcal M}_{\rho,7/8}(\zeta)$ such that $\overline{R}(Y) = 4^{-q-8} \rho$.  By Lemma~\ref{relharm critrad2 lemma} (with $Y$ in place of $Z$) and $B^{\mathcal M}_{\overline{R}(Y)}(Y) \subset B^{\mathcal M}_{2\rho}(\exp_{Z_0}(\zeta))$ 
\begin{equation}\label{keyest eqn30} 
	\min\{C_e,C_h\} \cdot \mathbf{m}_0 \overline{R}^{2-2\delta_2} 
	\leq \frac{C}{\overline{R}^{n+2}} \int_{B^{\mathcal M}_{\overline{R}}(Y)} |N|^2 
	\leq \frac{C}{\overline{R}^{n+2}} \int_{B^{\mathcal M}_{2\rho}(\exp_{Z_0}(\zeta))} |N|^2 , 
\end{equation}
where $\overline{R} = \overline{R}(Y)$ and $C = C(n,m,q,M_0) \in (0,\infty)$ is a constant.  By applying the triangle inequality and Hypothesis~$(\dagger)$ as in Remark~\ref{graphical rmk}(2) and then using \eqref{graphical thm eqn4}, $\|\varphi\|_{L^2(B_1(0,P_{Z_0}))} \leq 4 \rho_0^{-n/2-\alpha} \|N\|_{L^2(B^{\mathcal M}_{\rho_0}(Z_0))} \leq C(n,m,q) \,\mathbf{m}_0^{1/2} \rho_0^{2-\alpha}$.  By the triangle inequality and $\alpha \geq 2$ 
\begin{align}\label{keyest eqn31} 
	\int_{B^{\mathcal M}_{2\rho}(\exp_{Z_0}(\zeta))} |N|^2 
	\leq\,& 2 \int_{B^{\mathcal M}_{2\rho}(\exp_{Z_0}(\zeta))} \mathcal{G}(N,\varphi\circ\exp_{Z_0}^{-1})^2 
		+ 4 \int_{B_{2\rho}(\zeta,P_{Z_0})} |\varphi|^2 
	\\ \leq\,& 2 \int_{B^{\mathcal M}_{2\rho}(\exp_{Z_0}(\zeta))} \mathcal{G}(N,\varphi\circ\exp_{Z_0}^{-1})^2 
		+ \frac{C \mathbf{m}_0 \rho^{n+2\alpha}}{\rho_0^{2\alpha-4}} \nonumber
	\\ \leq\,& 2 \int_{B^{\mathcal M}_{2\rho}(\exp_{Z_0}(\zeta))} \mathcal{G}(N,\varphi\circ\exp_{Z_0}^{-1})^2 
		+ C \mathbf{m}_0 \overline{R}^{n+4} , \nonumber 
\end{align}
where $C = C(n,m,q) \in (0,\infty)$ are constants.  Hence provided $C_e$ and $C_h$ are sufficiently large, by \eqref{keyest eqn30} and \eqref{keyest eqn31}
\begin{equation}\label{keyest eqn32} 
	\int_{B_{2\rho}(\exp_{Z_0}(\zeta))} \mathcal{G}(N,\varphi\circ\exp_{Z_0}^{-1})^2 \geq \mathbf{m}_0 \overline{R}^{n+4-2\delta_2} .
\end{equation}
Moreover, by \eqref{keyest eqn31} and \eqref{keyest eqn32} 
\begin{equation}\label{keyest eqn33} 
	\int_{B^{\mathcal M}_{2\rho}(\exp_{Z_0}(\zeta))} |N|^2 \leq C \int_{B_{2\rho}(\exp_{Z_0}(\zeta))} \mathcal{G}(N,\varphi\circ\exp_{Z_0}^{-1})^2 , 
\end{equation}
where $C = C(n,m,q) \in (0,\infty)$ is a constant.  By~\cite[Lemmas~7.12 and 7.16]{GilbargTrudinger} applied to $N \circ \exp_{Z_0}$ 
\begin{equation*}
	\int_{B^{\mathcal M}_{10\rho}(\exp_{Z_0}(\zeta))} \big| |N| - \mu \big|^2 
	\leq C \rho^2 \int_{B^{\mathcal M}_{10\rho}(\exp_{Z_0}(\zeta))} |DN|^2 
\end{equation*}
where $\mu = \frac{1}{\omega_n (2\rho)^n} \int_{B_{2\rho}(\zeta,P_{Z_0})} |N \circ \exp_{Z_0}^{-1}|$ and $C = C(n,m,q,M_0,N_0,C_e,C_h) \in (0,\infty)$ is a constant.  Hence 
\begin{align}\label{keyest eqn34} 
	\int_{B^{\mathcal M}_{10\rho}(\exp_{Z_0}(\zeta))} |N|^2 
	\leq\,& C |\mu|^2 \rho^n + C \rho^2 \int_{B^{\mathcal M}_{5s}(\exp_{Z_0}(\zeta))} |DN|^2 
	\\ \leq\,& C \int_{B^{\mathcal M}_{2\rho}(\exp_{Z_0}(\zeta))} |N|^2 
		+ C \rho^2 \int_{B^{\mathcal M}_{10\rho}(\exp_{Z_0}(\zeta))} |DN|^2 , \nonumber 
\end{align}
where $C = C(n,m,q,M_0,N_0,C_e,C_h) \in (0,\infty)$ are constants.  Since there exists $\xi \in P_{Z_0} \cap S^{\perp}$ such that $\rho = |\xi|$ and $Z = \exp_{Z_0}(\xi+\zeta) \in B^{\mathcal M}_{15\rho_0/16}(Z_0) \cap \partial U$, by \eqref{keyest eqn21} we have that $\overline{R}(\exp_{Z_0}(\zeta)) \leq 4^{-q-8} \rho < 2 \rho$.  Thus by \eqref{normal approx est1} and \eqref{keyest eqn32} 
\begin{equation}\label{keyest eqn35} 
	\rho^2 \int_{B^{\mathcal M}_{10\rho}(\exp_{Z_0}(\zeta))} |DN|^2 
	\leq C \mathbf{m}_0 \rho^{n+4-2\delta_2} 
	\leq C \int_{B_{2\rho}(\exp_{Z_0}(\zeta))} \mathcal{G}(N,\varphi\circ\exp_{Z_0}^{-1})^2 , 
\end{equation}
where $C = C(n,m,q,M_0,N_0,C_e,C_h) \in (0,\infty)$ are constants.  By \eqref{keyest eqn33}, \eqref{keyest eqn34}, and \eqref{keyest eqn35} 
\begin{equation}\label{keyest eqn36} 
	\int_{B^{\mathcal M}_{10\rho}(\exp_{Z_0}(\zeta))} |N|^2 + \rho^2 \int_{B^{\mathcal M}_{10\rho}(\exp_{Z_0}(\zeta))} |DN|^2 
	\leq C \int_{B_{2\rho}(\exp_{Z_0}(\zeta))} \mathcal{G}(N,\varphi\circ\exp_{Z_0}^{-1})^2 , 
\end{equation}
where $C = C(n,m,q,M_0,N_0,C_e,C_h) \in (0,\infty)$ is a constant.  By summing \eqref{keyest eqn36} over all $(\rho,\zeta) \in \mathcal{J}$ 
\begin{equation}\label{keyest eqn37} 
	\int_{B^{\mathcal M}_{15\rho_0/16}(Z_0) \cap \{ d(X) < r(X)/2 \}} (|N|^2 + r^2 |DN|^2) 
	\leq C \int_{B^{\mathcal M}_{\rho_0}(Z_0)} \mathcal{G}(N,\varphi\circ\exp_{Z_0}^{-1})^2 , 
\end{equation}
where $C = C(n,m,q,M_0,N_0,C_e,C_h) \in (0,\infty)$ is a constant.  By \eqref{keyest eqn37}  
\begin{align}\label{keyest eqn38} 
	&-\int_{\mathcal M} \left( \tfrac{1}{2} r^2 |DN|^2 - \sum_{l=1}^q \langle r \partial_{\widehat r} N_l - \alpha N_l , r \partial_{\widehat r} N_l \rangle 
		- \alpha^2 |N|^2 \right) \frac{1}{R} \,\psi'(R) \,(1-\chi) 
	\\&\leq \frac{C}{\rho_0^2} \int_{B^{\mathcal M}_{\rho_0}(Z_0)} \mathcal{G}(N,\varphi\circ\exp_{Z_0}^{-1})^2 \nonumber
\end{align}
and 
\begin{equation}\label{keyest eqn39} 
	\int_{\mathcal M} \sum_{l=1}^q |r \partial_{\widehat r} N_l - \alpha N_l|^2 \,(1-\chi) 
	\leq \frac{C}{\rho_0^2} \int_{B^{\mathcal M}_{\rho_0}(Z_0)} \mathcal{G}(N,\varphi\circ\exp_{Z_0}^{-1})^2 ,  
\end{equation}
where $C = C(n,m,q,M_0,N_0,C_e,C_h,\alpha) \in (0,\infty)$ are constants.  Therefore, by combining \eqref{keyest eqn20}, \eqref{keyest eqn25}, \eqref{keyest eqn26}, \eqref{keyest eqn38}, and \eqref{keyest eqn39}, we obtain \eqref{keyest concl}.
\end{proof}

Throughout the remainder of the paper, we shall assume that the constants $M_0,N_0,C_e,C_h,\varepsilon_2$ satisfy the following: 

\begin{hypothesis}\label{noncon constants hypothesis}{\rm
Let $M_0,N_0,C_e,C_h,\varepsilon_2$ be constants satisfying the requirements of Hypothesis~\ref{DeLSpa consts freq hypothesis} such that Lemma~\ref{relharm critrad2 lemma}, Theorem~\ref{graphical_thm}, Lemma~\ref{graphical_lemma} and Theorem~\ref{keyest thm} hold true.
}\end{hypothesis}

\begin{corollary}\label{radial noncon cor}
Let $\alpha = \ell_0/q_0 \geq 2$ for some relatively prime positive integers $\ell_0$ and $q_0 \leq q$.  Let $p \in \{1,2,\ldots, \lfloor q/q_0 \rfloor \}$ and $\sigma \in (0,1)$.  Let $\varsigma = \varsigma(m,n,q) \in (0,1)$ be as in Theorem~\ref{relharm separation main thm}.  Let $M_0,N_0,C_e,C_h,\varepsilon_2$ be as in Hypothesis~\ref{noncon constants hypothesis}.  There exists there exist $\varepsilon_0, \beta \in (0,1)$ depending only on $n,m,q,p,\alpha,M_0,N_0,C_e,C_h$ such that if Hypothesis~$(\star)$ (in Section~\ref{sec:graph_rep_results}) and Hypothesis~$(\dagger)$ hold true and if $\mathbf{I}_{N,Z_0}(0^+) \geq \alpha$, then 
\begin{equation}\label{radial noncon concl}
	\rho_0^{2\alpha-2-\sigma} \int_{B^{\mathcal M}_{3\rho_0/4}(Z_0)} \frac{\mathcal{G}(N,\varphi)^2}{R^{n+2\alpha-\sigma}} 
	\leq C E_{\varsigma/2}(N,\varphi,B^{\mathcal M}_{\rho_0}(Z_0))^2 ,
\end{equation}
where $R$ is the geodesic distance from $Z_0$ (as in~\eqref{normal coord eqn2}) and $C = C(n,m,q,M_0,N_0,C_e,C_h,\alpha,$ $\sigma) \in (0,\infty)$ is a constant.
\end{corollary}

\begin{proof}
Let $\psi : [0,\rho_0] \rightarrow \mathbb{R}$ be a smooth function such that $\psi(s) = 1$ if $0 \leq s \leq 3\rho_0/4$, $\psi(s) = 0$ if $7\rho_0/8 \leq s \leq 1$, and $-12/\rho_0 \leq \psi'(s) \leq 0$.  For each $\delta \in (0,1/2]$, let $\eta_{\delta} : [0,\rho_0] \rightarrow \mathbb{R}$ be a smooth function such that $\eta_{\delta}(s) = 0$ if $0 \leq s \leq \delta\rho_0/2$, $\eta_{\delta}(s) = 1$ if $\delta\rho_0 \leq s \leq \rho_0$, and $0 \leq \eta'_{\delta}(s) \leq \frac{3}{\delta\rho_0}$.  By the divergence theorem $\int_{\mathcal M} \op{div}_{\mathcal M} \xi \,d\mathcal{H}^n = 0$ with $\xi = R^{1-n-2\alpha+\sigma} \,\mathcal{G}(N,\varphi)^2 \,\psi(R)^2 \,\eta_{\delta}(R) \,\widehat{R}$ 
\begin{align}\label{radial noncon eqn1}
	&\sigma \int_{\mathcal M} \frac{\mathcal{G}(N,\varphi)^2}{R^{n+2\alpha-\sigma}} \,\psi(R)^2 \,\eta_{\delta}(R) 
	= -\int_{\mathcal M} \frac{\mathcal{G}(N,\varphi)^2}{R^{n-1+2\alpha-\sigma}} \,\bigg( \Delta_{\mathcal M} R - \frac{n-1}{R} \bigg) 
		\,\psi(R)^2 \,\eta_{\delta}(R) \\&\hspace{15mm} - 2 \int_{\mathcal M} \frac{\mathcal{G}(N,\varphi)^2}{R^{n-1+2\alpha-\sigma}} 
		\,\psi(R) \,\psi'(R) \,\eta_{\delta}(R) - \int_{\mathcal M} \frac{\mathcal{G}(N,\varphi)^2}{R^{n-1+2\alpha-\sigma}} \,\psi(R)^2 \,\eta'_{\delta}(R) 
		\nonumber \\&\hspace{15mm} - \int_{\mathcal M} R^{1-n+\sigma} \,\partial_{\hat R} \bigg(\frac{\mathcal{G}(N,\varphi)^2}{R^{2\alpha}} \bigg) 
		\,\psi(R)^2 \,\eta_{\delta}(R) \nonumber
\end{align}
To bound the first term on the right-hand side of \eqref{radial noncon eqn1}, let $(x_1,x_2,\ldots,x_n)$ be normal coordinates of $\mathcal{M}$ at $Z_0$ as in \eqref{normal coord eqn1} and let $g_{ij}$ be as in \eqref{normal coord eqn3}.  By \eqref{normal coord eqn4} 
\begin{equation}\label{radial noncon eqn2}
	\left| \Delta_{\mathcal M} R - \frac{n-1}{R} \right| \leq C \mathbf{m}_0^{1/2} 
\end{equation}
for some constant $C = C(n,m,q,\alpha,\sigma) \in (0,\infty)$.  The third term on the right-hand side of \eqref{radial noncon eqn1} is $\leq 0$.  To bound the fourth term on the right-hand side of \eqref{radial noncon eqn1}, by the homogeneity of $\varphi$
\begin{equation*}
	\left| \partial_{\hat R} \bigg(\frac{\mathcal{G}(N,\varphi)^2}{R^{2\alpha}} \bigg) \right| 
	= \left| \partial_{\hat R} \bigg( \mathcal{G}\bigg(\frac{N}{R^{\alpha}}, \frac{\varphi}{R^{\alpha}}\bigg)^2 \bigg) \right| 
	\leq \frac{2 \mathcal{G}(N,\varphi)}{R^{\alpha}} \left| \partial_{\hat R} \bigg(\frac{N}{R^{\alpha}}\bigg) \right| . 
\end{equation*}
Hence by Cauchy's inequality and Theorem~\ref{keyest thm} 
\begin{align}\label{radial noncon eqn3}
	&\int_{\mathcal M} R^{1-n+\sigma} \left| \partial_{\hat R} \bigg(\frac{\mathcal{G}(N,\varphi)^2}{R^{2\alpha}} \bigg) \right| \psi(R)^2 \,\eta_{\delta}(R) 
	\\ \leq\,& \frac{\sigma}{2} \int_{\mathcal M} \frac{\mathcal{G}(N,\varphi)^2}{R^{n+2\alpha-\sigma}} \,\psi(R)^2 \,\eta_{\delta}(R) 
		+ \frac{2}{\sigma} \,\rho_0^{\sigma} \int_{\mathcal M} R^{2-n} \left| \partial_{\hat R} \bigg(\frac{N}{R^{\alpha}}\bigg) \right|^2 \psi(R)^2 \nonumber 
	\\ \leq\,& \frac{\sigma}{2} \int_{\mathcal M} \frac{\mathcal{G}(N,\varphi)^2}{R^{n+2\alpha-\sigma}} \,\psi(R)^2 \,\eta_{\delta}(R) 
		+ \frac{C}{\sigma} \,\rho_0^{-2\alpha+2+\sigma} E_{\varsigma/2}(N,\varphi,B^{\mathcal M}_{\rho_0}(Z_0))^2 , \nonumber 
\end{align}
where $C = C(n,m,q,M_0,N_0,C_e,C_h,\alpha) \in (0,\infty)$ is a constant.  Hence using \eqref{radial noncon eqn2} and \eqref{radial noncon eqn3} to bound the right-hand side of \eqref{radial noncon eqn1} 
\begin{equation*}
	\bigg( \frac{\sigma}{2} - C\mathbf{m}_0^{1/2} \rho_0 \bigg) \int_{\mathcal M} \frac{\mathcal{G}(N,\varphi)^2}{R^{n+2\alpha-\sigma}} \,\psi(R)^2 
		\,\eta_{\delta}(R) \leq C \rho_0^{-2\alpha+2+\sigma} E_{\varsigma/2}(N,\varphi,B^{\mathcal M}_{\rho_0}(Z_0))^2 , 
\end{equation*}
where $C = C(n,m,q,M_0,N_0,C_e,C_h,\alpha,\sigma) \in (0,\infty)$ is a constant.  Assuming $\varepsilon_0$ is sufficiently small, we obtain \eqref{radial noncon concl}.
\end{proof}

\section{A priori estimates: Part~II}\label{sec:nonconcentration} 

Assume that Hypothesis~$(\star)$ and Hypothesis~$(\dagger)$ hold true.  In Theorem~\ref{nonconcentration thm} below we obtain several further integral estimates for $N$.  A key consequence of these estimates is Corollary~\ref{nonconcentration cor}, which rules out concentration of the $L^2$-excess of $T$ relative to $\varphi$ near the set $\mathcal{B}^{(\geq 2)}_q \cap \{ Z : \mathbf{I}_{N,Z}(0^+) \geq \alpha \}$.  The results below use adaptations of the corresponding arguments in~\cite{Sim93} and~\cite{Wic14}.

Recall that $\mathcal{M}$ is the center manifold of $T$, $Z_0 \in \mathcal{B}^{(\geq 2)}_q$, and $P_{Z_0} = {\rm T}_{Z_0} \mathcal{M}$ and let $Z = \exp_{Z_0}(\xi+\zeta)$ where $\xi \in P_{Z_0} \cap S^{\perp}$ and $\zeta \in S$.  Let $P_Z = {\rm T}_Z \mathcal{M}$.  Since $\|\bs{\varphi}\|_{C^3} \leq C \mathbf{m}_0^{1/2}$, there exists a rotation $\Gamma_Z$ of $\mathbb{R}^{n+m}$ such that $\Gamma_Z(P_{Z_0}) = P_Z$ and 
\begin{align}\label{noncon rotate}
	\|\Gamma_Z - I\| \leq C \mathbf{m}_0^{1/2} |Z - Z_0|
\end{align}
where $I$ is the identity map on $\mathbb{R}^{n+m}$ and $C = C(n,m,q,M_0,N_0,C_e,C_h) \in (0,\infty)$ is a constant.  Let $\widehat{\varphi} = \Gamma_Z \circ \varphi \circ \Gamma_Z^{-1} \in \Phi_{\alpha,q,p}(P_Z,\Gamma_Z(S))$.

\begin{theorem}\label{nonconcentration thm} 
Let $\alpha = \ell_0/q_0 \geq 2$ for some relatively prime positive integers $\ell_0$ and $q_0 \leq q$.  Let $p \in \{1,2,\ldots, \lfloor q/q_0 \rfloor \}$, $\tau \in (0,1/2]$, and $\sigma \in (0,1)$.  Let $\varsigma = \varsigma(m,n,q) \in (0,1)$ be as in Theorem~\ref{relharm separation main thm}.  Let $M_0,N_0,C_e,C_h,\varepsilon_2$ be as in Hypothesis~\ref{noncon constants hypothesis}.  There exist $\varepsilon_0, \beta \in (0,1)$ depending only on $n,m,q,p,M_0,N_0,C_e,C_h,\alpha,\tau$ such that the following holds true.  Suppose that Hypothesis~$(\star)$ (of Section~\ref{sec:graph_rep_results}) and Hypothesis~$(\dagger)$ (of Section~\ref{sec:key_est_sec}) hold true.  Let $Z \in \mathcal{B}^{(\geq 2)}_q \cap B^{\mathcal M}_{\rho_0/2}(Z_0)$ such that $\mathbf{I}_{N,Z}(0^+) \geq \alpha$ and let $Z = \exp_{Z_0}(\xi+\zeta)$ where $\xi \in P_{Z_0} \cap S^{\perp}$ and $\zeta \in S$.  Then 
\begin{align*}
	&(A) \quad \widehat{E} \,|\xi| \leq C \rho_0 \,E_{\varsigma/2}(N,\varphi,B^{\mathcal M}_{\rho_0}(Z_0)), \\ 
	&(B) \quad \frac{1}{\rho_0^{n+2}} \int_{B^{\mathcal M}_{\rho_0}(Z_0)} \mathcal{G}(N,\widehat{\varphi} \circ \exp_Z^{-1})^2 
		\leq C E_{\varsigma/2}(N,\varphi,B^{\mathcal M}_{\rho_0}(Z_0))^2 , 
\end{align*}
where $\widehat{E} = \rho_0^{-n/2-1} \|N\|_{L^2(B^{\mathcal M}_{\rho_0}(Z_0))}$ and $C = C(n,m,q,p,M_0,N_0,C_e,C_h,\alpha) \in (0,\infty)$ is a constant.  Moreover, 
\begin{align*}
	(C) \quad& \frac{1}{\rho_0^{2+\sigma}} \int_{B^{\mathcal M}_{3\rho_0/4}(Z_0)} \frac{\mathcal{G}(N,\varphi \circ \exp_{Z_0}^{-1})^2}{R_Z^{n-\sigma}} \,d\mathcal{H}^n 
		\leq C E_{\varsigma/2}(N,\varphi,B^{\mathcal M}_{\rho_0}(Z_0))^2 , \\
	(D) \quad& \rho_0^{2\alpha-2-\sigma} \int_{B^{\mathcal M}_{3\rho_0/4}(Z_0) \cap \{r > \tau\rho_0\}} \frac{\mathcal{G}(N(X), 
		\varphi(\exp_{Z_0}^{-1}(X)) - D_x \varphi(\exp_{Z_0}^{-1}(X)) \cdot \xi)^2}{R_Z^{n+2\alpha-\sigma}} \,d\mathcal{H}^n(X)
		\\&\leq C E_{\varsigma/2}(N,\varphi,B^{\mathcal M}_{\rho_0}(Z_0))^2 
\end{align*}
for some constant $C = C(n,m,q,p,M_0,N_0,C_e,C_h,\alpha,\sigma) \in (0,\infty)$.  Here for each $X \in B^{\mathcal M}_{\rho_0}(Z_0)$, we let $r(X) = \op{dist}(\exp_{Z_0}^{-1}(X),S)$ (as in \eqref{normal coord eqn2}) and we let $R_Z = R_Z(X)$ denote the geodesic distance of $X$ to $Z$ in $\mathcal{M}$.  We let 
\begin{equation*}
	\varphi(X) - D\varphi(X) \cdot \xi 
	= \sum_{j=1}^J \sum_{l=1}^{q_j} \llbracket \varphi_{j,l}(X) - D\varphi_{j,l}(X) \cdot \xi \rrbracket 
\end{equation*}
for each $X \in P_{Z_0}$, where $\varphi_{j,l}$ are locally defined by \eqref{varphi_localized}.
\end{theorem}

We will prove Theorem~\ref{nonconcentration thm} by induction on $p$, assuming that: 
\begin{enumerate}[topsep=0mm]
	\item[(H3)]  whenever $\widetilde{p} \in \{1,2,\ldots,p-1\}$, Theorem~\ref{nonconcentration thm} holds true with $\widetilde{p}$ in place of $p$.
\end{enumerate}
First we have a few preliminary lemmas, Lemma~\ref{nonconcentration3 lemma}, Lemma~\ref{nonconcentration4 lemma}, and Lemma~\ref{nonconcentration5 lemma}.

\begin{lemma}\label{nonconcentration3 lemma} 
Let $\alpha = \ell_0/q_0 \geq 2$ for some relatively prime positive integers $\ell_0$ and $q_0 \leq q$.  Let $p \in \{1,2,\ldots, \lfloor q/q_0 \rfloor \}$.  Let $M_0,N_0,C_e,C_h,\varepsilon_2$ be as in Hypothesis~\ref{noncon constants hypothesis}.  For every $\delta \in (0,1)$ there exist $\varepsilon_0, \beta \in (0,1)$ depending only on $n,m,q,p,M_0,N_0,C_e,C_h,\alpha,\delta$ such that the following holds true.  Suppose that Hypothesis~$(\star)$ (of Section~\ref{sec:graph_rep_results}) and Hypothesis~$(\dagger)$ (of Section~\ref{sec:key_est_sec}) holds true.  Let $Z \in \mathcal{B}^{(\geq 2)}_q \cap B^{\mathcal M}_{\rho_0/2}(Z_0)$ such that $\mathbf{I}_{N,Z}(0^+) \geq \alpha$ and let $Z = \exp_{Z_0}(\xi+\zeta)$ where $\xi \in P_{Z_0} \cap S^{\perp}$ and $\zeta \in S$.  Then 
\begin{equation}\label{noncon3 concl} 
	|\xi| \leq \delta \rho_0 .
\end{equation}
\end{lemma}

\begin{proof} 
We proceed by contradiction.  Fix $\delta \in (0,1)$ and $p \in \{1,2,\ldots,\lfloor q/q_0 \rfloor\}$.  For $\nu = 1,2,3,\ldots$, suppose that $\varepsilon_{\nu} \rightarrow 0^+$, $\beta_{\nu} \rightarrow 0^+$, $T^{(\nu)}$, $\mathcal{M}^{(\nu)}$, $N^{(\nu)}$, $\mathcal{K}_{\nu}$, $Z_{0,\nu}$, $\rho_{\nu}$, $P_{\nu}$, $S_{\nu}$, and $\varphi^{(\nu)}$ satisfy Hypothesis~$(\star)$ and Hypothesis~$(\dagger)$ in place of $\varepsilon_0,\beta,T,\mathcal{M},N.\mathcal{K},Z_0,\rho,P_{Z_0},S,\varphi$.  Let $Z_{\nu} \in \mathcal{B}^{(\geq 2)}_q(T^{(\nu)}) \cap B^{\mathcal M}_{\rho_{\nu}/2}(Z_{0,\nu})$ such that $\mathbf{I}_{N^{(\nu)},Z_{\nu}}(0^+) \geq \alpha$.  Let $Z_{\nu} = \mathbf{e}_{\nu}(\xi_{\nu} + \zeta_{\nu})$, where $\xi_{\nu} \in P_{\nu} \cap S_{\nu}^{\perp}$ and $\zeta_{\nu} \in S_{\nu}$, where $\mathbf{e}_{\nu}$ is the exponential map of $\mathcal{M}^{(\nu)}$ at $Z_{0,\nu}$.  Suppose that $|\xi_{\nu}| \geq \delta\rho_{\nu}$.  

Let $\mathfrak{q}_{\nu}$ be a rotation of $\mathbb{R}^{n+m}$ such that $\mathfrak{q}_{\nu}(P_{\nu}) = \mathbb{R}^n \times \{0\}$ and $\mathfrak{q}_{\nu}(S_{\nu}) = \{0\}^2 \times \mathbb{R}^{n-2} \times \{0\}^m$.  Rotate and rescale letting $\overline{T}^{(\nu)}$, $\overline{\mathcal M}^{(\nu)}$, and $\overline{N}^{(\nu)}$ be as in \eqref{centered blowup eqn1}, $\overline{\varphi}^{(\nu)}(X) = \rho_{\nu}^{\alpha-1} (\mathfrak{q}_{\nu} \circ \varphi^{(\nu)} \circ \mathfrak{q}_{\nu}^{-1})(X)$ for each $X \in \mathbb{R}^n$, and $(\overline{\xi}_{\nu},\overline{\zeta}_{\nu}) = \mathfrak{q}_{\nu}(\xi_{\nu} + \zeta_{\nu})/\rho_{\nu}$ where $\overline{\xi}_{\nu} \in \mathbb{R}^2$ and $\overline{\zeta}_{\nu} \in \mathbb{R}^{n-2}$.  Following the blow-up procedure from Subsection~\ref{sec:centered blowup sec}, let $\varphi^{(0)}$ be a coarse blow-up of $T^{(\nu)}$ relative to $\mathcal{M}^{(\nu)}$ in $B^{\mathcal{M}^{(\nu)}}_{2\rho_{0,\nu}}(Z_{0,\nu})$.  By Hypothesis~$(\dagger)$, $\overline{\varphi}^{(\nu)}/\widehat{E}_{\nu} \rightarrow \varphi^{(0)}$ uniformly on compact subsets of $\mathbb{R}^n$, where $\widehat{E}_{\nu} = \rho_{\nu}^{-n/2-1} \|N^{(\nu)}\|_{L^2(B^{\mathcal{M}^{(\nu)}}_{\rho_{\nu}}(Z_{0,\nu})}$.  Hence $\varphi^{(0)}$ is homogeneous degree $\alpha$ and translation invariant along $\{0\} \times \mathbb{R}^{n-2}$.  Note that by Remark~\ref{graphical rmk}(2), $\varphi^{(0)}$ is not identically zero on $\mathbb{R}^n$.  After passing to a subsequence, let $\overline{\xi}_{\nu} \rightarrow \overline{\xi}$ in $\mathbb{R}^2$ and $\overline{\zeta}_{\nu} \rightarrow \overline{\zeta}$ in $\mathbb{R}^{n-2}$ as $\nu\rightarrow\infty$.  By Lemma~\ref{centered semicont2 lemma}
\begin{equation*}
	\mathcal{N}_{\varphi^{(0)}}(\overline{\xi},\overline{\zeta}) \geq \limsup_{\nu \rightarrow \infty} \mathbf{I}_{N,Z_{\nu}}(0^+) \geq \alpha 
\end{equation*}
and thus $(\overline{\xi},\overline{\zeta})$ is in the spine of $\varphi^{(0)}$, i.e.~$\overline{\xi} = 0$.  Hence $|\xi_{\nu}|/\rho_{\nu} = |\overline{\xi}_{\nu}| \rightarrow 0$ as $\nu\rightarrow\infty$, contradicting $|\xi_{\nu}| \geq \delta\rho_{\nu}$.  
\end{proof}

\begin{lemma}\label{nonconcentration4 lemma} 
Let $\alpha = \ell_0/q_0 \geq 2$ for some relatively prime positive integers $\ell_0$ and $q_0 \leq q$.  Let $p \in \{1,2,\ldots, \lfloor q/q_0 \rfloor \}$ such that (H3) holds true.  Let $M_0,N_0,C_e,C_h,\varepsilon_2$ be as in Hypothesis~\ref{noncon constants hypothesis}.  For every $\delta \in (0,1)$ there exist $\varepsilon_0, \beta, \gamma_0 \in (0,1)$ depending only on $n,m,q,p,M_0,N_0,C_e,$ $C_h,\alpha,\delta$ such that the following holds true.  Suppose that Hypothesis~$(\star)$ and Hypothesis~$(\star\star)$ (of Section~\ref{sec:graph_rep_results}) hold true and that there exists $s \in \{1,2,\ldots,p-1\}$ such that 
\begin{equation}\label{noncon4 hyp} 
	\inf_{\varphi' \in \bigcup_{p'=0}^s \Phi_{\alpha,q,p'}} E_{\varsigma/2}(N,\varphi',B^{\mathcal M}_{\rho_0}(Z_0)) 
		\leq \gamma_0 \inf_{\varphi' \in \bigcup_{p'=0}^{s-1} \Phi_{\alpha,q,p'}} E_{\varsigma/2}(N,\varphi',B^{\mathcal M}_{\rho_0}(Z_0)) .
\end{equation}
Let $Z \in \mathcal{B}^{(\geq 2)}_q \cap B^{\mathcal M}_{\rho_0/2}(Z_0)$ such that $\mathbf{I}_{N,Z}(0^+) \geq \alpha$ and let $Z = \exp_{Z_0}(\xi+\zeta)$ where $\xi \in P_{Z_0} \cap S^{\perp}$ and $\zeta \in S$.  Then 
\begin{equation}\label{noncon4 concl} 
	\widehat{E} \,|\xi| \leq \delta \rho_0 \inf_{\varphi' \in \bigcup_{p'=0}^s \Phi_{\alpha,q,p'}} E_{\varsigma/2}(N,\varphi',B^{\mathcal M}_{\rho_0}(Z_0)) . 
\end{equation}
where $\widehat{E} = \rho_0^{-n/2-1} \|N\|_{L^2(B^{\mathcal M}_{\rho_0}(Z_0))}$.
\end{lemma}

\begin{proof}
Fix $p \in \{2,3,\ldots,\lfloor q/q_0 \rfloor\}$ such that (H3) holds true and fix $\delta \in (0,1)$ and $s \in \{1,2,\ldots,p-1\}$.  Suppose to the contrary that there exists $\varepsilon_{\nu} \rightarrow 0^+$, $\beta_{\nu} \rightarrow 0^+$, $T^{(\nu)}$, $\mathcal{M}^{(\nu)}$, $N^{(\nu)}$,$Z_{0,\nu}$, $\rho_{\nu}$, $P_{\nu}$, $S_{\nu}$, and $\varphi^{(\nu)}$ satisfying Hypothesis~$(\star)$, Hypothesis~$(\dagger)$, and \eqref{noncon4 hyp} but not \eqref{noncon4 concl} in place of $\varepsilon_0$, $\beta$, $T$, $\mathcal{M}$, $N$, $Z_0$, $\rho$, $P_{Z_0}$, $S$, and $\varphi$ respectively.  Let $\widetilde{\varphi}^{(\nu)} \in \Phi_{\alpha,q,s}$ such that 
\begin{equation*} 
	E_{\varsigma/2}(N^{(\nu)},\widetilde{\varphi}^{(\nu)},B^{\mathcal{M}^{(\nu)}}_{\rho_{\nu}}(Z_{0,\nu})) 
		\leq 2 \inf_{\varphi' \in \Phi_{\alpha,q,s}} E_{\varsigma/2}(N^{(\nu)},\varphi',B^{\mathcal{M}^{(\nu)}}_{\rho_{\nu}}(Z_{0,\nu})) .
\end{equation*}
We argue exactly as in~\cite[Lemma~8.5]{KrumWic2} by constructing a blow-up $w$ of $T^{(\nu)}$ relative to $\widetilde{\varphi}^{(\nu)}$ in $B^{\mathcal{M}^{(\nu)}}_{\rho_{\nu}}(Z_{0,\nu})$ as in Subsection~\ref{sec:blowup subsec} and then use Theorem~\ref{nonconcentration thm}(A)(D) to reach a contradiction.
\end{proof}

\begin{lemma}\label{nonconcentration5 lemma}
Let $\alpha = \ell_0/q_0 \geq 2$ for some relatively prime positive integers $\ell_0$ and $q_0 \leq q$.  Let $p \in \{1,2,\ldots, \lfloor q/q_0 \rfloor \}$ such that (H3) holds true.  Let $M_0,N_0,C_e,C_h,\varepsilon_2$ be as in Hypothesis~\ref{noncon constants hypothesis}.  For every $\delta \in (0,1)$ there exist $\varepsilon_0, \beta \in (0,1)$ depending only on $n,m,q,p,M_0,N_0,C_e,C_h,\alpha,\delta$ such that the following holds true.  Suppose that Hypothesis~$(\star)$ and Hypothesis~$(\star\star)$ (of Section~\ref{sec:graph_rep_results}) hold true.  Let $Z \in \mathcal{B}^{(\geq 2)}_q \cap B^{\mathcal M}_{\rho_0/2}(Z_0)$ such that $\mathbf{I}_{N,Z}(0^+) \geq \alpha$ and let $Z = \exp_{Z_0}(\xi+\zeta)$ where $\xi \in P_{Z_0} \cap S^{\perp}$ and $\zeta \in S$.  Then 
\begin{equation}\label{noncon5 concl} 
	\widehat{E} \,|\xi| \leq \delta \rho_0 \inf_{\varphi' \in \bigcup_{p'=0}^{p-1} \Phi_{\alpha,q,p'}} E_{\varsigma/2}(N,\varphi',B^{\mathcal M}_{\rho_0}(Z_0)) . 
\end{equation}
where $\widehat{E} = \rho_0^{-n/2-1} \|N\|_{L^2(B^{\mathcal M}_{\rho_0}(Z_0))}$.
\end{lemma}

\begin{proof}
Lemma~\ref{nonconcentration5 lemma} follows from Lemma~\ref{nonconcentration3 lemma} and Lemma~\ref{nonconcentration4 lemma} by arguing as in~\cite[Lemma~8.4]{KrumWic2} (for which Lemma~8.4 follows from Lemma~8.5).
\end{proof}

\begin{remark}\label{noncon triangle rmk}{\rm 
We make the following observations which we shall use to prove of Theorem~\ref{nonconcentration thm}.  

\noindent (1)  \textit{There exists $\varepsilon_0 \in (0,1)$ depending only on $n,m,q,M_0,N_0,C_e,C_h$ such that if $X \in B^{\mathcal M}_{\rho_0}(Z_0)$ with $|\xi| < r(X)/2$, then }
\begin{align}\label{noncon triangle1} 
	\mathcal{G}(N(X),\widehat{\varphi}(\exp_Z^{-1}(X))) \leq\,& \mathcal{G}(N(X),\varphi(\exp_{Z_0}^{-1}(X)) - D\varphi(\exp_{Z_0}^{-1}(X)) \cdot \xi)
	\\&+ C \widehat{E} \rho_0^{-1} |\xi|^2 + C \mathbf{m}_0^{1/2} \widehat{E} \rho_0^2 \nonumber 
\end{align}
\textit{for all $X \in B^{\mathcal M}_{\rho_0}(Z_0)$, where $\widehat{\varphi} = \Gamma_Z \circ \varphi \circ \Gamma_Z^{-1}$, $\widehat{E} = \rho_0^{-n/2-1} \|N\|_{L^2(B^{\mathcal M}_{\rho_0}(Z_0))}$, and $C = C(n,m,q,\alpha,$ $M_0,N_0,C_e,C_h) \in (0,\infty)$ is a constant.}

By \eqref{exp map est}
\begin{equation}\label{noncon triangle eqn1} 
	|{\rm exp}_Z(X) - (Z+X)| \leq C \mathbf{m}_0^{1/2} |X|^2 
\end{equation}
for all $X \in B_1(0,P_Z)$, where $C = C(n,m,q,M_0,N_0,C_e,C_h) \in (0,\infty)$ is a constant.  Thus 
\begin{equation}\label{noncon triangle eqn2}
	|{\rm exp}_Z^{-1}(X) - (X-Z)| \leq C \mathbf{m}_0^{1/2} |X-Z|^2 
\end{equation}
for all $X \in B^{\mathcal M}_1(Z)$, where $C = C(n,m,q,M_0,N_0,C_e,C_h) \in (0,\infty)$ is a constant.  By similar reasoning, \eqref{noncon triangle eqn1} and \eqref{noncon triangle eqn2} hold true with $Z_0$ in place of $Z$.  By \eqref{noncon rotate} and \eqref{noncon triangle eqn2}
\begin{equation}\label{noncon triangle eqn3}
	|\Gamma_Z^{-1} \exp_Z^{-1}(X) - \exp_{Z_0}^{-1}(X) + \exp_{Z_0}^{-1}(Z)| \leq C \mathbf{m}_0^{1/2} \rho_0^2 
\end{equation}
for all $X \in B^{\mathcal M}_{\rho_0}(Z_0)$, where $C = C(n,m,q,M_0,N_0,C_e,C_h) \in (0,\infty)$ are constants.  Let $\varphi_{j,l}$ is as in \eqref{varphi_localized} with $B = B_{r(X)/2}(\exp_{Z_0}^{-1}(X),P_{Z_0})$ and let $\widehat{\varphi}_{j,l} = \Gamma_Z \circ \varphi_{j,l} \circ \Gamma_Z^{-1}$ on $B$.  Recall that $\exp_{Z_0}^{-1}(Z) = \xi+\zeta$ where $\xi \in P_{Z_0} \cap S^{\perp}$ and $\zeta \in S$.  By applying Taylor's theorem to $\varphi_{j,l}$ at the point $\exp_{Z_0}^{-1}(X)$ and using Remark~\ref{graphical rmk}(2), \eqref{noncon rotate}, and \eqref{noncon triangle eqn3}  
\begin{equation}\label{noncon triangle eqn4} 
	\widehat{\varphi}_{j,l}(\exp_Z^{-1}(X)) 
	= \varphi_{j,l}(\exp_{Z_0}^{-1}(X)) - D\varphi_{j,l}(\exp_{Z_0}^{-1}(X)) \cdot \xi + \mathcal{R}_{j,l}(X) ,  
\end{equation}
where $|\mathcal{R}_{j,l}(X)| \leq C \widehat{E} \rho_0^{-1} |\xi|^2 + C \mathbf{m}_0^{1/2} \widehat{E} \rho_0^2$ for some constant $C = C(n,m,q,\alpha,M_0,N_0,C_e,C_h,\alpha) \in (0,\infty)$.  It now follows from \eqref{noncon triangle eqn4} that \eqref{noncon triangle1} holds true.

\noindent (2)  \textit{For each $\tau \in (0,1/2)$ there exist $\varepsilon_0,\beta \in (0,1)$ depending only on $n,m,q,M_0,N_0,C_e,C_h,\tau$ such that if Hypothesis~$(\star\star)$ holds true and if $X \in B^{\mathcal M}_{\rho_0}(Z_0)$ with $|\xi| < r(X)/2$ and $r(X) \geq \tau/2$, then }
\begin{align}\label{noncon triangle2} 
	|D\varphi(\exp_{Z_0}^{-1}(X)) \cdot \xi| \leq\,& \mathcal{G}(N(X),\varphi(\exp_{Z_0}^{-1}(X))) + \mathcal{G}(N(X),\widehat{\varphi}(\exp_Z^{-1}(X)) 
		\\&+ C \widehat{E} \rho_0^{-1} |\xi|^2 + C \mathbf{m}_0^{1/2} \widehat{E} \rho_0^2 , \nonumber 
\end{align}
\textit{where $C = C(n,m,q,M_0,N_0,C_e,C_h,\alpha) \in (0,\infty)$ is a constant.}

Let $\delta \in (0,1)$.  Assuming that $\varepsilon_0,\beta \in (0,1)$ are suitably small (depending on $n,m,q,M_0,N_0,$ $C_e,C_h,\alpha,\delta$), by Remark~\ref{graphical rmk}(2), \eqref{noncon5 concl}, and \eqref{graphical conclA1} 
\begin{align*}
	|D\varphi(\exp_{Z_0}^{-1}(X)) \cdot \xi| \leq\,& C \widehat{E} |\xi| \leq C \delta \inf_{\mathbb{S}^{n+m-1} \cap P \cap S^{\perp}} \op{sep} \varphi 
		\leq C \delta r(X)^{-\alpha} \op{sep} \varphi(\exp_{Z_0}^{-1}(X)) \\ \leq\,& C \delta \tau^{-\alpha} \op{sep} \varphi(\exp_{Z_0}^{-1}(X)) , 
\end{align*}
where $C = C(n,m,q,M_0,N_0,C_e,C_h) \in (0,\infty)$ are constants.  Hence choosing $\delta = \delta(n,m,q,M_0,N_0,$ $C_e,C_h,\alpha,\tau)$ sufficiently small, it follows that $\mathcal{G}(\varphi(X),\varphi(X) - D\varphi(X) \cdot \xi) = |D\varphi(X) \cdot \xi|$.  Hence by \eqref{noncon triangle eqn4}, it follows that \eqref{noncon triangle2} holds true.

\noindent (3)  By Remark~\ref{graphical rmk}(2) and the homogeneity of $\varphi$, $|\varphi| \leq C\widehat{E} \rho_0$ and $\op{Lip} \varphi \leq C \widehat{E}$ on $B_{\rho_0}(0,P_{Z_0})$.  Using this together with \eqref{noncon rotate} and \eqref{noncon triangle eqn3} 
\begin{align}\label{noncon triangle3}
	& |\mathcal{G}(N(X),\varphi(\exp_{Z_0}^{-1}(X))) - \mathcal{G}(N(X),\widehat{\varphi}(\exp_Z^{-1}(X)))|
	\\ \leq\,& \mathcal{G}(\varphi(\exp_{Z_0}^{-1}(X)), \widehat{\varphi}(\exp_Z^{-1}(X))) 
	\leq C \widehat{E} |\xi| + C \mathbf{m}_0^{1/2} \widehat{E} \rho_0^2 \nonumber 
\end{align}
for all $X \in B^{\mathcal M}_{\rho_0}(Z_0)$, where $C = C(n,m,q,M_0,N_0,C_e,C_h,\alpha) \in (0,\infty)$ is a constant.

\noindent (4)  \textit{For each $\delta \in (0,1)$ there exist $\varepsilon_0,\beta \in (0,1)$ depending only on $n,m,q,M_0,N_0,C_e,C_h,\delta$ such that  such that if Hypothesis~$(\star)$ and Hypothesis~$(\dagger)$ hold true, then } 
\begin{equation}\label{noncon triangle eqn5} 
	E_{\varsigma/2}(N,\widehat{\varphi},B^{\mathcal M}_{\rho_0/4}(Z)) \leq C (\beta + \delta) \,\rho_0^{-n/2-1} \|N\|_{L^2(B^{\mathcal M}_{\rho_0/4}(Z))} , 
\end{equation}
\textit{where $C = C(n,m,q,M_0,N_0,C_e,C_h,\alpha) \in (0,\infty)$ is a constant.}

By \eqref{noncon triangle3}, Hypothesis~$(\dagger)$, and Lemma~\ref{nonconcentration3 lemma} 
\begin{align}\label{noncon triangle eqn6} 
	E_{\varsigma/2}(N,\widehat{\varphi},B^{\mathcal M}_{\rho_0/4}(Z)) 
	\leq\,& C E_{\varsigma/2}(N,\varphi,B^{\mathcal M}_{\rho_0}(Z_0)) + C \widehat{E} |\xi|/\rho_0 
	\\ \leq\,& C (\beta + \delta) \,\widehat{E} 
	\leq C (\beta + \delta) \,\rho_0^{\alpha-1} \|\varphi\|_{L^2(B_1(0,P_{Z_0}))} \nonumber 
\end{align}
where $\widehat{E} = \rho_0^{-n/2-1} \|N\|_{L^2(B^{\mathcal M}_{\rho_0}(Z_0))}$ and $C = C(n,m,q,M_0,N_0,C_e,C_h,\alpha) \in (0,\infty)$ are constants.  On the other hand, by the triangle inequality and Hypothesis~$(\dagger)$ 
\begin{align}\label{noncon triangle eqn7} 
	\rho_0^{\alpha-1} \|\varphi\|_{L^2(B_1(0,P_{Z_0}))} 
	=\,& \rho_0^{\alpha-1} \|\widehat{\varphi}\|_{L^2(B_1(0,P_Z))}
	\\ \leq\,& C E_{\varsigma/2}(N,\widehat{\varphi},B^{\mathcal M}_{\rho_0/4}(Z)) 
		+ C \rho_0^{-n/2-1} \|N\|_{L^2(B^{\mathcal M}_{\rho_0/4}(Z))} , \nonumber
\end{align}
where $C = C(n,m,q,M_0,N_0,C_e,C_h,\alpha) \in (0,\infty)$ is a constant.  Combining \eqref{noncon triangle eqn6} and \eqref{noncon triangle eqn7}
\begin{equation}\label{noncon triangle eqn8} 
	E_{\varsigma/2}(N,\widehat{\varphi},B^{\mathcal M}_{\rho_0/4}(Z)) 
	\leq C (\beta_0 + \delta) \,\big(E_{\varsigma/2}(N,\widehat{\varphi},B^{\mathcal M}_{\rho_0/4}(Z)) 
		+ \rho_0^{-n/2-1} \|N\|_{L^2(B^{\mathcal M}_{\rho_0/4}(Z))}\big) ,  
\end{equation}
where $C = C(n,m,q,M_0,N_0,C_e,C_h,\alpha) \in (0,\infty)$ is a constant.  Assuming $\beta$ and $\delta$ are small enough that $C (\beta + \delta) < 1/2$ (where $C$ is as in \eqref{noncon triangle eqn8}), we deduce that \eqref{noncon triangle eqn5} holds true.
}\end{remark}

\begin{proof}[Proof of Theorem~\ref{nonconcentration thm}]  
Suppose that $p \in \{1,2,\ldots,\lfloor q/q_0 \rfloor\}$ such that either $p = 1$ or $p > 1$ and (H4) holds true.  Assume that for suitably small constants $\varepsilon_0,\beta \in (0,1)$ depending only on $n,m,q,p,M_0,N_0,C_e,C_h,\alpha,\sigma$, Hypothesis~$(\star)$ and Hypothesis~$(\dagger)$ hold true.  Let $Z \in \mathcal{B}^{(\geq 2)}_q \cap B^{\mathcal M}_{\rho_0/2}(Z_0)$ such that $\mathbf{I}_{N,Z}(0^+) \geq \alpha$ and let $Z = \exp_{Z_0}(\xi+\zeta)$ where $\xi \in P_{Z_0} \cap S^{\perp}$ and $\zeta \in S$.

\noindent\textit{Proof of Conclusion~(A).}  To prove conclusion~(A), it suffices to assume that for some $\beta_0 = \beta_0(n,m,q,p,M_0,N_0,C_e,C_h,\alpha) \in (0,1)$ suitably small, Hypothesis~$(\star\star)$ holds true.  If instead Hypothesis~$(\star\star)$ is false, there exists $\widetilde{p} \in \{0,1,2,\ldots,p-1\}$ and $\widetilde{\varphi} \in \Phi_{\alpha,q,\widetilde{p}}(P_{Z_0},S)$ such that 
\begin{gather}
	\label{keyest hyp rmk eqn1} E_{\varsigma/2}(N,\widetilde{\varphi},B^{\mathcal M}_{\rho_0}(Z_0)) 
		\leq \beta_0 \inf_{\varphi' \in \bigcup_{p'=0}^{\widetilde{p}-1} \Phi_{\alpha,q,p'}(P_{Z_0},S)} 
		E_{\varsigma/2}(N,\varphi',B^{\mathcal M}_{\rho_0}(Z_0)) , \\
	\label{keyest hyp rmk eqn2} E_{\varsigma/2}(N,\widetilde{\varphi},B^{\mathcal M}_{\rho_0}(Z_0)) 
		\leq 2^p \beta_0^{1-p} \inf_{\varphi' \in \bigcup_{p'=0}^{p-1} \Phi_{\alpha,q,p'}(P_{Z_0},S)} 
		E_{\varsigma/2}(N,\varphi',B^{\mathcal M}_{\rho_0}(Z_0)) .
\end{gather}
By Hypothesis~$(\dagger)$ and \eqref{keyest hyp rmk eqn2} 
\begin{align}\label{keyest hyp rmk eqn3}
	E_{\varsigma/2}(N,\widetilde{\varphi},B^{\mathcal M}_{\rho_0}(Z_0)) 
	\leq 2^{p_0} \beta_0^{-p_0} E_{\varsigma/2}(N,\varphi,B^{\mathcal M}_{\rho_0}(Z_0))
	\leq 2^{p_0} \beta_0^{-p_0} \beta \rho_0^{-n/2-1} \|N\|_{L^2(B^{\mathcal M}_{\rho_0}(Z_0))} .  
\end{align}
Since $\widetilde{p} = 0$ precisely when $\widetilde{\varphi} = q \llbracket 0 \rrbracket$, it follows from \eqref{keyest hyp rmk eqn3} that $\widetilde{p} > 0$.  In light of \eqref{keyest hyp rmk eqn1}, \eqref{keyest concl} holds true with $\widetilde{\varphi}$ in place of $\varphi$.  Noting that the left-hand side of \eqref{keyest concl} is independent of $\varphi$, it follows from \eqref{keyest hyp rmk eqn2} that \eqref{keyest concl} holds true with the original $\varphi$.

We claim that:  \textit{for every $\vartheta \in (0,1/8]$ there exists $\varepsilon_0, \beta_0 \in (0,1)$ depending on $n,m,q,p,M_0,$ $N_0,C_e,C_h,\alpha,\vartheta$ such that if Hypothesis~$(\star)$ and Hypothesis~$(\star\star)$ hold true and if $Z \in \mathcal{B}^{(\geq 2)}_q \cap B^{\mathcal M}_{\rho_0/2}(Z_0)$ with $\mathbf{I}_{N,Z}(0^+) \geq \alpha$ and $Z = \exp_{Z_0}(\xi+\zeta)$ where $\xi \in P_{Z_0} \cap S^{\perp}$ and $\zeta \in S$, then 
\begin{equation}\label{nonconA eqn1}
	c \widehat{E}^2 |\xi|^2 (\vartheta\rho_0)^n \leq \int_{B^{\mathcal M}_{\vartheta\rho_0/4}(Z) \cap \{r > \vartheta\rho_0/32\}} 
		|D\varphi(\exp_{Z_0}^{-1}(X)) \cdot \xi|^2 \,d\mathcal{H}^n(X) ,
\end{equation}
where $c = c(n,m,q,p,M_0,N_0,C_e,C_h,\alpha) > 0$ is a constant (independent of $\rho$).}  To see this, suppose to the contrary that for some $\vartheta \in (0,1/8]$ and for $\nu = 1,2,3,\ldots$ there exists $\varepsilon_{\nu} \rightarrow 0^+$, $\beta_{\nu} \rightarrow 0^+$, $T^{(\nu)}$, $\mathcal{M}^{(\nu)}$, $N^{(\nu)}$, $\mathcal{K}_{\nu}$, $Z_{0,\nu}$, $\rho_{\nu}$, $P_{\nu}$, $S_{\nu}$, and $\varphi^{(\nu)}$ satisfying Hypothesis~$(\star)$ and Hypothesis~$(\star\star)$ in place of $\varepsilon_0$, $\beta_0$, $T$, $\mathcal{M}$, $N$, $\mathcal{K}$, $Z_0$, $\rho_0$, $P_{Z_0}$, $S$, and $\varphi$ and there exists $Z_{\nu} \in \mathcal{B}^{(\geq 2)}_q(T^{(\nu)}) \cap B^{\mathcal{M}^{(\nu)}}_{\rho_{\nu}/2}(Z_{0,\nu})$ with $\mathbf{I}_{N^{(\nu)},Z_{\nu}}(0^+) \geq \alpha$ given by $Z_{\nu} = \mathbf{e}_{\nu}(\xi_{\nu} + \zeta_{\nu})$ for unique $\xi_{\nu} \in P_{\nu} \cap S_{\nu}^{\perp}$ and $\zeta_{\nu} \in S_{\nu}$ such that  
\begin{equation}\label{nonconA eqn2} 
	\int_{B^{\mathcal{M}^{(\nu)}}_{\vartheta\rho_{\nu}/4}(Z_{\nu}) \cap \{r > \vartheta\rho_{\nu}/32\}} 
		|D\varphi^{(\nu)}(\mathbf{e}_{\nu}^{-1}(X)) \cdot \xi_{\nu}|^2 d\mathcal{H}^n(X)
	\leq \frac{1}{\nu} \widehat{E}_{\nu}^2 |\xi_{\nu}|^2 (\vartheta\rho_{\nu})^n ,
\end{equation}
where $\mathbf{e}_{\nu}$ is the exponential map of $\mathcal{M}^{(\nu)}$ at $Z_{0,\nu}$ and $\widehat{E}_{\nu} = \rho_{\nu}^{-n/2-1} \|N^{(\nu)}\|_{L^2(B^{\mathcal{M}^{(\nu)}}_{\rho_{\nu}}(Z_{0,\nu}))}$.  Let $\mathfrak{q}_{\nu}$ be a rotation of $\mathbb{R}^{n+m}$ such that $\mathfrak{q}_{\nu}(P_{\nu}) = \mathbb{R}^n \times \{0\}$ and $\mathfrak{q}_{\nu}(S_{\nu}) = \{0\}^2 \times \mathbb{R}^{n-2} \times \{0\}^m$.  Rotate and rescale letting $\overline{T}^{(\nu)}$, $\overline{\mathcal M}^{(\nu)}$, and $\overline{N}^{(\nu)}$ be as in \eqref{centered blowup eqn1}, $\overline{\varphi}^{(\nu)}(X) = \rho_{\nu}^{\alpha-1} (\mathfrak{q}_{\nu} \circ \varphi^{(\nu)} \circ \mathfrak{q}_{\nu}^{-1})(X)$ for each $X \in \mathbb{R}^n$, and $(\overline{\xi}_{\nu},\overline{\zeta}_{\nu}) = \mathfrak{q}_{\nu}(\xi_{\nu} + \zeta_{\nu})/\rho_{\nu}$ where $\overline{\xi}_{\nu} \in \mathbb{R}^2$ and $\overline{\zeta}_{\nu} \in \mathbb{R}^{n-2}$.  By Remark~\ref{graphical rmk}(2), after passing to a subsequence, there is a non-zero function $\varphi^{(0)} \in \Phi_{\alpha,q}$ such that $\overline{\varphi}^{(\nu)}/\widehat{E}_{\nu} \rightarrow \varphi^{(0)}$ uniformly on compact subsets of $\mathbb{R}^n$ as $\nu \rightarrow \infty$.  By Lemma~\ref{nonconcentration3 lemma}, $\overline{\xi}_{\nu} \rightarrow 0$ in $\mathbb{R}^2$.  After passing to a subsequence, let $\overline{\xi}_{\nu}/|\overline{\xi}_{\nu}| \rightarrow \widehat{\xi}$ in $\mathbb{R}^2$ and $\overline{\zeta}_{\nu} \rightarrow \overline{\zeta}$ in $\mathbb{R}^{n-2}$ as $\nu \rightarrow \infty$.  By rotating and rescaling, dividing both sides of \eqref{nonconA eqn2} by $\widehat{E}_{\nu}^2 |\xi_{\nu}|^2$, and letting $\nu \rightarrow \infty$  
\begin{equation*}
	D\varphi^{(0)}_{j,l}(X) \cdot \widehat{\xi} = 0 \text{ for all } X \in B_{\vartheta/4}(0,\zeta) \cap \{r > \vartheta/32\} .
\end{equation*}
Hence $\varphi^{(0)}$ is translation invariant in the $\widehat{\xi}$ direction.  This contradicts $\varphi^{(0)}$ being a non-zero function in $\Phi_{\alpha,q}$. 

Now let $\vartheta \in (0,1/8]$ to be later determined.  Recalling that Hypothesis~$(\star\star)$ holds true, we can use \eqref{nonconA eqn1} and \eqref{noncon triangle2} with $\tau = \vartheta\rho_0/32$ to obtain 
\begin{align}\label{nonconA eqn3}
	\frac{\widehat{E}^2 |\xi|^2 \vartheta^n}{\rho_0^2} 
	\leq\,& C E_{\varsigma/2}(N,\varphi,B^{\mathcal M}_{\rho_0}(Z_0))^2 
		+ \frac{C}{\rho_0^{n+2}} \int_{B^{\mathcal M}_{\vartheta\rho_0}(Z)} \mathcal{G}(N,\widehat{\varphi} \circ \exp_Z^{-1})^2 
		\\&+ \frac{C \widehat{E}^2 |\xi|^4 \vartheta^n}{\rho_0^4} + \frac{C}{\rho_0^{n+2}} \int_{B^{\mathcal M}_{\vartheta\rho_0/4}(Z) \cap \{r \leq 2|\xi|\}} 
		|D\varphi(\exp_{Z_0}^{-1}(X)) \cdot \xi|^2 \,d\mathcal{H}^n(X) , \nonumber
\end{align}
where $C = C(n,m,q,p,M_0,N_0,C_e,C_h,\alpha) \in (0,\infty)$ is a constant.  To bound the second term on the right-hand side of \eqref{nonconA eqn3}, recall Remark~\ref{noncon triangle rmk}(4) and use Corollary~\ref{radial noncon cor} and \eqref{noncon triangle3} to obtain 
\begin{align*}
	\rho_0^{2\alpha-2-1/q} \int_{B^{\mathcal M}_{\rho_0/8}(Z)} \frac{\mathcal{G}(N,\widehat{\varphi} \circ \exp_Z^{-1})^2}{R_Z^{n+2\alpha-1/q}} 
	\leq\,& C E_{\varsigma/2}(N,\widehat{\varphi},B^{\mathcal M}_{\rho_0/4}(Z))^2 
	\\ \leq\,& C E_{\varsigma/2}(N,\varphi,B^{\mathcal M}_{\rho_0}(Z_0))^2 + \frac{C \widehat{E}^2 |\xi|^2}{\rho_0^2}
\end{align*}
where $C = C(n,m,q,p,M_0,N_0,C_e,C_h,\alpha) \in (0,\infty)$ are constants.  In particular, 
\begin{equation}\label{nonconA eqn5}
	\frac{1}{\rho_0^{n+2}} \int_{B^{\mathcal M}_{\vartheta\rho_0}(Z)} \mathcal{G}(N,\widehat{\varphi} \circ \exp_Z^{-1})^2 
	\leq C E_{\varsigma/2}(N,\varphi,B^{\mathcal M}_{\rho_0}(Z_0))^2 + \frac{C \vartheta^{n+2\alpha-1/q} \widehat{E}^2 |\xi|^2}{\rho_0^2}
\end{equation}
where $C = C(n,m,q,p,M_0,N_0,C_e,C_h,\alpha) \in (0,\infty)$ are constants.  To bound the last term on the right-hand side of \eqref{nonconA eqn3}, using Remark~\ref{graphical rmk}(2) and the fact that $\varphi$ is homogeneous degree $\alpha$ and translation invariant along $S$ 
\begin{align}\label{nonconA eqn6}
	\int_{B^{\mathcal M}_{\vartheta\rho_0/4}(Z) \cap \{r \leq 2|\xi|\}} |D\varphi(\exp_{Z_0}^{-1}(X)) \cdot \xi|^2 \,d\mathcal{H}^n(X)
	\leq\,& C \widehat{E}^2 \vartheta^{n-2} \rho_0^{n-2\alpha} \int_{B^2_{2|\xi|}(0)} |x|^{2\alpha-2} |\xi|^2 \,dx 
	\\ \leq\,& C \widehat{E}^2 |\xi|^{2\alpha+2} \vartheta^{n-2} \rho_0^{n-2\alpha} \nonumber 
\end{align}
where $C = C(n,m,q,p,M_0,N_0,C_e,C_h,\alpha) \in (0,\infty)$ are constants.  Thus bounding terms on the right-hand side of \eqref{nonconA eqn3} using \eqref{nonconA eqn5}, \eqref{nonconA eqn6}, and the definition of $E_{\varsigma/2}(N,\varphi,B^{\mathcal M}_{\rho_0}(Z_0))$ 
\begin{equation}\label{nonconA eqn7}
	\frac{\widehat{E}^2 |\xi|^2 \vartheta^n}{\rho_0^2} 
	\leq C E_{\varsigma/2}(N,\varphi,B^{\mathcal M}_{\rho_0}(Z_0))^2 
		+ C \frac{\widehat{E}^2 |\xi|^2 \vartheta^n}{\rho_0^2} \left( \vartheta^{2\alpha-1/q} + \frac{|\xi|^2}{\rho_0^2} 
		+ \frac{|\xi|^{2\alpha}}{\vartheta^2 \rho_0^{2\alpha}} \right)  
\end{equation}
where $C = C(n,m,q,p,M_0,N_0,C_e,C_h,\alpha) \in (0,\infty)$ is a constant.  Let $\delta \in (0,1)$.  Assuming that $\varepsilon_0$ and $\beta_0$ are sufficiently small (depending on $\delta$), by Lemma~\ref{nonconcentration3 lemma} we have that $|\xi| \leq \delta\rho_0$.  Thus choosing $\vartheta = \vartheta(n,m,q,p,M_0,N_0,C_e,C_h,\alpha)$ and $\delta = \delta(n,m,q,p,M_0,N_0,C_e,C_h,\alpha,\vartheta) \in (0,\rho/2]$ small enough that $C (\vartheta^{2-1/q} + \delta^2 + \delta^{2\alpha} \vartheta^{-2}) < 1/2$ (where $C$ is as in \eqref{nonconA eqn7}), we deduce that conclusion~(A) holds true.

\noindent\textit{Proof of Conclusion~(B).}  Conclusion~(B) follows directly from \eqref{noncon triangle3} and conclusion~(A).

\noindent\textit{Proof of Conclusion~(C).}  By \eqref{noncon triangle3}, 
\begin{align}\label{nonconB eqn1}
	&\frac{1}{\rho_0^{2+\sigma}} \int_{B^{\mathcal M}_{3\rho_0/4}(Z_0)} \frac{\mathcal{G}(N,\varphi\circ\exp_{Z_0}^{-1})^2}{R_Z^{n-\sigma}} 
	\leq \frac{8^{n-\sigma}}{\rho_0^{n+2}} \int_{B^{\mathcal M}_{3\rho_0/4}(Z_0)} \mathcal{G}(N,\varphi\circ\exp_{Z_0}^{-1})^2
	\\&\hspace{15mm} + \frac{2}{\rho_0^{2+\sigma}} \int_{B^{\mathcal M}_{\rho_0/8}(Z)} \frac{\mathcal{G}(N,\widehat{\varphi}\circ\exp_Z^{-1})^2}{
	R_Z^{n-\sigma}} + \frac{C \widehat{E}^2}{\rho_0^{2+\sigma}} \int_{B^{\mathcal M}_{3\rho_0/4}(Z_0)} \frac{|\xi|^2 + \mathbf{m}_0 \rho_0^4}{
	R_Z^{n-\sigma}} \nonumber
\end{align}
where $C = C(n,m,q,p,M_0,N_0,C_e,C_h,\alpha,\sigma) \in (0,\infty)$ is a constant.  To bound the second term on the right-hand side of \eqref{nonconB eqn1}, recall Remark~\ref{noncon triangle rmk}(4) and use Corollary~\ref{radial noncon cor} and conclusion~(B) to obtain 
\begin{align}\label{nonconB eqn2}
	\rho_0^{2\alpha-2-\sigma} \int_{B^{\mathcal M}_{\rho_0/8}(Z)} \frac{\mathcal{G}(N,\widehat{\varphi} \circ \exp_Z^{-1})^2}{R_Z^{n+2\alpha-\sigma}} 
	\leq\,& C E_{\varsigma/2}(N,\widehat{\varphi},B^{\mathcal M}_{\rho_0/4}(Z))^2 
	\\ \leq\,& C E_{\varsigma/2}(N,\varphi,B^{\mathcal M}_{\rho_0}(Z_0))^2 \nonumber 
\end{align}
where $C = C(n,m,q,p,M_0,N_0,C_e,C_h,\alpha,\sigma) \in (0,\infty)$ are constants.  To bound the third term on the right-hand side of \eqref{nonconB eqn1}, we have that $\int_{B^{\mathcal M}_{3\rho_0/4}(Z_0)} R_Z^{-n+\sigma} \leq C \rho_0^{\sigma}$ and thus by conclusion~(A) and the definition of $E_{\varsigma/2}(N,\varphi,B^{\mathcal M}_{\rho_0}(Z_0))$ 
\begin{align}\label{nonconB eqn3}
	\frac{\widehat{E}^2}{\rho_0^{2+\sigma}} \int_{B^{\mathcal M}_{3\rho_0/4}(Z_0)} \frac{|\xi|^2 + \mathbf{m}_0 \rho_0^4}{R_Z^{n-\sigma}} 
	\leq C \widehat{E}^2 \left( \frac{|\xi|^2}{\rho_0^2} + \mathbf{m}_0 \rho_0^2 \right) 
	\leq C E_{\varsigma/2}(N,\varphi,B^{\mathcal M}_{\rho_0}(Z_0))^2 , 
\end{align}	
where $C = C(n,m,q,p,M_0,N_0,C_e,C_h,\alpha,\sigma) \in (0,\infty)$ are constants.  Combining \eqref{nonconB eqn1}, \eqref{nonconB eqn2}, and \eqref{nonconB eqn3} gives us conclusion~(C).

\noindent\textit{Proof of Conclusion~(D).}  By \eqref{noncon triangle eqn4} 
\begin{align}\label{nonconD eqn1}
	&\int_{B^{\mathcal M}_{3\rho_0/4}(Z_0) \cap \{r > \tau\rho_0\}} \frac{\mathcal{G}(N(X), 
		\varphi(\exp_{Z_0}^{-1}(X)) - D\varphi(\exp_{Z_0}^{-1}(X)) \cdot \xi)^2}{R_Z^{n+2\alpha-\sigma}} \,d\mathcal{H}^n(X)
	\\ \leq\,& \int_{B^{\mathcal M}_{\rho_0/8}(Z)} \frac{\mathcal{G}(N,\widehat{\varphi} \circ \exp_Z^{-1})^2}{
		R_Z^{n+2\alpha-\sigma}} 
	+ C \widehat{E}^2 \int_{B^{\mathcal M}_{3\rho_0/4}(Z_0) \cap \{r > \tau\rho_0\}} \frac{\rho_0^{-2} |\xi|^4 + \mathbf{m}_0 \rho_0^4}{
		R_Z^{n+2\alpha-\sigma}} \,d\mathcal{H}^n(X) \nonumber 
	\\&+ \frac{C}{\rho_0^{n+2\alpha-\sigma}} \int_{B^{\mathcal M}_{3\rho_0/4}(Z_0) \setminus B^{\mathcal M}_{\rho_0/8}(Z)} 
		(\mathcal{G}(N,\varphi)^2 + \widehat{E}^2 |\xi|^2) , \nonumber
\end{align}
where $C = C(n,m,q,\alpha,\sigma) \in (0,\infty)$ are constants.  To bound the first term on the right-hand side of \eqref{nonconD eqn1}, we use \eqref{nonconB eqn2}.  To bound the second term on the right-hand side of \eqref{nonconD eqn1}, let $\delta \in (0,1)$.  Assuming that $\varepsilon_0$ and $\beta$ are sufficiently small (depending on $\delta$), by Lemma~\ref{nonconcentration3 lemma} we have that $|\xi| \leq \delta\rho_0$.  We have that $\int_{B^{\mathcal M}_{3\rho_0/4}(Z_0) \cap \{r > \tau\rho_0\}} R_Z^{-n-2\alpha+\sigma} \leq C (\tau\rho_0)^{-2\alpha+\sigma}$ and thus by conclusion~(A) and the definition of $E_{\varsigma/2}(N,\varphi,B^{\mathcal M}_{\rho_0}(Z_0))$ 
\begin{align*}
	\widehat{E}^2 \int_{B^{\mathcal M}_{3\rho_0/4}(Z_0) \cap \{r > \tau\rho_0\}} \frac{\rho_0^{-2} |\xi|^4 + \mathbf{m}_0 \rho_0^4}{
		R_Z^{n+2\alpha-\sigma}} \,d\mathcal{H}^n(X)
	\leq C \widehat{E}^2 (\delta^2 |\xi|^2 + \mathbf{m}_0 \rho_0^4) \,(\tau\rho_0)^{-2\alpha+\sigma} , 
\end{align*}
where $C = C(n,m,q,\alpha,\sigma) \in (0,\infty)$ is a constant.  Assuming that $\delta$ and $\varepsilon_0$ are small enough that $(\delta^2 + \mathbf{m}_0 \rho_0^2) \,\tau^{-2\alpha+\sigma} < 1$ 
\begin{align}\label{nonconD eqn2}
	\widehat{E}^2 \int_{B^{\mathcal M}_{3\rho_0/4}(Z_0) \cap \{r > \tau\rho_0\}} \frac{\rho_0^{-2} |\xi|^4 + \mathbf{m}_0 \rho_0^4}{
		R_Z^{n+2\alpha-\sigma}} \,d\mathcal{H}^n(X)
	\leq C \rho_0^{-2\alpha+2+\sigma} E_{\varsigma/2}(N,\varphi,B^{\mathcal M}_{\rho_0}(Z_0))^2 , 
\end{align}
where $C = C(n,m,q,\alpha,\sigma) \in (0,\infty)$ is a constant.  The last term on the right-hand side of \eqref{nonconD eqn1} is readily bounded using conclusion~(A) and the definition of $E_{\varsigma/2}(N,\varphi,B^{\mathcal M}_{\rho_0}(Z_0))$.  By combining \eqref{nonconD eqn1}, \eqref{nonconD eqn2}, and conclusion~(A) we obtain conclusion~(D).
\end{proof}

\begin{corollary}\label{nonconcentration cor} 
Let $\alpha = \ell_0/q_0 \geq 2$ for some relatively prime positive integers $\ell_0$ and $q_0 \leq q$.  Let $p \in \{1,2,\ldots, \lfloor q/q_0 \rfloor \}$.  Let $M_0,N_0,C_e,C_h,\varepsilon_2$ be as in Hypothesis~\ref{noncon constants hypothesis}.  For every $\delta \in (0,1/4)$ and $\sigma \in (0,1)$ there exists $\varepsilon_0, \beta \in (0,1)$ depending only on $n,m,q,p,M_0,N_0,C_e,C_h,\alpha,\delta,\sigma$ such that if Hypothesis~$(\star)$ (of Section~\ref{sec:graph_rep_results}) and Hypothesis~$(\dagger)$ (of Section~\ref{sec:key_est_sec}) hold, and if 
\begin{equation*}
	B^{\mathcal M}_{\delta\rho_0}(\exp_{Z_0}(\zeta)) \cap B^{\mathcal M}_{\rho_0/2}(Z_0) 
		\cap \{ X \in \mathcal{B}^{(\geq 2)}_q : \mathbf{I}_{N,X}(0^+) \geq \alpha \} \neq \emptyset 
\end{equation*}
for all $\zeta \in S \cap B_{\rho_0/2}(0,P_{Z_0})$, then 
\begin{equation}\label{noncon cor concl}
	\frac{1}{\rho_0^{n+\sigma}} \int_{B^{\mathcal M}_{\rho_0/2}(Z_0)} \frac{\mathcal{G}(N,\varphi)^2}{r_{\delta}^{2-\sigma}} 
	\leq C E_{\varsigma/2}(N,\varphi,B^{\mathcal M}_{\rho_0}(Z_0))^2 , 
\end{equation}
where $r(X) = \op{dist}(\exp_{Z_0}^{-1}(X),S)$ (as in \eqref{normal coord eqn2}), $r_{\delta} = \max\{r,\delta\}$, and $C = C(n,m,q,p,M_0,N_0,C_e,C_h,$ $\alpha,\sigma) \in (0,\infty)$ is a constant (independent of $\delta$).
\end{corollary}

\begin{proof}
This follows from Theorem~\ref{nonconcentration thm}(C) by arguing as in the proof of~\cite[Corollary~3.2]{Sim93}.
\end{proof}

\section{Fine blow-up class}\label{sec:fine blowups}

\subsection{Definition of the fine blow-up class and notation}  Throughout this section, let $\alpha = \ell_0/q_0 \geq 2$ for some relatively prime positive integers $\ell_0$ and $q_0 \leq q$.  Let $M_0,N_0,C_e,C_h,\varepsilon_2$ be as in Hypothesis~\ref{noncon constants hypothesis}.  Let $\varsigma = \varsigma(n,m,q) \in (0,1)$ is as in Theorem~\ref{relharm separation3 thm}.  For $\nu = 1,2,3,\ldots$, let $\varepsilon_{\nu} \rightarrow 0^+$.  In Definition~\ref{blowup class defn} below, we define $\mathfrak{B}$ to be the class of all fine blow-ups $w = (w_j)$ of a sequence of $n$-dimensional locally area minimizing rectifiable currents $T^{(\nu)}$ relative to a sequence of homogeneous harmonic $q$-valued functions $\varphi^{(\nu)} \in \Phi_{\alpha,q}(P_{\nu},S_{\nu})$ as constructed in Subsection~\ref{sec:blowup subsec}.  The blow-up class $\mathfrak{B}$ and the proof of its fundamental properties is based on~\cite[Sections~9--11]{KrumWic2} and~\cite[Sections~9--10]{KrumWicc} and we provide details for completion.

\begin{definition}\label{domain class defn}
$\mathfrak{D}$ is the set of all ordered collections $(m_j,\varphi^{(0)}_j,\varphi^{(\infty)}_j)_{1 \leq j \leq J}$ where: 
\begin{enumerate}[itemsep=3mm,topsep=0mm]
	\item[{\rm (I)}]  $J$, $m_1,m_2,\ldots,m_J$, and $q_1,q_2,\ldots,q_J$ are positive integers such that $\sum_{j=1}^J m_j q_j = q$; 
	\item[{\rm (II)}]  for each $j \in \{1,2,\ldots,J\}$, $\varphi^{(0)}_j : \mathbb{R}^n \rightarrow \mathcal{A}_{q_j}(\mathbb{R}^m)$ and $\varphi^{(0)}_j : \mathbb{R}^n \rightarrow \mathcal{A}_{q_j}(\mathbb{R}^m)$ are $q_j$-valued functions such that either: 
	\begin{enumerate}[itemsep=3mm,topsep=3mm]
		\item[{\rm (a)}]  $q_j = q_0$ and $\varphi^{(0)}$ and $\varphi^{(\infty)}$ are $q_0$-valued functions defined by 
		\begin{equation*}
			\varphi^{(0)}_j(x,y) = \op{Re}(c^{(0)}_j (x_1+ix_2)^{\alpha}), \quad \varphi^{(\infty)}_j(x,y) = \op{Re}(c^{(0)}_j (x_1+ix_2)^{\alpha})
		\end{equation*}
		for all $x = (x_1,x_2) \in \mathbb{R}^2$ and $y \in \mathbb{R}^{n-2}$, where $c^{(0)}_j, c^{(\infty)}_j \in \mathbb{C}^m$ with $|c^{(\infty)}_j| = 1$ and $c^{(0)}_j = \gamma_j c^{(\infty)}_j$ for some $\gamma_j \in [0,\infty)$;
		
		\item[{\rm (b)}]  $q_j = 1$ and $\varphi^{(0)}_j : \mathbb{R}^n \rightarrow \mathcal{A}_{q_0}(\mathbb{R}^m)$ and $\varphi^{(\infty)}_j : \mathbb{R}^n \rightarrow \mathcal{A}_{q_0}(\mathbb{R}^m)$ are single-valued functions defined by $\varphi^{(0)}_j(x,y) = \varphi^{(\infty)}_j(x,y) = 0$ for all $(x,y) \in \mathbb{R}^n$.
	\end{enumerate}
\end{enumerate}
\end{definition}

\begin{definition}\label{blowup class defn} 
$\mathfrak{B}$ is the set of all ordered collections $w = (w_j)_{1 \leq j \leq J}$ of functions 
$$w_j : \op{graph}\varphi^{(\infty)}_j |_{B_{3/4}(0) \setminus \{0\} \times \mathbb{R}^{n-2}} \rightarrow \mathcal{A}_{m_j}(\mathbb{R}^m)$$
such that $w = (w_j)$ is the fine blow-up of a sequence $(T^{(\nu)})_{\nu=1}^{\infty}$ relative to a sequence $(\varphi^{(\nu)})_{\nu=1}^{\infty}$ in $B^{\mathcal{M}^{(\nu)}}_{\rho_{\nu}}(Z_{\nu})$ by the excess $E_{\nu} = E_{\varsigma/2,\varepsilon_{\nu}}(N^{(\nu)},\varphi^{(\nu)},B^{\mathcal{M}^{(\nu)}}_{\rho_{\nu}}(Z_{\nu}))$ (as in~\eqref{graphical E defn}) as constructed in Subsection~\ref{sec:blowup subsec}, where for some sequence $(\varepsilon_{\nu})_{\nu=1}^{\infty}$, $(\beta_{\nu})_{\nu=1}^{\infty}$, and $(\delta_{\nu})_{\nu=1}^{\infty}$ tending to zero: 
\begin{enumerate}[itemsep=3mm,topsep=0mm]
	\item[(A)]  $T^{(\nu)}$ is an $n$-dimensional locally area minimizing rectifaible current of $\mathbf{C}_{11\sqrt{n}/2}(0)$, $\mathbf{m}_{0,\nu} > 0$, $\mathcal{M}^{(\nu)}$ is the center manifold of $T^{(\nu)}$, $N^{(\nu)} : \mathcal{M}^{(\nu)} \rightarrow ({\rm T}\mathcal{M}^{(\nu)})^{\perp}$ is the normal approximation map of $T^{(\nu)}$, and $\mathcal{K}_{\nu} \subseteq \mathcal{M}^{(\nu)}$ is a closed set such that Hypothesis~\ref{center mnfd hypothesis} and Hypothesis~\ref{M and N hypothesis} hold true with $T^{(\nu)}, \mathbf{m}_{0,\nu}, \mathcal{M}^{(\nu)}, N^{(\nu)}, \mathcal{K}_{\nu}$ in place of $T, \mathbf{m}_0, \mathcal{M},N,\mathcal{K}$; 
	
	\item[(B)]  	$Z_{\nu} \in \mathcal{B}^{(\geq 2)}_q(T^{(\nu)})$ with $\mathbf{I}_{N^{(\nu)},Z_{\nu}}(0^+) \geq \alpha$, $\rho_{\nu} \in (0,1]$ with $\mathbf{m}_{0,\nu} \rho_{\nu} < \varepsilon_{\nu}$, $P_{\nu} = {\rm T}_{Z_{\nu}} \mathcal{M}^{(\nu)}$, and $S_{\nu}$ is an $(n-2)$-dimensional linear subspace of $P_{\nu}$; 
	
	\item[(C)]  $p \in \{1,2,\ldots,\lfloor q/q_0 \rceil\}$ (independent of $\nu$) and $\varphi^{(\nu)} \in \Phi_{\alpha,q,p}(P_{\nu},S_{\nu})$; 
	
	\item[(D)]  	following the procedure from Subsection~\ref{sec:blowup subsec}, we apply a rotation $\mathfrak{q}_{\nu}$ is a rotation of $\mathbb{R}^{n+m}$ such that $\mathfrak{q}_{\nu}(P_{\nu}) = \mathbb{R}^n \times \{0\}$ and $\mathfrak{q}_{\nu}(S_{\nu}) = \{0\}^2 \times \mathbb{R}^{n-2} \times \{0\}^m$; 
		
	\item[(E)]  
	\begin{equation}\label{blowups hypothesis2}
		E_{\varsigma/2,\varepsilon_{\nu}}(N^{(\nu)},\varphi^{(\nu)},B^{\mathcal{M}^{(\nu)}}_{\rho_{\nu}}(Z_{\nu})) 
		\leq \beta_{\nu} \inf_{\varphi' \in \bigcup_{p'=0}^{p-1} \Phi_{\alpha,q,p'}(P_{\nu},S_{\nu})} 
			E_{\varsigma/2,\varepsilon_{\nu}}(N^{(\nu)},\varphi',B^{\mathcal{M}^{(\nu)}}_{\rho_{\nu}}(Z_{\nu}));
	\end{equation}
	
	\item[(F)]  for every $\zeta \in S_{\nu} \cap B_{\rho_{\nu}/2}(0,P_{\nu})$ 
	\begin{equation}\label{blowups no gaps}
		B^{\mathcal{M}^{(\nu)}}_{\delta_{\nu}\rho_{\nu}}(\mathbf{e}_{\nu}(\zeta)) \cap B^{\mathcal{M}^{(\nu)}}_{\rho_{\nu}/2}(Z_{\nu}) \cap 
			\{ X \in \mathcal{B}^{(\geq 2)}_q(T^{(\nu)}) : \mathbf{I}_{N^{(\nu)},X}(0^+) \geq \alpha \} \neq \emptyset , 
	\end{equation}
	where $\mathbf{e}_{\nu} : B_3(0,P_{\nu}) \rightarrow \mathcal{M}^{(\nu)}$ is the exponential map of $\mathcal{M}^{(\nu)}$ at $Z_{\nu}$.
\end{enumerate}
\end{definition}

\begin{remark}\label{blowup on domain rmk}{\rm 
Let $w = (w_j) \in \mathfrak{B}$ be the blow-up of $T^{(\nu)}$ relative to a sequence $\varphi^{(\nu)}$ in $B^{\mathcal{M}^{(\nu)}}_{\rho_{\nu}}(Z_{\nu})$ by the excess $E_{\nu} = E_{\varsigma/2,\varepsilon_{\nu}}(N^{(\nu)},\varphi^{(\nu)},B^{\mathcal{M}^{(\nu)}}_{\rho_{\nu}}(Z_{\nu}))$ (as in Definition~\ref{blowup class defn}).  Express $\varphi^{(\nu)} = \sum_{j=1}^J m_j \varphi^{(\nu)}_j$ for some positive integers $J,m_1,m_2,\ldots,m_J$ and distinct functions $\varphi^{(\nu)}_j : \mathbb{R}^n \rightarrow \mathcal{A}_{q_j}(\mathbb{R}^m)$ (as in Definition~\ref{Phi-p_defn}).  As was discussed in Subsection~\ref{sec:blowup subsec}, after passing to a subsequence for each $w = (w_j) \in \mathfrak{B}$ there exists $(m_j,\varphi^{(0)}_j,\varphi^{(\infty)}_j) \in \mathfrak{D}$ such that:
\begin{enumerate}[itemsep=3mm,topsep=0mm]
	\item[(I)]  for each $j \in \{1,\ldots,J\}$, $m_j$ is the multiplicity of the component $\varphi^{(\nu)}_j$ in accordance with Definition~\ref{Phi-p_defn} (with $\varphi^{(\nu)}$ in place of $\varphi$); 
	
	\item[(II)]  for each $j \in \{1,\ldots,J\}$, precisely one of the following holds true: 
	\begin{enumerate}[itemsep=3mm,topsep=3mm]
		\item[(a)]  $q_j = q_0$, $\varphi^{(\nu)}_j$ are non-zero functions for all $\nu = 1,2,3,\ldots$, and $\varphi^{(\nu)}_j, \varphi^{(0)}_j, \varphi^{(\infty)}_j$ are $q_0$-valued functions given by 
		\begin{gather*}
			\rho_{\nu}^{\alpha-1} (\mathfrak{q}_{\nu} \circ \varphi^{(\nu)}_j \circ \mathfrak{q}_{\nu}^{-1})(x,y) = \op{Re}(c^{(\nu)}_j(x_1+ix_2)^{\alpha}), \\ 
			\varphi^{(0)}_j(x,y) = \op{Re}(c^{(0)}_j(x_1+ix_2)^{\alpha}), \quad 
			\varphi^{(\infty)}_j(x,y) = \op{Re}(c^{(\infty)}_j(x_1+ix_2)^{\alpha})
		\end{gather*}
		for all $x = (x_1,x_2) \in \mathbb{R}^2$ and $y \in \mathbb{R}^{n-2}$, where $c^{(\nu)} \in \mathbb{C}^m \setminus \{0\}$ and 
		\begin{equation*}
			c^{(0)}_j = \lim_{\nu\rightarrow\infty} \frac{c^{(\nu)}_j}{\widehat{E}_{\nu}}, \quad 
			c^{(\infty)}_j = \lim_{\nu\rightarrow\infty} \frac{c^{(\nu)}_j}{|c^{(\nu)}_j|} ,
		\end{equation*}
		where $\widehat{E}_{\nu} = \rho_{\nu}^{-n/2-1} \|N^{(\nu)}\|_{L^2(B^{\mathcal{M}^{(\nu)}}_{\rho_{\nu}}(Z_{\nu}))}$; 
		
		\item[(b)]  $q_j = 1$, $\varphi^{(\nu)}_j : \mathbb{R}^n \rightarrow \mathbb{R}^m$ are identically zero for all $\nu = 1,2,3,\ldots$, and $\varphi^{(0)}_j : \mathbb{R}^n \rightarrow \mathbb{R}^m$ and $\varphi^{(\infty)}_j : \mathbb{R}^n \rightarrow \mathbb{R}^m$ are also identically zero.
	\end{enumerate}
\end{enumerate}
By Remark~\ref{graphical rmk}(2), 
\begin{equation}\label{blowup on domain rmk eqn1}
	\frac{1}{C} \leq \sum_{j=1}^J \|\varphi^{(0)}_j\|_{L^2(B_1(0))} \leq C 
\end{equation}
for some constant $C = C(n,m,q,\alpha) \in (0,\infty)$.  Note that even in case~(a) where $\varphi^{(\nu)}_j$ are non-zero functions for all $\nu = 1,2,3,\ldots$, we might have that $\varphi^{(0)}_j$ is identically zero.  However, by \eqref{blowup on domain rmk eqn1}, we know that $\varphi^{(0)}_j$ is non-zero for some $j \in \{1,2,\ldots,J\}$.
}\end{remark}

\begin{definition}\label{blowup on domain defn}
Given $(m_j,\varphi^{(0)}_j,\varphi^{(\infty)}_j) \in \mathfrak{D}$, we let $\mathfrak{B}(m_j,\varphi^{(0)}_j,\varphi^{(\infty)}_j)$ be the set of all blow-ups $w = (w_j) \in \mathfrak{B}$ associated with $(m_j,\varphi^{(0)}_j,\varphi^{(\infty)}_j)$ as in Remark~\ref{blowup on domain rmk}.  (Note that here and throughout this section, we write $(m_j,\varphi^{(0)}_j,\varphi^{(\infty)}_j)$ to mean $(m_j,\varphi^{(0)}_j,\varphi^{(\infty)}_j)_{1 \leq j \leq J}$).
\end{definition}

We have, of course, that $\mathfrak{B} = \bigcup_{(m_j,\varphi^{(0)}_j,\varphi^{(\infty)}_j) \in \mathfrak{D}} \mathfrak{B}(m_j,\varphi^{(0)}_j,\varphi^{(\infty)}_j)$.

In the next definition, we identify a general point $x = (x_1,x_2) \in \mathbb{R}^2$ with $re^{i\theta}$, where $(r,\theta)$ ($r \geq 0$, $0 \leq \theta < 2\pi$) are polar coordinates of $(x_1,x_2)$.  Recall that whenever $\varphi^{(\infty)}_j$ is non-zero, $(re^{i\theta},y,\op{Re}(c^{(\infty)}_j r^{\alpha} e^{i\alpha\theta}))$ determines a unique point on $\op{graph} \varphi^{(\infty)}_j$ for each $r \geq 0$ and $\theta \in [0,2\pi q_0)$. 

\begin{definition}\label{homogeneous class defn}
Given $(m_j,\varphi^{(0)}_j,\varphi^{(\infty)}_j) \in \mathfrak{D}$, we let $\widetilde{\mathfrak L}(m_j,\varphi^{(0)}_j,\varphi^{(\infty)}_j)$ be the set of all ordered collections $\psi = (\psi_j)_{1 \leq j \leq J}$ of functions 
$$\psi_j : \op{graph}\varphi^{(\infty)}_j |_{\mathbb{R}^n \setminus \{0\} \times \mathbb{R}^{n-2}} \rightarrow \mathcal{A}_{m_j}(\mathbb{R}^m)$$
such that: 
\begin{enumerate}[itemsep=3mm,topsep=0mm]
	\item[(a)]  if $q_j = q_0$ and $\varphi^{(\infty)}_j$ is non-zero, then 
	\begin{equation}\label{homogeneous class eqn1}
		\psi_j(re^{i\theta},y,\op{Re}(c^{(\infty)}_j r^{\alpha} e^{i\alpha\theta})) 
			= \sum_{h=1}^{m_j} \llbracket \op{Re}(a_{j,h} r^{\alpha} e^{i\alpha\theta} + \alpha c^{(0)}_j r^{\alpha-1} e^{i(\alpha-1)\theta} (b \cdot y)) \rrbracket 
	\end{equation}
	for some $a_{j,h} \in \mathbb{C}^m$ and $b \in \mathbb{C}^{n-2}$, with $b$ independent of $j$ and $h$; 
	
	\item[(b)]  if $q_j = 1$ and $\varphi^{(0)}_j$ and $\varphi^{(\infty)}_j$ are both identically zero, then 
	\begin{equation}\label{homogeneous class eqn2}
		\psi_j(re^{i\theta},y,0) = (m_j - q_0 N_j) \llbracket 0 \rrbracket 
			+ \sum_{h=1}^{N_j} \sum_{l=0}^{q_0-1} \llbracket \op{Re}(a_{j,h} r^{\alpha} e^{i\alpha(\theta+2\pi l)}) \rrbracket 
	\end{equation}
	for some integer $N_j \geq 0$ and $a_{j,l} \in \mathbb{C}^m \setminus \{0\}$. 
\end{enumerate}
We let $\mathfrak{L}(m_j,\varphi^{(0)}_j,\varphi^{(\infty)}_j)$ be the set of all ordered collections of functions $\psi = (\psi_j)_{1 \leq j \leq J}$ (as above) such that $\psi_j$ is given by \eqref{homogeneous class eqn1} and \eqref{homogeneous class eqn2} with $b = 0$.  We let $\widetilde{\mathfrak L} = \bigcup_{(m_j,\varphi^{(0)}_j,\varphi^{(\infty)}_j) \in \mathfrak{D}} \widetilde{\mathfrak L}(m_j,\varphi^{(0)}_j,\varphi^{(\infty)}_j)$ and $\mathfrak{L} = \bigcup_{(m_j,\varphi^{(0)}_j,\varphi^{(\infty)}_j) \in \mathfrak{D}} \mathfrak{L}(m_j,\varphi^{(0)}_j,\varphi^{(\infty)}_j)$. 
\end{definition}

\begin{remark}\label{L scale varphi0 rmk}
If $\psi \in \widetilde{\mathfrak L}(m_j,\varphi^{(0)}_j,\varphi^{(\infty)}_j)$, then $\psi \in \widetilde{\mathfrak L}(m_j,c\varphi^{(0)}_j,\varphi^{(\infty)}_j)$ for all $c \in (0,\infty)$.  In particular, the expressions \eqref{homogeneous class eqn1} and \eqref{homogeneous class eqn2} for $\psi_j$ are unchanged if we replace $\varphi^{(0)}$ and $b$ with $c\varphi^{(0)}$ and $b/c$.
\end{remark}

Let $(m_j,\varphi^{(0)}_j,\varphi^{(\infty)}_j) \in \mathfrak{D}$.  Following the above notation, for each open ball $B \subset\subset B_{7/8}(0) \setminus \{0\} \times \mathbb{R}^{n-2}$, there are single-valued harmonic functions $\varphi^{(0)}_{j,l} : B \rightarrow \mathbb{R}^m$ and $\varphi^{(\infty)}_{j,l} : B \rightarrow \mathbb{R}^m$ such that 
\begin{equation}\label{varphi0 localized}
	\varphi^{(0)}_j(x,y) = \sum_{l=1}^{q_j} \llbracket \varphi^{(0)}_{j,l}(x,y) \rrbracket , \quad 
	\varphi^{(\infty)}_j(x,y) = \sum_{l=1}^{q_j} \llbracket \varphi^{(\infty)}_{j,l}(x,y) \rrbracket 
\end{equation}
for all $(x,y) \in B$.  We shall adopt the convention that $\varphi^{(0)}_{j,l}(x,y) = \gamma_j \varphi^{(\infty)}_{j,l}(x,y)$ for all $(x,y) \in B$, where $\gamma_j \in [0,\infty)$ is the constant such that $c^{(0)}_j = \gamma_j c^{(\infty)}_j$ (see Definition~\ref{domain class defn} and note that $\gamma_j = \lim_{\nu\rightarrow\infty} |c^{(\nu)}_j|/\widehat{E}_{\nu}$).  Moreover, if $(\varphi^{(\nu)})$ satisfies the requirements of Remark~\ref{blowup on domain rmk}, we shall adopt the convention that 
\begin{equation}\label{varphi nu localized}
	\varphi^{(\nu)}_j(X) = \sum_{l=1}^{q_j} \llbracket \varphi^{(\nu)}_{j,l}(X) \rrbracket 
\end{equation}
for each $X \in \mathfrak{q}_{\nu}^{-1}(B)$ (as in \eqref{varphi_localized}) where $\varphi^{(\nu)}_{j,l} : \mathfrak{q}_{\nu}^{-1}(B) \rightarrow \mathbb{R}^m$ are single-valued harmonic functions such that if $j \in \{1,2,\ldots,J\}$ such that $\varphi^{(\nu)}_j$ is non-zero for all $\nu$, then $\widehat{E}_{\nu}^{-1} \rho_{\nu}^{\alpha-1} (\mathfrak{q}_{\nu} \circ \varphi^{(\nu)}_{j,l} \circ \mathfrak{q}_{\nu}^{-1}) \rightarrow \varphi^{(0)}_j$ and $|c^{(\nu)}_j|^{-1} \rho_{\nu}^{\alpha-1} (\mathfrak{q}_{\nu} \circ \varphi^{(\nu)}_{j,l} \circ \mathfrak{q}_{\nu}^{-1}) \rightarrow \varphi^{(\infty)}_j$ uniformly on $B$ as $\nu\rightarrow\infty$.  For each $w = (w_j) \in \mathfrak{B}(m_j,\varphi^{(0)}_j,\varphi^{(\infty)}_j)$, there exists $m_j$-valued functions $w_{j,l} : B \rightarrow \mathcal{A}_{m_j}(\mathbb{R}^m)$ given by  
\begin{equation}\label{w localized}
	w_{j,l}(x,y) = w_j(x,y,\varphi^{(\infty)}_{j,l}(x,y)) = \sum_{h=1}^{m_j} \llbracket w_{j,l,h}(x,y) \rrbracket 
\end{equation}
for all $(x,y) \in B$.  Since $w = (w_j)$ is component-wise Dirichlet energy minimizing (as in Definition~\ref{v comp minimizing defn}), each function $w_{j,l}$ is locally Dirichlet energy minimizing in $B$.  Moreover, if $B$ does not intersect the $(n-2)$-dimensional singular set of $w_{j,l}$, each $w_{j,l,h}$ can be chosen to be smooth single-valued harmonic functions $w_{j,l,h} : B \rightarrow \mathbb{R}^m$.  For each Lebesgue measurable subset $\Omega \subset B_1(0)$ we let 
\begin{equation*}
	\|w\|_{L^2(\Omega)}^2 = \int_{\Omega} \sum_{j=1}^J \sum_{l=1}^{q_j} \sum_{h=1}^{m_j} |w_{j,l,h}(x,y)|^2 \,dx \,dy , 
\end{equation*}
where $w_{j,l,h}$ are locally defined by \eqref{w localized}.  Similarly, For each $\psi = (\psi_j) \in \widetilde{\mathfrak L}(m_j,\varphi^{(0)}_j,\varphi^{(\infty)}_j)$, there exists $m_j$-valued functions $\psi_{j,l} : B \rightarrow \mathcal{A}_{m_j}(\mathbb{R}^m)$ and single-valued harmonic functions $\psi_{j,l,h} : B \rightarrow \mathbb{R}^m$ such that 
\begin{equation}\label{psi localized}
	\psi_{j,l}(x,y) = \psi_j(x,y,\varphi^{(\infty)}_{j,l}(x,y)) = \sum_{h=1}^{m_j} \llbracket \psi_{j,l,h}(x,y) \rrbracket 
\end{equation}
for all $(x,y) \in B$.  For each Lebesgue measurable subset $\Omega \subset B_1(0)$ we let 
\begin{equation*}
	\|\psi\|_{L^2(\Omega)}^2 = \int_{\Omega} \sum_{j=1}^J \sum_{l=1}^{q_j} \sum_{h=1}^{m_j} |\psi_{j,l,h}(x,y)|^2 \,dx \,dy , 
\end{equation*}
where $\psi_{j,l,h}$ are locally defined by \eqref{psi localized}.

\subsection{Main estimates for fine blow-ups}  Below we introduce the main properties elements of the blow-up classes $\mathfrak{B}$ and $\widetilde{\mathfrak L}$.  In Lemma~\ref{integrability lemma}, we prove an ``integrability condition'' that each $\psi \in \widetilde{\mathfrak L}(m_j,\varphi^{(0)}_j,\varphi^{(\infty)}_j)$ arises as a blow-up of some sequence $\widetilde{\varphi}^{(\nu)} \in \Phi_{\alpha,q}$ relative to $\varphi^{(\nu)}$.  In Lemma~\ref{blowup norms lemma}, we establish convergence of the excess of $T^{(\nu)}$ relative to $\varphi^{(\nu)}$ to the $L^2$-norm of the blow-up $w \in \mathfrak{B}$, which in part follows from the non-concentration estimate from Corollary~\ref{nonconcentration cor}.  Finally, in Lemma~\ref{blowup estimates lemma}, we establish integral estimates for elements $w$ of the blow-up class $\mathfrak{B}$ which are direct analogues of estimates from Theorem~\ref{keyest thm}, Corollary~\ref{nonconcentration cor}, and Theorem~\ref{nonconcentration thm}.  Using the estimates from Lemma~\ref{blowup estimates lemma}, we can show that the asymptotic analysis blow-up from~\cite{KrumWic2} readily applies to the present setting.

\begin{lemma}\label{integrability lemma} 
Let $(m_j,\varphi^{(0)}_j,\varphi^{(\infty)}_j) \in \mathfrak{D}$ and let $\psi = (\psi_j) \in \widetilde{\mathfrak L}(m_j,\varphi^{(0)}_j,\varphi^{(\infty)}_j)$.  For each $\nu = 1,2,3,\ldots$ let $\varphi^{(\nu)} \in \Phi_{\alpha,q,p}(P_{\nu},S_{\nu})$ satisfies the requirements of Remark~\ref{blowup on domain rmk}, and in particular $\varphi^{(\nu)}/\widehat{E}_{\nu} \rightarrow \varphi^{(0)}$ for some sequence $(\widehat{E}_{\nu})$ of positive numbers with $\widehat{E}_{\nu} \rightarrow 0^+$.  Let $(E_{\nu})$ be a sequence of positive numbers with $E_{\nu} \rightarrow 0^+$ and $E_{\nu} \leq \beta_{\nu} \widehat{E}_{\nu}$, where $\beta_{\nu} \rightarrow 0$.  For each positive integer $\nu$, there exists a rotation $\Gamma_{\nu}$ of $P_{\nu}$ and $\widetilde{\varphi}^{(\nu)} \in \Phi_{\alpha,q}(P_{\nu}, \Gamma_{\nu}(S_{\nu}))$ (so that $\widetilde{\varphi}^{(\nu)} \circ \Gamma_{\nu}^{-1} \in \Phi_{\alpha,q}(P_{\nu},S_{\nu})$) such that $\|\Gamma_{\nu} - I_{P_{\nu}}\| \leq C(n,q,\alpha) \,E_{\nu} \widehat{E}_{\nu}^{-1} \|\psi\|_{L^2(B_1(0))}$ (where $I_{P_{\nu}}$ is the identity map on $P_{\nu}$) and:
\begin{enumerate}[itemsep=3mm,topsep=0mm]
	\item[(i)]  for any ball $B \subset \{ (x,y) \in \mathbb{R}^n : |x| \geq 4 E_{\nu} |b|/\widehat{E}_{\nu} \}$ 
	\begin{align}\label{integrability concl1}
		&\rho_{\nu}^{\alpha-1} (\mathfrak{q}_{\nu} \circ \widetilde{\varphi}^{(\nu)} \circ \mathfrak{q}_{\nu}^{-1})(x,y) 
		\\=\,& \sum_{j=1}^J \sum_{l=1}^{q_j} \sum_{h=1}^{m_j} \llbracket \rho_{\nu}^{\alpha-1} (\mathfrak{q}_{\nu} \circ 
			\varphi^{(\nu)}_{j,l} \circ \mathfrak{q}_{\nu}^{-1})(x,y) + E_{\nu} \psi_{j,l,h}(x,y) + \mathcal{R}_{j,l,h}(x,y) \rrbracket \nonumber 
	\end{align}
	for all $(x,y) \in B$ where $\mathcal{R}_{j,l,h}(x,y) \in \mathbb{R}^m$ with 
	\begin{equation}\label{integrability concl2}
		|\mathcal{R}_{j,l,h}(x,y)| \leq C E_{\nu} |b| \,|c^{(\nu)}_j/\widehat{E}_{\nu}  - c^{(0)}_j| \,|x|^{\alpha-1} 
			+ C \beta_{\nu} E_{\nu} |a_{j,h}| \,|b| \,|x|^{\alpha-1} + C \beta_{\nu} E_{\nu} |b|^2 |x|^{\alpha-2} 
	\end{equation}
	if $\varphi^{(\nu)}_j$ is non-zero and 
	\begin{equation}\label{integrability concl3}
		|\mathcal{R}_{j,l,h}(x,y)| \leq C \beta_{\nu} E_{\nu} |a_{j,h}| \,|b| \,|x|^{\alpha-1} + C \beta_{\nu} E_{\nu} |b|^2 |x|^{\alpha-2} 
	\end{equation}
	if $\varphi^{(\nu)}_j$ is identically zero, where $a_{j,h}$ and $b$ are as in Definition~\ref{homogeneous class defn}, $\varphi^{(\nu)}_{j,l}$ as in \eqref{varphi nu localized}, $\psi_{j,l,h}$ is as in \eqref{psi localized}, and $C = C(n,m,q,\alpha) \in (0,\infty)$ is a constant; 
	
	\item[(ii)]  for $\tau \in (0,1/4)$ and $\nu$ sufficiently large 
	\begin{equation}\label{integrability concl4}
		\rho_{\nu}^{2\alpha-2} \int_{B_{1/2}(0,P_{\nu}) \cap \{ X : \op{dist}(X,S_{\nu}) < \tau\}} \mathcal{G}(\widetilde{\varphi}^{(\nu)},\varphi^{(\nu)})^2 
		\leq C \|\psi\|_{L^2(B_1(0))}^2 \tau^{2\alpha} E_{\nu}^2 ,
	\end{equation}
	where $C = C(n,m,q,\alpha) \in (0,\infty)$ is a constant; 	
	
	\item[(iii)]  
	\begin{equation}\label{integrability concl5}
		\int_{B_{1/2}(0)} \sum_{j=1}^J \sum_{l=1}^{q_j} \sum_{h=1}^{m_j} |\psi_{j,l,h}|^2 
		= \lim_{\nu\rightarrow\infty} \frac{\rho_{\nu}^{2\alpha-2}}{E_{\nu}^2} \int_{B_{1/2}(0,P_{\nu})} 
			\mathcal{G}(\widetilde{\varphi}^{(\nu)},\varphi^{(\nu)})^2 . 
	\end{equation}
\end{enumerate}
\end{lemma}

\begin{proof} 
After rotating and rescaling as to replace $\varphi^{(\nu)}$ with $\rho_{\nu}^{\alpha-1} (\mathfrak{q}_{\nu} \circ \varphi^{(\nu)} \circ \mathfrak{q}_{\nu}^{-1})$ and $\widetilde{\varphi}^{(\nu)}$ with $\rho_{\nu}^{\alpha-1} (\mathfrak{q}_{\nu} \circ \widetilde{\varphi}^{(\nu)} \circ \mathfrak{q}_{\nu}^{-1})$, we argue exactly as in~\cite[Lemma~8.8]{KrumWicc}, which is a modification of~\cite[Lemma~9.7]{KrumWic2}.
\end{proof}

\begin{lemma}\label{blowup norms lemma}
Let $(m_j,\varphi^{(0)}_j,\varphi^{(\infty)}_j) \in \mathfrak{D}$.  Let $w = (w_j) \in \mathfrak{B}(m_j,\varphi^{(0)}_j,\varphi^{(\infty)}_j)$ be the blow-up of $T^{(\nu)}$ relative to a sequence $\varphi^{(\nu)}$ in $B^{\mathcal{M}^{(\nu)}}_{\rho_{\nu}}(Z_{\nu})$ by the excess $E_{\nu} = E_{\varsigma/2,,\varepsilon_{\nu}}(N^{(\nu)},\varphi^{(\nu)},B^{\mathcal{M}^{(\nu)}}_{\rho_{\nu}}(Z_{\nu}))$ (as in Definition~\ref{blowup class defn}).  Let $\psi = (\psi_j) \in \widetilde{\mathfrak L}(m_j,\varphi^{(0)}_j,\varphi^{(\infty)}_j)$ and for $\nu = 1,2,3,\ldots$ let $\widetilde{\varphi}^{(\nu)} \in \Phi_{\alpha,q}$ correspond to $\psi$ as in Lemma~\ref{integrability lemma} with $\widehat{E}_{\nu} = \rho_{\nu}^{-n/2-1} \|N^{(\nu)}\|_{L^2(B^{\mathcal{M}^{(\nu)}}_{\rho_{\nu}}(Z_{\nu}))}$ and $E_{\nu} = E_{\varsigma/2,\varepsilon_{\nu}}(N^{(\nu)},\varphi^{(\nu)},B^{\mathcal{M}^{(\nu)}}_{\rho_{\nu}}(Z_{\nu}))$.  Then for every $\vartheta \in (0,1/2]$ 
\begin{equation}\label{blowup norms concl1}
	\lim_{\nu\rightarrow\infty} \frac{E_{\varsigma/2,1}(N^{(\nu)},\widetilde{\varphi}^{(\nu)},B^{\mathcal{M}^{(\nu)}}_{\vartheta\rho_{\nu}}(Z_{\nu}))^2}{E_{\nu}^2} 
	= \frac{1}{\vartheta^{n+2}} \int_{B_{\vartheta}(0)} \sum_{j=1}^J \sum_{l=1}^{q_j} \mathcal{G}(w_{j,l},\psi_{j,l})^2 .
\end{equation}
In particular, 
\begin{equation}\label{blowup norms concl2}
	\lim_{\nu\rightarrow\infty} \frac{E_{\varsigma/2,1}(N^{(\nu)},\varphi^{(\nu)},B^{\mathcal{M}^{(\nu)}}_{\vartheta\rho_{\nu}}(Z_{\nu}))^2}{E_{\nu}^2} 
	= \frac{1}{\vartheta^{n+2}} \int_{B_{\vartheta}(0)} \sum_{j=1}^J \sum_{l=1}^{q_j} |w_{j,l}|^2 
\end{equation}
and so $\|w\|_{L^2(B_{1/2}(0))} \leq \omega_n^{1/2} (1/2)^{n/2+1}$.
\end{lemma}

\begin{proof}
By Corollary~\ref{nonconcentration cor}, for every $\tau \in (0,1/4)$ and $\nu$ sufficiently large (depending on $\tau$) 
\begin{equation}\label{blowup norms eqn1}
	\int_{B^{\mathcal{M}^{(\nu)}}_{\rho_{\nu}/2}(Z_{\nu}) \cap \{r_{\nu} < \tau\rho_{\nu}\}} \mathcal{G}(N^{(\nu)},\varphi^{(\nu)})^2 
	\leq C \tau^{1/q} E_{\nu}^2 
\end{equation}
for some constant $C = C(n,m,q,M_0,N_0,C_e,C_h,\alpha) \in (0,\infty)$.  In view of this, \eqref{blowup norms concl1} follows from $\varepsilon_{\nu} \rightarrow 0$, the definition of the blow-up $w$, and \eqref{integrability concl1} together with the estimates \eqref{integrability concl2} and \eqref{integrability concl3}.  \eqref{blowup norms concl2} follows by setting $\psi = 0$ in \eqref{blowup norms concl2}.
\end{proof}

\begin{remark}\label{rescale blowup rmk}
Let $(m_j,\varphi^{(0)}_j,\varphi^{(\infty)}_j) \in \mathfrak{D}$ and let $w = (w_j) \in \mathfrak{B}(m_j,\varphi^{(0)}_j,\varphi^{(\infty)}_j)$ be the blow-up of $T^{(\nu)}$ relative to $\varphi^{(\nu)}$ in $B^{\mathcal{M}^{(\nu)}}_{\rho_{\nu}}(Z_{\nu})$ by the excess $E_{\nu} = E_{\varsigma/2,\varepsilon_{\nu}}(N^{(\nu)},\varphi^{(\nu)},B^{\mathcal{M}^{(\nu)}}_{\rho_{\nu}}(Z_{\nu}))$ as in Definition~\ref{blowup class defn}.  Given $\vartheta \in (0,1/2]$ and $w_j \neq m_j \llbracket 0 \rrbracket$ on $\op{graph}\varphi^{(\infty)}_j |_{B_{\vartheta}(0)}$ for some $j \in \{1,2,\ldots,J\}$, let $\widetilde{w} = (\widetilde{w}_j)_{1 \leq j \leq J}$ such that $\widetilde{w}_j : \op{graph}\varphi^{(\infty)}_j |_{B_1(0)} \rightarrow \mathcal{A}_{m_j}(\mathbb{R}^m)$ is given by 
\begin{equation*}
	\widetilde{w}_j(x,y,\varphi^{(\infty)}_{j,l}(x,y)) = \frac{w_j(\rho x, \rho y, \varphi^{(\infty)}_{j,l}(\rho x, \rho y))}{\rho^{-n/2} \|w\|_{L^2(B_{\rho}(0))}} 
\end{equation*}
for all $(x,y) \in B_1(0)$, $j \in \{1,\ldots,J\}$, and $l \in \{1,\ldots,q_j\}$, where $\varphi^{(\infty)}_{j,l}$ is as in \eqref{varphi0 localized}.  Then $\widetilde{w} \in \mathfrak{B}$.  In fact, $E_{\varsigma/2,\varepsilon_{\nu}}(N^{(\nu)},\varphi^{(\nu)},B^{\mathcal{M}^{(\nu)}}_{\vartheta\rho_{\nu}}(Z_{\nu})) > 0$ for all sufficiently large $\nu$.  After passing to a subsequence, $\widetilde{w}$ is a fine blow-up of $T^{(\nu)}$ relative to $\vartheta^{\alpha-1} \varphi^{(\nu)}$ by the excess $E_{\nu} = E_{\varsigma/2,\varepsilon_{\nu}}(N^{(\nu)},\varphi^{(\nu)},B^{\mathcal{M}^{(\nu)}}_{\vartheta\rho_{\nu}}(Z_{\nu}))$.  $\widetilde{w} \in \mathfrak{B}(m_j,c\varphi^{(0)}_j,\varphi^{(\infty)}_j)$ where $c = \lim_{\nu\rightarrow\infty} \vartheta^{n/2+\alpha} \|N^{(\nu)}\|_{L^2(B^{\mathcal{M}^{(\nu)}}_{\rho_{\nu}}(Z_{\nu}))}/\|N^{(\nu)}\|_{L^2(B^{\mathcal{M}^{(\nu)}}_{\vartheta\rho_{\nu}}(Z_{\nu}))}$.  (By Remark~\ref{L scale varphi0 rmk}, $\widetilde{\mathfrak L}(m_j,c\varphi^{(0)}_j,\varphi^{(\infty)}_j) = \widetilde{\mathfrak L}(m_j,\varphi^{(0)}_j,\varphi^{(\infty)}_j)$.

{\rm  To see this, observe that by the triangle inequality, Remark~\ref{graphical rmk}(2), and \eqref{blowups hypothesis2} for all sufficiently large $\nu$ 
\begin{align*}
	\rho_{\nu}^{-n/2-1} \|N^{(\nu)}\|_{L^2(B^{\mathcal{M}^{(\nu)}}_{\vartheta\rho_{\nu}}(Z_{\nu}))} 
	\geq\,& \tfrac{1}{2} \vartheta^{n/2+\alpha} \|\varphi^{(\nu)}\|_{L^2(B_1(0,P_{\nu}))} 
		- E_{\varsigma/2,1}(N^{(\nu)},\varphi^{(\nu)},B^{\mathcal{M}^{(\nu)}}_{\vartheta\rho_{\nu}}(Z_{\nu})) 
	\\ \geq\,& \tfrac{1}{4} \vartheta^{n/2+\alpha} \rho_{\nu}^{-n/2-1} \|N^{(\nu)}\|_{L^2(B^{\mathcal{M}^{(\nu)}}_{\rho_{\nu}}(Z_{\nu}))} .
\end{align*}
Hence after passing to a subsequence, $c = \lim_{\nu\rightarrow\infty} \vartheta^{n/2+\alpha} \|N^{(\nu)}\|_{L^2(B^{\mathcal{M}^{(\nu)}}_{\rho_{\nu}}(Z_{\nu}))}/\|N^{(\nu)}\|_{L^2(B^{\mathcal{M}^{(\nu)}}_{\vartheta\rho_{\nu}}(Z_{\nu}))}$ exists.  By Lemma~\ref{blowup norms lemma}, $E_{\varsigma/2,\varepsilon_{\nu}}(N^{(\nu)},\varphi^{(\nu)},B^{\mathcal{M}^{(\nu)}}_{\vartheta\rho_{\nu}}(Z_{\nu})) > 0$ for all sufficiently large $\nu$.  Clearly conditions~(A)--(D) of Definition~\ref{blowup class defn} hold true with $\vartheta\rho_{\nu}$ in place of $\rho_{\nu}$.  Condition~(F) of Definition~\ref{blowup class defn} holds true with $\delta/\vartheta$ and $\vartheta\rho_{\nu}$ in place of $\delta$ and $\rho_{\nu}$.  To show condition~(E) of Definition~\ref{blowup class defn}, let $p' \in \{0,1,2,\ldots,p-1\}$ and each $\varphi' \in \Phi_{\alpha,q,p'}(P_{\nu},S_{\nu})$.  By the triangle inequality and \eqref{blowups hypothesis2} 
\begin{equation}\label{rescale blowup rmk eqn1}
	E_{\varsigma/2,\varepsilon_{\nu}}(N^{(\nu)},\varphi',B^{\mathcal{M}^{(\nu)}}_{\rho_{\nu}}(Z_{\nu})) 
	\leq 4 \rho_{\nu}^{\alpha-1} \|\mathcal{G}(\varphi',\varphi^{(\nu)})\|_{L^2(B_1(0,P_{\nu}))} 
\end{equation}
and 
\begin{align}\label{rescale blowup rmk eqn2}
	\tfrac{1}{2} (\vartheta\rho_{\nu})^{\alpha-1} \|\mathcal{G}(\varphi',\varphi^{(\nu)})\|_{L^2(B_1(0,P_{\nu}))} 
	\leq\,& E_{\varsigma/2,\varepsilon_{\nu}}(N^{(\nu)},\varphi^{(\nu)},B^{\mathcal{M}^{(\nu)}}_{\vartheta\rho_{\nu}}(Z_{\nu})) 
		\\&+ E_{\varsigma/2,\varepsilon_{\nu}}(N^{(\nu)},\varphi',B^{\mathcal{M}^{(\nu)}}_{\vartheta\rho_{\nu}}(Z_{\nu})) . \nonumber 
\end{align}
By \eqref{blowups hypothesis2}, \eqref{rescale blowup rmk eqn1}, and \eqref{rescale blowup rmk eqn2} 
\begin{align*}
	E_{\varsigma/2,\varepsilon_{\nu}}(N^{(\nu)},\varphi^{(\nu)},B^{\mathcal{M}^{(\nu)}}_{\vartheta\rho_{\nu}}(Z_{\nu})) 
	\leq\,& \vartheta^{-n/2-1} E_{\varsigma/2,\varepsilon_{\nu}}(N^{(\nu)},\varphi^{(\nu)},B^{\mathcal{M}^{(\nu)}}_{\rho_{\nu}}(Z_{\nu})) 
	\\ \leq\,& \beta_{\nu} \vartheta^{-n/2-1} E_{\varsigma/2,\varepsilon_{\nu}}(N^{(\nu)},\varphi',B^{\mathcal{M}^{(\nu)}}_{\rho_{\nu}}(Z_{\nu})) 
	\\ \leq\,& 4 \beta_{\nu} \vartheta^{-n/2-1} \rho_{\nu}^{\alpha-1} \|\mathcal{G}(\varphi',\varphi^{(\nu)})\|_{L^2(B_1(0,P_{\nu}))} 
	\\ \leq\,& 8 \beta_{\nu} \vartheta^{-n/2-\alpha} \rho_{\nu}^{\alpha-1} (E_{\varsigma/2,\varepsilon_{\nu}}(N^{(\nu)},\varphi^{(\nu)},
		B^{\mathcal{M}^{(\nu)}}_{\vartheta\rho_{\nu}}(Z_{\nu})) \\&+ E_{\varsigma/2,\varepsilon_{\nu}}(N^{(\nu)},\varphi',
		B^{\mathcal{M}^{(\nu)}}_{\vartheta\rho_{\nu}}(Z_{\nu}))) .
\end{align*}
Provided $\nu$ is large enough that $8 \beta_{\nu} \vartheta^{-n/2-\alpha} < 1/2$, we conclude that 
\begin{equation*}
	E_{\varsigma/2,\varepsilon_{\nu}}(N^{(\nu)},\varphi^{(\nu)},B^{\mathcal{M}^{(\nu)}}_{\vartheta\rho_{\nu}}(Z_{\nu})) 
	\leq 8\beta_{\nu} \vartheta^{-n/2-\alpha}  \inf_{\varphi' \in \bigcup_{p'=0}^{p-1} \Phi_{\alpha,q,p'}(P_{\nu},S_{\nu})} 
			E_{\varsigma/2,\varepsilon_{\nu}}(N^{(\nu)},\varphi',B^{\mathcal{M}^{(\nu)}}_{\vartheta\rho_{\nu}}(Z_{\nu}))
\end{equation*}
as in condition~(E) of Definition~\ref{blowup class defn} with $8\beta_{\nu} \vartheta^{-n/2-\alpha}$ and $\vartheta\rho_{\nu}$ in place of $\beta_{\nu}$ and $\rho_{\nu}$.  Following the blow-up procedure from Subsection~\ref{sec:blowup subsec}, it is readily verified that $\widetilde{w}$ is a fine blow-up of $T^{(\nu)}$ relative to $\vartheta^{\alpha-1} \varphi^{(\nu)}$ by the excess $E_{\nu} = E_{\varsigma/2,\varepsilon_{\nu}}(N^{(\nu)},\varphi^{(\nu)},B^{\mathcal{M}^{(\nu)}}_{\vartheta\rho_{\nu}}(Z_{\nu}))$.
}\end{remark}

\begin{lemma}\label{blowup estimates lemma}
Let $(m_j,\varphi^{(0)}_j,\varphi^{(\infty)}_j) \in \mathfrak{D}$ and $w = (w_j) \in \mathfrak{B}(m_j,\varphi^{(0)}_j,\varphi^{(\infty)}_j)$.  The following hold true: 
\begin{enumerate}[itemsep=3mm,topsep=0mm]
	\item[(i)]  For each $\psi \in \widetilde{\mathfrak L}(m_j,\varphi^{(0)}_j,\varphi^{(\infty)}_j)$ and each $\rho \in (0,1/2]$  
	\begin{gather}
		\label{blowup estimates concl1} \int_{B_{\rho/2}(0)} \sum_{j=1}^J \sum_{l=1}^{q_j} R^{2-n} \left|\frac{\partial (w_{j,l}/R^{\alpha})}{\partial R}\right|^2 
			\leq C \rho^{-n-2\alpha} \int_{B_{\rho}(0)} \sum_{j=1}^J \sum_{l=1}^{q_j} \mathcal{G}(w_{j,l},\psi_{j,l})^2  , \\
		\label{blowup estimates concl2} \int_{B_{\rho/2}(0)} \sum_{j=1}^J \sum_{l=1}^{q_j} \frac{\mathcal{G}(w_{j,l},\psi_{j,l})^2}{r^{2+1/q}} 
			\leq C \rho^{-2-1/q} \int_{B_{\rho}(0)} \sum_{j=1}^J \sum_{l=1}^{q_j} \mathcal{G}(w_{j,l},\psi_{j,l})^2  , 
	\end{gather}
	where $w_{j,l}$ is as in \eqref{w localized}, $\psi_{j,l}$ is as in \eqref{psi localized}, $R(X) = |X|$ and $r(X) = |x|$ for each $X = (x,y) \in \mathbb{R}^n$, and $C = C(n,m,q,M_0,N_0,C_e,C_h,\alpha) \in (0,\infty)$ is a constant. 
	
	\item[(ii)]  	
	\begin{equation}\label{blowup estimates concl3} 
		\int_{B_{1/2}(0)} \sum_{j=1}^J \sum_{l=1}^{q_j} \sum_{h=1}^{m_j} |D_y w_{j,l}|^2 \leq C  ,
	\end{equation}
	where we let $X = (x,y)$ for $x \in \mathbb{R}^2$ and $y \in \mathbb{R}^{n-2}$ and $C = C(n,m,q,M_0,N_0,C_e,C_h,\alpha) \in (0,\infty)$ is a constant. 
		
	\item[(iii)]  	There exists functions $\kappa : B^{n-2}_{1/4}(0) \rightarrow \mathbb{R}^2$ such that for each $\sigma \in (0,2/q)$ we have that $\|\kappa\|_{C^{1-\sigma/2}(B_{1/2}(0))} \leq C$ and 
	\begin{equation}\label{blowup estimates concl4} 
		\int_{B_{1/2}(0)} \sum_{j=1}^J \sum_{l=1}^{q_j} \sum_{h=1}^{m_j} \frac{|w_{j,l,h}(x,y) + D_x \varphi^{(0)}_{j,l}(x,y) \cdot \kappa(\zeta)|^2
			}{|(x,y) - (0,\zeta)|^{n+2\alpha-\sigma}} \,dx \,dy \leq C 
	\end{equation}
	for all $\zeta \in B^{n-2}_{1/4}(0)$, where $w_{j,l,h}$ is as in \eqref{w localized}, $\varphi^{(0)}_{j,l}$ is as in \eqref{varphi0 localized}, and $C = C(n,m,q,M_0,N_0,$ $C_e,C_h,\alpha,\sigma) \in (0,\infty)$ is a constant.\vspace{3mm}
	
	Moreover, $\kappa(\zeta)$ is unique subject to the condition that for some $\sigma \in (0,2/q)$ 
	\begin{equation*} 
		\int_{B_{1/2}(0)} \sum_{j=1}^J \sum_{l=1}^{q_j} \sum_{h=1}^{m_j} \frac{|w_{j,l,h}(x,y) + D_x \varphi^{(0)}_{j,l}(x,y) \cdot \kappa(\zeta)|^2}{
			|(x,y) - (0,\zeta)|^{n+2\alpha-\sigma}} \,dx \,dy < \infty .
	\end{equation*}
\end{enumerate}
\end{lemma}

\begin{proof} 
Let $w = (w_j)$ be the blow-up of $T^{(\nu)}$ relative to $\varphi^{(\nu)}$ in $B^{\mathcal{M}^{(\nu)}}_{\rho_{\nu}}(Z_{\nu})$ by the excess $E_{\nu} = E_{\varsigma/2,\varepsilon_{\nu}}(N^{(\nu)},\varphi^{(\nu)},B^{\mathcal{M}^{(\nu)}}_{\rho_{\nu}}(Z_{\nu}))$ as in Definition~\ref{blowup class defn}.  Let $\mathfrak{q}_{\nu}$ be the corresponding rotation of $\mathbb{R}^{n+m}$ from Definition~\ref{blowup class defn}.  Let $\widehat{E}_{\nu} = \rho_{\nu}^{-n/2-1} \|N^{(\nu)}\|_{L^2(B^{\mathcal{M}^{(\nu)}}_{\rho_{\nu}}(Z_{\nu}))}$.  Let $\widetilde{\varphi}^{(\nu)}$ be the function corresponding to $\psi$, $\varphi^{(\nu)}$, $E_{\nu}$, and $\widehat{E}_{\nu}$ as in Lemma~\ref{rescale blowup rmk}.  Rotate and rescale by letting $\overline{T}^{(\nu)}$, $\overline{\mathcal M}^{(\nu)}$, $\overline{N}^{(\nu)}$, and $\overline{\varphi}^{(\nu)}$ be as in \eqref{blowup eqn1}.  Let $\overline{\mathbf p}_{\nu}$ be the closest point projection map onto $\overline{\mathcal M}^{(\nu)}$ and let $\overline{\mathbf e}_{\nu}$ be the exponential map of $\overline{\mathcal M}^{(\nu)}$ at the origin.  Let $v^{(\nu)}_j : \mathbf{G}_{\overline{\varphi}^{(\nu)}_j \circ \overline{\mathbf e}_{\nu}^{-1}} \cap \mathbf{C}^{\overline{\mathcal M}^{(\nu)}}_{3/4}(0) \cap \{r_{\nu} > \tau_{\nu}\} \rightarrow \mathcal{A}_{m_j}(\mathbb{R}^{n+m})$ be as in \eqref{blowup eqn6} and \eqref{blowup eqn7}, where $r_{\nu}(X) = \op{dist}(\overline{\mathbf e}_{\nu}^{-1}(\overline{\mathbf p}_{\nu}(X)),\{0\} \times \mathbb{R}^{n-2})$ for each $X \in \mathbf{C}^{\overline{\mathcal M}^{(\nu)}}_{3/4}(0)$ and $\tau_{\nu} \rightarrow 0$.  Recall that $\mathbf{I}_{N^{(\nu)},Z_{\nu}}(0^+) \geq \alpha$.  By Theorem~\ref{keyest thm} (with $\mathcal{M}^{(\nu)},N^{(\nu)},\varphi^{(\nu)},Z_{\nu},\rho\rho_{\nu}$ in place of $\mathcal{M},N,\varphi,Z_0,\rho_0$), \eqref{blowup eqn6}, and the homogeneity of $\varphi^{(\nu)}$ 
\begin{equation}\label{blowup estimates eqn1}
	\int_{B^{\overline{\mathcal M}^{(\nu)}}_{1/2}(0) \cap \{r_{\nu} > \tau_{\nu}\}} 
		R_{\nu}^{2-n} \,\bigg| \partial_{\hat R_{\nu}} \bigg(\frac{v^{(\nu)}}{R_{\nu}^{\alpha}}\bigg) \bigg|^2 
	\leq C E_{\varsigma/2,\varepsilon_{\nu}}(N^{(\nu)},\varphi^{(\nu)},B^{\mathcal{M}^{(\nu)}}_{\rho \rho_{\nu}}(Z_{\nu}))^2 , 
\end{equation}
where $R_{\nu}(X)$ is the geodesic distance of $X$ from the origin in $\overline{\mathcal M}^{(\nu)}$, $\partial_{\hat R_{\nu}}$ is the derivative in the direction of the gradient of $R_{\nu}$, and $C = C(n,m,q,M_0,N_0,C_e,C_h,\alpha) \in (0,\infty)$ is a constant.  Recall that $\overline{\mathcal M}^{(\nu)} \rightarrow \mathbb{R}^n \times \{0\}$ in the $C^2$-topology in $\mathbf{C}_1(0)$ and in particular \eqref{exp map est} and \eqref{normal coord eqn4} hold true.  Recall from Subsection~\ref{sec:blowup subsec} that $(w_j^{(\nu)})$ is as in \eqref{blowup eqn10} and \eqref{blowup eqn11} and satisfies $w_j^{(\nu)} \rightarrow w_j$ and $|Dw_j^{(\nu)}| \rightarrow |Dw_j|$ in $L^2(\op{spt} \mathbf{G}_{\varphi^{(\infty)}_j} \cap \Omega' \times \mathbb{R}^m)$ for each $\Omega' \subset\subset \Omega \setminus \{0\}\times\mathbb{R}^{n-2}$.  In particular, letting $L = L(m,q) \geq 1$ be an integer and $\bs{\xi} : \mathcal{A}_q(\mathbb{R}^m) \rightarrow \mathbb{R}^L$ be White's bi-Lipschitz injection (\cite[Definition~1.1(6) and Theorem~1.2]{Almgren}, \cite[p.~706]{Cha88}, \cite[Corollary~2.2]{DeLSpaDirMin}), $\bs{\xi} \circ w_j^{(\nu)} \rightarrow \bs{\xi} \circ w_j$ in $W^{1,2}(\op{spt} \mathbf{G}_{\varphi^{(\infty)}_j} \cap \Omega' \times \mathbb{R}^m)$ for each $\Omega' \subset\subset \Omega \setminus \{0\}\times\mathbb{R}^{n-2}$.  Using all this together with Lemma~\ref{blowup norms lemma}, by dividing both sides of \eqref{blowup estimates eqn1} by $E_{\nu}^2$ and letting $\nu\rightarrow\infty$, we deduce that \eqref{blowup estimates concl1} holds true.  By similar reasoning using Lemma~\ref{rescale blowup rmk}, it follows from Corollary~\ref{nonconcentration cor} that \eqref{blowup estimates concl2} holds true and it follows from Theorem~\ref{keyest thm} that \eqref{blowup estimates concl3} holds true.

To see~(iii), let $\zeta \in B^{n-2}_{1/4}(0)$.  By \eqref{blowups no gaps}, there exists $X_{\nu} \in \mathcal{B}^{(\geq 2)}_q(T^{(\nu)}) \cap B^{\mathcal{M}^{(\nu)}}_{\delta_{\nu}\rho_{\nu}}(\mathbf{e}_{\nu}(\mathfrak{q}_{\nu}^{-1}(0,\rho_{\nu} \zeta)))$ such that $\mathbf{I}_{N^{(\nu)},X_{\nu}}(0^+) \geq \alpha$.  Let $X_{\nu} = \mathbf{e}_{\nu}(\xi_{\nu}+\zeta_{\nu})$ where $\xi_{\nu} \in P_{\nu} \cap S_{\nu}^{\perp}$ and $\zeta_{\nu} \in S_{\nu}$.  By Theorem~\ref{nonconcentration thm}(A), $\widehat{E}_{\nu} |\xi_{\nu}| \leq C E_{\nu} \rho_{\nu}$ and thus after passing to a subsequence $\widehat{E}_{\nu} \mathfrak{q}_{\nu}(\xi_{\nu})/(E_{\nu} \rho_{\nu})$ converges to some $\kappa(\zeta) \in \mathbb{R}^2$.  Now arguing as above, it follows from Theorem~\ref{nonconcentration thm}(D) that \eqref{blowup estimates concl4}.  The uniqueness of $\kappa(\zeta)$ and the bound $\|\kappa\|_{C^{1-\sigma/2}(B_{1/2}(0))} \leq C$ follow by arguing as in the proof of~\cite[Lemma~9.9(c)]{KrumWic2}.
\end{proof}

\subsection{Classification of homogeneous degree $\alpha$ blow-ups}  

\begin{theorem}\label{classify in L thm}
Let $w = (w_j) \in \mathfrak{B}$ be homogeneous degree $\alpha$ in the sense that 
\begin{equation*}
	(x,y) \,\mapsto\, \sum_{l=1}^{q_j} \sum_{h=1}^{m_j} \llbracket w_{j,l,h}(x,y) \rrbracket 
\end{equation*}
is a homogeneous degree $\alpha$, $\mathcal{A}_{m_j q_j}(\mathbb{R}^m)$-valued function of $(x,y) \in B_{1/2}(0) \setminus \{0\} \times \mathbb{R}^{n-2}$ for all $j \in \{1,2,\ldots,J\}$, where $w_{j,l,h}$ is as in \eqref{w localized}.  Then $w \in \widetilde{\mathfrak L}$.
\end{theorem}

\begin{proof}  
Recall from Lemma~\ref{blowup estimates lemma} that there exists a function $\kappa : B^{n-2}_{1/4}(0) \rightarrow \mathbb{R}^2$ such that \eqref{blowup estimates concl4} holds true for all $\zeta \in B^{n-2}_{1/4}(0)$ and $\sigma \in (0,1/q)$. 

\noindent\textit{Step~1.  Proof when $\kappa = 0$.}  In the case that \eqref{blowup estimates concl4} holds true with $\kappa = 0$, we can show that $w_j$ is of the form of $\psi_j$ in \eqref{homogeneous class eqn1} and \eqref{homogeneous class eqn2} with $b = 0$, i.e.~$w \in \mathfrak{L}$.  This follows by arguing exactly as in Step~1 of the proof of~\cite[Theorem~10.1]{KrumWic2}.  However, like in~\cite[Theorem~8.12]{KrumWicc}, since $\alpha > 1$ the inner variational formula (10.6) of~\cite[Theorem~10.1]{KrumWic2} follows from the cut-off function argument used to prove the outer variational formula (10.5) of~\cite[Theorem~10.1]{KrumWic2} and inner variational formula (10.6) of~\cite[Theorem~10.1]{KrumWic2}.

\noindent\textit{Step~2.  Proof when $\kappa \neq 0$.}  Arguing exactly as in Step~3 of the proof of~\cite[Theorem~8.12]{KrumWicc}, we can show that \eqref{blowup estimates concl4} holds true with $\kappa(\zeta) = (b_1 \cdot \zeta, b_2 \cdot \zeta)$ for all $\zeta \in B^{n-2}_{1/16}(0)$, where $b_1,b_2 \in \mathbb{R}^{n-2}$.  Let 
\begin{equation*}
	B = \left[\begin{array}{ccc} 
		0 & 0 & b_1 \\
		0 & 0 & b_2 \\
		b_1^T & b_2^T & 0
	\end{array}\right]
\end{equation*}
where we represent $b_1,b_2$ as row vectors.  Note that 
\begin{equation}\label{classify in L eqn1}
	\|\Gamma_{\nu} - I\| \leq 2 E_{\nu} (|b_1|+|b_2|)/\widehat{E}_{\nu} .
\end{equation}
Let $\widetilde{\Gamma}_{\nu}$ be a rotation of $\mathbb{R}^n$ given by $\widetilde{\Gamma}_{\nu} (x',x'') = (e^{-E_{\nu} B/\widehat{E}_{\nu}} x', x'')$ for each $x' \in \mathbb{R}^n$ and $x'' \in \mathbb{R}^m$.  Let $\Gamma_{\nu} = \mathfrak{q}_{\nu}^{-1} \circ \widetilde{\Gamma}_{\nu} \circ \mathfrak{q}_{\nu}$ as a rotation of $\mathbb{R}^{n+m}$ which fixes vectors in $P_{\nu}^{\perp}$.  Following~\cite[Theorem~10.1]{KrumWic2}, we would like to construct the blow-up $w^{\star}$ of $\Gamma_{\nu\#} T^{(\nu)}$ relative to $\varphi^{(\nu)}$.  However, in the construction of the center manifold we fixed the plane $\mathbb{R}^n \times \{0\}$ and $\Gamma_{\nu}$ might rotate $\mathbb{R}^n \times \{0\}$.  Instead, we blow-up $T^{(\nu)}$ relative to $\varphi^{(\nu)} \circ \Gamma_{\nu}$.  Observe that $\varphi^{(\nu)} \circ \Gamma_{\nu} \in \Phi_{\alpha,q,p}(P_{\nu},\Gamma^{-1}_{\nu}(S_{\nu}))$.  By the triangle inequality, Remark~\ref{graphical rmk}(2), and \eqref{classify in L eqn1} 
\begin{equation}\label{classify in L eqn2}
	E_{\varsigma/2,\varepsilon_{\nu}}(N^{(\nu)},\varphi^{(\nu)} \circ \Gamma_{\nu}, B^{\mathcal{M}^{(\nu)}}_{\rho_{\nu}}(Z_{\nu})) \leq C E_{\nu} , 
\end{equation}
where $C = C(n,m,q,M_0,N_0,C_e,C_h,\alpha) \in (0,\infty)$ is a constant.  Clearly condition~(A)--(D) of Definition~\ref{blowup class defn} hold true with $\varphi^{(\nu)} \circ \Gamma_{\nu}$ and $\Gamma_{\nu}^{-1}(S_{\nu}))$ in place of $\varphi^{(\nu)}$ and $S_{\nu}$.  By the computation leading to (10.30) of~\cite[Theorem~10.1]{KrumWic2}, it follows using Remark~\ref{graphical rmk}(1), \eqref{blowups hypothesis2}, \eqref{classify in L eqn1}, and \eqref{classify in L eqn2} that condition~(E) of Definition~\ref{blowup class defn} holds true with $C\beta_{\nu}$, $\varphi^{(\nu)} \circ \Gamma_{\nu}$, and $\Gamma_{\nu}^{-1}(S_{\nu}))$ in place of $\beta_{\nu}$, $\varphi^{(\nu)}$, and $S_{\nu}$.  By \eqref{blowups hypothesis2} and \eqref{classify in L eqn1}, condition~(F) of Definition~\ref{blowup class defn} holds true with $\delta_{\nu} + 2 \beta_{\nu} \,(|b_1|+|b_2|)$ and $\Gamma_{\nu}^{-1}(S_{\nu}))$ in place of $\delta_{\nu}$ and $S_{\nu}$.  Let $w^{\star} = (w^{\star}_j)$ be the blow-up of $T^{(\nu)}$ relative to $\varphi^{(\nu)} \circ \Gamma_{\nu}$ by excess $E_{\varsigma/2,\varepsilon_{\nu}}(N^{(\nu)},\varphi^{(\nu)} \circ \Gamma_{\nu},B^{\mathcal{M}^{(\nu)}}_{\rho_{\nu}}(Z_{\nu}))$.  By Taylor's theorem, for each $\tau \in (0,1/2]$ 
\begin{equation*}
	\varphi^{(\nu)}(\Gamma_{\nu} X) 
	= \sum_{j=1}^J \sum_{l=1}^{q_j} \llbracket \varphi^{(\nu)}_{j,l}(X) - D\varphi^{(\nu)}_{j,l}(X) 
		\cdot E_{\nu} \mathfrak{q}_{\nu}^{-1} (b_1 \cdot y,b_2 \cdot y,0)/\widehat{E}_{\nu} + \mathcal{R}^{(\nu)}_{j,l}(X) \rrbracket 
\end{equation*}
for each $X \in B_1(0,P_{\nu})$ with $\op{dist}(X,S_{\nu}) \geq \tau$, where $X = \mathfrak{q}_{\nu}^{-1}(x,y)$ for $x \in \mathbb{R}^2$ and $y \in \mathbb{R}^{n-2}$, $\varphi^{(\nu)}_{j,l}$ as in \eqref{varphi nu localized}, and $E_{\nu}^{-1} \mathcal{R}^{(\nu)}_{j,l}(X) \rightarrow 0$ as $\nu\rightarrow\infty$.  Hence it follows that 
\begin{equation}\label{classify in L eqn3}
	c^{\star} w^{\star}_j(X,\varphi^{(\infty)}_{j,l}(X)) = \sum_{h=1}^{m_j} \llbracket w_{j,l,h}(X) 
		+ D\varphi^{(0)}_{j,l}(X) \cdot (b_1 \cdot y, b_2 \cdot y, 0) \rrbracket
\end{equation}
for each $X \in B_{1/2}(0)$ if $\varphi^{(\infty)}_j$ is non-zero, where $\varphi^{(0)}_{j,l}$ and $\varphi^{(\infty)}_{j,l}$ are as in \eqref{varphi0 localized} and $w_{j,l,h}$ is as in \eqref{w localized}, and $c^{\star} w^{\star}_j = w_j$ if $\varphi^{(\infty)}_j$ is identically zero.  Here $c^{\star} = \lim_{\nu\rightarrow\infty} E_{\nu}^{-1} E_{\varsigma/2,\varepsilon_{\nu}}(N^{(\nu)},\varphi^{(\nu)} \circ \Gamma_{\nu}, B^{\mathcal{M}^{(\nu)}}_{\rho_{\nu}}(Z_{\nu}))$, which after passing to a subsequence exists by \eqref{classify in L eqn2}.  If $c^{\star} = 0$, then it follows directly from \eqref{classify in L eqn3} that $w \in \widetilde{\mathfrak L}$.  If on the other hand $c^{\star} = 0$, then it follows from Step~1 with $w^{\star}$ in place of $w$ that $w^{\star} \in \mathfrak{L}$ and thus by \eqref{classify in L eqn3} we have that $w \in \widetilde{\mathfrak L}$.
\end{proof}

\subsection{Asymptotic decay of fine blow-ups}\label{sec:fine blowups decay}  The main result of this section is a decay estimate for the fine blow-ups, Theorem~\ref{blowup decay thm}.  Our proof of Theorem~\ref{blowup decay thm} employs the well-known hole-filling technique, in a manner similar to its use in~\cite{Sim93}, and will be based on the following lemma.

\begin{lemma}\label{reverseHS lemma}
Let $(m_j,\varphi^{(0)}_j,\varphi^{(\infty)}_j) \in \mathfrak{D}$, $w = (w_j) \in \mathfrak{B}(m_j,\varphi^{(0)}_j,\varphi^{(\infty)}_j)$, and $\psi \in \widetilde{\mathfrak L}(m_j,\varphi^{(0)}_j,\varphi^{(\infty)}_j)$ such that 
\begin{equation}\label{reverseHS hyp1}
	\int_{B_{1/2}(0)} \sum_{j=1}^J \sum_{l=1}^{q_j} \mathcal{G}(w_{j,l},\psi_{j,l})^2 
	\leq 2 \inf_{\psi' \in \widetilde{\mathfrak L}(m_j,\varphi^{(0)}_j,\varphi^{(\infty)}_j)} 
		\int_{B_{1/2}(0)} \sum_{j=1}^J \sum_{l=1}^{q_j} \mathcal{G}(w_{j,l},\psi'_{j,l})^2 ,
\end{equation}
where $w_{j,l}$ is as in \eqref{w localized}, $\psi_{j,l}$ is as in \eqref{psi localized}, and $\psi'_{j,l}(x,y) = \psi'_j(x,y,\varphi^{(\infty)}_{j,l}(x,y))$ for each $(x,y) \in B_{1/2}(0)$ where $\varphi^{(\infty)}_{j,l}$ is as in \eqref{varphi0 localized} (as in \eqref{psi localized}).  Then 
\begin{equation}\label{reverseHS concl}
	\int_{B_{1/2}(0)} \sum_{j=1}^J \sum_{l=1}^{q_j} \mathcal{G}(w_{j,l},\psi_{j,l})^2 
	\leq C \int_{B_{1/2}(0) \setminus B_{1/8}(0)} \sum_{j=1}^J \sum_{l=1}^{q_j} \left|\frac{\partial (w_{j,l}/R^{\alpha})}{\partial R}\right|^2 ,
\end{equation}
where $R(x,y) = |(x,y)|$ for each $(x,y) \in \mathbb{R}^n$ and $C = C(n,m,q,M_0,N_0,C_e,C_h,\alpha) \in (0,\infty)$ is a constant.
\end{lemma}

Since our setting is higher multiplicity and the functions involved are multi-valued (in contrast with the multiplicity one setting of~\cite{Sim93}), we prove Lemma~\ref{reverseHS lemma} following an approach introduced in~\cite{KrumWic2} in which we first prove Lemma~\ref{reverseHS lemma0} below.  Let $(m_j,\varphi^{(0)}_j,\varphi^{(\infty)}_j) \in \mathfrak{D}$ and $\psi \in \mathfrak{L}(m_j,\varphi^{(0)}_j,\varphi^{(\infty)}_j)$.  Given $j \in \{1,2,\ldots,J\}$ such that $\varphi^{(\infty)}_j$ is non-zero, we know that $\psi_j$ is given by \eqref{homogeneous class eqn1} for some $a_{j,h} \in \mathbb{C}^m$ and $b \in \mathbb{C}^{n-2}$, with $b$ independent of $j$ and $l$.  Let $d_j(\psi)$ denote the number of distinct values of $a_{j,h}$, i.e.~the number of distinct components $\op{Re}(a_{j,h} r^{\alpha} e^{i\theta\theta} + \alpha c^{(0)}_j r^{\alpha-1} e^{i(\alpha-1)\theta} (b \cdot y))$ as in \eqref{homogeneous class eqn1}.  Similarly, given $j \in \{1,2,\ldots,J\}$ such that $\varphi^{(0)}_j$ and $\varphi^{(0)}_j$ are both identically zero, $\psi_j$ is given by \eqref{homogeneous class eqn2} for some positive integer $N$ and some $a_{j,h} \in \mathbb{C}^m$.  Let $d_j(\psi)$ denote the number of distinct values of $a_{j,h}$, i.e.~the number of distinct non-zero components $\op{Re}(a_{j,h} r^{\alpha} e^{i\theta\theta})$ as in \eqref{homogeneous class eqn2}.  Let 
\begin{equation*}
	d(\psi) = \sum_{j=1}^J d_j(\psi). 
\end{equation*}
Note that $1 \leq d_j(\psi) \leq m_j$ if $\varphi^{(\infty)}_j$ is non-zero and $0 \leq d_j(\psi) \leq \lfloor m_j/q_0 \rfloor$ if $\varphi^{(0)}_j$ and $\varphi^{(\infty)}_j$ are both identically zero.  It follows that $p \leq d(\psi) \leq \sum_{j=1}^J m_j q_j/q_0 = q/q_0$; that is, $p \leq d(\psi) \leq \lfloor q/q_0 \rfloor$.  (Here $p$ is the number of non-zero components of $\varphi^{(\nu)}$ as in Definition~\ref{blowup class defn}, and thus $p$ is equal to the number of $j \in \{1,2,\ldots,J\}$ for which $\varphi^{(\infty)}_j$ is non-zero.)  For $s \in \{p,p+1,\ldots,\lfloor q/q_0 \rfloor\}$, let 
\begin{equation*}
	\mathfrak{L}_s(m_j,\varphi^{(0)}_j,\varphi^{(\infty)}_j) = \{ \psi \in \mathfrak{L}(m_j,\varphi^{(0)}_j,\varphi^{(\infty)}_j) : d(\psi) = s \} 
\end{equation*}
and note that $\mathfrak{L}(m_j,\varphi^{(0)}_j,\varphi^{(\infty)}_j) = \bigcup_{s=p}^{\lfloor q/q_0 \rfloor} \mathfrak{L}_s(m_j,\varphi^{(0)}_j,\varphi^{(\infty)}_j)$.

\begin{lemma}\label{reverseHS lemma0}
For every $\overline{M} \in [1,\infty)$ there exists $\overline{\beta} \in (0,1)$ and $\overline{C} \in (0,\infty)$ depending only on $n,m,q,M_0,N_0,C_e,C_h,\alpha$ such that the following holds true.  Let $(m_j,\varphi^{(0)}_j,\varphi^{(\infty)}_j) \in \mathfrak{D}$ and let $p$ equal the number of $j \in \{1,2,\ldots,J\}$ for which $\varphi^{(\infty)}_j$ is non-zero.  Let $s \in \{p,p+1,\ldots,\lfloor q/q_0 \rfloor\}$, $w = (w_j) \in \mathfrak{B}(m_j,\varphi^{(0)}_j,\varphi^{(\infty)}_j)$, and $\psi \in \mathfrak{L}_s(m_j,\varphi^{(0)}_j,\varphi^{(\infty)}_j)$ such that 
\begin{equation}\label{reverseHS hyp2}
	\int_{B_{1/2}(0)} \sum_{j=1}^J \sum_{l=1}^{q_j} \mathcal{G}(w_{j,l},\psi_{j,l})^2 
	\leq 2 \overline{M}^2 \inf_{\psi' \in \widetilde{\mathfrak L}(m_j,\varphi^{(0)}_j,\varphi^{(\infty)}_j)} 
		\int_{B_{1/2}(0)} \sum_{j=1}^J \sum_{l=1}^{q_j} \mathcal{G}(w_{j,l},\psi'_{j,l})^2 ,
\end{equation}
and either (i) $s = p$ or (ii) $s > p$ and 
\begin{equation}\label{reverseHS hyp3}
	\int_{B_{1/2}(0)} \sum_{j=1}^J \sum_{l=1}^{q_j} \mathcal{G}(w_{j,l},\psi_{j,l})^2 
	\leq \overline{\beta}^2 \inf_{\psi' \in \bigcup_{s'=p}^{s-1} \mathfrak{L}_{s'}(m_j,\varphi^{(0)}_j,\varphi^{(\infty)}_j)} 
		\int_{B_{1/2}(0)} \sum_{j=1}^J \sum_{l=1}^{q_j} \mathcal{G}(w_{j,l},\psi'_{j,l})^2 
\end{equation}
(where $w_{j,l}$, $\psi_{j,l}$, and $\psi'_{j,l}$ are as in Lemma~\ref{reverseHS lemma}).  Then \eqref{reverseHS concl} holds true with $C = \overline{C}$.
\end{lemma}

\begin{proof} 
The proof is nearly identical to that of~\cite[Lemma~5.10]{KrumWicb} and~\cite[Lemma~8.13]{KrumWicc} with obvious modifications; we refer to the reader to~\cite{KrumWicb} and~\cite{KrumWicc} for details.  Note that since the fine excess $E_{\varsigma,\varepsilon}(N,\varphi,B^{\mathcal M}_{\rho_0}(Z_0))$ is defined in \eqref{graphical E defn} in terms of the $L^2$ multi-valued distance between $N$ and $\varphi$, Steps~2 and 3 of the proof of~\cite[Lemma~8.13]{KrumWicc} can be replaced by a straight-forward triangle inequality argument, see the argument leading to (11.10) of~\cite{KrumWic2}.
\end{proof}

\begin{proof}[Proof of Lemma~\ref{reverseHS lemma}]
Lemma~\ref{reverseHS lemma} follows from Lemma~\ref{reverseHS lemma0} by arguing as in~\cite[Lemma~11.1]{KrumWic2} (for which Lemma~11.1 follows from Lemma~11.5).
\end{proof}

\begin{lemma}\label{hole filling lemma}
Let $(m_j,\varphi^{(0)}_j,\varphi^{(\infty)}_j) \in \mathfrak{D}$ and $w = (w_j) \in \mathfrak{B}(m_j,\varphi^{(0)}_j,\varphi^{(\infty)}_j)$.  For each $\vartheta \in (0,1/8]$, let $\psi^{(\vartheta)} = (\psi^{(\vartheta)}_j) \in \widetilde{\mathfrak L}(m_j,\varphi^{(0)}_j,\varphi^{(\infty)}_j)$ such that 
\begin{equation}\label{hole filling proj}
	\int_{B_{\vartheta}(0)} \sum_{j=1}^J \sum_{l=1}^{q_j} \mathcal{G}(w_{j,l},\psi^{(\vartheta)}_{j,l})^2 
	\leq 2 \inf_{\psi' \in \widetilde{\mathfrak L}(m_j,\varphi^{(0)}_j,\varphi^{(\infty)}_j)} 
		\int_{B_{\vartheta}(0)} \sum_{j=1}^J \sum_{l=1}^{q_j} \mathcal{G}(w_{j,l},\psi'_{j,l})^2
\end{equation}
where $w_{j,l}$ is as in \eqref{w localized} and $\psi^{(\vartheta)}_{j,l}(x,y) = \psi^{(\vartheta)}_j(x,y,\varphi^{(\infty)}_{j,l}(x,y))$ and $\psi'_{j,l}(x,y) = \psi'_j(x,y,\varphi^{(\infty)}_{j,l}(x,y))$ for each $(x,y) \in B_{1/2}(0)$ where $\varphi^{(\infty)}_{j,l}$ is as in \eqref{varphi0 localized} (as in \eqref{psi localized}).  Then for each $\vartheta_1,\vartheta_2$ with $0 < \vartheta_1 \leq \vartheta_2 \leq 1/8$ 
\begin{equation}\label{hole filling concl}
	\vartheta_1^{-n-2\alpha} \int_{B_{\vartheta_1}(0)} \sum_{j=1}^J \sum_{l=1}^{q_j} \mathcal{G}(w_{j,l},\psi^{(\vartheta_1)}_{j,l})^2 
	\leq C \bigg(\frac{\vartheta_1}{\vartheta_2}\bigg)^{2\mu} \vartheta_2^{-n-2\alpha} \int_{B_{\vartheta_2}(0)} \sum_{j=1}^J \sum_{l=1}^{q_j} 
		\mathcal{G}(w_{j,l},\psi^{(\vartheta_2)}_{j,l})^2 ,
\end{equation}
where $\mu = \mu(n,m,q,\alpha) \in (0,1)$ and $C = C(n,m,q,M_0,N_0,C_e,C_h,\alpha) \in (0,\infty)$ are constants (independent of $\vartheta_1$ and $\vartheta_2$).
\end{lemma}

\begin{proof}
Without loss of generality, assume that $w$ is not identically zero on $B_{\varsigma_1}(0)$.  Noting that \eqref{hole filling concl} obviously holds true if $\vartheta_1 \geq \vartheta_2/32$, we may assume that $\vartheta_1 < \vartheta_2/32$.  Let $\rho \in [2\vartheta_1,1]$.  By \eqref{blowup estimates concl1} with $\psi^{(\rho)}$ in place of $\psi$ 
\begin{equation}\label{hole filling eqn1} 
	\int_{B_{\rho/2}(0)} \sum_{j=1}^J \sum_{l=1}^{q_j} R^{2-n} \left|\frac{\partial (w_{j,l}/R^{\alpha})}{\partial R}\right|^2 
	\leq C \rho^{-n-2\alpha} \int_{B_{\rho}(0)} \sum_{j=1}^J \sum_{l=1}^{q_j} \mathcal{G}(w_{j,l},\psi^{(\rho)}_{j,l})^2  , 
\end{equation}
where $C = C(n,m,q,M_0,N_0,C_e,C_h,\alpha) \in (0,\infty)$ is a constant.  Let $w^{(\rho)} = (w^{(\rho)}_1,\ldots,w^{(\rho)}_J)$ so that $w^{(\rho)}_j : \op{graph}\varphi^{(\infty)}_j \rightarrow \mathcal{A}_{m_j}(\mathbb{R}^m)$ is given by 
\begin{equation*}
	w^{(\rho)}_j(x,y,\varphi^{(\infty)}_{j,l}(x,y)) = \frac{w_j(\rho x, \rho y, \varphi^{(\infty)}_{j,l}(\rho x, \rho y))}{\rho^{-n/2} \|w\|_{L^2(B_{\rho}(0))}} 
\end{equation*}
for all $(x,y) \in B_1(0)$, $j \in \{1,\ldots,J\}$, and $l \in \{1,\ldots,q_j\}$, where $\varphi^{(\infty)}_{j,l}$ is as in \eqref{varphi0 localized}.  By Remark~\ref{rescale blowup rmk}, since $w$ is not identically zero in $B_{\rho}(0)$, $w^{(\rho)} \in \mathfrak{B}(m_j,c\varphi^{(0)}_j,\varphi^{(\infty)}_j)$, where we let $c = \lim_{\nu\rightarrow\infty} \vartheta^{n/2+\alpha} \|N^{(\nu)}\|_{L^2(B^{\mathcal{M}^{(\nu)}}_{\rho_{\nu}}(Z_{\nu}))}/\|N^{(\nu)}\|_{L^2(B^{\mathcal{M}^{(\nu)}}_{\vartheta\rho_{\nu}}(Z_{\nu}))} \in (0,\infty)$.  Recall from Remark~\ref{L scale varphi0 rmk} that $\psi^{(\rho)} \in \widetilde{\mathfrak L}(m_j,c\varphi^{(0)}_j,\varphi^{(\infty)}_j)$.  By Lemma~\ref{reverseHS lemma} with $w^{(\rho)}$ in place of $w$ 
\begin{equation}\label{hole filling eqn2} 
	\rho^{-n-2\alpha} \int_{B_{\rho/2}(0)} \sum_{j=1}^J \sum_{l=1}^{q_j} \mathcal{G}(w_{j,l},\psi^{(\rho/2)}_{j,l})^2 
	\leq C \int_{B_{\rho/2}(0) \setminus B_{\rho/8}(0)} \sum_{j=1}^J \sum_{l=1}^{q_j} R^{2-n} \left|\frac{\partial (w_{j,l}/R^{\alpha})}{\partial R}\right|^2 , 
\end{equation}
where $C = C(n,m,q,M_0,N_0,C_e,C_h,\alpha) \in (0,\infty)$ is a constant.  By combining \eqref{hole filling eqn1} and \eqref{hole filling eqn2} 
\begin{equation}\label{hole filling eqn3} 
	\int_{B_{\rho/8}(0)} \sum_{j=1}^J \sum_{l=1}^{q_j} R^{2-n} \left|\frac{\partial (w_{j,l}/R^{\alpha})}{\partial R}\right|^2 
	\leq C_0 \int_{B_{\rho/2}(0) \setminus B_{\rho/8}(0)} \sum_{j=1}^J \sum_{l=1}^{q_j} R^{2-n} \left|\frac{\partial (w_{j,l}/R^{\alpha})}{\partial R}\right|^2 , 
\end{equation}
where $C_0 = C_0(n,m,q,M_0,N_0,C_e,C_h,\alpha) \in (0,\infty)$ is a constant.  Adding $C_0$ times the left-hand side of \eqref{hole filling eqn3} to both sides of \eqref{hole filling eqn3} 
\begin{equation}\label{hole filling eqn4} 
	\int_{B_{\rho/8}(0)} \sum_{j=1}^J \sum_{l=1}^{q_j} R^{2-n} \left|\frac{\partial (w_{j,l}/R^{\alpha})}{\partial R}\right|^2 
	\leq \gamma \int_{B_{\rho/2}(0)} \sum_{j=1}^J \sum_{l=1}^{q_j} R^{2-n} \left|\frac{\partial (w_{j,l}/R^{\alpha})}{\partial R}\right|^2 
\end{equation}
for all $\rho \in [\vartheta_1/4,1]$, where $\gamma = C_0/(1+C_0) \in (0,1)$.  Let $N$ be a positive integer such that $2^{-2N-3} \vartheta_2 < \vartheta_1 \leq 2^{-2N-1} \vartheta_2$.  By iteratively applying \eqref{hole filling eqn4} with $\rho = 4^{-i} \vartheta_2$ for $i = 0,1,2,\ldots,N-1$ 
\begin{equation}\label{hole filling eqn5} 
	\int_{B_{\vartheta_1}(0)} \sum_{j=1}^J \sum_{l=1}^{q_j} R^{2-n} \left|\frac{\partial (w_{j,l}/R^{\alpha})}{\partial R}\right|^2 
	\leq C\bigg(\frac{\vartheta_1}{\vartheta_2}\bigg)^{2\mu} \int_{B_{\vartheta_2/2}(0) \setminus B_{\vartheta_2/8}(0)} \sum_{j=1}^J \sum_{l=1}^{q_j} 
		R^{2-n} \left|\frac{\partial (w_{j,l}/R^{\alpha})}{\partial R}\right|^2 
\end{equation}
where $\mu = -\log\gamma/\log 16$ and $C = C(n,m,q,M_0,N_0,C_e,C_h,\alpha) \in (0,\infty)$ is a constant.  Combining \eqref{hole filling eqn2} with $\rho = 2\vartheta_1$, \eqref{hole filling eqn5}, and \eqref{hole filling eqn1} with $\rho = \vartheta_2$, we obtain \eqref{hole filling concl}.
\end{proof}

\begin{theorem}\label{blowup decay thm} 
Let $(m_j,\varphi^{(0)}_j,\varphi^{(\infty)}_j) \in \mathfrak{D}$ and $w = (w_j) \in \mathfrak{B}(m_j,\varphi^{(0)}_j,\varphi^{(\infty)}_j)$.  There exists a unique $\psi = (\psi_j) \in \widetilde{\mathfrak L}(m_j,\varphi^{(0)}_j,\varphi^{(\infty)}_j)$ such that 
\begin{equation}\label{blowup decay concl} 
	\frac{1}{\rho^{n+2\alpha}} \int_{B_{\rho}(0)} \sum_{j=1}^J \sum_{l=1}^{q_j} \mathcal{G}(w_{j,l},\psi_{j,l})^2 
	\leq C \rho^{2\mu} \int_{B_{1/2}(0)} \sum_{j=1}^J \sum_{l=1}^{q_j} |w_{j,l}|^2 
\end{equation}
where $w_{j,l}$ is as in \eqref{w localized} and $\psi_{j,l}$ is as in \eqref{psi localized}, and $\mu \in (0,1)$ and $C \in (0,\infty)$ are constants depending only on $n,m,q,M_0,N_0,C_e,C_h,\alpha$.  Additionally, $\|\psi\|_{L^2(B_1(0))} \leq C$ for some constant $C = C(n,m,q,M_0,N_0,C_e,C_h,\alpha) \in (0,\infty)$.
\end{theorem}

\begin{proof}
Immediate consequence of Lemma~\ref{hole filling lemma}; see~\cite[Theorem~9.4]{KrumWicc} for details.
\end{proof}

\section{Decay of fine excess of area-minimizing currents relative to homogeneous Dirichlet energy minimizing multi-valued functions}\label{sec:fine excess decay} 

We now our main excess decay lemma, Theorem~\ref{main excess decay thm}, which plays a key role in the proof of the uniqueness of tangent functions.  First in Lemma~\ref{main excess decay lemma1}, we prove the excess decay lemma assuming that $\varphi$ and $T$ satisfy Hypothesis~$(\star)$ and Hypothesis~$(\star\star)$ (of Section~\ref{sec:graph_rep_results}).  Using an argument of~\cite{Wic14}, in Lemma~\ref{main excess decay lemma2} and in Theorem~\ref{main excess decay thm}, we replace Hypothesis~$(\star\star)$ with the weaker hypothesis \eqref{main excess decay2 hyp} that the $T$ is much closer to $\varphi$ (in $L^2$) than $T$ is to the center manifold $\mathcal{M}$ provided we allow conclusion~(B) (including the improvement of fine excess) to hold at one of a fixed number of smaller scales $\vartheta_1,\vartheta_2,\ldots,\vartheta_{\lfloor q/q_0 \rfloor}$.  For the purposes of obtaining uniform decay estimates for $T$ relative to a tangent function in the next section, having a finite number of scales to choose from at each step of iterating Theorem~\ref{main excess decay thm} is as good as using a single scale at every step.

\begin{lemma}\label{main excess decay lemma1}
Let $\alpha = \ell_0/q_0 \geq 2$ for some relatively prime positive integers $\ell_0$ and $q_0 \leq q$.  Let $\varsigma = \varsigma(m,n,q) \in (0,1)$ be as in Theorem~\ref{relharm separation main thm}.  Let $M_0,N_0,C_e,C_h,\varepsilon_2$ be as in Hypothesis~\ref{noncon constants hypothesis}.  For each $\vartheta \in (0,1/4)$ there exist $\overline{\varepsilon}, \overline{\beta} \in (0,1)$ and $\overline{\delta} \in (0,1/2)$ depending only on $n,m,q,p,M_0,N_0,C_e,C_h,\alpha,\vartheta$ such that if Hypothesis~$(\star)$ (of Section~\ref{sec:graph_rep_results})   holds true with $\overline{\varepsilon}$ in place of $\varepsilon_0$, $\mathbf{I}_{N,Z_0}(0^+) \geq \alpha$, and 
\begin{equation}\label{main excess decay1 hyp}
	E_{\varsigma/2,\overline{\varepsilon}}(N,\varphi,B^{\mathcal M}_{\rho_0}(Z_0)) 
	\leq \overline{\beta} \inf_{\varphi' \in \bigcup_{p'=0}^{p-1} \Phi_{\alpha,q,p'}(P_{Z_0},S)} 
		E_{\varsigma/2,\overline{\varepsilon}}(N,\varphi',B^{\mathcal M}_{\rho_0}(Z_0)) , 
\end{equation}
(cf.~Hypothesis~$(\star\star)$), then either:
\begin{enumerate}[itemsep=2mm,topsep=0mm]
	\item[(A)]  for some $\zeta \in S \cap B_{\rho_0/2}(0,P_{Z_0})$ 
	\begin{equation*} 
		B^{\mathcal M}_{\overline{\delta}\rho_0}(\exp_{Z_0}(\zeta)) \cap B^{\mathcal M}_{\rho_0/2}(Z_0) 
			\cap \{ X \in \mathcal{B}^{(\geq 2)}_q : \mathbf{I}_{N,X}(0^+) \geq \alpha \} = \emptyset 
	\end{equation*}
	where $\exp_{Z_0}$ is the exponential map of $\mathcal{M}$ at $Z_0$ or 
	
	\item[(B)]  there exists a rotation $\Gamma$ of $P_{Z_0}$ and $\widetilde{\varphi} \in \Phi_{\alpha,q}(P_{Z_0},\Gamma(S))$ such that: 
	\begin{gather} 
		\label{main excess decay1 concl1} \|\Gamma - I_{P_{Z_0}}\| \leq \frac{\overline{\kappa} E_{\varsigma/2,\overline{\varepsilon}}(N,\varphi,
			B^{\mathcal M}_{\rho_0}(Z_0))}{\rho_0^{-n/2-1} \|N\|_{L^2(B^{\mathcal M}_{\rho_0}(Z_0))})} , \\
		\label{main excess decay1 concl2} \rho_0^{\alpha-1} \|\mathcal{G}(\varphi,\widetilde{\varphi})\|_{L^2(B_1(0,P_{Z_0}))} 
			\leq \overline{C}_0 E_{\varsigma/2,\overline{\varepsilon}}(N,\varphi,B^{\mathcal M}_{\rho_0}(Z_0)) , \\
		\label{main excess decay1 concl3} E_{\varsigma/2,\overline{\varepsilon}}(N,\widetilde{\varphi},B^{\mathcal M}_{\vartheta\rho_0}(Z_0)) 
			\leq \overline{\nu} \vartheta^{\alpha-1+\mu} E_{\varsigma/2,\overline{\varepsilon}}(N,\varphi,B^{\mathcal M}_{\rho_0}(Z_0)) .
	\end{gather}
	Here $I_{P_{Z_0}}$ is the identity map on $P_{Z_0}$ and $\mu \in (0,1)$ and $\overline{\kappa},\overline{C}_0,\overline{\nu} \in (0,\infty)$ are constants depending only on $n,m,q,M_0,N_0,C_e,C_h,\alpha$.
\end{enumerate}
\end{lemma} 

\begin{proof}
Fix $\vartheta \in (0,1/4)$.  For $\nu = 1,2,3,\ldots$ let $\varepsilon_{\nu},\beta_{\nu},\delta_{\nu} \in (0,1)$ with $\varepsilon_{\nu}\rightarrow 0$, $\beta_{\nu}\rightarrow 0$, and $\delta_{\nu}\rightarrow 0$ as $\nu\rightarrow\infty$ and let $T^{(\nu)}, \mathcal{M}^{(\nu)}, N^{(\nu)}, \varphi^{(\nu)}, Z_{\nu}, \rho_{\nu}, P_{\nu}, S_{\nu}$ such that Hypothesis~$(\star)$ and \eqref{main excess decay1 hyp} hold true with $\varepsilon_{\nu},\beta_{\nu},T^{(\nu)}, \mathcal{M}^{(\nu)}, N^{(\nu)}, \varphi^{(\nu)}, Z_{\nu}, \rho_{\nu}, P_{\nu}, S_{\nu}$ in place of $\varepsilon_0,\beta_0,T, \mathcal{M}, N, \varphi, Z_0, \rho_0, P_{Z_0}, S$ and $\mathbf{I}_{N^{(\nu)},Z_{\nu}}(0^+) \geq \alpha$ but conclusion~(A) does not hold true with $\delta_{\nu}, T^{(\nu)}, \mathcal{M}^{(\nu)}, N^{(\nu)}, Z_{\nu}, \rho_{\nu}, P_{\nu}, S_{\nu}$ in place of $\overline{\delta},T, \mathcal{M}, N, Z_0, \rho_0, P_{Z_0}, S$.  Note that after passing to a subsequence, we may assume that the number $p$ of distinct non-zero components of $\varphi^{(\nu)}$ is independent of $\nu$.  To prove the lemma, it suffices to show that for infinitely many $\nu$ there exists an orthogonal rotation $\Gamma_{\nu}$ of $\mathbb{R}^n$ and $\widetilde{\varphi}^{(\nu)} \in \Phi_{\alpha,q}(P_{\nu},S_{\nu})$ such that conclusion~(B) holds true with $T^{(\nu)}, \mathcal{M}^{(\nu)}, N^{(\nu)}, \Gamma_{\nu}, \varphi^{(\nu)}, \widetilde{\varphi}^{(\nu)}, Z_{\nu}, \rho_{\nu}, P_{\nu},S_{\nu}$ in place of $T, \mathcal{M}, N, \Gamma, \varphi, \widetilde{\varphi}, Z_0, \rho_0, P_{Z_0}, S$ and with fixed constants $\mu = \mu(n,m,q) \in (0,1)$ and $\overline{\kappa},\overline{\nu},\overline{C}_0 \in (0,\infty)$ depending only on $n,m,q,M_0,N_0,C_e,C_h,\alpha$.

After passing to a subsequence, let $w = (w_j) \in \mathfrak{B}$ be the fine blow-up of $(T^{(\nu)})$ relative to $(\varphi^{(\nu)})$ (as in Subsection~\ref{sec:blowup subsec} and Definition~\ref{blowup class defn}).  Let $w \in \mathfrak{B}(m_j,\varphi^{(0)},\varphi^{(\infty)}_j)$, where $(m_j,\varphi^{(0)},\varphi^{(\infty)}_j) \in \mathfrak{D}$ is as in Remark~\ref{blowup on domain rmk}.  By Theorem~\ref{blowup decay thm}, there exists (a unique) $\psi = (\psi_j) \in \widetilde{\mathfrak L}(m_j,\varphi^{(0)},\varphi^{(\infty)}_j)$ such that \eqref{blowup decay concl} holds true and $\|\psi\|_{L^2(B_1(0))} \leq C$, where $\mu \in (0,1)$ and $C \in (0,\infty)$ are constants depending only on $n,m,q,M_0,N_0,C_e,C_h,\alpha$.  Let $\Gamma_{\nu}$ and $\widetilde{\varphi}^{(\nu)} \in \Phi_{\alpha,q}(P_{\nu},\Gamma_{\nu}(S_{\nu}))$ be as in Lemma~\ref{integrability lemma} corresponding to $\psi$ and $(\varphi^{(\nu)})$.  By Lemma~\ref{integrability lemma} and $\|\psi\|_{L^2(B_1(0))} \leq C$, we immediately have that \eqref{main excess decay1 concl1} and \eqref{main excess decay1 concl2} hold true with $\mathcal{M}^{(\nu)},N^{(\nu)},\Gamma_{\nu},\varphi^{(\nu)},\widetilde{\varphi}^{(\nu)},Z_{\nu},\rho_{\nu},P_{\nu}$ in place of $\mathcal{M},N,\Gamma,\varphi,\widetilde{\varphi},Z_0,\rho_0,P_{Z_0}$ and with constants $\overline{\kappa},\overline{C}_0 \in (0,\infty)$ depending only on $n,m,q,M_0,N_0,C_e,C_h,\alpha$.  By Lemma~\ref{blowup norms lemma} and \eqref{blowup decay concl} with $\rho = \vartheta$, \eqref{main excess decay1 concl3} hold true with $\mathcal{M}^{(\nu)},N^{(\nu)},\varphi^{(\nu)},\widetilde{\varphi}^{(\nu)},Z_{\nu},\rho_{\nu}$ in place of $\mathcal{M},N,\varphi,\widetilde{\varphi},Z_0,\rho_0$ and for some constant $\overline{\nu} \in (0,\infty)$ depending only on $n,m,q,M_0,N_0,C_e,C_h,\alpha$.  
\end{proof}

\begin{lemma}\label{main excess decay lemma2}
Let $\alpha = \ell_0/q_0 \geq 2$ for some relatively prime positive integers $\ell_0$ and $q_0 \leq q$.  Let $\varsigma = \varsigma(m,n,q) \in (0,1)$ be as in Theorem~\ref{relharm separation main thm}.  Let $M_0,N_0,C_e,C_h,\varepsilon_2$ be as in Hypothesis~\ref{noncon constants hypothesis}.  For each $p \in \{1,2,\ldots,\lfloor q/q_0 \rfloor\}$ and $\vartheta_1,\vartheta_2,\ldots,\vartheta_p \in (0,1/4)$ there exists $\varepsilon^{(p)}, \beta^{(p)} \in (0,1)$ and $\delta^{(p)} \in (0,1/2)$ depending only on $n,m,q,M_0,N_0,C_e,C_h,\alpha,\vartheta_1,\ldots,\vartheta_p$ such that if Hypothesis~$(\star)$ (of Section~\ref{sec:graph_rep_results}) holds true with $\varepsilon^{(p)}$ in place of $\varepsilon_0$, $\mathbf{I}_{N,Z_0}(0^+) \geq \alpha$, and 
\begin{equation}\label{main excess decay2 hyp}
	E_{\varsigma/2,\varepsilon^{(p)}}(N,\varphi,B^{\mathcal M}_{\rho_0}(Z_0)) \leq \beta^{(p)} \rho_0^{-n/2-1} \|N\|_{L^2(B^{\mathcal M}_{\rho_0}(Z_0))} , 
\end{equation}
then either:  
\begin{enumerate}[itemsep=2mm,topsep=0mm]
	\item[(A)]  for some $\zeta \in S \cap B_{\rho_0/2}(0,P_{Z_0}),$ 
	\begin{equation*} 
		B^{\mathcal M}_{\delta^{(p)}\rho_0}(\exp_{Z_0}(\zeta)) \cap B^{\mathcal M}_{\rho_0/2}(Z_0) 
			\cap \{ X \in \mathcal{B}^{(\geq 2)}_q : \mathbf{I}_{N,X}(0^+) \geq \alpha \} = \emptyset ; \;or,  
	\end{equation*}
	
	\item[(B)]  there exists a rotation $\Gamma$ of $P_{Z_0}$, $\widetilde{\varphi} \in \Phi_{\alpha,q}(P_{Z_0},\Gamma(S))$, and $j \in \{1,2,\ldots,p\}$ such that: 
	\begin{gather*}
		\|\Gamma - I_{P_{Z_0}}\| \leq \frac{\kappa^{(p)} E_{\varsigma/2,\varepsilon^{(p)}}(N,\varphi,B^{\mathcal M}_{\rho_0}(Z_0))}{
			\rho_0^{-n/2-1} \|N\|_{L^2(B^{\mathcal M}_{\rho_0}(Z_0))})} , \\
		\rho_0^{\alpha-1} \|\mathcal{G}(\varphi,\widetilde{\varphi})\|_{L^2(B_1(0,P_{Z_0}))} 
			\leq C_0^{(p)} E_{\varsigma/2,\varepsilon^{(p)}}(N,\varphi,B^{\mathcal M}_{\rho_0}(Z_0)) , \; and\\
		E_{\varsigma/2,\varepsilon^{(p)}}(N,\widetilde{\varphi},B^{\mathcal M}_{\vartheta_j\rho_0}(Z_0)) 
			\leq \nu_j^{(p)} \vartheta_j^{\alpha-1+\mu} E_{\varsigma/2,\varepsilon^{(p)}}(N,\varphi,B^{\mathcal M}_{\rho_0}(Z_0)).
	\end{gather*}
	Here $\mu = \mu(n,m,q,M_0,N_0,C_e,C_h,\alpha) \in (0,1)$ is as in Lemma~\ref{main excess decay lemma1}; $\kappa^{(p)}, C_0^{(p)} \in (0,\infty)$ depend only on $n,m,q,M_0,N_0,C_e,C_h,\alpha$ if $p = 1$ and on $n,m,q,M_0,N_0,C_e,C_h,\alpha,\vartheta_1,\ldots,\vartheta_{p-1}$ if $p > 1$; $\nu_j^{(p)} \in (0,\infty)$ depends only on $n,m,q,M_0,N_0,C_e,C_h,\alpha$ if $j = 1$ and on $n,m,q,M_0,$ $N_0,C_e,C_h,\alpha,\vartheta_1,\ldots,\vartheta_{j-1}$ if $j > 1$.
\end{enumerate}
\end{lemma}

\begin{proof} 
Throughout this proof, we let $\overline{\varepsilon} = \overline{\varepsilon}(\vartheta)$, $\overline{\beta} = \overline{\beta}(\vartheta)$, $\overline{\delta} = \overline{\delta}(\vartheta)$, $\overline{\kappa}$, $\overline{C}_0$, and $\overline{\nu}$ be as in Lemma~\ref{main excess decay lemma1}, omitting the dependence on $n,m,q,M_0,N_0,C_e,C_h,\alpha$.  Similarly, $\varepsilon^{(p)} = \varepsilon^{(p)}(\vartheta_1,\ldots,\vartheta_p)$, $\beta^{(p)} = \beta^{(p)}(\vartheta_1,\ldots,\vartheta_p)$, $\delta^{(p)} = \delta^{(p)}(\vartheta_1,\ldots,\vartheta_p)$, $\kappa^{(p)} = \kappa^{(p)}(\vartheta_1,\ldots,\vartheta_{p-1})$, $C_0^{(p)} = C_0^{(p)}(\vartheta_1,\ldots,\vartheta_{p-1})$, and $\nu_j^{(p)} = \nu_j^{(p)}(\vartheta_1,\ldots,\vartheta_{j-1})$ are as in Lemma~\ref{main excess decay lemma2}.

We proceed by induction on $p$.  In the base case $p = 1$, set $\varepsilon^{(1)} = \overline{\varepsilon}(\vartheta_1)$, $\beta^{(1)} = \overline{\beta}(\vartheta_1)$, $\delta^{(1)} = \overline{\delta}(\vartheta_1)$, $\kappa^{(1)} = \overline{\kappa}$, $C_0^{(1)} = \overline{C}_0$, $\nu_1^{(1)} = \overline{\nu}$.  Note that this proves lemma in the case $q_0 = q$; henceforth, we shall assume that $q_0 < q$.

Suppose that $p_0 \in \{2,3,\ldots,\lfloor q/q_0 \rfloor\}$ such that Lemma~\ref{main excess decay lemma1} holds true whenever $p \in \{1,2,\ldots,p_0-1\}$.  Let $\vartheta_1,\vartheta_2,\ldots,\vartheta_{p_0-1} \in (0,1/4)$.  By induction, we obtain parameters $\varepsilon^{(p)},\beta^{(p)},\delta^{(p)},\kappa^{(p)},C_0^{(p)},$ $C_2^{(p)},\nu^{(p)}_j$ as in Lemma~\ref{main excess decay lemma2} for each $p \in \{1,2,\ldots,p_0-1\}$.  We claim that Lemma~\ref{main excess decay lemma2} holds true with $p = p_0$ and with the constants 
\begin{align*}
	\varepsilon^{(p_0)} =\,& \min\{\overline{\varepsilon}(\vartheta_1)\} \cup \{\varepsilon^{(p)}(\vartheta_2,\ldots,\vartheta_p) 
		: p = 1,\ldots,p_0-1 \} , \\
	\beta^{(p_0)} =\,& \min\{\overline{\beta}(\vartheta_1)\} \cup \{ 2^{-p_0} \overline{\beta}(\vartheta_1)^{p_0} 
		\beta^{(p)}(\vartheta_2,\ldots,\vartheta_p) : p = 1,\ldots,p_0-1 \} , \\
	\delta^{(p_0)} =\,& \min\{\overline{\delta}(\vartheta_1)\} \cup \{\delta^{(p)}(\vartheta_2,\ldots,\vartheta_p) : p = 1,\ldots,p_0-1 \} , \\
	\kappa^{(p_0)} =\,& \max\{\overline{\kappa}\} \cup \{ 2^{p_0} \overline{\beta}(\vartheta_1)^{-p_0} 
		\kappa^{(p)}(\vartheta_2,\ldots,\vartheta_p) : p = 1,\ldots,p_0-1 \} , \\
	C_0^{(p_0)} =\,& \max\{\overline{C}_0\} \cup \{ \overline{C} + 2^{p_0} \overline{\beta}(\vartheta_1)^{-p_0} 
		(\overline{C} + C_0^{(p)}(\vartheta_2,\ldots,\vartheta_p)) : p = 1,\ldots,p_0-1 \} , \\
	\nu_j^{(p_0)} =\,& \max\{\overline{\nu}\} \cup \{ 2^{p_0} \overline{\beta}(\vartheta_1)^{-p_0} \nu_j^{(p)}(\vartheta_1,\ldots,\vartheta_{j-1}) 
		: p = 1,\ldots,p_0-1 \} ,
\end{align*}
where $\overline{\beta} = \overline{\beta}(n,m,q,\alpha,\vartheta_1)$ is as in Lemma~\ref{main excess decay lemma1}, $\beta_0$ is as in Theorem~\ref{graphical_thm} with $\tau = 1/16$, and $\overline{C} = \overline{C}(n,\alpha)$ is as in \eqref{main excess decay2 eqn6} below.  To see this, let $T, \mathcal{M}, N, \varphi, Z_0, \rho_0, P_{Z_0}, S$ satisfy Hypothesis~$(\star)$ and \eqref{main excess decay2 hyp} with $p_0, \varepsilon^{(p_0)}, \beta^{(p_0)}$ in place of $p,\varepsilon_0,\beta_0$
\begin{equation}\label{main excess decay2 eqn1}
	E_{\varsigma/2}(N,\varphi,B^{\mathcal M}_{\rho_0}(Z_0)) \leq \beta^{(p_0)} \rho_0^{-n/2-1} \|N\|_{L^2(B^{\mathcal M}_{\rho_0}(Z_0))} .
\end{equation}
If additionally \eqref{main excess decay1 hyp} holds true with $p = p_0$, $\overline{\varepsilon} = \overline{\varepsilon}(\vartheta_1)$, and $\overline{\beta} = \overline{\beta}(\vartheta_1)$, then the conclusion of the present lemma follows directly from Lemma~\ref{main excess decay lemma1} with $\vartheta = \vartheta_1$.  If on the other hand 
\begin{equation}\label{main excess decay2 eqn2}
	E_{\varsigma/2,\overline{\varepsilon}(\vartheta_1)}(N,\varphi,B^{\mathcal M}_{\rho_0}(Z_0)) 
	\geq \overline{\beta}(\vartheta_1) \inf_{\varphi' \in \bigcup_{p'=1}^{p_0-1} \Phi_{\alpha,q,p'}} 
		E_{\varsigma/2,\overline{\varepsilon}(\vartheta_1)}(N,\varphi',B^{\mathcal M}_{\rho_0}(Z_0)) ,
\end{equation}
then we can choose an integer $\overline{p} \in \{0,1,2,\ldots,p_0-1\}$ and $\overline{\varphi} \in \Phi_{\alpha,q,\overline{p}}$ such that 
\begin{equation}\label{main excess decay2 eqn3}
	E_{\varsigma/2,\overline{\varepsilon}(\vartheta_1)}(N,\overline{\varphi},B^{\mathcal M}_{\rho_0}(Z_0)) 
	\leq 2^{p_0} \overline{\beta}(\vartheta_1)^{1-p_0} \inf_{\varphi' \in \bigcup_{p'=0}^{p_0-1} \Phi_{\alpha,q,p'}} 
		E_{\varsigma/2,\overline{\varepsilon}(\vartheta_1)}(N,\varphi',B^{\mathcal M}_{\rho_0}(Z_0)) 
\end{equation}
and either $\overline{p} = 0$ or $\overline{p} > 0$ and 
\begin{equation}\label{main excess decay2 eqn4} 
	E_{\varsigma/2,\overline{\varepsilon}(\vartheta_1)}(N,\overline{\varphi},B^{\mathcal M}_{\rho_0}(Z_0)) 
	\leq \overline{\beta}(\vartheta_1) \inf_{\varphi' \in \bigcup_{p'=0}^{\overline{p}-1} \Phi_{\alpha,q,p'}} 
		E_{\varsigma/2,\overline{\varepsilon}(\vartheta_1)}(N,\varphi',B^{\mathcal M}_{\rho_0}(Z_0)) .
\end{equation}
By \eqref{main excess decay2 eqn2} and \eqref{main excess decay2 eqn3}
\begin{equation}\label{main excess decay2 eqn5} 
	E_{\varsigma/2,\overline{\varepsilon}(\vartheta_1)}(N,\overline{\varphi},B^{\mathcal M}_{\rho_0}(Z_0)) 
	\leq 2^{p_0} \overline{\beta}(\vartheta_1)^{-p_0} E_{\varsigma/2,\overline{\varepsilon}(\vartheta_1)}(N,\varphi,B^{\mathcal M}_{\rho_0}(Z_0)) . 
\end{equation}
Provided $\varepsilon^{(p_0)}$ is sufficiently small, by the triangle inequality  
\begin{equation}\label{main excess decay2 eqn6} 
	\rho_0^{\alpha-1} \|\mathcal{G}(\varphi,\overline{\varphi})\|_{L^2(B_1(0,P_{Z_0}))} 
	\leq \overline{C} E_{\varsigma/2,\overline{\varepsilon}(\vartheta_1)}(N,\varphi,B^{\mathcal M}_{\rho_0}(Z_0)) 
		+ \overline{C} E_{\varsigma/2,\overline{\varepsilon}(\vartheta_1)}(N,\overline{\varphi},B^{\mathcal M}_{\rho_0}(Z_0)) 
\end{equation}
for some constant $\overline{C} = \overline{C}(n,\alpha) \in (0,\infty)$.  By \eqref{main excess decay2 eqn5}, \eqref{main excess decay2 eqn1}, and the definition of $\beta^{(p_0)}$ 
\begin{align}\label{main excess decay2 eqn7} 
	E_{\varsigma/2,\overline{\varepsilon}(\vartheta_1)}(N,\overline{\varphi},B^{\mathcal M}_{\rho_0}(Z_0)) 
	\leq\,& 2^{p_0} \overline{\beta}(\vartheta_1)^{-p_0} E_{\varsigma/2,\overline{\varepsilon}(\vartheta_1)}(N,\varphi,B^{\mathcal M}_{\rho_0}(Z_0))
	\\ \leq\,& 2^{p_0} \overline{\beta}(\vartheta_1)^{-p_0} \beta^{(p_0)} \rho_0^{-n/2-1} \|N\|_{L^2(B^{\mathcal M}_{\rho_0}(Z_0))} \nonumber
	\\ \leq\,& \beta^{(\overline{p})} \rho_0^{-n/2-1} \|N\|_{L^2(B^{\mathcal M}_{\rho_0}(Z_0))} . \nonumber 
\end{align}
Noting that $\beta^{(\overline{p})} \in (0,1)$, we must have that $p_0 \geq 1$; otherwise, $p_0 = 0$ implies that $\overline{\varphi} = q \llbracket 0 \rrbracket$ contradicting \eqref{main excess decay2 eqn7}.  In light of \eqref{main excess decay2 eqn7}, we can apply the induction hypothesis to deduce that the conclusion of Lemma~\ref{main excess decay lemma2} holds true with $p = \overline{p}$ and $\overline{\varphi}$ in place of $\varphi$.  It follows using \eqref{main excess decay2 eqn5} and \eqref{main excess decay2 eqn6} that the conclusion of Lemma~\ref{main excess decay lemma2} holds true.
\end{proof}

\begin{theorem}\label{main excess decay thm}
Let $\alpha = \ell_0/q_0 \geq 2$ for some relatively prime positive integers $\ell_0$ and $q_0 \leq q$.  Let $\varsigma = \varsigma(m,n,q) \in (0,1)$ be as in Theorem~\ref{relharm separation main thm}.  Let $M_0,N_0,C_e,C_h,\varepsilon_2$ be as in Hypothesis~\ref{noncon constants hypothesis}.  For each $\vartheta_1,\vartheta_2,\ldots,\vartheta_{\lceil q/q_0 \rceil} \in (0,1/4)$ there exists $\varepsilon,\beta \in (0,1)$ and $\delta \in (0,1/2)$ depending only on $n,m,q,M_0,N_0,C_e,C_h,\alpha,\vartheta_1,\ldots,\vartheta_{\lceil q/q_0 \rceil}$ such that if Hypothesis~$(\star)$ (of Section~\ref{sec:graph_rep_results}) holds true with $\varepsilon$ in place of $\varepsilon_0$, $\mathbf{I}_{N,Z_0}(0^+) \geq \alpha$, and 
\begin{equation}\label{main excess decay3 hyp}
	E_{\varsigma/2,\varepsilon}(N,\varphi,B^{\mathcal M}_{\rho_0}(Z_0)) \leq \beta \rho_0^{-n/2-1} \|N\|_{L^2(B^{\mathcal M}_{\rho_0}(Z_0))} , 
\end{equation}
then either:  
\begin{enumerate}[itemsep=2mm,topsep=0mm]
	\item[(A)]  for some $\zeta \in S \cap B_{\rho_0/2}(0,P_{Z_0}),$ 
	\begin{equation*} 
		B^{\mathcal M}_{\delta\rho_0}(\exp_{Z_0}(\zeta)) \cap B^{\mathcal M}_{\rho_0/2}(Z_0) 
			\cap \{ X \in \mathcal{B}^{(\geq 2)}_q : \mathbf{I}_{N,X}(0^+) \geq \alpha \} = \emptyset ; \; or, 
	\end{equation*}
	
	\item[(B)]  there exists a rotation $\Gamma$ of $P_{Z_0}$, $\widetilde{\varphi} \in \Phi_{\alpha,q}(P_{Z_0},\Gamma(S))$, and $j \in \{1,2,\ldots,p\}$ such that: 
	\begin{gather} 
		\label{main excess decay3 concl1} \|\Gamma - I_{P_{Z_0}}\| \leq \frac{\kappa \,E_{\varsigma/2,\varepsilon}(N,\varphi,
			B^{\mathcal M}_{\rho_0}(Z_0))}{\rho_0^{-n/2-1} \|N\|_{L^2(B^{\mathcal M}_{\rho_0}(Z_0))})} , \\
		\label{main excess decay3 concl2} \rho_0^{\alpha-1} \|\mathcal{G}(\varphi,\widetilde{\varphi})\|_{L^2(B_1(0,P_{Z_0}))} 
			\leq C_0 E_{\varsigma/2,\varepsilon}(N,\varphi,B^{\mathcal M}_{\rho_0}(Z_0)) , \; and\\
		\label{main excess decay3 concl3} E_{\varsigma/2,\varepsilon}(N,\widetilde{\varphi},B^{\mathcal M}_{\vartheta_j\rho_0}(Z_0)) 
			\leq \nu_j \vartheta_j^{\alpha-1+\mu} E_{\varsigma/2,\varepsilon}(N,\varphi,B^{\mathcal M}_{\rho_0}(Z_0)) .
	\end{gather}
	Here $\mu = \mu(n,m,q,M_0,N_0,C_e,C_h,\alpha) \in (0,1)$ is as in Lemma~\ref{main excess decay lemma1}; $\kappa, C_0 \in (0,\infty)$ depend only on $n,m,q,M_0,N_0,C_e,C_h,\alpha$ if $p = 1$ and on $n,m,q,M_0,N_0,C_e,C_h,\alpha,\vartheta_1,\ldots,\vartheta_{\lceil q/q_0 \rceil-1}$ if $p > 1$; $\nu_j \in (0,\infty)$ depends only on $n,m,q,M_0,N_0,C_e,C_h,\alpha$ if $j = 1$ and on $n,m,q,M_0,$ $N_0,C_e,C_h,\alpha,\vartheta_1,\ldots,\vartheta_{j-1}$ if $j > 1$.
\end{enumerate}
\end{theorem}

\begin{proof}
As in the proof of Lemma~\ref{main excess decay lemma2}, we let $\varepsilon^{(p)} = \varepsilon^{(p)}(\vartheta_1,\ldots,\vartheta_p)$, $\beta^{(p)} = \beta^{(p)}(\vartheta_1,\ldots,\vartheta_p)$, $\delta^{(p)} = \delta^{(p)}(\vartheta_1,\ldots,\vartheta_p)$, $\kappa^{(p)} = \kappa^{(p)}(\vartheta_1,\ldots,\vartheta_{p-1})$, $C_0^{(p)} = C_0^{(p)}(\vartheta_1,\ldots,\vartheta_{p-1})$, and $\nu_j^{(p)} = \nu_j^{(p)}(\vartheta_1,\ldots,\vartheta_{j-1})$ be as in Lemma~\ref{main excess decay lemma2}, omitting the dependence on $n,m,q,M_0,N_0,C_e,C_h,\alpha$.  Set 
\begin{align*}
	\varepsilon =\,& \min\{\varepsilon^{(p)}(\vartheta_1,\ldots,\vartheta_p) : p = 1,2,\ldots,\lfloor q/q_0 \rfloor\} , \\
	\beta =\,& \min\{\beta^{(p)}(\vartheta_1,\ldots,\vartheta_p) : p = 1,2,\ldots,\lfloor q/q_0 \rfloor\} , \\
	\delta =\,& \min\{\delta^{(p)}(\vartheta_1,\ldots,\vartheta_p) : p = 1,2,\ldots,\lfloor q/q_0 \rfloor\} 
\end{align*}
Set 
\begin{align*}
	\kappa =\,& \max\{\kappa^{(1)}\} \cup \{\kappa^{(p)}(\vartheta_1,\ldots,\vartheta_{p-1}) : p = 2,3,\ldots,\lfloor q/q_0 \rfloor\} , \\
	C_0 =\,& \max\{C_0^{(1)}\} \cup \{C_0^{(p)}(\vartheta_1,\ldots,\vartheta_{p-1}) : p = 2,3,\ldots,\lfloor q/q_0 \rfloor\} 
\end{align*}
if $q_0 < q$ and $\kappa = \kappa^{(1)}$ and $C_0 = C_0^{(1)}$ if $q_0 = q$.  Set 
\begin{equation*}
	\nu_j = \max\{\nu_j^{(p)}(\vartheta_1,\ldots,\vartheta_{j-1}) : p = 1,2,3,\ldots,\lfloor q/q_0 \rfloor\} 
\end{equation*}
if $j > 1$ and $\nu_1 = \max\{\nu_1^{(p)} : p = 1,2,3,\ldots,\lfloor q/q_0 \rfloor\}$.  With this choice of parameters, the conclusion of the present theorem follows directly from Lemma~\ref{main excess decay lemma2}.
\end{proof}

\section{Proofs of the main theorems}\label{sec:proof of main results}

We begin with the following two results Theorem~\ref{SandT thm} and Theorem~\ref{precise rectifiability thm}. 

\begin{theorem}\label{SandT thm}
Let $\alpha = \ell_0/q_0 \geq 2$ for some relatively prime positive integers $\ell_0$ and $q_0 \leq q$.  Let $\varsigma = \varsigma(m,n,q) \in (0,1)$ be as in Theorem~\ref{relharm separation main thm}.  Let $M_0,N_0,C_e,C_h,\varepsilon_2$ be as in Hypothesis~\ref{noncon constants hypothesis}.  There exists $\varepsilon, \beta, \eta, \overline{\mu} \in (0,1)$ depending only on $n,m,q,M_0,N_0,C_e,C_h,\alpha$ such that the following holds true.  Suppose that Hypothesis~$(\star)$ (of Section~\ref{sec:graph_rep_results}) and Hypothesis~$(\dagger)$ (of Section~\ref{sec:key_est_sec}) hold true with $\varepsilon,\beta$ in place of $\varepsilon_0,\beta_0$.  Then 
\begin{equation}\label{SandT concl1}
	\{ Z \in \mathcal{B}^{(\geq 2)}_q \cap B^{\mathcal M}_{\rho_0/2}(Z_0) : \mathbf{I}_{N,Z}(0^+) \geq \alpha \} = \mathcal{S} \cup \mathcal{T} 
\end{equation}
where $\mathcal{S} \subseteq L$ for some properly embedded $(n-2)$-dimensional $C^{1,\overline{\mu}}$-submanifold $L$ of $B^{\mathcal M}_{\rho_0/2}(Z_0)$ with $\mathcal{H}^{n-2}(L) \leq 2\omega_{n-2} (\rho_0/2)^{n-2}$ and $\mathcal{T} \subseteq \bigcup_{j=1}^{\infty} B^{\mathcal M}_{\rho_j}(Z_j)$ for a countable family of geodesic balls $\{B^{\mathcal M}_{\rho_j}(Z_j)\}_{j=1}^{\infty}$ with $\sum_{j=1}^{\infty} \rho_j^{n-2} \leq (1-\eta) \,\rho_0^{n-2}$.  Moreover, for $\mathcal{H}^{n-2}$-a.e.~$Z \in \mathcal{S}$, there exists a unique non-zero, homogeneous degree $\alpha$, locally Dirichlet energy minimizing function $\varphi^{(Z)} : P_Z \rightarrow \mathcal{A}_q(P_Z^{\perp})$ such that $\omega_n^{1/2} \|\varphi^{(Z)}\|_{L^2(B_1(0,P_Z))}^{-1} \,\varphi^{(Z)}$ is the unique tangent function to $T$ relative to the center manifold $\mathcal{M}$ at $Z$, $\dim\op{spine}\varphi^{(Z)} = n-2$, and 
\begin{equation}\label{SandT concl2}
	\int_{B^{\mathcal M}_{\rho}(Z)} \mathcal{G}(N,\varphi \circ \exp_{Z}^{-1})^2 \leq C \rho^{n+2\alpha+2\overline{\mu}} 
\end{equation}
for all $\rho \in (0,\rho_0/4]$, where $P_Z = {\rm T}_Z \mathcal{M}$ and $C = C(n,m,q,\alpha) \in (0,\infty)$ is a constant.
\end{theorem}

\begin{proof} 
For $j = 1,2,\ldots,\lfloor q/q_0 \rfloor$ inductively choose $\vartheta_j = \vartheta_j(n,m,q,M_0,N_0,C_e,C_h,\alpha) \in (0,1/4)$ such that $\vartheta_{j+1} \leq \vartheta_j/8$ if $1 \leq j < \lfloor q/q_0 \rfloor$ and $\nu_j \vartheta_j^{\mu/2} < 1$ for all $j$, where $\mu = \mu(n,m,q,M_0,N_0,C_e,C_h,\alpha)$, $\nu_1 = \nu_1(n,m,q,M_0,N_0,C_e,C_h,\alpha)$, and $\nu_j = \nu_1(n,m,q,\alpha,M_0,N_0,C_e,C_h,\vartheta_1,\ldots,\vartheta_{j-1})$ for $2 \leq j \leq \lfloor q/q_0 \rfloor$ are as in Theorem~\ref{main excess decay thm}.  Let $\delta = \delta(n,m,q,M_0,N_0,C_e,C_h,\alpha,\vartheta_1,\vartheta_2,\ldots,\vartheta_{\lfloor q/q_0 \rfloor})$ be as in Theorem~\ref{main excess decay thm}.  Let $\varepsilon = \varepsilon(n,m,q,M_0,N_0,C_e,C_h,\alpha) \in (0,1)$ and $\beta = \beta(n,m,q,M_0,N_0,C_e,C_h,\alpha) \in (0,1)$ be suitable small constants to be later determined.  

Set 
\begin{equation*}
	\mathcal{B}^{\star}_q = \{Z \in \mathcal{B}_q : \mathbf{I}_{N,Z}(0^+) \geq \alpha \} .
\end{equation*}
Suppose that $\mathcal{B}^{\star}_q \cap B^{\mathcal M}_{\delta\vartheta_{\lfloor q/q_0 \rfloor}\rho_0/2}(\exp_{Z_0}(\zeta)) = \emptyset$ for some $\zeta \in S \cap B_{\rho_0/2}(0,P_{Z_0})$.  Then by Theorem~\ref{nonconcentration thm}(A), \eqref{exp map est}, and Hypothesis~$(\dagger)$ 
\begin{equation}\label{SandT eqn1}
	\mathcal{B}^{\star}_q \cap B^{\mathcal M}_{\rho_0/2}(Z_0) 
	\subset \{ X \in \mathbb{R}^{n+m} : \op{dist}(X, Z_0 + S) < C \beta \rho_0 + C \mathbf{m}_0^{1/2} \rho_0^2 \} , 
\end{equation}
where $C = C(n,m,q,M_0,N_0,C_e,C_h,\alpha) \in (0,\infty)$.  Hence assuming $\varepsilon$ and $\beta$ are sufficiently small, the conclusion of the theorem trivially holds true with $\mathcal{S} = \emptyset$ and $\mathcal{T} = \mathcal{B}^{\star}_q \cap B^{\mathcal M}_{\rho_0}(Z_0)$.  Thus we may assume that $\mathcal{B}^{\star}_q \cap B^{\mathcal M}_{\delta\vartheta_{\lfloor q/q_0 \rfloor}\rho_0/2}(\exp_{Z_0}(\zeta)) \neq \emptyset$ for all $\zeta \in S \cap \mathbf{B}_{\rho_0/2}(0)$.

Let $Z \in \mathcal{B}^{\star}_q \cap B^{\mathcal M}_{\rho_0/2}(Z_0)$.  Let $\Gamma_Z$ be a rotation of $\mathbb{R}^{n+m}$ such that $\Gamma_Z(P_{Z_0}) = P_Z$ and \eqref{noncon rotate} holds true.  We claim that:  \textit{there exists $k_{\star}(Z) \in \{1,2,3,\ldots\} \cup \{\infty\}$ such that, setting $s_0 = \rho_0/4$, $\widetilde{\Gamma}_0 = \Gamma_Z$, and $\varphi_0 = \Gamma_Z \circ \varphi \circ \Gamma_Z^{-1}$:}
\begin{enumerate}[itemsep=2mm,topsep=0mm]
	\item[(i)]  \textit{for each positive integer $i < k_{\star}(Z)$ there exists $j(i) \in \{1,2,\ldots,\lceil q/q_0 \rceil\}$, a radius $s_i$ given by $s_i = \vartheta_{j(i)} s_{i-1}$, and rotation $\widetilde{\Gamma}_i$ of $\mathbb{R}^{n+m}$ with $\widetilde{\Gamma}_i(P_{Z_0}) = P_Z$ and $\widetilde{\Gamma}_i |_{P_{Z_0}} = \Gamma_Z |_{P_{Z_0}}$, and $\varphi_i \in \Phi_{\alpha,q}(P_Z,\widetilde{\Gamma}_i(S))$ such that 
	\begin{gather}
		\label{SandT eqn2} \|\widetilde{\Gamma}_{i-1} - \widetilde{\Gamma}_i\| \leq \frac{C E_{\varsigma/2,\varepsilon}(N,\varphi_{i-1},
			B^{\mathcal M}_{s_i}(Z))}{s_i^{-n/2-1} \|N\|_{L^2(B^{\mathcal M}_{s_i}(Z))})} , \\
		\label{SandT eqn3} s_i^{\alpha-1} \|\mathcal{G}(\varphi_{i-1},\varphi_i)\|_{L^2(B_1(0,P_Z))} 
			\leq C E_{\varsigma/2,\varepsilon}(N,\varphi_{i-1},B^{\mathcal M}_{s_i}(Z)) , \\
		\label{SandT eqn4} E_{\varsigma/2,\varepsilon}(N,\varphi_i,B^{\mathcal M}_{s_i}(Z)) 
			\leq \vartheta_{j(i)}^{\alpha-1+\mu/2} E_{\varsigma/2,\varepsilon}(N,\varphi_{i-1},B^{\mathcal M}_{s_{i-1}}(Z)) , 
	\end{gather}
	where $C = C(n,m,q,M_0,N_0,C_e,C_h,\alpha) \in (0,\infty)$;}
	
	\item[(ii)]  \textit{either $k_{\star}(Z) = \infty$ or $k_{\star}(Z) < \infty$ and 
	\begin{equation}\label{SandT eqn5} 
		\mathcal{B}^{\star}_q \cap B^{\mathcal M}_{s_{k_{\star}(Z)}/2}(Z) \cap B^{\mathcal M}_{\delta s_{k_{\star}(Z)}/2}(\exp_Z(\zeta)) = \emptyset 
			\text{ for some } \zeta \in \Gamma_Z(S) \cap \mathbf{B}_{s_{k_{\star}(Z)}/2}(0).
	\end{equation}}
\end{enumerate}
To see this, first observe that by \eqref{noncon triangle3} and Theorem~\ref{nonconcentration thm}(A) 
\begin{equation}\label{SandT eqn6} 
	E_{\varsigma/2,\varepsilon}(N,\varphi_0,B^{\mathcal M}_{\rho_0/4}(Z)) 
	\leq C E_{\varsigma/2,\varepsilon}(N,\varphi,B^{\mathcal M}_{\rho_0}(Z_0)) + C \widehat{E} |\xi|/\rho_0  
	\leq C E_{\varsigma/2,\varepsilon}(N,\varphi,B^{\mathcal M}_{\rho_0}(Z_0)) , 
\end{equation}
where $\widehat{E} = \rho_0^{-n/2-1} \|N\|_{L^2(B^{\mathcal M}_{\rho_0}(Z_0))}$ and $C = C(n,m,q,M_0,N_0,C_e,C_h,\alpha) \in (0,\infty)$ are constants.  Now we prove the claim via induction on $i$.  In the base case $i = 1$, recall that $\varphi_0 = \Gamma_Z^{-1} \circ \varphi \circ \Gamma_Z$.  By the triangle inequality, \eqref{SandT eqn6}, Hypothesis~$(\dagger)$, and Remark~\ref{graphical rmk}(2) 
\begin{align*}
	E_{\varsigma/2,\varepsilon}(N,\varphi_0,B^{\mathcal M}_{\rho_0/4}(Z)) 
	\leq\,& C E_{\varsigma/2,\varepsilon}(N,\varphi,B^{\mathcal M}_{\rho_0}(Z_0)) 
	\leq C \beta \rho_0^{-n/2-1} \|N\|_{L^2(B^{\mathcal M}_{\rho_0}(Z_0))} 
	\\ \leq\,& C \beta \rho_0^{\alpha-1} \|\varphi\|_{L^2(B_1(0,P_{Z_0}))} 
	= C \beta \rho_0^{\alpha-1} \|\varphi_0\|_{L^2(B_1(0,P_Z))} 
	\\ \leq\,& C \beta E_{\varsigma/2,\varepsilon}(N,\varphi_0,B^{\mathcal M}_{\rho_0/4}(Z)) 
		+ C \beta \rho_0^{-n/2-1} \|N\|_{L^2(B^{\mathcal M}_{\rho_0/4}(Z))} 
\end{align*}
where $C = C(n,m,q,M_0,N_0,C_e,C_h,\alpha) \in (0,\infty)$ are constants.  Hence provided $\beta$ is small enough that $C\beta < 1/2$ 
\begin{equation}\label{SandT eqn7} 
	E_{\varsigma/2,\varepsilon}(N,\varphi_0,B^{\mathcal M}_{\rho_0/4}(Z)) \leq 2C \beta \rho_0^{-n/2-1} \|N\|_{L^2(B^{\mathcal M}_{\rho_0/4}(Z))} .
\end{equation}
By assumption, conclusion~(A) of Theorem~\ref{main excess decay thm} does not hold true with $Z,\rho_0/4$ in place of $Z_0,\rho_0$.  Hence by \eqref{SandT eqn7}, conclusion~(B) of Theorem~\ref{main excess decay thm} must hold true with $Z,\varphi_0,\rho_0/2$ in place of $Z_0,\varphi,\rho_0$; that is, there exists $j(1) \in \{1,2,\ldots,\lceil q/q_0 \rceil\}$, a rotation $\widetilde{\Gamma}_1$, and $\varphi_1 \in \Phi_{\alpha,q}(P_Z,\widetilde{\Gamma}_1(S))$ such that $\widetilde{\Gamma}_1(P_{Z_0}) = P_Z$, $\widetilde{\Gamma}_1 |_{P_{Z_0}} = \Gamma_Z |_{P_{Z_0}}$, and \eqref{SandT eqn2}--\eqref{SandT eqn4} holds true.  For the inductive step, suppose that for some positive integer $i_0$, for all $i \in \{1,2,\ldots,i_0-1\}$ we have found $j(i) \in \{1,2,\ldots,\lceil q/q_0 \rceil\}$, a rotation $\widetilde{\Gamma}_i$ of $\mathbb{R}^{n+m}$, and $\varphi_i \in \Phi_{\alpha,q}(P_{Z_0},\widetilde{\Gamma}_i(S))$ such that $\widetilde{\Gamma}_i(P_{Z_0}) = P_Z$, $\widetilde{\Gamma}_i |_{P_{Z_0}} = \Gamma_Z |_{P_{Z_0}}$, and \eqref{SandT eqn2}--\eqref{SandT eqn4} holds true.  By iteratively applying \eqref{SandT eqn4} and using \eqref{SandT eqn6} 
\begin{equation}\label{SandT eqn8} 
	E_{\varsigma/2,\varepsilon}(N,\varphi_i,B^{\mathcal M}_{s_i}(Z)) 
			\leq C \Big(\frac{s_i}{\rho_0}\Big)^{\alpha-1+\mu/2} E_{\varsigma/2,\varepsilon}(N,\varphi,B^{\mathcal M}_{\rho_0}(Z_0))
\end{equation}
for all $i \in \{0,1,2,\ldots,i_0-1\}$, where $C = C(n,m,q,M_0,N_0,C_e,C_h,\alpha) \in (0,\infty)$ is a constant.  By \eqref{SandT eqn3} and \eqref{SandT eqn8}
\begin{equation}\label{SandT eqn9} 
	\rho_0^{\alpha-1} \|\mathcal{G}(\varphi_{i-1},\varphi_i)\|_{L^2(B_1(0,P_Z))}
			\leq C \Big(\frac{s_i}{\rho_0}\Big)^{\mu/2} E_{\varsigma/2,\varepsilon}(N,\varphi,B^{\mathcal M}_{\rho_0}(Z_0))
\end{equation}
for all $i \in \{0,1,2,\ldots,i_0-1\}$, where $C = C(n,m,q,M_0,N_0,C_e,C_h,\alpha) \in (0,\infty)$ is a constant.  Noting that $s_l \leq \vartheta_1^{l-i} s_i$ and using the triangle inequality and \eqref{SandT eqn9}  
\begin{align}\label{SandT eqn10} 
	\rho_0^{\alpha-1} \|\mathcal{G}(\varphi_i,\varphi_j)\|_{L^2(B_1(0,P_Z))} 
	\leq\,& \sum_{l=i}^j \rho_0^{\alpha-1} \|\mathcal{G}(\varphi_{l-1},\varphi_l)\|_{L^2(B_1(0,P_Z))} 
	\\ \leq\,& \sum_{l=i}^l C \vartheta_1^{(l-i) \mu/2} \Big(\frac{s_i}{\rho_0}\Big)^{\mu/2} 
		E_{\varsigma/2,\varepsilon}(N,\varphi,B^{\mathcal M}_{\rho_0}(Z_0)) \nonumber 
	\\ \leq\,& C \Big(\frac{s_i}{\rho_0}\Big)^{\mu/2} E_{\varsigma/2,\varepsilon}(N,\varphi,B^{\mathcal M}_{\rho_0}(Z_0)) \nonumber 
\end{align}
for all $0 \leq i < j \leq i_0-1$ and in particular 
\begin{align}\label{SandT eqn11} 
	\rho_0^{\alpha-1} \|\mathcal{G}(\varphi_i,\varphi_0)\|_{L^2(B_1(0,P_Z))} 
	\leq C E_{\varsigma/2,\varepsilon}(N,\varphi,B^{\mathcal M}_{\rho_0}(Z_0)) 
\end{align}
for all $i \in \{0,1,2,\ldots,i_0-1\}$, where $C = C(n,m,q,M_0,N_0,C_e,C_h,\alpha) \in (0,\infty)$ are constants.  Thus by the triangle inequality, Remark~\ref{graphical rmk}(2), Hypothesis~$(\dagger)$, and \eqref{SandT eqn11} 
\begin{equation}\label{SandT eqn12} 
	\frac{1}{C} \rho_0^{-n/2-\alpha} \|N\|_{L^2(B^{\mathcal M}_{\rho_0}(Z_0))} \leq \|\varphi_i\|_{L^2(B_1(0,P_Z))} 
		\leq C \rho_0^{-n/2-\alpha} \|N\|_{L^2(B^{\mathcal M}_{\rho_0}(Z_0))}
\end{equation}
for all $i \in \{0,1,2,\ldots,i_0-1\}$, where $C = C(n,m,q,M_0,N_0,C_e,C_h,\alpha) \in (1,\infty)$ is a constant.  By the triangle inequality, \eqref{SandT eqn8}, Hypothesis~$(\dagger)$, and \eqref{SandT eqn12} 
\begin{align*}
	E_{\varsigma/2,\varepsilon}(N,\varphi_i,B^{\mathcal M}_{s_i}(Z)) 
	\leq\,& C \Big(\frac{s_i}{\rho_0}\Big)^{\alpha-1+\mu/2} E_{\varsigma/2,\varepsilon}(N,\varphi,B^{\mathcal M}_{\rho_0}(Z_0)) 
	\\ \leq\,& C \beta \Big(\frac{s_i}{\rho_0}\Big)^{\alpha-1+\mu/2} \rho_0^{-n/2-1} \|N\|_{L^2(B^{\mathcal M}_{\rho_0}(Z_0))}
	\\ \leq\,& C \beta \Big(\frac{s_i}{\rho_0}\Big)^{\mu/2} s_i^{\alpha-1} \|\varphi_i\|_{L^2(B_1(0,P_Z))} 
	\\ \leq\,& C \beta \Big(\frac{s_i}{\rho_0}\Big)^{\mu/2} \big( s_i^{-n/2-1} \|N\|_{L^2(B^{\mathcal M}_{s_i}(Z))}
		+ C \beta E_{\varsigma/2,\varepsilon}(N,\varphi_i,B^{\mathcal M}_{s_i}(Z)) \big)
\end{align*}
for all $i \in \{0,1,2,\ldots,i_0-1\}$, where $C = C(n,m,q,M_0,N_0,C_e,C_h,\alpha) \in (0,\infty)$ are constants.  Provided $\beta$ is small enough that $C\beta < 1/2$ 
\begin{equation}\label{SandT eqn13} 
	E_{\varsigma/2,\varepsilon}(N,\varphi_i,B^{\mathcal M}_{s_i}(Z)) 
	\leq C \beta \Big(\frac{s_i}{\rho_0}\Big)^{\mu/2} s_i^{-n/2-1} \|N\|_{L^2(B^{\mathcal M}_{s_i}(Z))} 
\end{equation}
for all $i \in \{0,1,2,\ldots,i_0-1\}$, where $C = C(n,m,q,M_0,N_0,C_e,C_h,\alpha) \in (0,\infty)$ is a constant.  By \eqref{SandT eqn2} and \eqref{SandT eqn13} 
\begin{equation}\label{SandT eqn14} 
	\|\widetilde{\Gamma}_{i-1} - \widetilde{\Gamma}_i\| \leq C \beta \Big(\frac{s_i}{\rho_0}\Big)^{\mu/2} 
\end{equation}
for all $i \in \{0,1,2,\ldots,i_0-1\}$, where $C = C(n,m,q,M_0,N_0,C_e,C_h,\alpha) \in (1,\infty)$ are constants.  Noting that $s_i \leq \vartheta_1^i s_0$ and using the triangle inequality and \eqref{SandT eqn14} as above
\begin{equation}\label{SandT eqn15} 
	\|\widetilde{\Gamma}_i - \widetilde{\Gamma}_j\| \leq C \beta \Big(\frac{s_i}{\rho_0}\Big)^{\mu/2} 
\end{equation}
for all $0 \leq i < j \leq i_0-1$.  In particular, $\|\widetilde{\Gamma}_i - \Gamma_Z\| \leq C\beta$ and thus by \eqref{noncon rotate}
\begin{equation}\label{SandT eqn16} 
	\|\widetilde{\Gamma}_i - I\| \leq C \beta + C \mathbf{m}_0^{1/2} |Z - Z_0|
\end{equation}
for all $i \in \{0,1,2,\ldots,i_0-1\}$, where $C = C(n,m,q,M_0,N_0,C_e,C_h,\alpha) \in (0,\infty)$ are constants.   In light of \eqref{SandT eqn13}, we can apply Theorem~\ref{main excess decay thm} with $\varphi_{i_0-1},Z,s_{i_0-1}$ in place of $\varphi,Z_0,\rho_0$.  If conclusion~(A) of Theorem~\ref{main excess decay thm} holds true with $Z,s_{i_0-1}$ in place of $Z_0,\rho_0$, then we stop and set $k_{\star}(Z) = i_0$ as assuming $\beta$ is sufficiently small by \eqref{SandT eqn16} we know that \eqref{SandT eqn5} holds true.  Otherwise, conclusion~(B) of Theorem~\ref{main excess decay thm}  true with $\varphi_{i_0-1},Z,s_{i_0-1}$ in place of $\varphi,Z_0,\rho_0$; that is, there exists $j(i_0) \in \{1,2,\ldots,\lceil q/q_0 \rceil\}$, a rotation $\widetilde{\Gamma}_{i_0}$, and $\varphi_{i_0} \in \Phi_{\alpha,q}(P_Z,\Gamma_{i_0}(S))$ such that \eqref{SandT eqn2}--\eqref{SandT eqn4} holds true with $i = i_0$.  

Having proved the claim, for each positive integer $k$ we let 
\begin{equation*}
	\Upsilon_k = \{ Z \in \mathcal{B}^{\star}_q \cap B^{\mathcal M}_{\rho_0/2}(Z_0) : k_{\star}(Z) = k \} 
\end{equation*}
and we let $\Upsilon_{\infty} = \{ Z \in \mathcal{B}^{\star}_q \cap B^{\mathcal M}_{\rho_0/2}(Z_0) : k_{\star}(Z) = \infty \}$.  Note that $\mathcal{B}^{\star}_q \cap B^{\mathcal M}_{\rho_0}(Z_0) = \bigcup_{1 \leq k < \infty} \Upsilon_k \cup \Upsilon_{\infty}$.  We shall show that the conclusion of the theorem holds true with $\mathcal{S} = \Upsilon_{\infty}$ and $\mathcal{T} = \bigcup_{1 \leq k < \infty} \Upsilon_k$

Let $Z \in \Upsilon_{\infty}$.  By \eqref{SandT eqn15}, $(\widetilde{\Gamma}_i)_{i=0}^{\infty}$ is a Cauchy sequence of linear automorphism of $\mathbb{R}^{n+m}$ and thus there exists a rotation $\widetilde{\Gamma}_Z$ of $\mathbb{R}^{n+m}$ such that $\|\widetilde{\Gamma}_i - \widetilde{\Gamma}_Z\| \rightarrow 0$ as $i\rightarrow\infty$.  Similarly, \eqref{SandT eqn10}, $(\varphi_i)_{i=0}^{\infty}$ is a Cauchy sequence in $L^2(B_1(0,P_Z),\mathcal{A}_q(\mathbb{R}^{n+m}))$ and thus there exists $\varphi^{(Z)} \in \Phi_{\alpha,q}(P_Z,\widetilde{\Gamma}_Z(S))$ such that $\varphi_i \rightarrow \varphi^{(Z)}$ in $L^2(B_1(0,P_Z),\mathcal{A}_q(\mathbb{R}^{n+m}))$.  Letting $i\rightarrow\infty$ in \eqref{SandT eqn16}
\begin{equation*} 
	\|\widetilde{\Gamma}_Z - I\| \leq C \beta + C \mathbf{m}_0^{1/2} |Z - Z_0| , 
\end{equation*}
where $C = C(n,m,q,M_0,N_0,C_e,C_h,\alpha) \in (0,\infty)$ is a constant.   Hence $\op{spine}\varphi^{(Z)} = \widetilde{\Gamma}_Z(S)$ satisfies 
\begin{equation}\label{SandT eqn17} 
	\op{dist}_{\mathcal H}(\op{spine}\varphi^{(Z)} \cap \mathbf{B}_1(0), S \cap \mathbf{B}_1(0)) \leq C \beta + C \mathbf{m}_0^{1/2} |Z - Z_0| , 
\end{equation}
where $C = C(n,m,q,M_0,N_0,C_e,C_h,\alpha) \in (0,\infty)$ is a constant.   Letting $j\rightarrow\infty$ in \eqref{SandT eqn12}
\begin{equation}\label{SandT eqn18} 
	\frac{1}{C} \rho_0^{-n/2-\alpha} \|N\|_{L^2(B^{\mathcal M}_{\rho_0}(Z_0))} \leq \|\varphi^{(Z)}\|_{L^2(B_1(0,P_Z))} 
		\leq C \rho_0^{-n/2-\alpha} \|N\|_{L^2(B^{\mathcal M}_{\rho_0}(Z_0))} , 
\end{equation}
where $C = C(n,m,q,M_0,N_0,C_e,C_h,\alpha) \in (1,\infty)$ is a constant.  Letting $j\rightarrow\infty$ in \eqref{SandT eqn10} 
\begin{align}\label{SandT eqn19} 
	\rho_0^{\alpha-1} \|\mathcal{G}(\varphi_i,\varphi^{(Z)})\|_{L^2(B_1(0,P_Z))} 
	\leq C \Big(\frac{s_i}{\rho_0}\Big)^{\mu/2} E_{\varsigma/2,\varepsilon}(N,\varphi,B^{\mathcal M}_{\rho_0}(Z_0)) 
\end{align}
for all $i \in \{0,1,2,\ldots\}$, where $C = C(n,m,q,M_0,N_0,C_e,C_h,\alpha) \in (0,\infty)$ is a constant.  By the triangle inequality, \eqref{SandT eqn8}, and \eqref{SandT eqn19}  
\begin{equation}\label{SandT eqn20} 
	E_{\varsigma/2,\varepsilon}(N,\varphi^{(Z)},B^{\mathcal M}_{s_i}(Z)) 
		\leq C \Big(\frac{s_i}{\rho_0}\Big)^{\alpha-1+\mu/2} E_{\varsigma/2,\varepsilon}(N,\varphi,B^{\mathcal M}_{\rho_0}(Z_0))
\end{equation}
for all $i \in \{0,1,2,\ldots,i_0-1\}$, where $C = C(n,m,q,M_0,N_0,C_e,C_h,\alpha) \in (0,\infty)$ is a constant.  For each $\rho \in (0,\rho_0/4]$ there exists a (unique) integer $i \geq 0$ such that $s_{i+1} < \rho \leq s_i$ and thus by \eqref{SandT eqn20} 
\begin{equation}\label{SandT eqn21} 
	E_{\varsigma/2,\varepsilon}(N,\varphi^{(Z)},B^{\mathcal M}_{\rho}(Z)) 
		\leq C \Big(\frac{\rho}{\rho_0}\Big)^{\alpha-1+\mu/2} E_{\varsigma/2,\varepsilon}(N,\varphi,B^{\mathcal M}_{\rho_0}(Z_0)) 
\end{equation}
for some constant $C = C(n,m,q,M_0,N_0,C_e,C_h,\alpha) \in (0,\infty)$, proving that $\varphi^{(Z)}$ is the unique tangent function of $T$ relative to the center manifold $\mathcal{M}$ at $Z$ and that \eqref{SandT concl2} holds true.  Arguing as we did to prove \eqref{SandT eqn13} using \eqref{SandT eqn21} and \eqref{SandT eqn18} in place of \eqref{SandT eqn8} and \eqref{SandT eqn12} 
\begin{equation}\label{SandT eqn22} 
	E_{\varsigma/2,\varepsilon}(N,\varphi^{(Z)},B^{\mathcal M}_{\rho}(Z)) \leq C \beta \rho^{-n/2-1} \|N\|_{L^2(B^{\mathcal M}_{\rho}(Z))} 
\end{equation}
for all $\rho \in (0,\rho_0/4]$, where $C = C(n,m,q,M_0,N_0,C_e,C_h,\alpha) \in (0,\infty)$ is a constant.  In light of \eqref{SandT eqn22}, we can apply Theorem~\ref{nonconcentration thm}(A) together with \eqref{exp map est} 
\begin{equation}\label{SandT eqn23} 
	\mathcal{B}^{(\geq 2)}_q \cap B^{\mathcal M}_{\rho/2}(Z) 
	\subset \{ X \in \mathbb{R}^{n+m} : \op{dist}(X, Z + \op{spine} \varphi^{(Z)}) \leq C \beta \rho^{1+\mu/2} + C \mathbf{m}_0^{1/2} \rho^2 \} 
\end{equation}
for all $\rho \in (0,\rho_0/4]$, where $C = C(n,m,q,M_0,N_0,C_e,C_h,\alpha) \in (0,\infty)$ is a constant.

Let $Y,Z \in \Upsilon_{\infty}$ such that $|Y-Z| \leq \rho_0/32$.  In light of \eqref{SandT eqn22}, we can inductively apply Theorem~\ref{main excess decay thm} with $Y,P_Y,\op{spine}\varphi_Y$ in place of $Z_0,P_{Z_0},S$ and starting with $s_0 = 2 |Y-Z|$, $\widetilde{\Gamma}_0 = \Gamma_Z \circ \Gamma_Y^{-1}$, and $\varphi_0 = \Gamma_Z \circ \Gamma_Y^{-1} \circ \varphi_Y \circ \Gamma_Y \circ \Gamma_Z^{-1}$ to construct $j(i) \in \{1,2,\ldots,\lceil q/q_0 \rceil\}$, rotations $\widetilde{\Gamma}_i$, and $\varphi_i \in \Phi_{\alpha,q}(P_{Z_0},\widetilde{\Gamma}_i(S))$ such that \eqref{SandT eqn2}--\eqref{SandT eqn4} holds true for all integers $i \geq 0$.  Since $\varphi^{(Z)}$ is the unique tangent function of $T$ relative to the center manifold $\mathcal{M}$ at $Z$, $\varphi_i \rightarrow \varphi^{(Z)}$ in $L^2(B_1(0,P_Z),\mathcal{A}_q(\mathbb{R}^{n+m}))$ as $i \rightarrow \infty$.  In place of \eqref{SandT eqn17}, we get 
\begin{equation}\label{SandT eqn24} 
	\op{dist}_{\mathcal H}(\op{spine}\varphi^{(Z)} \cap \mathbf{B}_1(0), \op{spine}\varphi_Y \cap \mathbf{B}_1(0)) 
	\leq C \beta \,|Y - Z|^{\mu/2} + C \mathbf{m}_0^{1/2} |Y - Z| 
\end{equation}
where $C = C(n,m,q,M_0,N_0,C_e,C_h,\alpha) \in (0,\infty)$ are constants.  By \eqref{SandT eqn1}, \eqref{SandT eqn23}, and \eqref{SandT eqn24}, there exists a $C^{1,\mu/2}$-function $f : (Z_0 + S) \cap \mathbf{B}_{\rho_0/2}(Z_0) \rightarrow S^{\perp}$ such that $\Upsilon_{\infty} \cap B^{\mathcal M}_{\rho_0/2}(Z_0) \subseteq \op{graph} f$ and the $C^{1,\mu/2}$-norm of $f$ is $\leq C \beta + C \mathbf{m}_0^{1/2}$.

Suppose instead that $Z \in \Upsilon_k$ for some positive integer $k$.  By Theorem~\ref{nonconcentration thm}(A) and \eqref{SandT eqn13} 
\begin{equation*} 
	\mathcal{B}^{(\geq 2)}_q \cap B^{\mathcal M}_{s_i/2}(Z) 
	\subset \{ X \in \mathbb{R}^{n+m} : \op{dist}(X, Z + \Gamma_k(S)) \leq C \beta s_i (s_i/\rho_0)^{\mu/2} + C \mathbf{m}_0^{1/2} s_i^2 \} 
\end{equation*}
for all $i \in \{0,1,2,\ldots,k\}$, where $C = C(n,m,q,M_0,N_0,C_e,C_h,\alpha) \in (0,\infty)$ is a constant.  By \eqref{SandT eqn16}
\begin{equation}\label{SandT eqn25} 
	\mathcal{B}^{(\geq 2)}_q \cap B^{\mathcal M}_{s_i/2}(Z) 
	\subset \{ X \in \mathbb{R}^{n+m} : \op{dist}(X, Z + S) \leq C \beta s_i + C \mathbf{m}_0^{1/2} |Z - Z_0| \,s_i \} 
\end{equation}
for all $i \in \{0,1,2,\ldots,k\}$, where $C = C(n,m,q,M_0,N_0,C_e,C_h,\alpha) \in (0,\infty)$ is a constant.  For each $\rho \in [s_k,\rho_0/4]$ there exists a (unique) integer $i \in \{0,1,2,\ldots,k\}$ such that $s_{i+1} < \rho \leq s_i$ and thus by \eqref{SandT eqn25} 
\begin{equation}\label{SandT eqn26} 
	\mathcal{B}^{(\geq 2)}_q \cap B^{\mathcal M}_{\rho/2}(Z) 
	\subset \{ X \in \mathbb{R}^{n+m} : \op{dist}(X, Z + S) \leq C \beta \rho + C \mathbf{m}_0^{1/2} |Z - Z_0| \,\rho \} , 
\end{equation}
where $C = C(n,m,q,M_0,N_0,C_e,C_h,\alpha) \in (0,\infty)$ is a constant.  Assuming that $\varepsilon$ and $\beta$ are sufficiently small depending only on $n,m,q,M_0,N_0,C_e,C_h,\alpha,\delta$, it follows from \eqref{SandT eqn5} that 
\begin{equation}\label{SandT eqn27} 
	\mathcal{B}^{\star}_q \cap B^{\mathcal M}_{s_k/2}(Z) 
		\cap \mathbf{C}_{\delta s_k/4}(Z+\zeta,P_{Z_0}) = \emptyset \text{ for some } \zeta \in S \cap \mathbf{B}_{s_k/2}(0) .
\end{equation}
By projecting $\mathcal{B}^{\star}_q \cap B^{\mathcal M}_{\rho_0/2}(Z_0)$ onto $P_{Z_0}$ and then arguing exactly as in pages~642--643 of~\cite{Sim93} (using \eqref{SandT eqn26} and \eqref{SandT eqn27} in place of (12) and (13) on page~642 of~\cite{Sim93}), we obtain a covering of $\bigcup_{1 \leq k < \infty} \Upsilon_k$ by cylinders $\{\mathbf{C}_{\rho_j}(x_j,P_{Z_0})\}_{j=1}^{\infty}$ with $x_j \in B_{\rho_0}(Z_0,P_{Z_0})$ and $\sum_{j=1}^{\infty} \rho_j^{n-2} \leq (1-2\eta) \,\rho_0^{n-2}$.  Express $\mathcal{M}$ as the graph of a $C^{3,\kappa}$-function $\bs{\varphi}' : B_{\rho_0}(Z_0,P_{Z_0}) \rightarrow P_{Z_0}^{\perp}$.  Assuming $\varepsilon$ is suitably small, there exists geodesic balls $\{B^{\mathcal M}_{(1+C\mathbf{m}_0^{1/2}) \,\rho_j}(x_j + \bs{\varphi}'(x_j))\}_{j=1}^{\infty}$ covering of $\bigcup_{1 \leq k < \infty} \Upsilon_k$ such that $\sum_{j=1}^{\infty} \rho_j^{n-2} \leq (1-\eta) \,\rho_0^{n-2}$, where $C = C(n,m,q,M_0,N_0,C_e,C_h,\alpha) \in (0,\infty)$ is a constant.
\end{proof}

\begin{theorem}\label{precise rectifiability thm}
Let $q \geq 2$ be an integer.  Let $M_0,N_0,C_e,C_h,\varepsilon_2$ be as in Hypothesis~\ref{noncon constants hypothesis}.  Let $T$, $\mathcal{M}$, and $N$ be as in Hypothesis~\ref{M and N hypothesis} and let $\mathcal{B}^{(\geq 2)}_q$ be as in Definition~\ref{center mnfd branch defn}.  For each $\alpha \geq 2$ let $\mathcal{B}^{(\geq 2)}_{q,\alpha}$ be the set of all branch points $Z \in \mathcal{B}^{(\geq 2)}_q$ such that $\mathbf{I}_{N,Z}(0^+) = \alpha$ and $T$ has a tangent function $\varphi$ relative to $\mathcal{M}$ at $Z$ with $\dim\op{spine}\varphi = n-2$.  Then $\mathcal{B}^{(\geq 2)}_q$ is countably $(n-2)$-rectifiable.  Moreover: 
\begin{enumerate}[itemsep=2mm,topsep=0mm]
	\item[{\rm (a)}]  For any compact set $K \subset \mathbf{C}_{1/2}(0)$, we have that $K \cap \mathcal{B}^{(\geq 2)}_{q,\alpha} \neq \emptyset$ for only finitely many $\alpha \geq 2$. 
	
	\item[{\rm (b)}]  For every $\alpha \geq 2$ and $\mathcal{H}^{n-2}$-a.e.~$Z \in \mathcal{B}^{(\geq 2)}_{q,\alpha}$ there exists a unique tangent function $\varphi^{(Z)} : P_Z \rightarrow \mathcal{A}_q(P_Z^{\perp})$ of $T$ relative to $\mathcal{M}$ at $Z$ such that $\dim\op{spine}\varphi^{(Z)} = n-2$ and 
	\begin{equation}\label{precise rectifiability concl}
		\int_{B^{\mathcal M}_{\rho}(Z)} \mathcal{G}(N,\varphi^{(Z)}\circ\exp_Z^{-1})^2 \leq C_Z \rho^{n+2\alpha+2\mu} 
	\end{equation}
	for all $\rho \in (0,\rho_Z]$, where $\mu = \mu(n,m,q,M_0,N_0,C_e,C_h,\alpha) \in (0,1)$ is a constant and where $\rho_Z \in (0,1/4]$ and $C_Z \in (0,\infty)$ are constants depending on $T$ and $Z$. 
	
	\item[{\rm (c)}]  For every $\alpha \geq 2$ there is a relatively open set $V_{\alpha}$ of $\mathcal{M}$ such that $\mathcal{B}^{(\geq 2)}_{q,\alpha} \subset V_{\alpha}$ and $V_{\alpha} \cap \{X \in \mathcal{B}^{(\geq 2)}_q : \mathbf{I}_{N,X}(0^+) \geq \alpha \}$ has locally finite $\mathcal{H}^{n-2}$-measure, i.e.~for every $Y \in V_{\alpha} \cap \{X \in \mathcal{B}^{(\geq 2)}_q : \mathbf{I}_{N,X}(0^+) \geq \alpha \}$ there exists $\rho > 0$ such that $\mathcal{H}^{n-2}(B^{\mathcal M}_{\rho}(Y) \cap \{X \in \mathcal{B}^{(\geq 2)}_q : \mathbf{I}_{N,X}(0^+) \geq \alpha \}) < \infty$.  In particular, for each $\alpha \geq 2$, $\mathcal{B}^{(\geq 2)}_{q,\alpha}$ has locally finite $\mathcal{H}^{n-2}$-measure.
\end{enumerate}
\end{theorem}

\begin{proof} 
The proof is similar to that of~\cite[Theorem~$2'$]{Sim93}, so we shall sketch the proof for completion.  Let $T$, $\mathcal{M}$, and $N$ be as in Hypothesis~\ref{M and N hypothesis} and let $\mathcal{B}^{(\geq 2)}_q$ be as in Definition~\ref{center mnfd branch defn}.  By Theorem~\ref{centered stratification thm}, the set $S^{(n-3)}$ of all points $Z \in \mathcal{B}^{(\geq 2)}_q(T)$ such that every tangent function of $T$ relativity to $\mathcal{M}$ at $Z$ has spine dimension $\leq n-3$ satisfies $\dim_{\mathcal H} S^{(n-3)} \leq n-3$.  Hence it suffices to consider the set $\mathcal{B}^{(\geq 2)}_{q,\alpha}$ of all points $Z \in \mathcal{B}^{(\geq 2)}_q(T)$ at which at the planar frequency $= \alpha$ and least one tangent function has spine dimension $= n-2$. 

Let $Z_0 \in \mathcal{B}^{(\geq 2)}_{q,\alpha}$ and $\varphi$ be a tangent function of $T$ at $Z_0$ with $\dim \op{spine}\varphi = n-2$.  Let $\eta \in (0,1)$ to be later determined.  By the monotonicity formula for Almgren's frequency function Theorem~\ref{freq mono centered thm} and definition of tangent function, for every $\eta > 0$ there exists $\sigma \in (0,1]$ such that  
\begin{equation*}
	\mathbf{I}_{N,Z_0}(2\sigma/\delta(\eta)) < \alpha + \delta(\eta) ,  
\end{equation*}
where $\delta(\eta)$ is as in Lemma~\ref{lemma2_4}.  For each $\rho_{\star} \in (0,1/2)$, define the outer measure $\mu_{\rho_{\star}}$ on $B^{\mathcal M}_{\rho_{\star}}(Z_0)$ by 
\begin{equation*}
	\mu_{\rho_{\star}}(A) = \inf \sum_{i=1}^N \omega_{n-2} \sigma_i^{n-2} 
\end{equation*}
for every set $A \subseteq B^{\mathcal M}_{\rho_0}(Z_0)$, where the infimum is taken over all finite covers of $A$ by open geodesic balls $\mathbf{B}^{\mathcal M}_{\sigma_i}(Y_i)$ for $i = 1,2,\ldots,N$ with $Y_i \in B^{\mathcal M}_{\rho_0}(Z_0)$ and $\sigma_i \leq \rho_{\star}$.  Let $\mathcal{B}^+_{q,\alpha} = \{ X \in \mathcal{B}^{(\geq 2)}_q \cap \overline{B^{\mathcal M}_{\rho_0}(Z_0)} : \mathbf{I}_{N,X}(0^+) \geq \alpha \}$ and note that $\mathcal{B}^{(\geq 2)}_{q,\alpha} \subseteq \mathcal{B}^+_{q,\alpha}$.  Choose a cover of $\mathcal{B}^+_{q,\alpha}$ by a finite collection of open geodesic balls $B^{\mathcal M}_{\sigma_i}(Y_i)$ such that $Y_i \in B^{\mathcal M}_{\rho_0}(Z_0)$, $\sigma_i \leq \rho_{\star}$, and 
\begin{equation*}
	\sum_{i=1}^N \omega_{n-2} \sigma_i^{n-2} \leq \mu_{\rho_{\star}}(\mathcal{B}^+_{q,\alpha}) + 1
\end{equation*}
Remove the balls $B^{\mathcal M}_{\sigma_i}(Y_i)$ which do not intersect $\mathcal{B}^+_{q,\alpha}$ from the collection.  For each $i$, let $Z_i \in B^{\mathcal M}_{\sigma_i}(Y_i) \cap \mathcal{B}^+_{q,\alpha}$.  By Lemma~\ref{lemma2_4}, either there exists a non-zero, homogeneous degree $\alpha$, locally Dirichlet energy minimizing $q$-valued function $\varphi_i \in W^{1,2}_{\rm loc}(P_{Z_i},\mathcal{A}_q(P_{Z_i}^{\perp}))$ such that $\dim \op{spine}\varphi_i = n-2$ and 
\begin{equation}\label{precise rectifiability eqn1}
	\int_{B^{\mathcal M}_{2\sigma_i}(Z_i)} \mathcal{G}(N,\varphi_i)^2 \leq \eta^2 \int_{B^{\mathcal M}_{2\sigma_i}(Z_i)} |N|^2 
\end{equation}
or there exists an $(n-3)$-dimensional linear subspace $L_i$ of $P_{Z_i}$ such that 
\begin{equation}\label{precise rectifiability eqn2}
	\{ X \in \mathcal{B}^{(\geq 2)}_q \cap \overline{B^{\mathcal M}_{2\sigma_i}(Z_i)} : \mathbf{I}_{N,Z_i}(0^+) \geq \alpha \} 
		\subset \{ X \in \mathbb{R}^{n+m} : \op{dist}(X,Z_i+L_i) < 2\eta\sigma_i \} .
\end{equation}
Note that since $\sigma_i \leq \rho_{\star}$, \eqref{precise rectifiability eqn1} implies that 
\begin{equation}\label{precise rectifiability eqn3}
	E_{\varsigma/2,1}(N,\varphi_i,B^{\mathcal M}_{2\sigma_i}(Z_i)) 
	\leq (\eta + \rho_{\star}^{\varsigma/4}) \,(2\sigma_i)^{-n/2-1} \|N\|_{L^2(B^{\mathcal M}_{2\sigma_i}(Z_i))} . 
\end{equation}
Moreover, the degree of homogeneity $\alpha$ of a non-zero locally Dirichlet energy minimizing $q$-valued function must equal $\ell_0/q_0$ for some relatively prime positive integers $\ell_0,q_0$ with $q_0 \leq q$ and thus the set of all such $\alpha$ is discrete.  Assuming $\rho_{\star}$ is sufficiently small, the conclusion of the theorem now follows by arguing exactly as in~\cite{Sim93} by iteratively applying Theorem~\ref{SandT thm} using the fact that either \eqref{precise rectifiability eqn3} or \eqref{precise rectifiability eqn2} holds true.
\end{proof}

\begin{proof}[Proof of Theorem~\ref{B rectifiable thm} of the Introduction]  Fix constants $M_0,N_0,C_e,C_h,\varepsilon_2$ satisfying the requirements of Hypothesis~\ref{noncon constants hypothesis} (and including the requirements of Hypothesis~\ref{DeLSpa constants hypothesis}).  Let $T$, $\mathcal{M}$, and $N$ be as in Hypothesis~\ref{M and N hypothesis} and let $\mathcal{B}^{(\geq 2)}_q$ be as in Definition~\ref{center mnfd branch defn}.  Let $K$ be a compact subset of $\mathbf{C}_1(0)$ such that $\mathcal{B}^{(\geq 2)}_q \cap K \neq \emptyset$.  By Theorem~\ref{precise rectifiability thm}(a), there is a finite set $\{\alpha_1,\alpha_2,\ldots,\alpha_J\}$ such that $K \cap \mathcal{B}^{(\geq 2)}_q \cap \mathcal{B}^{(\geq 2)}_{q,\alpha_j} \neq \emptyset$ for all $j = 1,2,\ldots,J$.  By Theorem~\ref{precise rectifiability thm}(c), for each $j = 1,2,\ldots,J$ there exists an open set $V_{\alpha_j} \supset \mathcal{B}^{(\geq 2)}_{q,\alpha}$ such that $V_{\alpha_j} \cap \{X \in \mathcal{B}^{(\geq 2)}_q : \mathbf{I}_{N,X}(0^+) = \alpha \}$ has locally finite $\mathcal{H}^{n-2}$-measure in $V_{\alpha_j}$.  Set $\alpha_0 = 2$ and $\alpha_{J+1} = \infty$.  For $j = 0,1,2,\ldots,J$ let 
\begin{equation*}
	\Gamma_j = \{ X \in K \cap \mathcal{B}^{(\geq 2)}_q : \alpha_j \leq \mathbf{I}_{N,X}(0^+) < \alpha_{j+1} \} \cap V_{\alpha_j} 
\end{equation*}
so that $\Gamma_j$ has locally finite measure in $V_{\alpha_j}$ and let 
\begin{equation*}
	\widetilde{\Gamma}_j = \{ X \in K \cap \mathcal{B}^{(\geq 2)}_q : \alpha_j \leq \mathbf{I}_{N,X}(0^+) < \alpha_{j+1} \} \setminus V_{\alpha_j} .
\end{equation*}
Since $\widetilde{\Gamma}_j \subseteq S^{(n-3)}_{q}$, by Theorem~\ref{centered stratification thm} we have that $\dim_{\mathcal H} \widetilde{\Gamma}_j \leq n-3$.  By~\cite[Theorem~1.1 and Corollary~5.5]{KrumWica}, $\mathcal{B}_q^{(\geq 2)}$ is a relatively closed subset of $\mathbf{C}_1(0)$.  By Corollary~\ref{centered semicont1 cor}, $\mathbf{I}_{N,X}(0^+)$ is an upper semi-continuous function of $X \in \mathcal{B}^{(\geq 2)}_q$.  It follows that each of $\Gamma_j$ and $\widetilde{\Gamma}_j$ is the intersection of an open set and closed set, and hence is locally compact.  Clearly $\mathcal{B}^{(\geq 2)}_q \cap K = \bigcup_{j=0}^J (\Gamma_j \cup \widetilde{\Gamma}_j)$.
\end{proof}

\begin{proof}[Proof of Theorem~\ref{unique tangent function thm} of the Introduction] 
Fix constants $M_0,N_0,C_e,C_h,\varepsilon_2$ satisfying the requirements of Hypothesis~\ref{noncon constants hypothesis} (including the requirements of Hypothesis~\ref{DeLSpa constants hypothesis}).  Let $\mathcal{M}$ be the center manifold of $T$ and $N$ be the normal approximation map of $T$ (as in Hypothesis~\ref{M and N hypothesis}).  By Theorem \ref{precise rectifiability thm}(b), for $\mathcal{H}^{n-2}$-a.e.~$Z \in \mathcal{B}^{(\geq 2)}_q$ there exists a 
$\varphi^{(Z)}$ is the unique tangent function of $T$ at $Z$ such that $\dim\op{spine}\varphi^{(Z)} = n-2$ and \eqref{precise rectifiability concl} holds true.  Set $\alpha = \mathbf{I}_{N,Z}(0^+)$ and recall that by Lemma~\ref{planar vs alm freq lemma}, $\alpha \geq \mathcal{N}_{T,{\rm Pl}}(Z) \geq 2$.  Since $\varphi^{(Z)}$ is homogeneous degree $\alpha$ and $\dim\op{spine}\varphi^{(Z)} = n-2$, conclusions~(i) and (iii) holds true with $\mathcal{O}_T(Z) = \alpha$.  By approximating $T$ by the graph of $N$ and applying the area formula together with \eqref{normal approx est1}, \eqref{normal approx est3}, and Corollary~\ref{doubling centered cor}, \eqref{precise rectifiability concl} gives us 
\begin{align}\label{unique tangent function eqn1}
	&\frac{1}{\rho^n} \int_{\mathbf{C}_{\rho}(Z,P_Z)} 
		\op{dist}^2(X, Z + \op{spt} \mathbf{G}_{\varphi^{(Z)} \circ \exp_Z^{-1}}) \,d\|T\|(X) 
	\\&+ \sum_{j=1}^{N_Z} \frac{m^{(Z)}_j}{\rho^n} \int_{B_{\rho/2}(0,P_Z) \cap \{ r < \rho/16\}} 
		\sum_{l=1}^{q^{(Z)}_j} \op{dist}^2(\exp_Z(X) + \varphi^{(Z)}_{j,l}(X), \op{spt} T) \,d\mathcal{H}^n(X) \nonumber 
	\\&\leq C \rho^{2\alpha + 2\mu_1} \nonumber
\end{align}
for all $\rho \in (0,\rho_Z]$, where $r(X) = \op{dist}(X,\op{spine}\varphi^{(Z)})$, $\varphi^{(Z)}_j(X) = \sum_{l=1}^{q^{(Z)}_j} \llbracket \varphi^{(Z)}_{j,l}(X) \rrbracket$ for each $X \in P_Z$, $\mu_1 = \mu_1(n,m,q,M_0,N_0,C_e,C_h,\alpha) \in (0,1)$ is as in Theorem~\ref{precise rectifiability thm}, and $C \in (0,\infty)$ is a constant (depending on $T$ and $Z$ but independent of $\rho$).  Let $\psi^{(Z)} : B_2(0,P_Z) \rightarrow P_Z^{\perp}$ be the $C^{3,\kappa}$-function such that $\mathcal{M} \cap \mathbf{C}_2(Z,P_Z) = Z + \op{spt}\mathbf{G}_{\psi^{(Z)}}$ and note that $\psi^{(Z)}(0) = 0$, $D\psi^{(Z)}(0) = 0$, and $\|\psi^{(Z)}\|_{C^3(B_2(0,P_Z))} \leq C \mathbf{m}_0^{1/2}$, where $C = C(n,m,q,M_0,N_0,C_e,C_h) \in (0,\infty)$ is a constant.  For each $X \in P_Z$, let $\varphi^{(Z)}(X) = \sum_{l=1}^q \llbracket \varphi^{(Z)}_l(X) \rrbracket$.  Observe that if $X,X' \in B_2(0,P_Z)$ such that $\exp_Z(X) = Z+X'+\psi^{(Z)}(X')$, then it follows using Taylor's theorem that $|X-X'| \leq C \mathbf{m}_0^{1/2} |X|^2$ and thus 
\begin{align}\label{unique tangent function eqn2}
	\mathcal{G}(\varphi^{(Z)}(X), \varphi^{(Z)}(X')) 
	\leq\,& C \|\varphi^{(Z)}\|_{L^2(B_1(0,P_Z))} |X|^{\alpha-1} |X-X'| 
	\\ \leq\,& C \mathbf{m}_0^{1/2} \|\varphi^{(Z)}\|_{L^2(B_1(0,P_Z))} |X|^{\alpha+1} , \nonumber
\end{align}
where $C = C(n,m,q,M_0,N_0,C_e,C_h,\alpha) \in (0,\infty)$ are constants.  Hence it follows from \eqref{unique tangent function eqn1} that \eqref{unique tangent function concl1} holds true, proving conclusion~(iv). 

Now we want to show that  conclusion~(v) holds true.  By arguing as in~\cite[Lemma~5.1]{DeLSpa0}, provided $\mathbf{m}_0 \rho_Z$ is sufficiently small there exists a Lipschitz $q$-valued function $u : B_{\rho_Z}(0,P_Z) \rightarrow \mathcal{A}_q(P_Z^{\perp})$ and a closed set $K \subset B_{\rho_Z}(0,P_Z)$ such that 
\begin{equation}\label{unique tangent function eqn3}
    \mathbf{G}_u = (\eta_{Z,1\#} \mathbf{G}_N) \llcorner \mathbf{C}_{\rho_Z}(0,P_Z)
\end{equation}
and by \eqref{normal approx est1}, \eqref{normal approx est3}, and Corollary~\ref{doubling centered cor} 
\begin{gather}\label{unique tangent function eqn4}
	(\eta_{Z,1\#} T) \llcorner \pi_{P_Z}^{-1}(K) = \mathbf{G}_u \llcorner \pi_{P_Z}^{-1}(K) , \\
    \sup_{B_{\rho}(0,P_Z)} |u| \leq C \mathbf{m}_0^{1/(2n)} \rho^{1+\beta_2} , \quad 
    \op{Lip} u \leq C \mathbf{m}_0^{\gamma_2} \rho^{\gamma_2} , \nonumber \\
	\mathcal{H}^n(B_{\rho}(0,P_Z) \setminus K) + \|\eta_{Z,1\#} T - \mathbf{G}_u\|(\mathbf{C}_{\rho}(0,P_Z)) \leq C \mathbf{m}_0^{1+3\gamma_2/4} \rho^{n+2\alpha-2+3\gamma_2/4}  \nonumber 
\end{gather}
for all $\rho \in (0,\rho_Z]$, where $C = C(n,m,q,M_0,N_0,C_e,C_h) \in (0,\infty)$ is a constant.  For each $X \in B_{\rho_Z}(0,P_Z)$, let $u(X) = \sum_{l=1}^q \llbracket u_l(X) \rrbracket$ and let $u_a(X) = \frac{1}{q} \sum_{l=1}^q u_l(X)$ be the single-valued average of $u$ (as in~\cite[Subsection~3.5]{KrumWica}).  By \eqref{precise rectifiability concl} and Corollary~\ref{doubling centered cor}, Hypothesis~$(\star\star)$ holds true with $\Gamma_{Z\#} \eta_{Z,\rho\#} T$ and $\rho^{\alpha-1} (\Gamma_Z \circ \varphi^{(Z)} \circ \Gamma_Z^{-1})$ in place of $T$ and $\varphi$ for all sufficiently small $\rho > 0$, where $\Gamma_Z$ is an orthogonal rotation of $\mathbb{R}^{n+m}$ such that $\Gamma_Z(P_Z) = \mathbb{R}^n \times \{0\}$ and $\Gamma_Z(\op{spine}\varphi^{(Z)}) = \{0\}^2 \times \mathbb{R}^{n-2} \times \{0\}^2$.  Thus by Theorem~\ref{graphical_thm} and Remark~\ref{graphical beta0 rmk} 
\begin{align}\label{unique tangent function eqn10}
    &\op{dist}_{\mathcal H}((\op{spt} T - Z) \cap \mathbf{C}_{\rho}(Z,P_Z) \cap \{r > \tau(\rho)\rho \}, \\&\hspace{1in} \op{spt} \mathbf{G}_{\varphi^{(Z)}} \cap \mathbf{C}_{\rho}(Z,P_Z) \cap \{r > \tau(\rho)\rho \}) \leq C\rho^{\alpha+\mu_1} \nonumber 
\end{align}
for all sufficiently small $\rho > 0$, where $\mu_1$ is as in Theorem~\ref{precise rectifiability thm} (as above) and $C \in (0,\infty)$ is a constant (depending on $T$ and $Z$ but independent of $\rho$) and where $\tau(\rho) = \beta_{\star} \rho^{\frac{\mu_1}{n+2\alpha}}$ for $\beta_{\star} = \beta_{\star}(n,m,q) \in (0,1)$ sufficiently small.  Hence using the area formula, \eqref{unique tangent function eqn4}, and \eqref{unique tangent function concl1}, it follows from \eqref{unique tangent function eqn10} that, after possibly replacing $\rho_Z$ with a smaller positive number, 
\begin{align}\label{unique tangent function eqn11}
    \frac{1}{\rho^n} \int_{B_{\rho}(0,P_Z) \cap \{r > \tau(\rho)\rho\}} \mathcal{G}(u,\psi^{(Z)}+\varphi^{(Z)})^2 \leq C \rho^{2\alpha+\min\{2\mu_1,3\gamma_2/4\}}
\end{align}
for all $\rho \in (0,\rho_Z]$, where $C \in (0,\infty)$ are constants.  By \eqref{unique tangent function concl1} and the homogeneity of $\varphi^{(Z)}$ 
\begin{equation*}
    \frac{1}{\rho^{n+2\alpha}} \int_{\mathbf{C}_{\rho}(Z,P_Z) \cap \{r \leq \tau(\rho)\rho\}} \op{dist}^2(X,Z+\op{spt}\mathbf{G}_{\psi^{(Z)}}) \,d\|T\|(X) \leq C \rho^{2\mu_1} + C \tau(\rho)^{2\alpha+2} 
    \leq C \rho^{\frac{(2\alpha+2)\mu_1}{n+2\alpha}}
\end{equation*}
for all $\rho \in (0,\rho_Z]$, where $C \in (0,\infty)$ is a constant.  Hence using the area formula and \eqref{unique tangent function eqn4} 
\begin{equation*}
    \frac{1}{\rho^{n+2\alpha}} \int_{B_{\rho}(0,P_Z) \cap \{r \leq \tau(\rho)\rho\}} \mathcal{G}(u,q\llbracket\psi^{(Z)}\rrbracket)^2 \leq C \rho^{{\rm min}\left\{\frac{(2\alpha+2)\mu_1}{n+2\alpha},\frac{3\gamma_2}{4}\right\}} 
\end{equation*}
for all $\rho \in (0,\rho_Z]$, where $C \in (0,\infty)$ is a constant.  Again using the homogeneity of $\varphi^{(Z)}$ 
\begin{equation}\label{unique tangent function eqn12}
    \frac{1}{\rho^{n+2\alpha}} \int_{B_{\rho}(0,P_Z) \cap \{r \leq \tau(\rho)\rho\}} \mathcal{G}(u,\psi^{(Z)}+\varphi^{(Z)})^2 \leq C \rho^{{\rm min}\left\{\frac{(2\alpha+2)\mu_1}{n+2\alpha},\frac{3\gamma_2}{4}\right\}} 
\end{equation}
for all $\rho \in (0,\rho_Z]$, where $C \in (0,\infty)$ is a constant.  Combining \eqref{unique tangent function eqn11} and \eqref{unique tangent function eqn12}
\begin{equation}\label{unique tangent function eqn13}
    \frac{1}{\rho^{n+2\alpha}} \int_{B_{\rho}(0,P_Z)} \mathcal{G}(u,\psi^{(Z)}+\varphi^{(Z)})^2 \leq C \rho^{{\rm min}\left\{\frac{(2\alpha+2)\mu_1}{n+2\alpha},\frac{3\gamma_2}{4}\right\}} 
\end{equation}
for all $\rho \in (0,\rho_Z]$, where $C \in (0,\infty)$ is a constant.

If $\mathcal{N}_{T,{\rm Pl}}(Z) = 2$, then by Taylor's theorem there exists a quadratic polynomial $h^{(Z)}$ such that 
\begin{equation}\label{unique tangent function eqn14}
    |\psi^{(Z)}(X) - h^{(Z)}(X)| \leq C |X|^3 
\end{equation}
for all $X \in B_{\rho_Z}(0,P_Z)$ and in particular \eqref{unique tangent function concl2} holds true.  By \eqref{unique tangent function eqn11} we have that $\rho^{-2} u_a(\rho X) \rightarrow h^{(Z)}(X)$ in the $L^2$-topology on compact subsets of $P_Z$ and thus $h^{(Z)}$ is harmonic on $P_Z$, proving conclusion~(v).  If instead $\mathcal{N}_{T,{\rm Pl}}(Z) > 2$, then by~\cite[Theorem~1.6]{KrumWicc} there exists a non-zero, homogeneous degree $\mathcal{N}_{T,{\rm Pl}}(Z)$, locally Dirichlet energy minimizing $q$-valued function $\chi^{(Z)} : P_Z \rightarrow \mathcal{A}_q(P_Z^{\perp})$ such that $\chi^{(Z)}$ is the unique tangent function obtained as the blow-up of $\eta_{Z,\rho\#} T$ relative to $P_Z$ by excess $\widehat{E}(T,P_Z,\mathbf{C}_{\rho}(Z,P_Z))$ as $\rho \rightarrow 0^+$.  Moreover, 
\begin{equation}\label{unique tangent function eqn15}
	\frac{1}{\rho^n} \int_{\mathbf{C}_{\rho}(Z,P_Z)} 
		\op{dist}^2(X, Z + \op{spt} \mathbf{G}_{\chi^{(Z)}}) \,d\|T\|(X) \leq C \rho^{2\mathcal{N}_{T,{\rm Pl}}(Z) + 2\mu_2} 
\end{equation}
for all sufficiently small $\rho > 0$ and there exists $\widetilde{\tau} : (0,\rho_Z] \rightarrow (0,1/2)$ such that 
\begin{align}\label{unique tangent function eqn16}
	&\op{dist}_{\mathcal H}((\op{spt} T - Z) \cap \mathbf{C}_{\rho}(Z,P_Z) \cap \{\widetilde{r} > \widetilde{\tau}(\rho)\rho \}, \\&\hspace{1in} \op{spt}\mathbf{G}_{\chi^{(Z)}} \cap \mathbf{C}_{\rho}(Z,P_Z) \cap \{\widetilde{r} > \widetilde{\tau}(\rho)\rho \}) \leq C \rho^{\mathcal{N}_{T,{\rm Pl}}(Z) + \mu_2} \nonumber 
\end{align}
for all sufficiently small $\rho > 0$, where $\widetilde{r}(X) = \op{dist}(X,\op{spine}\chi^{(Z)})$ and where $\mu_2 = \mu_2(n,m,q,\mathcal{N}_{T,{\rm Pl}}(Z)) \in (0,1)$ and $C \in (0,\infty)$ are constants (independent of $\tau(\rho)$ and $\rho$).  
We observe that~\cite[Theorem~5.6]{KrumWicc} requires the hypotheses of~\cite[Lemma~5.10]{KrumWicc} to hold true with $\eta_{0,1/4\#} T$ and $4^{1-\alpha} \varphi$ in place of $T$ and $\varphi$ and further requires that $\widehat{E}(T,P_0,A_{\rho,1}(\zeta)) \leq C(n) \,\tau^{-(n+2)/2} \varepsilon_0$ to be sufficiently small whenever $\rho \geq \tau/4$ and $\zeta \in \mathbb{R}^{n-2}$ with $\rho^2 + |\zeta|^2 \leq (3+\gamma)^2/16$.  Hence~\cite[Theorem~5.6]{KrumWicc} holds true with $\varepsilon_0 = \varepsilon_{\star} \tau^{(n+2)/2}$ where $\varepsilon_{\star},\beta_0 \in (0,1)$ are sufficiently small constants depending only on $n,m,q,p,\alpha,\gamma$ but independent of $\tau$.  Using this together with~\cite[Theorem~3.2(e)]{KrumWicc} and the use of~\cite[Theorem~5.6]{KrumWicc} to prove the uniqueness of tangent functions in~\cite[Theorem~1.6]{KrumWicc}, we deduce that \eqref{unique tangent function eqn16} holds true with $\widetilde{\tau}(\rho) = \widetilde{\varepsilon}_{\star} \rho^{\frac{2\mathcal{N}_{T,{\rm Pl}}(Z)}{n+2}}$ for some $\widetilde{\varepsilon}_{\star} = \widetilde{\varepsilon}_{\star}(n,m,q,\mathcal{N}_{T,{\rm Pl}}(Z)) \in (0,1)$.  (Note for later purposes in the proof that when $\mathcal{N}_{T,{\rm Pl}}(Z) = 2$, it follows from \eqref{unique tangent function concl1}, \eqref{unique tangent function eqn10} and \eqref{unique tangent function eqn14} that \eqref{unique tangent function eqn15} and \eqref{unique tangent function eqn16} hold true with $\chi^{(Z)} = h^{(Z)} + \varphi^{(Z)}$ when $\alpha = 2$ and $\chi^{(Z)} = q \llbracket h^{(Z)} \rrbracket$ when $\alpha > 2$.  Moreover, by \eqref{unique tangent function eqn15} and~\cite[Corollary~5.3]{KrumWica}, $\chi^{(Z)}$ must be non-zero.  Thus the subsequent discussion applies when $\mathcal{N}_{T,{\rm Pl}}(Z) = 2$.)  We may require that $\rho_Z$ is chosen so that \eqref{unique tangent function concl1}, \eqref{unique tangent function eqn11}, \eqref{unique tangent function eqn15} and \eqref{unique tangent function eqn16} hold true for all $\rho \in (0,\rho_Z]$.  By the area formula, \eqref{unique tangent function eqn4} and \eqref{unique tangent function eqn16}
\begin{equation}\label{unique tangent function eqn17}
    \frac{1}{\rho^n} \int_{B_{\rho}(0,P_Z) \cap \{\widetilde{r} > \widetilde{\tau}(\rho)\rho\}} \mathcal{G}(u,\chi^{(Z)})^2 \leq C \rho^{2\mathcal{N}_{T,{\rm Pl}}(Z) + \min\{2\mu_2,3\gamma_2/4\}} 
\end{equation}
for all $\rho \in (0,\rho_Z]$, where $C \in (0,\infty)$ is a constant. Let $h^{(Z)}$ denote the single-valued average of $\chi^{(Z)}$ and $\chi^{(Z)}_f$ denote the average-free part of $\chi^{(Z)}$ (as in~\cite[Subsection~3.5]{KrumWica}).  Note that $h^{(Z)}$ is a homogeneous degree $\mathcal{N}_{T,{\rm Pl}}(Z)$, harmonic single-valued function and $\chi^{(Z)}_f$ is a homogeneous degree $\mathcal{N}_{T,{\rm Pl}}(Z)$, locally Dirichlet energy minimizing $q$-valued function.  Using the argument leading to \eqref{unique tangent function eqn12} using $h^{(Z)}$ and $\chi^{(Z)}_f$ in place of $\psi^{(Z)}$ and $\varphi^{(Z)}$ respectively, it follows from the area formula, \eqref{unique tangent function eqn4}, \eqref{unique tangent function eqn15}, \eqref{unique tangent function eqn16} and the homogeneity of $\chi^{(Z)}$ that 
\begin{equation}\label{unique tangent function eqn18}
    \frac{1}{\rho^{n+2\mathcal{N}_{T,{\rm Pl}}(Z)}} \int_{B_{\rho}(0,P_Z) \cap \{\widetilde{r} \leq \widetilde{\tau}(\rho)\rho\}} \mathcal{G}(u,\chi^{(Z)})^2 \leq C \rho^{2\mu_3} 
\end{equation}
for all $\rho \in (0,\rho_Z]$, where $\mu_3 = {\rm min}\left\{\mu,\frac{(\mathcal{N}_{T,{\rm Pl}}(Z)+1)\,\mathcal{N}_{T,{\rm Pl}}(Z)}{n+2}, \frac{3\gamma_2}{8}\right\}$ and $C \in (0,\infty)$ is a constant.  Combining \eqref{unique tangent function eqn17} and \eqref{unique tangent function eqn18}
\begin{equation}\label{unique tangent function eqn19}
    \frac{1}{\rho^{n+2\mathcal{N}_{T,{\rm Pl}}(Z)}} \int_{B_{\rho}(0,P_Z)} \mathcal{G}(u,\chi^{(Z)})^2 \leq C \rho^{2\mu_3} 
\end{equation}
for all $\rho \in (0,\rho_Z]$, where $C \in (0,\infty)$ is a constant.  Recall that by Lemma~\ref{planar vs alm freq lemma}, $\mathcal{N}_{T,{\rm Pl}}(Z) \leq \alpha$.  Thus by the triangle inequality, \eqref{unique tangent function eqn13} and \eqref{unique tangent function eqn19}
\begin{equation}\label{unique tangent function eqn20}
	\frac{1}{\rho^n} \int_{B_{\rho}(0,P_Z)} \mathcal{G}(\chi^{(Z)},\psi^{(Z)}+\varphi^{(Z)})^2 \leq C \rho^{2\mathcal{N}_{T,{\rm Pl}}(Z)+2\mu} 
\end{equation}
for all $\rho \in (0,\rho_Z]$, where $\mu = \mu(n,m,q,M_0,N_0,C_e,C_h,\mathcal{N}_{T,{\rm Pl}}(Z),\alpha) \in (0,1)$ and $C \in (0,\infty)$ are constants.  \eqref{unique tangent function eqn20} implies that \eqref{unique tangent function concl2} holds true and 
\begin{equation}
    \label{unique tangent function eqn22} \frac{1}{\rho^n} \int_{B_{\rho}(0,P_Z)} \mathcal{G}(\chi^{(Z)}_f,\varphi^{(Z)})^2 \leq C \rho^{2\mathcal{N}_{T,{\rm Pl}}(Z)+2\mu} 
\end{equation}
for all $\rho \in (0,\rho_Z]$.  This completes the proof of conclusion~(v).

To see conclusion~(vi), first recall that by Lemma~\ref{planar vs alm freq lemma} we already know that $\mathcal{N}_{T,{\rm Pl}}(Z) \leq \alpha$.   Observe that by letting $\rho\rightarrow 0^+$ in \eqref{unique tangent function eqn22}, if $\mathcal{N}_{T,{\rm Pl}}(Z) = \alpha$, then $\chi^{(Z)}_f \equiv \varphi^{(Z)}$ and in particular $\chi^{(Z)}_f$ is non-zero.  If on the other hand $\mathcal{N}_{T,{\rm Pl}}(Z) < \alpha$, then $\chi^{(Z)}_f$ is identically zero so that $\chi^{(Z)} = q \llbracket h^{(Z)} \rrbracket$ on $P_Z$ and $\mathcal{N}_{T,{\rm Pl}}(Z)$ is an integer.  Thus not only is $\mathcal{N}_{T,{\rm Pl}}(Z) < \alpha$, but by conclusion~(i) we have that $\mathcal{N}_{T,{\rm Pl}}(Z) + 1/q \leq \alpha$.  By comparing \eqref{unique tangent function concl2} and \eqref{unique tangent function concl3}, \eqref{unique tangent function concl3} holds true for some $\mu' > 0$ if and only if $\chi^{(Z)} = q \llbracket h^{(Z)} \rrbracket$ on $P_Z$, which occurs precisely when $\mathcal{N}_{T,{\rm Pl}}(Z) < \alpha$.

Finally, to show conclusion~(vii), suppose that there exists $\widetilde{\alpha} > 0$, a $C^3$ single-valued function $\widetilde{\psi}^{(Z)} : P_Z \rightarrow P_Z^{\perp}$, a homogeneous degree $\mathcal{N}_{T,{\rm Pl}}(Z)$, harmonic single-valued function $\widetilde{h}^{(Z)} : P_Z \rightarrow P_Z^{\perp}$ and an average-free, homogeneous degree degree $\widetilde{\alpha}$, locally Dirichlet energy minimizing $q$-valued $\widetilde{\varphi}^{(Z)} : P_Z \rightarrow \mathcal{A}_q(P_Z^{\perp})$ such that \eqref{unique tangent function concl1} and \eqref{unique tangent function concl2} hold true with $\widetilde{\alpha}$, $\widetilde{\psi}^{(Z)}$, $\widetilde{h}^{(Z)}$ and $\widetilde{\varphi}^{(Z)}$ in place of $\mathcal{O}_T(Z)$, $\psi^{(Z)}$, $h^{(Z)}$ and $\varphi^{(Z)}$ respectively.  Then \eqref{unique tangent function eqn13} holds true, and similarly it follows from \eqref{unique tangent function concl1} with $\widetilde{\psi}^{(Z)}, \widetilde{\varphi}^{(Z)},\widetilde{\alpha}$ in place of $\psi^{(Z)},\varphi^{(Z)},\alpha$ that 
\begin{equation}\label{unique tangent function eqn25} 
    \frac{1}{\rho^n} \int_{B_{\rho}(0,P_Z)} \mathcal{G}(u,\widetilde{\psi}^{(Z)}+\widetilde{\varphi}^{(Z)})^2 \leq C \rho^{2\min\{\alpha,\widetilde{\alpha}\}+2\mu} 
\end{equation}
for all sufficiently small $\rho > 0$, where $\mu \in (0,1)$ and $C \in (0,\infty)$ are constants.  By \eqref{unique tangent function eqn13} and \eqref{unique tangent function eqn25} 
\begin{equation*} 
    \frac{1}{\rho^n} \int_{B_{\rho}(0,P_Z)} \mathcal{G}(\psi^{(Z)}+\varphi^{(Z)},\widetilde{\psi}^{(Z)}+\widetilde{\varphi}^{(Z)})^2 \leq C \rho^{2\min\{\alpha,\widetilde{\alpha}\}+2\mu} 
\end{equation*}
for all sufficiently small $\rho > 0$, where $\mu \in (0,1)$ and $C \in (0,\infty)$ are constants.  Hence 
\begin{align}
    \label{unique tangent function eqn26} \frac{1}{\rho^n} \int_{B_{\rho}(0,P_Z)} |\psi^{(Z)} - \widetilde{\psi}^{(Z)}|^2 \leq\,& C \rho^{2\min\{\alpha,\widetilde{\alpha}\}+2\mu} , \\
    \label{unique tangent function eqn27} \frac{1}{\rho^n} \int_{B_{\rho}(0,P_Z)} \mathcal{G}(\varphi^{(Z)},\widetilde{\varphi}^{(Z)})^2 \leq\,& C \rho^{2\min\{\alpha,\widetilde{\alpha}\}+2\mu} 
\end{align}
for all sufficiently small $\rho > 0$, where $\mu \in (0,1)$ and $C \in (0,\infty)$ are constants.  Since $\varphi^{(Z)}$ and $\widetilde{\varphi}^{(Z)}$ are both non-zero and homogeneous, \eqref{unique tangent function eqn27} implies that $\alpha = \widetilde{\alpha}$ and $\varphi^{(Z)} = \widetilde{\varphi}^{(Z)}$ on $P_Z$.  By $\alpha = \widetilde{\alpha}$ and \eqref{unique tangent function eqn26}, \eqref{unique tangent function concl4} holds true.  By \eqref{unique tangent function concl4} and by \eqref{unique tangent function concl2} (with $\psi^{(Z)}, h^{(Z)}$ and with $\widetilde{\psi}^{(Z)},\widetilde{h}^{(Z)}$ in place of $\psi^{(Z)}, h^{(Z)}$), $h^{(Z)} = \widetilde{h}^{(Z)}$ on $P_Z$.
\end{proof}

\begin{proof}[Proof of Theorem~\ref{topology thm} of the Introduction] 
Fix constants $M_0,N_0,C_e,C_h,\varepsilon_2$ satisfying the requirements of Hypothesis~\ref{noncon constants hypothesis} (including the requirements of Hypothesis~\ref{DeLSpa constants hypothesis}).  Applying the blow-up argument from~\cite[Theorem~1.9]{KrumWicc} using \eqref{topology hyp1}, \eqref{topology hyp2}, the coarse blow-up procedure from Subsection~\ref{sec:centered blowup sec}, and Lemma~\ref{centered semicont2 lemma}, it follows that for every $\beta \in (0,1)$ there exists $\varepsilon = \varepsilon(n,m,q,M_0,N_0,C_e,C_h,\alpha,\lambda,\beta) \in (0,1)$ such that if $T, \mathcal{M}, N, Z_0, \rho_0, L$ satisfy the hypothesis of Theorem~\ref{topology thm}, then there exists a non-zero, homogeneous degree $\alpha$, locally Dirichlet energy minimizing function $\varphi \in W^{1,2}(P_{Z_0},\mathcal{A}_q(P_{Z_0}^{\perp})$ such that $\dim\op{spine}\varphi = n-2$, 
\begin{gather}
	\label{topology eqn1} \int_{B^{\mathcal M}_{\lambda\rho_0}(Z_0)} \mathcal{G}(N,\varphi\circ\exp_{Z_0}^{-1})^2 
		< \beta^2 \int_{B^{\mathcal M}_{\lambda\rho_0}(Z_0)} |N|^2 , \\
	\label{topology eqn02} \{ X \in \mathcal{B}_q^{(\geq 2)} \cap \mathbf{I}_{N,X}(0^+) \geq \alpha \} \subseteq \exp_{Z_0}(\{ X' : \op{dist}(X',\op{spine}\varphi) < \beta \} , \\
    \label{topology eqn2} \forall Y \in L \;\; \exists! X \in\op{spine}\varphi \text{ such that } 2\lambda |X| \leq |Y| \text{ and } \pi_L(X) = Y . 
\end{gather}

Let $\vartheta \in (0,1/4)$ such that $\overline{\nu} \vartheta^{\mu/2} < 1/2$, where $\mu$ and $\overline{\nu}$ are as in Lemma~\ref{main excess decay lemma1} (depending on $n,m,q,M_0,N_0,C_e,C_,\alpha$).  Let $Z \in \mathcal{B}_q^{(\geq 2)} \cap \mathbf{C}_{\lambda\rho_0/2}(Z_0)$ and set $S = \op{spine}\varphi$.  In light of \eqref{topology eqn1}, we can iteratively apply Lemma~\ref{main excess decay lemma1} (in place of Theorem~\ref{main excess decay thm}) as in the proof of Theorem~\ref{SandT thm} with $\lambda\rho_0$ in place of $\rho_0$.  In particular, let $\Gamma_Z$ be a rotation of $\mathbb{R}^{n+m}$ such that $\Gamma_Z(P_{Z_0}) = P_Z$ and \eqref{noncon rotate} holds true.  Set $\widetilde{\Gamma}_0 = \Gamma_Z$ and $\varphi_0 = \Gamma_Z \circ \varphi \circ \Gamma_Z^{-1}$ and set $s_k = \vartheta^k \lambda\rho_0/4$  for each integer $i \geq 0$.  We claim that for every $i \in \{1,2,3,\ldots\}$ there exists a rotation $\widetilde{\Gamma}_i$ of $\mathbb{R}^{n+m}$ with $\widetilde{\Gamma}_i(P_{Z_0}) = P_Z$ and $\widetilde{\Gamma}_i |_{P_{Z_0}} = \Gamma_Z |_{P_{Z_0}}$ and $\varphi_i \in \Phi_{\alpha,q}(P_Z,\widetilde{\Gamma}_i(S))$ such that \eqref{SandT eqn2}--\eqref{SandT eqn4} hold true with $\vartheta$ in place of $\vartheta_{j(i)}$.  Having found $\widetilde{\Gamma}_i$ and $\varphi_i$ for each $i \in \{0,1,\ldots,k\}$ for some integer $k \geq 0$, to apply Lemma~\ref{main excess decay lemma1} to find $\widetilde{\Gamma}_{k+1}$ and $\varphi_{K+1}$, we must verify that every $y_0 \in \Gamma_Z(S) \cap \mathbf{B}_{s_k/2}(0)$ 
\begin{equation*}
	\mathcal{B}_q^{\star} \cap B^{\mathcal M}_{\delta s_k/2}(\exp_Z(y_0)) \neq \emptyset . 
\end{equation*}
By \eqref{topology hyp2}, there exists $X \in \mathcal{B}_q^{\star} \cap B^{\mathcal M}_{\rho_0/4}(Z_0)$ such that $\pi_L(X-Z_0) = \pi_L(\exp_Z(y_0)-Z_0)$.  Let $X = \exp_Z(x+y)$ where $x \in P_Z \cap (\Gamma_Z(S))^{\perp}$ and $y \in \Gamma_Z(S)$.  It follows using Theorem~\ref{nonconcentration thm}(A) that 
\begin{equation}\label{topology eqn04}
    |x| \leq C (\beta + \mathbf{m}_0^{1/2} \rho_0) \max\{|y|,s_k\} 
\end{equation}
where $C = C(n,m,q,M_0,N_0,C_e,C_h,\alpha) \in (0,\infty)$ is a constant.  In particular, if $X \not\in B^{\mathcal M}_{\lambda\rho_0/8}(Z)$, then \eqref{topology eqn04} follows from \eqref{topology eqn02} and \eqref{noncon rotate}.  If $X \in B^{\mathcal M}_{s_i/2}(Z) \setminus B^{\mathcal M}_{s_{i+1}/2}(Z)$ for some $i \in \{0,1,2,\ldots,k-1\}$, then by Theorem~\ref{nonconcentration thm}(A) with $Z$, $s_i$, and $s_i^{\alpha-1} \varphi_i$ in place of $Z_0$, $\rho_0$, and $\varphi$ and by \eqref{SandT eqn13}, $|x| \leq C\beta s_i \leq 2\vartheta^{-1} C\beta |X|$ and thus \eqref{topology eqn04} holds true.  Finally, if $X \in B^{\mathcal M}_{s_k/2}(Z) \setminus B^{\mathcal M}_{s_i/2}(Z)$, then similarly by Theorem~\ref{nonconcentration thm}(A) with $Z$, $s_k$, and $s_k^{\alpha-1} \varphi_i$ in place of $Z_0$, $\rho_0$, and $\varphi$ and by \eqref{SandT eqn13}, $|x| \leq C\beta s_k$ as in \eqref{topology eqn04}.  By \eqref{topology eqn04}, \eqref{noncon rotate}, and \eqref{exp map est}
\begin{align*}
	|\pi_L(\Gamma_Z^{-1}(y - y_0))| 
	\leq\,& |\pi_L(x+y - y_0)| + |x| + \|\Gamma_Z - I\| \,|y - y_0| 
	\\ \leq\,& |\pi_L(\exp_Z(x+y) - \exp_Z(y_0))| + C (\beta + \mathbf{m}_0^{1/2} \rho_0) \max\{|y|,s_k\}
	\\ \leq\,& C (\beta + \mathbf{m}_0^{1/2} \rho_0) \max\{|y|,s_k\} , 
\end{align*}
where $C = C(n,m,q,M_0,N_0,C_e,C_h,\alpha) \in (0,\infty)$ are constants.  Hence by \eqref{topology eqn2}
\begin{equation}\label{topology eqn4}
	|y - y_0| \leq 2\lambda C (\beta + \mathbf{m}_0^{1/2} \rho_0) \max\{|y|,s_k\} .
\end{equation}
Note that we must have that $|y| \leq s_k$; otherwise, assuming $\varepsilon$ and $\beta$ are sufficiently small, it follows from \eqref{topology eqn4} that if $|y| > s_k$, then $|y| \leq |y - y_0| + |y_0| \leq |y|/2 + |y_0|$ so that $|y| \leq 2 |y_0| < s_k$, giving us a contradiction.  Hence by $|x| \leq C \beta$ and \eqref{topology eqn4}, provided $\varepsilon$ and $\beta$ are sufficiently small 
\begin{equation*}
	|X -  \exp_Z(y_0)| = |\exp_Z(x+y) - \exp_Z(y_0)| \leq C |x+y - y_0| \leq C (\beta + \mathbf{m}_0^{1/2} \rho_0) \,s_k < \delta s_k/2 , 
\end{equation*}
where $C = C(n,m,q,M_0,N_0,C_e,C_h,\alpha) \in (0,\infty)$ are constants.

Now we have that $\Upsilon_{\infty} = \mathcal{B}^{\star}_q \cap B^{\mathcal M}_{\lambda\rho_0/2}(Z_0)$.  In particular, conclusion~(a) of the theorem holds true as $\Gamma(\eta_{Z_0,1}(\op{sing} T)) \cap \mathbf{C}_{\lambda\rho_0/4}(0) = \{ (g(y), y, h(y)) : y \in B^{n-2}_{\lambda\rho_0/4}(0) \} \cap \mathbf{C}_{\lambda/4}(0)$ for some rotation $\Gamma$ of $\mathbb{R}^{n+m}$ with $\Gamma(P_{Z_0}) = \mathbb{R}^n \times \{0\}$ and some functions $g \in C^{1,\mu/2}(B^{n-2}_{\lambda\rho_0/4}(0),\mathbb{R}^2)$ and $h \in C^{1,\mu/2}(B^{n-2}_{\lambda\rho_0/4}(0),\mathbb{R}^m)$ such that $g(0) = 0$, $Dg(0) = 0$, $h(0) = 0$, $Dh(0) = 0$, and $[Dg]_{\mu/2,B^{n-2}_{\lambda\rho_0/4}(0)} + [Dh]_{\mu/2,B^{n-2}_{\lambda\rho_0/4}(0)} < C\beta$, where $C = C(n,m,q,M_0,N_0,C_e,C_h,\alpha) \in (0,\infty)$ is a constant.  Moreover, $T$ has a unique tangent function $\varphi^{(Z)}$ at each $Z \in \mathcal{B}^{\star}_q \cap B^{\mathcal M}_{\rho_0/2}(Z_0)$.  Let $Z \in \mathcal{B}^{\star}_q \cap B^{\mathcal M}_{\rho_0/2}(Z_0)$ and $S_Z = \op{spine}\varphi^{(Z)}$.  Let $\sigma \in (0,\lambda\rho_0/8]$, $\rho \geq \sigma^{1+\mu/(8n)}/4$, and $\zeta \in S_Z$ with $\rho^2 + |\zeta|^2 \leq (7\sigma/8)^2$.  For each $\kappa \in (0,1]$, let 
\begin{equation*}
	A^Z_{\rho,\kappa}(\zeta) = \{ x+y : x \in P_Z \cap S_Z^{\perp},\, y \in S_Z,\, (|x| - \rho)^2 + |y - \zeta|^2 < (\kappa\rho/64)^2 \} . 
\end{equation*}
By \eqref{SandT eqn21}, for each $\sigma \in (0,\lambda\rho_0/8]$, $\rho \geq \sigma^{1+\mu/(8n)}/4$, and $\zeta \in S$ with $\rho^2 + |\zeta|^2 \leq (7\sigma/8)^2$ 
\begin{align*}
	\rho \,E_{\varsigma/2,\varepsilon}(N,\varphi^{(Z)},\exp_Z(A^Z_{\rho,1}(\zeta)))
	\leq\,& C (\rho/\sigma)^{-n/2} \sigma \,E_{\varsigma/2,\varepsilon}(N,\varphi^{(Z)},B^{\mathcal M}_{\sigma}(Z_0)) 
	\\ \leq\,& C \rho^{-n/2} \sigma^{n/2+\alpha+\mu/2} E_{\varsigma/2,\varepsilon}(N,\varphi,B^{\mathcal M}_{\lambda\rho_0}(Z_0))
	\\ \leq\,& C \sigma^{\alpha+7\mu/16} E_{\varsigma/2,\varepsilon}(N,\varphi,B^{\mathcal M}_{\lambda\rho_0}(Z_0)) , 
\end{align*}
where $C = C(n,m,q,M_0,N_0,C_e,C_h,\alpha) \in (0,\infty)$ are constants.  Thus by Lemma~\ref{graphical_lemma}(B) 
\begin{equation*}
	\sup_{X \in \op{spt} T \cap \mathbf{p}^{-1}(\exp_Z(A^Z_{\rho,3/4}(\zeta)))} \op{dist}(X,\op{spt}\mathbf{G}_{\varphi^{(Z)} \circ \exp_Z^{-1}}) 
	\leq C \sigma^{\alpha+7\mu/16} E_{\varsigma/2,\varepsilon}(N,\varphi,B^{\mathcal M}_{\lambda\rho_0}(Z_0)) , 
\end{equation*}
where $C = C(n,m,q,M_0,N_0,C_e,C_h,\alpha) \in (0,\infty)$ is a constant.  Let $\psi^{(Z)} : B_2(0,P_Z) \rightarrow P_Z^{\perp}$ be a $C^3$-function such that $\mathcal{M} \cap \mathbf{C}_2(0,P_Z) = Z + \op{spt}\mathbf{G}_{\psi^{(Z)}}$, $\psi^{(Z)}(0) = 0$, $D\psi^{(Z)}(0) = 0$, and $\|\psi^{(Z)}\|_{C^3(B_2(0,P_Z))} \leq C \mathbf{m}_0^{1/2}$, where $C = C(n,m,q,M_0,N_0,C_e,C_h) \in (0,\infty)$ is a constant.  By \eqref{exp map est} (see the proof Theorem~\ref{unique tangent function thm}) 
\begin{equation*}
	\sup_{X \in \op{spt} T \cap \pi_{P_Z}^{-1}(A^Z_{\rho,5/8}(\zeta))} \op{dist}(X,\op{spt}\mathbf{G}_{\psi^{(Z)} + \varphi^{(Z)}}) 
	\leq C \sigma^{\alpha+7\mu/16} E_{\varsigma/2,\varepsilon}(N,\varphi,B^{\mathcal M}_{\lambda\rho_0}(Z_0)) , 
\end{equation*}
where $C = C(n,m,q,M_0,N_0,C_e,C_h,\alpha) \in (0,\infty)$ is a constant.  Since the only component of $\varphi^{(Z)}$ has multiplicity one component, by locally applying the Allard regularity theorem (\cite{Allard}), there exists $\widetilde{u}^{(Z,\rho,\zeta)} : A^Z_{\rho,1/2}(\zeta) \rightarrow \mathcal{A}_q(P_Z^{\perp})$ and $\widetilde{v}^{(Z,\rho,\zeta)} : \op{spt} \mathbf{G}_{\varphi^{(Z)}} \cap \pi_{P_Z}^{-1}(A^Z_{\rho,1}(\zeta)) \rightarrow P_Z^{\perp}$ such that for each open ball $B \subset A^Z_{\rho,1/2}(\zeta)$ 
\begin{equation}\label{topology eqn5} 
	\widetilde{u}^{(Z,\rho,\zeta)}(X) = \sum_{l=1}^q \llbracket \psi^{(Z)}(X) + \varphi^{(Z)}_l(X) + \widetilde{v}^{(Z,\rho,\zeta)}_l(X) \rrbracket 
\end{equation}
for each $X \in B$, where $\varphi^{(Z)}(X) = \sum_{l=1}^q \llbracket \varphi^{(Z)}_l(X) \rrbracket$ for single-valued harmonic functions $\varphi^{(Z)}_l : B \rightarrow P_Z^{\perp}$ and $\widetilde{v}^{(Z,\rho,\zeta)}_l : B \rightarrow P_Z^{\perp}$ is a $C^{3,\kappa}$-function given by $\widetilde{v}^{(Z,\rho,\zeta)}_l(X) = \widetilde{v}^{(Z,\rho,\zeta)}(X,\varphi^{(Z)}_l(X))$ for each $X \in B$, and 
\begin{gather}
	\label{topology eqn6} (\eta_{Z,1\#} T) \llcorner \pi_{P_Z}^{-1}(A^Z_{\rho,1/2}(\zeta)) = \mathbf{G}_{\widetilde{u}^{(Z,\rho,\zeta)}} , \\
	\label{topology eqn7} \sup_{X \in A^Z_{\rho,1/2}(\zeta)} \sum_{l=1}^q |\widetilde{v}^{(Z,\rho,\zeta)}_l(X)| 
		\leq C \sigma^{\alpha+7\mu/16} E_{\varsigma/2,\varepsilon}(N,\varphi,B^{\mathcal M}_{\lambda\rho_0}(Z_0)) , 
\end{gather}
where $C = C(n,m,q,M_0,N_0,C_e,C_h,\alpha) \in (0,\infty)$ is a constant.  

Let $B \subset A_{\rho,1/2}(\zeta)$ be an arbitrary ball of radius $\rho/256$.  In $B$, $\widetilde{u}^{(Z,\rho,\zeta)}_l = \psi^{(Z)} + \varphi^{(Z)}_l + \widetilde{v}^{(Z,\rho,\zeta)}_l$ is a solution to the minimal surface system, which in terms of coordinates with respect to an orthonormal basis for $P_Z$ can be expressed as 
\begin{equation}\label{topology eqn8}
	D_i \bigg( \sqrt{G(D\widetilde{u}^{(Z,\rho,\zeta)}_l)} \,G^{ij}(D\widetilde{u}^{(Z,\rho,\zeta)}_l) \,D_j \widetilde{u}^{(Z,\rho,\zeta)}_l \bigg) = 0 \text{ in } B , 
\end{equation}
where $G_{ij}(D\widetilde{u}^{(Z,\rho,\zeta)}_l) = \delta_{ij} + D_i \widetilde{u}^{(Z,\rho,\zeta)}_l \cdot D_j \widetilde{u}^{(Z,\rho,\zeta)}_l$ for all $1 \leq i,j \leq n$, $[G^{ij}(D\widetilde{u}^{(Z,\rho,\zeta)}_l)]$ is the inverse matrix of $[G_{ij}(D\widetilde{u}^{(Z,\rho,\zeta)}_l)]$, and $G(D\widetilde{u}^{(Z,\rho,\zeta)}_l)$ is the determinant of $[G_{ij}(D\widetilde{u}^{(Z,\rho,\zeta)}_l)]$.  By $\mathcal{N}_{T,{\rm Pl}}(Z) \geq 2$, \eqref{fast decay eqn3},~\cite[Corollary~5.3]{KrumWica}, and standard estimates for the minimal surface system 
\begin{equation}\label{topology eqn9}
	\rho^j \sup_{X \in A^Z_{\rho,7/16}(\zeta)} |D^j \widetilde{u}^{(Z,\rho,\zeta)}(X)| \leq C \mathbf{m}_0^{1/2} \sigma^2
\end{equation}
for each $j \in \{0,1,2,3,4\}$, where $C = C(n,m,q,M_0,N_0,C_e,C_h,\alpha) \in (0,\infty)$ is a constant.  By \eqref{topology eqn8} 
\begin{equation}\label{topology eqn10a}
	D_i \big( a^{ij}_{\kappa\lambda}(D\widetilde{u}^{(Z,\rho,\zeta)}_l,\widetilde{u}^{(Z,\rho,\zeta)}_{l'}) \,D_j (\widetilde{u}^{(Z,\rho,\zeta)}_l - \widetilde{u}^{(Z,\rho,\zeta)}_{l'})^{\lambda} \big) = 0 \text{ in } B 
\end{equation}
for all $1 \leq l < l' \leq q$.  Here we let $(\widetilde{u}^{(Z,\rho,\zeta)}_l)^{\lambda}$ denote the $\lambda$-th coordinate of $\widetilde{u}^{(Z,\rho,\zeta)}_l$ with respect to some orthonormal basis for $P_Z^{\perp}$.  We let $a^{ij}_{\kappa\lambda}(P,Q) = \int_0^1 D_{P_j^{\lambda}} A_{\kappa}^i(t \,P + (1-t) \,Q) \,dt$ for each $m \times n$ matrix $P$ and $Q$, where $A^i_{\kappa}(P) = \sqrt{G(P)} \,G^{ij}(P) \,P_i^{\kappa}$ for each $m \times n$ matrix $P = [P_i^{\kappa}]$.  Hence by the Schauder estimates, Remark~\ref{graphical rmk}(2), \eqref{topology eqn1}, \eqref{topology eqn7}, and \eqref{topology eqn9} 
\begin{equation}\label{topology eqn11} 
	\rho^j \sup_{A^Z_{\rho,3/8}(\zeta)} |D^j (\widetilde{u}^{(Z,\rho,\zeta)}_l - \widetilde{u}^{(Z,\rho,\zeta)}_{l'})| \leq C \widehat{E} \sigma^{\alpha}
\end{equation}
for all $1 \leq l < l' \leq q$ and all $j \in \{0,1,2,3\}$, where $\widehat{E} = \rho_0^{-n/2-\alpha} \|N\|_{L^2(B^{\mathcal M}_{\rho_0}(Z_0))}$ and $C = C(n,m,q,M_0,N_0,C_e,C_h,\alpha) \in (0,\infty)$ is a constant.  Since each $\varphi^{(Z)}_l$ is harmonic on $B$ and by \eqref{topology eqn10a}, 
\begin{equation*}
	\Delta (\widetilde{v}^{(Z,\rho,\zeta)}_l - \widetilde{v}^{(Z,\rho,\zeta)}_{l'})^{\kappa} = -D_i \big( (a^{ij}_{\kappa\lambda} - \delta_{ij} \delta_{\kappa\lambda}) \,D_j (\widetilde{u}^{(Z,\rho,\zeta)}_l - \widetilde{u}^{(Z,\rho,\zeta)}_{l'})^{\lambda} \big) \text{ in } B 
\end{equation*}
for all $1 \leq l < l' \leq q$.  By \eqref{topology eqn9} and \eqref{topology eqn11}
\begin{equation*}
	\rho^{2+j} \sup_{A^Z_{\rho,3/8}(\zeta)} |D^j \Delta (\widetilde{v}^{(Z,\rho,\zeta)}_l - \widetilde{v}^{(Z,\rho,\zeta)}_{l'})| \leq C \mathbf{m}_0 \widehat{E} \rho^{-2} \sigma^{4+\alpha} \leq C \mathbf{m}_0 \widehat{E} \sigma^{\alpha+7/4}
\end{equation*}
for all $1 \leq l < l' \leq q$ and all $j \in \{0,1\}$, where $C = C(n,m,q,M_0,N_0,C_e,C_h,\alpha) \in (0,\infty)$ is a constant.  Thus by the Schauder estimates, \eqref{topology eqn9}, and \eqref{topology eqn11}
\begin{align}\label{topology eqn12}
	\sup_{A^Z_{\rho,5/16}(\zeta)} |D^j (\widetilde{v}^{(Z,\rho,\zeta)}_l - \widetilde{v}^{(Z,\rho,\zeta)}_{l'})| \leq\,& C \rho^{-j} \sigma^{\alpha+7\mu/16} E_{\varsigma/2,\varepsilon}(N,\varphi,B^{\mathcal M}_{\lambda\rho_0}(Z_0)) \\ \leq\,& C \sigma^{\alpha+\mu/8-j} E_{\varsigma/2,\varepsilon}(N,\varphi,B^{\mathcal M}_{\lambda\rho_0}(Z_0)) \nonumber
\end{align}
for all $1 \leq l < l' \leq q$ and all $j \in \{0,1,2\}$, where $C = C(n,m,q,M_0,N_0,C_e,C_h,\alpha) \in (0,\infty)$ is a constant.

Let $Z \in \mathcal{B}^{\star}_q \cap B^{\mathcal M}_{\rho_0/2}(Z_0)$ and $\rho \in (0,7\lambda\rho_0/64]$ and set $\rho = 8\sigma/7$ and $\zeta = 0$.  In light of \eqref{topology eqn6}, there exists $u^{(Z,\rho,0)} : A_{\rho,1/8}(0) \rightarrow \mathcal{A}_q(\mathbb{R}^m)$ such that 
\begin{align}\label{topology eqn010}
	&(\Gamma_{\#} \eta_{Z_0,1\#} T) \llcorner ((\pi(\Gamma(Z-Z_0)) + A_{\rho,1/8}(0)) \times \mathbb{R}^m) 
    \\&\hspace{0.65in} = \eta^{-1}_{\Gamma(Z-Z_0)\#} ((\Gamma_{\#} \mathbf{G}_{\widetilde{u}^{(Z,\rho,0)}}) \llcorner (A_{\rho,1/8}(0) \times \mathbb{R}^m))
    = \mathbf{G}_{u^{(Z,\rho,0)}} , \nonumber
\end{align}
where $\pi$ is the orthogonal projection map onto $\mathbb{R}^n \times \{0\}$.  (More generally, if $\sigma \in (0,\lambda\rho_0/8]$, $\rho \geq \sigma^{1+\mu/(8n)}/4$, and $\zeta \in S_Z$ with $\rho^2 + |\zeta|^2 \leq (7\sigma/8)^2$, then we can similarly express $\eta^{-1}_{\Gamma(Z-Z_0)\#} \Gamma_{\#} \mathbf{G}_{\widetilde{u}^{(Z,\rho,\zeta)}}$ as the graph of a $q$-valued function $u^{(Z,\rho,\zeta)}$ over the domain $\pi(\Gamma(Z-Z_0+\zeta)) + A_{\rho,1/8}(0)$.)  By gluing together the graphs of of the functions $u^{(Z,\rho,0)}$, we can construct a $q$-valued function $u : B_{\lambda\rho_0/16}(0) \rightarrow \mathcal{A}_q(\mathbb{R}^m)$ such that 
\begin{equation*}
	(\Gamma_{\#} \eta_{Z_0,1\#} T) \llcorner \mathbf{C}_{\lambda\rho_0/16}(0) = \op{graph} u . 
\end{equation*}

Let $\tilde{\psi} = \Gamma \circ \psi^{(Z)} \circ \Gamma^{-1}$ so that $\tilde{\psi} : B_2(0) \rightarrow \mathbb{R}^m$ is the $C^{3,\kappa}$-function whose graph is $\Gamma(\mathcal{M} - Z_0) \cap \mathbf{C}_2(0)$.  Recall that $\mathcal{B}^{\star}_q \subset \mathcal{M}$ and $T$ is tangent to $\mathcal{M}$ at each point of $\mathcal{B}^{\star}_q$.  It follows from \eqref{topology eqn9} that 
\begin{align}\label{topology eqn10}
    &|u(x,y) - \widetilde{\psi}(g(\zeta_0),\zeta_0) - D\widetilde{\psi}(g(\zeta_0),\zeta_0) \cdot (x-g(\zeta_0),y-\zeta_0)| \\&+ |x - g(\zeta_0)| \,|Du(x,y) - D\widetilde{\psi}(g(\zeta_0),\zeta_0)| + |x - g(\zeta_0)|^2 |D^2 u(x,y)| \leq C \mathbf{m}_0^{1/2} |x-g(\zeta_0)|^2 \nonumber
\end{align}
for all $(x,y) \in B_{\lambda\rho_0/16}(0)$ and $\zeta_0 \in B^{n-2}_{\lambda\rho_0/16}(0)$ with $|y-\zeta_0| < |x-g(\zeta_0)|/128$, where $C = C(n,m,q,M_0,$ $N_0,C_e,C_h,\alpha) \in (0,\infty)$ is a constant.  It follows from \eqref{topology eqn10} that $u \in C^{1,1}(B_2(0),\mathcal{A}_q(\mathbb{R}^m))$.  To see that $u$ is differentiable at $(g(y_0),y_0)$, where $y_0 \in B^{n-2}_{\lambda\rho_0/16}(0)$, let $(x,y) \in B_{\lambda\rho_0/16}(0)$ with $|(x,y) - (y_0,g(y_0))| < \lambda\rho_0/32$.  If $|y-\zeta_0| < |x-g(\zeta_0)|/128$, then by \eqref{topology eqn10} 
\begin{equation*}
    |u(x,y) - \widetilde{\psi}(g(y_0),y_0) - D\widetilde{\psi}(g(y_0),y_0) \cdot (x-g(y_0),y-y_0)| \leq C \mathbf{m}_0^{1/2} |x-g(y_0)|^2 ,
\end{equation*}
where $C = C(n,m,q,M_0,$ $N_0,C_e,C_h,\alpha) \in (0,\infty)$ is a constant.  If instead $|y-\zeta_0| \geq |x-g(\zeta_0)|/128$, then by $\|\widetilde{\psi}\|_{C^{3,\kappa}(B_2(0))} \leq C\mathbf{m}_0^{1/2}$ and \eqref{topology eqn10}
\begin{align*}
    &|u(x,y) - \widetilde{\psi}(g(y_0),y_0) - D\widetilde{\psi}(g(y_0),y_0) \cdot (x-g(y_0),y-y_0)|
    \\ \leq\,& |u(x,y) - \widetilde{\psi}(g(y),y) - D\widetilde{\psi}(g(y),y) \cdot (x-g(y),0)| \\&+ |\widetilde{\psi}(g(y),y) - \widetilde{\psi}(g(y_0),y_0) - D\widetilde{\psi}(g(y_0),y_0) \cdot (g(y)-g(y_0),y-y_0)| \\&+ |D\widetilde{\psi}(g(y),y) - D\widetilde{\psi}(g(y_0),y_0)| \,|x - g(y)|
    \\ \leq\,& C\mathbf{m}_0^{1/2} |x-g(y)|^2 + C\mathbf{m}_0^{1/2} |y-y_0|^2 + C\mathbf{m}_0^{1/2} |y-y_0| \,|x-g(y)|
    \\ \leq\,& C\mathbf{m}_0^{1/2} |(x-g(y_0),y-y_0)|^2 
\end{align*}
where $C = C(n,m,q,M_0,$ $N_0,C_e,C_h,\alpha) \in (0,\infty)$ are constants.  Therefore, $u$ is differentiable at $(g(y_0),y_0)$ with $u(g(y_0),y_0) = q \llbracket \widetilde{\psi}(g(y_0),y_0) \rrbracket$ and $Du(g(y_0),y_0) = q \llbracket D\widetilde{\psi}(g(y_0),y_0) \rrbracket$.  To bound the Lipschitz constant of $Du$, let $(x,y),(x',y') \in B_{1/16}(0)$ with $|(x,y) - (x',y')| < 1/32$.  Assume that $|x' - g(y')| \leq |x - g(y)|$.  If $|(x,y) - (x',y')| < |x-g(y)|/4$, then by \eqref{topology eqn10}
\begin{equation*}
    \mathcal{G}(Du(x,y),Du(x',y')) \leq \sup_{B_{2d}(x,y)} |D^2 u| \,|(x,y) - (x',y')| \leq C\mathbf{m}_0^{1/2} |(x,y) - (x',y')| , 
\end{equation*}
where $C = C(n,m,q,M_0,$ $N_0,C_e,C_h,\alpha) \in (0,\infty)$ is a constant.  If instead $|(x,y) - (x',y')| \geq |x-g(y)|/4$, then by $\|\widetilde{\psi}\|_{C^{3,\kappa}(B_2(0))} \leq C\mathbf{m}_0^{1/2}$ and \eqref{topology eqn10} 
\begin{align*}
    \mathcal{G}(Du(x,y),Du(x',y')) \leq\,& |Du(x,y) - D\widetilde{\psi}(g(y),y)| + |D\widetilde{\psi}(g(y),y) - D\widetilde{\psi}(g(y'),y')| \\&+ |D\widetilde{\psi}(g(y'),y') - Du(x',y')|
    \\ \leq\,& C\mathbf{m}_0^{1/2} |x - g(y)| + C\mathbf{m}_0^{1/2} |y - y'| + C\mathbf{m}_0^{1/2} |x' - g(y')|
    \\ \leq\,& C\mathbf{m}_0^{1/2} |(x,y) - (x',y')| , 
\end{align*}
where $C = C(n,m,q,M_0,$ $N_0,C_e,C_h,\alpha) \in (0,\infty)$ are constants.

By the structure of $\widetilde{u}^{(Z,0,\rho)}$ from \eqref{topology eqn5} (with $\zeta = 0$), there exists $f \in C^{1,\mu/2}(B^2_1(0) \times B^{n-2}_1(0),\mathbb{R}^m)$ such that 
\begin{equation*}
	u(g(y)+x,y) = \sum_{z^q = x} \llbracket f(z,y) \rrbracket 
\end{equation*}
for all $x \in B^2_1(0)$ and $y \in B^{n-2}_1(0)$ so that conclusion~(b) holds true.  By \eqref{topology eqn1}, Remark~\ref{graphical rmk}(2), Theorem~\ref{graphical_thm}(A), and \eqref{topology eqn7}, $f(\nu z,y) \neq f(z,y)$ for every $q^{\rm th}$ root of unity $\nu$, $z \in \mathbf{C}$ with $|z| < 1$, and $y \in B^{n-2}_1(0)$ and thus the mapping $(z,y) \mapsto (g(y) + z^q, y, f(z,y))$ is injective.  Let $\op{Re}(c (x_1+ix_2)^{p/q})$ be the unique tangent function (as a blow-up) of $\Gamma_{\#} \eta_{Z_0,1\#} T$ relative to $\Gamma(\eta_{Z_0,1}(\mathcal{M}))$ at the origin.  Note that in the special case $Z = Z_0$, $u = \Gamma \circ \widetilde{u}^{(0,\rho,\zeta)} \circ \Gamma^{-1}$ on $B_{\lambda\rho_0/16}(0) \cap A_{\rho,1/8}(\Gamma(\zeta))$ for all $\sigma \in (\lambda\rho_0/8]$, $\rho \geq \sigma^{1+\mu/(8n)}/4$, and $\zeta \in S_{Z_0}$ with $\rho^2 + |\zeta|^2 \leq (7\sigma/8)^2$ and $B_{\lambda\rho_0/16}(0) \cap A_{\rho,1/8}(\Gamma(\zeta)) \neq \emptyset$.  (Here we note that $\Gamma(S_{Z_0}) = \{0\}^2 \times \mathbb{R}^{n-2} \times \{0\}^m$ and thereby identify $\Gamma(\zeta)$ with the corresponding point on $\mathbb{R}^{n-2}$ for each $\zeta \in S_{Z_0}$.)  Hence setting $\Omega = \{ (x,y) \in B_{\lambda\rho_0/16}(0) : |x| \geq |(x,y)|^{1+\mu/(8n)} \}$, by \eqref{topology eqn5} (with $Z = 0$, $\rho = |x|$, $\sigma = 2 |(x,y)|$, and $\Gamma(\zeta) = y$) there exists $v : \op{graph}\varphi^{(0)} \cap \Omega \times \mathbb{R}^m \rightarrow \mathcal{A}_q(\mathbb{R}^m)$ such that for each open ball $B \subset \Omega$ 
\begin{equation*} 
	u(x,y) = \sum_{l=1}^q \llbracket \varphi^{(0)}_l(x,y) + v_l(x,y) \rrbracket \nonumber 
\end{equation*}
for each $(x,y) \in B$, where $\varphi^{(0)}(X) = \sum_{l=1}^q \llbracket \varphi^{(0)}_l(x,y) \rrbracket$ for single-valued harmonic functions $\varphi^{(0)}_l : B \rightarrow \mathbb{R}^m$ and $v_l : B \rightarrow \mathbb{R}^m$ is a smooth function given by $v_l(x,y) = v(x,y,\varphi_l(x,y))$ for each $(x,y) \in B$.  By \eqref{topology eqn12} (with $Z = 0$, $\rho = |x|$, $\sigma = 2 |(x,y)|$, and $\Gamma(\zeta) = y$) 
\begin{equation*} 
	\sum_{l=1}^q |D^j v_l(x,y) - D^j v_{l'}(x,y)| \leq C |(x,y)|^{\alpha+\mu/8-j} E_{\varsigma/2,\varepsilon}(T,\varphi,\mathbf{C}_1(0))  , 
\end{equation*}
for all $(x,y) \in \Omega$, $1 \leq l < l' \leq q$, and $j \in \{0,1\}$, where $C = C(n,m,q,\alpha) \in (0,\infty)$ is a constant.  Therefore, applying the chain rule as in~\cite[Theorem~1.9]{KrumWicc}, conclusion~(d) holds true (with $a = (1 - \nu^p) \,c$) whenever $z \in B^2_{(\lambda\rho_0/32)^{1/q}}(0)$, $y \in B^{n-2}_{\lambda\rho_0/32}(0)$ such that setting $x = g(y) + z^q$ we have that $|x| \geq |(x,y)|^{1+\mu/(8n)}$.  Suppose on the other hand that $(x,y) \in B_{\lambda\rho_0/16}(0)$ and $z \in \mathbb{C}$ such that $|x| < |(x,y)|^{1+\mu/(8n)}$ and $x=g(y)+z^q$.  Let $Z = (g(y),y,h(y))$ and $Y \in A^Z_{\rho,1/8}(0)$ such that $(x,y) \in _{\rho/512}(Y,P_Z) \subset A^Z_{\rho,1/8}(0)$.  Recall that in the open ball $B = B_{\rho/512}(Y,P_Z) \subset A^Z_{\rho,1/8}(0)$, $\widetilde{u}^{(Z,\rho,0)} = \sum_{l=1}^q \llbracket \widetilde{u}^{(Z,\rho,0)}_l \rrbracket$ for some smooth single-valued functions $\widetilde{u}^{(Z,\rho,0)}_l : B \rightarrow P_Z^{\perp}$.  For each $l \in \{1,2,\ldots,q\}$, there exists a single-valued function $u^{(Z,\rho,0)}_l : B_{\rho/1024}(\pi(\Gamma(Y))) \rightarrow \mathbb{R}^m$ such that 
\begin{equation*}
    \eta^{-1}_{\Gamma(Z-Z_0)\#} ((\Gamma_{\#} \mathbf{G}_{\widetilde{u}^{(Z,\rho,0)}_l}) \llcorner (B_{\rho/1024}(\pi(\Gamma(Y))) \times \mathbb{R}^m))
    = \mathbf{G}_{u^{(Z,\rho,0)}_l}
\end{equation*}
(cf.~\eqref{topology eqn010}).  By \eqref{topology eqn11} with $Z = (g(y),y,h(y))$, $\rho = |x|$, $\sigma = 8\rho/7$, and $\zeta = 0$ 
\begin{align*} 
	|u_l(x,y) - u_{l'}(x,y)| \leq\,& C \sup_{A^Z_{\rho,3/8}(0)} |\widetilde{u}^{(Z,\rho,0)}_l - \widetilde{u}^{(Z,\rho,0)}_{l'}| 
    \\ \leq\,& C |(x,y)|^{\alpha+\mu/8-j} E_{\varsigma/2,\varepsilon}(T,\varphi,\mathbf{C}_1(0)) , 
\end{align*}
for all $1 \leq l < l' \leq q$, where $C = C(n,m,q,\alpha) \in (0,\infty)$ are constants.  Notice that if $X,X' \in B_{\rho/512}(Y,P_Z)$ such that $(x,y,u_l(x,y)) = \Gamma(Z - Z_0 + X + \widetilde{u}^{(Z,\rho,0)}_l(X))$ and $(x,y,u_{l'}(x,y)) = \Gamma(Z - Z_0 + X' + \widetilde{u}^{(Z,\rho,0)}_{l'}(X'))$, then $|X - X'| \leq |u_l(x,y) - u_{l'}(x,y)|$ and thus by \eqref{topology eqn9} and \eqref{topology eqn11} with $Z = (g(y),y,h(y))$, $\rho = |x|$, $\sigma = 8\rho/7$, and $\zeta = 0$
\begin{align*} 
	|D^j u_l(x,y) - D^j u_{l'}(x,y)| 
    \leq\,& C |D\widetilde{u}^{(Z,\rho,0)}_l(X) - D\widetilde{u}^{(Z,\rho,0)}_{l'}(X')| 
    \\ \leq\,& C |D\widetilde{u}^{(Z,\rho,0)}_l(X) - D\widetilde{u}^{(Z,\rho,0)}_{l'}(X)| + C \op{Lip} \widetilde{u}^{(Z,\rho,0)}_{l'} |X- X'| 
    \\ \leq\,& C |D\widetilde{u}^{(Z,\rho,0)}_l(X) - D\widetilde{u}^{(Z,\rho,0)}_{l'}(X)| + C \op{Lip} \widetilde{u}^{(Z,\rho,0)}_{l'} |u_l(x,y) - u_{l'}(x,y)|
    \\ \leq\,& C |(x,y)|^{\alpha+\mu/8-1} E_{\varsigma/2,\varepsilon}(T,\varphi,\mathbf{C}_1(0))  , 
\end{align*}
for all $1 \leq l < l' \leq q$, where $C = C(n,m,q,\alpha) \in (0,\infty)$ are constants.  Thus again using the chain rule as in~\cite[Theorem~1.9]{KrumWicc}, we complete the proof of conclusion~(d).  Conclusion~(c) follows from conclusion~(d) with $\psi(x,y) = u_a(g(y)+x,y) = \sum_{z^q = x} f(z,y)$ for each $(x,y) \in B_{\lambda\rho_0/32}(0)$, where $u_a$ is the single-valued average of $u$.  (Note that the graph of $\psi$ is not necessarily the center manifold $\mathcal{M}$ of $T$.)
\end{proof}

\appendix
\section{Change of coordinates}\label{sec:appendix sec}

Let $\mathcal{M}$ be a $C^3$ $n$-dimensional submanifold of $\mathbb{R}^{n+m}$.  For each $X \in \mathcal{M}$, we let ${\rm T}_X \mathcal{M}$ denote the tangent plane to $\mathcal{M}$ at $X$ and we let ${\rm T}_X \mathcal{M}^{\perp}$ denote its orthogonal complement in $\mathbb{R}^{n+m}$.  We let $A_{\mathcal M}$ denote the second fundamental form of $\mathcal{M}$ so that $A_{\mathcal{M},X}(\vec\sigma,\vec\tau) = \pi_{{\rm T}_X \mathcal{M}^{\perp}}(D_{\vec\sigma} \vec\tau)$ for each $X \in \mathcal{M}$ and tangent vector fields $\vec\sigma,\vec\tau$ on $\mathcal{M}$.  We let $H_{\mathcal M}$ denote the mean curvature of $\mathcal{M}$.  We let $B^{\mathcal M}_{\rho}(Z)$ denote the geodesic ball in $\mathcal{M}$ with center $Z \in \mathcal{M}$ and radius $\rho > 0$.

Let $\Omega \subseteq \mathcal{M}$ be an open subset.  Let $\varphi \in C^2(\Omega,\mathbb{R}^{n+m})$ be an arbitrary vector field and write $\varphi(X) = \sum_{k=1}^{n+m} \varphi^k(X) \,e_k$ for all $X \in \Omega$, where $\varphi^k \in C^2(\Omega)$ are real-valued functions.  We equip the trivial bundle $\Omega \times \mathbb{R}^{n+m}$ over $\Omega$ with a connection $D$ given by $D_{\vec \tau}\varphi(X) = \sum_{k=1}^{n+m} D_{\vec \tau}\varphi^k(X) \,e_k$ for all $X \in \Omega$ and $\vec\tau \in {\rm T}_X \mathcal{M}$, where $D_{\vec \tau}\varphi^k(X)$ is the derivative of $\varphi^k$ in the direction $\vec\tau$.  We define the Hessian of $\varphi$ as the second-order absolute derivative of $D\varphi$ as a section of ${\rm T}^* \mathcal{M} \times \mathbb{R}^{n+m}$ over $\mathcal{M}$ (where ${\rm T}^* \mathcal{M}$ is the cotangent bundle of $\mathcal{M}$); that is, the Hessian of $\varphi$ is given by 
\begin{equation}\label{app hessian defn}
	({\rm Hess}\, \varphi)_X(\vec\sigma,\vec\tau) = D_{\vec\sigma} (D_{\vec\tau} \varphi)(X) 
		- D_{\pi_{{\rm T}_X\mathcal{M}}(D_{\vec\sigma} \vec\tau)} \varphi(X) 
\end{equation}
for each $X \in \Omega$ and $C^1$ tangent vector fields $\vec\sigma,\vec\tau$ on $\mathcal{M}$.  Notice that ${\rm Hess} \varphi$ is a $C^1(\Omega)$-bilinear map on $C^1$ tangent vector fields and thus for each $X \in \Omega$ we can regard $({\rm Hess} \varphi)_X$ is a bilinear form on ${\rm T}_X \mathcal{M}$.  
We define the Laplacian of $\varphi$ by 
\begin{equation}\label{app laplace defn}
	\Delta_{\mathcal M} \varphi(X) = \sum_{i=1}^n ({\rm Hess}\, \varphi)_X(\vec\tau_i,\vec\tau_i) 
\end{equation}
for each $X \in \Omega$, where $\vec\tau_1,\vec\tau_2,\ldots,\vec\tau_n$ is an orthonormal basis for ${\rm T}_X \mathcal{M}$.

\begin{lemma}\label{app lemma}
There exists $\varepsilon = \varepsilon(n,m) > 0$ such that the following holds true.  Let $Z \in \mathbb{R}^{n+m}$ and $\rho > 0$.  Let $\mathcal{M}$ and $\mathcal{N}$ be $C^3$ $n$-dimensional submanifolds of $\mathbf{C}_{\rho}(Z,{\rm T}_Z\mathcal{M})$ such that 
\begin{gather}\label{app hyp1}
	Z \in \mathcal{M} \cap \mathcal{N}, \quad 
	\op{dist}_{\mathcal H}({\rm T}_Z \mathcal{M} \cap \mathbf{B}_1(0), {\rm T}_Z \mathcal{N} \cap \mathbf{B}_1(0)) < \varepsilon , \\
	\rho \|A_{\mathcal M}\|_{C^0(\mathcal{M})} < \varepsilon, \quad \rho \|A_{\mathcal N}\|_{C^0(\mathcal{N})} < \varepsilon . \nonumber 
\end{gather}
Let $\varphi \in C^2(B^{\mathcal M}_{\rho}(Z),\mathbb{R}^{n+m})$ such that 
\begin{equation*}
	\frac{1}{\rho} \|\varphi\|_{C^0(B^{\mathcal M}_{\rho}(Z))} + \|D\varphi\|_{C^0(B^{\mathcal M}_{\rho}(Z))} < \varepsilon . 
\end{equation*}
Then there exists a unique function $\psi \in C^2(B^{\mathcal N}_{3\rho/4}(Z),{\rm T}\mathcal{N}^{\perp})$ such that 
\begin{equation}\label{app concl1}
	\mathbf{G}_{\psi} = \mathbf{G}_{\varphi} \llcorner \mathbf{C}^{\mathcal N}_{3\rho/4}(Z) .
\end{equation}
Moreover, letting $\xi : B^{\mathcal N}_{3\rho/4}(Z) \rightarrow B^{\mathcal M}_{\rho}(Z)$ be the unique injective $C^2$-function such that 
\begin{equation}\label{app concl2}
	X + \psi(X) = \xi(X) + \varphi(\xi(X)) 
\end{equation}
for each $X \in B^{\mathcal N}_{3\rho/4}(Z)$, we have that 
\begin{align}
	\label{app concl3} \sup_{X \in B^{\mathcal N}_{3\rho/4}(Z)} |D\xi(X) -& I_{{\rm T}_X \mathcal{N}}| \leq C \sup_{X \in B^{\mathcal N}_{3\rho/4}(Z)} 
		\|\pi_{{\rm T}_{\xi(X)} \mathcal{M}} - \pi_{{\rm T}_X \mathcal{N}}\| \\&+ \|A_{\mathcal N}\|_{C^0(B^{\mathcal N}_{3\rho/4}(Z))} 
		\|\psi\|_{C^0(B^{\mathcal N}_{3\rho/4}(Z))} \nonumber \\&+ \sup_{Y \in B^{\mathcal M}_{\rho}(Z),\, \vec\tau \in {\rm T}_Y \mathcal{M} \cap 
		\partial\mathbf{B}_1(0)} |\pi_{{\rm T}_Y \mathcal{M}}(D_{\vec\tau} \varphi(Y))| , \nonumber \\
	\label{app concl4} \|D\psi\|_{C^0(B^{\mathcal N}_{3\rho/4}(Z))} \leq\,& C \sup_{X \in B^{\mathcal N}_{3\rho/4}(Z)} \|\pi_{{\rm T}_{\xi(X)} \mathcal{M}} 
		- \pi_{{\rm T}_X \mathcal{N}}\| \\&+ \|A_{\mathcal N}\|_{C^0(B^{\mathcal N}_{3\rho/4}(Z))} \|\psi\|_{C^0(B^{\mathcal N}_{3\rho/4}(Z))} 
		+ C \|D\varphi\|_{C^0(B^{\mathcal M}_{\rho}(Z))} , \nonumber \\
	\label{app concl5} \|{\rm Hess}\,\psi\|_{C^0(B^{\mathcal N}_{3\rho/4}(Z))} 
		\leq\,& C\sup_{X \in B^{\mathcal N}_{3\rho/4}(Z), \,|\vec v| = |\vec w| = 1}
			 |A_{\mathcal{M},\xi(X)}(D_{\vec v} \xi(X),D_{\vec w} \xi(X)) - A_{\mathcal{N},X}(\vec v,\vec w)| 
		\\&+ C\sup_{X \in B^{\mathcal N}_{3\rho/4}(Z)} \|\pi_{{\rm T}_{\xi(X)} \mathcal{M}} - \pi_{{\rm T}_X \mathcal{N}}\| 
		\,\|A_{\mathcal N}\|_{C^0(B^{\mathcal N}_{3\rho/4}(Z))} \nonumber\\&+ C \|DA_{\mathcal N}\|_{C^0(B^{\mathcal N}_{3\rho/4}(Z))} 
		\|\psi\|_{C^0(B^{\mathcal N}_{3\rho/4}(Z))} \nonumber \\&+ C \|A_{\mathcal N}\|_{C^0(B^{\mathcal N}_{3\rho/4}(Z))} 
		\|D\psi\|_{C^0(B^{\mathcal N}_{3\rho/4}(Z))} + C \|{\rm Hess}\,\varphi\|_{C^0(B^{\mathcal M}_{\rho}(Z))} , \nonumber \\
	\label{app concl6} \|\Delta_{\mathcal N} \psi\|_{C^0(B^{\mathcal N}_{3\rho/4}(Z))} 
		\leq\,& C \sup_{X \in B^{\mathcal N}_{3\rho/4}(Z)} |H_{\mathcal M}(\xi(X)) - H_{\mathcal N}(X)| 
		\\&+ C \sup_{X \in B^{\mathcal N}_{3\rho/4}(Z)} \|\pi_{{\rm T}_{\xi(X)} \mathcal{M}} - \pi_{{\rm T}_X \mathcal{N}}\| 
		\,\|A_{\mathcal N}\|_{C^0(B^{\mathcal N}_{3\rho/4}(Z))} \nonumber
		\\&+ C \|DA_{\mathcal N}\|_{C^0(B^{\mathcal N}_{3\rho/4}(Z))} \|\psi\|_{C^0(B^{\mathcal N}_{3\rho/4}(Z))} \nonumber
		\\&+ C \|A_{\mathcal N}\|_{C^0(B^{\mathcal N}_{3\rho/4}(Z))} \|D\psi\|_{C^0(B^{\mathcal N}_{3\rho/4}(Z))} 
		+ C \|\Delta_{\mathcal M} \varphi\|_{C^0(B^{\mathcal M}_{\rho}(Z))} \nonumber 
		\\&+ C \|A_{\mathcal M}\|_{C^0(B^{\mathcal M}_{\rho}(Z))} \sup_{X \in B^{\mathcal N}_{3\rho/4}(Z)} |D\xi(X) - I_{{\rm T}_X \mathcal{N}}| \nonumber 
		\\&+ C \|{\rm Hess}\,\varphi\|_{C^0(B^{\mathcal M}_{\rho}(Z))} \sup_{X \in B^{\mathcal N}_{3\rho/4}(Z)} |D\xi(X) - I_{{\rm T}_X \mathcal{N}}| , \nonumber
\end{align}
where $I_{{\rm T}_X \mathcal{N}}$ is the identity map on ${\rm T}_X \mathcal{N}$ for each $X \in \mathcal{N}$ and $C = C(n,m) \in (0,\infty)$ is a constant.
\end{lemma}

\begin{proof}
By translating and rescaling we may assume that $Z = 0$ and $\rho = 1$.  By the Implicit Function Theorem, there exists a (unique) function $\psi \in C^2(B^{\mathcal N}_{3/4}(0),{\rm T}\mathcal{N}^{\perp})$ such that $\|\psi\|_{C^1(B^{\mathcal N}_{3/4}(0))} < C(n,m) \,\varepsilon$ and \eqref{app concl1} holds true.  Let $\xi : B^{\mathcal N}_{3/4}(0) \rightarrow B^{\mathcal M}_1(0)$ be the unique injective $C^2$-function such that \eqref{app concl2} holds true. 

Differentiating \eqref{app concl2} 
\begin{equation}\label{app eqn1}
	\vec v + D_{\vec v} \psi(X) = D_{\vec v} \xi(X) + D_{D_{\vec v} \xi(X)} \varphi(\xi(X))
\end{equation}
for all $X \in B^{\mathcal N}_{3/4}(0)$ and $\vec v \in {\rm T}_X \mathcal{N}$.  Using the definition of the second fundamental form $A_{\mathcal N}$ 
\begin{equation}\label{app eqn2}
	\vec u \cdot D_{\vec v} \psi(X) = -\psi(X) \cdot A_{\mathcal{N},X}(\vec u,\vec v)
\end{equation}
for all $X \in B^{\mathcal N}_{3/4}(0)$ and $\vec u,\vec v \in {\rm T}_X \mathcal{N}$.  Notice that $D_{\vec v} \xi(X) \in {\rm T}_{\xi(X)} \mathcal{M}$ and recall that $\|D\psi\|_{C^0(B^{\mathcal N}_{3/4}(0))} < C(n,m) \,\varepsilon$.  Hence by projecting both sides of \eqref{app eqn1} onto ${\rm T}_{\xi(X)} \mathcal{M}$ 
\begin{align*}
	|D_{\vec v} \xi(X) - \vec v| \leq\,& \|\pi_{{\rm T}_{\xi(X)} \mathcal{M}} - \pi_{{\rm T}_X \mathcal{N}}\| \,(|\vec v| + |D_{\vec v} \psi(X)|) 
		+ |\pi_{{\rm T}_X \mathcal{N}}(D_{\vec v} \psi(X))| \\&+ |\pi_{{\rm T}_{\xi(X)} \mathcal{M}}(D_{D_{\vec v} \xi(X)} \varphi(\xi(X)))| 
	\\ \leq\,& C \|\pi_{{\rm T}_{\xi(X)} \mathcal{M}} - \pi_{{\rm T}_X \mathcal{N}}\| + |A_{\mathcal N,X}| \,|\psi(X)| 
		+ |\pi_{{\rm T}_{\xi(X)} \mathcal{M}}(D_{D_{\vec v} \xi(X)} \varphi(\xi(X)))| 
\end{align*}
for all $X \in B^{\mathcal N}_{3/4}(0)$ and $\vec v \in {\rm T}_X \mathcal{N}$ with $|v| = 1$, where $C = C(n,m) \in (0,\infty)$ is a constant.  In particular, \eqref{app concl3} holds true.  Subtracting $\vec v$ from both sides of \eqref{app eqn1} and using \eqref{app concl3} to bound $D_{\vec v} \xi(X) - \vec v$ gives us \eqref{app concl4}.

Differenting \eqref{app eqn1} 
\begin{align}\label{app eqn3}
	A_{\mathcal N,X}(\vec v,\vec w) + (\op{Hess} \psi)_X(\vec v,\vec w) =\,& (\op{Hess} \xi)_X(\vec v,\vec w) 
		+ (\op{Hess} \varphi)_{\xi(X)}(D_{\vec v} \xi(X),D_{\vec w} \xi(X)) 
		\\&+ D_{\pi_{{\rm T}_{\xi(X)}\mathcal{M}}((\op{Hess} \xi)_X(\vec v,\vec w))} \varphi(\xi(x)) \nonumber 
\end{align}
for each $X \in B^{\mathcal N}_{3/4}(0)$ and $\vec v,\vec w \in {\rm T}_X \mathcal{N}$.  Noting that $D_{\pi_{{\rm T}_X \mathcal{N}}(D_{\vec v} \vec w)} \xi(X) \in {\rm T}_{\xi(X)} \mathcal{M}$ and using the chain rule and the definition of the second fundamental form $A_{\mathcal N}$ 
\begin{align}\label{app eqn4}
	\pi_{{\rm T}_{\xi(X)} \mathcal{M}^{\perp}}(({\rm Hess}\,\xi)_X(\vec v,\vec w)) 
	=\,& \pi_{{\rm T}_{\xi(X)} \mathcal{M}^{\perp}}(D_{\vec v} D_{\vec w} \xi(X)) 
	\\=\,& \pi_{{\rm T}_{\xi(X)} \mathcal{M}^{\perp}}(D_{D_{\vec v} \xi(X)} (D_{\vec w} \xi \circ \xi^{-1})(\xi(X))) \nonumber
	\\=\,& A_{\mathcal{M},\xi(X)}(D_{\vec v} \xi(X),D_{\vec w} \xi(X)) \nonumber
\end{align} 
for each $X \in B^{\mathcal N}_{3/4}(0)$ and tangent vector fields $\vec v,\vec w$ on $\mathcal{N}$.  Using \eqref{app eqn2}, the definition of the second fundamental form $A_{\mathcal N}$, and the Codazzi equation 
\begin{align}\label{app eqn5}
	\vec u \cdot ({\rm Hess}\,\psi)(\vec v,\vec w) 
	=\,& \vec u \cdot D_{\vec v} (\pi_{{\rm T}\mathcal{N}}(D_{\vec w} \psi)) + \vec u \cdot D_{\vec v} (\pi_{{\rm T}\mathcal{N}^{\perp}}(D_{\vec w} \psi)) 
	 	- \vec u \cdot D_{\pi_{{\rm T}\mathcal{N}}(D_{\vec v} \vec w)} \psi 
	\\=\,& -\sum_{i=1}^n \vec u \cdot D_{\vec v} (\psi \cdot A_{\mathcal N}(\vec w,\vec\tau_i) \,\vec\tau_i) \nonumber 
		- D_{\vec w} \psi \cdot A_{\mathcal N}(\vec u,\vec v) \\&+ \psi \cdot A_{\mathcal N}(\vec u, \pi_{{\rm T}\mathcal{N}}(D_{\vec v} \vec w)) \nonumber 
	\\=\,& -\psi \cdot D_{\vec u} A_{\mathcal N}(\vec v,\vec w) - D_{\vec v} \psi \cdot A_{\mathcal N}(\vec u,\vec w) 
		- D_{\vec w} \psi \cdot A_{\mathcal N}(\vec u,\vec v) \nonumber 
\end{align}
at each $X \in B^{\mathcal N}_{3/4}(0)$ and for all tangent vector fields on $\vec u,\vec v,\vec w$ on $\mathcal{N}$, where $\vec\tau_1,\vec\tau_2,\ldots,\vec\tau_n$ is an orthonormal tangent frame field on $\mathcal{N}$ with $D_{\vec\tau_i} \vec\tau_j \cdot \vec\tau_k = 0$ at $X$ for all $i,j,k \in \{1,2,\ldots,n\}$.  Assuming $\varepsilon$ is sufficiently small, by projecting both sides of \eqref{app eqn3} onto ${\rm T}_{\xi(X)} \mathcal{M}$ and using $\|D\varphi\|_{C^0(B^{\mathcal M}_1(0))} < \varepsilon$ and \eqref{app eqn5} 
\begin{align}\label{app eqn6} 
	&\tfrac{1}{2} \big| \pi_{{\rm T}_{\xi(X)} \mathcal{M}}(({\rm Hess}\,\xi)_X(\vec v,\vec w)) \big| 
	\\ \leq\,& \|\pi_{{\rm T}_{\xi(X)} \mathcal{M}} - \pi_{{\rm T}_X \mathcal{N}}\| \,\big( |A_{\mathcal{N},X}| + |(\op{Hess} \psi)_X(\vec v,\vec w)| \big) \nonumber
		\\&+ |\pi_{{\rm T}_X \mathcal{N}}((\op{Hess} \psi)_X(\vec v,\vec w)| + |({\rm Hess}\,\varphi)_{\xi(X)}(D_{\vec v} \xi(X),D_{\vec w} \xi(X))| \nonumber 
	\\ \leq\,& \|\pi_{{\rm T}_{\xi(X)} \mathcal{M}} - \pi_{{\rm T}_X \mathcal{N}}\| \,\big( |A_{\mathcal{N},X}| + |(\op{Hess} \psi)_X(\vec v,\vec w)| \big) \nonumber
		\\&+ |\psi(X)| \,|DA_{\mathcal{N},X}| + 2 |D\psi(X)| \,|A_{\mathcal{N},X}| + |({\rm Hess}\,\varphi)_{\xi(X)}(D_{\vec v} \xi(X),D_{\vec w} \xi(X))| 
		\nonumber 
\end{align}
for all $X \in B^{\mathcal N}_{3/4}(0)$ and $\vec v,\vec w \in {\rm T}_X \mathcal{N}$ with $|\vec v| = |\vec w| = 1$.  By \eqref{app eqn4} and \eqref{app eqn6} 
\begin{align}\label{app eqn7} 
	&\big| ({\rm Hess}\,\xi)_X(\vec v,\vec w) - A_{\mathcal{N},X}(\vec v,\vec w) \big| 
	\\ \leq\,& \big| \pi_{{\rm T}_{\xi(X)} \mathcal{M}^{\perp}}(({\rm Hess}\,\xi)_X(\vec v,\vec w)) - A_{\mathcal{N},X}(\vec v,\vec w) \big| 
		+ \big| \pi_{{\rm T}_{\xi(X)} \mathcal{M}}(({\rm Hess}\,\xi)_X(\vec v,\vec w)) \big| \nonumber 
	\\ \leq\,& \big| A_{\mathcal{M},\xi(X)}(D_{\vec v} \xi(X),D_{\vec w} \xi(X)) - A_{\mathcal{N},X}(\vec v,\vec w) \big| \nonumber 
		\\&+ 2 \|\pi_{{\rm T}_{\xi(X)} \mathcal{M}} - \pi_{{\rm T}_X \mathcal{N}}\| \,\big( |A_{\mathcal{N},X}| + |(\op{Hess} \psi)_X(\vec v,\vec w)| \big) 
		\nonumber \\&+ 2 |\psi(X)| \,|DA_{\mathcal{N},X}| + 4 |D\psi(X)| \,|A_{\mathcal{N},X}| 
		+ 2 |({\rm Hess}\,\varphi)_{\xi(X)}(D_{\vec v} \xi(X),D_{\vec w} \xi(X))| \nonumber 
\end{align}
for all $X \in B^{\mathcal N}_{3/4}(0)$ and $v,w \in {\rm T}_X \mathcal{N}$ with $|\vec v| = |\vec w| = 1$.  Recall that $\|D\varphi\|_{C^0(B^{\mathcal M}_1(0))} < \varepsilon$ and observe that by \eqref{app hyp1}, $\|\pi_{{\rm T}_{\xi(X)} \mathcal{M}} - \pi_{{\rm T}_X \mathcal{N}}\| < C(n,m) \,\varepsilon$.  Hence assuming $\varepsilon$ is sufficiently small, we can subtract $A_{\mathcal{N},X}(\vec v,\vec w)$ from both sides of \eqref{app eqn3} and use \eqref{app eqn6} and \eqref{app eqn7} to obtain 
\begin{align}\label{app eqn8} 
	|({\rm Hess}\,\psi)_X(\vec v,\vec w)| 
	\leq\,& C |A_{\mathcal{M},\xi(X)}(D_{\vec v} \xi(X),D_{\vec w} \xi(X)) - A_{\mathcal{N},X}(\vec v,\vec w)| 
		\\&+ C \|\pi_{{\rm T}_{\xi(X)} \mathcal{M}} - \pi_{{\rm T}_X \mathcal{N}}\| \,|A_{\mathcal{N},X}| 
		+ C |\psi(X)| \,|DA_{\mathcal{N},X}| \nonumber\\&+ C |D\psi(X)| \,|A_{\mathcal{N},X}| 
		+ C |({\rm Hess}\,\varphi)_{\xi(X)}(D_{\vec v} \xi(X), D_{\vec w} \xi(X))| \nonumber 
\end{align}
for all $X \in B^{\mathcal N}_{3/4}(0)$ and $\vec v,\vec w \in {\rm T}_X \mathcal{N}$, where $C = C(n,m) \in (0,\infty)$ is a constant.  In particular, it follows using \eqref{app concl3} that \eqref{app concl5} holds true.

Let $\vec\tau_1,\vec\tau_2,\ldots,\vec\tau_n$ is an orthonormal tangent frame field on $\mathcal{N}$.  By the Gram-Schmidt process, there is an orthonormal tangent frame field $\vec\sigma_1,\vec\sigma_2,\ldots,\vec\sigma_n$ on $\mathcal{M}$ such that 
\begin{equation}\label{app eqn9} 
	|D_{\vec\tau_i} \xi(X) - \vec\sigma_i(X)| \leq C |D\xi(X) - I_{{\rm T}_X \mathcal{N}}| 
\end{equation}
for all $X B^{\mathcal N}_{3/4}(0)$ and $i \in \{1,2,\ldots,n\}$, where $C = C(n,m) \in (0,\infty)$ is a constant.  (Note that $|D\xi(X) - I_{{\rm T}_X \mathcal{N}}|$ can be bounded using \eqref{app concl3}.)  By \eqref{app eqn8} with $\vec v = \vec w = \vec\tau_i(X)$ and \eqref{app eqn9}, \eqref{app concl6} holds true.  
\end{proof}

\begin{lemma}\label{app lemma2}
There exists $\varepsilon = \varepsilon(n,m) > 0$ such that the following holds true.  Let $Z \in \mathbb{R}^{n+m}$ and $\rho > 0$.  Let $\mathcal{M}$ and $\mathcal{N}$ be $C^3$ $n$-dimensional submanifolds of $\mathbf{C}_{\rho}(Z,{\rm T}_Z\mathcal{M})$ such that \eqref{app hyp1} holds true.  Let $\varphi_1,\varphi_2 \in C^2(B^{\mathcal M}_{\rho}(Z),\mathbb{R}^{n+m})$ such that for some constant $\kappa \in (0,\infty)$
\begin{gather}
    \label{app2 hyp1} \frac{1}{\rho} \|\varphi_i\|_{C^0(B^{\mathcal M}_{\rho}(Z))} + \|D\varphi_i\|_{C^0(B^{\mathcal M}_{\rho}(Z))} < \varepsilon \text{ for all } i = 1,2 , \\
    \label{app2 hyp2} |\pi_{{\rm T}_X \mathcal{M}}(\varphi_1(X) - \varphi_2(X))|
    \leq \tfrac{1}{16} |\varphi_1(X) - \varphi_2(X)| \text{ for all } X \in B^{\mathcal M}_{\rho}(Z) , \\
    \label{app2 hyp3} \rho \sup_{B^{\mathcal M}_{\rho}(Z)} |D\varphi_1 - D\varphi_2| \leq \kappa \inf_{B^{\mathcal M}_{\rho}(Z)} |\varphi_1 - \varphi_2| .
\end{gather}
Let $\psi_1,\psi_2 \in C^2(B^{\mathcal N}_{3\rho/4}(Z),{\rm T}\mathcal{N}^{\perp})$ such that $\mathbf{G}_{\psi_i} = \mathbf{G}_{\varphi_i} \llcorner \mathbf{C}^{\mathcal N}_{3\rho/4}(Z)$ for all $i = 1,2$.  Then 
\begin{equation}\label{app2 concl}
    \rho \sup_{B^{\mathcal N}_{3\rho/4}(Z)} |D\psi_1 - D\psi_2| \leq 40\kappa \inf_{B^{\mathcal N}_{3\rho/4}(Z)} |\psi_1 - \psi_2| . 
\end{equation}
\end{lemma}

\begin{proof}
Assume that $Z = 0$ and $\rho = 1$.  For $i = 1,2$, let $\xi_i : B^{\mathcal N}_{3/4}(0) \rightarrow B^{\mathcal M}_1(0)$ be the unique injective $C^2$-function such that 
\begin{equation}\label{app2 eqn1}
	X + \psi_i(X) = \xi_i(X) + \varphi_i(\xi_i(X)) 
\end{equation}
for all $X \in B^{\mathcal N}_{3/4}(0)$ and $i = 1,2$.  By subtracting \eqref{app2 eqn1} for $i = 1,2$ 
\begin{equation}\label{app2 eqn2}
	\psi_1(X) - \psi_2(X) = \xi_1(X) - \xi_2(X) + \varphi_1(\xi_1(X)) - \varphi_2(\xi_2(X)) 
\end{equation}
for all $X \in B^{\mathcal N}_{3/4}(0)$.  Recall that $\psi_i(X),\psi_2(X) \in {\rm T}_X \mathcal{N}^{\perp}$ and $\|D\varphi_2\|_{C^0(B^{\mathcal M}_1(0))} < \varepsilon$.  Note that by \eqref{app hyp1}, $\|\pi_{{\rm T}_{\xi_1(X)} \mathcal{M}} - \pi_{{\rm T}_X \mathcal{N}}\| \leq C(n,m) \,\varepsilon$.  Hence assuming $\varepsilon$ is sufficiently small, by projecting both sides of \eqref{app2 eqn2} onto ${\rm T}_{\xi_1(X)} \mathcal{M}$ and using \eqref{app2 hyp2}
\begin{align}\label{app2 eqn3}
    \tfrac{1}{4} |\xi_1(X) - \xi_2(X)| 
    \leq\,& \tfrac{1}{2} |\pi_{{\rm T}_{\xi_1(X)} \mathcal{M}}(\xi_1(X) - \xi_2(X))| 
    \\ \leq\,& \|\pi_{{\rm T}_{\xi_1(X)} \mathcal{M}} - \pi_{{\rm T}_X \mathcal{N}}\| \,|\psi_1(X) - \psi_2(X)| \nonumber\\&+ |\pi_{{\rm T}_{\xi_1(X)} \mathcal{M}}(\varphi_1(\xi_1(X)) - \varphi_2(\xi_2(X)))| \nonumber 
    \\ \leq\,& C \varepsilon \,|\psi_1(X) - \psi_2(X)| + \tfrac{1}{16} |\varphi_1(\xi_1(X)) - \varphi_2(\xi_2(X))| \nonumber 
\end{align}
for all $X \in B^{\mathcal N}_{3/4}(0)$, where $C = C(n,m) \in (0,\infty)$ is a constsant.  Hence assuming $\varepsilon$ is sufficiently small, by \eqref{app2 eqn2}
\begin{align*}
    |\varphi_1(\xi_1(X)) - \varphi_2(\xi_1(X))| 
    \leq\,& |\psi_1(X) - \psi_2(X)| + 2 |\xi_1(X) - \xi_2(X)| 
    \\ \leq\,& (1+4C\varepsilon) \,|\psi_1(X) - \psi_2(X)| + \tfrac{1}{2} |\varphi_1(\xi_1(X)) - \varphi_2(\xi_2(X))|
\end{align*}
or 
\begin{equation}\label{app2 eqn4}
    |\varphi_1(\xi_1(X)) - \varphi_2(\xi_1(X))| 
    \leq 2 (1+4C\varepsilon) \,|\psi_1(X) - \psi_2(X)|  
\end{equation}
for all $X \in B^{\mathcal N}_{3/4}(0)$.

Differentiating \eqref{app2 eqn2}
\begin{equation}\label{app2 eqn5}
	D_{\vec v} \psi_1(X) - D_{\vec v} \psi_2(X) = D_{\vec v} \xi_1(X) - D_{\vec v} \xi_2(X) + D_{D_{\vec v} \xi_1(X)} \varphi_1(\xi_1(X)) - D_{D_{\vec v} \xi_2(X)} \varphi_2(\xi_2(X)) 
\end{equation}
for all $X \in B^{\mathcal N}_{3/4}(0)$ and all $\vec v \in {\rm T}_X \mathcal{N}$.  Recall that \eqref{app eqn2} holds true with $\psi_i$ in place of $\psi$ for all $i = 1,2$ and that $\|D\varphi_i\|_{C^0(B^{\mathcal M}_1(0))} < \varepsilon$ for all $i = 1,2$.  Note that by \eqref{app concl3}, $|D_{\vec v} \xi_i(X)| \leq 1 + C(n,m)\,\varepsilon$ for all $i = 1,2$.  Hence assuming by projecting \eqref{app2 eqn5} onto ${\rm T}_{\xi_1(X)} \mathcal{M}$ 
\begin{align}\label{app2 eqn6}
    &\tfrac{1}{2} |\pi_{{\rm T}_{\xi_1(X)} \mathcal{M}}(D_{\vec v} \xi_1(X) - D_{\vec v} \xi_2(X))| \leq |A_{\mathcal N,X}| \,|\psi_1(X) - \psi_2(X)| \\&\hspace{15mm} + 2 |D\varphi_1(\xi_1(X)) - D\varphi_2(\xi_1(X))| + 2 \varepsilon  |\xi_1(X) - \xi_2(X)| \nonumber 
\end{align}
for all $X \in B^{\mathcal N}_{3/4}(0)$ and all $\vec v \in {\rm T}_X \mathcal{N}$.  Note that by \eqref{app hyp1}, $\|\pi_{{\rm T}_{\xi_1(X)} \mathcal{M}} - \pi_{{\rm T}_{\xi_2(X)} \mathcal{M}}\| \leq C(n,m) \,\varepsilon \,|\xi_1(X) - \xi_2(X)|$.  Since $D_{\vec v} \xi_i(X) \in {\rm T}_{\xi_i(X)} \mathcal{M}$ for all $i = 1,2$, using \eqref{app2 eqn6}
\begin{align}\label{app2 eqn7}
    |D_{\vec v} \xi_1(X) - D_{\vec v} \xi_2(X)| 
    \leq\,& |\pi_{{\rm T}_{\xi_1(X)} \mathcal{M}}(D_{\vec v} \xi_1(X) - D_{\vec v} \xi_2(X))| \\&+ \|\pi_{{\rm T}_{\xi_1(X)} \mathcal{M}} - \pi_{{\rm T}_{\xi_2(X)} \mathcal{M}}\| \,|D_{\vec v} \xi_2(X)| \nonumber 
    \\ \leq\,& |A_{\mathcal N,X}| \,|\psi_1(X) - \psi_2(X)| + 2 |D\varphi_1(\xi_1(X)) - D\varphi_2(\xi_1(X))| \nonumber \\&+ C \varepsilon  |\xi_1(X) - \xi_2(X)| \nonumber
\end{align}
for all $X \in B^{\mathcal N}_{3/4}(0)$ and all $\vec v \in {\rm T}_X \mathcal{N}$, where $C = C(n,m) \in (0,\infty)$ is a constant.  Therefore, by bounding the right-hand side of \eqref{app2 eqn5} using \eqref{app hyp1}, $\|D\varphi_2\|_{C^0(B^{\mathcal M}_1(0))} < \varepsilon$, and \eqref{app2 eqn7} 
\begin{align*}
	|D_{\vec v} \psi_1(X) - D_{\vec v} \psi_2(X)| \leq\,& 
    2\varepsilon |\psi_1(X) - \psi_2(X)| + 5 |D\varphi_1(\xi_1(X)) - D\varphi_2(\xi_1(X))| \\&+ C \varepsilon  |\xi_1(X) - \xi_2(X)| \nonumber 
\end{align*}
for all $X \in B^{\mathcal N}_{3/4}(0)$ and all $\vec v \in {\rm T}_X \mathcal{N}$, where $C = C(n,m) \in (0,\infty)$ is a constant.  By \eqref{app2 eqn3} and \eqref{app2 eqn4}
\begin{equation}\label{app2 eqn8}
	|D_{\vec v} \psi_1(X) - D_{\vec v} \psi_2(X)| \leq
    C\varepsilon |\psi_1(X) - \psi_2(X)| + 5 |D\varphi_1(\xi_1(X)) - D\varphi_2(\xi_1(X))|  
\end{equation}
for all $X \in B^{\mathcal N}_{3/4}(0)$ and all $\vec v \in {\rm T}_X \mathcal{N}$, where $C = C(n,m) \in (0,\infty)$ is a constant.  By combining \eqref{app2 hyp3}, \eqref{app2 eqn4}, and \eqref{app2 eqn8}
\begin{equation}\label{app2 eqn9}
    \sup_{B^{\mathcal N}_{3/4}(0)} |D\psi_1 - D\psi_2| \leq C\varepsilon \sup_{B^{\mathcal N}_{3/4}(0)} |\psi_1 - \psi_2| + 10 (1+C\varepsilon) \,\kappa \inf_{B^{\mathcal N}_{3/4}(0)} |\psi_1 - \psi_2| , 
\end{equation}
where $C = C(n,m) \in (0,\infty)$ is a constant.  Assuming $\varepsilon$ is sufficiently small and using the interpolation inequality 
\begin{equation*}
    \sup_{B^{\mathcal N}_{3/4}(0)} |\psi_1 - \psi_2| 
    \leq \inf_{B^{\mathcal N}_{3/4}(0)} |\psi_1 - \psi_2| + C \sup_{B^{\mathcal N}_{3/4}(0)} |D\psi_1 - D\psi_2| , 
\end{equation*}
where $C = C(n,m) \in (0,\infty)$ is a constant, \eqref{app2 concl} follows from \eqref{app2 eqn9}.
\end{proof}

 \bigskip
\hskip-.2in\vbox{\hsize3in\obeylines\parskip -1pt 
  \small 
Brian Krummel
School of Mathematics \& Statistics 
University of Melbourne
Parkville,VIC  3010, Australia
\vspace{4pt}
{\tt brian.krummel@unimelb.edu.au}} 
\vbox{\hsize3in
\obeylines 
\parskip-1pt 
\small 
Neshan Wickramasekera
DPMMS 
University of Cambridge 
Cambridge CB3 0WB, United Kingdom
\vspace{4pt}
{\tt N.Wickramasekera@dpmms.cam.ac.uk}
}

\end{document}